\documentclass[final]{amsart} 
\usepackage[utf8]{inputenc}
\usepackage[english]{babel}
\usepackage{amssymb,amsthm}
\usepackage{mathtools}
\usepackage[final,kerning=true,spacing=true]{microtype}
\usepackage[inline,shortlabels]{enumitem}
\usepackage{mathrsfs}
\usepackage[bottom=1in]{geometry}
\usepackage{tikz-cd}
\usepackage{calc}
\usepackage{stmaryrd}
\usepackage{aliascnt}
\usepackage{mathbbol} 
\usepackage{appendix}
\usepackage[final,pagebackref,hypertexnames=false]{hyperref}
\usepackage{csquotes}
\usepackage[capitalize]{cleveref}
\usepackage{autonum}
\usepackage[inline]{showlabels}
\AtBeginDocument{%
\let\oldref\ref
\renewcommand{\ref}[1]{\mbox{\oldref{#1}}}}

\renewcommand*{\backref}[1]{}
\renewcommand*{\backrefalt}[4]{%
\ifnum#1=1 %
\ifnum#3=1 %
Cited on p. %
\else
Cited on p. %
\fi
\else
Cited on pp. %
\fi
#2,\par
}

\DeclareSymbolFontAlphabet{\mathbb}{AMSb}
\DeclareSymbolFontAlphabet{\mathbbl}{bbold}
\newcommand{\bbDelta}{\mathbbl{\Delta}}
\newcommand{\bbSigma}{\mathbbl{\Sigma}}
\newcommand{\bbLambda}{\mathbbl{\Lambda}}

\newcommand{\NN}{\mathbb{N}}
\newcommand{\ZZ}{\mathbb{Z}}
\newcommand{\QQ}{\mathbb{Q}}
\newcommand{\RR}{\mathbb{R}}
\newcommand{\CC}{\mathbb{C}}
\newcommand{\FF}{\mathbb{F}}

\newcommand{\LM}{\mathsf{LM}}
\newcommand{\Ass}{\mathsf{Ass}}
\newcommand{\Comm}{\mathsf{Comm}}

\newcommand{\K}{\mathrm{K}}
\newcommand{\D}{\mathrm{D}}
\newcommand{\R}{\mathrm{R}}

\renewcommand{\L}{\mathrm{L}}
\renewcommand{\H}{\mathrm{H}}

\newcommand{\co}{\mathrm{co}}
\newcommand{\op}{\mathrm{op}}

\newcommand{\lax}{\mathrm{lax}}
\newcommand{\rev}{{\mathrm{rev}}}
\newcommand{\const}{\mathrm{const}}
\newcommand{\cocart}{\mathrm{cocart}}
\newcommand{\cart}{\mathrm{cart}}
\newcommand{\rfib}{\mathrm{rfib}}
\newcommand{\lfib}{\mathrm{lfib}}
\newcommand{\inert}{\mathrm{int}}
\newcommand{\ordinary}{\mathrm{ord}}

\newcommand{\gsindex}{0}  

\newcommand{\act}{\mathrm{act}}

\newcommand{\dbl}{\mathrm{dbl}} 

\newcommand{\homotopy}{{\mathsf{h}}} 

\newcommand{\one}{\mathbf{1}}

\newcommand{\blank}{{\mathchoice%
{\raisebox{-1.5pt}{$-$}} 
{\raisebox{-1.5pt}{$-$}} 
{\scalebox{.7}{\raisebox{-1.5pt}{$-$}} } 
{\scalebox{.5}{\raisebox{-1.5pt}{$-$}} } 
}}

\newcommand{\loccit}{\textit{loc.\,cit.}}
\newcommand{\opcit}{\textit{op.\,cit.}}

\renewcommand{\setminus}{\smallsetminus}
\renewcommand{\emptyset}{\varnothing}

\newcommand{\Fil}{\mathcal{F}\kern-1.5pt{\mathit{il}} } 

\newcommand{\Yo}{\mathsf{Y}} 
\newcommand{\Horn}{\Lambda}
\newcommand{\Hecke}{\mathcal{H}} 
\newcommand{\Inv}{\mathrm{Inv}}

\newcommand{\too}{\longrightarrow}
\newcommand{\To}{\Rightarrow}

\newcommand{\ol}[1]{\overline{#1}}
\newcommand{\ul}[1]{\underline{#1}}

\newcommand{\set}[2]{\left\{#1\,\middle|\,#2\right\}}
\newcommand{\wt}[1]{\widetilde{#1}}
\newcommand{\wh}[1]{\widehat{#1}}
\newcommand{\cat}[1]{\mathcal{#1}}
\newcommand{\shv}[1]{\mathcal{#1}}
\newcommand{\lie}[1]{\mathfrak{#1}}
\newcommand{\gs}[1]{\wt{\cat{#1}}} 

\newcommand{\Ab}{\mathsf{Ab}}
\newcommand{\Ani}{\mathsf{Ani}}
\newcommand{\Cat}{\mathsf{Cat}}
\newcommand{\TwoCat}{{\mathsf{Cat}_2}}
\newcommand{\TwoFun}{{\Fun_2}}
\newcommand{\CSegAni}{\mathsf{CSegAni}}
\newcommand{\qCat}{\mathsf{qCat}}
\newcommand{\Set}{\mathsf{Set}}
\newcommand{\sSet}{\mathsf{sSet}}

\newcommand{\ProFin}{\mathsf{ProFin}}
\newcommand{\Kan}{\mathsf{Kan}}
\newcommand{\Top}{\mathsf{Top}}
\newcommand{\Fin}{\mathsf{Fin}}

\newcommand{\GeomSetup}{\mathsf{GeomSetup}}
\newcommand{\Grp}{\mathsf{Grp}}
\newcommand{\Sp}{\mathsf{Sp}}
\newcommand{\CHaus}{\mathsf{CHaus}} 

\newcommand{\longtwoheadrightarrow}%
{\longrightarrow\hspace{-1.2em}\rightarrow\hspace{.2em}}
\renewcommand{\twoheadrightarrow}%
{\rightarrow\hspace{-1.2em}\rightarrow\hspace{.2em}}
\newcommand{\longhookrightarrow}%
{\lhook\joinrel\relbar\joinrel\rightarrow}

\newcommand{\comp}{\mathbin{\circ}}
\newcommand{\isect}{\mathbin{\cap}}

\newcommand{\union}{\mathbin{\cup}}
\newcommand{\bigunion}{\bigcup}
\newcommand{\dunion}{\mathbin{\sqcup}}
\newcommand{\bigdunion}{\bigsqcup}
\newcommand{\tensor}{\otimes}
\newcommand{\bigtensor}{\bigotimes}
\newcommand{\dsum}{\oplus}
\newcommand{\isom}{\cong}
\newcommand{\injto}{\mathrel{\hookrightarrow}}
\newcommand{\injfrom}{\mathrel{\hookleftarrow}}
\newcommand{\surjto}{\mathrel{\twoheadrightarrow}}

\newcommand{\from}{\mathbin{\leftarrow}}
\newcommand{\xto}[1]{\mathbin{\xrightarrow{#1}}} 
\newcommand{\xlongto}[1]{\mathbin{\xrightarrow{\ \ #1\ \ }}} 
\newcommand{\xfrom}[1]{\mathbin{\xleftarrow{#1}}} 
\newcommand{\xinjto}[1]{\mathbin{\xhookrightarrow{#1}}} 
\newcommand{\xinjfrom}[1]{\mathbin{\xhookleftarrow{#1}}} 
\newcommand{\isoto}{\xto{\smash{\raisebox{-.4ex}{\ensuremath{\scriptstyle\sim}}}}}

\newcommand{\xtofrom}[1]{\mathbin{\xleftrightarrow{#1}}}
\newcommand{\isotofrom}{\xtofrom\sim}

\newcommand{\epito}{\mathbin{\to\hspace{-.8em}\to}}

\DeclareMathOperator*{\colim}{colim}
\DeclareMathOperator*{\dotimes}{\otimes}

\DeclareMathOperator{\Alg}{Alg}

\DeclareMathOperator{\CAlg}{CAlg}
\DeclareMathOperator{\cd}{cd} 
\DeclareMathOperator{\Ch}{Ch} 

\DeclareMathOperator{\CMon}{CMon} 
\DeclareMathOperator{\CoCorr}{CoCorr} 
\DeclareMathOperator{\Cond}{Cond} 
\DeclareMathOperator{\Cont}{\mathcal{C}}
\DeclareMathOperator{\coev}{coev}

\DeclareMathOperator{\conj}{conj}
\DeclareMathOperator{\Corr}{Corr}

\DeclareMathOperator{\DSuave}{SD} 
\DeclareMathOperator{\DPrim}{PD} 
\DeclareMathOperator{\Desc}{Desc} 

\DeclareMathOperator{\End}{End}
\DeclareMathOperator{\Enr}{Enr} 
\DeclareMathOperator{\ev}{ev}
\DeclareMathOperator{\Ext}{Ext}
\DeclareMathOperator{\fib}{fib}
\DeclareMathOperator{\Fun}{Fun}
\DeclareMathOperator{\GL}{GL}
\DeclareMathOperator{\Gpd}{Gpd}

\DeclareMathOperator{\Hom}{Hom}
\DeclareMathOperator{\Hyp}{Hyp}
\DeclareMathOperator{\HypShv}{HypShv}
\DeclareMathOperator{\id}{id}
\DeclareMathOperator{\iHom}{\ul{Hom}}

\DeclareMathOperator{\Ind}{Ind}
\DeclareMathOperator{\cInd}{c-Ind}

\DeclareMathOperator{\Inf}{Inf}

\DeclareMathOperator{\Ker}{Ker}
\DeclareMathOperator{\LAct}{LAct}

\DeclareMathOperator{\LMod}{LMod}
\DeclareMathOperator{\Mod}{Mod}
\DeclareMathOperator{\Mon}{Mon} 

\DeclareMathOperator{\Mul}{Mul}

\DeclareMathOperator{\Nerve}{N}
\DeclareMathOperator{\RezkNerve}{N^r}

\DeclareMathOperator{\Path}{Path}
\DeclareMathOperator{\pr}{pr}
\DeclareMathOperator{\Pres}{Pr} 
\DeclareMathOperator{\Prim}{Prim}
\newcommand{\PrL}[1][]{\operatorname{Pr}^L_{#1}}
\newcommand{\PrR}[1][]{\operatorname{Pr}^R_{#1}}
\DeclareMathOperator{\Pro}{Pro}
\DeclareMathOperator{\QCoh}{QCoh}
\DeclareMathOperator{\Rep}{Rep}
\DeclareMathOperator{\hatRep}{\wh{Rep}}

\DeclareMathOperator{\Res}{Res}
\DeclareMathOperator{\RHom}{RHom}
\DeclareMathOperator{\RInd}{RInd}
\DeclareMathOperator{\Shv}{Shv} 
\DeclareMathOperator{\Sing}{Sing}

\DeclareMathOperator{\Suave}{Suave}

\DeclareMathOperator{\Supp}{Supp}

\DeclareMathOperator{\Str}{Str}
\DeclareMathOperator{\Tot}{Tot}
\DeclareMathOperator{\Un}{Un}

\DeclareMathOperator{\Op}{Op} 
\DeclareMathOperator{\GOp}{GOp} 
\DeclareMathOperator{\Quiv}{Quiv} 
\DeclareMathOperator{\PreEnr}{PreEnr} 
\newcommand{\PSh}[1][]{\cat P_{#1}} 
\DeclareMathOperator{\Tw}{Tw}
\newcommand{\enrHom}[1][]{\operatorname{Hom}^{#1}}
\newcommand{\enrFun}[1][]{\operatorname{Fun}^{#1}}

\newcommand{\Einfty}{\mathbb E_\infty}
\newcommand{\fine}[1]{\textit{$#1$-fine}}
\newcommand{\ufine}{\textit{fine}}
\newcommand{\et}{\mathrm{et}}
\newcommand{\nuc}{\mathrm{nuc}}
\newcommand{\solid}{{\scalebox{0.5}{$\square$}}}

\newcommand{\noloc}{\nobreak\mskip6mu plus1mu\mathpunct{}\nonscript\mkern-\thinmuskip{:}\mskip2mu\relax} 

\newcommand{\CorrHomInline}[5]{
    #1 \xfrom{#2} #3 \xto{#4} #5
}

\newcommand{\CorrHomDisplay}[5]{
    \begin{tikzcd}[column sep=tiny,row sep=small,ampersand replacement=\&]
        \& #3 \ar[dl,"#2"'] \ar[dr,"#4"] \\
        #1 \& \& #5
    \end{tikzcd}
}

\newcommand{\CorrHom}[5]{\mathchoice
    {\CorrHomDisplay{#1}{#2}{#3}{#4}{#5}}
    {\CorrHomInline{#1}{#2}{#3}{#4}{#5}}
    {\CorrHomInline{#1}{#2}{#3}{#4}{#5}}
    {\CorrHomInline{#1}{#2}{#3}{#4}{#5}}
}

\newcommand{\abs}[1]{|#1|}

\swapnumbers
\newaliascnt{main}{subsubsection}
\newtheorem{proposition}[main]{Proposition}
\newtheorem{theorem}[main]{Theorem}
\newtheorem{lemma}[main]{Lemma}
\newtheorem{corollary}[main]{Corollary}

\theoremstyle{definition}
\newtheorem{definition}[main]{Definition}
\newtheorem*{claim*}{Claim}
\newtheorem*{fact*}{Fact}
\newtheorem{notation}[main]{Notation}
\newtheorem{construction}[main]{Construction}
\newtheorem{convention}[main]{Convention}

\newtheorem{example}[main]{Example}
\newtheorem{examples}[main]{Examples}
\newtheorem{remark}[main]{Remark}
\newtheorem{remarks}[main]{Remarks}

\newtheorem{principle}[main]{Principle}

\theoremstyle{remark}
\newtheorem*{remark*}{Remark}
\newtheorem*{example*}{Example}
\newtheorem*{notation*}{Notation}

\numberwithin{equation}{subsection}
\makeatletter
\renewcommand{\@secnumfont}{\bfseries}
\def\section{\@startsection{section}{1}%
  \z@{.7\linespacing\@plus\linespacing}{.5\linespacing}%
  {\large\normalfont\bfseries\centering}}
\let\c@equation\c@subsubsection

\makeatother

\tikzcdset{proof/.style={sep=normal,font=\tiny,labels={font=\tiny}} } 

\newenvironment{hideproof}{\noindent\textit{Proof.}}{\hfill\qedsymbol}
\newif\ifhideproofs
\hideproofstrue 
\ifhideproofs
	\usepackage{environ}
	\NewEnviron{hide}{}
	
\fi

\newcounter{hideeq}[main]
{\tag{\themain.\thehideeq}\endequation}

\newlist{thmenum}{enumerate}{1}
\setlist[thmenum]{label=(\roman*), ref=\thetheorem.(\roman*)}
\crefalias{thmenumi}{theorem}

\newlist{propenum}{enumerate}{1}
\setlist[propenum]{label=(\roman*), ref=\theproposition.(\roman*)}
\crefalias{propenumi}{proposition}

\newlist{corenum}{enumerate}{1}
\setlist[corenum]{label=(\roman*), ref=\thecorollary.(\roman*)}
\crefalias{corenumi}{corollary}

\newlist{lemenum}{enumerate}{1}
\setlist[lemenum]{label=(\roman*), ref=\thelemma.(\roman*)}
\crefalias{lemenumi}{lemma}

\newlist{exampleenum}{enumerate}{1}
\setlist[exampleenum]{label=(\alph*), ref=\theexample.(\alph*)}
\crefalias{exampleenumi}{example}

\newlist{remarksenum}{enumerate}{1}
\setlist[remarksenum]{label=(\roman*), ref=\theremarks.(\roman*)}
\crefalias{remarksenumi}{remark}

\newlist{defenum}{enumerate}{1}
\setlist[defenum]{label=(\alph*), ref=\thedefinition.(\alph*)}
\crefalias{defenumi}{definition}

\title{6-Functor Formalisms and Smooth Representations}
\author{Claudius Heyer}
\address{Institut für Mathematik, Universität Paderborn, Warburger Straße 100, D-33100 Paderborn, Germany}
\email{\href{mailto:heyer@math.uni-paderborn.de}{heyer@math.uni-paderborn.de}}

\author{Lucas Mann}
\thanks{The project was funded by the Deutsche
Forschungsgemeinschaft (DFG, German Research Foundation) – Project-ID 427320536
– SFB 1442, as well as under Germany's Excellence Strategy EXC 2044 390685587,
Mathematics Münster: Dynamics–Geometry–Structure. The first author was partially supported by the German Research Foundation SFB-TRR 358/1 2023 -- 491392403. The second author was partially supported by ERC Consolidator Grant 770936:NewtonStrat.}
\address{Mathematisches Institut, Universität Münster, Einsteinstraße 62, D-48149 Münster, Germany}
\email{\href{mailto:mann.lucas@uni-muenster.de}{mann.lucas@uni-muenster.de}}


\begin{document}
\begin{abstract}
The purpose of this article is threefold: Firstly, we propose some enhancements to the existing definition of 6-functor formalisms. Secondly, we systematically study the category of kernels, which is a certain 2-category attached to every 6-functor formalism. It provides powerful new insights into the internal structure of the 6-functor formalism and allows to abstractly define important finiteness conditions, recovering well-known examples from the literature. Finally, we apply our methods to the theory of smooth representations of $p$-adic Lie groups and, as an application, construct a canonical anti-involution on derived Hecke algebras generalizing results of Schneider--Sorensen. In an appendix we provide the necessary background on $\infty$-categories, higher algebra, enriched $\infty$-categories and $(\infty,2)$-categories. Among others we prove several new results on adjunctions in an $(\infty,2)$-category and in particular show that passing to the adjoint morphism is a functorial operation.
\end{abstract}
\maketitle
\tableofcontents

\section{Introduction}

Cohomology theories provide a powerful and extremely successful way to understand geometric objects, like manifolds, schemes or rigid varieties, by associating algebraic invariants with them and studying their properties. Some of the most well-known examples include singular cohomology of topological spaces as well as étale and coherent cohomology of schemes. In each of these cases one has a base ring $\Lambda$ and then associates with every \enquote{space} $X$ its cohomology $\Gamma(X, \Lambda) \in \D(\Lambda)$, which is a complex of $\Lambda$-modules. This cohomology is expected to satisfy certain properties, like the Künneth formula or Poincaré duality, which often require a lot of effort to prove with classical methods.

A 6-functor formalism is a way to formalize the notion of a cohomology theory and allows one to simplify and streamline the proofs of the basic properties, as well as enlarge the toolbox one works with. The basic idea is to replace the cohomology modules $\Gamma(X, \Lambda)$ by much more structure. Namely, fixing a category $\cat C$ of geometric objects (which we assume to have finite limits), a 6-functor formalism $\D$ on $\cat C$ captures the following data:
\begin{itemize}
	\item For every space $X \in \cat C$ there is a category $\D(X)$ of \emph{sheaves on $X$}.

	\item On every $\D(X)$ there is a tensor product $\tensor$ which is closed, i.e.\ admits an internal hom $\iHom(M, N)$ for any two sheaves $M, N \in \D(X)$.

	\item For every map $f\colon Y \to X$ in $\cat C$ there is a \emph{pullback functor} $f^*\colon \D(X) \to \D(Y)$ which admits a right adjoint $f_*\colon \D(Y) \to \D(X)$, called the \emph{pushforward}.

	\item For every map $f\colon Y \to X$ in $\cat C$ there is an \emph{exceptional pushforward functor} $f_!\colon \D(Y) \to \D(X)$. It admits a right adjoint $f^!\colon \D(X) \to \D(Y)$, called the \emph{exceptional pullback}.
\end{itemize}
If $*$ denotes the final object in $\cat C$ then one usually has $\D(*) = \D(\Lambda)$. In this case we can define the \emph{cohomology} of an object $X \in \cat C$ as $\Gamma(X, \Lambda) = p_* \one$, where $p\colon X \to *$ is the projection and $\one \in \D(X)$ is the tensor unit. We can also define the \emph{cohomology with compact support} of $X$ as $\Gamma_c(X, \Lambda) \coloneqq p_! \one$. Thus one can see the functors $f_*$ and $f_!$ as a relative version of cohomology and cohomology with compact support, and one can see $\D(X)$ as a relative version of $\D(\Lambda)$.

As explained above, a 6-functor formalism captures the six functors $\tensor, \iHom, f^*, f_*, f_!, f^!$. For this collection of functors to be useful, we need to know that they satisfy additional compatibilities:
\begin{itemize}
	\item $f^*$ commutes with $\tensor$.

	\item $f_!$ and $f^*$ commute in the following sense: Suppose we are given a cartesian square
	\begin{equation}\begin{tikzcd}
		Y' \arrow[r,"g'"] \arrow[d,"f'"'] & Y \arrow[d,"f"] \\
		X' \arrow[r,"g"'] & X
	\end{tikzcd}\end{equation}
	in $\cat C$, i.e.\ such that $Y' = Y \times_X X'$. Then there is an isomorphism of functors $g^! f_* \isom f'_* g'^!$. This property is usually called \emph{proper base-change}.

	\item $f_!$ is compatible with $\tensor$ in the following sense: Given a map $f\colon Y \to X$ and objects $M \in \D(X)$ and $N \in \D(Y)$, there is an isomorphism $M \tensor f_! N \isom f_!(f^* M \tensor N)$. This isomorphism is usually called the \emph{projection formula}.
\end{itemize}
Using these compatibilities one can deduce many formal results about the cohomology theories one started with. We provide two examples:

\begin{example}[Künneth formula]
Given two objects $X, Y \in \cat C$ there is an isomorphism
\begin{align}
	\Gamma_c(X \times Y, \Lambda) \isom \Gamma_c(X, \Lambda) \tensor \Gamma_c(Y, \Lambda)
\end{align}
in $\D(\Lambda)$, known as the \emph{Künneth formula}. It can be proved as follows: Consider the cartesian diagram
\begin{equation}\begin{tikzcd}
	X \times Y \arrow[r,"f_X"] \arrow[d,"f_Y",swap] \arrow[dr,"p"] & Y \arrow[d,"p_Y"] \\
	X \arrow[r,"p_X"'] & *
\end{tikzcd}\end{equation}
Then we compute
\begin{align}
	&\Gamma_c(X \times Y, \Lambda) = p_! \one \isom p_{Y!} f_{X!} \one \isom p_{Y!} f_{X!} f_Y^* \one \isom p_{Y!} p_Y^* p_{X!} \one = p_{Y!}(p_Y^* p_{X!} \one \tensor \one) \isom\\&\qquad\isom p_{X!} \one \tensor p_{Y!} \one = \Gamma_c(X, \Lambda) \tensor \Gamma_c(Y, \Lambda)
\end{align}
using proper base-change and the projection formula.
\end{example}

\begin{example}[Poincaré duality]
Given a \enquote{smooth} space $X \in \cat C$ (e.g.\ a manifold in the topological setting or a smooth scheme in the algebo-geometric setting) there is an equivalence
\begin{align}
	\Gamma(X, \omega_X) \isom \iHom(\Gamma_c(X, \Lambda), \Lambda)
\end{align}
in $\D(\Lambda)$, where we define $\omega_X \coloneqq p^! \one$ and $\Gamma(X, \omega_X) \coloneqq p_* \omega_X$ with $p\colon X \to *$ the projection. This equivalence is known as \emph{Poincaré duality}. To prove it, we note that for every $M \in \D(\Lambda)$ we have
\begin{align}
	&\Hom(M, \Gamma(X, \omega_X)) = \Hom(M, p_* p^! \one) = \Hom(p_! p^* M, \one) \isom \Hom(M \tensor p_! \one, \one) =\\&\qquad= \Hom(M, \iHom(\Gamma_c(X, \Lambda), \Lambda)),
\end{align}
so the claim follows from the Yoneda lemma. If $X$ is an orientable manifold of some dimension $n$ and $\D$ is the 6-functor formalism of sheaves on topological spaces, then the choice of an orientation of $X$ defines an isomorphism $\omega_X \isom \one[n]$. If furthermore $\Lambda$ is a field (so that there are no higher $\Ext$-groups in $\D(\Lambda)$) then the above identity reduces to
\begin{align}
	H^i(X, \Lambda) \isom H_c^{n-i}(X, \Lambda)^\vee,
\end{align}
where $H^i_{(c)}(X, \Lambda) \coloneqq H^i \Gamma_{(c)}(X, \Lambda)$. This recovers the classical formulation of Poincaré duality.
\end{example}

The purpose of this paper is to lay the foundations for an abstract theory of 6-functor formalisms that allows one to formalize constructions and results in different kinds of geometry as well as introduce new ideas and definitions. More precisely, the goals of this paper are threefold:
\begin{enumerate}[(a)]
	\item We introduce an abstract definition of 6-functor formalisms that captures all of the data explained above. This is more subtle than one might think at first, because one certainly wants all the compatibility isomorphisms to be \enquote{natural} in an appropriate sense.

	\item We introduce the 2-category of kernels attached to a 6-functor formalism, which we believe to be the next step in the process of understanding cohomology theories. Using the category of kernels we can abstractly define many useful notions, like étale, smooth and proper maps, as well as derive their basic properties. As an application we derive classical results in algebraic topology from scratch.

	\item We apply the theory of 6-functor formalisms to the study of smooth representations and observe that many constructions and results in representation theory have clean interpretations in terms of a 6-functor formalism. We use this observation to construct a canonical anti-involution on the derived Hecke algebra.
\end{enumerate}
Part (a) is covered in \cref{sec:Corr,sec:6ff} and mainly collects existing results from the literature, with a few enhancements. Part (b) is covered in \cref{sec:kerncat} and is the heart of this paper; while the category of kernels was considered in the literature before, we provide new results and rigorous conceptual proofs of existing results and definitions using the abstract theory of 2-categories. Part (c) is covered in \cref{sec:reptheory}. In addition, we have a rather long appendix which develops the necessary categorical machine and should be of independent interest. In the following subsections of the introduction we provide more information on each of the three parts above.

\subsection{Categorical background}

This paper makes heavy use of category theory, in fact in the appendix we provide many results in abstract category theory which we require for the study of 6-functor formalisms and could not find in the literature. We work with $\infty$-categories throughout because they provide an elegant way to formulate our ideas and in particular solve many naturality issues that occur in classical approaches to 6-functor formalisms. Among others, $\infty$-categories provide the following benefits:
\begin{itemize}
	\item Using $\infty$-categories allows a good treatment of stacks, which are just sheaves in $\infty$-land. Namely, in the classical literature (for example in \cite{stacks-project}) stacks on a site $\cat C$ are defined using certain fibrations $\cat E \to \cat C$. By the straightening/unstraightening correspondence (see \cref{rslt:Straightening-Unstraightening}), these fibrations correspond to $(2,1)$-functors $\cat C^\op \to \Gpd$ satisfying a certain descent condition, where $\Gpd$ is the $(2,1)$-category\footnote{Recall that a $(2,1)$-category is a $2$-category where all $2$-morphisms are isomorphisms; in the case of $\Gpd$ the $2$-morphisms are natural isomorphisms of functors.} of groupoids. Dealing with $(2,1)$-functors and $(2,1)$-categories by hand is cumbersome and hence avoided usually, but in the $\infty$-categorical framework they are just a special case of $\infty$-functors and $\infty$-categories.

	\item We will often use descent techniques to reduce questions about sheaves on a space to questions about sheaves on nice coverings of the space. Naturally the categories of sheaves that we deal with are derived categories and gluing derived categories along covers is not possible with the language of ordinary categories---too much information is lost when passing to the homotopy category of the category of complexes. With the $\infty$-categorical framework, gluing of derived categories can be defined very naturally, see \cref{sec:sheaves}.

	\item In the construction of the 2-category of kernels, $\infty$-categorical techniques are very helpful. In the setting of ordinary 2-categories it is tempting to make constructions by explicitly specifying objects and morphisms and ignoring higher compatibilities, but this easily leads to sloppy proofs. For example, defining the 2-category of kernels by hand seems easy enough, but verifying that it is indeed a 2-category is harder and much more tedious than one might think at first (see the beginning of \cite[\S2.2]{Zavyalov.2023}).

	\item In our application we want to construct an anti-involution on the derived Hecke algebra. Working with derived algebras over an arbitrary ring $\Lambda$ is hard without $\infty$-categorical techniques, and even if $\Lambda$ is a field (where one ends up with differential graded algebras), $\infty$-categorical techniques are extremely helpful and allow one to avoid nasty computations.
\end{itemize}
We refer the reader to \cref{sec:cat} for a quick introduction to $\infty$-categories. Since we believe that $\infty$-categories are the \enquote{correct} version of category theory, we will drop the \enquote{$\infty$} in this paper. For example, we refer to $\infty$-categories as categories and denote the $\infty$-category of $\infty$-categories by $\Cat$. Similarly, we refer to $(\infty,2)$-categories as $2$-categories and so on. The main results of this paper should be applicable without a deep understanding of higher category theory.

The appendix covers many results in abstract category theory which we have not found in the literature and believe to be useful outside of the scope of this paper. Many of these results are straightforward and so we will not explain them here. However, let us highlight some parts of the appendix which we are particularly proud of. The reader who is only interested in 6-functor formalisms is invited to skip the rest of this subsection.

\subsubsection*{Enriched categories}

In \cref{sec:enr} we provide an introduction to enriched category theory based on previous treatments in the literature. Given a category $\cat V$ with tensor product, a $\cat V$-enriched category is a category $\cat C$ where for any two objects $X, Y \in \cat C$ the maps $X \to Y$ form an object $\enrHom[\cat V]_{\cat C}(X, Y)$ of $\cat V$. Examples of enriched categories are ubiquitous in ordinary category theory and are often used without explicit mentioning, e.g.\ every abelian category is enriched in the category of abelian groups (i.e.\ homomorphism sets in an abelian category are abelian groups) or the category of (ordinary) $\Lambda$-valued sheaves on a space is enriched in the category of $\Lambda$-modules, for any ring $\Lambda$. Similar examples exist in $\infty$-category theory, e.g.\ a stable category (see \cite[\S1]{HA}) is naturally enriched in the category of spectra and the (derived) category of $\Lambda$-valued sheaves on a space is enriched in $\D(\Lambda)$; in fact, the $\RHom$ functor secretly encodes the $\D(\Lambda)$-enrichment.

While there are many useful examples, defining enrichment in the $\infty$-setting is harder than in ordinary category theory, as there is an infinite amount of coherences to keep track of. There are several approaches in the literature, some of which we find hard to understand for non-experts. In \cref{sec:enr} we explain these constructions and how they relate in order to (hopefully) make enriched category theory available to a wider audience.

\subsubsection*{Cartesian actions}

Suppose we are given a monoidal category $\cat V$ (i.e.\ a category equipped with an associative tensor product with unit) such that the tensor unit $\one$ is a final object of $\cat V$ and for every two objects $V, W \in \cat V$ the induced maps $V \tensor W \to V$ and $V \tensor W \to W$ induce an isomorphism $V \tensor W = V \times W$. We call such a monoidal category a \emph{cartesian monoidal category}. It is then natural to believe that the monoidal structure on $\cat V$ is uniquely determined, i.e.\ all higher associators are fixed. In the case that $\cat V$ is a \emph{symmetric} monoidal category, i.e.\ $\tensor$ is symmetric in both arguments (which is an additional structure in form of commutator homotopies!), this follows easily from results by Lurie, see \cite[\S2.4.1]{HA}. In the non-symmetric case we did not find a similar result in the literature, so we prove one. Before we state the result, let us indicate another situation where a similar result holds: Suppose $\cat V$ acts on a category $\cat C$, i.e.\ we have a functor $\tensor\colon \cat V \times \cat C \to \cat C$ together with higher homotopies exhibiting associativity. Assume furthermore that $\cat C$ has a final object $*$. We call the action of $\cat V$ on $\cat C$ a \emph{cartesian left action} if for any $V \in \cat V$ and $X \in \cat C$ the natural map $V \tensor X \isoto (V \tensor *) \times X$ is an isomorphism. Then the main results of \cref{sec:alg.cartact} can be summarized as follows:

\begin{proposition}
Let $\cat V$ be a cartesian monoidal category.
\begin{propenum}
	\item The monoidal structure on $\cat V$ is uniquely determined.
	\item \label{rslt:intro-uniqueness-of-cartesian-left-action} Suppose $\cat V$ acts on a category $\cat C$ which has a final object $*$. If this action is a cartesian left action, then it is uniquely determined by the induced functor $\cat V \to \cat C$, $V \mapsto V \tensor *$.
\end{propenum}
\end{proposition}
\begin{proof}
See \cref{rslt:cartesian-monoidal-structure-unique} for (i) and \cref{rslt:uniqueness-of-cartesian-left-actions} for (ii).
\end{proof}

The main reason why we are interested in this result is that it often provides an easy proof that two enrichments are the same (using that enrichments can be induced by left actions, see \cref{def:Lurie-enriched-category}).

\begin{example}
Take $\cat V = \Cat$, so that $\cat V$-enriched categories are the same as 2-categories. Given a category $\cat C$ one can form the 2-category $\Fun(\cat C, \Cat)$, where the 2-categorical structure is induced by the natural action of $\Cat$ on the functor category (see \cref{rslt:2-yoneda}). By the straightening/unstraightening correspondence the underlying 1-category is equivalent to the category $\Cat^\cocart_{/\cat C}$ of cocartesian fibrations over $\cat C$. This category also admits a natural 2-categorical enhancement (see \cref{rslt:2-categorical-enhancement-of-Cat-over-S}) and it is natural to ask if the equivalence
\begin{align}
	\Fun(\cat C, \Cat) = \Cat^\cocart_{/\cat C}
\end{align}
enhances to an equivalence of 2-categories. But by \cref{rslt:intro-uniqueness-of-cartesian-left-action} this amounts to checking that the natural functor from $\Cat$ agrees on both sides, which follows easily from the functoriality of the straightening/unstraightening correspondence.
\end{example}

\subsubsection*{Passing to adjoint morphisms}

In every 2-category $\cat C$ there is a notion for a morphism $f\colon Y \to X$ to be left adjoint to a morphism $g\colon X \to Y$: This is the case if there are a unit $\eta\colon \id_Y \to gf$ and a counit $\varepsilon\colon fg \to \id_X$ satisfying the triangle identities (see \cref{def:adjoint-morphisms}). It is then natural to ask in what sense the association $f \mapsto g$ is functorial in $f$. We answer this question as follows:

\begin{theorem} \label{rslt:intro-passing-to-adjoints}
Let $\cat C$ be a 2-category and let $\cat C^L, \cat C^R \subseteq \cat C$ be the sub-2-categories with the same objects as $\cat C$ but where we only allow left adjoint (resp.\ right adjoint) morphisms. Then there is an isomorphism
\begin{align}
	\cat C^L = (\cat C^R)^{\co,\op}
\end{align}
which acts as the identity on objects and sends a left adjoint morphism to a corresponding right adjoint morphism.
\end{theorem}
\begin{proof}
See \cref{rslt:passing-to-adjoints}.
\end{proof}

Here $(\blank)^\co$ associates with a 2-category the 2-category where all 2-morphisms are inverted, and $(\blank)^\op$ inverts all 1-morphisms. Our proof of \cref{rslt:intro-passing-to-adjoints} is self-contained and does not depend on existing literature beyond the definition of 2-categories. To our knowledge, the statement has not appeared in this generality before.

\subsection{Definition and construction of 6-functor formalisms} \label{sec:intro.def6ff}

A 6-functor formalism encodes a lot of information: Not only does it define the six functors $\tensor$, $\iHom$, $f^*$, $f_*$, $f_!$ and $f^!$, but it also contains isomorphisms exhibiting proper base-change and the projection formula, as well as higher homotopies exhibiting the naturality of these isomorphisms. Finding a good representation of these data is not an easy task, and we rely on previous work by Gaitsgory--Rozenblyum based on a suggestion by Lurie (see \cite[Part~III]{Gaitsgory-Rozenblyum:Vol1}) and by Liu--Zheng (see \cite{Liu-Zheng.2012}). There are more recent references based on these ideas, see for example \cite[Appendix~A.5]{Mann.2022a}, \cite{Scholze:Six-Functor-Formalism}, \cite{Khan:Weaves}, \cite{Kujper:6ff} to name a few. In this paper we recall the definitions and basic constructions of 6-functor formalisms, following mostly the exposition in \cite[Appendix~A.5]{Mann.2022a} with a few enhancements. 

\begin{remark}
Currently our main construction result relies on subtle unpublished results by Liu--Zheng \cite{Liu-Zheng:Gluing-Nerves} which are very model dependent and hard to understand. In work in progress, we are currently developing new (and model independent) foundations for the construction of 6-functor formalisms.
\end{remark}

Without further ado, let us introduce the definition of 6-functor formalisms. Suppose we are given a \emph{geometric setup} $(\cat C, E)$, i.e.\ a pair consisting of a category $\cat C$ of \enquote{spaces} and a collection of \emph{exceptional morphisms} $E$ in $\cat C$ satisfying some basic conditions (see \cref{defn:geometric-setup}).\footnote{We warn the reader that we changed the definition of geometric setups compared to the original one in \cite[Definition~A.5.1]{Mann.2022a}: We additionally require $E$ to be stable under taking diagonals; see \cref{rmk:geometric-setup-morphisms}.} We want to define the notion of a 6-functor formalism on $\cat C$, where the functors $f_!$ and $f^!$ are only defined for maps $f$ which lie in $E$ (in practice one often cannot define the $!$-functors for all maps $f$). With $(\cat C, E)$ we can associate the following category:

\begin{definition}
Let $(\cat C, E)$ be a geometric setup. We define the \emph{category of correspondences} $\Corr(\cat C, E)$ to be the following category:
\begin{itemize}
	\item The objects of $\Corr(\cat C, E)$ are those of $\cat C$.
	
	\item For every two objects $X, Y \in \cat C$, the homomorphisms from $X$ to $Y$ in $\Corr(\cat C, E)$ are given by diagrams
	\begin{align}
		\CorrHom{X}{g}{X'}{f}{Y}
	\end{align}
	in $\cat C$ such that $f \in E$.
	
	\item The composition in $\Corr(\cat C, E)$ is given by
	\begin{align}
		\left[ \CorrHom{Y}{g'}{Y'}{f'}{Z} \right] \comp \left[ \CorrHom{X}{g}{X'}{f}{Y} \right] = \begin{tikzcd}[column sep=tiny,row sep=small,ampersand replacement=\&]
			\& \& X' \times_Y Y' \ar[dl] \ar[dr] \\
	        \& X' \ar[dl,"g"] \ar[dr,"f"] \& \& Y' \ar[dl,"g'"] \ar[dr,"f'"] \\
	       	X \& \& Y \& \& Z
	    \end{tikzcd}
	\end{align}
\end{itemize}
The category $\Corr(\cat C, E)$ is equipped with a symmetric monoidal structure (i.e.\ a tensor product) given by the formula $X \tensor Y \coloneqq X \times Y$, where on the right-hand side we denote the product in $\cat C$.\footnote{For simplicity we implicitly assume that $\cat C$ admits finite products. This is not strictly necessary for the definition of 6-functor formalisms and the main text works without this assumption.}
\end{definition}

We refer the reader to \cref{sec:Corr} for the precise definitions of geometric setups (see \cref{subsec:geometric-setups}), the correspondence category (see \cref{def:Correspondence-Category}) and its symmetric monoidal structure (see \cref{def:Corr-operad}). Note that even if $\cat C$ is an ordinary category, $\Corr(\cat C, E)$ is a $(2, 1)$-category, because the composition is given by pullbacks, which are only defined up to isomorphism.

With the category of correspondences at hand, we can now provide a formal definition of a 6-functor formalism on $(\cat C, E)$. The first observation is that only three of the six functors need to be encoded, namely $\tensor$, $f^*$ and $f_!$. Namely, the other three functors are defined to be right adjoints of these first three and hence automatically inherit all desired compatibilities.

\begin{definition}
Let $(\cat C, E)$ be a geometric setup. A \emph{3-functor formalism} on $(\cat C, E)$ is a lax symmetric monoidal functor
\begin{align}
	\D\colon \Corr(\cat C, E) \to \Cat, \qquad X \mapsto \D(X),
\end{align}
where we equip $\Cat$ with the tensor product given by product.
\end{definition}

This definition is taken straight from \cref{def:3-functor-formalism}. We observe that a 3-functor formalism encodes the following data:
\begin{itemize}
	\item For every $X \in \cat C$ there is a category $\D(X)$.
	\item For every map $g\colon Y \to X$ in $\cat C$ there is a correspondence $\CorrHom{X}{g}{Y}{\id}{Y}$ and the functor $\D$ provides an associated functor $g^*\colon \D(X) \to \D(Y)$.
	\item For every map $f\colon Y \to X$ in $E$ there is a correspondence $\CorrHom{Y}{\id}{Y}{f}{X}$ and the functor $\D$ provide an associated functor $f_!\colon \D(Y) \to \D(X)$.
	\item Given $X \in \cat C$, the lax symmetric monoidal structure on $\D$ defines a natural map $\D(X) \times \D(X) \to \D(X \times X)$. By composing this map with the pullback $\Delta^*\colon \D(X \times X) \to \D(X)$ along the diagonal $\Delta$ of $X$, we obtain a functor $\tensor\colon \D(X) \times \D(X) \to \D(X)$.
\end{itemize}
Thus we indeed obtain the three functors $\tensor$, $f^*$ and $f_!$. Note that a general correspondence $\CorrHom{X}{g}{X'}{f}{Y}$ is sent by $\D$ to the functor $f_! g^*\colon \D(X) \to \D(Y)$. Together with the definition of the composition of morphisms in $\Corr(\cat C, E)$ one easily obtains proper base-change (which is thus encoded in the composition functoriality of the functor $\D$!). One similarly obtains the projection formula by looking closer at the lax symmetric monoidal structure on $\D$. We refer the reader to \cref{rslt:compatibilities-of-3ff} for details. It now makes sense to reintroduce the remaining three functors:

\begin{definition}
A \emph{6-functor formalism} is a 3-functor formalism $\D$ such that the right adjoints $\iHom$, $f_*$ and $f^!$ exist.
\end{definition}

We refer the reader to \cref{sec:6ff.def} for a more elaborate version of that definition as well as a list of additional compatibilities enjoyed by $\iHom$, $f_*$ and $f^!$. We emphasize that one of the observations in this paper is that the three functors $\iHom$, $f_*$ and $f^!$ are not required for most definitions and constructions in this paper---even notions like cohomologically smooth maps do not require the existence of $f^!$-functors. We refer the reader to \cref{sec:intro.kerncat} for more information.

We have discussed how to properly define 6-functor formalisms. For this definition to be of any use, we need a way to construct them. Here we recall the following result from \cite[Proposition~A.5.10]{Mann.2022a} which is heavily based on results by Liu--Zheng \cite{Liu-Zheng:Gluing-Nerves,Liu-Zheng.2012}. The idea of the construction is as follows. It is usually not hard to construct a functor $\D\colon \cat C^\op \to \CMon$, where $\CMon$ is the category of symmetric monoidal categories, i.e.\ this functor $\D$ encodes the functors $\tensor$ and $f^*$. In order to add $f_!$ into the mix, we write every map $f \in E$ as a composition $f = g \comp j$ of an \enquote{open immersion} $j$ and a \enquote{proper} map $g$. Then $f_! = g_! \comp j_!$, where $j_!$ is left adjoint to $j^*$ and $g_!$ is right adjoint to $g^*$. The precise definition of such a \emph{suitable decomposition} of $E$ into classes $I, P \subseteq E$ of \enquote{open immersions} and \enquote{proper} maps is given in \cref{def:suitable-decomposition}. With this definition at hand, we can formulate the following construction result:

\begin{proposition} \label{rslt:intro-construct-3ff}
Let $(\cat C, E)$ be a geometric setup and $I, P \subseteq E$ be a suitable decomposition. Suppose $\D\colon \cat C^\op \to \CMon$ is a functor with the following properties:
\begin{enumerate}[(a)]
	\item For every $[j\colon U \to X] \in I$ the functor $j^*$ admits a left adjoint $j_!\colon \D(U) \to \D(X)$ which satisfies base-change and the projection formula.
	\item For every $[g\colon Y \to X] \in P$ the functor $g^*$ admits a right adjoint $g_*\colon \D(Y) \to \D(X)$ which satisfies base-change and the projection formula.
	\item For every cartesian square
	\begin{equation}\begin{tikzcd}
		V \arrow[r,"j'"] \arrow[d,"g'"] & Y \arrow[d,"g"]\\
		U \arrow[r,"j"] & X
	\end{tikzcd}\end{equation}
	in $\cat C$ such that $j \in I$ and $g \in P$, the natural transformation $j_! g'_* \isoto g_* j'_!$ is an isomorphism of functors $\D(V) \to \D(X)$.
\end{enumerate}
Then $\D$ extends to a 3-functor formalism $\D\colon \Corr(\cat C, E) \to \Cat$ such that $f_! = g_* j_!$ for every decomposition $f = gj$ with $j \in I$ and $g \in P$.
\end{proposition}

We refer the reader to \cref{rslt:construct-3ff-from-suitable-decomp} for a more elaborate version of the result. As mentioned above this result appeared in the literature before but we include it here due to its high significance for the construction of 6-functor formalisms. We remark that in (a) and (b) of \cref{rslt:intro-construct-3ff} both base-change and projection formula are a \emph{condition} rather than additional structure because in both cases the $!$-functor is given as an adjoint of the pullback functor. One of the new results in this paper is that condition (c) is often automatic:

\begin{lemma}
In the setting of \cref{rslt:construct-3ff-from-suitable-decomp} assume that (a) and (b) are satisfied and that all $j \in I$ are monomorphisms. Then (c) is automatically satisfied.
\end{lemma}

For the proof and a more general version of this result we refer the reader to \cref{rslt:cond-c-satisfied-on-tensor-product}, see specifically \cref{rslt:cond-c-auto-for-fully-faithful-j-shriek,rmk:cond-c-auto-for-monomorphism-j}.

With \cref{rslt:intro-construct-3ff} we can construct the functor $f_!$ whenever $f$ admits a decomposition into an open immersion and a proper morphism. While this is often satisfied in controlled geometric situations, one of the observations in \cite{Liu-Zheng.2012,Gulotta-Hansen-Weinstein:Enhanced-6FF-on-vStacks} and \cite[\S A.5]{Mann.2022a} is that one can extend the definition of $f_!$ for a lot more maps $f$. In fact, we have the following result based on \cite[Theorem~4.20]{Scholze:Six-Functor-Formalism}, which composes the extension results from \cite[\S A.5]{Mann.2022a} in a clever way. In the following we call a 3-functor formalism $\D$ \emph{sheafy} if the functor $\D\colon \cat C^\op \to \Cat$ defines a sheaf of categories for some predefined site on $\cat C$ (see \cref{sec:sheaves} for the definition of sheaves). Moreover, we call $\D$ \emph{presentable} if it factors over the category of presentable categories (this is usually satisfied in practice, see \cref{sec:cat.pres} for the definition of presentable categories).

\begin{theorem}
Let $(\cat C, E)$ be a geometric setup with $\cat C$ a site. Let $\D\colon \Corr(\cat C, E) \to \PrL$ be a sheafy presentable 6-functor formalism on $(\cat C, E)$. Then there is a collection of edges $E'$ on $\cat X \coloneqq \Shv(\cat C)$ with the following properties:
\begin{thmenum}
	\item $(\cat X, E')$ is a geometric setup and $\D$ extends uniquely to a sheafy presentable 6-functor formalism $\D'$ on $(\cat X, E')$.
	\item $E'$ is $*$-local on the target and $!$-local on source and target.
	\item $E'$ is tame, i.e.\ every map in $E'$ with target in $\cat C$ is $!$-locally on the source in $E$.
\end{thmenum}
\end{theorem}

We refer the reader to \cref{rslt:extend-6ff-to-stacks-and-stacky-maps} for a more elaborate and more precise formulation of the result. The upshot is that a sheafy presentable 6-functor formalism on some site $\cat C$ can always be extended uniquely to a sheafy presentable 6-functor formalism on the topos $\cat X = \Shv(\cat C)$ of sheaves on $\cat C$. In more classical terminology, $\cat X$ is the category of ($\infty$-)stacks on $\cat C$. Moreover, the class $E'$ of exceptional edges for this extension contains all edges which are representable in $E$ as well as all edges which \emph{$!$-locally} on source or target lie in $E$. Here a $!$-cover is a family $(U_i \to X)_i$ along which $\D^!$ admits universal descent, where $\D^!$ is the functor $\D$ with $f^!$-functors as transition maps. An important example of $!$-covers are suave covers (see \cref{rslt:!-descent-along-suave-map}), so one should expect a robust definition of $!$-functors for all \enquote{Artin stacks} on $\cat C$. There are other interesting examples of $!$-covers however, as we will see in \cref{sec:intro.reptheory}.

\begin{example} \label{ex:intro-6ff-on-CondAni}
Let $\Cond(\Ani)$ be the category of \emph{condensed anima}, i.e.\ the topos of sheaves on the site of ($\kappa$-small) profinite sets (see \cref{def:condensed-anima}). It contains the category of ($\kappa$-small) locally compact Hausdorff spaces as a full subcategory, but also allows stacky phenomena, e.g.\ one can form the classifying stack $*/G$ of a locally profinite group $G$. Fix moreover some ring $\Lambda$. With the above constructions we can define a sheafy presentable 6-functor formalism
\begin{align}
	\D(\blank,\Lambda)\colon \Corr(\Cond(\Ani), \fine{\Lambda}) \to \Cat,
\end{align}
where $\D(X,\Lambda)$ is roughly the (derived) category of sheaves of $\Lambda$-modules on the condensed anima $X$ (see \cref{def:6ff-on-Cond-Ani}). The class $\fine{\Lambda}$ of \emph{$\Lambda$-fine} maps in $\Cond(\Ani)$ contains all maps between nice enough topological spaces as well as maps between classifying stacks of nice enough locally profinite groups; see \cref{ex:intro-basic-properties-of-6ff-on-CondAni,sec:intro.reptheory} for more information.
\end{example}

\subsection{The category of kernels} \label{sec:intro.kerncat}

The focus of this paper is the rigorous treatment of the \emph{2-category of kernels} $\cat K_\D$ associated with a 6-functor formalism $\D$, which we believe to be an extremely powerful tool for studying $\D$---in fact we believe that at some point one should replace $\D$ by $\cat K_\D$. The category of kernels was first defined in \cite[\S IV.2.3.3]{Fargues-Scholze:Geometrization} where it was the key ingredient in the proof of Verdier duality for universally locally acyclic sheaves. The definition is based on a related (but different) construction by Lu--Zheng \cite{Lu-Zheng:ULA} which is why Fargues--Scholze gave the 2-category the name \emph{Lu--Zheng category} and denote it $\operatorname{LZ}_\D$. We chose to instead use the name \enquote{category of kernels}, the meaning of which will be explained below.

The definition in \cite{Fargues-Scholze:Geometrization} is rather ad-hoc (and only defines an \emph{ordinary} 2-category), but a more conceptual definition using enriched category theory was provided in \cite[Proposition~2.2.6]{Zavyalov.2023}. The category of kernels also appeared in \cite{Mann.2022b,Scholze:Six-Functor-Formalism} where it is used to define notions like cohomologically smooth and proper maps in a 6-functor formalism and provide some basic properties. In this paper we provide a comprehensive study of the category of kernels, with conceptual and easy arguments for all of the results (as well as many results that have not appeared in the literature before). All of our arguments work on the $(\infty,2)$-categorical level in contrast to previous results in the literature.

Let us finally come to the definition of the category of kernels. The following definition is taken from \cref{def:kern-cat}:

\begin{definition}
Let $\D$ be a 3-functor formalism on some geometric setup $(\cat C, E)$ and fix an object $S \in \cat C$. The \emph{2-category of kernels} associated with $\D$ and $S$ is the 2-category
\begin{align}
	\cat K_{\D,S} \in \TwoCat
\end{align}
given by transferring the self-enrichment of $\Corr((\cat C_E)_{/S})$ along $\D$. It has the following explicit description:
\begin{itemize}
	\item The objects of $\cat K_{\D,S}$ are those of $(\cat C_E)_{/S}$, i.e.\ objects $X \in \cat C$ together with a map $X \to S$ in $E$.
	\item Given two objects $X \to S$ and $Y \to S$ in $\cat K_{\D,S}$, the category of morphisms $Y \to X$ in $\cat K_{\D,S}$ is given by
	\begin{align}
		\Fun_S(Y, X) \coloneqq \Fun_{\cat K_{\D,S}}(Y, X) = \D(X \times_S Y).
	\end{align}
	\item Given morphisms $M\colon Y \to X$ and $N\colon Z \to Y$ in $\cat K_{\D,S}$, represented by objects $M \in \D(X \times_S Y)$ and $N \in \D(Y \times_S Z)$, the composition is given by
	\begin{align}
		M \comp N = \pi_{13!}(\pi_{12}^* M \tensor \pi_{23}^* N) \in \D(X \times_S Z),
	\end{align}
	where $\pi_{ij}$ denote the projections on $X \times_S Y \times_S Z$.
\end{itemize}
\end{definition}

Our first result is that for every $S \in \cat C$ the category of kernels $\cat K_{\D,S}$ fits into the following diagram of (2-)categories and (2-)functors:

\begin{proposition}
Let $\D$ be a 3-functor formalism on some geometric setup $(\cat C, E)$ and let $S$ be an object in $\cat C$. Then there is a diagram of 2-functors
\begin{equation}\begin{tikzcd}
	(\cat C_E)_{/S} \arrow[dr]\\
	& \Corr((\cat C_E)_{/S}) \arrow[r,"\Phi_{\D,S}"] & \cat K_{\D,S} \arrow[r,"\Psi_{\D,S}"] & \Cat\\
	(\cat C_E)_{/S}^\op \arrow[ur]	
\end{tikzcd}\end{equation}
such that $\Psi_{\D,S} \comp \Phi_{\D,S} = \D$.
\end{proposition}
\begin{proof}
See \cref{rslt:functors-from-and-to-kercat} and the explanations below.
\end{proof}

Here the diagonal functors are the canonical embeddings used in the definition of $f^*$ and $f_!$ in \cref{sec:intro.def6ff}. The functor $\Phi_{\D,S}$ is the identity on objects and acts on morphisms by
\begin{align}
	[\CorrHomInline{X}{g}{Z}{f}{Y}] \mapsto [(f,g)_! \one \in \D(Y \times_S X)],
\end{align}
where $(f,g)\colon Z \to Y \times_S X$ is the obvious map. Finally, the 2-functor $\Psi_{\D,S}$ is given by $\Fun_S(S,\blank)$. Explicitly, it sends an object $X \in \cat K_{\D,S}$ to $\D(X)$ and it sends a morphism $M \in \Fun_S(Y, X) = \D(X \times_S Y)$ to the functor
\begin{align}
	\pi_{1!}(M \tensor \pi_2^*)\colon \D(Y) \to \D(X).
\end{align}
Therefore, the morphism $M$ in $\cat K_{\D,S}$ acts as a \emph{kernel} for a functor $\D(Y) \to \D(X)$, hence the name \emph{category of kernels}.

\begin{remark}
See also \cite[Proposition~3.3.21]{Lurie:Cobordism} for a slightly different perspective on $\cat K_{\D,S}$, characterizing it as the unique symmetric monoidal 2-category with an essentially surjective symmetric monoidal map from $\Corr((\cat C_E)_{/S})$ such that the above factorization of $\D$ exists.
\end{remark}

The main power of the category of kernels comes from the fact that it allows to conceptualize many important constructions in a 3-functor formalism. This provides a new and much simpler view on these constructions and therefore leads to simple and elegant arguments for things that otherwise end up as a huge technical nightmare. Let us make the following core definitions (see \cref{def:adm-and-coadm-sheaves,def:suave-and-prim-map}):

\begin{definition}
Let $\D$ be a 3-functor formalism on some geometric setup $(\cat C, E)$, let $f\colon X \to S$ be a map in $E$ and let $P \in \D(X)$.
\begin{defenum}
	\item $P$ is called \emph{$f$-suave} if it is a left adjoint morphism $X \to S$ in $\cat K_{\D,S}$, where we observe that $\Fun_S(X, S) = \D(X)$. We denote the associated right adjoint morphism by $\DSuave_f(P) \in \D(X)$.

	\item $P$ is called \emph{$f$-prim} if it is a right adjoint morphism $X \to S$. We denote the associated left adjoint morphism by $\DPrim_f(P) \in \D(X)$.
\end{defenum}
\end{definition}

\begin{definition}
Let $\D$ be a 3-functor formalism on some geometric setup $(\cat C, E)$ and let $f\colon Y \to X$ be a map in $E$.
\begin{defenum}
	\item \label{def:intro-suave-map} We say that $f$ is \emph{$\D$-suave} if it is sent to a left adjoint morphism under the functor $(\cat C_E)_{/X} \to \cat K_{\D,X}$ from above or, equivalently, if $\one \in \D(Y)$ is $f$-suave. In this case we call $\omega_f \coloneqq \DSuave_f(\one) \in \D(Y)$ the \emph{dualizing complex} of $f$. We say that $f$ is \emph{$\D$-smooth} if additionally $\omega_f$ is invertible.

	\item We say that $f$ is \emph{$\D$-prim} if it is sent to a right adjoint morphism in $\cat K_{\D,X}$. In this case we call $\delta_f \coloneqq \DPrim_f(\one) \in \D(Y)$ the \emph{codualizing complex} of $f$
\end{defenum}
\end{definition}

While the above definitions look very simple, they encode important finiteness conditions in many different settings. For example, $f$-suave objects correspond to universally acyclic sheaves in the étale setting (which is the original motivation for introducing the category of kernels in \cite{Fargues-Scholze:Geometrization}), to admissible representations when applied to a classifying stack (see \cref{rslt:intro-suave-duality-on-Rep-G} below) and to coherent sheaves in the quasi-coherent setting (this will be shown in upcoming work of the second author and David Hansen). The $f$-prim objects are more mysterious and closely related to compact objects. For example, on the classifying stack of a locally profinite group they correspond exactly to compact representations, which is the source of our main application to representation theory (see \cref{sec:intro.reptheory}). Primness was first considered in \cite[Definition~7.12]{Mann.2022b}.

\begin{remark}
The name \enquote{prim} was suggested by David Hansen. It is supposed to confer the idea of something related to \enquote{properness}, but still different (we define proper maps below). The name \enquote{suave} was suggested by Peter Scholze as a word close to \enquote{smooth}.
\end{remark}

The definition of $f$-suave and $f$-prim objects in terms of adjoint morphisms is very abstract. In \cref{rslt:criterion-for-suave-object,rslt:criterion-for-prim-object} we recall more concrete criteria for checking suaveness and primness, at least in the presence of all \emph{six} functors. For example, $P \in \D(X)$ is suave along $f\colon X \to S$ if and only if the natural map
\begin{align}
	\pi_1^* \iHom(P, f^! \one) \tensor \pi_2^* P \isoto \iHom(\pi_1^* P, \pi_2^! P)
\end{align}
is an isomorphism in $\D(X \times_S X)$, where $\pi_i\colon X \times_S X \to X$ are the two projections. In this case we also have $\DSuave_f(P) = \iHom(P, f^! \one)$; in particular the dualizing complex of a $\D$-suave map is given by $\omega_f = f^! \one$. There is a similar, but less useful, formula for primness and prim duality. These criteria follow immediately from a simple general criterion for adjunctions in an arbitrary 2-category, which we develop in the appendix (see \cref{rslt:pointwise-criterion-for-adjunction}).

From basic properties of adjunctions in 2-categories, we can deduce many properties of suave and prim objects. We list some of them here and refer the reader to \cref{sec:kerncat.suave-prim-obj,sec:kerncat.suave-prim-map} for the full list.
\begin{enumerate}[(1)]
	\item Suave objects are stable under suave duality and for all $f$-suave objects $P$ there is a natural isomorphism $\DSuave_f(\DSuave_f(P)) = P$. This follows from uniqueness of adjoints and $\cat K_{\D,S} \isom \cat K_{\D,S}^\op$, see \cref{rslt:SD-PD-are-self-inverse-equivalences}. The same is true for prim in place of suave.

	\item Suaveness and primness are local on the target, see \cref{rslt:suave-and-prim-obj-is-local-on-target,rslt:suave-and-prim-map-is-local-on-target}. Under mild assumptions on $\D$, suaveness is suave local on the source, see \cref{rslt:suave-is-local-on-source,rslt:suave-map-is-local-on-source}. Primness is usually not local on the source, but see \cref{rmk:how-to-apply-locality-of-primness-on-source,rslt:check-prim-object-on-!-descendable-cover} for important exceptions.

	\item Suave objects are stable under suave pullback and prim exceptional pushforward, and the same is true with the roles of suave and prim swapped; see \cref{rslt:suave-pullback-and-prim-pushforward-preservations}.

	\item \label{rslt:intro-twist-of-shriek-and-star} If $f$ is $\D$-suave then there are natural isomorphisms $f^! = f^* \tensor \omega_f$ and $f^* = \iHom(\omega_f, f^!)$. If $f$ is $\D$-prim then there are natural isomorphisms $f_* = f_!(\delta_f \tensor \blank)$ and $f_! = f_* \iHom(\delta_f, \blank)$. See \cref{rslt:suave-and-prim-maps-induce-twist-of-shriek-functors} for a proof, and see \cref{rslt:suave-and-prim-obj-induce-twist-of-shriek-functors} for a more general version for suave and prim objects.

	\item If $\D$ is stable, i.e.\ all $\D(X)$ are stable categories (roughly, triangulated) and $\tensor$, $f^*$ and $f_!$ are exact, then suave and prim objects are stable under (co)fibers and retracts.
\end{enumerate}
From \ref{rslt:intro-twist-of-shriek-and-star} we see that suave and prim objects provide a close relationship between $!$- and $*$-functors. We also see that $\D$-suaveness is closely related to the classical definition of cohomological smoothness, where one requires the identity $f^! = f^* \tensor \omega_f$ to hold after any base-change (and additionally assumes $\omega_f$ to be invertible, see \cref{rmk:cohomologically-smooth-maps} for a discussion).

In addition to $\D$-suave and $\D$-prim maps one can also consider the more restrictive notions of $\D$-étale and $\D$-proper maps, which are discussed in \cref{sec:kerncat.etale-proper} and where we have $f^! = f^*$ and $f_! = f_*$, respectively.

As mentioned above, we provide short and conceptual (and oftentimes new) proofs for all the statements above that rely on general properties of 2-categories. Perhaps the most subtle argument is the proof that suave and prim objects are local on the target. To see this we use the functoriality of the 2-category of kernels:

\begin{theorem} \label{rslt:intro-functoriality-of-kercat}
Let $\D$ be a 3-functor formalism on some geometric setup $(\cat C, E)$. Then the assignment $S \mapsto \cat K_{\D,S}$ upgrades to a lax symmetric monoidal 2-functor
\begin{align}
	\cat K_{\D,(\blank)}\colon \Corr(\cat C, E) \to \TwoCat, \qquad S \mapsto \cat K_{\D,S}.
\end{align}
\end{theorem}

In other words, the 2-category of kernels forms itself a 3-functor formalism, but one categorical level higher. We refer the reader to \cref{rslt:functoriality-of-kerncat} for the proof and for a more precise statement, where we make the functoriality explicit.

\begin{example} \label{ex:intro-basic-properties-of-6ff-on-CondAni}
Using the results on the category of kernels, we further study the 6-functor formalism $\D(\blank,\Lambda)$ on condensed anima from \cref{ex:intro-6ff-on-CondAni}, see \cref{sec:kerncat.examples}. We show \emph{from scratch} (without using any results from algebraic topology) that all maps between locally \emph{finite dimensional} (see \cref{def:finite-dimensional-CHaus-space}) locally compact Hausdorff spaces are $\Lambda$-fine and moreover that the 6-functor formalism behaves as expected: $f_!$ computes (relative) cohomology with compact support, proper maps are $\D$-proper, open immersions are $\D$-étale and topological manifolds are $\D$-suave. In particular we deduce Poincaré duality for the cohomology of manifolds (see \cref{rslt:Poincare-duality-on-manifolds}), although we only show that the dualizing complex is locally free of the correct rank and do not explicitly compare it to the orientation sheaf---this is postponed to a future paper, where we will work out a general technique to compute dualizing complexes. Using our results on the category of kernels, these claims reduce easily to the following two computations:
\begin{enumerate}[(a)]
	\item The unit interval $I  = [0,1]$ is $\D$-proper (and in particular $\Lambda$-fine) over $*$, see \cref{rslt:unit-interval-computations}. Here the main idea is to cover $I$ by the Cantor set (which is profinite) and do an explicit computation of the descent data.

	\item The real line $\RR$ is $\D$-suave, see \cref{rslt:6ff-on-manifolds}. We show this by explicitly computing the required adjunction in $\cat K_{\D,*}$ according to \cref{def:intro-suave-map}, following a similar argument in \cite[Proposition~5.10]{Scholze:Six-Functor-Formalism}. We spend a considerable amount of work to showing that certain maps are the expected ones.
\end{enumerate}
The above results provide the classical 6-functor formalism in algebraic topology, computing cohomology and cohomology with compact support. The main difference to many other sources (see e.g.\ \cite{Volpe:6ff-in-Topology}) is that we spend only little effort \emph{constructing} the six functors---they exist for all $\Lambda$-fine maps. Instead, the non-trivial computations are shifted towards \emph{proving} that certain maps are $\Lambda$-fine, and here one can often use general stability properties of $\Lambda$-fine maps to reduce to easy base cases. In summary, we have shifted the problem from \emph{providing a structure} to \emph{proving a property}. One can push this approach to six functors on topological spaces much further (see e.g.\ \cite[Lecture~VII]{Scholze:Six-Functor-Formalism}), but we do not pursue this in this paper.
\end{example}

\begin{example}
At the end of \cref{sec:kerncat.examples} we study a different kind of $\Lambda$-fine maps: All truncated maps between anima are $\Lambda$-fine, and the $f_!$-functors compute (relative) homology in this case. The truncatedness assumption is an artifact of our current construction of the 6-functor formalism and can hopefully be removed in the future. We refer the reader to \cref{ex:primness-of-anima} for a fun interplay between the six functors on topological spaces and on anima.
\end{example}

\begin{remark} \label{rmk:presentable-kerncat}
It has been observed by Aoki, Scholze and the second author that there is a more enhanced version of the 2-category of kernels which allows one to better deal with infinite stacks. The idea is as follows: Given a presentable 6-functor formalism $\D$ on a geometric setup $(\cat C, E)$ with $E = all$, the 2-category $\cat K_\D$ is not only a 2-category, but a $\PrL$-enriched category (we ignore set-theoretic considerations that go into the definition of this notion). One can therefore pass to the $\PrL$-enriched Yoneda category
\begin{align}
	\cat K'_{\D} \coloneqq \PSh[\PrL](\cat K_{\D}) = \Fun_{\PrL}(\cat K_{\D}, \PrL),
\end{align}
where the right-hand side denotes the category of $\PrL$-enriched functors $\cat K_\D \to \PrL$. By abstract nonsense, the functor $\cat C \to \K_\D$ extends uniquely to a colimit-preserving functor $\PSh(\cat C) \to \cat K'_\D$, so we can view any prestack on $\cat C$ as an object in $\cat K'_\D$. Given prestacks $X = \varinjlim_i X_i$ and $Y = \varinjlim_j Y_j$ with $X_i, Y_j \in \cat C$, we have
\begin{align}
	\Fun_{\cat K'_\D}(Y, X) = \varprojlim_j \varinjlim_i \D(X_i \times Y_j),
\end{align}
where the transition maps are given by pullbacks along $Y_{j'} \to Y_j$ and exceptional pushforwards along $X_{i'} \to X_i$. As a special case we get
\begin{align}
	\D^*(X) \coloneqq \Fun_{\cat K'_\D}(X, *) = \varprojlim_i \D(X_i), \qquad \D^!(X) \coloneqq \Fun_{\cat K'_\D}(*, X) = \varinjlim_i \D(X_i).
\end{align}
Then $\D^*$ is a theory which admits all $*$-functors and $\D^!$ is a theory which admits $!$-functors. If $X$ is \enquote{small}, i.e.\ self-dual in $\cat K_\D$, then $\D^*(X) = \D^!(X)$ and we have both $*$-functors and $!$-functors on $X$. The two sheaf categories $\D^*$ and $\D^!$ have appeared already in the literature (see e.g.\ \cite[\S2]{Raskin:Infinite-Dimensional}) but the 2-category $\cat K'_\D$ makes them appear naturally and ties them together. For example, if $X = \varinjlim X_i$ is a colimit of $\D$-suave stacks then by (the dual version of) \cref{rslt:limits-and-adjunctions} the map $X \to *$ is left adjoint in $\cat K'_\D$ and its right adjoint is given by an object $\omega_X \in \D^!(X)$. We will pursue these ideas in a future paper.
\end{remark}

\begin{remark}
Other 2-categories associated with a 6-functor formalism have appeared in the literature, and we briefly sketch their relation to the 2-category of kernels. In the following we fix a 6-functor formalism $\D$ on a geometric setup $(\cat C, E)$ and for simplicity we assume $E = all$. We furthermore fix an object $S \in \cat C$ and assume that $\cat C$ has finite products (as well as fiber products).
\begin{itemize}
	\item In \cite[Construction~2.6]{Lu-Zheng:ULA} Lu--Zheng define the \emph{2-category of cohomological correspondences} $\CoCorr_{\D,S}$, which they denote by $\cat C_S$. An object in $\CoCorr_{\D,S}$ is a pair $(X, M)$ consisting of a map $X \to S$ and an object $M \in \D(X)$. A morphism $(X, M) \to (Y, N)$ in $\CoCorr_{\D,S}$ consists of a correspondence $X \xfrom{f} Z \xto{g} Y$ together with a map $f^* M \to g^! N$ in $\D(Z)$. A 2-morphism consists of a $\D$-proper map $Z \to Z'$ of correspondences over $X$ and $Y$ which is compatible with the maps $f^* M \to g^! N$ and $f'^* M \to g'^! N$ in a suitable sense.

	Lu--Zheng use $\CoCorr_{\D,S}$ to provide an elegant characterization of ULA étale sheaves: A sheaf $M \in \D(X)$ is ULA over $S$ if and only if the pair $(X, M)$ is dualizable in $\CoCorr_{\D,S}$ (see \cite[Theorem~2.16]{Lu-Zheng:ULA}). This remarkable observation is what lead Fargues--Scholze to their version of the 2-category, which is the one we call the 2-category of kernels now.

	There is a more conceptual construction of $\CoCorr_{\D,S}$: As observed in \cite[Remark~2.7]{Lu-Zheng:ULA} one can construct $\CoCorr_{\D,S}$ as the symmetric monoidal 2-unstraightening of the 2-functor $\D\colon \Corr'(\cat C_{/S}) \to \Cat$, where $\Corr'$ denotes the 2-category of correspondences where we allow $\D$-proper maps as above as 2-morphisms. Of course some effort is required to make this construction rigorous in the $\infty$-categorical setting; we do not pursue this further.

	\item In \cite{Hansen-Scholze:Relative-Perversity} Hansen--Scholze consider the (co)lax (co)slice 2-category $\cat C'_{\D,S} \coloneqq (\cat K_{\D,S})_{/S}$. This is a symmetric monoidal 2-category whose objects are pairs $(X, M)$ consisting of a map $X \to S$ in $\cat C$ and a sheaf $M \in \D(X)$ (viewed as a morphism $X \to S$ in $\cat K_{\D,S}$).	By abstract nonsense the (co)slice is obtained as the 2-unstraightening of the 2-functor $\cat K_{\D,S} \to \Cat$, $X \mapsto \Hom_{\cat K_{\D,S}}(X, S) = \D(X)$. Again, we ignore the details on how to make this rigorous.
\end{itemize}
The above 2-categories fit in the following commuting diagram:
\begin{equation}\begin{tikzcd}
	\CoCorr_{\D,S} \arrow[r] \arrow[d] & \cat C'_{\D,S} \arrow[d] \\
	\Corr'(\cat C_{/S}) \arrow[r,"\Phi_{\D,S}"] & \cat K_{\D,S} \arrow[r,"\Psi_{\D,S}"] & \Cat
\end{tikzcd}\end{equation}
Here the vertical maps are cocartesian fibrations (in a 2-categorical sense) and the lower horizontal maps provide the corresponding functors to $\Cat$. In particular we expect the above square to be cartesian.
\end{remark}

\subsection{Applications to representation theory} \label{sec:intro.reptheory}

A perhaps somewhat surprising application of the 6-functor formalism from \cref{ex:intro-6ff-on-CondAni} is to smooth representation theory of locally profinite groups, by applying it to a classifying stack. This idea was brought to the second author's attention in \cite{Gulotta-Hansen-Weinstein:Enhanced-6FF-on-vStacks} and he explored it further in \cite{Hansen-Mann:Mod-p-Stacky-6FF,Mann.2022b}. In this paper we present this idea in full detail and use the category of kernels to push it further than before (while also fixing some gaps in existing proofs). We believe this approach to be especially powerful in natural characteristic (e.g.\ for smooth $\FF_p$-representations of $\GL_n(\QQ_p)$), where classical representation theoretic methods often fail. But even in the setting of $\ell$-adic smooth representation theory our constructions provide a new view on some classical definitions.

Let us now explain the promised application to representation theory in more detail. In the following, for a ring $\Lambda$ and a locally profinite group $G$ we denote by $\Rep_\Lambda(G)$ the \emph{derived} category of smooth representations of $G$ on $\Lambda$-modules, i.e.\ continuous $G$-representations on $\Lambda$-modules, where we implicitly equip the module with the discrete topology. We also denote by $\hatRep_\Lambda(G)$ the left-completion of $\Rep_\Lambda(G)$ (note that this does not change the full subcategory of left-bounded complexes). We can now state our first result (for the definition of the classifying stack $*/G$ we refer the reader to \cref{sec:alg.grp}):

\begin{proposition} \label{rslt:intro-identification-of-Rep-with-D-*-G}
Let $G$ be a locally profinite group and $\Lambda$ a (commutative) ring. There is a natural t-exact equivalence of categories
\begin{align}
	\D(*/G, \Lambda) = \hatRep_\Lambda(G),
\end{align}
compatibly with tensor products. Moreover, given a homomorphism $\varphi\colon H \to G$ of locally profinite groups with associated map $f\colon */H \to */G$ of classifying stacks, we have the following description of the adjunction $f^*\colon \D(*/G, \Lambda) \rightleftarrows \D(*/H, \Lambda) \noloc f_*$ in terms of smooth representations:
\begin{propenum}
	\item $f^*$ is the (derived) functor of taking a $G$-representation to the $H$-representation on the same underlying module, where $H$ acts via $\varphi$. This functor is often called \emph{restriction} or \emph{inflation}, depending on whether $\varphi$ is injective or surjective.

	\item If $\varphi$ is the inclusion of a closed subgroup, then $f_*$ is the right derived functor $\RInd_H^G$ of smooth induction. If $\varphi$ is a topological quotient map with kernel $U$ then $f_*$ is the right derived functor of taking $U$-invariants.
\end{propenum}
\end{proposition}
\begin{proof}
For the first statement see \cref{rslt:identification-of-sheaves-and-representations}. The hard part is to show that $\D(*/G, \Lambda)$ is the left-completed derived category of its heart, for which we develop some formalism of \enquote{derived descent from abelian descent} in \cref{sec:reptheory.descent} that may be of independent interest. The identification of $f^*$ and $f_*$ is straightforward, see \cref{rslt:*-functors-smooth-representations}.
\end{proof}

As an example, the pullback along the projection $*/G \to *$ sends a $\Lambda$-module to the trivial $G$-representation on it, while the pushforward along this map computes $G$-cohomology. \Cref{rslt:intro-identification-of-Rep-with-D-*-G} is an incarnation of the well-known paradigm that a category of representations can be realized as a category of sheaves on a classifying stack (this is valid for all kinds of representations). By its own this identification does not buy us much, but things get more interesting when we add $!$-functors. For the following result, we say that $G$ has \emph{finite $\Lambda$-cohomological dimension} if the group cohomology functor $(\blank)^G\colon \Rep_\Lambda(G) \to \D(\Lambda)$ has finite cohomological dimension, and we say that $G$ has \emph{locally finite $\Lambda$-cohomological dimension} if it has an open subgroup of finite $\Lambda$-cohomological dimension.

\begin{proposition} \label{rslt:intro-shriek-functors-for-classifying-stack}
Let $G$ be a locally profinite group and $\Lambda$ a ring such that $G$ has locally finite $\Lambda$-cohomological dimension. Then $\Rep_\Lambda(G)$ is left-complete, i.e.\ we have
\begin{align}
	\D(*/G,\Lambda) = \hatRep_\Lambda(G) = \Rep_\Lambda(G).
\end{align}
Moreover, if $H$ is another such group and $\varphi\colon H \to G$ is a homomorphism, then the associated map $f\colon */H \to */G$ is $\Lambda$-fine and hence we have an adjunction $f_!\colon \D(*/H,\Lambda) \rightleftarrows \D(*/G,\Lambda) \noloc f^!$. Moreover, the following is true:
\begin{propenum}
	\item Suppose $\varphi$ is the inclusion $H \subset G$ of a closed subgroup. If $H$ is open then $f$ is $\D$-étale and if $G/H$ is compact then $f$ is $\D$-proper. In any case, $f_!$ is the compact induction functor $\cInd_H^G$.

	\item Suppose that $H$ is compact and has finite $\Lambda$-cohomological dimension. Then $*/H \to *$ is $\D$-proper.

	\item \label{rslt:intro-class-stack-suave-if-Poincare} Suppose that $\varphi$ is a topological quotient map with kernel $U$ and assume that $U$ is \emph{locally $\Lambda$-Poincaré}, i.e.\ for an open profinite subgroup $K \subseteq U$ with finite $\Lambda$-cohomological dimension, the $K$-cohomology functor $(\blank)^K\colon\Rep_\Lambda(K) \to \D(\Lambda)$ preserves dualizable objects. Then $f$ is $\D$-suave.
\end{propenum}
\end{proposition}
\begin{proof}
For the left-completeness statement see \cref{rslt:Rep-G-left-complete}; the crucial observation is that under the assumptions on $G$ and $\Lambda$, products in $\Rep_\Lambda(G)$ have finite cohomological dimension. For the $\Lambda$-fineness claim, see \cref{rslt:maps-of-class-stacks-are-fine}. The claim easily reduces to the case $\varphi\colon H \to *$ such that $H$ has finite $\Lambda$-cohomological dimension. In this case we analyze the descent along the map $* \to */H$ and in particular show that it has universal $!$-descent using the primness criterion for $!$-descent \cref{rslt:!-descent-along-prim-morphism}. By the locality criterion for properness on the source (see \cref{rslt:check-properness-on-!-descendable-cover}) we deduce that $*/H \to *$ is $\D$-proper; this also proves (ii).

The étaleness and properness claims in (i) follow from locality on the target and $*/H \times_{*/G} * = G/H$. For (iii) see \cref{rslt:suave-equiv-to-locally-Poincare}.	
\end{proof}

The condition of being locally $\Lambda$-Poincaré can be verified in many examples of interest, for a fixed prime $p$:
\begin{itemize}
	\item If $\Lambda$ is a $\ZZ[1/p]$-algebra and $G$ is locally pro-$p$, then $G$ is locally $\Lambda$-Poincaré (see \cref{example:dualizing-complex-away-from-p}.
	\item If $\Lambda$ is a $\ZZ/p^n$-algebra and $G$ is a $p$-adic Lie group, then $G$ is locally $\Lambda$-Poincaré (see \cref{example:dualizing-complex-mod-p}).
\end{itemize}
In both cases we have an explicit description of the dualizing complex, which is in particular invertible (i.e.\ given by a character of $G$). Note that by \cref{rslt:intro-class-stack-suave-if-Poincare} in both cases we immediately deduce Poincaré duality for $G$-cohomology. We point out that this does not really give a new proof of Poincaré duality, because the proofs in the above examples rely on results in classical representation theory which are close to Poincaré duality already. Note that the functor $f_!$ also recovers the derived coinvariants functor constructed by the first author in the mod $p$ case (see \cref{rmk:Heyer-coinvariants-functor}).

We have defined and explained all six functors on classifying stacks and have seen that they recover many important constructions from representation theory (some of which are not easy to perform by hand in the mod $p$ case). The next result shows that suave and prim objects also have a tight relation to classical representation theoretic notions:

\begin{proposition}
Let $G$ be a locally profinite group and $\Lambda$ a ring such that $G$ has locally finite $\Lambda$-cohomological dimension. Let $f\colon */G \to *$ denote the structure map and let $V \in \D(*/G,\Lambda)$ be a representation.
\begin{propenum}
	\item \label{rslt:intro-suave-duality-on-Rep-G}$V$ is $f$-suave if and only if it is \emph{admissible}, i.e.\ $V^K \in \D(\Lambda)$ is perfect for small enough compact open subgroups $K \subseteq G$ (see \cref{def:admissible-representation}). In this case the $f$-suave dual is given by
	\begin{align}
		\DSuave_f(V) = \iHom(V, f^! \one).
	\end{align}

	\item \label{rslt:intro-prim-duality-on-Rep-G} $V$ is $f$-prim if and only if it is compact. In this case the $f$-prim dual is given by
	\begin{align}
		\DPrim_f(V) = \RHom_G(V, \Cont_c(G, \Lambda)),
	\end{align}
	where the $\RHom_G$ on the right denotes derived $G$-equivariant maps for the right translation action of $G$ on $\Cont_c(G,\Lambda)$. If $K \subseteq G$ is a compact open subgroup with finite $\Lambda$-cohomological dimension and $V$ is a dualizable smooth $K$-representation then
	\begin{align}
		\DPrim_f(\cInd_K^G V) = \cInd_K^G V^\vee.
	\end{align}
\end{propenum}
\end{proposition}
\begin{proof}
This is an easy application of the general results on suave and prim objects, see \cref{rslt:identify-suave-and-prim-representations}. We explain the proof of the formula for prim duality of $\cInd_K^G V$ here, as we find it quite enlightening. Consider the maps $g\colon */K \to *$ and $h\colon */K \to */G$, so that by \cref{rslt:intro-shriek-functors-for-classifying-stack} $g$ is $\D$-proper and $h$ is $\D$-étale. Since also the diagonal of $g$ is proper (its fiber is $K$), by \cref{rslt:relation-between-suave-prim-and-dualizable} the dualizable $K$-representations are exactly the $g$-prim ones, and $g$-prim duality is just the usual duality. Now by \cref{rslt:suave-pullback-and-prim-pushforward-preservations} $h_! = \cInd_K^G$ preserves primness and commutes with prim duality, giving us the claim.
\end{proof}

If $p$ is a prime, $G$ is locally pro-$p$ and $\Lambda$ is a $\ZZ[1/p]$-algebra, then the prim duality functor has been studied classically, under the name of \enquote{cohomological duality} (see \cite[\S5.1]{Bernstein:Lectures} or \cite[\S III]{Schneider-Stuhler.1997}). There is a version of this duality in the geometrization of the local Langlands correspondence (see \cite[\S V.5]{Fargues-Scholze:Geometrization}), which is another incarnation of prim duality (this will be part of a joint work by the second author and Hansen). In the case that $\Lambda$ is an $\FF_p$-algebra, a version of Bernstein's cohomological duality is less understood in the literature, but appeared in \cite{Emerton-Gee-Hellmann:IHES}, see the explanations following Remark~6.1.6. Prim duality provides a conceptual unifying approach to all of these constructions.

Using prim duality we obtain new results for derived Hecke algebras. In \cite{Schneider-Sorensen.2023b} the authors construct a canonical anti-involution on the pro-$p$ Iwahori--Hecke $\Ext$-algebra in analogy with a similar anti-involution for (underived) Hecke algebras. However, it remained unclear whether the construction lifts to the level of dg-algebras, cf.\ \cite[Remark~2.5]{Schneider-Sorensen.2023b}. We provide a positive answer:

\begin{theorem} \label{rslt:intro-antiinvolution}
Let $\Lambda$ be a field of characteristic $p > 0$, $G$ a $p$-adic Lie group and $I \subseteq G$ a $p$-torsionfree compact open subgroup. Denote by
\begin{align}
	\Hecke_I^\bullet \coloneqq \RHom_G(\cInd_K^G \one, \cInd_K^G \one) \in \Alg(\D(\Lambda))
\end{align}
the \emph{derived Hecke algebra} of $G$. Then there is a canonical anti-involution
\begin{align}
	\Inv\colon (\Hecke_I^\bullet)^\op \isoto \Hecke_I^\bullet,
\end{align}
which coincides with the one in \cite{Schneider-Sorensen.2023b} after passage to the cohomology algebras.
\end{theorem}
\begin{proof}
By \cref{rslt:intro-prim-duality-on-Rep-G} prim duality induces an equivalence $(\Rep_\Lambda(G)^\omega)^\op \isoto \Rep_\Lambda(G)^\omega$ sending $\cInd_K^G \one$ to itself. Evaluating this equivalence on the ($\Lambda$-enriched) endomorphisms of $\cInd_K^G$ provides the result. See \cref{rslt:anti-involution-comparison-with-Schneider-Sorensen} for details and for the comparison with the construction in \cite{Schneider-Sorensen.2023b}.
\end{proof}

We remark that the restriction on $\Lambda$ being a field of characteristic $p$ is not necessary in \cref{rslt:intro-antiinvolution}. We only require that $G$ is locally $\Lambda$-Poincaré and $I$ is compact with finite $\Lambda$-cohomological dimension. In particular the same technique constructs an anti-involution on the Hecke algebra away from $p$, which coincides with the classical construction.

\subsection{Notation and conventions}

As mentioned in the introduction, this paper is rooted in higher categories and so we decided to drop the \enquote{$\infty$} from the notation, i.e.\ we refer to $\infty$-categories as categories and to $(\infty,2)$-categories as 2-categories. We refer the reader to the appendix for an introduction to these notions. In the same spirit, we will implicitly work in the derived world throughout, i.e.\ we will denote by $f^*$, $f_*$, $\Gamma(X,\blank)$, etc. the derived versions of these functors and we will refer to complexes of modules/representations simply as modules/representations. In contrast, by a \enquote{ring} we will mean an ordinary ring (to the $\infty$-categorical versions we will refer to as $\Einfty$-rings and animated rings). In order to indicate that a derived object lies in the heart, we will call it \emph{static}, e.g.\ a static $\Lambda$-module is just an ordinary $\Lambda$-module (for a ring $\Lambda$). Moreover, we use the following notation:
\begin{enumerate}[label=--]
\item $\bbDelta$ is the simplex category. The objects are finite, non-empty, linearly ordered sets together with order preserving maps. Every object of $\bbDelta$ is isomorphic to $[n] \coloneqq \{0\le 1\le \dotsb\le n\}$, for some $n\ge0$ (see \cref{sec:quasicategories}).

\item $\Ani$ is the category of anima (see \cref{ex:category-of-anima}).

\item For a category $\cat C$ we denote by $\PSh(\cat C) = \Fun(\cat C^{\op}, \Ani)$ the category of presheaves on $\cat C$.

\item $\Cat$ is the (2-)category of categories (see \cref{ex:category-of-categories}, and \cref{ex:2-categorical-enhancement-of-Cat}), $\Cat_2$ is the (2\nobreakdash-)category of 2-categories (see \cref{def:2-category-of-2-categories}).

\item For a 2-category $\cat C$ with objects $X,Y$ we denote by $\Fun_{\cat C}(X,Y)$ the category of morphisms $X\to Y$ in $\cat C$ (see \cref{def:Hom-category-in-2-Cat}).

\item $\Fin_*$ is the category of finite pointed sets. For $n\ge0$, we write $\langle n\rangle^{\circ} \coloneqq \{1,2,\dotsc,n\}$ and $\langle n\rangle \coloneqq \{*\}\sqcup \langle n\rangle^{\circ}$ (see \cref{def:category-of-finite-pointed-sets}).

\item $\Comm^\tensor$ is the commutative operad, $\Ass^\tensor$ is the associative $\infty$-operad and $\LM^\tensor$ is the operad of left modules (see \cref{def:basic-operads}).

\item Denote $\Cat^\times$ the Cartesian operad with underlying category $\Cat$ (see \cref{ex:cartesian-monoidal-structure}). 

\item $\CMon$ is the category of symmetric monoidal categories and symmetric monoidal functors, $\Mon$ is the category of monoidal categories and monoidal functors, see \cref{def:monoidal-categories}.

\item We denote by $\CAlg$ the category of $\Einfty$-ring spectra, which contains classical (commutative) rings as a full subcategory. For an $\Einfty$-ring $\Lambda$ we denote by $\Mod_{\Lambda} = \D(\Lambda)$ the (derived) category of $\Lambda$-modules (see \cref{ex:Einfty-rings-and-modules}).

\item For a monoidal category $\cat V$ we denote by $\Enr_{\cat V}$ the category of $\cat V$-enriched categories (see \cref{def:category-of-Lurie-enriched-categories}).
\end{enumerate}

\subsection{Acknowledgements}

We thank Thomas Nikolaus, David Hansen, Juan Esteban Rodríguez Camargo, Eugen Hellmann, Peter Schneider, Claus Sorensen, Zhixiang Wu and Fabian Hebestreit for enlightening discussions during the preparation of this article. 
We thank Peter Scholze, David Hansen, João Lourenço and Juan Esteban Rodríguez Camargo for many helpful corrections and comments on a first draft.
Some of the ideas in this paper have already appeared in joint work of the second author with David Hansen and in Peter Scholze's lecture notes on 6-functor formalisms. We are particularly grateful to David Hansen, whose preliminary work on the Liu--Zheng papers started the whole idea of abstract 6-functor formalisms, and discussions with whom have lead to powerful new criteria for detecting suave and prim sheaves.
We furthermore thank Peter Scholze and David Hansen for suggesting the names \enquote{suave} and \enquote{prim}.

Some of the results of this article were presented at the \enquote{Workshop on Shimura varieties, representation theory and related topics} held in October 2024 in Tokyo, and C.H.\ wants to thank the organizers, particularly Noriyuki Abe for the invitation.

\section{The category of correspondences} \label{sec:Corr}

In this section we will introduce the main ingredient for 6-functor formalisms: the category of correspondences. We start in \cref{subsec:geometric-setups} with a discussion of geometric setups. From each geometric setup $(\cat C, \cat C_\gsindex)$ we will then construct in \cref{subsec:definition-of-Corr} the category of correspondences $\Corr(\cat C, \cat C_\gsindex)$. If $\cat C$ admits finite products, then $\Corr(\cat C, \cat C_\gsindex)$ can be endowed with the structure of a closed symmetric monoidal category, which is the content of \cref{subsec:monoidal-structure-on-Corr,subsec:self-enrichment-of-Corr}. For a detailed discussion of the category of correspondences in the more general context of adequate triples we refer to \cite{Haugseng-Hebestreit.Spans}.

\subsection{Geometric setups}\label{subsec:geometric-setups}
In this subsection we recall the definition of geometric setups $(\cat C, \cat C_{\gsindex})$ and give first examples. 

\begin{definition}\label{defn:geometric-setup}
A \emph{geometric setup} is a pair $(\cat C, \cat C_{\gsindex})$ consisting of a category $\cat C$ and a wide subcategory $\cat C_{\gsindex} \subseteq \cat C$ satisfying the following properties:
\begin{defenum}
\item\label{defn:geometric-setup-base-change} $\cat C_{\gsindex}$ is closed under pullbacks along morphisms in $\cat C$,
\item\label{defn:geometric-setup-fiber-products} $\cat C_{\gsindex}$ admits and the inclusion $\cat C_{\gsindex} \subseteq \cat C$ preserves fiber products.
\end{defenum}
We denote $\GeomSetup$ the (non-full) subcategory of $\Fun([1],\Cat)$ classifying geometric setups $\cat C_\gsindex \subseteq \cat C$, where
\[
\Hom_{\GeomSetup}((\cat C,\cat C_\gsindex), (\cat D,\cat D_\gsindex)) \subseteq \bigl(\Fun(\cat C, \cat D) \times_{\Fun(\cat C_\gsindex,\cat D)} \Fun(\cat C_\gsindex, \cat D_\gsindex)\bigr)^\simeq
\]
is the full subanima spanned by those pairs $(f,f_\gsindex)$ such that
$f$ preserves the pullbacks of $X\xrightarrow{g} Y \from Z$ in $\cat C$ with $g$ in $\cat C_\gsindex$.
Note that the forgetful map $\Hom_{\GeomSetup}((\cat C,\cat C_\gsindex), (\cat D,\cat D_\gsindex)) \to \Fun(\cat C, \cat D)^\simeq$ is a full subanima, meaning that the pair $(f,f_\gsindex)$ is already determined by $f$. Note also that $f_0$ automatically preserves fiber products.
\end{definition}

\begin{remark}\label{rmk:geometric-setup-morphisms}
Observe that $\cat C_{\gsindex}$ contains all objects and isomorphisms in $\cat C$. Hence, the datum of a geometric setup is equivalent to the datum of a category $\cat C$ together with a (homotopy) class of edges $E$ in $\cat C$ which 
\begin{enumerate}[label=(\roman*)]
\item contains all isomorphisms, 
\item is closed under composition and pullbacks along edges in $\cat C$, and
\item for every $Y\to X$ in $E$, the diagonal $Y\to Y\times_XY$ lies in $E$ (see \cref{rslt:right-cancellative-fiber-products} below for the equivalence of this condition with (b) in \cref{defn:geometric-setup}).
\end{enumerate}
This point of view is taken in \cite{Mann.2022a} and \cite{Scholze:Six-Functor-Formalism}, where the condition (iii) is dropped. The reasons why we include (iii) in the definition of a geometric setup are that
\begin{enumerate*}
\item it is often necessary in order to construct 6-functor formalisms, see \cref{sec:6ff.construct},
\item it will be needed throughout \cref{sec:kerncat}, where we discuss the $2$-category of kernels and its functoriality properties, and
\item it is basically always satisfied in practice anyway. 
\end{enumerate*}
However, we point out that (iii) is not essential for the results in \cref{sec:Corr}.
\end{remark}

We make the following convention:

\begin{convention}\label{convention:geometric-setup-morphisms}
We will, whenever it is convenient, denote a geometric setup by $(\cat C, E)$, where $E$ is a collection of morphisms of $\cat C$ satisfying the conditions in \cref{rmk:geometric-setup-morphisms}. In this case, we denote $\cat C_E$ the wide subcategory of $\cat C$ with morphisms $E$.
\end{convention}

\begin{example}
Let $\cat C$ be a category.
\begin{exampleenum}
\item The pair $(\cat C, \cat C^\simeq)$ is a geometric setup, where $\cat C^\simeq$ is the underlying anima of $\cat C$.

\item If $\cat C$ admits fiber products, then $(\cat C, \cat C)$ is a geometric setup.

\item If $(\cat C, \cat C_{\gsindex})$ is a geometric setup, then so too is $(\cat C_{\gsindex}, \cat C_{\gsindex})$ and the map $(\cat C_{\gsindex}, \cat C_{\gsindex}) \to (\cat C, \cat C_{\gsindex})$ is a morphism of geometric setups.
\end{exampleenum}
\end{example}

\begin{lemma}\label{rslt:right-cancellative-fiber-products}
Let $\cat C$ be a category and $E$ a collection of morphisms in $\cat C$ satisfying (i) and (ii) in \cref{rmk:geometric-setup-morphisms}. The following conditions are equivalent:
\begin{lemenum}
\item $\cat C_E$ satisfies condition (b) in \cref{defn:geometric-setup}, i.e.\ $\cat C_E$ admits and the inclusion $\cat C_E \subseteq \cat C$ preserves fiber products;

\item \label{rslt:E-closed-under-diagonals} $E$ satisfies (iii) in \cref{rmk:geometric-setup-morphisms}, i.e.\ for every morphism $Y\to X$ in $E$, the diagonal morphism $Y\to Y\times_XY$ lies in $E$;

\item \label{rslt:E-is-cancellative} $E$ is \emph{right cancellative}, i.e.\ given composable morphisms $Z\xto{g} Y \xto{f} X$ such that $f, f\comp g \in E$, then $g\in E$.
\end{lemenum}
\end{lemma}
\begin{proof}
The implication (i) $\To$ (ii) is obvious.

We now prove that (ii) implies (iii). So let $Z\xto{g} Y \xto{f} X$ be composable morphisms such that $f, f\comp g \in E$. We need to show $g\in E$. Since $g$ is the composition of $Z \xto{(g,\id)} Y\times_XZ \xrightarrow{\pr_1} Y$, it suffices to show $(g,\id), \pr_1 \in E$. But note that $\pr_1$ is the pullback of $f\comp g \in E$ along $g$, so that $\pr_1 \in E$. Similarly, $(g,\id)$ is the pullback of the map $\id\times g \colon Y\times_XZ \to Y\times_XY$ along the diagonal map $Y \to Y\times_XY$, which lies in $E$ by (ii) and the assumption that $f\in E$. Therefore, we have $(g,\id) \in E$, which finishes the proof that $g\in E$.

Finally, we show that (iii) implies (i). Given two morphisms $Y\to X \from Z$ in $\cat C_E$, the projections $\pr_1\colon Y\times_XZ \to Y$ and $Y\times_XZ \to Z$ lie in $E$. In order to verify that $Y\times_XZ$ is also a fiber product in $\cat C_E$, consider a commutative diagram with solid arrows
\begin{equation}
\begin{tikzcd}
W \ar[ddr, bend right, "u"'] \ar[rrd, bend left] \ar[dr,dashed, "\exists! h"] \\[-1em]
& Y\times_XZ \ar[r] \ar[d,"\pr_1"'] & Z \ar[d] \\
& Y \ar[r] & X
\end{tikzcd}
\end{equation}
in $\cat C_E$. There exists a unique map $h\colon W \to Y\times_XZ$ in $\cat C$ making the whole diagram commutative. It suffices to prove $h\in E$. Since $E$ is right cancellative and $u, \pr_1 \in E$, we deduce that $h$ lies in $E$. This shows that $Y\times_XZ$ is a fiber product in $\cat C_E$.
\end{proof}

\begin{lemma}\label{rslt:GeomSetup-has-limits}
$\GeomSetup$ is complete and the functor $\GeomSetup \to \Fun([1],\Cat)$ preserves small limits.
\end{lemma}
\begin{proof}
Let $\cat I$ be a small category and $F\colon \cat I \to \GeomSetup$, $i\mapsto (\cat C_i, \cat C_{i,\gsindex})$ a diagram. Put $\cat C \coloneqq \lim_i \cat C_i$ and $\cat C_\gsindex \coloneqq \lim_i \cat C_{i,\gsindex}$ in $\Cat$. We need to show that $(\cat C, \cat C_\gsindex)$ is a geometric setup and satisfies the universal property of a limit in $\GeomSetup$.

First note that $\cat C_\gsindex$ is a wide subcategory of $\cat C$ by \cref{rslt:limit-of-subcategories}. (The objects are the same, since each $\cat C_{i,\gsindex}$ is a wide subcategory of $\cat C_i$ containing all isomorphisms.) Consider a diagram $X \xto{g} Y \from Z$ with $g\in \cat C_\gsindex$, and denote its image in $\cat C_i$ by $X_i \xto{g_i} Y_i \from Z_i$. Note that each $g_i$ lies in $\cat C_{i,\gsindex}$ and for any edge $i\to j$ in $\cat I$, the map $\cat C_i\to \cat C_j$ sends $X_i\times_{Y_i}Z_i$ to $X_j\times_{Y_j}Z_j$, since $(\cat C_i,\cat C_{i,\gsindex}) \to (\cat C_j, \cat C_{j,\gsindex})$ is a morphism of geometric setups. The pullback diagrams 
\begin{equation}
\begin{tikzcd}
X_i\times_{Y_i} Z_i \ar[d] \ar[r] \ar[dr,phantom, very near start, "\lrcorner"] & 
Z_i \ar[d] \\
X_i \ar[r] & Y_i
\end{tikzcd}
\quad
\text{assemble into a pullback diagram}
\quad
\begin{tikzcd}
X\times_YZ \ar[d] \ar[r] \ar[dr,phantom, very near start, "\lrcorner"] & Z \ar[d] \\
X \ar[r] & Y.
\end{tikzcd}
\end{equation}
Thus, $U \in \cat C$ is a pullback of $X\to Y \from Z$ if and only if its image in $\cat C_i$ is a pullback of $X_i\to Y_i \from Z_i$, for all $i$. 

We deduce that the projections $(\cat C, \cat C_\gsindex) \to (\cat C_i, \cat C_{i,\gsindex})$ are morphisms of geometric setups, and that for any map $\const_{(\cat A, \cat A_\gsindex)} \to F$ in $\Fun(\cat I, \GeomSetup)$, the induced map $(\cat A, \cat A_\gsindex) \to (\cat C,\cat C_\gsindex)$ (in $\Fun([1],\Cat)$) is a morphism of geometric setups. Therefore, $(\cat C, \cat C_\gsindex)$ is the limit of $F$ in $\GeomSetup$.
\end{proof}

\subsection{Definition of the correspondence category}\label{subsec:definition-of-Corr} 
In this section we associate with each geometric setup $(\cat C, \cat C_\gsindex)$ the category $\Corr(\cat C, \cat C_\gsindex)$ of correspondences. The correspondence category has appeared in the literature under several names. It was introduced by Barwick in \cite{Barwick-Rognes.2013}, where it is called the \emph{$Q$-construction}, and further studied in the series of papers \cite{Barwick.2017, Barwick-Glasman-Shah.2019} under the name \emph{effective Burnside category}. Our exposition was partly inspired by the very readable treatment in \cite{Haugseng.2017a} (who uses the more common name \emph{category of spans}).

\begin{definition} \label{def:Sigma-n-and-Lambda-n}
\begin{defenum}
\item We denote $\bbSigma^n$ the poset of pairs $(i,j)$ with $0\le i\le j\le n$, where the partial ordering is given by $(i,j) \le (k,l)$ if $i\le k$ and $l\le j$. We obtain a cosimplicial category
\[
\bbSigma^\bullet \colon \bbDelta \to \Cat_{\ordinary} \subseteq \Cat.
\]

\item Let $\bbSigma_2^n \subseteq \bbSigma^n$ be the subposet with the same elements but with partial ordering given by $(i,j)\le (k,l)$ if $i\le k$ and $j=l$. This defines a cosimplicial category $\bbSigma^\bullet_2 \colon \bbDelta \to \Cat$ and a natural functor $\bbSigma_2^\bullet \to \bbSigma^\bullet$.

\item We denote $\bbLambda^n \subseteq \bbSigma^n$ the full subcategory consisting of those pairs $(i,j)$ with $j-i \le 1$. Note that the functor $\bbSigma^n \to \bbSigma^m$ induced by $\phi \colon [n] \to [m]$ does not restrict to $\bbLambda^n$ unless $\phi$ is the inclusion of an interval. We thus obtain a functor $\bbLambda^\bullet\colon \bbDelta_{\inert} \to \Cat$ and a natural embedding $\bbLambda^\bullet \injto \bbSigma^\bullet|_{\bbDelta_{\inert}}$, where $\bbDelta_\inert \subseteq \bbDelta$ is the wide subcategory with interval inclusions as morphisms.
\end{defenum}
\end{definition}

\begin{figure}[ht]
\begin{center}
\begin{tikzcd}[sep = small]
& & & (0,3) \ar[dl,dashed] \ar[dr,dashed] \\
& & (0,2) \ar[dl,dashed] \ar[dr,dashed] & & (1,3) \ar[dl,dashed] \ar[dr,dashed] \\
& (0,1) \ar[dl] \ar[dr] & & (1,2) \ar[dl] \ar[dr] & & (2,3) \ar[dl] \ar[dr] \\
(0,0) & & (1,1) & & (2,2) & & (3,3)
\end{tikzcd}
\end{center}
\caption{The diagram depicts $\bbSigma^n$ (all arrows) and $\bbLambda^n$ (solid arrows) for $n=3$. The poset $\bbSigma_2^3$ contains only the arrows pointing to the right.}
\label{figure:Sigma_Lambda}
\end{figure}
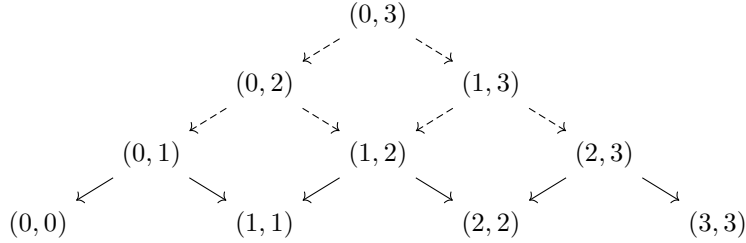

\begin{example}\label{ex:geometric-setups}
Let $n\ge0$.
\begin{exampleenum}
\item The pair $(\bbSigma^n, \bbSigma_2^n)$ is a geometric setup which is functorial in $[n]\in \bbDelta$. We obtain a cosimplicial object $(\bbSigma^\bullet, \bbSigma_2^\bullet) \colon \bbDelta \to \GeomSetup$. 

\item The pair $([n]^\op, [n]^\simeq)$ is a geometric setup which is functorial in $[n]\in \bbDelta$. We obtain a cosimplicial object $([\blank]^\op, [\blank]^\simeq) \colon \bbDelta \to \GeomSetup$. The projection $(i,j) \mapsto j$ induces a map $(\bbSigma^\bullet, \bbSigma_2^\bullet) \to ([\blank]^\op, [\blank]^\simeq)$ of cosimplicial geometric setups.

\item The pair $([n], [n])$ is a geometric setup which is functorial in $[n]\in \bbDelta$. We obtain a cosimplicial object $([\blank], [\blank]) \colon \bbDelta \to \GeomSetup$. The projection $(i,j)\mapsto i$ induces a map $(\bbSigma^\bullet, \bbSigma_2^\bullet)\to ([\blank], [\blank])$ of cosimplicial geometric setups.
\end{exampleenum}
\end{example}

\begin{definition}
The composite
\begin{align}
\bbDelta^\op \times \GeomSetup \xrightarrow{(\bbSigma^\bullet, \bbSigma_2^\bullet) \times \id} \GeomSetup^\op \times \GeomSetup \xrightarrow{\Hom} \Ani
\end{align}
defines a functor
\begin{align}\label{eq:functor-GeomSetup-to-sAni}
\Corr\colon \GeomSetup \too \Fun(\bbDelta^\op, \Ani)
\end{align}
sending a geometric setup $(\cat C,\cat C_\gsindex)$ to the simplicial anima $[n] \mapsto \Hom_{\GeomSetup}((\bbSigma^n,\bbSigma_2^n), (\cat C, \cat C_\gsindex))$. See also \cref{def:Correspondence-Category} below for a more explicit description, which also explains the notation.
\end{definition}

Our goal in the following is to show that $\Corr$ factors through the full subcategory of complete Segal anima and hence defines a functor to $\Cat$, see \cref{sec:Segal-anima}. This can be checked on objects, so let us fix a geometric setup $(\cat C,\cat C_\gsindex)$. It will be helpful to describe a morphism $(\bbSigma^n, \bbSigma_2^n) \to (\cat C, \cat C_\gsindex)$ as a condition on the underlying functor $\bbSigma^n\to \cat C$.

\begin{definition}\label{def:cartesian-map}
Let $\cat A$ be a category. A functor $f\colon \bbSigma^n \to \cat A$ is called \emph{cartesian} if $f$ is a pointwise right Kan extension of its restriction to $\bbLambda^n$. We denote $\Fun^{\cart}(\bbSigma^n,\cat A)$ the full subcategory of cartesian functors and observe that the restriction map $\Fun^{\cart}(\bbSigma^n, \cat A) \to \Fun(\bbLambda^n,\cat A)$ is fully faithful (by \cite[Proposition~4.3.2.15]{HTT}).
\end{definition}

\begin{remark}\label{rmk:cartesian-map}
Since right Kan extensions can be computed in stages, the following assertions for a functor $f\colon \bbSigma^n\to \cat C$ are equivalent:
\begin{remarksenum}
\item\label{rmk:cartesian-map-i} $f$ is cartesian;

\item\label{rmk:cartesian-map-ii} each square
\begin{equation}
\begin{tikzcd}
f(i,j) \ar[d] \ar[r] \ar[dr,phantom,very near start, "\lrcorner"] & f(l,j) \ar[d] \\
f(i,k) \ar[r] & f(l,k)
\end{tikzcd}
\end{equation}
is a pullback in $\cat C$, for all $i\le l\le k\le j$;

\item\label{rmk:cartesian-map-iii} each square as in (ii) with $(l,k) = (i+1,j-1)$ is a pullback in $\cat C$.
\end{remarksenum}
The fact that (iii) implies (ii) is an easy application of the pasting law, see~\cite[\href{https://kerodon.net/tag/03FZ}{Tag 03FZ}]{kerodon}.

We argue that (iii) is equivalent to (i): For any $1\le a \le n$ and $0\le b \le n-a$, let $\bbSigma^{n,a,b} \subseteq \bbSigma^n$ be the full subcategory of $(i,j)$ such that either $j-i < a$ or $j-a = i \le b$. We observe:
\begin{itemize}
\item $\bbSigma^n = \bbSigma^{n,n,0}$ and $\bbLambda^n = \bbSigma^{n,1,n-1}$;

\item $\bbSigma^{n,a,b} = \bbSigma^{n,a,b-1} \sqcup \{(b,a+b)\}$ for all $0\le b\le n-a$, where we make the convention $\bbSigma^{n,a,-1}\coloneqq \bbSigma^{n,a-1,n-a+1}$.
\end{itemize}
Now note that a functor $f\colon \bbSigma^{n,a,b} \to \cat C$ is a right Kan extension of $f|_{\bbSigma^{n,a,b-1}}$ if and only if
\begin{align}\label{eq:cartesian-map-pullback}
f(b,a+b) &\isoto \lim_{(i,j) \in \bbSigma^{n,a,b-1}_{(b,a+b)/}} f(i,j).
\end{align}
We observe $\bbSigma^{n,a,b-1}_{(b,a+b)/} \simeq \bbSigma^{a,a-1,1}$ and that the inclusion $\Horn^2_2= \{(0,a-1), (1,a-1), (1,a)\} \injto \bbSigma^{a,a-1,1}$ is cofinal as follows immediately from Quillen's Theorem~A \cite[\href{https://kerodon.net/tag/02P0}{Tag 02P0}]{kerodon} and the observation that
\[
\{(0,a-1), (1,a-1), (1,a)\}_{/(i,j)} \simeq \begin{cases}
*, & \text{if $i=0$ or $j=a$,}\\
\Horn^2_2, & \text{else,}
\end{cases}
\]
is contractible (since $\Horn^2_2$ contains a terminal element)
for all $(i,j) \in \bbSigma^{a,a-1,1}$. Hence, \eqref{eq:cartesian-map-pullback} holds if and only if $f(b,a+b)$ is the pullback of the diagram $f(b,a+b-1) \to f(b+1,a+b-1) \from f(b+1,a+b)$.
\end{remark}

\begin{definition}
Let $n\ge0$.
\begin{defenum}
\item We denote by $\Fun^{\cart, \cat C_\gsindex}(\bbSigma^n, \cat C)^\simeq \subseteq \Fun(\bbSigma^n,\cat C)^\simeq$ the full subanima of cartesian functors $f \colon \bbSigma^n \to \cat C$ such that all $f(i_1,j) \to f(i_2,j)$ with $i_1\le i_2$ lie in $\cat C_{\gsindex}$. 

\item We denote $\Fun^{\cat C_\gsindex}(\bbLambda^n,\cat C)^\simeq\subseteq \Fun(\bbLambda^n,\cat C)^\simeq$ the full subanima of functors $f$ such that all maps $f(i-1,i) \to f(i,i)$ lie in $\cat C_\gsindex$.
\end{defenum}
\end{definition}

\begin{lemma}\label{rslt:Hom_Corr-identifications}
For any $n\ge0$, the restriction maps
\[
\Hom((\bbSigma^n, \bbSigma_2^n), (\cat C, \cat C_\gsindex)) \isoto \Fun^{\cart,\cat C_\gsindex}(\bbSigma^n, \cat C)^\simeq \isoto \Fun^{\cat C_\gsindex}(\bbLambda^n,\cat C)^\simeq
\]
are isomorphisms of anima. 
\end{lemma}
\begin{proof}
The assertion follows easily from \cref{rmk:cartesian-map}.
\end{proof}

\begin{proposition}\label{rslt:Corr-complete-Segal-anima}
The functor $\Corr\colon \GeomSetup\to \Fun(\bbDelta^\op, \Ani)$ factors through the full subcategory of complete Segal anima.
\end{proposition}
\begin{proof}
Let $(\cat C, \cat C_\gsindex)$ be a geometric setup. 
We first show that $\Corr(\cat C, \cat C_\gsindex)$ is a Segal anima. Recall that a simplicial anima $X \in \Fun(\bbDelta^\op, \Ani)$ is a Segal anima if for each $n\ge2$ the restriction map
\begin{align}\label{eq:Segal-anima}
\Hom(\Delta^n, X) \isoto \Hom(I_n, X) \simeq X_1 \times_{X_0} X_1 \times_{X_0}\dotsb \times_{X_0} X_1
\end{align}
is an isomorphism, where $I_n \subseteq \Delta^n$ is the $n$-th spine and we view $I_n$ and $\Delta^n$ as simplicial anima via the embedding $\Set\subseteq \Ani$.

By Lemma~\ref{rslt:Hom_Corr-identifications} we have $\Hom(\Delta^n, \Corr(\cat{C}, \cat{C}_\gsindex)) = \Fun^{\cart, \cat C_\gsindex}(\bbSigma^n, \cat C)^\simeq$ and $\Hom(I_n, \Corr(\cat C, \cat C_\gsindex)) = \Fun^{\cat C_\gsindex}(\bbLambda^1,\cat C)^\simeq \times_{\cat C} \dotsb \times_{\cat C} \Fun^{\cat C_\gsindex}(\bbLambda^1,\cat C)^\simeq = \Fun^{\cat C_\gsindex}(\bbLambda^n, \cat C)^\simeq$, where the last identification uses the easy fact that $\bbLambda^n = \bbLambda^1 \dunion_{\bbLambda^0} \dotsb \dunion_{\bbLambda^0} \bbLambda^1$ in $\Cat$, see~\cite[Proposition~5.13]{Haugseng.2017a}. This means that we have to check that the restriction map
\[
\Fun^{\cart, \cat C_{\gsindex}}(\bbSigma^n, \cat C)^\simeq \to \Fun^{\cat{C}_\gsindex}(\bbLambda^n, \cat C)^\simeq
\]
is an isomorphism, which is \cref{rslt:Hom_Corr-identifications}. Hence, $\Corr(\cat C, \cat C_\gsindex)$ is a Segal anima.

The completeness is proved as in \cite[Proposition~8.1]{Haugseng.2017a}.
\end{proof}

\begin{definition}\label{def:Correspondence-Category}
Let $(\cat C, \cat C_\gsindex)$ be a geometric setup. The complete Segal anima $\Corr(\cat C, \cat C_\gsindex)$ can be identified with a category (see \cref{rslt:Segal-anima-categories}) which we call the \emph{category of correspondences associated with $(\cat C, \cat C_\gsindex)$}. If $\cat C_\gsindex = \cat C$, we simply write $\Corr(\cat C)$ instead of $\Corr(\cat C, \cat C)$.

In concrete terms $\Corr(\cat C, \cat C_\gsindex)$ can be described as follows:
\begin{defenum}
\item the objects are the objects of $\cat C$;

\item a morphism from $X$ to $Y$ is a correspondence
\begin{equation}
\begin{tikzcd}[column sep=small]
& V \ar[dl,"f"'] \ar[dr,"g"] \\
X & & Y,
\end{tikzcd}
\end{equation}
where $g$ lies in $\cat C_\gsindex$;

\item the composite of two correspondences $[X\xleftarrow{f} V\xrightarrow{g}Y]$ and $[Y\xleftarrow{h} W \xrightarrow{k} Z]$ is given by
\begin{equation}
\begin{tikzcd}[column sep=small]
& & U \ar[dl,"h'"'] \ar[dr,"g'"] \\
& V \ar[dl,"f"'] \ar[dr,"g"] & & W \ar[dl,"h"'] \ar[dr,"k"] \\
X & & Y & & Z
\end{tikzcd}
\end{equation}
where $U$ is a pullback of $V\xrightarrow{g} Y \xleftarrow{h} W$. Note that with $g$ also $g'$ lies in $\cat C_\gsindex$ and hence so does the composite $k\comp g'$.
\end{defenum}
\end{definition}

\begin{lemma}\label{rslt:Corr-preserves-limits}
The functor $\Corr\colon \GeomSetup \to \Cat$ preserves limits.
\end{lemma}
\begin{proof}
By \cref{rslt:GeomSetup-has-limits}, $\GeomSetup$ is complete. Note that $\Corr \colon \GeomSetup \to \Fun(\bbDelta^\op, \Ani)$ preserves limits, because limits in $\Fun(\bbDelta^\op, \Ani)$ can be computed pointwise and the functor $\GeomSetup \to \Ani$, $\Hom_{\GeomSetup}((\bbSigma^n, \bbSigma_2^n), \blank)$, clearly commutes with limits for all $n\ge0$. Moreover, the embedding $\CSegAni \subseteq \Fun(\bbDelta^\op,\Ani)$ of complete Segal anima into simplicial anima preserves and reflects limits since it is a right adjoint. Finally, the equivalence $\CSegAni \isoto \Cat$ trivially preserves limits.
\end{proof}

\begin{remark}
\begin{remarksenum}
\item\label{rmk:embedding-into-Corr(C)} The morphisms $(\bbSigma^\bullet, \bbSigma_2^\bullet) \to ([\blank],[\blank])$ and $(\bbSigma^\bullet, \bbSigma_2^\bullet) \to ([\blank]^\op, [\blank]^\simeq)$ of cosimplicial geometric setups from \cref{ex:geometric-setups} induce functors
\[
\cat C_\gsindex \to \Corr(\cat C, \cat C_\gsindex)
\qquad\text{and}\qquad
\cat C^\op \to \Corr(\cat C,\cat C_\gsindex),
\]
which are functorial in $(\cat C, \cat C_\gsindex) \in \GeomSetup$.
The first functor sends a morphism $X\to Y$ in $\cat C_\gsindex$ to $X\xleftarrow{\id} X\to Y$, while the second sends a morphism $X\to Y$ in $\cat C$ to $Y\from X \xrightarrow{\id} X$. 

\item\label{rmk:Corr(C)-self-dual} The non-trivial automorphism $\bbSigma^\bullet \isoto (\bbSigma^\bullet)^\op$ given by $(i,j) \mapsto (j,i)$ induces an isomorphism
\[
\Corr(\cat C_\gsindex) \isoto \Corr(\cat C_\gsindex)^\op,
\]
which is functorial in $(\cat C, \cat C_\gsindex) \in \GeomSetup$.
\end{remarksenum}
\end{remark}

\begin{example}\label{ex:Corr(C-Iso)}
Let $\cat C$ be a category with underlying anima $\cat C^\simeq$. Then $(\cat C, \cat C^\simeq)$ is a geometric setup and the functor
\begin{align}
\cat C^\op \xto{\sim} \Corr(\cat C, \cat C^\simeq)
\end{align}
from \cref{rmk:embedding-into-Corr(C)} is an isomorphism in $\Cat$. Indeed, it is clearly essentially surjective, and the fully faithfulness follows easily from the fact that the anima $\cat C^\simeq_{X/}$ is contractible for each $X\in \cat C$ (see~\cite[\href{https://kerodon.net/tag/018Y}{Tag 018Y}]{kerodon}). 
\end{example}

\begin{proposition}\label{rslt:Corr-cocartesian-fibration}
Let $F\colon (\cat E, \cat E_\gsindex)\to (\cat C, \cat C_\gsindex)$ be a morphism of geometric setups. Suppose the following conditions are satisfied:
\begin{enumerate}[(a)]
\item the underlying functor $F \colon \cat E \to \cat C$ is a cartesian fibration;
\item $F$ admits cocartesian lifts of edges in $\cat C_{\gsindex}$ which lie in $\cat E_{\gsindex}$, i.e.\ for every $X \in \cat E$ and every edge $\ol{f}\colon F(X) \to \ol Y$ in $\cat C_{\gsindex}$ there exists an $F$-cocartesian edge $f\colon X\to Y$ in $\cat E_{\gsindex}$ with $F(f) = \ol{f}$;
\item let $\sigma\colon [1]\times [1]\to \cat E$ be a commutative square such that the restrictions $\sigma|_{\{i\}\times [1]}$ are $F$-co\-car\-te\-sian edges of $\cat E_\gsindex$ and $F\comp\sigma$ is a pullback square in $\cat C$. Then $\sigma$ is a pullback square.
\end{enumerate}
Then the induced functor $\Corr(\cat E,\cat E_\gsindex) \to \Corr(\cat C, \cat C_\gsindex)$ is a cocartesian fibration.
\end{proposition}
\begin{proof}
We call a diagram $\bbSigma^1\to \cat E$ \emph{special} if it sends $(0,1) \to (0,0)$ to an $F$-cartesian and $(0,1)\to (1,1)$ to an $F$-cocartesian morphism in $\cat E$. We call a diagram $\bbSigma^2\to \cat E$ \emph{special} if its restriction along $d_2\colon \bbSigma^1\to \bbSigma^2$ is special and it sends $(0,2)\to (1,2)$ to an $F$-cocartesian map. We will prove that the diagram 
\begin{equation}
\begin{tikzcd}
\Hom'\bigl((\bbSigma^2,\bbSigma^2_2), (\cat E, \cat E_\gsindex)\bigr) \ar[d] \ar[r,"d^*_1"] &
\Hom\bigl((\bbSigma^1, \bbSigma^1_2), (\cat E, \cat E_\gsindex)\bigr) \ar[d] \\
\Hom\bigl((\bbSigma^2,\bbSigma^2_2), (\cat C,\cat C_\gsindex)\bigr) \ar[r,"d^*_1"'] &
\Hom\bigl((\bbSigma^1, \bbSigma^1_2), (\cat C, \cat C_\gsindex)\bigr)
\end{tikzcd}
\end{equation}
is a pullback of anima, where $\Hom'((\bbSigma^2,\bbSigma^2_2), (\cat E, \cat E_\gsindex)) \subseteq \Hom((\bbSigma^2,\bbSigma^2_2), (\cat E, \cat E_\gsindex))$ denotes the full subanima spanned by the special diagrams. In view of \cref{rslt:Corr-complete-Segal-anima} this implies that every morphism in $\Corr(\cat C, \cat C_\gsindex)$ admits a special lift and that every special morphism in $\Corr(\cat E, \cat E_\gsindex)$ is $\Corr(F)$-cocartesian, proving the claim.

In order to show that the above diagram is indeed a pullback diagram, it suffices to prove that restriction along $d_1\colon \bbSigma^1 \to \bbSigma^2$ induces an isomorphism
\begin{align}
\Fun'(\bbSigma^2, \cat E) \isoto \Fun(\bbSigma^1,\cat E) \times_{\Fun(\bbSigma^1,\cat C)} \Fun(\bbSigma^2,\cat C),
\end{align}
where $\Fun'(\bbSigma^2,\cat E) \subseteq \Fun(\bbSigma^2,\cat E)$ is the full subcategory spanned by the special diagrams: Using the assumptions on $F$ it is then easily checked that this isomorphism restricts to an isomorphism as in the claim.

We introduce auxiliary subposets $\bbSigma^1 \subseteq P_1 \subseteq P_2 \subseteq P_3 \subseteq \bbSigma^2$ given by
\begin{align}
P_1 &= 
\left\{\begin{tikzcd}[ampersand replacement=\&, sep=small, cramped]
\& \&[-.7em] \bullet \ar[ddll, bend right] \ar[ddrr,bend left] \&[-.7em] \\
\& \bullet \ar[dl,"a"] \& \& \phantom{\bullet}\\[-.7em]
\bullet \& \& \& \& \bullet
\end{tikzcd}\right\},
&
P_2 &=
\left\{\begin{tikzcd}[ampersand replacement=\&, sep=small, cramped]
\& \&[-.7em] \bullet \ar[ddll, bend right] \ar[ddrr,bend left] \ar[dl] \&[-.7em] \\
\& \bullet \ar[dl,"a"] \& \& \phantom{\bullet}\\[-.7em]
\bullet \& \& \& \& \bullet
\end{tikzcd}\right\},
&
P_3 &=
\left\{\begin{tikzcd}[ampersand replacement=\&, sep=small, cramped]
\& \&[-.7em] \bullet \ar[ddll, bend right] \ar[ddrr,bend left] \ar[dl] \ar[dr, "c"'] \&[-.7em] \\
\& \bullet \ar[dl,"a"] \ar[dr,"b"'] \& \& \bullet \ar[dl] \\[-.7em]
\bullet \& \& \bullet \& \& \bullet
\end{tikzcd}\right\}.
\end{align}
We then compute
\begin{align}
\Fun'(\bbSigma^2,\cat E) &= \Fun'(P_3, \cat E)\times_{\Fun(P_3,\cat C)} \Fun(\bbSigma^2,\cat C) \\
&= \Fun'(P_2,\cat E) \times_{\Fun(P_2, \cat C)} \Fun(\bbSigma^2,\cat C) \\
&= \Fun'(P_1, \cat E) \times_{\Fun(P_1,\cat C)} \Fun(\bbSigma^2, \cat C) \\
&= \Fun(\bbSigma^1, \cat E) \times_{\Fun(\bbSigma^1, \cat C)} \Fun(\bbSigma^2,\cat C),
\end{align}
where $\Fun'(P_i,\cat E) \subseteq \Fun(P_i,\cat E)$ denote the full subcategories spanned by those functors which send $a$ to an $F$-cartesian morphism and, in case $i=3$, the maps $b$, $c$ to an $F$-cocartesian morphism.
The first isomorphism uses $\bbSigma^2 = P_3 \coprod_{\Horn^2_0} [2]$ and \cref{rslt:lifting-property-of-cocartesian-fibrations}, the second isomorphism uses $P_3 = P_2 \coprod_{[1]} [1]^2$ and \cref{rslt:exponentiation-of-cocartesian-fibrations}, the third isomorphism uses $P_2 = P_1 \coprod_{\Horn^2_2} [2]$ and (the dual of) \cref{rslt:lifting-property-of-cocartesian-fibrations}, and the last isomorphism uses $P_1 = \bbSigma^1 \coprod_{[0]} [1]$ and the definition of cartesian fibration (\cref{def:cocartesian-fibration}).
\end{proof}

\subsection{The symmetric monoidal structure}\label{subsec:monoidal-structure-on-Corr} 
Let $(\cat C, \cat C_\gsindex)$ be a geometric setup. If $\cat C$ admits finite products, then $\Corr(\cat C, \cat C_\gsindex)$ can be endowed with a symmetric monoidal structure which we will discuss next.

Recall the cocartesian operad $(\cat C^\op)^\amalg \to \Fin_*$ from \cref{ex:cocartesian-monoidal-structure}. For each $\langle n\rangle \in \Fin_*$, we identify the fiber $(\cat C^\op)^{\amalg, \op}_{\langle n\rangle}$ with the $n$-fold product $\cat C^n$. Hence, an edge $f\colon X_1\oplus\dotsb\oplus X_n\to Y_1\oplus\dotsb \oplus Y_n$ in the fiber is given by morphisms $f_i\colon X_i\to Y_i$ in $\cat C$ for all $1\le i\le n$.


\begin{lemma}\label{rslt:Corr^tensor-explicit}
Let $(\cat C, \cat C_\gsindex)$ be a geometric setup and put
\begin{align}
\cat C_\gsindex^- \coloneqq (\cat C_\gsindex^\op)^{\amalg,\op} \times_{\Fin_*^\op} \Fin_*^\simeq.
\end{align}
Then $((\cat C^\op)^{\amalg,\op}, \cat C_\gsindex^-)$ is a geometric setup.
\end{lemma}
\begin{proof}
We have to check that $\cat C_\gsindex^-$ is stable under pullbacks and that $\cat C_{\gsindex}^-$ admits and $\cat C_{\gsindex}^- \to (\cat C^\op)^{\amalg,\op}$ preserves fiber products. So let $f\colon X_1\oplus \dotsb \oplus X_n\to Y_1\oplus\dotsb\oplus Y_n$ be an edge in $\cat C_\gsindex^-$, which without loss of generality we may assume to lift the identity of $\langle n\rangle$. Let $g\colon Z_1\oplus\dotsb\oplus Z_m \to Y_1\oplus\dotsb\oplus Y_n$ be an arbitrary edge in $(\cat C^\op)^{\amalg,\op}$ lifting a map $\alpha\colon \langle n\rangle \to \langle m\rangle$; then $g$ is given by maps $g_i\colon Z_{\alpha(i)} \to Y_i$ in $\cat C$, for each $1\le i\le n$ with $\alpha(i) \in \langle m\rangle^\circ$.

For any $1\le j\le m$, let $W_j$ be the iterated pullback over $Z_j$ of the family $\{X_i\times_{Y_i}Z_j\}_{i\in \alpha^{-1}(j)}$. Since $(\cat C, \cat C_\gsindex)$ is a geometric setup, it is clear that $W_j$ exists and the projection $W_j\to Z_j$ lies in $\cat C_\gsindex$, and hence the induced map $f' \colon W_1\oplus\dotsb \oplus W_m \to Z_1\oplus \dotsb\oplus Z_m$ lies in $\cat C_\gsindex^-$. By the construction, $f'$ is a pullback of $f$ along $g$. The other statements are immediate.
Hence, $((\cat C^\op)^{\amalg,\op}, \cat C_\gsindex^-)$ is a geometric setup.
\end{proof}

\begin{definition}\label{def:Corr-operad}
Let $(\cat C, \cat C_\gsindex)$ be a geometric setup. We put 
\begin{align}
\Corr(\cat C, \cat C_\gsindex)^\tensor \coloneqq \Corr\bigl((\cat C^\op)^{\amalg,\op}, \cat C_{\gsindex}^-\bigr).
\end{align}
\end{definition}

Explicitly, an object in $\Corr(\cat C, \cat C_\gsindex)^\tensor$ is a tuple $X_\bullet = (X_1, \dots, X_n)$ of objects in $\cat C$ and a morphism of such tuples is given by a correspondence $X_\bullet \from Y_\bullet \to Z_\bullet$ such that $Y_\bullet$ and $Z_\bullet$ are tuples of the same size and the map $Y_\bullet \to Z_\bullet$ is given by a tuple of maps $Y_i \to Z_i$ in $E$.

\begin{proposition}\label{rslt:Corr-tensor-is-an-operad}
Let $(\cat C, \cat C_\gsindex)$ be a geometric setup. Then $\Corr(\cat C, \cat C_\gsindex)^\tensor$ is an operad with underlying category $\Corr(\cat C, \cat C_\gsindex)$. If $\cat C$ admits finite products, then $\Corr(\cat C, \cat C_\gsindex)^{\tensor}$ is a symmetric monoidal category.
\end{proposition}
\begin{proof}
Denote $q\colon \Corr(\cat C,\cat C_\gsindex)^\tensor \to \Fin_*$ the projection. We verify the axioms of an operad, see \cref{def:operad}.

Let $\alpha\colon \langle n\rangle \to \langle m\rangle$ be an inert edge in $\Fin_*$ and let $X\in \Corr(\cat C, \cat C_\gsindex)^{\tensor}_{\langle n\rangle}$. We need to find a $q$-cocartesian lift $X\to Y$ of $\alpha$ in $\Corr(\mathcal C, \mathcal C_\gsindex)^\tensor$. Let $f\colon X\to Y$ be a cocartesian lift of $\alpha$ in $(\cat C^\op)^\amalg$. Then $\CorrHom{X}{f^\op}{Y}{\id}Y$ is the desired $q$-cocartesian lift of $\alpha$. The same argument shows that $q$ is a cocartesian fibration if $\cat C$ admits finite products so that $(\cat{C}^\op)^{\amalg}\to \Fin_*$ is a cocartesian fibration. 

If $X_1,\dotsc,X_n \in \Corr(\cat C, \cat C_\gsindex)^\tensor_{\langle 1\rangle}$ are arbitrary, we need to find $X \in \Corr(\cat C, \cat C_\gsindex)^\tensor_{\langle n\rangle}$ with $X = X_1\oplus\dotsb\oplus X_n$; in other words, we need to find $X$ together with $q$-cocartesian lifts $X\to X_i$ of the Segal maps $\rho_i\colon \langle n\rangle \to \langle 1\rangle$, which are given by sending $i\mapsto 1$ and everything else to the base point. Since $(\cat C^\op)^\amalg$ is an operad, we find $X \in (\cat C^\op)^\amalg_{\langle n\rangle}$ with $X = X_1\oplus\dotsb\oplus X_n$, which then also holds true in $\Corr(\cat C,\cat C_\gsindex)^\tensor$. 

Finally, let $\alpha \colon \langle n\rangle \to \langle m\rangle$ in $\Fin_*$. Given $X, Y \in \Corr(\cat C, \cat C_\gsindex)^\tensor$ living over $\langle n\rangle$ and $\langle m\rangle$, respectively, we need to check that the map
\[
\Hom^{\alpha}(X,Y) \to \prod_{i=1}^m \Hom^{\rho_i\alpha}(X, Y_i),
\]
induced by the cocartesian lifts $Y\to Y_i$ of $\rho_i$,
is an isomorphism of anima. We compute
\begin{align}
\Hom^\alpha(X,Y) 
&= \Bigl[(\cat{C}^{\op,\amalg})_{X/_{\alpha}} \times_{\cat C^{\op,\amalg}} (\cat{C}_\gsindex^{-,\op})_{Y/_{\id}} \Bigr]^{\simeq} \\
&= \Bigl[(\cat C^{\op,\amalg}_{\langle m\rangle})_{X/_{\alpha}} \times_{\cat C^{\op,\amalg}_{\langle m\rangle}} \prod_{i=1}^m (\cat C_\gsindex^\op)_{Y_{i}/} \Bigr]^{\simeq} \\
&= \Bigl[\prod_{i=1}^m (\cat C^{\op,\amalg})_{X/_{\rho_i\alpha}} \times_{\cat{C}^{\op,\amalg}} (\cat C_\gsindex^\op)_{Y_{i}/} \Bigr]^{\simeq} \\
&= \prod_{i=1}^m \Hom^{\rho_{i}\alpha}(X, Y_{i}),
\end{align}
where, for example, $(\cat C^{\op,\amalg})_{X/_{\alpha}}$ is the full subcategory of $(\cat C^{\op,\amalg})_{X/}$ of those maps $X\to X'$ which lift $\alpha$. Hence, $\Corr(\cat C, \cat C_\gsindex)^\tensor$ is an operad.
\end{proof}

\begin{example}
Let $(\cat C, \cat C_\gsindex)$ be a geometric setup such that $\cat C$ admits finite products. The symmetric monoidal structure on $\Corr(\cat C, \cat C_\gsindex)$ is given by the cartesian product in $\cat C$. Indeed, the proof of \cref{rslt:Corr-tensor-is-an-operad} shows that for two objects $X,Y \in \Corr(\cat C, \cat C_\gsindex)$, the cocartesian lift of the unique map $\alpha\colon \langle 2\rangle \to \langle 1\rangle$ with $\alpha(1) = \alpha(2) = 1$ is given by a  correspondence of the form
\begin{align}
\CorrHom{X\oplus Y}{}{X \times Y}{\id}{X\times Y.}
\end{align}
\end{example}

\begin{proposition}\label{rslt:functoriality-of-Corr^tensor}
The assignment $(\cat C, \cat C_\gsindex) \mapsto \Corr(\cat C, \cat C_\gsindex)^\tensor$ enhances to a functor
\[
\GeomSetup \to \Op
\]
into the category $\Op$ of operads. In particular, every morphism $(\cat C, \cat C_\gsindex) \to (\cat D, \cat D_\gsindex)$, where $\cat C$ and $\cat D$ admit finite products, induces a lax symmetric monoidal functor $\Corr(\cat C, \cat C_\gsindex)^\tensor \to \Corr(\cat D, \cat D_\gsindex)^\tensor$.
\end{proposition}
\begin{proof}
It is clear from \cite[Construction~2.4.3.1]{HA} that the assignment $\cat C \mapsto (\cat C^\op)^{\amalg,\op}$ determines a functor $\Cat \to \Cat_{/\Fin_*^\op}$. The induced map $\Fun([1], \Cat)\to \Fun([1], \Cat_{/\Fin_*^\op}) \cong \Fun([1], \Cat)_{/\id_{\Fin_*^\op}}$ restricts to a functor
\[
\GeomSetup \to \GeomSetup_{/(\Fin_*^\op, \Fin_*^\op)}.
\]
The composite with $\GeomSetup_{/(\Fin_*^\op, \Fin_*^\op)} \xrightarrow{\Corr} \Cat_{/\Corr(\Fin_*^\op)} \to \Cat_{/\Fin_*}$, where the second functor is given by pullback along the functor $\Fin_* \to \Corr(\Fin_*^\op)$ from \cref{rmk:embedding-into-Corr(C)}, then factors through $\Op \subseteq \Cat_{/\Fin_*}$ by \cref{rslt:Corr-tensor-is-an-operad}. The remaining statements are now obvious.
\end{proof}

\begin{example}
Let $f\colon (\cat C, \cat C_\gsindex) \to (\cat D, \cat D_\gsindex)$ be a morphism of geometric setups, where $\cat C$ and $\cat D$ admit finite products. Concretely, the lax symmetric monoidal structure on $\Corr(\cat C, \cat C_\gsindex)^\tensor \to \Corr(\cat D, \cat D_\gsindex)^\tensor$ is induced by the correspondences
\[
\
\CorrHom{\prod_{i=1}^nf(X_i)}{}{f(\prod_{i=1}^nX_i)}{\id}{f(\prod_{i=1}^nX_i).}
\]
In particular, we see that if $f$ preserves finite products, then $\Corr(\cat C, \cat C_\gsindex)^\tensor \to \Corr(\cat D, \cat D_\gsindex)^\tensor$ is even symmetric monoidal.
\end{example}

For a category $\cat C$ admitting finite products, we recall the Cartesian operad $\cat C^\times \to \Fin_*$ from \cref{ex:cartesian-monoidal-structure}.

\begin{proposition}\label{rslt:embedding-into-Corr^tensor}
Let $(\cat C, \cat C_\gsindex)$ be a geometric setup. 
\begin{propenum}
\item The functor $\cat C^\op \to \Corr(\cat C, \cat C_\gsindex)$ enhances to a morphism of operads
\begin{align}
(\cat C^\op)^\amalg \to \Corr(\cat C, \cat C_\gsindex)^\tensor,
\end{align}
which is a symmetric monoidal functor if $\cat C$ admits finite products.

\item Suppose that $\cat C$ admits finite products and $*$ is also terminal in $\cat C_{\gsindex}$. Then the functor $\cat C_\gsindex \to \Corr(\cat C, \cat C_\gsindex)$ enhances to a symmetric monoidal functor
\begin{align}
\cat C_\gsindex^\times \to \Corr(\cat C, \cat C_\gsindex)^\tensor.
\end{align}
\end{propenum}
\end{proposition}
\begin{proof}
The map $(\cat C^\op)^\amalg \to \Corr((\cat C^\op)^{\amalg,\op}, \cat C_\gsindex^-) = \Corr(\cat C, \cat C_\gsindex)^\tensor$ is provided by \cref{rmk:embedding-into-Corr(C)}, which establishes (i).

We now prove (ii). The cocartesian fibration $\cat C_\gsindex^\times \to \Fin_*$ is equivalent to the dual fibration $(\cat C_\gsindex^{\op,\amalg,\op})^\vee \to \Fin_*$ defined in \cite[Definition~3.5]{Barwick-Glasman-Nardin.2018}; see \cite[Theorem~1.7]{Barwick-Glasman-Nardin.2018} and \cite[Corollary~2.4.1.8]{HA}. The construction of the dual fibration provides a canonical symmetric monoidal functor
\begin{align}\label{eq:dual-fibration-to-Corr}
(\cat C_\gsindex^{\op,\amalg,\op})^\vee \to \Corr(\cat C_\gsindex)^\tensor.
\end{align}
From the definition of geometric setups and our assumptions on $\cat C$ and $\cat C_{\gsindex}$ it follows that $\cat C_\gsindex \to \cat C$ preserves finite limits. Hence, $\Corr(\cat C_\gsindex)^\tensor \to \Corr(\cat C, \cat C_\gsindex)^\tensor$ is a symmetric monoidal functor by \cref{rslt:functoriality-of-Corr^tensor}, and precomposing with \eqref{eq:dual-fibration-to-Corr} yields the desired map.
\end{proof}

We record the following result on duals in the correspondence category, which has already been observed in various other places, see e.g.\ \cite[2.18]{Barwick-Glasman-Shah.2019}, \cite[Corollary~12.5]{Haugseng.2017a}, \cite[Lemma~2.2.5]{Zavyalov.2023}.

\begin{proposition}\label{rslt:duals-in-Corr}
Let $(\cat C, E)$ be a geometric setup and assume that $\cat C$ admits finite products. Let $X\in \cat C$ such that $X\to *$ lies in $E$. Then $X$ is dualizable in $\Corr(\cat C, E)^\tensor$ with dual $X^\vee = X$. 
\end{proposition}
\begin{proof}
For the notion of dualizable object we refer to \cref{def:dualizable-object}. Observe that with $X\to *$ also the diagonal $\Delta\colon X\to X\times X$ lies in $E$. We claim that the evaluation and coevaluation maps are given by the correspondences $\ev_X = [\CorrHom{X\times X}{\Delta}{X}{}{*}]$ and $i_X = [\CorrHom{*}{}{X}{\Delta}{X\times X}]$, respectively. We need to check that $(\ev_X\times\id_X)\comp(\id_X\times i_X) = \id_X$ and $(\id_X\times \ev_X) \comp (i_X\times\id_X) = \id_X$. We only verify the first identity, since the other one is similar. The claim reduces to showing that the square in the diagram
\begin{equation}
\begin{tikzcd}[column sep=small]
& & X \ar[dl,"\Delta"'] \ar[dr,"\Delta"] \\
& X\times X \ar[dl,"\pr_1"'] \ar[dr,"\id\times\Delta"'] & & X\times X \ar[dl,"\Delta\times \id"] \ar[dr,"\pr_2"] \\
X & & X\times X\times X & & X
\end{tikzcd}
\end{equation}
is a pullback in $\cat C$. We contemplate the diagram
\begin{equation}
\begin{tikzcd}
X \ar[d,"\Delta"'] \ar[r,"\Delta"] & X\times X \ar[d,"\id\times\Delta"'] \ar[r,"\pr_2"] & X \ar[d,"\Delta"] \\
X\times X \ar[r,"\Delta\times\id"'] & X\times X\times X \ar[d,"\pr_1"'] \ar[r,"\pr_{2,3}"'] & X\times X \ar[d] \\
& X \ar[r] & *.
\end{tikzcd}
\end{equation}
The lower right square and the outer right rectangle are pullbacks. It follows from the pasting law (\cite[\href{https://kerodon.net/tag/03FZ}{Tag 03FZ}]{kerodon}) that the top right square is a pullback. Since the composites of the horizontal arrows are each the identity, the outer left rectangle is a pullback. By the pasting law again, the top left square is a pullback. This finishes the proof.
\end{proof}

\begin{remark}
If $(\cat C, E)$ is a geometric setup such that $\cat C$ admits finite limits and $*$ is also terminal in $\cat C_E$, then \cref{rslt:duals-in-Corr} shows that every object of $\Corr(\cat C, E)^\tensor$ is self-dual.
\end{remark}

\subsection{Self-enrichment}\label{subsec:self-enrichment-of-Corr} 
We fix a category $\cat C$ which admits finite limits so that $(\cat C, \cat C)$ is a geometric setup. In this section we will show that $\Corr(\cat C)$ is enriched over itself. We start with the following consequence of \cref{rslt:duals-in-Corr}:

\begin{proposition}\label{rslt:Corr-closed-symmetric-monoidal}
The category $\Corr(\cat C)$ is closed symmetric monoidal. The internal hom is given by $\iHom(X,Y) = Y\times X$ for each $X,Y\in \Corr(\cat C)$.
\end{proposition}
\begin{proof}
\cref{rslt:duals-in-Corr} shows that every object in $\Corr(\cat C)$ is self-dual. The computation of the internal hom now follows from \cref{rslt:properties-of-dualizable-objects}.
\end{proof}

\begin{corollary}\label{rslt:Corr-self-enriched-and-self-dual}
The category $\Corr(\cat C)$ is enriched over itself, and $X\mapsto X$ defines a $\Corr(\cat C)$-enriched equivalence $\Corr(\cat C)^\op \isoto \Corr(\cat C)$.
\end{corollary}
\begin{proof}
By \cref{rslt:Corr-closed-symmetric-monoidal}, every object $X \in \Corr(\cat C)$ is dualizable with $X^\vee = X$. Now the assertion follows from \cref{rslt:symmetric-monoidal-self-duality}.
\end{proof}

\begin{example}
We describe $\Corr(\cat C)$ as a $\Corr(\cat C)$-enriched category:
\begin{itemize}
\item The objects are the objects of $\cat C$.
\item The enriched hom is given by $\enrHom[\Corr(\cat C)](X,Y) = \iHom(X,Y) = Y\times X$, and the evaluation map is given by the correspondence
\[
\CorrHom{(Y\times X)\times X}{\id_Y\times\Delta_X}{Y\times X}{\pr_Y}{Y.}
\]
\item The composition $\iHom(Y,Z) \times \iHom(X,Y) \to \iHom(X,Z)$ is given by the correspondence
\[
\CorrHom{(Z\times Y)\times (Y\times X)}{\id_Z\times\Delta_Y\times\id_X}{Z\times Y\times X}{\pr_{Z,X}}{Z\times X.}
\]
\item The unit $*\to \iHom(X,X)$ is given by $i_X = [\CorrHom{*}{}{X}{\Delta_X}{X\times X}]$.
\end{itemize}
\end{example}





\subsection{A 2-categorical enhancement on geometric setups}
The results of this section will not strictly be necessary in this paper and hence may be skipped on first reading. In the following we will answer the following question: In what sense is the construction $(\cat C, \cat C_{\gsindex}) \mapsto \Corr(\cat C, \cat C_{\gsindex})$ 2-functorial, i.e.\ which kind of \enquote{natural transformations} between morphisms of geometric setups will induce natural transformations between the associated functors on the correspondence categories? To this end, we will enrich the category $\GeomSetup$ by introducing a 2-category $\cat G$ as follows:
\begin{itemize}
\item The objects of $\cat G$ are geometric setups $\gs{C} = (\cat C, \cat C_{\gsindex})$.
\item The $1$-morphisms of $\cat G$ are morphisms of geometric setups.
\item Let $F,F'\colon \gs{C} \to \gs{D}$ be two $1$-morphisms in $\cat G$. A $2$-morphism from $F$ to $F'$ is a correspondence $\CorrHom{F}{\alpha}{F''}{\beta}{F'}$, where $\alpha$ and $\beta$ are natural transformations for which certain naturality squares are cartesian; see \cref{def:internal-Hom-in-GeomSetup} for the precise condition.
\end{itemize}

We will then show that $\Corr(\blank)$ upgrades to a 2-functor $\cat G \to \TwoCat$ (see \cref{rslt:Corr-as-a-2-functor} below). In order to construct $\cat G$, recall from \cref{rslt:GeomSetup-has-limits} that $\GeomSetup$ admits all small limits, which can be computed in $\Fun([1],\Cat)$. In particular, $\GeomSetup$ is equipped with the cartesian monoidal structure from \cref{ex:cartesian-monoidal-structure}. As a first step, we will show that the monoidal structure on $\GeomSetup$ is closed, i.e.\ that internal hom functors exist. A similar construction was independently performed in \cite[Lemma~2.5]{Haugseng-Hebestreit.Spans}.

\begin{definition}\label{def:internal-Hom-in-GeomSetup}
Let $\gs{C} = (\cat C, \cat C_{\gsindex})$ and $\gs{D} = (\cat D, \cat D_{\gsindex})$ be geometric setups.
\begin{defenum}
\item We denote $\Fun^1(\gs{C}, \gs{D})$ the (non-full) subcategory of $\Fun(\cat C, \cat D)$ whose objects are the morphisms of geometric setups and whose morphisms are natural transformations $\alpha\colon F\to F'$ which satisfy the following condition: For every morphism $f\colon Y\to X$ in $\cat C_{\gsindex}$ the following commutative diagram in $\cat D$ is cartesian:
\begin{equation}\label{eq:internal-Hom-in-GeomSetup}
\begin{tikzcd}
F(Y) \ar[d,"F(f)"'] \ar[r,"\alpha_Y"] \ar[dr,phantom,very near start,"\lrcorner"] & F'(Y) \ar[d,"F'(f)"] \\
F(X) \ar[r,"\alpha_X"'] & F'(X).
\end{tikzcd}
\end{equation}

\item We denote $\Fun^2(\gs{C}, \gs{D})$ the wide subcategory of $\Fun^1(\gs{C}, \gs{D})$ whose morphisms are the natural transformations $\alpha\colon F\to F'$ which satisfy the following conditions: For every $X\in \cat C$ the morphism $\alpha_X$ lies in $\cat D_{\gsindex}$, and for every morphism $f\colon Y\to X$ in $\cat C$ the diagram \eqref{eq:internal-Hom-in-GeomSetup} is cartesian in $\cat D$.
\end{defenum}
\end{definition}

\begin{proposition}\label{rslt:internal-Hom-is-a-GeomSetup}
Let $\gs{C}, \gs{D}$ be geometric setups. The pair
\[
\Fun_{\GeomSetup}\bigl(\gs{C}, \gs{D}\bigr) \coloneqq \bigl(\Fun^1(\gs{C}, \gs{D}), \Fun^2(\gs{C}, \gs{D})\bigr)
\]
is a geometric setup.
\end{proposition}

The proposition is an immediate consequence of the following technical result:

\begin{lemma}\label{rst:Fun1-has-fiber-products}
Let $\gs{C} = (\cat C, \cat C_{\gsindex})$ and $\gs{D} = (\cat D, \cat D_{\gsindex})$ be geometric setups, let $\alpha\colon F'\to F$ be an edge in $\Fun^1(\gs{C}, \gs{D})$ and let $\beta\colon F''\to F$ be an edge in $\Fun^2(\gs{C}, \gs{D})$. We denote $H\coloneqq F'\times_FF''$ the fiber product in $\Fun(\cat C, \cat D)$. 
\begin{lemenum}
\item The functor $H \colon \cat C \to \cat D$ is a morphism of geometric setups.
\item The projection $p''\colon H\to F''$ lies in $\Fun^1(\gs{C}, \gs{D})$.
\item The projection $p'\colon H\to F'$ lies in $\Fun^2(\gs{C}, \gs{D})$.
\item $H$ is also the fiber product in $\Fun^1(\gs{C}, \gs{D})$.
\item If $\alpha$ lies in $\Fun^2(\gs{C}, \gs{D})$, then $H$ is the fiber product in $\Fun^2(\gs{C}, \gs{D})$.
\end{lemenum}
\end{lemma}
\begin{proof}
We first prove the following claim:
\begin{claim*}
\begin{enumerate}[label=(\alph*)]
\item For every morphism $f\colon Y\to X$ in $\cat C_{\gsindex}$ the following diagram in $\cat D$ is cartesian:
\begin{equation}
\begin{tikzcd}
H(Y) \ar[d,"H(f)"'] \ar[r,"p''_Y"] \ar[dr,phantom, very near start, "\lrcorner"] & F''(Y) \ar[d,"F''(f)"]
\\
H(X) \ar[r,"p''_X"'] & F''(X).
\end{tikzcd}
\end{equation}

\item For every morphism $f\colon Y\to X$ in $\cat C$ the following diagram in $\cat D$ is cartesian:
\begin{equation}
\begin{tikzcd}
H(Y) \ar[d,"H(f)"'] \ar[r,"p'_Y"] \ar[dr,phantom, very near start, "\lrcorner"] & F'(Y) \ar[d,"F'(f)"] 
\\
H(X) \ar[r,"p'_X"'] & F'(X).
\end{tikzcd}
\end{equation}
\end{enumerate}
\end{claim*}
\begin{proof}[Proof of the claim]
We will only prove (a); the proof of (b) is analogous. We have the following commutative diagrams, where the indicated squares are cartesian by the hypotheses on $H$, $\alpha$ and $f$:
\begin{equation}
\begin{tikzcd}
H(Y) \ar[d,"p'_Y"'] \ar[r,"p''_Y"] \ar[dr,phantom, very near start, "\lrcorner"] & F''(Y) \ar[d,"\beta_Y"]
&[5em]
H(Y) \ar[d,"H(f)"'] \ar[r,"p''_Y"] & F''(Y) \ar[d,"F''(f)"]
\\
F'(Y) \ar[d,"F'(f)"'] \ar[r,"\alpha_Y"] \ar[dr,phantom, very near start, "\lrcorner"] & F(Y) \ar[d,"F(f)"] 
&
H(X) \ar[d,"p'_X"'] \ar[r, "p''_X"] \ar[dr,phantom, very near start, "\lrcorner"] & F''(X) \ar[d,"\beta_X"]
\\
F'(X) \ar[r,"\alpha_X"'] & F(X)
& 
F'(X) \ar[r,"\alpha_X"'] & F(X).
\end{tikzcd}
\end{equation}
The outer diagrams coincide by the naturality of $p''$ and $\beta$. It follows that the upper right square is cartesian, which proves the claim.
\end{proof}

We prove (i). Claim (a) shows that $H \colon \cat C\to \cat D$ sends $\cat C_{\gsindex}$ into $\cat D_{\gsindex}$. Pullbacks of the diagrams $X \xto{f} S \from Y$, where $f$ lies in $\cat C_{\gsindex}$, are preserved by $F, F', F''$ and hence also by $H$, since fiber products commute with fiber products. Therefore, $H$ is a morphism of geometric setups.

Now, (ii) and (iii) are a reformulation of claims (a) and (b), respectively.

We prove (iv). Let $\varphi'\colon G\to F'$ and $\varphi''\colon G\to F''$ be morphisms in $\Fun^1(\gs{C}, \gs{D})$ such that $\alpha\comp\varphi' = \beta \comp \varphi''$ (by which we mean that the obvious diagram commutes). Let $\psi\colon G\to H$ be the unique natural transformation satisfying $p'\comp\psi = \varphi'$ and $p''\comp \psi = \varphi''$. We will show that $\psi$ lies in $\Fun^1(\gs{C}, \gs{D})$. For any morphism $f\colon Y\to X$ in $\cat C_{\gsindex}$ we have a commutative diagram
\begin{equation}
\begin{tikzcd}
G(Y) \ar[d,"G(f)"'] \ar[r,"\psi_Y"] & H(Y) \ar[d,"H(f)"'] \ar[r,"p''_Y"] \ar[dr,phantom, very near start, "\lrcorner"] & F''(Y) \ar[d,"F''(f)"]
\\
G(X) \ar[r,"\psi_X"'] & H(X) \ar[r,"p''_X"'] & F''(X),
\end{tikzcd}
\end{equation}
where the right square is cartesian since $p''$ lies in $\Fun^1(\gs{C}, \gs{D})$. Note that the outer diagram is cartesian, since $p''\comp \psi = \varphi''$ lies in $\Fun^1(\gs{C}, \gs{D})$. It follows that the left square is cartesian, which proves (iv).

The proof of (v) is analogous to (iv). Note that, with the notation in (iv), the fact that $\psi_X$ lies in $\cat D_{\gsindex}$ is immediate from the fact that the morphisms in $\cat D_{\gsindex}$ are right cancellative, see~\cref{rslt:right-cancellative-fiber-products}.
\end{proof}

\begin{proposition}\label{rslt:GeomSetup-is-closed}
The cartesian symmetric monoidal structure on $\GeomSetup$ is closed.
\end{proposition}
\begin{proof}
We will prove the assertion by showing that for any geometric setup $\gs{C} = (\cat C, \cat C_{\gsindex})$, the endofunctor $\Fun(\gs{C}, \blank)$ on $\GeomSetup$ constructed in \cref{rslt:internal-Hom-is-a-GeomSetup} is right adjoint to $\blank \times \gs{C}$. Concretely, we will define for any geometric setup $\gs{D}$ natural maps
\begin{align}
\ev\colon \Fun(\gs{C}, \gs{D})\times \gs{C} \to \gs{D}
\qquad\text{and}\qquad
\coev\colon \gs{D} \to \Fun(\gs{C}, \gs{D}\times \gs{C})
\end{align}
in $\GeomSetup$ such that the following triangles commute:
\begin{equation}
\begin{tikzcd}
\gs{D}\times \gs{C}
\ar[r,"\coev\times\id"] \ar[dr, equals] &[2em]
\Fun(\gs{C}, \gs{D}\times \gs{C}) \times \gs{C} \ar[d,"\ev"]
&[3em]
\Fun(\gs{C}, \gs{D}) 
\ar[dr,equals] \ar[r,"\coev"] &
\Fun\bigl(\gs{C}, \Fun(\gs{C}, \gs{D})\times \gs{C}\bigr) \ar[d,"{\Fun(\gs{C}, \ev)}"]
\\
& \gs{D}\times \gs{C}
& &
\Fun(\gs{C}, \gs{D}).
\end{tikzcd}
\end{equation}

We start with the coevaluation. Under the identification $\Fun(\cat C, \cat D\times \cat C) = \Fun(\cat C, \cat D) \times \Fun(\cat C, \cat D)$ we consider the functor $F= (\delta, \const_{\id_{\cat C}}) \colon \cat D\to \Fun(\cat C, \cat D\times \cat C)$, where $\delta$ is the diagonal functor defined by precomposing with $\cat C\to *$. Note that for every $X\in \cat D$ the functor $F(X)$ coincides with $(\const_X, \id_{\cat C})$. It is thus immediate that $F(X)\colon \gs{C} \to \gs{D}\times \gs{C}$ is a morphism in $\GeomSetup$ for all $X\in \cat D$, and that in fact $F$ identifies with a morphism $\coev\colon \gs{D} \to \Fun(\gs{C}, \gs{D}\times \gs{C})$.

We now construct the evaluation. Let $\ev\colon \Fun^1(\gs{C}, \gs{D}) \times \cat C \to \cat D$ be the restriction of the usual evaluation map. We claim that $\ev$ defines a morphism of geometric setups. Observe that $\ev$ clearly sends $\Fun^2(\gs{C}, \gs{D}) \times \cat C_{\gsindex}$ into $\cat D_{\gsindex}$, so it remains to verify that it preserves certain pullbacks. Consider diagrams $F'\to F \xfrom{\beta}F''$ be a diagram in $\Fun^1(\gs{C}, \gs{D})$ and $X\to S \xfrom{f}Y$ a diagram in $\cat C$, where $\beta$ lies in $\Fun^2(\gs{C}, \gs{D})$ and $f$ lies in $\cat C_{\gsindex}$. Denote $H = F'\times_FF''$ and $Z = X\times_SY$ the respective pullbacks. We need to show that
\begin{equation}
\begin{tikzcd}
H(Z) \ar[d] \ar[r] & F''(Y) \ar[d,"\beta\star f"] \\
F'(X) \ar[r] & F(S)
\end{tikzcd}
\end{equation}
is a pullback diagram in $\cat D$. We consider the following commutative diagram:
\begin{equation}
\begin{tikzcd}[sep=small]
F''(Z) \ar[dd] \ar[rr] \ar[dr]
& & 
F''(X) \ar[dr] \ar[dd]
\\
& F''(Y) \ar[rr,crossing over] 
& &
F''(S) \ar[dd]
\\
F(Z) \ar[rr] \ar[dr]
& &
F(X) \ar[from=dd] \ar[dr]
\\
& F(Y) \ar[from=uu, crossing over] \ar[rr,crossing over]
& & 
F(S)
\\
F'(Z) \ar[uu] \ar[rr] \ar[dr]
& &
F'(X) \ar[dr]
\\
& F'(Y) \ar[uu, crossing over] \ar[rr]
& &
F'(S) \ar[uu]
\end{tikzcd}
\end{equation}
The upper left face is cartesian since $\beta\colon F''\to F$ lies in $\Fun^2(\gs{C}, \gs{D})$. The middle horizontal face is cartesian since $F$ is a morphism of geometric setups and $f\colon Y\to S$ lies in $\cat C_{\gsindex}$. The bottom back face is cartesian because $Z\to X$ lies in $\cat C_{\gsindex}$ (being a pullback of $f$) and $F'\to F$ lies in $\Fun^1(\gs{C}, \gs{D})$. We now compute
\begin{align}
H(Z) &= F'(Z) \times_{F(Z)} F''(Z) \\
&= F'(Z) \times_{F(Z)} \bigl(F(Z)\times_{F(Y)} F''(Y)\bigr) \\
&= \bigl(F'(X) \times_{F(X)} F(Z)\bigr) \times_{F(Z)} \bigl(F(Z) \times_{F(Y)} F''(Y)\bigr) \\
&= F'(X) \times_{F(X)} F(Z) \times_{F(Y)} F''(Y) \\
&= F'(X) \times_{F(X)} \bigl(F(X) \times_{F(S)} F(Y)\bigr) \times_{F(Y)} F''(Y) \\
&= F'(X) \times_{F(S)} F''(Y),
\end{align}
where the second, third and fifth identification uses the cartesianness of the upper left, bottom back and middle horizontal face, respectively. This shows that $\ev$ is a morphism of geometric setups.

It remains to show that the triangle identities hold. For the left triangle we need to check that the composite of the horizontal arrows in the diagram
\begin{equation}
\begin{tikzcd}[column sep=3.5em]
\cat D\times \cat C \ar[dr] \ar[r,"\coev\times\id"] \ar[dr,"\coev\times\id"']
&
\Fun^1(\gs{C}, \gs{D}\times \gs{C})\times \cat C \ar[d] \ar[r,"\ev"]
&
\cat D\times \cat C \\
& \Fun(\cat C, \cat D\times \cat C)\times \cat C \ar[ur,"\ev"']
\end{tikzcd}
\end{equation}
is the identity, which is clear from the fact that this is true for the composite of the oblique arrows. For the right triangle we need to show that the composite, say $\varphi$, of the top horizontal arrows in the diagram
\begin{equation}
\begin{tikzcd}[column sep=4em]
\Fun^1(\gs{C}, \gs{D}) \ar[d,"\iota"'] \ar[r,"\coev"]
&
\Fun^1\bigl(\gs{C}, \Fun(\gs{C}, \gs{D})\times \gs{C}\bigr)
\ar[r,"{\Fun^1(\gs{C},\ev)}"] \ar[d]
&
\Fun^1(\gs{C}, \gs{D}) \ar[d,"\iota"]
\\
\Fun(\cat C, \cat D) \ar[dr,"\coev"'] \ar[r] & \Fun\bigl(\cat C, \Fun^1(\gs{C}, \gs{D}) \times \cat C\bigr) \ar[d] \ar[r] & \Fun(\cat C, \cat D)
\\
& \Fun\bigl(\cat C, \Fun(\cat C, \cat D)\times \cat C\bigr) \ar[ur, "{\Fun(\cat C, \ev)}"']
\end{tikzcd}
\end{equation}
is the identity. But the composite of the oblique arrows is the identity, and hence the commutativity of the diagram shows $\iota\comp \varphi = \iota$. But since $\iota \colon \Fun^1(\gs{C}, \gs{D}) \to \Fun(\cat C, \cat D)$ is a subcategory, hence a monomorphism by \cref{rslt:subcategories-are-monomorphisms}, we deduce $\varphi = \id$ as desired.
\end{proof}

\begin{definition}
We define the $2$-categorical enhancement of $\GeomSetup$ as follows: By \cref{rslt:GeomSetup-is-closed} and \cref{ex:monoidal-category-enriched-in-itself} the category $\GeomSetup$ is enriched over itself. Transferring the enrichment (\cref{def:transfer-of-enrichment}) along the monoidal functor $\Corr\colon \GeomSetup \to \Cat$ (see \cref{rslt:Corr-preserves-limits}) yields the $2$-category
\[
\cat G \coloneqq \tau_{\Corr}(\GeomSetup),
\]
which was described on objects, morphisms and 2-morphisms at the beginning of this section.
\end{definition}

\begin{proposition}\label{rslt:Corr-as-a-2-functor}
The functor $\Corr\colon \GeomSetup \to \Cat$ enhances to a $2$-functor
\begin{align}
\Corr \colon \cat G \to \Cat_2.
\end{align}
Concretely, it sends:
\begin{itemize}
\item a geometric setup $(\cat C, \cat C_{\gsindex})$ to $\Corr(\cat C, \cat C_{\gsindex})$, 
\item a morphism $F\in \Fun((\cat C, \cat C_{\gsindex}), (\cat D, \cat D_{\gsindex}))$ to the functor
\begin{align}
\Corr(F) \colon \Corr(\cat C, \cat C_{\gsindex}) \to \Corr(\cat D, \cat D_{\gsindex}),
\end{align}
\item a morphism $\CorrHom{F}{\alpha}{F''}{\beta}{F'}$ in $\Corr(\Fun((\cat C, \cat C_{\gsindex}), (\cat D, \cat D_{\gsindex})))$ to the natural transformation $\Corr(F) \to \Corr(F')$ given by $\CorrHom{F(X)}{\alpha_{X}}{F''(X)}{\beta_{X}}{F'(X)}$ for every $X \in \Corr(\cat C, \cat C_{\gsindex})$. 
\end{itemize}
\end{proposition}
\begin{proof}
Apply \cref{rslt:functor-from-transfer-of-self-enrichment} to the monoidal functor $\Corr\colon \GeomSetup \to \Cat$.
\end{proof}

\section{6-functor formalisms} \label{sec:6ff}

We now come to the central object of study in this paper: 6-functor formalisms. With the preparations from \cref{sec:Corr} it is now easy to define 6-functor formalisms, which we will do in \cref{sec:6ff.def3ff,sec:6ff.def}, following the first half of \cite[\S A.5]{Mann.2022a}. In \cref{sec:6ff.construct,sec:6ff.extend} we recall the results from the second half of \cite[\S A.5]{Mann.2022a} and from \cite[Appendix to Lecture~IV]{Scholze:Six-Functor-Formalism} on how to construct 6-functor formalisms in practice. We will then use these tools to provide some examples of 6-functor formalisms in \cref{sec:6ff.example}. Most importantly, we will construct the 6-functor formalism of $\Lambda$-valued sheaves on condensed anima, which will play an important role in our applications to representation theory in \cref{sec:reptheory}.

\subsection{3-functor formalisms} \label{sec:6ff.def3ff}

Recall our \cref{convention:geometric-setup-morphisms} on the notation of geometric setups.
As explained in the introduction, a 6-functor formalism $\D$ on a geometric setup $(\cat C, E)$ consists roughly of the following data:
\begin{itemize}
	\item For every $X \in \mathcal C$ there is a category $\D(X)$ of \enquote{sheaves on $X$}. This category is closed symmetric monoidal, i.e.\ comes equipped with a tensor product operation $\tensor$ admitting an internal hom functor $\iHom$.
	
	\item For every morphism $f\colon Y \to X$ in $\mathcal C$ there is a pullback functor $f^*\colon \D(X) \to \D(Y)$. Moreover, $f^*$ admits a right adjoint $f_*\colon \D(Y) \to \D(X)$, the pushforward functor.

	\item For every morphism $f\colon Y \to X$ in $E$ there is an exceptional pushforward functor $f_!\colon \D(Y) \to \D(X)$ which admits a right adjoint $f^!\colon \D(X) \to \D(Y)$, the exceptional pullback.
\end{itemize}
These data are required to satisfy some compatibilities: Everything is functorial, all pullback functors are symmetric monoidal (i.e.\ commute with the tensor product), and proper base-change and projection formula hold. These compatibilities result in additional data one needs to keep track of, resulting in an additional infinite amount of higher homotopies that need to be encoded in $\D$.

In the present subsection we will focus on the three functors $\tensor$, $f^*$ and $f_!$. A very elegant way to encode these three functors is given by the following definition due to Lurie:

\begin{definition} \label{def:3-functor-formalism}
Let $(\cat C, E)$ be a geometric setup. A \emph{3-functor formalism} on $(\cat C, E)$ is a map of operads
\begin{align}
	\D\colon \Corr(\cat C, E)^\tensor \to \Cat^\times,
\end{align}
where $\Corr(\cat C, E)^\tensor$ denotes the operad from \cref{def:Corr-operad} and $\Cat^\times$ denotes the cartesian operad associated with the category of categories (see \cref{ex:cartesian-monoidal-structure}). We say that $\D$ is \emph{stable} or \emph{presentable} if it factors over the category of stable categories or over the category of presentable categories (with the Lurie tensor product) respectively.
\end{definition}

\begin{remark} \label{rmk:3ff-is-lax-symm-mon-functor-if-C-has-finite-limits}
If $\cat C$ admits finite limits then $\Corr(\cat C, E)^\tensor$ is a symmetric monoidal category (by \cref{rslt:Corr-tensor-is-an-operad}) and hence a 3-functor formalism on $(\cat C, E)$ is a lax symmetric monoidal functor $\Corr(\cat C, E) \to \Cat$.
\end{remark}

\begin{remark}
In \cite[Definition~A.5.6]{Mann.2022a} we used the name \enquote{pre-6-functor formalism} to refer to 3-functor formalisms. However, we chose to abandon this name and rather use the terminology from \cite[Definition~2.4]{Scholze:Six-Functor-Formalism} because a surprising amount of constructions can already be performed for 3-functor formalisms, so they are a very natural notion to consider.
\end{remark}

Given a 3-functor formalism $\D$ one can easily extract the categories $\D(X)$ and the functors $\tensor$, $f^*$ and $f_!$ from it (see below for a more explicit explanation).

\begin{definition} \label{def:3-functors-from-3ff}
Let $\D\colon \Corr(\cat C, E)^\tensor \to \Cat^\times$ be a given 3-functor formalism on a geometric setup $(\cat C, E)$.
\begin{defenum}
	\item \label{def:tensor-product-from-3ff} By precomposing $\D$ with the map of operads $(\cat C^\op)^\amalg \to \Corr(\cat C, E)^\tensor$ from \cref{rslt:embedding-into-Corr^tensor}, we obtain a map of operads $(\cat C^\op)^\amalg \to \Cat^\times$, which by \cite[Theorem 2.4.3.18]{HA} is equivalent to a functor $\D^*\colon \cat C^\op \to \CMon$, where $\CMon$ denotes the category of symmetric monoidal categories. In particular, for every $X \in \cat C$ we obtain a symmetric monoidal structure on $\D(X)$. We denote the associated tensor operation by
	\begin{align}
		\blank \tensor \blank\colon \D(X) \times \D(X) \to \D(X).
	\end{align}

	\item For every morphism $f\colon Y \to X$ we denote
	\begin{align}
		f^* \coloneqq D^*(f)\colon \D(X) \to \D(Y)
	\end{align}
	and call it the \emph{pullback along $f$}.

	\item By precomposing $\D$ along the functor $\cat C_E \to \Corr(\cat C, E)$ from \cref{rmk:embedding-into-Corr(C)} we obtain a functor $\D_!\colon \cat C_E \to \Cat$. For every $f\colon Y \to X$ in $E$ we denote
	\begin{align}
		f_! \coloneqq \D_!(f)\colon \D(Y) \to \D(X)
	\end{align}
	and call it the \emph{exceptional pushforward along $f$}.
\end{defenum}
\end{definition}

Let us spell out more explicitly what happens in \cref{def:3-functors-from-3ff}. Since the objects of $\Corr(\cat C, E)$ are the same as those of $\cat C$, we get a category $\D(X)$ for every object $X \in \cat C$. Moreover, for every morphism $f\colon Y \to X$ we obtain the morphism $[\CorrHom{X}{f}{Y}{\id}{Y}]$ in $\Corr(\cat C, E)$, which is sent by $\D$ to the functor $f^*\colon \D(X) \to \D(Y)$. Similarly, if $g\colon Y \to Z$ is a morphism in $E$ then we get a morphism $[\CorrHom{Y}{\id}{Y}{g}{Z}]$ in $\Corr(\cat C, E)$, resulting in the functor $g_!\colon \D(Y) \to \D(Z)$ via $\D$. Now note that in $\Corr(\cat C, E)$ we have
\begin{align}
	\left[\CorrHom{Y}{\id}{Y}{g}{Z}\right] \comp \left[\CorrHom{X}{f}{Y}{\id}{Y}\right] = \left[\CorrHom{X}{f}{Y}{g}{Z}\right],
\end{align}
which results in the following description of $\D$:

\begin{remark} \label{rmk:3ff-applied-to-general-morphism}
Given a 3-functor formalism $\D\colon \Corr(\cat C, E)^\tensor \to \Cat^\times$, the functor $\D$ acts on objects as $X \mapsto \D(X)$ and on morphisms as
\begin{align}
	\D\left(\CorrHom{X}{f}{Y}{g}{Z}\right) = g_! f^*\colon \D(X) \to \D(Y)
\end{align}
for all $f\colon Y \to X$ in $\cat C$ and $g\colon Y \to Z$ in $E$.
\end{remark}

It is sometimes useful to know that the 3-functors $\tensor$, $f^*$ and $f_!$ preserve certain limits or colimits, so in the following we provide a convenient terminology for this case:

\begin{definition} \label{def:3ff-compatible-with-lim-and-colim}
Let $I$ be a category and $\D$ a 3-functor formalism on some geometric setup $(\cat C, E)$. We say that $\D$ is \emph{compatible with $I$-indexed colimits} if the following is satisfied:
\begin{enumerate}[(i)]
	\item For all $X \in \cat C$, $\D(X)$ has $I$-indexed colimits and $\tensor$ preserves them in each argument.
	\item For all $f \in E$, $f^*$ and $f_!$ preserve $I$-indexed colimits.
\end{enumerate}
We similarly define the compatibility of $\D$ with $I$-indexed limits.
\end{definition}

\begin{remark}
A 3-functor formalism $\D$ on a geometric setup $(\cat C, E)$ is stable if and only if all $\D(X)$ are stable categories and $\D$ is compatible with finite limits (equivalently, finite colimits). A similar criterion exists for presentable 3-functor formalisms and is discussed in \cref{rslt:6ff-with-presentable-cat}. Note that every presentable 3-functor formalism is compatible with all small colimits.
\end{remark}

We now have a good understanding of how a 3-functor formalism encodes pullback and exceptional pushforward. It remains to explain how the tensor product is encoded, i.e.\ unpack \cref{def:tensor-product-from-3ff}. For simplicity, let us assume that $\cat C$ admits finite limits, so that $\Corr(\cat C, E)^\tensor$ is a symmetric monoidal category and $\D$ is a lax symmetric monoidal functor $\D\colon \Corr(\cat C, E) \to \Cat$ (see \cref{rmk:3ff-is-lax-symm-mon-functor-if-C-has-finite-limits}). This lax symmetric monoidal structure provides us with natural maps
\begin{align}
	\blank \boxtimes \blank\colon \D(X) \times \D(Y) \to \D(X \times Y)
\end{align}
for all $X, Y \in \cat C$. Then the tensor product on $\D(X)$ is defined as the composition
\begin{align}
	\blank \tensor \blank\colon \D(X) \times \D(X) \xlongto{\blank \boxtimes \blank} \D(X \times X) \xlongto{\Delta^*} \D(X),
\end{align}
where $\Delta\colon X \to X \times X$ denotes the diagonal map.

We have learned how a 3-functor formalism, as defined in \cref{def:3-functor-formalism}, encodes the three functors $\tensor$, $f^*$ and $f_!$. The following result shows that it also encodes the expected compatibilities of these three functors:

\begin{proposition} \label{rslt:compatibilities-of-3ff}
Let $\D\colon \Corr(\cat C, E)^\tensor \to \Cat^\times$ be a 3-functor formalism on a geometric setup $(\cat C, E)$. Then the three functors $\tensor$, $f^*$ and $f_!$ satisfy the following compatibilities:
\begin{propenum}
	\item \label{rslt:functoriality-of-pullback-except-pushforward} (Functoriality) The functors $f^*$ and $f_!$ are natural in $f$, i.e.\ for composable morphisms $f$ and $g$ in $\cat C$ there are natural isomorphisms
	\begin{align}
		(g \comp f)^* = f^* \comp g^*, \qquad (g \comp f)_! = g_! \comp f_!,
	\end{align}
	where the second one only applies if $f, g \in E$.

	\item (Compatibility of pullback and tensor) For every morphism $f\colon Y \to X$ in $\cat C$ the functor $f^*$ is symmetric monoidal, i.e.\ for all $M, N \in \D(X)$ there is a natural isomorphism
	\begin{align}
		f^*(M \tensor N) = f^* M \tensor f^* N.
	\end{align}

	\item \label{rslt:proper-base-change} (Proper base-change) For every cartesian square
	\begin{equation}\begin{tikzcd}
		Y' \arrow[r,"g'"] \arrow[d,"f'"] & Y \arrow[d,"f"] \\
		X' \arrow[r,"g"] & X
	\end{tikzcd}\end{equation}
	in $\cat C$ with $f \in E$ there is a natural isomorphism
	\begin{align}
		g^* f_! = f'_! g'^*
	\end{align}
	of functors $\D(Y) \to \D(X')$.

	\item \label{rslt:projection-formula} (Projection formula) For every $f\colon Y \to X$ in $E$ and all $M \in \D(X)$ and $N \in \D(Y)$ there is a natural isomorphism
	\begin{align}
		M \tensor f_! N = f_!(f^* M \tensor N).
	\end{align}
\end{propenum}
\end{proposition}
\begin{proof}
This is proved in \cite[Proposition~A.5.8]{Mann.2022a}. For the convenience of the reader we explain the proof of (i), (ii) and (iii).

Parts (i) and (ii) are obvious from the construction, as the functors $\D^*$ and $\D_!$ from \cref{def:3-functors-from-3ff} provide the desired compatibilities. To prove (iii), let the cartesian square be given as in the claim. Then we have
\begin{align}
	\left[\CorrHom{X}{g}{X'}{\id}{X'}\right] \comp \left[\CorrHom{Y}{\id}{Y}{f}{X}\right] = \left[\CorrHom{Y}{g'}{Y'}{f'}{X'}\right]
\end{align}
in $\Corr(\cat C, E)$ (more precisely, this isomorphism of maps is equivalently specified by the homotopy witnessing the pullback diagram in the claim). Applying $\D$ and using \cref{rmk:3ff-applied-to-general-morphism} we obtain an isomorphism of functors
\begin{align}
	(\id_! g^*) \comp (f_! \id^*) = f'_! g'^*,
\end{align}
which is exactly what we wanted. Part (iv) can be deduced in a similar way by applying $\D$ to a suitable diagram in $\Corr(\cat C, E)^\tensor$.
\end{proof}

One way to understand the projection formula is to say that for every map $f\colon Y \to X$ in $E$, the functor $f_!\colon \D(Y) \to \D(X)$ is $\D(X)$-linear, where we equip $\D(Y)$ with the $\D(X)$-linear structure induced by the symmetric monoidal functor $f^*\colon \D(X) \to \D(Y)$. In fact, the whole 3-functor formalism is automatically linear over the base, in the following sense:

\begin{lemma} \label{rslt:3ff-linear-over-base}
Let $\D\colon \Corr(\cat C, E)^{\tensor} \to \Cat^{\times}$ be a 3-functor formalism. Then $\D$ canonically enhances to a lax symmetric monoidal functor
\begin{align}
\D\colon \Corr(\cat C, E)^{\tensor} \to \LMod_{\D(*)}(\Cat)^{\tensor}.
\end{align}
\end{lemma}
\begin{proof}
Note that the unit $*$ in $\Corr(\cat C, E)^{\tensor}$ enhances to a commutative algebra by \cref{ex:trivial-algebra} and thus $\D(*)$ is a symmetric monoidal category (in fact, this is true for all objects of $\cat C$ by \cref{def:tensor-product-from-3ff}). By \cref{rslt:extension-of-monoidal-functor-to-modules} we obtain an operad map
\begin{align}
\LMod_*\bigl(\Corr(\cat C,E)\bigr)^{\tensor} \to \LMod_{\D(*)}(\Cat)^{\tensor}.
\end{align}
Precomposing with an inverse of the equivalence $\LMod_*(\Corr(\cat C,E))^{\tensor} \isoto \Corr(\cat C, E)^{\tensor}$ from \cref{ex:modules-over-trivial-algebra} yields the assertion.
\end{proof}

\subsection{6-functor formalisms} \label{sec:6ff.def}

In the previous subsection we defined 3-functor formalisms, encoding the functors $\tensor$, $f^*$ and $f_!$. We will now introduce the remaining three functors, finally arriving at the notion of a full 6-functor formalism. Fortunately, the three functors $\iHom$, $f_*$ and $f^!$ do not require any extra data to be defined, so a 6-functor formalism is just a special form of 3-functor formalism:

\begin{definition} \label{def:6-functor-formalism}
Let $(\cat C, E)$ be a geometric setup. A \emph{6-functor formalism} on $(\cat C, E)$ is a 3-functor formalism $\D\colon \Corr(\cat C, E)^\tensor \to \Cat^\times$ with the following properties:
\begin{enumerate}[(i)]
	\item For every $X \in \cat C$ the symmetric monoidal category $\D(X)$ is closed, i.e.\ for every $M \in \D(X)$ the functor $\blank \tensor M\colon \D(X) \to \D(X)$ admits a right adjoint $\iHom(M, \blank)$ (see \cref{def:internal-hom}). We call
	\begin{align}
		\iHom(\blank, \blank)\colon \D(X)^\op \times \D(X) \to \D(X)
	\end{align}
	the \emph{internal hom} functor.

	\item For every morphism $f\colon Y \to X$ in $\cat C$ the functor $f^*$ admits a right adjoint
	\begin{align}
		f_*\colon \D(Y) \to \D(X),
	\end{align}
	which we call the \emph{pushforward} functor.

	\item For every morphism $f\colon Y \to X$ in $E$ the functor $f_!$ admits a right adjoint
	\begin{align}
		f^!\colon \D(X) \to \D(Y),
	\end{align}
	which we call the \emph{exceptional pullback} functor.
\end{enumerate}
\end{definition}

In a 6-functor formalism, we get additional compatibilities of the new three functors, both among each other and in combination with the three functors from the 3-functor formalism. For the convenience of the reader, we list the commonly used properties:

\begin{proposition} \label{rslt:compatibilities-of-6ff}
Let $\D\colon \Corr(\cat C, E)^\tensor \to \Cat^\times$ be a 6-functor formalism on a geometric setup $(\cat C, E)$. In addition to \cref{rslt:compatibilities-of-3ff}, the six functors $\tensor$, $\iHom$, $f^*$, $f_*$, $f_!$ and $f^!$ satisfy the following compatibilities:
\begin{propenum}
	\item \label{rslt:functoriality-of-pushforward-and-except-pullback} (Functoriality) The functors $f_*$ and $f^!$ are natural in $f$, i.e.\ for composable morphisms $f$ and $g$ in $\cat C$ there are natural isomorphisms
	\begin{align}
		(g \comp f)_* = g_* \comp f_*, \qquad (g \comp f)^! = f^! \comp g^!,
	\end{align}
	where the second one only applies if $f, g \in E$.

	\item For every cartesian square
	\begin{equation}\begin{tikzcd}
		Y' \arrow[r,"g'"] \arrow[d,"f'"] & Y \arrow[d,"f"] \\
		X' \arrow[r,"g"] & X
	\end{tikzcd}\end{equation}
	in $\cat C$ with $g \in E$ there is a natural isomorphism
	\begin{align}
		g^! f_* = f'_* g'^!
	\end{align}
	of functors $\D(Y) \to \D(X')$.

	\item \label{rslt:enriched-adjunction-of-shriek-functors} For every $f\colon Y \to X$ in $E$ and all $M \in \D(X)$ and $N \in \D(Y)$ there is a natural isomorphism
	\begin{align}
		\iHom(f_! N, M) = f_* \iHom(N, f^! M).
	\end{align}

	\item For every $f\colon Y \to X$ in $E$ and all $M, N \in \D(X)$ there is a natural isomorphism
	\begin{align}
		f^! \iHom(N, M) = \iHom(f^* N, f^! M).
	\end{align}
\end{propenum}
\end{proposition}
\begin{proof}
Part (i) follows formally from \cref{rslt:functoriality-of-pullback-except-pushforward} using the fact that forming adjoints is natural (see \cite[\href{https://kerodon.net/tag/02ES}{Tag 02ES}]{kerodon}). Part (ii) follows formally from proper base-change (see \cref{rslt:proper-base-change}) by passing to right adjoints.

We now prove (iii), so let everything be given as in the claim. Then for every $P \in \D(X)$ we make the following computation using \cref{rslt:compatibilities-of-3ff}:
\begin{align}
	& \Hom(P, \iHom(f_! N, M)) = \Hom(P \tensor f_! N, M) = \Hom(f_!(f^* P \tensor N), M) =\\&\qquad= \Hom(f^* P \tensor N, f^! M) = \Hom(f^* P, \iHom(N, f^! M)) = \Hom(P, f_* \iHom(N, f^! M)),
\end{align}
hence we can conclude by Yoneda.

Part (iv) can be proved similarly to (iii) and is left as an exercise for the reader (see e.g.\ \cite[Proposition~23.3.(ii)]{Scholze.2018} for an analogous computation).
\end{proof}

\begin{definition}
Let $\D\colon \Corr(\cat C, E)^\tensor \to \Cat^\times$ be a 6-functor formalism. Using the functoriality of $f_*$ and $f^!$ from \cref{rslt:functoriality-of-pushforward-and-except-pullback} we make the following definitions:
\begin{defenum}
	\item We denote by $\D_*\colon \cat C \to \Cat$ the functor obtained from $\D^*$ by passing to right adjoints. Hence $\D_*(X) = \D(X)$ for $X \in \cat C$ and $\D_*(f) = f_*$ for morphisms $f$ in $\cat C$.

	\item We denote by $\D^!\colon \cat C_E^\op \to \Cat$ the functor obtained from $\D_!$ by passing to right adjoints. Hence $\D^!(X) = \D(X)$ for $X \in \cat C$ and $\D^!(f) = f^!$ for morphisms $f$ in $E$.
\end{defenum}
\end{definition}


\begin{remark}
The reader who is familiar with classical theories of 6-functor formalisms may wonder why the list of compatibilities in \cref{rslt:compatibilities-of-6ff} lacks any form of Poincaré duality, i.e.\ a relation between $f^*$ and $f^!$ for smooth maps $f$. Since we are working in a purely abstract setting, we do not have a notion of smooth maps and can therefore not hope to get Poincaré duality as a property of $\D$. However, we will see below (see \cref{rmk:cohomologically-smooth-maps}) that in every 6-functor formalism (even in a 3-functor formalism) there is a good notion of \enquote{cohomologically smooth maps} satisfying Poincaré duality. 
\end{remark}

In practice one often guarantees the existence of $\iHom$, $f_*$ and $f^!$ by employing the adjoint functor theorem, which constructs right adjoint functors for colimit-preserving functors between presentable categories (see \cref{rslt:adjoint-functor-theorem}). The following result solidifies this philosophy. For that result, recall that $\PrL$ denotes the category of presentable categories with colimit-preserving functors, which is equipped with the Lurie tensor product $\tensor$ (see \cref{sec:cat.pres}).

\begin{lemma} \label{rslt:6ff-with-presentable-cat}
Let $\D\colon \Corr(\cat C, E)^\tensor \to \Cat^\times$ be a 3-functor formalism. Then the following are equivalent:
\begin{lemenum}
\item $\D$ is presentable, i.e.\ factors over $(\PrL)^\tensor$.

\item The functors $\tensor$, $f^*$ and $f_!$ preserve all small colimits and $\D(X)$ is presentable for every $X \in \cat C$.
\end{lemenum}
Moreover, in this case $\D$ is automatically a 6-functor formalism and factors over $\LMod_{\D(*)}(\PrL)^\tensor$. In particular, the categories $\D(X)$, all 6 functors, the base-change and projection formulas etc.\ are enriched over $\D(*)$.
\end{lemma}
\begin{proof}
It is clear that (i) implies (ii) and \cite[Remark~4.8.1.9]{HA} shows that (ii) implies (i). By the adjoint functor theorem (see \cref{rslt:adjoint-functor-theorem}) every presentable 3-functor formalism is automatically a 6-functor formalism. The claim about $\LMod_{\D(*)}(\PrL)^\tensor$ follows in the same way as in \cref{rslt:3ff-linear-over-base}. For the claims about the enrichment, see \cref{ex:enrichment-of-linear-category,ex:lax-linear-equals-enriched}.
\end{proof}

\subsection{Construction of 6-functor formalisms} \label{sec:6ff.construct}

In the previous subsections we introduced 6-functor formalisms as a neat way of organizing the six functors $\tensor$, $\iHom$, $f^*$, $f_*$, $f_!$ and $f^!$ together with an infinite amount of higher homotopies exhibiting various compatibilities between these functors. In order for the notion of 6-functor formalisms to be useful we need some tools to construct them, as it is basically impossible to specify all of the higher homotopies by hand. The present subsection provides our main tool of constructing a 6-functor formalism \enquote{out of nowhere}, while the next subsection deals with extending an existing 6-functor formalism to a larger class of objects and morphisms.

Let us now come to the promised construction result. First of all, it is enough to provide a tool for constructing 3-functor formalisms, as a 6-functor formalism is just a 3-functor formalism with extra properties. In practice, one is often in the following situation: We are given a geometric setup $(\cat C, E)$ and a functor $\D\colon \cat C^\op \to \CMon$, associating with every object $X \in \cat C$ a symmetric monoidal category $\D(X)$ of \enquote{sheaves on $X$} and for every morphism $f\colon Y \to X$ in $\cat C$ a symmetric monoidal pullback functor $f^*\colon \D(X) \to \D(X)$. Moreover, the class $E$ decomposes into a class $P \subseteq E$ of \enquote{proper} maps and a class $I \subseteq E$ of \enquote{local isomorphisms} (roughly open immersions), so that every map in $E$ can be written as a composition of a proper map and a local isomorphism. Now we want to define a 3-functor formalism on $(\cat C, E)$ such that for every $f \in P$ the functor $f_!$ is right adjoint to $f^*$ and for every $j \in I$ the functor $j_!$ is left adjoint to $j^*$. Since every map in $E$ can be decomposed into maps in $P$ and $I$, these requirements should define the desired 3-functor formalism. The main result of this subsection asserts that such a 3-functor formalism indeed exists (under mild assumptions on the data above).

In order to state the promised existence result, let us first make precise what we mean by \enquote{$E$ decomposes into $P$ and $I$}.

\begin{definition} \label{def:right-cancellative}
Let $E$ be a collection of morphisms in a category $\cat C$. We say that $E$ is \emph{right cancellative} if for all morphisms $f\colon Y \to X$ and $g\colon Z \to Y$ in $\cat C$ with $f, f \comp g \in E$, we have $g \in E$.
\end{definition}

\begin{definition}[{cf. \cite[Definition~A.5.9]{Mann.2022a}}] \label{def:suitable-decomposition}
Let $(\cat C, E)$ be a geometric setup such that $\cat C$ admits pullbacks. A \emph{suitable decomposition} of $E$ is a pair $I, P \subseteq E$ satisfying the following properties:
\begin{enumerate}[(i)]
	\item $I$ and $P$ contain all identity morphisms and are stable under composition and pullback (i.e.\ they encode wide subcategories $\cat C_I, \cat C_P \subseteq \cat C_E$ that are stable under pullbacks in $\cat C$).

	\item Every $f \in E$ can be written as $f = p \comp i$ for some $p \in P$ and some $i \in I$.
	
	\item $I$ and $P$ are right cancellative in $\cat C$.

	\item Every morphism $f \in I \isect P$ is $n$-truncated for some $n \ge -2$ (possibly depending on $f$). Note that this is always satisfied if $\cat C$ is an ordinary category.
\end{enumerate}
\end{definition}

While conditions (i) and (ii) of \cref{def:suitable-decomposition} are very natural, conditions (iii) and (iv) may look weird at first glance and it is less clear why we need them. We give an explanation for this below. Let us first state the promised construction result:

\begin{proposition} \label{rslt:construct-3ff-from-suitable-decomp}
Let $(\cat C, E)$ be a geometric setup such that $\cat C$ admits pullbacks, let $I, P \subseteq E$ be a suitable decomposition of $E$ and let $\D\colon \cat C^\op \to \CMon$ be a functor. Assume that the following conditions are satisfied:
\begin{enumerate}[(a)]
	\item For every map $j\colon U \to X$ in $I$ the functor $j^*$ admits a left adjoint $j_!\colon \D(U) \to \D(X)$ satisfying base-change and the projection formula, i.e.:
	\begin{itemize}
		\item For every map $g\colon X' \to X$ in $\cat C$ with pullbacks $j'\colon U' \to X'$ and $g'\colon U' \to U$, the natural map $j'_! g'^* \isoto g^* j_!$ is an isomorphism of functors $\D(U) \to \D(X')$.

		\item For all $M \in \D(X)$ and $N \in \D(U)$ the natural map $j_!(j^* M \tensor N) \isoto M \tensor j_! N$ is an isomorphism.
	\end{itemize}

	\item For every map $f\colon Y \to X$ in $P$ the functor $f^*$ admits a right adjoint $f_*\colon \D(Y) \to \D(X)$ satisfying base-change and the projection formula, i.e.:
	\begin{itemize}
		\item For every map $g\colon X' \to X$ in $\cat C$ with pullbacks $f'\colon Y' \to X'$ and $g'\colon Y' \to Y$, the natural map $g^* f_* \isoto f'_* g'^*$ is an isomorphism of functors $\D(Y) \to \D(X')$.

		\item For all $M \in \D(X)$ and $N \in \D(Y)$ the natural map $M \tensor f_* N \isoto f_*(f^* M \tensor N)$ is an isomorphism.
	\end{itemize}

	\item For every cartesian diagram
	\begin{equation}\begin{tikzcd}
		U' \arrow[r,"j'"] \arrow[d,"f'"] & X' \arrow[d,"f"]\\
		U \arrow[r,"j"] & X
	\end{tikzcd}\end{equation}
	in $\cat C$ such that $j \in I$ and $f \in P$, the natural map $j_! f'_* \isoto f_* j'_!$ is an isomorphism of functors $\D(U') \to \D(X)$.
\end{enumerate}
Then $\D$ extends to a 3-functor formalism $\D\colon \Corr(\cat C, E)^\tensor \to \Cat$ such that for every $j \in I$ the functor $j_!$ is left adjoint to $j^*$ and for every $f \in P$ the functor $f_!$ is right adjoint to $f^*$.
\end{proposition}
\begin{proof}
This is shown in \cite[Proposition~A.5.10]{Mann.2022a} relying on powerful results by Liu--Zheng. We do not repeat the proof here, but below we provide some intuition for the conditions appearing in the claim.
\end{proof}

Let us explain \cref{rslt:construct-3ff-from-suitable-decomp}, especially why all the appearing conditions are required. It is clear that conditions (a) and (b) are necessary if one wants $j_!$ resp.\ $f_*$ to play the role of the exceptional pushforward (because the exceptional pushforward satisfies base-change and the projection formula by \cref{rslt:compatibilities-of-3ff}). Note that since in $j_!$ and $f_*$ are adjoint functors to $j^*$ and $f^*$, respectively, the base-change and the projection formula are merely \emph{conditions} and not extra data, because there are natural maps between both sides of the claimed identities. For example, the natural map $g^* f_* \to f'_* g'^*$ in (b) comes via adjunction from the map $f'^* g^* f_* = g'^* f^* f_* \to g'^*$, which in turn is induced by the counit map $f^* f_* \to \id$ from the adjunction of $f^*$ and $f_*$.

It is clear that condition (c) in \cref{rslt:construct-3ff-from-suitable-decomp} is necessary in order for the association $f \mapsto f_!$ to be functorial. Note that it is crucial in (c) that the provided diagram is a \emph{pullback} diagram in order to get a natural map $j_! f'_* \to f_* j'_!$. We need this natural map to exist in order for the desired equivalence of $j_! f'_*$ and $f_* j'_!$ being a condition rather than an extra datum. Note however that the conclusion of \cref{rslt:construct-3ff-from-suitable-decomp} implies that a posteriori there is a natural isomorphism of functors $j_! f'_* = f_* j'_!$ even if the diagram is not a pullback (i.e.\ just a commuting diagram with $j, j' \in I$ and $f, f' \in P$). Moreover, it is a priori not clear why different decompositions $f \comp j = f' \comp j'$ yield a natural isomorphism of functors $f_* \comp j_! = f'_* \comp j'_!$, although this also follows from the conclusion of \cref{rslt:construct-3ff-from-suitable-decomp}. It is natural to wonder how these isomorphisms are constructed naturally from the data provided in \cref{rslt:construct-3ff-from-suitable-decomp} -- this is why we need conditions (iii) and (iv) of \cref{def:suitable-decomposition}! For a detailed explanation how (iii) and (iv) allow us to construct the claimed compatibilities, we refer the reader to \cite[Constructions~4.4,~4.5]{Scholze:Six-Functor-Formalism}.

It turns out that condition (c) of \cref{rslt:construct-3ff-from-suitable-decomp} is very often automatic. Indeed, we have the following result:

\begin{lemma} \label{rslt:cond-c-satisfied-on-tensor-product}
In the setting of \cref{rslt:construct-3ff-from-suitable-decomp}, assume that conditions (a) and (b) are satisfied. Then for every diagram as in (c) and every $M \in \D(U')$ that can be written as $M = f'^* N \tensor j'^* N'$ for some $N \in \D(U)$ and $N' \in D(X')$, the natural map
\begin{align}
	j_! f'_* M \isoto f_* j'_! M
\end{align}
is an isomorphism.
\end{lemma}
\begin{proof}
It is not hard to see that both sides of the claimed isomorphism are indeed isomorphic, by the following computation involving (a) and (b):
\begin{align}
	j_! f'_* M = j_! f'_* (f'^* N \tensor j'^* N') = j_! (N \tensor f'_* j'^* N') = j_! (N \tensor j^* f_* N') = j_! N \tensor f_* N'.
\end{align}
A similar computation shows that $f_* j'_! M = j_! N \tensor f_* N'$, proving the desired isomorphism. However, from this computation it is not clear that the isomorphism we constructed is the one in the claim. For example, the computation made use of the projection formula for all four functors, whereas the claimed isomorphism only used the base-change property. In the following we provide a more conceptual view on the computation in order to verify that it does indeed construct the correct isomorphism.

We will employ higher algebra in the category $\Cat$ equipped with the cartesian symmetric monoidal structure, see \cref{ex:algebra-and-modules-in-Cat} for the basic properties. By assumption $f^*$ and $j^*$ are symmetric monoidal and hence define maps of algebras in $\Cat$. Their pushout is given by $\D(X') \tensor_{\D(X)} \D(U)$, so we obtain the following diagram, where the solid part commutes:
\begin{equation}\begin{tikzcd}[column sep=large,row sep=huge]
	\D(U') \arrow[drr,bend left,dashed,"f'_*"] \arrow[ddr,bend right=45,dashed,"j'_!",swap]\\
	& \D(X') \tensor_{\D(X)} \D(U) \arrow[ul,"h"] & \D(U) \arrow[l,"f^* \tensor \id"] \arrow[d,bend left,dashed,"j_!"] \arrow[ull,"f'^*",swap]\\
	& \D(X') \arrow[u,"\id \tensor j^*",swap] \arrow[r,dashed,bend right,"f_*"] \arrow[uul,"j'^*"] & \D(X) \arrow[l,"f^*",swap] \arrow[u,"j^*"]
\end{tikzcd}\end{equation}
The object $M = f'^* N \tensor j'^* N'$ lies in the essential image of $h$, so it is enough to show that the natural map $j_! f'_* h \isoto f_* j'_! h$ is an isomorphism. Note that by \cref{rslt:monoidal-right-adjoint,rslt:monoidal-left-adjoints} the projection formulas in (a) and (b) imply that $f_*$ and $j_!$ are $\D(X)$-linear, $f'_*$ is $\D(U)$-linear and $j'_!$ is $\D(X')$-linear and that the adjunctions to the pullback functors lift to adjunctions of functors of $\LM$-operads. We now claim that there is a natural isomorphism
\begin{align}
	f_* \tensor \id \isoto f'_* h.
\end{align}
Indeed, by the usual adjunctions of change-of-algebra (see \cref{ex:algebra-and-modules-in-Cat}) we have a natural isomorphism
\begin{align}
	\Hom_{\D(U)}(\D(X') \tensor_{\D(X)} \D(U), \D(U)) = \Hom_{\D(X)}(\D(X'), \D(U)),
\end{align}
where $\Hom_{\cat V}$ denotes the Hom anima in the category of $\cat V$-linear categories. The map from left to right is given by $\blank \comp (\id \tensor j^*)$ and hence it sends $f_* \tensor \id \mapsto j^* f_*$ and $f'_* h \mapsto f'_* j'^*$. The claimed isomorphism thus reduces to the base-change isomorphism $j^* f_* \isoto f'_* j'^*$. Note that there is a different way of constructing a morphism $f_* \tensor \id \to f'_* h$ by using the adjunctions of $(f^*, f_*)$ and $(f'^*, f'_*)$ and the isomorphism $f'^* = h (f^* \tensor \id)$. One checks that this is the same isomorphism as the one above (note that the used base-change isomorphism comes from the same adjunctions!). We similarly obtain a natural isomorphism
\begin{align}
	j'_! h \isoto \id \tensor j_!
\end{align}
coming from the base-change isomorphism $j'_! f'^* \isoto f^* j_!$ which is compatible with the obvious morphism coming from the adjunctions of $(j_!, j^*)$ and $(j'_!, j'^*)$. Note furthermore that there is an obvious isomorphism
\begin{align}
	j_! (f_* \tensor \id) = f_* (\id \tensor j_!)
\end{align}
(by writing all appearing categories as a tensor product of two categories over $\D(X)$). One checks that this isomorphism is compatible with the natural map from left to right coming out of the adjunctions of $(j_!, j^*)$ and $(f^*, f_*)$. By combining the above three isomorphisms we arrive at the claimed isomorphism $j_! f'_* h \isoto f_* j'_! h$.
\end{proof}

\begin{corollary} \label{rslt:cond-c-auto-for-fully-faithful-j-shriek}
In the setting of \cref{rslt:construct-3ff-from-suitable-decomp}, assume that conditions (a) and (b) are satisfied and assume furthermore that for all $j \in I$ the functor $j_!$ is fully faithful. Then condition (c) of \cref{rslt:construct-3ff-from-suitable-decomp} is automatically satisfied.
\end{corollary}
\begin{proof}
Suppose we are given a diagram as in condition (c) and let $M \in \D(U')$ be given. Then the full faithfulness of $j'_!$ implies that $M = f'^* \one \tensor j'^*(j'_! M)$, so \cref{rslt:cond-c-satisfied-on-tensor-product} implies that the natural map $j_! f'_* M \isoto f_* j'_! M$ is an isomorphism. As $M$ was arbitrary, we obtain the desired isomorphism of functors.
\end{proof}

\begin{remark} \label{rmk:cond-c-auto-for-monomorphism-j}
In the setting of \cref{rslt:cond-c-auto-for-fully-faithful-j-shriek}, note that $j_!$ is automatically fully faithful if $j\colon U \injto X$ is a monomorphism in $\cat C$. This follows from base-change and $U \times_X U = U$. In particular, in \cite[Remark~4.2]{Scholze:Six-Functor-Formalism} the existence of the \enquote{closed complement} is not necessary.
\end{remark}

\subsection{Extension of 6-functor formalisms} \label{sec:6ff.extend}

In the previous subsection we learned how to construct a 6-functor formalism on a geometric setup $(\cat C, E)$ by decomposing the maps in $E$ into \enquote{open immersions} and \enquote{proper maps}. In practice, such a decomposition may not always be possible for every map in $E$, only locally on source and target. By recalling and extending the results in the second part of \cite[\S A.5]{Mann.2022a} we present two powerful extension results, allowing one to extend a 6-functor formalism on a site $\cat C$ to the category of stacks on $\cat C$. We refer the reader to \cref{sec:sheaves} for an introduction to the language of sites, sheaves and descent for higher categories, which will be used freely in the following.


Let us come to our first extension result. It deals with the following situation: Suppose we have managed to construct a 6-functor formalism $\D$ on a site $\cat C$ and some nice class of maps $E$ in $\cat C$. We now want to extend $\D$ to a 6-functor formalism on all stacks on $\cat C$. The following result shows that this is basically always possible.

\begin{definition}
Let $(\cat C, E)$ be a geometric setup such that $\cat C$ is a site. A 3-functor formalism $\D$ on $(\cat C, E)$ is called \emph{sheafy} if the functor $\D^*\colon \cat C^\op \to \Cat$ is a sheaf of categories. A similar definition applies if $\cat C$ is a topos instead.
\end{definition}

\begin{proposition} \label{rslt:extend-3ff-to-stacks}
Let $(\cat C, E)$ be a geometric setup such that $\cat C$ is a subcanonical site. Let
\begin{itemize}
	\item $\cat X = \Shv(\cat C)$ be the topos of sheaves of anima on $\cat C$,
	\item $E'$ be the collection of those edges $f'\colon Y' \to X'$ in $\cat X$ such that for every map $g\colon X \to X'$ from an object $X \in \cat C$ the pullback of $f'$ along $g$ lies in $E$.
\end{itemize}
Then:
\begin{propenum}
 	\item The natural map $(\cat C, E) \to (\cat X, E')$ is a morphism of geometric setups.

	\item Every sheafy 3-functor formalism $\D$ on $(\cat C, E)$ extends uniquely to a sheafy 3-functor formalism $\D'$ on $(\cat X, E')$.
\end{propenum} 
A similar result is true for hypersheaves if one assumes that $\cat C$ is hyper-subcanonical.
\end{proposition}
\begin{proof}
We first show (i). Clearly $E'$ is stable under composition and pullback in $\cat X$. To show that $(\cat X, E')$ is a geometric setup it remains to verify that $E'$ is right cancellative (cf. \cref{rslt:right-cancellative-fiber-products}). Suppose we are given maps $f\colon Y \to X$ and $g\colon Z \to Y$ in $\cat X$ such that $f$ and $f \comp g$ lie in $E'$. We need to show that $g$ lies in $E'$, so pick an object $Y' \in \cat C$ together with a map $Y' \to Y$. We obtain the following diagram in $\cat X$, where the upper square is cartesian:
\begin{equation}\begin{tikzcd}
	Z \times_Y Y' \arrow[r] \arrow[d] & Z \times_X Y' \arrow[d] \\
	Y \times_Y Y' \arrow[r] & Y \times_X Y' \arrow[d]\\
	& X \times_X Y'
\end{tikzcd}\end{equation}
The left vertical map is the pullback of $g$ along $Y' \to Y$, so we need to show that it lies in $E$ (and in particular that $Z \times_Y Y' \in \cat C$). The bottom right vertical map is the pullback of $f$ along $Y' \to X$, so from $f \in E'$ we deduce that this map lies in $E$ and in particular that $Y \times_X Y' \in \cat C$. The composition of the right vertical maps is the pullback of $f \comp g$ along $Y' \to X$ and hence lies in $E$ as well (and thus $Z \times_X Y' \in \cat C$). Since $E$ is right cancellative we deduce that the upper right vertical map lies in $E$ and hence the left vertical map lies in $E'$. Hence, as $Y' = Y \times_Y Y' \in \cat C$, the left vertical map must lie in $E$ and we are done. This finishes the proof that $(\cat X, E')$ is a geometric setup.

Note that the natural functor $\cat C \injto \cat X$ is fully faithful because $\cat C$ is subcanonical. In order to verify that this functor defines a morphism of geometric setups, we need to see that it preserves pullbacks of maps in $E$; this is clear. We have thus finished the proof of (i).
	
It remains to prove (ii), for which we refer to \cite[\S A.5]{Mann.2022a}. Namely, by \cref{rslt:sheaves-on-topos-equiv-sheaves-on-site} the categories of sheaves of categories on $\cat C$ and on $\cat X$ coincide. Moreover, the functor from sheaves on $\cat C$ to sheaves on $\cat X$ is given by right Kan extension along the embedding $\cat C^\op \injto \cat X^\op$. Thus the claim reduces to \cite[Proposition~A.5.16]{Mann.2022a}. In the case of hypersheaves we can argue similarly (using \cite[Lemma~A.3.6]{Mann.2022a}).
\end{proof}

\begin{remark} \label{rmk:extend-6ff-to-stacks}
In the setting of \cref{rslt:extend-3ff-to-stacks} if $\D$ is a presentable 6-functor formalism, then the same is true for $\D'$. Indeed, by \cref{rslt:6ff-with-presentable-cat} the claim reduces immediately to the fact that the category $\PrL$ of presentable categories is stable under limits in $\Cat$ (and hence the right Kan extension of $\D$ lands in $\PrL$), which is shown in \cite[Proposition~5.5.3.13]{HTT}.
\end{remark}

\begin{remark}
Given a (subcanonical) site $\cat C$, the classical category of stacks on $\cat C$ can be identified with the category of $1$-truncated sheaves of anima on $\cat C$ and in particular is a full subcategory of $\Shv(\cat C)$. Thus \cref{rslt:extend-3ff-to-stacks} implies that a \enquote{sheafy} 3-functor formalism on $\cat C$ extends uniquely to a \enquote{sheafy} 3-functor formalism on the category of all stacks on $\cat C$ (and similarly for 6-functor formalisms as in \cref{rmk:extend-6ff-to-stacks}).
\end{remark}

\begin{examples}
If $\cat C$ does not come equipped with a site, here are two canonical choices of sites on $\cat C$ one can make:
\begin{exampleenum}
	\item The trivial site, where the only covers are of the form $\cat C_{/U}$ for $U \in \cat C$. In this case the sheafiness condition on $\D$ is vacuous and the topos $\cat X = \Shv(\cat C)$ identifies with the category of \emph{all} functors $\cat C^\op \to \Ani$.

	\item The \emph{$\D$-topology}, which is the site on $\cat C$ where covers are canonical covers that have \enquote{universal $*$- and $!$-descent} (these notions are defined below). This site is used in \cite[Definition~4.14]{Scholze:Six-Functor-Formalism}. Note that being canonical and having universal $*$-descent are the minimal requirements to apply \cref{rslt:extend-3ff-to-stacks}, while universal $!$-descent is not strictly necessary here, but can be very useful in applications.
\end{exampleenum}
\end{examples}

Let us spell out more explicitly what \cref{rslt:extend-3ff-to-stacks} does. Every $X' \in \cat X$ can be written as a colimit $X' = \varinjlim_{i\in I} X_i$ of representable objects $X_i \in \cat C$. For example, if $X'$ admits a cover $X_0 \surjto X'$ by some $X_0 \in \cat C$ such that all fiber products $X_0 \times_{X'} \dots \times_{X'} X_0$ are in $\cat C$, then we can write $X' = \varinjlim_{n \in \bbDelta^\op} X_n$, where $X_\bullet$ is the Čech nerve of $X_0 \to X'$. Since $\D'^*$ is supposed to be a sheaf of categories on $\cat X$, we are forced to define
\begin{align}
	\D'(X') = \varprojlim_i \D^*(X_i).
\end{align}
In the example of the Čech nerve above, this defines $\D'(X')$ as the category of descent data for $\D$ along the Čech nerve $X_\bullet$. We have explained the construction of the functor $\D'^*\colon \cat X^\op \to \Cat$. In order to see how the 3-functor formalism $\D'$ is constructed, it remains to see how the exceptional pushforward functors are defined. Let $f'\colon Y' \to X'$ be any map in $E'$ and write $X' = \varinjlim_i X_i$ as above. We define $Y_i := X_i \times_{X'} Y'$, so that $Y' = \varinjlim_i Y_i$ by general properties of sheaves. Moreover, each $Y_i$ lies in $\cat C$ and the induced map $f_i\colon Y_i \to X_i$ lies in $E$, so that $f_{i!}\colon \D(Y_i) \to \D(X_i)$ exists. We can now define the desired functor $f'_!$ componentwise, i.e.\ as
\begin{align}
	f'_! = f_{\bullet!}\colon \D(Y') = \varprojlim_i \D^*(Y_i) \to \varprojlim_i \D^*(X_i) = \D(X').
\end{align}
Note that in order for this definition to make sense, it is crucial to know that the functors $f_{i!}$ satisfy base-change with the transition functors in the above limits, which is guaranteed by proper base-change. The reader is invited to write down explicitly how $f'_!$ acts on descent data in the Čech nerve example.

We now come to the second extension result. So far we have explained how to extend a 6-functor formalism to stacks, but the allowed class $E$ of exceptional edges is quite restrictive: It consists only of maps which are representable in $\cat C$, i.e.\ no \enquote{stacky} maps are allowed. The next results show how to extend $E$. These results will make use of certain descent properties of the 6-functor formalism, so let us define these properties first.

\begin{definition}\label{def:D*-and-D!-covers}
Let $\D$ be a 3-functor formalism on a geometric setup $(\cat C, E)$.
\begin{defenum}
	\item We say that a sieve $\cat U$ in $\cat C$ is a \emph{(universal) $\D^*$-cover} if $\D^*$ descends (universally) along $\cat U$.
	\item Assume that $\D$ is a 6-functor formalism. We say that a sieve $\cat U \subseteq (\cat C_E)_{/U}$ in $\cat C_E$ is a \emph{small $\D^!$-cover} if it is generated by a small family of maps $(U_i \to U)_i$ and $\D^!$ descends along $\cat U$. We say that $\cat U$ is a small \emph{universal} $\D^!$-cover if for every map $V \to U$ in $\cat C$, the family $(U_i \times_U V \to V)_i$ generates a small $\D^!$-cover.
\end{defenum}
If $\D$ is clear from context we often abbreviate by saying $*$-cover and $!$-cover.
\end{definition}

\begin{remarks}
\begin{remarksenum}
	\item While universal $*$-covers often arise directly from the construction (e.g.\ as in \cref{rslt:extend-3ff-to-stacks}), it is less clear how to show that a sieve is a $!$-cover. We will provide tools for that in \cref{sec:kerncat.excdescent} below. For example, in \cref{rslt:!-descent-along-suave-map} we show that every cohomologically suave cover is a $!$-cover; in particular, open covers are usually $!$-covers.

	\item In order to check that a sieve is a $*$-cover (or $!$-cover), it is enough to find a universal $*$-cover (or $!$-cover) inside of it, see \cref{rslt:descent-along-sieve-from-subsieve}.
\end{remarksenum}
\end{remarks}

We can now present the main result for extending $E$. We start with the following technical result, which relies on left and right Kan extensions.

\begin{proposition} \label{rslt:Kan-extend-E}
Let $\D\colon \Corr(\cat C, E)^\tensor \to \Cat^\times$ be a 3-functor formalism on a geometric setup $(\cat C, E)$. Let $E' \supseteq E$ be another collection of morphisms in $\cat C$ such that $(\cat C, E')$ is a geometric setup and one of the following two conditions is satisfied:
\begin{propenum}
	\item \label{rslt:extend-E-on-the-target} For every $f\colon Y \to X$ in $E'$ there is a universal $\D^*$-cover $\cat U \subseteq \cat C_{/X}$ on $X$ such that the pullback of $f$ along every map $X' \to X$ in $\cat U$ lies in $E$.

	\item \label{rslt:extend-E-on-the-source} $\D$ is presentable such that for every $f\colon Y \to X$ in $E'$ there is a small universal $\D^!$-cover $\cat U \subseteq (\cat C_E)_{/Y}$ on $Y$ such that the composition of $f$ with every map $Y' \to Y$ in $\cat U$ lies in $E$.
\end{propenum}
Then $\D$ extends uniquely to a 3-functor formalism on $(\cat C, E')$.
\end{proposition}
\begin{proof}
We first prove (ii). By \cite[Lemma~A.5.11]{Mann.2022a} the desired extension of $\D$ exists uniquely if for every $f\colon Y \to X$ in $E'$, $\D^!$ descends along the sieve $\cat U \subset (\cat C_E)_{/Y}$ consisting of those maps $g\colon Y' \to Y$ such that $f\comp g \in E$. But by \cref{rslt:descent-along-sieve-from-subsieve} this reduces immediately to the assumption. Note that the smallness assumption on the $!$-cover is necessary to identify the limits in $\PrL$ and $\Cat$ (using \cref{rslt:descent-data-for-generated-sieve}); here \cite[Lemma~A.5.11]{Mann.2022a} is slightly inaccurate.

We now prove (i), for which we argue in a similar way as in \cite[Proposition~A.5.16]{Mann.2022a}. By functoriality of $\Corr$ there is a natural functor $\varphi\colon \Corr(\cat C, E)^\tensor \to \Corr(\cat C, E')^\tensor$. Moreover, by \cite[Proposition~2.4.1.7]{HA} we can identify 3-functor formalisms with lax cartesian structures on the above two operads. Thus in order to get the desired extension of $\D$, we perform a right Kan extension of $\D\colon \Corr(\cat C, E)^\tensor \to \Cat$ along $\varphi$. Let us denote the resulting functor by $\D'\colon \Corr(\cat C, E')^\tensor \to \Cat$. We need to show that for every $X_\bullet \in \Corr(\cat C, E)^\tensor$, the natural functor $\D'(X_\bullet) \isoto \D(X_\bullet)$ is an equivalence. We denote
\begin{align}
	\cat K := \Corr(\cat C, E)^\tensor_{X_\bullet/} = \Corr(\cat C, E)^\tensor \times_{\Corr(\cat C, E')^\tensor} \Corr(\cat C, E')^\tensor_{X_\bullet/}.
\end{align} 
Then $\D'(X_\bullet) = \varprojlim_{Y_\bullet \in \cat K} \D(Y_\bullet)$. The objects of $\cat K$ are diagrams
\begin{align}
	\CorrHom{(X_i)_{1\le i \le n}}{}{(Y'_j)_{1\le i \le m}}{}{(Y_j)_{1\le i\le m}},
\end{align}
such that the map to the right lies over the identity $\langle m \rangle \to \langle m \rangle$ in $\Fin_*$ and all the maps $Y'_j \to Y_j$ lie in $E'$. Morphisms in $\mathcal K$ are of similar shape, but here in addition we require containment in $E$ instead of $E'$. Let now $\cat K' \subseteq \cat K$ be the full subcategory spanned by those objects where each $Y'_j \to Y_j$ lies in $E$. Then $\cat K' = \Corr(\cat C, E)^\tensor_{X_\bullet/}$, which has the initial object $X_\bullet$ and hence $\varprojlim_{Y_\bullet \in \cat K'} \D(Y_\bullet) = \D(X_\bullet)$. It is therefore enough to show that $\D$ is the right Kan extension of itself along the embedding $\cat K' \injto \cat K$, as then the limits of $\D$ over $\cat K$ and $\cat K'$ agree. To prove this Kan extension property, let an object $f = [\CorrHom{X_\bullet}{}{Y'_\bullet}{}{Y_\bullet}] \in \cat K'$ be given. Then we need to show that
\begin{align}
	\varprojlim_{Z_\bullet \in \cat K'_{f/}} \D(Z_\bullet) = \D(Y_\bullet).
\end{align}
Note that $\cat K'_{f/}$ consists of morphisms $\CorrHom{Y_\bullet}{}{Z'_\bullet}{}{Z_\bullet}$ such that in the composition diagram
\begin{equation}\begin{tikzcd}[column sep=small, row sep=small]
	&& (Z''_i)_{1\le i \le k} \arrow[dr] \arrow[dl]\\
	& (Y'_i)_{1\le i \le m} \arrow[dr] \arrow[dl] && (Z'_i)_{1\le i \le k} \arrow[dr] \arrow[dl] \\
	(X_i)_{1\le i \le n} && (Y_i)_{1\le i \le m} && (Z_i)_{1\le i \le k}
\end{tikzcd}\end{equation}
all the maps $Z'_i \to Z_i$ and $Z''_i \to Z_i$ lie in $E$. By right cancellativeness, also the maps $Z''_i \to Z'_i$ lie in $E$. Let us denote by $\cat K'' \subseteq \cat K'_{f/}$ the full subcategory spanned by those objects $\CorrHom{Y_\bullet}{}{Z'_\bullet}{}{Z_\bullet}$ such that the map $Z'_\bullet \to Z_\bullet$ is degenerate. Then the inclusion $\cat K''^\op \injto (\cat K'_{f/})^\op$ is cofinal, by a similar argument as in the proof of \cite[Proposition~A.5.16]{Mann.2022a} (the crucial observation is that for an object $\CorrHom{Y_\bullet}{}{Z'_\bullet}{}{Z_\bullet}$ in $\cat K'_{f/}$ as above, also the morphism $\CorrHom{Y_\bullet}{}{Z'_\bullet}{\id}{Z'_\bullet}$ is in $\cat K'_{f/}$). We let $\cat K''' \subseteq \cat K''$ be the full subcategory where $Z_\bullet \to Y_\bullet$ lies over the identity in $\Fin_*$ and for $i = 1, \dots, m$ we let $\cat K'''_i \subseteq \cat K''$ be the full subcategory where $Z \to Y_\bullet$ lies over the map $\langle m \rangle \to \langle 1 \rangle$ that sends everything to $*$ but $i$. Then $\cat K''' \injto \cat K'_{f/}$ is cofinal and hence
\begin{align}
	\varprojlim_{Z_\bullet \in \cat K'_{f/}} \D(Z_\bullet) = \varprojlim_{Z_\bullet \in \cat K'''} \D(Z_\bullet) = \prod_{i=1}^m \varprojlim_{Z \in \cat K'''_i} \D(Z).
\end{align}
Now observe that $\cat K'''_i$ can be identified with the full subcategory of $\cat C_{/Y_i}$ consisting of those maps $Z \to Y_i$ such that the pullback of $Y'_i \to Y_i$ along $Z \to Y_i$ lies in $E$. By our assumption on $\D$ and \cref{rslt:descent-along-sieve-from-subsieve} we deduce that $\D$ descends along the sieve $\cat K'''_i$ and hence altogether we deduce
\begin{align}
	\varprojlim_{Z_\bullet \in \cat K'_{f/}} \D(Z_\bullet) = \prod_{i=1}^m \varprojlim_{Z \in \cat K'''_i} \D(Z) = \prod_{i=1}^m \D(Y_i) = \D(Y_\bullet),
\end{align}
as desired. 
\end{proof}

\begin{remark}
It is surprisingly hard to apply \cref{rslt:extend-E-on-the-target}, as it is not easy to construct a suitable class $E'$ out of $E$. Namely, if one just takes $E'$ to consist of all those maps that $\D^*$-locally on the target lie in $E$, then $E'$ may not be stable under composition.
\end{remark}

Let us unpack the construction in \cref{rslt:Kan-extend-E}. We first discuss the extension on the target, i.e.\ \cref{rslt:extend-E-on-the-target}. Given a map $f\colon Y \to X$ in $E'$, pick a $\D^*$-cover $(X_i \to X)_{i\in I}$ such that all the pulled back maps $Y_i \to X_i$ lie in $E$. Then the functor $f_!$ is defined as the composition
\begin{align}
	f_!\colon \D(Y) \to \varprojlim_i \D(Y_i) \xto{f_{\bullet!}} \varprojlim_i \D(X_i) = \D(X),
\end{align}
similarly to what happens in \cref{rslt:extend-3ff-to-stacks}. Note that we crucially use the assumption that $E$ is right cancellative (which appears in the definition of geometric setups, see \cref{rslt:right-cancellative-fiber-products}) in order to ensure that $f_!$ is functorial in $f$. Indeed, suppose we have another morphism $g\colon Z \to Y$ in $E'$. Then locally on $X$, both $f$ and $f \comp g$ lie in $E$. By right cancellativeness, also $g$ lies in $E$ (locally on $X$) and the desired identity $(f \comp g)_! = f_! \comp g_!$ comes from $\D$.

Let us now explain the construction in \cref{rslt:extend-E-on-the-source}. Suppose we are given a map $f\colon Y \to X$ in $E'$. Pick a small $!$-cover $\cat U \subseteq \cat (C_E)_{/Y}$ such that for all $Y' \to Y$ in $\cat U$, the composed map $Y' \to X$ lies in $E$. From the descent $\D^!(Y) = \varprojlim_{Y' \in \cat U} \D^!(Y')$ we deduce that for every $M \in \D(Y)$ we have $M = \varinjlim_{Y' \in \cat U} g_{Y'!} g_{Y'}^! M$, where $g_{Y'}\colon Y' \to Y$ denotes the obvious map. We can thus define
\begin{align}
	f_! M = \varinjlim_{Y' \in \cat U} (f \comp g_{Y'})_! g_{Y'}^! M.
\end{align}
Note that the smallness assumption on $\cat U$ guarantees that this colimit exists in $\D(X)$.

We now combine all of the above results in order to deduce a powerful result on extensions of 6-functor formalisms from sites to stacks and exceptional pushforwards along stacky morphisms. The following result is in the spirit of \cite[Theorem~4.20]{Scholze:Six-Functor-Formalism}, but in a slightly different formulation and with a few corrections. The reader should see the following result as the most common use case of the above extension results, but in specific situations it may be necessary to use the extension results directly.

\begin{definition}
Let $\D$ be a 6-functor formalism on a geometric setup $(\cat C, E)$, let $S$ be a collection of maps in $\cat C$ and let $f\colon Y \to X$ be a fixed map.
\begin{defenum}
	\item We say that $f$ is \emph{$\D^!$-locally on the source in $S$} if there is a universal small $\D^!$-cover of $Y$ such that the composition of $f$ with every element in the cover lies in $S$.

	\item We say that $f$ is \emph{$\D^!$-locally on the target in $S$} if there is a universal small $\D^!$-cover of $X$ such that the pullback of $f$ along every element in the cover lies in $S$.
\end{defenum}
If $\D$ is clear from context we simply write \enquote{$!$-local} in place of $\D^!$-local.
\end{definition}

\begin{theorem} \label{rslt:extend-6ff-to-stacks-and-stacky-maps}
Let $\D$ be a sheafy presentable 6-functor formalism on a geometric setup $(\cat C, E)$, where $\cat C$ is a subcanonical site. Then there is a collection of edges $E'$ in $\cat X \coloneqq \Shv(\cat C)$ with the following properties:
\begin{thmenum}
	\item The inclusion $\cat C \injto \cat X$ defines a morphism of geometric setups $(\cat C, E) \to (\cat X, E')$ and $\D$ extends uniquely to a sheafy 6-functor formalism $\D'$ on $(\cat X, E')$.

	\item $E'$ is $*$-local on the target: Let $f\colon Y \to X$ be a map in $\cat X$ whose pullback to every object in $\cat C$ lies in $E'$; then $f$ lies in $E'$.

	\item $E'$ is $!$-local: Let $f\colon Y \to X$ be a map in $\cat X$ which is $!$-locally on source or target in $E'$; then $f$ lies in $E'$.

	\item $E'$ is tame: Every map $f\colon Y \to X$ in $E'$ with $X \in \cat C$ is $!$-locally on the source in $E$.
\end{thmenum}
Moreover, there is a minimal choice of $E'$. The same result holds for hypercomplete sheaves in place of sheaves if one assumes that $\cat C$ is hyper-subcanonical.
\end{theorem}
\begin{proof}
It is clear that the collection of possible $E'$ is stable under intersection and hence admits a minimal choice. In the following we prove existence of $E'$.

Let $E_0$ be the collection of morphisms in $\cat X$ whose pullback to every $X \in \cat C$ lies in $E$. By \cref{rslt:extend-3ff-to-stacks} $\D$ extends uniquely to a sheafy 3-functor formalism $\D_0$ on $(\cat X, E_0)$, which by \cref{rmk:extend-6ff-to-stacks} is a 6-functor formalism such that all $\D_0(X)$ are presentable. We now proceed as in the proof of \cite[Theorem~4.20]{Scholze:Six-Functor-Formalism}. Let us denote by $A$ the class of all collections of edges $E'$ that satisfy (i) and (iv). Note that $A$ is stable under filtered unions and it is non-empty because it contains $E_0$. We first show the following claim:
\begin{itemize}
	\item[($*$)] Every $E' \in A$ can be enlarged to some $E'_! \in A$ which satisfies (iii). 
\end{itemize}
To prove this claim, fix some $E' \in A$. Let $E''$ be the collection of maps $f\colon Y \to X$ in $\cat X$ which are $\D'^!$-locally on the source in $E'$. We claim that $E''$ lies in $A$, i.e.\ $E''$ satisfies (i) and (iv). It is clear that $E''$ is stable under pullback in $\cat X$. Let us check that $E''$ is stable under composition, so let $f\colon Y \to X$ and $g\colon Z \to Y$ be two maps in $E''$. Pick a small family $(Y_i \to Y)_{i \in I}$ in $E'$ which is a $\D'^!$-cover after every pullback in $\cat X$ and such that all $Y_i \to X$ lie in $E'$. For each $i$ we denote by $g_i\colon Z_i \to Y_i$ the pullback of $g$ along $Y_i \to Y$. Then $g_i$ lies in $E''$ and hence we can find a small family $(Z_{ij} \to Z_i)_{j\in J_i}$ in $E$ (with $!$-descent as before) such that the composed maps $Z_{ij} \to Y_i$ lie in $E'$. Then all maps $Z_{ij} \to X$ lie in $E'$. Moreover, $(Z_{ij} \to Z)_{ij}$ is a small family of maps in $E$ which has $\D'^!$-descent after every pullback in $\cat X$ (by \cref{rslt:descent-along-sieve-from-subsieve}). This proves that $Z \to X$ lies in $E''$, as desired. One checks similarly that $E''$ is right cancellative and hence that $(\cat X, E'')$ is a geometric setup. The fact that $\D'$ extends uniquely to a 6-functor formalism on $(\cat X, E'')$ follows immediately from \cref{rslt:extend-E-on-the-source}. This proves that $E''$ satisfies (i). It is easy to see that it also satisfies (iv) and hence $E'' \in A$. By applying the same construction to $E''$ in place of $E'$, repeating this procedure and taking the union of all these classes, we arrive at the desired $E'_! \in A$, proving ($*$). We now prove the following claim:
\begin{itemize}
	\item[($**$)] Every $E' \in A$ satisfying (iii) can be enlarged to some $E'_* \in A$ satisfying (ii). 
\end{itemize}
To prove this claim we fix some $E' \in A$ satisfying (iii). Let $E'_*$ be the collection of maps $f\colon Y \to X$ whose pullback to every object in $\cat C$ lies in $E'$. We claim that $E'_* \in A$, i.e.\ $E'_*$ satisfies (i) and (iv). To prove this, we first show the following property of $E'_*$: If $f\colon Y \to X$ is in $E'_*$ such that $X$ admits an $E'$-map to some object in $\cat C$ then $f$ lies in $E'$. Namely, since $E'$ satisfies (iii), it is enough to show that $!$-locally on the target, $f$ lies in $E'$. But by (iv), $X$ lies $!$-locally in $\cat C$, so we conclude by the definition of $E'_*$. It follows easily that $E'_*$ is stable under composition, and it is clear that it is stable under pullback. By the same argument as in the proof of \cref{rslt:extend-3ff-to-stacks} we deduce that $E'_*$ is right cancellative and therefore $(\cat X, E'_*)$ is a geometric setup. By \cref{rslt:extend-E-on-the-target}, $\D'$ extends uniquely to $(\cat X, E'_*)$, hence $E'_*$ satisfies (i). It is clear that $E'_*$ satisfies (iv) and thus $E'_* \in A$, proving ($**$).

From claims ($*$) and ($**$) it is now easy to verify the existence of some $E'$ satisfying (i)--(iv). Namely, start with $E_0 \in A$, then define $E_{n+1} \coloneqq E_{n!*}$ for $n \ge 0$ and take $E'$ to be the union of all $E_n$.
\end{proof}

\begin{remark} \label{rmk:corrections-to-Scholze-Thm-4-20}
\Cref{rslt:extend-6ff-to-stacks-and-stacky-maps} is very similar to \cite[Theorem~4.20]{Scholze:Six-Functor-Formalism}, but with two small corrections:
\begin{itemize}
	\item We require that $E'$ is right cancellative, so that $\cat X_{E'}$ has pullbacks and the functor $\cat X_{E'} \to \cat X$ commutes with these pullbacks. This is in order to ensure that $!$-covers behave as expected (e.g.\ \cref{rslt:descent-data-for-generated-sieve} applies). In particular this is required to apply \cite[Proposition~A.5.14]{Mann.2022a} in the way it is claimed in \cite[Theorem~4.20]{Scholze:Six-Functor-Formalism}. While we include right cancellativeness for necessity, it is also a very useful property to have in practice.

	\item We stick to covering sieves and covering families instead of reducing to maps of the form $\bigdunion_i U_i \to V$ for $U_i \to V$ in $E'$. This is because in order to show that $E'$ can be enlarged to contain such a map, one needs to know a priori that for all $i_0 \in I$ the map $j_{i_0}\colon U_{i_0} \to \bigdunion_i U_i$ lies in $E'$, which is not clear to us. In particular we do not see why one can always enlarge $E'$ to be stable under disjoint unions (as defined in \cite[Definition~4.18.(1)]{Scholze:Six-Functor-Formalism}). Note that $j_{i_0}$ is $*$-locally on the target in $E'$, so one may be tempted to apply \cref{rslt:extend-E-on-the-target}; however, we were not able to come up with a good class of maps containing $j_{i_0}$ and being stable under composition.
\end{itemize}
\end{remark}

As discussed in \cref{rmk:corrections-to-Scholze-Thm-4-20} we do not know how to guarantee that maps of the form $X \injto X \dunion Y$ lie in $E'$. However, if we know that they do, we can derive the following properties:

\begin{lemma} \label{rslt:extend-6ff-to-disjoint-unions}
In the setting of \cref{rslt:extend-6ff-to-stacks-and-stacky-maps}, assume that the map $* \to * \dunion *$ lies in $E'$, where $*$ denotes the final object of $\cat X$ and the map is given by inclusion of the first component. Then:
\begin{lemenum}
	\item For every collection of objects $(X_i)_i$ in $\cat X$ with disjoint union $X = \bigdunion_i X_i$, every inclusion $X_i \injto X$ lies in $E'$ and the family of maps $(X_i \to X)_i$ is a universal $!$-cover.

	\item For every family of maps $(Y_i \to X)_i$ in $E'$, also $\bigdunion_i Y_i \to X$ is in $E'$.

	\item For every family of maps $(Y_i \to X_i)_i$ in $E'$, also $\bigdunion_i Y_i \to \bigdunion_i X_i$ is in $E'$.
\end{lemenum}
\end{lemma}
\begin{proof}
Claims (ii) and (iii) follow immediately from (i) using $!$-locality of $E'$. To prove (i), let $(X_i)_i$ be given. Note that for fixed $i$ we can write $\bigdunion_i X_i = X_i \dunion X'$ for $X' = \bigdunion_{j\ne i} X_j$ and the inclusion map $f_i\colon X_i \injto X$ is the pullback of the map $* \to * \dunion *$ along the projection $X_i \dunion X' \to * \dunion *$. Since $E'$ is stable under pullbacks and $* \to * \dunion *$ lies in $E'$ by assumption, we deduce that all $f_i$ lie in $E'$. Moreover, by sheafiness of $\D'^*$ we know that there is an isomorphism $\D'(X) \isoto \prod_i \D'(X_i)$ induced by the maps $f_i^*$. We need to show that this isomorphism is also induced by the maps $f_i^!$, for which it is enough to show that $f_i^! \isom f_i^*$ for each $i$. This further reduces to showing that for each $i$ the functor $f_{i!}\colon \D'(X_i) \to \D'(X)$ is isomorphic to the functor $\D'(X_i) \to \prod_j \D'(X_j)$ sending $M \mapsto (M_j)_j$ with $M_j = \emptyset$ (the initial object) if $j \ne i$ and $M_i = M$. But note that by general properties of topoi, for $i \ne j$ we have the following pullback diagrams in $\cat X$:
\begin{equation}
	\begin{tikzcd}
		X_i \arrow[r,"\id"] \arrow[d,"\id"] & X_i \arrow[d,"f_i"] \\
		X_i \arrow[r,"f_i"] & X
	\end{tikzcd}
	\qquad
	\begin{tikzcd}
		\emptyset \arrow[r] \arrow[d] & X_i \arrow[d,"f_i"] \\
		X_j \arrow[r,"f_j"] & X
	\end{tikzcd}
\end{equation}
By proper base-change there are natural isomorphisms $f_i^*(f_{i!} M) \isom M$ and $f_j^*(f_{i!} M) = \emptyset$, as desired.
\end{proof}

\subsection{Examples} \label{sec:6ff.example}

In the previous subsections we provide several powerful tools to construct 6-functor formalisms. We now discuss several applications of this machine by providing examples of 6-functor formalisms. In this paper we are particularly interested in the 6-functor formalism of $\Lambda$-valued sheaves on condensed anima, hence we focus on the construction of this 6-functor formalism. A condensed anima is roughly a \enquote{stack of topological spaces} (see \cref{rslt:embed-Top-into-Cond-Ani} below) and our 6-functor formalism will associate the category of $\Lambda$-valued sheaves on it. By restricting the 6-functor formalism to topological spaces, we will recover classical results from algebraic topology, see \cref{sec:kerncat.examples} below. By restricting the 6-functor formalism to classifying stacks of locally profinite groups we will obtain results about representation theory, see \cref{sec:reptheory}. At the end of this subsection we will also sketch the construction of other examples of 6-functor formalisms, with appropriate references to the literature.

Without further ado, let us come to the construction of the 6-functor formalism for $\Lambda$-valued sheaves on condensed anima. In the following we use some basic results from condensed mathematics introduced by Clausen--Scholze in \cite{Clausen-Scholze:Condensed-Mathematics,Clausen-Scholze:Analytic-Geometry}. We refer the reader to \cite[\S1, \S2]{Clausen-Scholze:Condensed-Mathematics} for an introduction to condensed methods and to \cite[\S11]{Clausen-Scholze:Analytic-Geometry} for some basic properties and constructions with condensed anima. In the following we recall the definition of condensed anima.

\begin{definition}
\begin{defenum}
	\item We denote by $\ProFin \coloneqq \Pro(\Fin)$ the $\Pro$-category of the category of finite sets. This is an ordinary category whose objects are cofiltered diagrams $S = (S_i)_{i \in I}$ of finite sets $S_i$ and whose morphisms are given by
	\begin{align}
		\Hom_{\ProFin}((T_j)_{j \in J}, (S_i)_{i \in I}) = \varprojlim_{i \in I} \varinjlim_{j \in J} \Hom_\Fin(T_j, S_i).
	\end{align}
	The objects of $\ProFin$ are called \emph{profinite sets}.

	\item The inclusion $\Fin \injto \Set$ induces a functor $\alpha\colon \ProFin \to \Set, (S_i)_i \mapsto \varprojlim_i S_i$. We call $\alpha(S)$ the \emph{underlying set} of a profinite set $S$. Given an uncountable strong limit cardinal $\kappa$, we say that a profinite set $S$ is \emph{$\kappa$-small} if its underlying set is a $\kappa$-small set. We denote
	\begin{align}
		\ProFin_\kappa \subset \ProFin
	\end{align}
	the full subcategory spanned by the $\kappa$-small profinite sets.
\end{defenum}
\end{definition}

\begin{remark} \label{rmk:topology-on-profinite-sets}
The assignment $(S_i)_i \mapsto \varprojlim_i S_i$ defines a fully faithful functor from $\ProFin$ to the category of topological spaces whose essential image is given by the totally disconnected compact Hausdorff spaces. The fully faithfulness is a standard exercise in topology (using the above description of $\Hom$-sets in $\ProFin$ and \cite[\href{https://stacks.math.columbia.edu/tag/08ZZ}{Tag 08ZZ}]{stacks-project}) and the essential image is described in \cite[\href{https://stacks.math.columbia.edu/tag/08ZY}{Tag 08ZY}]{stacks-project}.
\end{remark}

We record the following result which we could not find in the literature. It will be used in the computation of the cohomology of the unit interval, but is interesting on its own.

\begin{lemma}\label{rslt:profinite-completion}
The inclusion $\ProFin \injto \CHaus$ into the category of compact Hausdorff spaces admits a left adjoint $\pi_0$, whose underlying set is the set of connected components. Moreover, $\pi_0$ commutes with cofiltered limits.
\end{lemma}
\begin{proof}
The existence of $\pi_0$ follows from \cite[\href{https://stacks.math.columbia.edu/tag/08ZL}{Tag 08ZL}]{stacks-project} together with \cite[\href{https://stacks.math.columbia.edu/tag/0900}{Tag 0900}]{stacks-project}. Concretely, given a compact Hausdorff space $X$, the profinite space $\pi_0(X)$ can be described as follows: Let $I$ be the filtered poset whose elements are the equivalence relations ${\sim_i} \subseteq X\times X$ such that the quotient space $X/\sim_i$ is finite and discrete. The ordering is given by reverse inclusion. Then $\pi_0(X) \coloneqq \varprojlim_{i} X/{\sim_i}$ has the required universal property.

We now prove that $\pi_0$ preserves cofiltered limits. Let $X = \varprojlim_{i\in I} X_i$ be a cofiltered limit in $\CHaus$ (in view of \cite[\href{https://kerodon.net/tag/02QA}{Tag 02QA}]{kerodon} we may assume for convenience that $I$ is a directed partially ordered set). Since $S \coloneqq \{0,1\}$ is a cocompact cogenerator for $\ProFin$ (it is a cogenerator for $\Fin$, hence also for $\ProFin$, and clearly every finite set is cocompact), we may check that $\pi_0(X) \to \varprojlim_i \pi_0(X_i)$ is an isomorphism after applying $\Hom(\blank, S)$. By adjointness it therefore suffices to show that the natural map
\begin{align}\label{eq:profinite-completion}
\varinjlim_i \Hom(X_i, S) \to \Hom(X,S)
\end{align}
is bijective. Denote $p_i\colon X\to X_i$ and $p_{ij}\colon X_j \to X_i$ the projections and recall the following standard facts:
\begin{enumerate}[(i)]
\item If all $X_i$ are non-empty, then $X$ is non-empty.
\item $p_i(X) = \bigcap_{j\ge i} p_{ij}(X_j)$ for every $i\in I$.
\item Every compact open subset of $X$ is of the form $p_i^{-1}(U_i)$ for some $i$ and some open subset $U_i \subseteq X_i$.
\end{enumerate}
Parts (i) and (iii) are \cite[\href{https://stacks.math.columbia.edu/tag/0A2R}{Tag 0A2R}]{stacks-project} and \cite[\href{https://stacks.math.columbia.edu/tag/0A2P}{Tag 0A2P}]{stacks-project}, respectively. In (ii), the inclusion \enquote{$\subseteq$} is clear; for the other inclusion, pick a point $x$ in the intersection. Then each $p_{ij}^{-1}(x)$ is non-empty and compact, hence $\varprojlim_{j\ge i} p_{ij}^{-1}(x) \neq \emptyset$ by (i), which shows $x \in p_i(X)$.

We now prove that the map \eqref{eq:profinite-completion} is surjective, so let $f\in \Hom(X,S)$. Then $f$ corresponds to a partition $X = U \dunion V$ into clopen subsets. By (iii) and the fact that $I$ is cofiltered, we find $i\in I$ and open subsets $U_i, V_i \subseteq X_i$ such that $U = p_i^{-1}(U_i)$ and $V = p_i^{-1}(V_i)$. Note that $U_i$ and $V_i$ are necessarily disjoint. Denoting $Z_i$ the closed complement of $U_i \dunion V_i$, we deduce from $Z_i\cap p_i(X) = \emptyset$ and (ii) that $\bigcap_{j\ge i}Z_i \cap p_{ij}(X_j) = \emptyset$. Since $Z_i$ is compact and $I$ is cofiltered, we find $j\ge i$ such that $Z_i \cap p_{ij}(X_j) = \emptyset$. Now the open subsets $U_j\coloneqq p_{ij}^{-1}(U_i)$ and $V_j\coloneqq p_{ij}^{-1}(V_i)$ satisfy $X_j = U_j \dunion V_j$, $p_j(U) \subseteq U_j$ and $p_j(V) \subseteq V_j$. But this means that $f$ factors over $X_j$ proving surjectivity.

We next show injectivity. Let $f, g \in \Hom(X_i, S)$ for some $i\in I$, whose images in $\Hom(X,S)$ agree. Again, $f$ and $g$ correspond to clopen partitions $U^f_i \dunion V^f_i = X_i = U^g_i \dunion V^g_i$, such that $U\coloneqq p_i^{-1}(U^f_i) = p_i^{-1}(U^g_i)$ and $V\coloneqq p_i^{-1}(V^f_i) = p_i^{-1}(V^g_i)$ constitute a partition of $X$. Denoting $Z_i\coloneqq (U^f_i \cap V^g_i) \union (V^f_i \cap U^g_i)$, we have $Z_i \cap p_i(X) = \emptyset$, so that arguing as above we find $j\ge i$ with $Z_i \cap p_{ij}(X_j) = \emptyset$. But this means $p_{ij}^{-1}(U_i^f) = p_{ij}^{-1}(U_i^g)$ and $p_{ij}^{-1}(V_i^f) = p_{ij}^{-1}(V_i^g)$. In other words we have proved $fp_{ij} = gp_{ij}$, i.e.\ the map \eqref{eq:profinite-completion} is injective.
\end{proof}

\begin{definition} \label{def:condensed-anima}
We implicitly fix an uncountable strong limit cardinal $\kappa$. We equip $\ProFin_\kappa$ with the site where covers are those sieves $(\ProFin_\kappa)_{/S}$ that contain a finite family of maps $(S_n \to S)_n$ which are jointly surjective on the underlying sets. We denote by
\begin{align}
	\Cond(\Ani) \coloneqq \HypShv(\ProFin_\kappa)
\end{align}
the category of $\Ani$-valued hypersheaves on $\ProFin_\kappa$. The objects of $\Cond(\Ani)$ are called the \emph{condensed anima}.
\end{definition}

\begin{remark}
Fixing the strong limit cardinal $\kappa$ in \cref{def:condensed-anima} may seem unnatural and is not strictly necessary. It can be avoided by taking the colimit of $\Cond(\Ani)$ over all possible choices of $\kappa$ (cf. \cite[Appendix to Lecture II]{Clausen-Scholze:Condensed-Mathematics}). However, the restriction to fixed $\kappa$ simplifies some aspects of the theory (e.g.\ because $\ProFin_\kappa$ is a site while $\ProFin$ is not) and in practice it does not seem to make a difference. 
\end{remark}

Like for any topos, there are certain finiteness (or \enquote{coherence}) conditions that one can impose on condensed anima, which we recall now:

\begin{definition}
Fix $\kappa$ as in \cref{def:condensed-anima} and let $X$ be a condensed anima.
\begin{defenum}
	\item We say that a map $Y \surjto X$ of condensed anima is \emph{injective} (resp.\ \emph{surjective}) if it is a monomorphism (resp.\ effective epimorphism). We say that a family $(U_i \to X)_i$ of maps of condensed sets is a \emph{covering} if the induced map $\bigdunion_i U_i \to X$ is surjective.

	\item We say that $X$ is \emph{quasicompact} if every covering family $(Y_i \to X)_i$ admits a finite subcovering.

	\item \label{def:quasiseparated-condensed-anima} We say that $X$ is \emph{quasiseparated} if for any two maps $Y, Z \to X$ with $Y$ and $Z$ quasicompact condensed anima, also $Y \times_X Z$ is quasicompact.

	\item We say that $X$ is \emph{qcqs} if it is quasicompact and quasiseparated.
\end{defenum}
\end{definition}

A map $Y \to X$ is injective if and only if for all ($\kappa$-small) profinite sets $S$ the map $f\colon Y(S) \to X(S)$ is a monomorphism in $\Ani$, i.e.\ induces an embedding of connected components. A map $Y \to X$ is surjective if and only if for every ($\kappa$-small) profinite set $S$ and every section $s \in X(S)$ there are a ($\kappa$-small) profinite set $T$ with a surjective map $T \to S$ and a section $t \in Y(T)$ such that $f(t) = s|_T$ (see \cite[Lemma~A.3.9]{Mann.2022a}). A condensed anima $X$ is quasicompact if and only if there exists a surjective map $S \surjto X$ from some ($\kappa$-small) profinite set $S$.

A condensed anima can intuitively be thought of as a homotopy type $X$ whose homotopy groups $\pi_i(X, x)$ are equipped with a topological structure. In particular, $0$-truncated condensed anima (i.e.\ condensed \emph{sets}) are closely related to topological spaces:

\begin{proposition} \label{rslt:embed-Top-into-Cond-Ani}
Fix $\kappa$ as in \cref{def:condensed-anima}. The assignment $X \mapsto [S \mapsto \Hom(S, X)]$ defines a fully faithful embedding
\begin{align}
	\{ \text{$\kappa$-compactly generated topological spaces} \} \injto \Cond(\Set) \subset \Cond(\Ani)
\end{align}
which preserves finite products and induces an equivalence
\begin{align}
	\{ \text{$\kappa$-small compact Hausdorff spaces} \} \isoto \{ \text{qcqs condensed sets} \}.
\end{align}
Here $\Hom(S, X)$ denotes the set of continuous maps from $S$ to $X$, where we identify $S$ with the topological space from \cref{rmk:topology-on-profinite-sets}.
\end{proposition}
\begin{proof}
For the first part, see \cite[Proposition~1.7]{Clausen-Scholze:Condensed-Mathematics}. We sketch the argument: One first checks that the assignment $X \mapsto \underline X \coloneqq [S \mapsto \Hom(S, X)]$ does indeed define a functor, i.e.\ $\underline X$ is a sheaf on $\ProFin_\kappa$; this statement works also for non-compactly generated spaces $X$. Now given two $\kappa$-compactly generated topological spaces $X$ and $Y$ and their associated condensed sets $\underline X$ and $\underline Y$, we need to show that $\Hom(X, Y) = \Hom(\underline X, \underline Y)$. By definition, a map $\underline X \to \underline Y$ is a natural collection of maps $\Hom(S, X) \to \Hom(S, Y)$ for all $\kappa$-small profinite sets $S$. Taking $S = *$ we see that a map $\underline X \to \underline Y$ has an associated underlying map of \emph{sets} $X \to Y$. This easily implies that the assignment $X \mapsto \underline X$ is faithful. For fullness one uses the $\kappa$-compact generation, because it implies that a map $X \to Y$ is continuous as soon as the composition $S \to X \to Y$ is continuous for all maps $S \to X$ from a $\kappa$-small profinite set $S$.

For the identification of compact Hausdorff spaces and qcqs condensed sets, see \cite[Theorem~2.16]{Clausen-Scholze:Condensed-Mathematics}. The idea is to write a compact Hausdorff space as the quotient of a profinite set by a closed equivalence relation.
\end{proof}

In the following we will implicitly view every $\kappa$-compactly generated topological space as a condensed set and in particular as a condensed anima. As a special case, every $\kappa$-small profinite set is seen as a condensed anima (via the Yoneda embedding). In this paper we will be mostly interested in topological spaces that can be written as a disjoint union of profinite sets (most importantly, locally profinite groups). From now on we will mostly drop the explicit mentioning of $\kappa$.

\begin{example}
Since the embedding in \cref{rslt:embed-Top-into-Cond-Ani} commutes with finite products, it identifies ($\kappa$-)compactly generated topological groups with group objects in $\Cond(\Ani)$. In particular every locally profinite group $G$ can be seen as a group in $\Cond(\Ani)$ (i.e.\ an algebra for the cartesian symmetric monoidal structure). If $G$ acts on a condensed set $X$ then we can form the quotient $X/G$ in $\Cond(\Ani)$. As we will see in \cref{sec:reptheory} such quotients have a tight relationship to the theory of smooth representations of $G$.
\end{example}

Having introduced the category of condensed anima, we now come to the construction of the promised 6-functor formalism of $\Lambda$-valued sheaves. With the results of the previous subsections, this is straightforward to do. Let us first define the category of $\Lambda$-valued sheaves:

\begin{definition} \label{def:Einfty-ring-associated-to-profin-set}
Let $\Lambda$ be an $\Einfty$-ring (see \cref{ex:Einfty-rings-and-modules}). Then right Kan extension along the inclusion $\{ * \} \injto \Fin$ produces a functor $\Fin^\op \to \CAlg$ sending $S \mapsto \Lambda(S) = \prod_{x\in S} \Lambda$. By using the universal property of the $\Ind$-category we extend this to a functor
\begin{align}
	\ProFin^\op \to \CAlg, \qquad S = (S_i)_i \mapsto \Lambda(S) \coloneqq \varinjlim_i \Lambda(S_i).
\end{align}
We view $\Lambda(S)$ as the ring of continuous maps $S \to \Lambda$, where $\Lambda$ is equipped with the discrete topology (if $\Lambda$ is an ordinary ring, then this identification is literally true). By composing the functor $S \mapsto \Lambda(S)$ with $\Mod_{(\blank)}$ we obtain the functor
\begin{align}
	\D(\blank, \Lambda)\colon \ProFin^\op \to \Cat, \qquad S \mapsto \D(S, \Lambda) \coloneqq \Mod_{\Lambda(S)}.
\end{align}
\end{definition}

If $\Lambda$ is an ordinary ring (viewed as an $\Einfty$-ring) and $S = \varprojlim_i S_i$ is a profinite set, then $\Lambda(S) = \varinjlim_i \Lambda(S_i)$ can canonically be identified with the ring of continuous maps $S \to \Lambda$. Hence $\Lambda(S)$ agrees with the global sections of the constant sheaf $\Lambda$ on the topological space $S$. The next result shows that this identification extends to an equivalence of $\D(S,\Lambda)$ with the category of $\Lambda$-valued sheaves on $S$ (see \cref{rslt:identify-6ff-with-sheaves-on-top-spaces} for a similar result applied to locally compact Hausdorff spaces instead of profinite sets).

\begin{proposition} \label{rslt:identify-sheaves-on-profinite-sets-with-D}
Let $\Lambda$ be an $\Einfty$-ring. Then for every profinite set $S$ there is a natural equivalence (in $S$)
\begin{align}
	\D(S,\Lambda) = \Shv(S,\Mod_\Lambda),
\end{align}
where the right-hand side denotes the category of $\Mod_\Lambda$-valued sheaves (cf. \cref{def:sheaf}) on the site of open subsets of $S$.
\end{proposition}
\begin{proof}
For every map $f\colon Y \to X$ of topological spaces we get an induced pushforward functor $f_*\colon \Shv(Y,\Mod_\Lambda) \to \Shv(X,\Mod_\Lambda)$ given by $(f_* \mathcal M)(U) = \mathcal M(f^{-1}(U))$: This is clear on the level of presheaves and one checks easily that $f_*$ preserves the sheaf property. Using \cref{rslt:sheafification-exists} and the adjoint functor theorem, we obtain a left adjoint $f^{-1}\colon \Shv(X,\Mod_\Lambda) \to \Shv(Y,\Mod_\Lambda)$. This construction is clearly functorial in $f$ and so we obtain a functor
\begin{align}
	\ProFin^\op \to \PrL, \qquad S \mapsto \Shv(S,\Mod_\Lambda).
\end{align}
We want to show that this functor is isomorphic to the functor $\D(\blank,\Lambda)$. This is clearly true when restricted to finite sets $S$. Moreover, since $\Lambda(\blank)$ preserves filtered colimits, the same is true for $\D(\blank,\Lambda)$ (for colimits in $\PrL$), see \cref{ex:Einfty-rings-and-modules}. By the universal property of $\ProFin^\op = \Ind(\Fin^\op)$ (see \cref{rslt:universal-property-of-Ind-kappa}) we obtain a natural transformation
\begin{align}
	\alpha\colon \D(\blank,\Lambda) \to \Shv(\blank,\Mod_\Lambda)
\end{align}
of functors $\ProFin^\op \to \PrL$.

Now fix a profinite set $S$ and let $\beta_S\colon \Shv(S,\Mod_\Lambda) \to \D(S,\Lambda)$ be the right adjoint of $\alpha_S$. We first show that $\beta_S$ is conservative. By construction it is compatible with pushforwards along $p_S\colon S \to *$, so we reduce to showing that $p_{S*}\colon \Shv(S,\Mod_\Lambda) \to \Shv(*,\Mod_\Lambda) = \Mod_\Lambda$ is conservative. In other words we need to see that if a given sheaf $\mathcal M \in \Shv(S,\Mod_\Lambda)$ satisfies $\mathcal M(S) = 0$ then $\mathcal M = 0$. But if $U \subseteq S$ is any clopen subset with complement $V$ then $0 = \mathcal M(S) = \mathcal M(U) \dsum \mathcal M(V)$, forcing $\mathcal M(U) = 0$. For a general open subset $U \subseteq S$ we can compute $\mathcal M(U)$ as a cofiltered limit of the values of $\mathcal M$ on clopen subsets of $U$, forcing $\mathcal M(U) = 0$ in general.

To show that $\alpha_S$ is an equivalence, it is now enough to show that it is fully faithful, i.e.\ that for every $M \in \D(S,\Lambda)$ the unit $M \isoto \beta_S \alpha_S M$ is an isomorphism. We first show that $\beta_S$ preserves colimits. This can be checked after applying the forgetful functor, which reduces the claim to showing that $p_{S*}$ preserves colimits, i.e.\ that the functor $\Shv(S,\Mod_\Lambda) \to \Mod_\Lambda$, $\mathcal M \mapsto \mathcal M(S)$ preserves colimits. Since clopen subsets of $S$ form a basis for the topology on $S$ we can identify $\Shv(S,\Mod_\Lambda)$ with the category of sheaves on the site of clopen subsets of $S$ (rather than all open subsets), cf. \cite[Proposition A.3.11]{Mann.2022a}. But here the sheaf condition simply amounts to $\mathcal M(W) = \mathcal M(U) \times \mathcal M(V)$ for any two disjoint clopen subsets $U, V \subseteq S$ with $W = U \dunion V$. Clearly this condition is stable under colimits of presheaves, proving that $\beta_S$ indeed preserves filtered colimits.

Since $\D(S,\Lambda) = \Mod_{\Lambda(S)}$ is generated under colimits by $\Lambda(S)$, the general isomorphism $M \isoto \beta_S \alpha_S M$ reduces to the case $M = \Lambda(S)$. By compatibility of $\alpha$ with pullbacks and $\beta_S$ with pushforwards, we deduce that the underlying $\Lambda$-module of $\beta_S \alpha_S \Lambda(S)$ is $p_{S*} p_S^{-1} \Lambda$. We are thus reduced to showing that the natural map $\Lambda(S) \isoto (p_S^{-1}\Lambda)(S)$ is an isomorphism of $\Lambda$-modules. This map is obtained from the canonical maps $\Lambda(S_i) = (p_{S_i}^{-1}\Lambda)(S_i) \to (p_S^{-1}\Lambda)(S)$ via passing to the colimit. Now define the sheaf $\Lambda_S \in \Shv(S,\Mod_\Lambda)$ by $\Lambda_S(U) = \Lambda(U)$ for every clopen subset $U \subseteq S$ (note that $U$ is again profinite; one checks easily that this assignment defines a sheaf on $S$). By adjunctions we obtain a natural map $p_S^{-1}\Lambda \to \Lambda_S$ such that the composition $\Lambda(S) \to (p_S^{-1}\Lambda)(S) \to \Lambda_S(S)$ is the identity (check the latter condition for the maps from each $\Lambda(S_i)$ separately). It is therefore enough to show that the map $p_S^{-1}\Lambda \isoto \Lambda_S$ is an isomorphism of sheaves on $S$.

For every point $x \in S$ we denote by $i_x\colon * \to S$ the induced map. We claim that the family of functors $i_x^{-1}\colon \Shv(S,\Mod_\Lambda) \to \Mod_\Lambda$, parametrized over all $x \in S$, is conservative. To see this, let $\mathcal M \in \Shv(S,\Mod_\Lambda)$ be such that $i_x^{-1} \mathcal M = 0$ for all $x$. Note that $i_x^{-1}$ is computed via a left Kan extension on presheaf categories (as there is no sheafification on $*$), hence $i_x^{-1} \mathcal M = \varinjlim_{x \in U \subseteq S} \mathcal M(U)$, where $U$ ranges over all clopen neighborhoods of $x$ in $S$. Fix a clopen subset $U \subseteq S$ and a map $f\colon \Lambda \to \mathcal M(U)$. By the usual stalk arguments we see that there is an open cover of $U$ on which $f$ becomes zero. W.l.o.g.\ this cover is finite and disjoint; but then $f$ must be zero, showing that $\mathcal M = 0$. Finally, observe that $i_x^{-1} (p_S^{-1} \Lambda) = i_x^{-1} \Lambda_S = \Lambda$ for all $x \in S$, proving the desired isomorphism of sheaves on $S$.
\end{proof}

While \cref{rslt:identify-sheaves-on-profinite-sets-with-D} provides good motivation for our definitions, we will almost nowhere need to use it directly (only the conservativity of stalks, which was part of the proof, will be used in \cref{rslt:cond-ani-stalks-are-conservative} below). We now extend the definition of $\D(X,\Lambda)$ to all condensed anima $X$, for which the following descent result is crucial: 

\begin{lemma} \label{rslt:Lambda-sheaves-have-hyperdescent-on-ProFin}
Let $\Lambda$ be an $\Einfty$-ring. Then:
\begin{lemenum}
\item For each map (resp.\ surjection) $T\to S$ of profinite sets, the induced map of $\Einfty$-rings $\Lambda(S)\to \Lambda(T)$ is flat (resp.\ faithfully flat).
\item $\D(\blank, \Lambda)\colon \ProFin^\op \to \Cat$ is a hypercomplete sheaf.
\end{lemenum}
\end{lemma}
\begin{proof}
We first prove (i). For the notion of faithful flatness we refer to \cite[Definition~B.6.1.1]{SAG}. Let $f\colon T\to S$ be a map of profinite sets. We can write $f$ as a cofiltered system of maps $(T_i \to S_i)_i$ between finite sets $T_i$ and $S_i$ (this is true in any $\Pro$-category by the dual version of \cite[Proposition~5.3.5.15]{HTT} and is easy to check for profinite sets by elementary means); we can easily ensure that all the maps $T_i \to S_i$ are surjective if $f$ is surjective. We deduce that the map $\Lambda(S) \to \Lambda(T)$ is the filtered colimit of the maps $\Lambda(S_i) \to \Lambda(T_i)$, so it is enough to show that each of these maps is (faithfully) flat. This follows immediately from writing $\Lambda(S_i) = \prod_{x\in S_i} \Lambda$ and $\Lambda(T_i) = \prod_{x \in T_i} \Lambda$ and the fact that flatness is local (in the Zariski topology) on the source and target.

We now prove (ii). Having established that covers in $\ProFin$ induce faithfully flat maps of rings, the desired descent of categories follows from faithfully flat hyperdescent (see \cite[Corollary~D.6.3.3]{SAG}). In the case that $\Lambda$ is an ordinary ring, one can reduce the claimed hyperdescent to classical faithfully flat descent via standard techniques (using \cite[Proposition~A.1.2]{Mann.2022a}).
\end{proof}

\begin{definition}
We denote
\begin{align}
	\D(\blank, \Lambda)\colon \Cond(\Ani)^\op \to \Cat, \qquad X \mapsto \D(X, \Lambda)
\end{align}
the hypercomplete sheaf of categories associated with \cref{rslt:Lambda-sheaves-have-hyperdescent-on-ProFin} via \cref{rslt:sheaves-on-hypertopos-equiv-hypersheaves-on-site}. The objects of $\D(X, \Lambda)$ are called the \emph{$\Lambda$-valued hypersheaves} on the condensed anima $X$.
\end{definition}

\begin{lemma}\label{rslt:Lambda-sheaves-t-structure}
Let $\Lambda$ be a connective $\Einfty$-ring.
The categories $\D(X,\Lambda)$ are presentable, stable and come equipped with a left-complete t-structure. Moreover, for any morphism $f\colon Y\to X$ of condensed anima, the pullback functor $f^*\colon \D(X,\Lambda) \to \D(Y,\Lambda)$ is t-exact and (co)connectivity can be checked after pullback along a cover.
\end{lemma}
\begin{proof}
Since the subcategory of presentable stable categories is closed under limits in $\Cat$, it is clear that the categories $\D(X,\Lambda)$ are presentable and stable. Fix a condensed anima $X$. To construct the t-structure on $\D(X, \Lambda)$, we define the full subcategory $\D(X, \Lambda)^{\le 0} \subseteq \D(X, \Lambda)$ of connective objects by saying that $M \in \D(X, \Lambda)$ is connective if its pullback to every profinite set $S$ is connective (for the standard t-structure on $\Mod_{\Lambda(S)}$). We have to show that this defines a left complete t-structure on $\D(X, \Lambda)$ whose coconnective objects are those that become coconnective after pullback to every profinite set. Then the last part of the claim is automatic, as (co)connectivity can be checked after applying a t-exact conservative functor.

The claim is clear if $X = S$ is profinite and recovers the standard t-structure on $\D(S,\Lambda) = \Mod_{\Lambda(S)}$, because pullbacks along maps of profinite sets are t-exact. For general $X$, write $X = \varinjlim_i S_i$ in $\Cond(\Ani)$, where all $S_i$ are profinite sets and w.l.o.g.\ all profinite sets appear in the colimit. For brevity we write $\cat C = \D(X, \Lambda)$ and $\cat C_i = \D(S_i,\Lambda)$, so that $\cat C= \varprojlim_i \cat C_i$. Then $\cat C^{\le 0}$ consists of exactly those objects whose image in each $\cat C_i$ lies in $\cat C_i^{\le 0}$, so that $\cat C^{\le 0} = \varprojlim_i \cat C_i^{\le 0}$ (see \cref{rslt:limit-of-subcategories}). Using \cref{rslt:termwise-adjoint-of-limit} the truncation functors $\tau^{\le0}\colon \cat C_i \to \cat C_i^{\le0}$ assemble into a right adjoint $\tau^{\le0}\colon \cat C \to \cat C^{\le0}$. As $\cat C^{\le 0}$ is closed under extensions, we deduce from (the dual of) \cite[Proposition~1.2.1.16]{HA} that $(\cat C^{\le 0}, \cat C)$ defines a t-structure on $\cat C$. One easily checks that
\begin{align}
	\cat C^{\ge1} \coloneqq \set{X\in \cat C}{\tau^{\le0}X = 0} = \lim_i \cat C_i^{\ge1}.
\end{align}
This implies that all pullbacks $\cat C \to \cat C_i$ are t-exact. It also implies that $\cat C$ is left-complete, since all $\cat C_i$ are left-complete and limits commute with limits.
\end{proof}

\begin{remark} \label{rslt:identify-6ff-with-sheaves-on-top-spaces}
Suppose $\Lambda$ is an ordinary ring and $X$ is a locally compact Hausdorff space (viewed as a condensed anima). Then $\D(X, \Lambda)$ is the left-completion of the derived category of the classical category of sheaves of $\Lambda$-modules on the topological space $X$. Indeed, using \cref{rslt:identify-sheaves-on-profinite-sets-with-D} one can construct a natural functor between both categories and then the statement reduces to a descent result for sheaves on locally compact Hausdorff spaces shown in \cite[Propositions~7.1,~7.3]{Scholze:Six-Functor-Formalism}. In particular, a left-bounded object in $\D(X, \Lambda)$ is the same as a left-bounded complex of sheaves of $\Lambda$-modules on $X$. If $X$ is paracompact and has finite covering dimension, then no left-completion is necessary (see \cite[Proposition~7.1]{Scholze:Six-Functor-Formalism}). We will only use these remarks as a motivation for our definitions, but we nowhere need to rely on these results directly.
\end{remark}

With the definition of $\D(X, \Lambda)$ at hand, we can finally introduce the 6-functor formalism on this category by combining \cref{rslt:construct-3ff-from-suitable-decomp,rslt:extend-6ff-to-stacks-and-stacky-maps}.

\begin{construction} \label{cstr:6ff-on-Cond-Ani}
Let $\Lambda$ be an $\Einfty$-ring. We let $E$ and $P$ denote the classes of all edges in $\ProFin_\kappa$ and we let $I$ be the class of isomorphisms in $\ProFin_\kappa$. Then $I, P$ is a suitable decomposition of $E$. Using the tensor products of modules, $\D(\blank, \Lambda)$ upgrades to a functor $\ProFin_\kappa^\op \to \CMon$ and one easily verifies the conditions of \cref{rslt:construct-3ff-from-suitable-decomp}; indeed, conditions (a) and (c) are essentially vacuous and condition (b) follows from standard properties of change-of-ring and forgetful functors of modules. By \cref{rslt:construct-3ff-from-suitable-decomp} we can thus upgrade $\D(\blank, \Lambda)$ to a 3-functor formalism, i.e.\ a lax symmetric monoidal functor
\begin{align}
	\D(\blank, \Lambda)\colon \Corr(\ProFin_\kappa) \to \Cat.
\end{align}
The site $\ProFin_\kappa$ is hyper-subcanonical (by \cref{rmk:subcanonical-ordinary-site-is-hyper-subcanonical}) and $\D(\blank, \Lambda)$ satisfies the conditions of \cref{rslt:extend-6ff-to-stacks-and-stacky-maps}. Therefore, this 6-functor formalism extends uniquely to a 6-functor formalism
\begin{align}
	\D(\blank, \Lambda)\colon \Corr(\Cond(\Ani), E') \to \Cat,
\end{align}
for a certain class of edges $E'$ satisfying the properties in \cref{rslt:extend-6ff-to-stacks-and-stacky-maps}.
\end{construction}

\begin{definition} \label{def:6ff-on-Cond-Ani}
Let $\Lambda$ be an $\Einfty$-ring. We say that a map $f\colon Y \to X$ in $\Cond(\Ani)$ is \emph{$\Lambda$-fine} if it belongs to the smallest class of edges $E'$ satisfying the properties of \cref{rslt:extend-6ff-to-stacks-and-stacky-maps} for the 6-functor formalism in \cref{cstr:6ff-on-Cond-Ani}. We denote by
\begin{align}
	\D(\blank, \Lambda)\colon \Corr(\Cond(\Ani), \fine{\Lambda}) \to \Cat
\end{align}
the 6-functor formalism constructed in \cref{cstr:6ff-on-Cond-Ani}.
\end{definition}

By construction, every map $f\colon T \to S$ of profinite sets is $\Lambda$-fine. The functor $f_*$ is the forgetful functor $\Mod_{\Lambda(T)} \to \Mod_{\Lambda(S)}$, the functor $f_!$ is isomorphic to $f_*$, the functor $f^*$ is given by the change-of-ring map $\blank \tensor_{\Lambda(S)} \Lambda(T)$ and the functor $f^!$ is given by $\iHom_{\Lambda(S)}(\Lambda(T), \blank)$ with the induced $\Lambda(T)$-module structure (for the computation of $f^!$ we use \cref{rslt:enriched-adjunction-of-shriek-functors} with $N = \Lambda(T)$). It is convenient to have the following generalization of $\Lambda(\blank)$:

\begin{definition}
Let $\Lambda$ be an $\Einfty$-ring and $X$ a condensed anima with projection $p_X\colon X \to *$.
\begin{defenum}
	\item We denote
	\begin{align}
		\Lambda(X) \coloneqq \Gamma(X, \Lambda) \coloneqq p_{X*} \one \in \Mod_\Lambda
	\end{align}
	and call it the \emph{$\Lambda$-valued cohomology of $X$}. We similarly denote $\Gamma(X,M) \coloneqq p_{X*} M$ for any $M \in \D(X,\Lambda)$.

	\item \label{def:Gamma-c} We say that $X$ is \emph{$\Lambda$-fine} if $p_X$ is so. In this case we denote
	\begin{align}
		\Gamma_c(X,\Lambda) \coloneqq p_{X!} \one \in \Mod_\Lambda.
	\end{align}
\end{defenum}
\end{definition}

\begin{lemma} \label{rslt:cohomology-is-sheafy}
The assignment $X \mapsto \Gamma(X,\Lambda)$ defines a hypercomplete sheaf $\Cond(\Ani)^\op \to \Mod_\Lambda$ and agrees with \cref{def:Einfty-ring-associated-to-profin-set} for profinite sets $X$.
\end{lemma}
\begin{proof}
We first construct the functor $\Gamma(\blank,\Lambda)$, which is a standard exercise with straightening. First note that there is a natural transformation $p_{\bullet}^*\colon \Mod_\Lambda \to \D(\blank,\Lambda)$ of functors $\Cond(\Ani)^\op \to \Cat$ given by $p_X^*\colon \Mod_\Lambda \to \D(X,\Lambda)$ for each condensed anima $X$. Namely, such a natural transformation is equivalent to a functor $\Mod_\Lambda \to \varprojlim_{X \in \Cond(\Ani)^\op} \D(X,\Lambda)$; but since $\Cond(\Ani)^\op$ has the initial object $*$, the limit on the right is simply given by $\D(*,\Lambda) = \Mod_\Lambda$. By letting $\cat E \to \Cond(\Ani)^\op$ denote the cocartesian fibration associated with $\D(\blank,\Lambda)$, we obtain a map of cocartesian fibrations $p^*\colon \Cond(\Ani)^\op \times \Mod_\Lambda \to \cat E$. By \cite[Proposition~7.3.2.6]{HA} the functor $p^*$ admits a right adjoint $p_*$ over $\Cond(\Ani)^\op$, which is fiberwise given by $p_{X*}$. We thus obtain the endofunctor $p_* p^*$ of $\Cond(\Ani)^\op \times \Mod_\Lambda$ over $\Cond(\Ani)^\op$. We can equivalently view it as the functor $\Cond(\Ani)^\op \times \Mod_\Lambda \to \Mod_\Lambda$ sending $(X, M) \mapsto p_{X*} p_X^* M$. Fixing $M = \Lambda$ provides the desired functor $\Gamma(\blank,\Lambda)$.

We now prove the sheafiness of $\Gamma(\blank,\Lambda)$, i.e.\ for a colimit diagram $X = \varinjlim_i X_i$ in $\Cond(\Ani)$ we need to see that the natural map
\begin{align}
	\Gamma(X,\Lambda) \isoto \varprojlim_i \Gamma(X_i,\Lambda)
\end{align}
is an isomorphism. Let us denote by $f_i\colon X_i \to X$ the natural map. Then it is enough to show that the natural map $\one \isoto \varprojlim_i f_{i*} \one$ is an isomorphism in $\D(X,\Lambda)$ (then the claimed isomorphism follows by applying $p_{X*}$). But this follows immediately from $\D(X,\Lambda) = \varprojlim_i \D(X_i,\Lambda)$ using \cref{rslt:limits-and-adjoints-in-Cat}.

To prove the last assertion, it remains to show that $\Gamma(S,\Lambda) = \Lambda(S)$ agrees with \cref{def:Einfty-ring-associated-to-profin-set} for profinite sets $S$. But this is clear from the definition $\D(S,\Lambda) = \Mod_{\Lambda(S)}$.
\end{proof}

To make a 6-functor formalism useful, one needs to find a large enough class of maps for which the $!$-functors are defined. In the case of $\D(\blank,\Lambda)$ this amounts to proving $\Lambda$-fineness for as many maps as possible. We start with a few basic examples:

\begin{lemma} \label{rslt:Lambda-fineness-and-disjoint-unions}
Let $\Lambda$ be an $\Einfty$-ring.
\begin{lemenum}
	\item For all condensed anima $X$ and $Y$ the inclusion $X \injto X \dunion Y$ is $\Lambda$-fine.

	\item If $(Y_i \to X)_i$ is a family of $\Lambda$-fine maps, then $\bigdunion_i Y_i \to X$ is a $\Lambda$-fine map.
\end{lemenum}
\end{lemma}
\begin{proof}
By \cref{rslt:extend-6ff-to-disjoint-unions} this reduces to showing that the map $* \to * \dunion *$ is $\Lambda$-fine, where $*$ is the set consisting of one object. But this is clear, since both $*$ and $* \dunion *$ are profinite sets and any map between profinite sets is $\Lambda$-fine by construction.
\end{proof}

From \cref{rslt:Lambda-fineness-and-disjoint-unions} we deduce that every map between disjoint unions of profinite sets is $\Lambda$-fine, in particular every homomorphism between locally profinite groups. Of course there are many more $\Lambda$-fine maps than these, but we have not yet developed the machinery to prove universal $!$-descent in a wide array of examples. We will therefore postpone the construction of more $\Lambda$-fine maps to \cref{sec:kerncat.examples} (for maps between topological spaces and between anima) and \cref{sec:reptheory} (for maps between stacks of locally profinite groups).

It is sometimes convenient to understand the behavior of the 6-functor formalism $\D(\blank,\Lambda)$ under change of $\Lambda$. This is solved by the following result, which tells us that $\D(\blank,\Lambda)$ is functorial in $\Lambda$.

\begin{definition}
A map $(f, \varphi)\colon (Y, \Lambda) \to (X, \Lambda')$ in $\Cond(\Ani) \times \CAlg^\op$ is called \emph{fine} if $\varphi\colon \Lambda \isoto \Lambda'$ is an isomorphism and $f\colon Y \to X$ is $\Lambda$-fine.
\end{definition}

\begin{proposition} \label{rslt:6ff-on-Cond-Ani-functorial-in-Lambda}
There is a 6-functor formalism
\begin{align}
	\D(\blank, \blank)\colon \Corr(\Cond(\Ani) \times \CAlg^\op, \ufine) \to \Cat
\end{align}
with the following properties:
\begin{propenum}
	\item For a fixed $\Einfty$-ring $\Lambda$, $\D(\blank,\Lambda)$ coincides with the 6-functor formalism from \cref{def:6ff-on-Cond-Ani}.

	\item Let $X \in \Cond(\Ani)$ and let $\varphi\colon \Lambda \to \Lambda'$ be a map of $\Einfty$-rings. Then there is a natural equivalence $\D(X,\Lambda') = \LMod_{\Lambda'}(\D(X, \Lambda))$ and we have
	\begin{align}
		\varphi^* \coloneqq (\id_X, \varphi)^* = \blank \tensor_\Lambda \Lambda'&\colon \D(X, \Lambda) \to \D(X, \Lambda')\\
		\varphi_* \coloneqq (\id_X, \varphi)_* = \mathrm{forget}&\colon \D(X, \Lambda') \to \D(X, \Lambda).
	\end{align}

	\item \label{rslt:compatibility-of-change-of-Lambda-with-6ff} In the setting of (ii), let $f\colon Y \to X$ be a $\Lambda$-fine map. Then $f$ is also $\Lambda'$-fine. Moreover, $\varphi^*$ commutes with $f^*$, $f_!$ and $f^!$, and $\varphi_*$ commutes with $f^*$, $f_!$ and $f^!$, via the natural maps.
\end{propenum}
\end{proposition}
\begin{proof}
We start with the functor $(\ProFin_\kappa)^{\op} \times \CAlg \to \CAlg$ given by $(S, \Lambda) \mapsto \Lambda(S)$. Composing this functor with \cref{rslt:functoriality-of-modules-over-algebras} defines the functor
\begin{align}
	\D(\blank,\blank)\colon (\ProFin_\kappa)^\op \times \CAlg \to \CMon, \qquad (S, \Lambda) \mapsto \D(S, \Lambda).
\end{align}
We now apply \cref{rslt:construct-3ff-from-suitable-decomp} to this functor, where we take $\cat C = \ProFin_\kappa \times \CAlg^\op$, let $E$ be the collection of maps $(f, \varphi)$ such that $\varphi$ is an isomorphism, let $P = E$ and let $I$ be the collection of isomorphisms. Using (the proof of) \cref{rslt:extend-3ff-to-stacks} and \cref{rmk:extend-6ff-to-stacks} we can extend $\D(\blank,\blank)$ to a 6-functor formalism
\begin{align}
	\D(\blank,\blank)\colon \Corr(\Cond(\Ani) \times \CAlg^\op, E) \to \Cat,
\end{align}
where $E$ is the collection of maps $(f, \varphi)$ such that $\varphi$ is an isomorphism and $f$ is representable in profinite sets. By construction, for fixed $\Lambda$ this 6-functor formalism coincides with the one from \cref{def:6ff-on-Cond-Ani}. We now prove that it satisfies (ii). 
Fix $X \in \Cond(\Ani)$ and a map $\varphi\colon \Lambda \to \Lambda'$ of $\Einfty$-rings. We first note that for any map $f\colon Y \to X$ in $\Cond(\Ani)$ the natural map $f^* \varphi_* \isoto \varphi_* f^*$ is an isomorphism of functors $\D(X, \Lambda') \to \D(Y, \Lambda)$. Indeed, by standard reductions this reduces to the case that $X$ and $Y$ are profinite sets, where it is clear. We deduce that $\varphi_*$ is conservative, preserves all small colimits and is linear over $\D(X,\Lambda)$ (by reduction to the case that $X$ is a profinite set, for linearity note that lax linearity is guaranteed by \cref{rslt:monoidal-right-adjoint}). Thus (ii) follows from Lurie's Barr--Beck theorem \cite[Theorem~4.7.3.5]{HA}---the only subtlety is to identify modules under the monad $\varphi_* \varphi^*$ with modules under $\Lambda'$, i.e.\ to identify the monad $\varphi_* \varphi^*$ with the monad $\blank \tensor_\Lambda \Lambda'$. To this end we can use \cite[Theorem~1.1]{Heine:Monadicity} to show that $\varphi_* \varphi^*$ is a monad in $\Mod_{\D(X,\Lambda)}(\PrL)$, i.e.\ an algebra in the endomorphism category of $\D(X,\Lambda)$ in $\Mod_{\D(X,\Lambda)}(\PrL)$; but this endomorphism category identifies with $\D(X,\Lambda)$ and hence the monad must be induced by an algebra in $\D(X,\Lambda)$.

Next we show that the above 6-functor formalism satisfies (iii) for maps $f \in E$. The commutation of $\varphi_*$ and $f^*$ was shown above. The commutation of $\varphi^*$ and $f_!$ is part of the 6-functor formalism, and we can also deduce the commutation of $\varphi_*$ and $f^!$ by passing to right adjoints. The commutation of $\varphi^*$ and $f^*$ is also encoded in the 6-functor formalism. Thus it only remains to show that $\varphi^*$ and $f^!$ commute via the natural map. But this reduces to the case that $X$ (and hence $Y$) is a profinite set, where it is clear.

We now extend the above 6-functor formalism from $E$ to $\ufine$ by proceeding as in the proof of \cref{rslt:extend-6ff-to-stacks-and-stacky-maps}. In each step we additionally require that property (iii) above is satisfied. Then everything is formal, but to show that we can extend to all of $\ufine$ we need to prove the following claim: If $\varphi\colon \Lambda \to \Lambda'$ is a map of $\Einfty$-rings and $(f_i\colon U_i \to X)_{i\in I}$ is a universal $\D^!(\blank,\Lambda)$-cover, then it is also a universal $\D^!(\blank,\Lambda')$-cover. To see this we note that by \cref{rslt:limits-and-adjoints-in-Cat} (and \cref{rslt:descent-data-for-generated-sieve-functor-on-Delta-I}) there is an adjunction
\begin{align}
	f_{i_\bullet!} \colon \varprojlim_{i_\bullet \in \bbDelta_I} \D^!(U_{i_\bullet}, \Lambda') &\rightleftarrows \D^!(X, \Lambda') \noloc f^!_{i_\bullet}\\
	(M_{i_\bullet})_{i_\bullet} &\mapsto \varinjlim_{i_\bullet \in \bbDelta_I^{\op}} f_{i_\bullet!} M_{i_\bullet}\\
	(f_{i_\bullet}^! N) & \mapsfrom N
\end{align}
and we need to show that this adjunction is an equivalence, i.e.\ that the unit and counit of the adjunction are isomorphisms. This can be checked after applying the conservative functor $\varphi_*$. But by (iii), $\varphi_*$ commutes with all the functors appearing above, hence the claim reduces easily to the case with $\Lambda$ in place of $\Lambda'$, where it is true by assumption.
\end{proof}

\begin{example} \label{rslt:functoriality-of-6ff-on-Cond-Ani-in-Lambda}
Let $\Lambda$ be an $\Einfty$-ring. By abuse of notation we denote by $\fine{\Lambda}$ the collection of morphisms $(f,\varphi)$ in $\Cond(\Ani) \times (\CAlg_{\Lambda/})^\op$ where $\varphi$ is an isomorphism and $f$ is $\Lambda$-fine. By \cref{rslt:compatibility-of-change-of-Lambda-with-6ff} we obtain a map of geometric setups
\begin{align}
	(\Cond(\Ani) \times (\CAlg_{\Lambda/})^\op, \fine{\Lambda}) \to (\Cond(\Ani) \times \CAlg^\op, \ufine).
\end{align}
Thus the 6-functor formalism from \cref{rslt:6ff-on-Cond-Ani-functorial-in-Lambda} together with \cref{rslt:Corr-preserves-limits} produces a lax monoidal functor
\begin{align}
	&\CAlg_{\Lambda/} \times \Corr(\Cond(\Ani), \fine{\Lambda}) = \Corr(\Cond(\Ani) \times (\CAlg_{\Lambda/}^\op), \fine{\Lambda}) \\&\qquad\to \Corr(\Cond(\Ani) \times \CAlg^\op, \ufine) \xto{\D(\blank,\blank)} \Cat.
\end{align}
In particular we obtain a functor
\begin{align}
	\CAlg_{\Lambda/} \to \Alg(\Corr(\Cond(\Ani), \fine{\Lambda}), \Cat), \qquad \Lambda' \mapsto \D(\blank,\Lambda').
\end{align}
which exhibits the 6-functor formalism $\D(\blank, \Lambda')$ as functorial in the $\Lambda$-algebra $\Lambda'$.
\end{example}

This is all we want to say about the 6-functor formalism of $\Lambda$-valued sheaves on condensed anima for now. In the following we provide examples of other 6-functor formalisms (with references to the literature) that can be constructed in a similar way.

\begin{examples}
\begin{exampleenum}
	\item (Étale sheaves on schemes.) Let $\Lambda$ be an ordinary torsion ring. By \cite[\S7]{Scholze:Six-Functor-Formalism} there is a 6-functor formalism $\D_\et(\blank, \Lambda)$ for the derived category of étale $\Lambda$-sheaves on qcqs schemes. Here the exceptional morphisms $E$ are the separated maps of \enquote{finite expansion} (this occurs as \enquote{$+$-finite type} in \cite[Definition~2.9.28]{Mann.2022a}), including all separated maps that are of finite type or integral. For a universally closed map $f\colon Y \to X$ in $E$, we have $f_! \isom f_*$ and for an open immersion $j\colon U \injto X$ we have $j^! \isom j^*$. In particular, this 6-functor formalism recovers the usual notion of étale cohomology with compact support. Note that the 6-functor formalism exists without any requirement on the torsion of $\Lambda$ being invertible on the schemes; this condition appears only when one talks about smoothness. Using \cref{rslt:extend-6ff-to-stacks-and-stacky-maps} one can extend the 6-functor formalism to stacks; this was the main goal of \cite{Liu-Zheng.2012}, whose methods we used in \cref{rslt:construct-3ff-from-suitable-decomp}.

	\item ($\ell$-adic sheaves on rigid varieties.) Fix primes $\ell \ne p$ and let $\Lambda$ be a nuclear $\ZZ_\ell$-algebra (i.e.\ a condensed animated ring which can be written as a filtered colimit of $\ZZ_\ell$-Banach algebras), for example $\Lambda \in \{ \ZZ_\ell, \QQ_\ell, \overline{\QQ_\ell}, \CC_\ell \}$ or $\Lambda$ is a $\ZZ/\ell^n$-algebra for some $n > 0$. In \cite{Mann.2022b} the second author defines a category $\D_\nuc(X, \Lambda)$ of \emph{nuclear $\Lambda$-sheaves} associated with every small v-stack on perfectoid spaces in characteristic $p$. This category admits a full 6-functor formalism and in the case that $\Lambda$ is discrete it coincides with the 6-functor formalism for étale $\Lambda$-modules as constructed in \cite{Scholze.2018}, except that with the techniques of the present paper we can enlarge the class $E$ of exceptional morphisms in order to include many stacky maps. Via \cite[\S15]{Scholze.2018} this in particular defines a 6-functor formalism for pro-étale $\Lambda$-cohomology on rigid-analytic varieties over non-archimedean fields over $\ZZ_p$. The first attempt at this extension of $E$ appeared in \cite{Gulotta-Hansen-Weinstein:Enhanced-6FF-on-vStacks}, whose content was a great inspiration to the second author.

	\item ($p$-adic sheaves on rigid varieties.) Fix a prime $p$ and an $\FF_p$-algebra $\Lambda$. In \cite{Mann.2022a} the second author constructs a category $\D(X, \Lambda)$ of \enquote{$\Lambda$-modules} on small v-stacks over $\ZZ_p$ and equips it with a full 6-functor formalism. The definition of $\D(X, \Lambda)$ is quite involved, but it contains the category of étale $\Lambda$-sheaves as a full subcategory, which allows one to deduce statements about classical étale $\Lambda$-cohomology on small v-stacks (and in particular on rigid-analytic varieties). In the process of constructing this 6-functor formalism, the second author developed the tools of constructing and extending abstract 6-functor formalisms, which now appear in a refined form in the previous subsections. In \cite{Hansen-Mann:Mod-p-Stacky-6FF} this $p$-adic 6-functor formalism is refined to include stacky maps and then studied on classifying stacks of locally profinite groups, in a similar fashion as in \cref{sec:reptheory} below.

	\item (Quasi-coherent sheaves.) One of the main results of \cite{Clausen-Scholze:Condensed-Mathematics} is the construction of a 6-functor formalism $\QCoh_\solid(\blank)$ for solid quasicoherent sheaves on schemes. Here the term \enquote{solid} means that the sheaves are equipped with a topological structure, which allows to define a left adjoint $j_!\colon \QCoh_\solid(U) \to \QCoh_\solid(X)$ to the pullback $j^*$ along an open immersion $j\colon U \injto X$ of schemes. In \cite[Proposition~2.9.31]{Mann.2022a} the construction of this 6-functor formalism is formalized and extended to qcqs discrete adic spaces. The class $E$ of exceptional morphisms consists of the separated maps of $+$-finite type (as in the above example for étale sheaves).
\end{exampleenum}
\end{examples}

With the tools provided in this paper, it should be straightforward to formalize many more 6-functor formalisms that appear in the literature and hence extend them to stacks and stacky maps. As we hope to show in this paper, this extension can be very fruitful and sometimes provide surprising new insights into existing constructions.

\section{The category of kernels} \label{sec:kerncat}

This section is devoted to constructing and studying the category of kernels associated with a 6-functor formalism. This category is a certain 2-category which encodes a surprising amount of useful geometric information in terms of simple 2-categorical data. This idea is due to Lu--Zheng and first appeared in \cite{Lu-Zheng:ULA} to study universally locally acyclic sheaves. It was then adapted by Fargues--Scholze in their geometrization of the local Langlands correspondence \cite{Fargues-Scholze:Geometrization} where they define a different 2-category but use it in an analogous way. In \cite{Mann.2022b,Scholze:Six-Functor-Formalism} some more results on this 2-category were produced and it became clear that this 2-category is a very powerful tool to study 6-functor formalisms in general. In the following we put all of these results on firm $\infty$-categorical grounds as well as providing some further criteria and compatibilities.

This section is structured as follows. In \cref{sec:kerncat.def} we start by defining the 2-category of kernels $\cat K_{\D,S}$ and study its basic properties. In \cref{sec:kerncat.functorial} we then show that $\cat K_{\D,S}$ is functorial in $\D$ and $S$ in an appropriate sense and in \cref{sec:kerncat.descent} we investigate how to interpret $*$-covers in terms of $\cat K_{\D,S}$. In \cref{sec:kerncat.suave-prim-obj} we then introduce suave and prim sheaves in a 6-functor formalism using adjunctions in $\cat K_\D$. This leads to the definition of smooth and proper maps in \cref{sec:kerncat.etale-proper} and various compatibilities between all of the notions introduced so far. In \cref{sec:kerncat.excdescent} we discuss some results for proving $!$-descent in the presence of smooth and proper maps. Finally, in \cref{sec:kerncat.examples} we then apply these results to the 6-functor formalisms from \cref{sec:6ff.example}.

The results in this section rely heavily on the material in \cref{sec:enr,sec:2cat} and hence it may be advisable for the reader to consult these sections beforehand.

\subsection{Basic definitions} \label{sec:kerncat.def}

We start by defining the 2-category of kernels associated with a 3-functor formalism $\D$. For its construction it is useful to single out a special class of geometric setups:

\begin{definition}
A geometric setup $(\cat C, E)$ is called \emph{full} if $E$ contains all edges in $\cat C$ and $\cat C$ admits finite products.
\end{definition}

\begin{remark}
If $(\cat C, E)$ is a geometric setup and $S \in \cat C$ is an object, then the geometric setup $((\cat C_E)_{/S}, \mathit{all})$ is full.
\end{remark}

We can now come to the definition of the category of kernels. We start with the formal definition and afterwards provide a more intuitive description.

\begin{definition} \label{def:kern-cat}
Let $\D$ be a 3-functor formalism on a geometric setup $(\cat C, E)$.
\begin{defenum}
	\item Suppose $(\cat C, E)$ is full. Then we define the \emph{2-category of kernels} associated with $\D$ as
	\begin{align}
		\cat K_\D \coloneqq \tau_\D(\Corr(\cat C)) \in \TwoCat.
	\end{align}
	Here $\tau_\D\colon \Enr_{\Corr(\cat C)} \to \Enr_{\Cat} = \TwoCat$ denotes the transfer of enrichment from \cref{def:transfer-of-enrichment} and we equip $\Corr(\cat C)$ with the natural self-enrichment from \cref{rslt:Corr-self-enriched-and-self-dual}.

	\item \label{def:relative-category-of-kernels} For every $S \in \cat C$ the natural functor $(\cat C_E)_{/S} \to \cat C$ induces a morphism of geometric setups $((\cat C_E)_{/S}, \mathit{all}) \to (\cat C, E)$ and hence a lax symmetric monoidal functor $\Corr((\cat C_E)_{/S}) \to \Corr(\cat C, E)$ by \cref{rslt:functoriality-of-Corr^tensor}. We denote by $\D_S$ the 3-functor formalism on $((\cat C_E)_{/S}, \mathit{all})$ induced by restricting along this functor and we call
	\begin{align}
		\cat K_{\D, S} \coloneqq \cat K_{\D_S} \in \TwoCat
	\end{align}
	the \emph{2-category of kernels} associated with $\D$ over $S$. We will often abbreviate
	\begin{align}
		\Fun_S(\blank, \blank) = \Fun_{\cat K_{\D,S}}(\blank, \blank)
	\end{align}
	to denote the category of homomorphisms in $\cat K_{\D,S}$, if $\D$ is clear from context.
\end{defenum}
\end{definition}

Fix a 3-functor formalism $\D$ on a geometric setup $(\cat C, E)$ and fix some $S \in \cat C$. Let us explain in more down-to-earth terms what the 2-category $\cat K_{\D,S}$ looks like:
\begin{itemize}
	\item The objects of $\cat K_{\D,S}$ are the morphisms $X \to S$ in $E$.

	\item Given two objects $X, Y \in \cat K_{\D,S}$, the category of homomorphisms from $Y$ to $X$ is given by
	\begin{align}
		\Fun_S(Y, X) = \D(X \times_S Y).
	\end{align}
	Indeed, by \cref{rslt:Corr-closed-symmetric-monoidal} we have $\iHom(Y, X) = X \times_S Y$ in $\Corr((\cat C_E)_{/S})$, so the above formula follows directly from the explicit description of the transfer of enrichment (see \cref{sec:enr.transfer}).

	\item Given objects $X, Y, Z \in \cat K_{\D,S}$ and morphisms $M\colon Y \to X$ and $N\colon Z \to Y$ in $\cat K_{\D,S}$, the composition $M \comp N$ is given by
	\begin{align}
		M \comp N = \pi_{13!}(\pi_{12}^* M \tensor \pi_{23}^* N),
	\end{align}
	where $\pi_{ij}$ denote the projections on $X \times_S Y \times_S Z$. Here we note that by the above description of morphisms in $\cat K_{\D,S}$ we have $M \in \D(X \times_S Y)$, $N \in \D(Y \times_S Z)$ and $M \comp N \in \D(X \times_S Z)$, so the above formula for $M \comp N$ makes sense.
\end{itemize}
From the above description of $\cat K_{\D,S}$ it seems clear that $\cat K_{\D,S}$ is isomorphic to its own opposite, because clearly $\D(Y \times_S X) = \D(X \times_S Y)$. Indeed, we have:

\begin{proposition} \label{rslt:kern-cat-is-self-dual}
Let $\D$ be a 3-functor formalism on a geometric setup $(\cat C, E)$. Then for every $S \in \cat C$ we have a natural isomorphism
\begin{align}
	\cat K_{\D,S}^\op = \cat K_{\D,S}
\end{align}
of 2-categories.
\end{proposition}
\begin{proof}
By \cref{rslt:Corr-self-enriched-and-self-dual} there is a natural equivalence $\Corr((\cat C_E)_{/S})^\op = \Corr((\cat C_E)_{/S})$ in the category of $\Corr((\cat C_E)_{/S})$-enriched categories. The claim therefore follows from the fact that transfer of enrichment commutes with taking opposite categories, see \cref{rslt:transfer-of-enrichment-commutes-with-op}.
\end{proof}

Another useful observation is that the correspondence category admits a natural functor to the category of kernels, as follows:

\begin{proposition} \label{rslt:functors-from-and-to-kercat}
Let $\D$ be a 3-functor formalism on a geometric setup $(\cat C, E)$. Then for every $S \in \cat C$ there are natural 2-functors
\begin{align}
	\Corr((\cat C_E)_{/S}) \xlongto{\Phi_{\D,S}} \cat K_{\D,S} \xlongto{\Psi_{\D,S}} \TwoCat
\end{align}
whose composition is the functor $\D$. They have the following properties:
\begin{propenum}
	\item The functor $\Phi_{\D,S}$ acts as the identity on objects. It sends a morphism $f = [\CorrHom{Y}{}{Z}{}{X}]$ to
	\begin{align}
		\Phi_{\D,S}(f) = f'_!(\one) \in \D(X \times_S Y) = \Fun_S(Y, X),
	\end{align}
	where $f'\colon Z \to X \times_S Y$ is the map induced by $f$.

	\item We have $\Psi_{\D,S} = \Fun_S(S, \blank)$. Concretely, this 2-functor sends $X \in \cat K_{\D,S}$ to $\D(X)$ and a morphism $M \in \Fun_S(Y, X)$ to the functor
	\begin{align}
		\Psi_{\D,S}(M) = \pi_{1!}(M \tensor \pi_2^*(\blank))\colon \D(Y) \to \D(X),
	\end{align}
	where $\pi_i$ denote the two projections on $X \times_S Y$.
\end{propenum}
\end{proposition}
\begin{proof}
The functor $\Phi_{\D,S}$ is the functor $F_\D$ from \cref{rslt:functor-on-underlying-categories-of-transfer-of-enrichment} (this is a 1-functor and in particular a 2-functor, since its source is a 1-category) and $\Psi_{\D,S}$ is the functor $G_\D$ from \cref{rslt:functor-from-transfer-of-self-enrichment}. The explicit descriptions in (i) and (ii) follow immediately from the explicit descriptions of $F_\D$ and $G_\D$ in their respective construction results. The fact that $\Psi_{\D,S} \comp \Phi_{\D,S} = \D$ is shown in \cref{rslt:factor-map-of-monoidal-categories-into-enriched-maps}.
\end{proof}

Associated with every 3-functor formalism $\D$ on a full geometric setup $(\mathcal C, \mathit{all})$, we have constructed the following diagram of (2-)functors, which summarizes the most important abstract constructions we have made so far:
\begin{equation}\begin{tikzcd}[row sep=tiny]
	\mathcal C \arrow[dr]\\
	& \Corr(\mathcal C) \arrow[r,"\Phi_\D"] & \mathcal K_\D \arrow[r,"\Psi_\D"] & \TwoCat\\
	\mathcal C^\op \arrow[ur]
\end{tikzcd}\end{equation}
The two diagonal functors on the left are symmetric monoidal (where we equip $\cat C$ and $\cat C^\op$ with the cartesian and cocartesian symmetric monoidal structures, respectively), while the composition $\D = \Psi_\D \comp \Phi_\D$ is lax symmetric monoidal. In future work we will discuss symmetric monoidal 2-categories and equip $\cat K_\D$ with such a structure; then $\Phi_\D$ is a symmetric monoidal functor and $\Psi_\D$ is a lax symmetric monoidal 2-functor.

\begin{remark}
The functor $\Psi_{\D,S}$ is the reason for our name \enquote{category of \emph{kernels}}. Indeed, functors of the form $\pi_{1!}(M \tensor \pi_2^*(\blank))\colon \D(Y) \to \D(X)$ are common in many geometric fields of mathematics, and the tensor factor $M$ is often called the \enquote{kernel} of this functor. Thus, $\mathcal K_{\D,S}$ is the category whose morphisms from $Y$ to $X$ are exactly the kernels for functors $\D(Y) \to \D(X)$.
\end{remark}

\begin{remark}\label{rmk:functor-from-kercat-to-LMod}
The functor $\Psi_{\D,S}$ factors through a 2-functor
\begin{align}
\cat K_{\D,S} \to \LMod_{\D(S)}(\Cat).
\end{align}
Indeed, we obtain a $\LMod_{\D(S)}(\Cat)$-enriched functor $\Psi_{\wt{\D},S} \colon \cat K_{\wt{\D},S} \to \LMod_{\D(S)}(\Cat)$ by applying \cref{rslt:functor-from-transfer-of-self-enrichment} to the lax symmetric monoidal functor $\alpha = \wt{\D}\colon \Corr((\cat C_E)_{/S})^{\tensor} \to \LMod_{\D(S)}(\Cat)^{\tensor}$ from \cref{rslt:3ff-linear-over-base}. By transferring the enrichment along the forgetful functor $\LMod_{\D(S)}(\Cat) \to \Cat$ we obtain the desired enhancement of $\Psi_{\D,S}$.
\end{remark}

\subsection{Functoriality} \label{sec:kerncat.functorial}

We have introduced the category of kernels $\cat K_{\D,S}$ attached to a fixed 3-functor formalism $\D$ and a fixed object $S \in \cat C$. The next step is to discuss the functoriality of $\cat K_{\D,S}$ when $\D$ and $S$ vary. This will be important when we study certain 2-categorical constructions in $\cat K_{\D,S}$, as these are preserved by 2-functors. Let us start with the functoriality in $\D$:

\begin{proposition}
Let $\D$ be a 3-functor formalism on a geometric setup $(\cat C, E)$.
\begin{propenum}
	\item \label{rslt:functoriality-of-kerncat-in-D} Let $\alpha\colon \D \to \D'$ be a morphism of 3-functor formalisms on $(\cat C, E)$, i.e.\ a natural transformation of functors $\Corr(\cat C, E)^\tensor \to \Cat^\times$. Then for every $S \in \cat C$, there is an induced 2-functor
	\begin{align}
		\varphi_\alpha\colon \cat K_{\D,S} \to \cat K_{\D',S}
	\end{align}
	which acts on objects as the identity and on morphism categories as
	\begin{align}
		\D(X \times_S Y) \xlongto{\alpha} \D'(X \times_S Y).
	\end{align}

	\item Let $f\colon (\cat C', E') \to (\cat C, E)$ be a morphism of geometric setups and denote by $f_*\D$ the 3-functor formalism on $(\cat C', E')$ given by precomposing $\D$ with the lax symmetric monoidal functor $\Corr(\cat C', E') \to \Corr(\cat C, E)$. Then for every $S' \in \cat C'$ there is an induced 2-functor
	\begin{align}
		\rho_f\colon \cat K_{f_*\D,S'} \to \cat K_{\D,f(S')}.
	\end{align}
	which acts on objects as $f$ and on morphism categories as
	\begin{align}
		\D(f(X' \times_{S'} Y')) \xlongto{(f_{X'Y'})_!} \D(f(X') \times_{f(S')} f(Y')),
	\end{align}
	where $f_{X'Y'}\colon f(X' \times_{S'} Y') \to f(X') \times_{f(S')} f(Y')$ is the natural map.
\end{propenum}
\end{proposition}
\begin{proof}
Part (i) follows from the 2-functoriality of the transfer of enrichment, see \cref{rslt:2-functoriality-of-transfer-of-enrichment}. To prove (ii), let everything be given as in the claim. To simplify notation we write $\cat C'$ for $\cat C'_{/S'}$ and $\cat C$ for $\cat C_{/f(S')}$. Denote $\alpha\colon \Corr(\cat C') \to \Corr(\cat C)$ the lax symmetric monoidal functor induced by $f$. Then by \cref{rslt:functor-from-transfer-of-self-enrichment} there is a natural $\Corr(\cat C)$-enriched functor $\tau_\alpha(\Corr(\cat C')) \to \Corr(\cat C)$. By applying the transfer of enrichment $\tau_\D$ along $\D$, we arrive at the desired 2-functor $\rho_f$. From the explicit description in \cref{rslt:functor-from-transfer-of-self-enrichment} we deduce the explicit description of $\rho_f$.
\end{proof}

It remains to study the functoriality of $\cat K_{\D,S}$ in $S$. It turns out that the assignment $S \mapsto \cat K_{\D,S}$ is itself some kind of \enquote{3-functor-formalism}, as we will see in \cref{rslt:functoriality-of-kerncat} below. Before proving that result, we need the following preparation:

\begin{lemma} \label{rslt:functoriality-of-Corr-in-Corr}
Let $(\cat C, E)$ be a geometric setup. Then there is a functor
\begin{align}
	\Corr(\cat C, E) \to \Enr_{\Corr(\cat C, E)}, \qquad S \mapsto \Corr((\cat C_E)_{/S}).
\end{align}
Here we view $\Corr((\cat C_E)_{/S})$ as $\Corr(\cat C, E)$-enriched by using transfer of enrichment along the lax symmetric monoidal functor
\begin{align}
	\Corr((\cat C_E)_{/S})^\tensor \to \Corr(\cat C, E)^\tensor
\end{align}
induced by the forgetful functor $(\cat C_E)_{/S} \to \cat C$.
\end{lemma}
\begin{proof}
There should be a way to obtain this result by applying a variant of \cref{rslt:construct-3ff-from-suitable-decomp}, where $I = E$ and $P$ is the collection of isomorphisms in $\cat C$. In fact, in this case one can use \cite[Theorem~1.1.1]{Stefanich:Correspondences}. We leave the details of this approach to the reader and instead perform an explicit construction of the desired functor.

By the definition of the cocartesian operad in \cite[Construction~2.4.3.1]{HA} we have
\begin{align}
	\Fun_{\Fin_*}(\cat C^\op \times \Fin_*, \cat C^{\op,\amalg}) = \Fun(\cat C^\op \times \Gamma^*, \cat C^\op),
\end{align}
where $\Gamma^*$ is the category with objects $(\langle n \rangle, i)$ for $\langle n \rangle \in \Fin_*$ and $i \in \{ 1, \dots, n \}$. In particular the projection $\cat C^\op \times \Gamma^* \to \cat C^\op$ induces a functor $\cat C^\op \times \Fin_* \to \cat C^{\op,\amalg}$ over $\Fin_*$, sending $(S, \langle n \rangle) \mapsto S^{\dsum n}$. Now let $\cat E' \subseteq \Fun([1], \cat C)$ be the full subcategory spanned by the edges in $E$ and let
\begin{align}
	\cat E \coloneqq \big(\cat E'^{\op,\amalg} \times_{\cat C^{\op,\amalg}} (\cat C^\op \times \Fin_*)\big)^\op,
\end{align}
where the functor $\cat E'^{\op,\amalg} \to \cat C^{\op,\amalg}$ is induced by the target map. The category $\cat E$ has the following more explicit description: An object in $\cat E$ is a finite family of maps $(X_i \to S)_i$ in $E$ (with fixed target) and a morphism in $\cat E$ is a commuting diagram
\begin{equation}\begin{tikzcd}
	(Y_j)_{1\le j\le m} \arrow[r] \arrow[d] & (X_i)_{1\le i\le n} \arrow[d]\\
	T \arrow[r] & S
\end{tikzcd}\end{equation}
where the upper horizontal map is a map in $\cat C^{\op,\amalg,\op}$, i.e.\ a collection of maps $Y_{\alpha(i)} \to X_i$ for some $\alpha\colon \langle n \rangle \to \langle m \rangle$, and the bottom horizontal map is a map in $\cat C$. The commutativity of the diagram means that it commutes for each fixed map $Y_{\alpha(i)} \to X_i$ at the top for which $\alpha(i) \ne *$. By construction there is a natural functor
\begin{align}
	t\colon \cat E \to \cat C \times \Fin_*^\op.
\end{align}
Let $E_{\cat E}$ be the collection of edges in $\cat E$ where (in the above description) $m = n$, $\alpha = \id$, and $T \to S$ is in $E$ (then automatically all the maps $Y_i \to X_i$ are in $E$ by cancellativeness, see \cref{rslt:E-is-cancellative}). A straightforward calculation shows that the pair $(\cat E, E_{\cat E})$ is a geometric setup, and $t$ extends to a morphism of geometric setups
\begin{align}
	t\colon (\cat E, E_{\cat E}) \to (\cat C, E) \times (\Fin_*^\op, \isom).
\end{align}
By applying $\Corr$ (and using that it commutes with products by \cref{rslt:Corr-preserves-limits}) we obtain a functor
\begin{align}
	p\colon \cat V^\tensor \coloneqq \Corr(\cat E, E_{\cat E}) \to \Corr(\cat C, E) \times \Fin_*.
\end{align}
We claim that $p$ is a $\Corr(\cat C, E)$-cocartesian family of operads, i.e.\ it satisfies the following conditions:
\begin{enumerate}[(a)]
	\item $p$ is a map of cocartesian fibrations over $\Corr(\cat C, E)$.
	\item The induced straightening $\Corr(\cat C, E) \to \Cat_{/\Fin_*}$ factors over $\Op$.
\end{enumerate}
To prove (a), we first verify that $\cat V^\tensor \to \Corr(\cat C, E)$ is a cocartesian fibration. By \cref{rslt:Corr-cocartesian-fibration} this boils down to understanding the behavior of $t'$-cartesian and $t'$-cocartesian edges, where $t'$ is the functor $t'\colon \cat E \to \cat C$. For a map $f\colon T \to S$ in $\cat C$ and some $(X_i)_i \to S$ in the fiber of $\cat E$ over $S$, the associated $t'$-cartesian lift is given by $(X_i \times_S T)_i \to (X_i)_i$. For a map $g\colon T \to S$ in $E$ and some $(Y_j)_j \to T$ in the fiber of $\cat E$ over $T$, the associated $t'$-cocartesian lift is given by the identity $(Y_j)_j \to (Y_j)_j$. We picture both the cartesian and the cocartesian case in the following diagrams:
\begin{equation}
	\begin{tikzcd}
		(X_i \times_S T)_i \arrow[d] \arrow[r,"\pr_1"] & (X_i)_i \arrow[d]\\
		T \arrow[r,"f"] & S
	\end{tikzcd}
	\qquad
	\begin{tikzcd}
		(X_i)_i \arrow[d] \arrow[r,"\id"] & (X_i)_i \arrow[d]\\
		T \arrow[r,"g"] & S
	\end{tikzcd}
\end{equation}
We leave it to the reader to verify these claims by reducing them to \cref{ex:cocartesian-fibration}. We deduce that condition (a) of \cref{rslt:Corr-cocartesian-fibration} is satisfied, and from the explicit description of cartesian and cocartesian edges we also deduce that conditions (b) and (c) are satisfied. This proves that $\cat V^\tensor \to \Corr(\cat C, E)$ is indeed a cocartesian fibration. Moreover, the explicit description of cocartesian morphisms in \cref{rslt:Corr-cocartesian-fibration} easily implies that $p$ is a map of cocartesian fibrations, proving claim (a) from above.

By straightening $p$ we obtain a functor $F\colon \Corr(\cat C, E) \to \Cat_{/\Fin_*}$. We now prove (b) from above, i.e.\ that $F$ factors through $\Op$. Fix an object $S \in \cat C$. By \cref{rslt:Corr-preserves-limits} we have
\begin{align}
	F(S) = \cat V^\tensor \times_{\Corr(\cat C, E)} \{ S \} = \Corr(\cat E \times_{\cat C} \{ S\}, E_{\cat E,S}),
\end{align}
where $E_{\cat E,S}$ is the collection of edges in $\cat E_S \coloneqq \cat E \times_{\cat C} \{ S \}$ which lie in $E_{\cat E}$. To show that $F(S)$ is an operad, it is enough to show the following stronger statement: There is a natural isomorphism $(\cat E_S, E_{\cat E,S}) = (((\cat C_E)_{/S})^{\op,\amalg,\op}, E^-)$ of geometric setups over $(\Fin_*, \isom)$, where $E^-$ is as in \cref{rslt:Corr^tensor-explicit}. Then by \cref{rslt:Corr^tensor-explicit} we deduce $F(S) = \Corr((\cat C_E)_{/S})^\tensor$, which is an operad. To construct the desired natural isomorphism of geometric setups, one easily reduces to constructing a natural isomorphism
\begin{align}
	\cat E_S^\op = ((\cat C_E)_{/S})^{\op,\amalg}.
\end{align}
Now note that for every category $\cat K$ over $\Fin_*$ we have
\begin{align}
	\Fun_{\Fin_*}(\cat K, \cat E_S^\op) &= \Fun_{\Fin_*}(\cat K, \cat E'^{\op,\amalg} \times_{\cat C^{\op,\amalg}} (\{ S \} \times \Fin_*))\\
	&= \Fun(\cat K \times_{\Fin_*} \Gamma^*, \cat E'^\op \times_{\cat C^\op} \{ S \})\\
	&= \Fun(\cat K \times_{\Fin_*} \Gamma^*, ((\cat C_E)_{/S})^\op)\\
	&= \Fun_{\Fin_*}(\cat K, ((\cat C_E)_{/S})^{\op,\amalg}),
\end{align}
so we conclude by Yoneda. This shows that $F(S)$ is an operad. Now it is easy to see that $F$ maps edges in $\Corr(\cat C, E)$ to operad maps and hence factors through $\Op$, proving (b). Altogether we have now proved the following claim:
\begin{itemize}
	\item[($*$)] The functor $p\colon \cat V^\tensor \to \Corr(\cat C, E) \times \Fin_*$ is a cocartesian family of operads and straightens to the functor
	\begin{align}
	 	F\colon \Corr(\cat C, E) \to \Op, \qquad S \mapsto \Corr((\cat C_E)_{/S})^\tensor.
	\end{align}
	The functor $F$ sends the morphism $\CorrHom{S}{f}{T}{\id}{T}$ to the functor $X \mapsto X \times_S T$ and a morphism $\CorrHom{T}{\id}{T}{g}{S}$ to the functor $(X \to T) \mapsto (X \to S)$ (on the respective correspondence categories).
\end{itemize}
With $F$ and $p$ at hand, we can now perform the construction of the desired functor. Namely, note that there is a natural commuting diagram
\begin{equation}\begin{tikzcd}
	\cat V^\tensor \arrow[rr,"\alpha"] \arrow[dr,"p",swap] && \Corr(\cat C, E) \times \Corr(\cat C, E)^\tensor \arrow[dl,"q"]\\
	& \Corr(\cat C, E) \times \Fin_*
\end{tikzcd}\end{equation}
of categories as follows. The functor $q$ is built as a product of the identity on $\Corr(\cat C, E)$ and the structure map $\Corr(\cat C, E)^\tensor \to \Fin_*$. The functor $\alpha$ is built from the functor $\cat V^\tensor \to \Corr(\cat C, E)$ as part of $p$ and the functor $\cat V^\tensor \to \Corr(\cat C, E)^\tensor$ induced by the source map $\cat E' \to \cat C$ by applying $(\blank)^{\op,\amalg}$ and $\Corr$; this is clearly a functor over $\Fin_*$, hence the above diagram commutes. Concretely, $\alpha$ sends an object $(X_i \to S)_i$ to the pair $(S, (X_i)_i)$. Clearly $\alpha$ is a morphism of $\Corr(\cat C, E)$-families of operads, so we can apply \cref{rslt:transfer-of-enrichment-in-families} with $\cat S = \Corr(\cat C, E)$, $F = F$ (which we can view as a functor to $\Enr$ by \cref{rslt:Corr-self-enriched-and-self-dual}), $\cat W^\tensor = \Corr(\cat C, E) \times \Corr(\cat C, E)^\tensor$ and $\alpha = \alpha$; here we implicitly pull back everything from $\Fin_*$ to $\Ass^\tensor$. By \cref{rslt:transfer-of-enrichment-in-families} we obtain a functor
\begin{align}
	G \coloneqq \tau_\alpha(F)\colon \Corr(\cat C, E) \to \Enr.
\end{align}
Moreover, by construction, the composition of $G$ with the projection $\Enr \to \Op_{/\Ass^\tensor}$ is the constant functor sending everything to $\Corr(\cat C, E)^\tensor$, so we can view $G$ as a functor
\begin{align}
	G\colon \Corr(\cat C, E) \to \Enr_{\Corr(\cat C, E)},
\end{align}
as desired.
\end{proof}

\begin{corollary} \label{rslt:functoriality-of-Corr-in-Corr-monoidal}
Let $(\cat C, E)$ be a geometric setup. Then the functor in \cref{rslt:functoriality-of-Corr-in-Corr} induces a lax symmetric monoidal functor
\begin{align}
	\Corr(\cat C, E)^\tensor \to \Cat^\times, \qquad S \mapsto \Corr((\cat C_E)_{/S}).
\end{align}
\end{corollary}
\begin{proof}
We apply \cref{rslt:functoriality-of-Corr-in-Corr} to the geometric setup $(\cat C', E') = (\cat C^{\op,\amalg,\op}, E^-)$ from \cref{rslt:Corr^tensor-explicit} in order to obtain a functor
\begin{align}
	\Corr(\cat C, E)^\tensor = \Corr(\cat C', E') \to \Enr_{\Corr(\cat C', E')} \to \Cat,
\end{align}
where the last functor in the above composition is the functor sending an enriched category to its underlying category. It is straightforward to check that the just constructed functor $\Corr(\cat C, E)^\tensor \to \Cat$ is a lax cartesian structure (in the sense of \cref{ex:cartesian-monoidal-structure}) and hence produces the desired lax symmetric monoidal functor.
\end{proof}

We can finally prove the promised functoriality of the 2-category of kernels $\cat K_{\D,S}$ in $S$. In fact, this is an easy consequence of \cref{rslt:functoriality-of-Corr-in-Corr}:

\begin{theorem} \label{rslt:functoriality-of-kerncat}
Let $\D$ be a 3-functor formalism on a geometric setup $(\cat C, E)$. Then the assignment $S \mapsto \cat K_{\D,S}$ upgrades to a lax symmetric monoidal functor
\begin{align}
	\cat K_{\D,(\blank)}\colon \Corr(\cat C, E)^\tensor \to \TwoCat^\times, \qquad S \mapsto \cat K_{\D,S}.
\end{align}
It has the following explicit description:
\begin{thmenum}
	\item For every $S \in \cat C$ the 2-category $\cat K_{\D,S}$ is the 2-category of kernels from \cref{def:relative-category-of-kernels}.

	\item \label{rslt:functoriality-of-kerncat-pullbacks} For a morphism $\CorrHom{S}{f}{T}{}{T}$ in $\Corr(\cat C, E)$ the induced 2-functor
	\begin{align}
		f^*\colon \cat K_{\D,S} \to \cat K_{\D,T}
	\end{align}
	acts on objects as $X \mapsto X \times_S T$ and on morphism categories as
	\begin{align}
		\D(X \times_S Y) \xlongto{f_{XY}^*} \D((X \times_S T) \times_T (Y \times_S T)),
	\end{align}
	where $f_{XY}$ is the obvious map induced by $f$.

	\item \label{rslt:functoriality-of-kerncat-pushforward} For a morphism $\CorrHom{T}{}{T}{g}{S}$ in $\Corr(\cat C, E)$, the induced 2-functor
	\begin{align}
		g_!\colon \cat K_{\D,T} \to \cat K_{\D,S}
	\end{align}
	acts on objects as the forgetful functor and on morphism categories as
	\begin{align}
		\D(X \times_T Y) \xto{g_{XY!}} \D(X \times_S Y),
	\end{align}
	where $g_{XY}$ is the obvious map (which lies in $E$ by \cref{rslt:E-closed-under-diagonals}, because it is a pullback of the diagonal of $g$).
\end{thmenum}
Moreover, the functor $\Phi_{\D,S}$ from \cref{rslt:functors-from-and-to-kercat} upgrades to a morphism
\begin{align}
	\Phi_{\D,(\blank)}\colon \Corr((\cat C_E)_{/(\blank)}) \to \cat K_{\D,(\blank)}
\end{align}
of lax symmetric monoidal functors $\Corr(\cat C, E)^\tensor \to \TwoCat^\times$, where the first functor is the one from \cref{rslt:functoriality-of-Corr-in-Corr-monoidal}.
\end{theorem}
\begin{proof}
We first prove everything without symmetric monoidal structures. Let $G\colon \Corr(\cat C, E) \to \Enr_{\Corr(\cat C, E)}$ be the functor from \cref{rslt:functoriality-of-Corr-in-Corr}. The desired functor $\cat K_{\D,(\blank)}\colon \Corr(\cat C, E) \to \TwoCat$ can then be obtained by composing $G$ with the transfer of enrichment $\tau_{\D}$. The explicit description in (ii) and (iii) follows from the construction and the description in \cref{rslt:transfer-of-enrichment-in-families}. It remains to prove the final part. To construct $\Phi_{\D,(\blank)}$, we note that by the definition of cocartesian fibrations in \cref{def:cocartesian-fibration}, applied to the cocartesian fibration $\Enr \to \Op_{/\Ass^\tensor}$, the morphism $\D\colon \Corr(\cat C, E)^\tensor \to \Cat^\times$ in $\Op_{/\Ass^\tensor}$ induces an equivalence
\begin{align}
\Fun_\D([1], \Enr)^\cocart \isoto \Enr_{\Corr(\cat C, E)},
\end{align}
identifying $\Corr(\cat C, E)$-enriched categories $\cat X$ with the associated cocartesian edge $\cat X \to \tau_\D(\cat X)$ in $\Enr$. By composing $G$ with the inverse of this equivalence we obtain a functor $\Corr(\cat C, E) \to \Fun([1], \Enr)$, equivalently a natural transformation of functors $\Corr(\cat C, E) \to \Enr$. By composing with the functor $\Enr \to \Cat$ sending an enriched category to its underlying category, we arrive at $\Phi_{\D,(\blank)}$.

We now upgrade everything to the symmetric monoidal setting. Recall from \cref{ex:cartesian-monoidal-structure} that the category of lax symmetric monoidal functors $\Corr(\cat C, E)^\tensor \to \TwoCat^\times$ is equivalent to the category of lax cartesian structures, i.e.\ functors $\Corr(\cat C, E)^\tensor \to \TwoCat$ satisfying certain compatibilities. Now let $(\cat C', E') = (\cat C^{\op,\amalg,\op}, E^-)$ be the geometric setup from \cref{rslt:Corr^tensor-explicit}. In order to get the desired symmetric monoidal upgrade of the constructions in the first part of the proof, it is enough to upgrade $\D$ to a 3-functor formalism $\D'\colon \Corr(\cat C', E')^\tensor \to \Cat^\times$ sending $S \dsum T \mapsto \D(S) \times \D(T)$, because then we can simply apply the first part of the proof to $\D'$ in place of $\D$. Constructing $\D'$ is equivalent to constructing a functor $\Corr(\cat C', E')^\tensor \to \Cat$ which is a lax cartesian structure. Since $\D$ has a similar description, it is enough to produce a (not lax symmetric monoidal!) functor $\Corr(\cat C', E')^\tensor \to \Corr(\cat C, E)^\tensor$. By definition of these operads, this construction reduces to constructing a functor $(\cat C^{\op,\amalg,\op})^{\op,\amalg,\op} \to \cat C^{\op,\amalg,\op}$ which induces a morphism of the associated geometric setups. Simplifying notation a bit, we are reduced to constructing a functor $\cat C^{\op,\amalg,\amalg} \to \cat C^{\op,\amalg}$. This functor is constructed in \cref{rslt:functor-from-C-amalg-to-C} (observe that $\cat C^{\op,\amalg}$ has finite coproducts given by $X \dunion Y = X \dsum Y$).
\end{proof}

It is natural to ask in what sense the functors $\Psi_{\D,S}\colon \cat K_{\D,S} \to \Cat$ are functorial in $S$. They should compose to a \emph{lax} natural transformation $\cat K_{\D,(\blank)} \to \Cat$ of functors $\Corr(\cat C, E) \to \TwoCat$, where we identify $\Cat$ with the associated constant functor. As we do not have developed the theory of lax natural transformations here, we content ourselves with the following special case of that functoriality:

\begin{proposition} \label{rslt:functoriality-of-Psi-D}
Let $\D$ be a 3-functor formalism on a geometric setup $(\cat C, E)$ and let $\CorrHom{T}{g}{T'}{f}{S}$ be a morphism in $\Corr(\cat C, E)$. Then there is a natural transformation
\begin{align}
	\sigma\colon \Psi_{\D,T} \to \Psi_{\D,S} \comp f_!g^*
\end{align}
of 2-functors $\cat K_{\D,T} \to \Cat$.
For each $X \in \cat K_{\D,T}$ the induced map
\begin{align}
	\sigma(X)\colon \Psi_{\D,T}(X) = \D(X) \to \D(X \times_T T') = \Psi_{\D,S}(f_! g^*(X))
\end{align}
is given by pullback along $X \times_T T' \to X$.
\end{proposition} 
\begin{proof}
We note that there is a diagram
\begin{equation}\begin{tikzcd}
	\Corr((\cat C_E)_{/T}) \arrow[d,"f_! g^*",swap] \arrow[dr,""{name=psi,below}]\\
	\Corr((\cat C_E)_{/S}) \arrow[r] \arrow[Rightarrow,from=psi] & \Corr(\cat C, E) \arrow[r,"\D"] & \Cat
\end{tikzcd}\end{equation}
of closed monoidal categories, lax monoidal functors and a natural transformation between them. Here the left-hand triangle comes from the functor $\alpha$ constructed in the proof of \cref{rslt:functoriality-of-Corr-in-Corr}. Now apply \cref{rslt:functoriality-of-functor-from-transfer-of-self-enrichment} to this diagram.
\end{proof}

\begin{remark}\label{rmk:kerncat-symmetric-monoidal-2-cat}
Let $\D$ be a 3-functor formalism on a geometric setup $(\cat C, E)$. It follows from \cref{rslt:functoriality-of-kerncat} that each $\cat K_{\D,S}$ is a symmetric monoidal 2-category with respect to $\blank\times_S\blank$.
\end{remark}

From the perspective of a 3-functor formalism, the functor $\cat K_{\D,(\blank)}$ is rather boring, as the exceptional pushforwards are completely determined by the pullbacks and vice-versa:

\begin{lemma} \label{rslt:adjunction-between-functors-of-kerncat}
Let $\D$ be a 3-functor formalism on a geometric setup $(\cat C, E)$ and let $f\colon T \to S$ be a morphism in $E$ with associated 2-functors
\begin{align}
	f_!\colon \cat K_{\D,T} \rightleftarrows \cat K_{\D,S} \noloc f^*
\end{align}
from \cref{rslt:functoriality-of-kerncat}. Then $f_!$ is left and right adjoint to $f^*$.
\end{lemma}
\begin{proof}
We can replace $\cat C$ by $\cat C_E$ to assume that $E = \mathit{all}$. Let us denote by $\D_0\colon \Corr(\cat C) \to \Cat^\times$ the 3-functor formalism from \cref{rslt:functoriality-of-Corr-in-Corr-monoidal}, i.e.\ $\D_0(S) = \Corr(\cat C_{/S})$, and let us define
\begin{align}
	\Corr_2(\cat C) \coloneqq \cat K_{\D_0} = \tau_{\D_0}(\Corr(\cat C)).
\end{align}
Concretely, $\Corr_2(\cat C)$ has the same objects and morphisms as $\Corr(\cat C)$, but allows additional 2-morphisms, which are themselves correspondences. This \enquote{iterated correspondence category} is studied in \cite{Haugseng.2017a} and it is shown in \cite[Lemma~12.3]{Haugseng.2017a} that all morphisms in $\Corr_2(\cat C)$ have both left and right adjoints, given by swapping the correspondence (to be precise, \cite{Haugseng.2017a} uses a different construction for $\Corr_2(\cat C)$, but the proof of the adjunctions works also in our setting). It is thus enough to show that the functor from \cref{rslt:functoriality-of-kerncat} upgrades to a 2-functor
\begin{align}
	\Corr_2(\cat C) \to \TwoCat, \qquad S \mapsto \cat K_{\D,S}.
\end{align}
From \cref{rslt:functoriality-of-kerncat} we obtain a morphism $\Psi_{\D,(\blank)}\colon \D_0 \to \cat K_{\D,(\blank)}$ of $\TwoCat$-valued 3-functor formalisms. In the following we note that that \cref{rslt:functoriality-of-kerncat-in-D} and \cref{rslt:functors-from-and-to-kercat} work analogously for $\TwoCat$-valued 3-functor formalisms, in which case they deal with 3-categories and 3-functors (i.e.\ $\TwoCat$-enriched categories and functors). By the analog of \cref{rslt:functoriality-of-kerncat-in-D} we obtain a 3-functor $\Corr_2(\cat C) \to \cat K_{\cat K_{\D,(\blank)}}$. By composing this functor with the 3-functor $\Psi_{\cat K_{\D,(\blank)}} = \Fun(*,\blank)\colon \cat K_{\cat K_{\D,(\blank)}} \to \TwoCat$ from the analog of \cref{rslt:functors-from-and-to-kercat} we arrive at the desired functor $\Corr_2(\cat C) \to \TwoCat$.
\end{proof}

To make \cref{rslt:adjunction-between-functors-of-kerncat} more explicit, let everything be given as in that result and pick objects $X \in \cat K_{\D,S}$ and $Y \in \cat K_{\D,T}$. Then by \cref{rslt:adjunctions-in-Enr-V} we get induced equivalences of categories
\begin{align}
	\Fun_{\cat K_{\D,S}}(f^* X, Y) = \Fun_{\cat K_{\D,T}}(X, f_! Y), \qquad \Fun_{\cat K_{\D,S}}(Y, f^* X) = \Fun_{\cat K_{\D,T}}(f_! Y, X).
\end{align}
Explicitly, $X$ and $Y$ correspond to maps $X \to S$ and $Y \to T$ in $E$, and the above equivalences amount to an equivalence
\begin{align}
	\D(Y \times_T (X \times_S T)) = \D(Y \times_S X).
\end{align}
By going through the construction of the adjunction, we note that this equivalence is induced by the obvious isomorphism $Y \times_T (X \times_S T) = Y \times_S X$.

\subsection{Descent} \label{sec:kerncat.descent}

In this subsection we investigate how descent of a 6-functor formalism translates into properties of the associated category of kernels. There are two such results one may consider: firstly, how descent for $\D$ affects the functor $S \mapsto \cat K_{\D,S}$, and secondly, how descent for $\D$ is reflected in the category $\cat K_{\D,S}$ for fixed $S$. Using the explicit description of homomorphism categories in $\cat K_{\D,S}$, both questions are easily answered. Let us start with the first case:

\begin{proposition} \label{rslt:*-descent-implies-weak-descent-for-K-D}
Let $\D$ be a 3-functor formalism on a geometric setup $(\cat C, E)$ and let $\cat U$ be a sieve over some $X \in \cat C$ such that $\D^*$ descends universally along $\cat U$. Then the 2-functor
\begin{align}
	\cat K_{\D,X} \injto \varprojlim_{U \in \cat U^\op} \cat K_{\D,U}
\end{align}
is fully faithful.
\end{proposition}
\begin{proof}
Fix $Y, Z \in \cat K_{\D,X}$, i.e.\ $Y$ and $Z$ are objects in $(\cat C_E)_{/X}$. Using the explicit computation of limits of 2-categories (see \cref{rslt:limits-of-enriched-categories}) and the explicit description of the functor $S \mapsto \cat K_{\D,S}$ we deduce that the claimed fully faithfulness amounts to saying that the natural functor
\begin{align}
	\D(Y \times_X Z) = \Fun_{\cat K_{\D,X}}(Z, Y) \isoto \varprojlim_{U\in \cat U^\op} \Fun_{\cat K_{\D,U}}(Z \times_X U, Y \times_X U) = \varprojlim_{U\in \cat U^\op} \D((Y \times_X Z) \times_X U)
\end{align}
is an equivalence (for all $Y$ and $Z$ as above). But this is the same as saying that $\D$ descends along $f^* \cat U$, where $f\colon Y \times_X Z \to X$ is the projection.
\end{proof}

\begin{remark}
We do not know if in \cref{rslt:*-descent-implies-weak-descent-for-K-D} it is even true that $\cat K_{\D,X} \isoto \varprojlim_U \cat K_{\D,U}$ is an equivalence of 2-categories, i.e.\ if this functor is essentially surjective. This should be true if $E = \mathit{all}$ and we replace $\cat K_{\D,X}$ by its 2-presentable version (see \cref{rmk:presentable-kerncat}), but we will not further pursue this approach here.
\end{remark}

The next result discusses how descent for $\D$ is reflected in the category of kernels with fixed base. Recall the notion of limits in a 2-category from \cref{def:limits-in-a-2-category}.

\begin{proposition} \label{rslt:*-descent-implies-limit-in-kerncat}
Let $\D$ be a 3-functor formalism on a geometric setup $(\cat C, E)$ with $E = \mathit{all}$ and let $\cat U$ be a sieve on some $X \in \cat C$ such that $\D^*$ descends universally along $\cat U$. Then the natural map
\begin{align}
	\varinjlim_{U\in \cat U} U \isoto X
\end{align}
is an isomorphism in $\cat K_\D$.
\end{proposition}
\begin{proof}
By definition we need to verify that for every $Z \in \cat C$ the natural map
\begin{align}
	\varprojlim_{U \in \cat U^\op} \D(Z \times U) = \varprojlim_{U \in \cat U^\op} \Fun_{\cat K_\D}(U, Z) \isoto \Fun_{\cat K_\D}(Z, X) = \D(Z \times X)
\end{align}
is an equivalence of categories. Here the transition functors are given by pullbacks along the obvious maps. But then the claimed equivalence of categories boils down to the descent of $\D^*$ along the sieve $f^* \cat U$, where $f\colon X\times Z \to X$ is the projection.
\end{proof}

\subsection{Suave and prim objects} \label{sec:kerncat.suave-prim-obj}

A fundamental 2-categorical notion is that of adjoint morphisms, which we review in \cref{sec:2cat.adj}. It was first observed in \cite[\S~IV.2]{Fargues-Scholze:Geometrization} (based on similar ideas of \cite{Lu-Zheng:ULA}) that adjoint morphisms in the category of kernels are a surprisingly fruitful notion. In the following we provide all the relevant definitions and basic properties of adjunctions in $\cat K_\D$. Most of these results have already appeared in some form in the literature, but many of the proofs are unsatisfying (or even incomplete) and assume statements about (classical) 2-categories that are not properly verified. In the following we fill all the gaps. Let us start with the core definition:

\begin{definition} \label{def:adm-and-coadm-sheaves}
Let $\D$ be a 3-functor formalism on a geometric setup $(\cat C, E)$, let $f\colon X \to S$ be a morphism in $E$ and fix $P \in \D(X)$.
\begin{defenum}
	\item We say that $P$ is \emph{$f$-suave} if, when viewed as a morphism $X \to S$ in $\cat K_{\D,S}$, it is a left adjoint morphism. We denote the associated right adjoint morphism $S \to X$ by $\DSuave_f(P) \in \D(X)$ and call it the \emph{$f$-suave dual} of $P$. We denote by $\Suave_f(X) \subseteq \D(X)$ the full subcategory spanned by the $f$-suave objects.
	\item We say that $P$ is \emph{$f$-prim} if, when viewed as a morphism $X \to S$ in $\cat K_{\D,S}$, it is a right adjoint morphism. We denote the associated left adjoint morphism $S \to X$ by $\DPrim_f(P) \in \D(X)$ and call it the \emph{$f$-prim dual} of $P$. We denote by $\Prim_f(X) \subseteq \D(X)$ the full subcategory spanned by the $f$-prim objects.
\end{defenum}
If $f$ is clear from context we may drop it from the notation or replace it by $S$.
\end{definition}

\begin{remark} \label{rmk:origin-of-suave-and-prim-objects}
The above definition of $f$-suave sheaves/objects first appeared in \cite[Theorem~IV.2.23(iii)]{Fargues-Scholze:Geometrization} based on ideas of \cite{Lu-Zheng:ULA}. In that reference, the definition only appeared in the context of étale sheaves on diamonds, in which case it coincides with the notion of $f$-ULA sheaves. In \cite{Mann.2022b} $f$-suaveness was introduced under the name \enquote{$f$-dualizable sheaves} (now in the slightly more general context of nuclear sheaves on diamonds) and many formal properties were deduced, the proofs of which are abstract enough to work in any 6-functor formalism. In \cite{Mann.2022b} the notion of $f$-prim sheaves was introduced for the first time (under the name \enquote{$f$-proper sheaves}). Both $f$-suave and $f$-prim sheaves were later studied in \cite{Scholze:Six-Functor-Formalism} (under the names of \enquote{$f$-smooth} and \enquote{$f$-proper} sheaves) in the context of abstract 6-functor formalisms and some basic properties were shown there. We chose to deviate slightly from the previous terminology in order to avoid clashes with smooth and proper morphisms, which are introduced in \cref{sec:kerncat.suave-prim-map,sec:kerncat.etale-proper} below. We thank Peter Scholze and David Hansen for suggesting the names \enquote{suave} and \enquote{prim}.
\end{remark}

To provide some further explanation, we observe that by construction of the category of kernels the morphism categories in $\cat K_{\D,S}$ are given by
\begin{align}
	\Fun_{\cat K_{\D,S}}(X, S) = \D(S \times_S X) = \D(X) = \D(X \times_S S) = \Fun_{\cat K_{\D,S}}(S, X),
\end{align}
i.e.\ we can indeed view objects $P \in \D(X)$ both as morphisms $X \to S$ and $S \to X$ in $\cat K_{\D,S}$. Suave sheaves take well-known shape in many classical 6-functor formalisms, e.g.\ in the form of admissible representations (see \cref{rslt:suave=admissible}) or universally locally acyclic (ULA) sheaves (see e.g.\ \cite[Theorem~IV.2.23]{Fargues-Scholze:Geometrization}). On the other hand, primness is somewhat harder to grasp, but closely related to being a compact object in $\D(X)$ (see \cref{rslt:relation-between-suave-prim-and-compact} below).

\begin{example} \label{rslt:suave-and-prim-for-identity-is-dualizable}
In \cref{def:adm-and-coadm-sheaves} assume that $X = S$ and $f = \id_S$. Then an object $P \in \D(S)$ is $f$-suave if and only if it is $f$-prim if and only if it is dualizable (cf. \cref{rslt:adj-in-*-V-same-as-dualizable}).
\end{example}

It is shown in \cref{rslt:passing-to-adjoints} that passing to adjoint morphisms in a 2-category is a functorial operation, which leads to the following observation:

\begin{lemma}\label{rslt:SD-PD-are-self-inverse-equivalences}
Let $\D$ be a 3-functor formalism on a geometric setup $(\cat C, E)$ and let $f\colon X \to S$ be a morphism in $E$.
\begin{lemenum}
	\item The assignment $P \mapsto \DSuave_f(P)$ defines a self-inverse equivalence
	\begin{align}
		\DSuave_f\colon \Suave_f(X)^\op \isoto \Suave_f(X).
	\end{align}

	\item The assignment $P \mapsto \DPrim_f(P)$ defines a self-inverse equivalence
	\begin{align}
	 	\DPrim_f\colon \Prim_f(X)^\op \isoto \Prim_f(X).
	\end{align}
\end{lemenum}
\end{lemma}
\begin{proof}
By \cref{rslt:passing-to-adjoints} there is an isomorphism $\cat K_{\D,S}^L \isom (\cat K_{\D,S}^R)^{\co,\op}$ which acts as the identity on objects and sends a left adjoint morphism to its corresponding right adjoint. In particular this isomorphism induces the desired equivalence
\begin{align}
	\DSuave_f\colon \Suave_f(X) = \Fun^L_{\cat K_{\D,S}}(X,S) \isoto \Fun^R_{\cat K_{\D,S}}(S,X)^\op = \Fun^L_{\cat K_{\D,S}}(X, S)^\op = \Suave_f(X)^\op,
\end{align}
where in the second to last step we used the equivalence $\cat K_{\D,S} = \cat K_{\D,S}^\op$ from \cref{rslt:kern-cat-is-self-dual}. It remains to prove that $\DSuave_f$ and $\DPrim_f$ are self-dual. By the construction of passage to the adjoint morphism in \cref{rslt:passing-to-adjoints-for-fixed-objects} we see that $\DSuave_f$ comes via the adjunction in \cref{rslt:pairing-of-categories} from the functor
\begin{align}
	\mu\colon \Suave_f(X) \times \Suave_f(X) \to \Ani, \qquad (P, Q) \mapsto \Hom(\id_S, P^\vee \comp Q),
\end{align}
where for a morphism $P\colon X \to S$ in $\cat K_{D,S}$ we denote by $P^\vee\colon S \to X$ the corresponding morphism via $\cat K_{\D,S} = \cat K_{\D,S}^\op$. Thus by \cref{rslt:pairing-of-categories} we obtain an adjunction
\begin{align}
	\DSuave_f\colon \Suave_f(X) \rightleftarrows \Suave_f(X)^\op \noloc \DSuave_f
\end{align}
which is also an equivalence, as noted above. The adjunction unit thus provides an equivalence $\id \isoto \DSuave_f \comp \DSuave_f$, as desired. The same argument works for $\DPrim_f$.
\end{proof}




In case that all six functors exist, there is a useful criterion for detecting suave and prim objects based on the pointwise criterion for adjunctions in \cref{rslt:pointwise-criterion-for-adjunction}. We also obtain an explicit formula for the suave and prim duals. Let us start with the suave case, where the following criterion first appeared in \cite[Theorem~VI.2.23.(ii)]{Fargues-Scholze:Geometrization}, based on \cite{Lu-Zheng:ULA}.

\begin{lemma} \label{rslt:criterion-for-suave-object}
Let $\D$ be a 6-functor formalism on a geometric setup $(\cat C, E)$, let $f\colon X \to S$ be a map in $E$ and $P \in \D(X)$. Then $P$ is $f$-suave if and only if the natural map
\begin{align}
	\pi_1^* \iHom(P, f^! \one) \tensor \pi_2^* P \isoto \iHom(\pi_1^* P, \pi_2^! P)
\end{align}
becomes an isomorphism after applying $\Hom(\one, \Delta^!(\blank))$. If this is the case, then
\begin{align}
	\DSuave_f(P) = \iHom(P, f^! \one),
\end{align}
and for all $Z\in (\cat C_E)_{/S}$ and all $M\in \D(Z)$ the natural map
\begin{align}
\pi_1^*\iHom(P,f^!\one) \tensor \pi_2^*M \isoto \iHom(\pi_1^*P, \pi_2^!M)
\end{align}
is an isomorphism in $\D(X\times_SZ)$.
\end{lemma}
\begin{proof}
We apply \cref{rslt:pointwise-criterion-for-adjunction} to the 2-category $\cat K_{\D,S}$ and the map $P\colon X \to S$. Fix some object $Z$ in $\cat K_{\D,S}$. Then the functor
\begin{align}
	P_*\colon \D(X \times_S Z) = \Fun_{\cat K_{\D,S}}(Z, X) \to \Fun_{\cat K_{\D,S}}(Z, S) = \D(Z)
\end{align}
is given by $\pi_{2!}(\pi_1^* P \tensor \blank)$, where $\pi_1$ and $\pi_2$ are the two projections from $X \times_S Z$. Consequently, $P_*$ admits a right adjoint $G_Z$ given by $G_Z(M) = \iHom(\pi_1^* P, \pi_2^!M)$. This shows that condition (a) of \cref{rslt:pointwise-criterion-for-adjunction} is always satisfied. Consequently, $P$ is $f$-suave if and only if condition (b) of \loccit{} is satisfied, i.e.\ the natural map $G_S(\id_S) \comp P \isoto G_X(P)$ becomes an isomorphism after applying $\Hom(\id_X,\blank)$. But by the above description of $G_Z$ this map is exactly the one in the claim, and since $\id_X$ is represented by $\Delta_! \one$ we deduce $\Hom(\id_X,\blank) = \Hom(\one, \Delta^!(\blank))$. This proves the first part of the claim. Moreover, if $P$ is $f$-suave then by \cref{rslt:pointwise-criterion-for-adjunction} the right adjoint is given by $G_S(\id_S) = \iHom(P, f^! \one)$ and the last asserted isomorphism holds, as desired.
\end{proof}

The next result provides a similar criterion as in \cref{rslt:criterion-for-suave-object} for $f$-prim objects. This criterion first appeared in \cite[Proposition~6.9]{Scholze:Six-Functor-Formalism}, but the proof of loc.\ cit.\ was left to the reader.

\begin{lemma} \label{rslt:criterion-for-prim-object}
Let $\D$ be a 6-functor formalism on a geometric setup $(\cat C, E)$, let $f\colon X \to S$ be a map in $E$ and $P \in \D(X)$. Then $P$ is $f$-prim if and only if the natural map
\begin{align}
	f_!(\pi_{2*} \iHom(\pi_1^* P, \Delta_! \one) \tensor P) \isoto f_* \iHom(P, P)
\end{align}
becomes an isomorphism after applying $\Hom(\one,\blank)$. Here $\pi_i\colon X \times_S X \to X$ denote the two projections and $\Delta\colon X \to X \times_S X$ is the diagonal. If $P$ is $f$-prim, then
\begin{align}
	\DPrim_f(P) = \pi_{2*} \iHom(\pi_1^* P, \Delta_! \one),
\end{align}
and for all $Z\in (\cat C_E)_{/S}$ and all $M\in \D(X\times_SZ)$ the natural map
\begin{align}
\pr_{2!}(\pr_1^*\pi_{2*}\iHom(\pi_1^*P,\Delta_!\one) \tensor M) \isoto \pr_{2*} \iHom(\pr_1^*P, M)
\end{align}
is an isomorphism in $\D(Z)$, where $\pr_1, \pr_2$ denote the two projections from $X\times_SZ$.
\end{lemma}
\begin{proof}
We argue very similarly to \cref{rslt:criterion-for-suave-object}, but now consider $P$ as a map $S \to X$ in $\cat K_{\D,S}$. For every $Z \in \cat K_{\D,S}$ the functor
\begin{align}
	P_*\colon \D(Z) = \Fun_{\cat K_{\D,S}}(Z, S) \to \Fun_{\cat K_{\D,S}}(Z, X) = \D(X \times_S Z)
\end{align}
is given by $\pr_1^* P \tensor \pr_2^*(\blank)$, where $\pr_1$ and $\pr_2$ are the two projections from $X \times_S Z$. This functor has a right adjoint $G_Z$ given by $G_Z(M) = \pr_{2*} \iHom(\pr_1^* P, M)$. Thus $P$ is $f$-prim if and only if the natural map $G_X(\id_X) \comp P \isoto G_S(P)$ is an isomorphism after applying $\Hom(\id_S,\blank)$. And if $P$ is $f$-prim, then \cref{rslt:pointwise-criterion-for-adjunction} shows that the natural map $G_X(\id_X) \comp M \isoto G_Z(M)$ is an isomorphism for all $Z \in \cat C$ and all $M\in \D(X\times_SZ)$.
By using the explicit formulas for composition in $\cat K_{\D,S}$ (see the remarks following \cref{def:kern-cat}) both isomorphisms reduce to the ones in the assertion. If $P$ is $f$-prim then its dual is given by $G_X(\id_X)$, which has the claimed formula.
\end{proof}

\begin{remark}
It is not really necessary to assume the existence of all six functors in \cref{rslt:criterion-for-suave-object,rslt:criterion-for-prim-object}. Instead, the same results apply for every 3-functor formalism as soon as all the functors appearing in the criteria (like $f^!$) exist.
\end{remark}

While the above criteria can be useful for detecting suave and prim objects in some cases, our experience suggests that one should usually try to avoid them. One of the main problems of the criteria is that in practice it is often hard to verify that the isomorphism one constructs is the natural one in the criteria. In the following we provide different tools for detecting suave and prim objects by reducing the question to simpler cases. We first observe that suaveness and primness are $*$-local on the target:

\begin{lemma} \label{rslt:suave-and-prim-obj-is-local-on-target}
Let $\D$ be a 3-functor formalism on a geometric setup $(\cat C, E)$, let $f\colon X \to S$ be a map in $E$ and let $P \in \D(X)$. For every map $g\colon S' \to S$ in $\cat C$ we denote $g'\colon X' \coloneqq X \times_S S' \to X$ and $f'\colon X' \to S'$ the base-changes. Then:
\begin{lemenum}
	\item If $P$ is $f$-suave then for every map $g\colon S' \to S$ in $\cat C$, $g'^* P$ is $f'$-suave.

	\item Let $\cat U$ be a sieve on $S$ along which $\D^*$ descends universally. If for every $g\colon U \to S$ in $\cat U$ the object $g'^* P$ is $f'$-suave, then $P$ is $f$-suave.
\end{lemenum}
The same result is true for primness in place of suaveness.
\end{lemma}
\begin{proof}
Part (i) follows from the fact that the 2-functor $g^*\colon \cat K_{\D,S} \to \cat K_{\D,S'}$ from \cref{rslt:functoriality-of-kerncat-pullbacks} sends $P$ to $g'^* P$ and preserves adjunctions (like every 2-functor). To prove (ii) let $\cat U$ be given as in the claim. Then by \cref{rslt:*-descent-implies-weak-descent-for-K-D} the 2-functor $\cat K_{\D,S} \injto \varprojlim_{U \in \cat U} \cat K_{\D,U}$ is fully faithful. Using \cref{rslt:functoriality-of-C-L-R} we deduce that a morphism in $\cat K_{\D,S}$ is a left adjoint if and only if its images in all $\cat K_{\D,U}$ are left adjoints. This immediately implies (ii).
\end{proof}

Having established locality on the \emph{target}, we now concern ourselves with locality on the \emph{source}. We first observe that suave and prim objects are stable under pullback along suave and prim maps, respectively, in the following sense (see also \cref{rslt:suave-pullback-and-prim-pushforward-preservations} for an important special case):

\begin{lemma} \label{rslt:composition-of-suave-and-prim-objects}
Let $\D$ be a 3-functor formalism on some geometric setup $(\cat C, E)$ and let $f\colon X \to S$ and $g\colon Y \to X$ be maps in $E$.
\begin{lemenum}
	\item \label{rslt:pullback-of-suave-and-prim-objects} Suppose that $P \in \D(X)$ is $f$-suave and $Q \in \D(Y)$ is $g$-suave. Then $g^* P \tensor Q$ is $fg$-suave and
	\begin{align}
		\DSuave_{fg}(g^* P \tensor Q) = g^* \DSuave_f(P) \tensor \DSuave_g(Q).
	\end{align}

	\item \label{rslt:pushforward-of-suave-and-prim-objects} Suppose that $P \in \D(Y)$ is $fg$-suave and $Q \in \D(Y)$ is $g$-prim. Then $g_!(Q \tensor P)$ is $f$-suave and
	\begin{align}
		\DSuave_f(g_!(Q \tensor P)) = g_!(\DPrim_g(Q) \tensor \DSuave_{fg}(P)).
	\end{align}
\end{lemenum}
The same is true with \enquote{suave} and \enquote{prim} swapped (and $\DPrim$ replaced by $\DSuave$).
\end{lemma}
\begin{proof}
We first prove (i). In \cref{rslt:functoriality-of-kerncat-pushforward} we constructed a 2-functor $f_!\colon \cat K_{\D,X} \to \cat K_{\D,S}$ sending
\begin{align}
	Q \in \Fun_{\cat K_{\D,X}}(Y, X) = \D(Y) \quad \mapsto \quad i_! Q \in \Fun_{\cat K_{\D,S}}(Y, X) = \D(X \times_S Y),
\end{align}
where $i\colon Y \to X \times_S Y$ is the map given by $(f, \id_Y)$. Since 2-functors preserve adjunctions, $i_! Q$ is a left adjoint morphism $Y \to X$ in $\cat K_{\D,S}$ and its right adjoint is given by $i_! \DSuave_g(Q)$ (viewed as a morphism $X \to Y$). One easily computes $P \comp i_! Q = g^* P \comp Q$ in $\cat K_{\D,S}$ (cf.\ the proof of \cite[Proposition~7.11]{Mann.2022b}) and similarly for the right adjoints. Thus (i) follows from the fact that adjoint morphisms are stable under composition.

To prove (ii) we argue similarly: $Q$ is a left adjoint morphism $X \to Y$ in $\cat K_{\D,X}$, hence $i_! Q$ is a left adjoint morphism $X \to Y$ in $\cat K_{\D,S}$. Moreover, one easily computes $P \comp i_! Q = g_!(Q \tensor P)$ as a morphism $X \to S$ in $\cat K_{\D,S}$. We again conclude by the observation that adjoints are stable under composition.
\end{proof}

The next result shows that suave and prim objects can be detected on a suave and prim $*$-cover of the source, respectively:

\begin{lemma} \label{rslt:suave-prim-is-local-on-source}
Let $\D$ be a 3-functor formalism on some geometric setup $(\cat C, E)$, let $f\colon X \to S$ be a map in $E$, let $P \in \D(X)$ and let $\cat U \subseteq (\cat C_E)_{/X}$ be a sieve along which $\D^*$ descends universally.
\begin{lemenum}
	\item \label{rslt:suave-is-local-on-source} Suppose that $\D(X)$ and $\D(S)$ admit $\cat U$-indexed colimits and $f_!(P\tensor \blank) \colon \D(X) \to \D(S)$ preserves them. If $\cat U$ is generated by maps $g\colon U \to X$ such that $\one$ is $g$-suave and $g^*P$ is $fg$-suave then $P$ is $f$-suave and
	\begin{align}
		\DSuave_f(P) = \varinjlim_{g \in \cat U} g_!\DSuave_{fg}(g^* P).
	\end{align}

	\item \label{rslt:prim-is-local-on-source} Suppose that $\D(X)$ and $\D(S)$ admit $\cat U^\op$-indexed limits and $f_!(P\tensor \blank) \colon \D(X) \to \D(S)$ preserves them. If $\cat U$ is generated by maps $g\colon U \to X$ such that $\one$ is $g$-prim and $g^*P$ is $fg$-prim then $P$ is $f$-prim and
	\begin{align}
		\DPrim_f(P) = \varprojlim_{g\in \cat U^\op}g_!\DPrim_{fg}(g^*P).
	\end{align}
\end{lemenum}
\end{lemma}
\begin{proof}
To prove (i) we note that by \cref{rslt:*-descent-implies-limit-in-kerncat} we have $X = \varinjlim_{U \in \cat U} U$ in $\cat K_{\D,S}$. Pick a generating family $(g_i\colon U_i \to X)_{i\in I}$ of $\cat U$ as in (i). Then by \cref{rslt:descent-data-for-generated-sieve} there is a cofinal functor $\bbDelta_I^\op \to \cat U$ whose image consists of maps of the form $U_{i_\bullet} \coloneqq U_{i_1} \times_X \dots \times_X U_{i_n} \to X$. By cofinality, colimits over $\cat U$ and over $\bbDelta_I^\op$ coincide and in particular we have $X = \varinjlim_{i_\bullet \in \bbDelta_I^\op} U_{i_\bullet}$. We now apply the ($(\co,\op)$-dual of the) criterion of \cref{rslt:limits-and-adjunctions}. For condition (a) of \loccit{} we need to show that for every $i_\bullet \in \bbDelta_I^\op$ with associated map $g_{i_\bullet}\colon U_{i_\bullet} \to X$, the object $g_{i_\bullet}^* P$ is $f g_{i_\bullet}$-suave. This follows easily by induction on the length $n$ of $i_\bullet$, the case $n = 1$ being true by assumption and the case $n > 1$ following from the fact that the projection $h\colon U_{i_1,\dots,i_n} \to U_{i_1,\dots,i_{n-1}}$ has the property that $\one$ is $h$-suave (by the case of $n = 1$ and \cref{rslt:suave-and-prim-obj-is-local-on-target}) and the stability of suaveness under suave pullback (see \cref{rslt:composition-of-suave-and-prim-objects}). For condition (b) of \cref{rslt:limits-and-adjunctions} we have to check that for every $Z \in \{X,S\}$ and every $Q \in \D(Z\times_SX)$ the colimit $\varinjlim_{g\in \cat U} Q\circ g_!\DSuave_{fg}(g^*P)$ exists in $\D(Z)$, and that for $Z = X$ the functor $f_!(P\tensor \blank) \colon \D(X)\to \D(S)$ preserves this colimit. But this is ensured by the assumptions. This proves (i).

For (ii), we argue in the same way as above using the $\op$-dual version of \cref{rslt:limits-and-adjunctions}.
\end{proof}

\begin{remarks}
We put ourselves in the setting of \cref{rslt:suave-prim-is-local-on-source}.
\begin{remarksenum}
	\item In practice the 3-functor formalism $\D$ is usually compatible with all small colimits (see \cref{def:3ff-compatible-with-lim-and-colim}), e.g.\ this is the case if $\D$ is presentable. Then \cref{rslt:suave-is-local-on-source} simply says that $P$ is suave over $S$ as soon as it is so on a universal $\D^*$-cover of $X$ in $\cat C_E$ by $\D$-suave maps (i.e.\ maps along which $\one$ is suave).

	\item The previous remark also applies dually to prim objects if $\D$ is compatible with enough limits. However, this is rarely the case, so usually primness is \emph{not} a local condition on the source. The next remark provides an important exception to this rule.

	\item \label{rmk:how-to-apply-locality-of-primness-on-source} Suppose $\D$ is stable and the sieve $\cat U$ on $X$ is generated by a finite family of monomorphisms $(U_i \injto X)_{i\in I}$ along which $\one$ is prim. Then $P$ is prim over $S$ as soon as it is so after pullback to each $U_i$. Indeed, by (the proof of) \cref{rslt:descent-data-for-generated-sieve-by-monomorphisms} there is a cofinal map $P_I^\op \to \cat U$, where $P_I$ is the category of non-empty subsets of $I$; thus a limit over $\cat U^\op$ is the same as a limit over the \emph{finite} category $P_I$ and hence preserved by $f_! (P \tensor \blank)$.
\end{remarksenum}
\end{remarks}

We have now shown that suaveness and primness is local on source and target. The next result concerns itself with the stability of suave and prim objects under limits and colimits:

\begin{lemma} \label{rslt:suave-prim-objects-stable-under-colim}
Let $\D$ be a 3-functor formalism on some geometric setup $(\cat C, E)$ and let $f\colon X \to S$ be a map in $E$.
\begin{lemenum}
	\item $f$-suave objects in $\D(X)$ are stable under retracts.
	\item Let $I$ be a category such that $\D$ is compatible with $I$-indexed colimits and $I^\op$-indexed limits. Then $f$-suave objects are stable under $I$-indexed colimits and $I^\op$-indexed limits.
\end{lemenum}
The same is true for prim objects.
\end{lemma}
\begin{proof}
This is a special case of \cref{rslt:Fun-L-stable-under-retracts,rslt:Fun-L-stable-under-colim-if-compatible-2-cat}.
\end{proof}

\begin{corollary}\label{rslt:suave-prim-stable-under-retracts}
Let $\D$ be a stable 3-functor formalism on some geometric setup $(\cat C, E)$. Then suave and prim objects are stable under retracts, fibers and cofibers.
\end{corollary}
\begin{proof}
This is a special case of \cref{rslt:suave-prim-objects-stable-under-colim}.
\end{proof}

The next result helps to detect suave objects when one has a sufficiently large class of prim objects (or vice-versa). It is one of the main tools we have for detecting suave objects.

\begin{lemma} \label{rslt:conservativity-criterion-for-suave-and-prim-objects}
Let $\D$ be a 3-functor formalism on some geometric setup $(\cat C, E)$, let $f\colon X \to S$ be a map in $E$ and let $(Q_i)_{i\in I}$ be a family of objects in $\D(X)$ such that:
\begin{itemize}
	\item[$(*)$] All $Q_i$ are $f$-prim and the family of functors
	\begin{align}
		\pi_{2!} (\pi_1^* \DPrim_f(Q_i) \tensor \blank)\colon \D(X \times_S X) \to \D(X)
	\end{align}
	is conservative.
\end{itemize}
Then an object $P \in \D(X)$ is $f$-suave if and only if $f_!(Q_i \tensor P)$ is dualizable for all $i$. The same result is true with \enquote{suave} and \enquote{prim} swapped everywhere.
\end{lemma}
\begin{proof}
This is an immediate application of \cref{rslt:conservativity-criterion-for-adjunction} in $\cat K_{\D,S}$, where the roles of $X$, $Y$ and $W_i$ are played by $S$, $X$ and $S$, and the roles of $f$, $\ell_i$ and $r_i$ are played by $P$, $Q_i$ and $\DPrim_f(Q_i)$. Now it only remains to observe that $f_!(Q_i \tensor P)$ is dualizable if and only if it is a left adjoint morphism $S \to S$ in $\cat K_{\D,S}$ (see \cref{rslt:suave-and-prim-for-identity-is-dualizable}).
\end{proof}


\begin{corollary} \label{rslt:prim-generation-criterion-for-suave-objects}
Let $\D$ be a 6-functor formalism on some geometric setup $(\cat C, E)$, let $f\colon X \to S$ be a map in $E$ and let $(Q_i)_{i\in I}$ be a family of objects in $\D(X)$. Assume that:
\begin{itemize}
	\item[$(*)$] All the $Q_i$ are $f$-prim and $\D(X\times_SX)$ is generated by the $\pi_1^*Q_i \tensor \pi_2^*Q_j$. 
\end{itemize}
Then an object $P \in \D(X)$ is $f$-suave if and only if $f_* \iHom(Q_i, P)$ is dualizable for all $Q_i$.
\end{corollary}
\begin{proof}
Note that the composite of $\pi_1^*\DPrim_f(Q_i) \colon X\times_SX \to X$ with $\DPrim_f(Q_j) \colon X\to S$ in $\cat K_{\D,S}$ is $\pi_1^*\DPrim_f(Q_i) \tensor \pi_2^*\DPrim_f(Q_j)$. By assumption ($*$), the induced functors
\begin{align}
(f\times f)_! (\pi_1^*\DPrim_f(Q_i) \tensor \pi_2^*\DPrim_f(Q_j) \tensor \blank) \cong (f\times f)_*\iHom(\pi_1^*Q_i \tensor \pi_2^*Q_j, \blank),
\end{align}
are jointly conservative; for the identification, see \cref{rslt:prim-obj-induce-twist-of-shriek-functors} below. But then the functors $\pi_{2!}(\pi_1\DPrim_f(Q_i) \tensor -)$ are jointly conservative. Thus, the criterion for $f$-suaveness of $P$ is an immediate consequence of \cref{rslt:conservativity-criterion-for-suave-and-prim-objects} using that $f_! (Q \tensor \blank) = f_* \iHom(\DPrim_f(Q), \blank)$.
\end{proof}

\begin{remark}
By \cref{rslt:relation-between-suave-and-prim-duality} below, in the setting of \cref{rslt:prim-generation-criterion-for-suave-objects} the suave dual of an $f$-suave object is characterized completely by the prim duality for the objects $Q_i$.
\end{remark}

We have now discussed plenty of stability results that help detect suave and prim objects. We will now turn to properties that suave and prim maps enjoy. We start with the observation that $f$-suave (resp.\ $f$-prim) objects provide a tight relation between $f^*$ and $f^!$ (resp.\ $f_!$ and $f_*$):

\begin{lemma} \label{rslt:suave-and-prim-obj-induce-twist-of-shriek-functors}
Let $\D$ be a 6-functor formalism on some geometric setup $(\cat C, E)$, let $f\colon X \to S$ be a map in $E$ and let $P \in \D(X)$.
\begin{lemenum}
	\item \label{rslt:suave-obj-induce-twist-of-shriek-functors} If $P$ is $f$-suave then for all $M \in \D(S)$ the natural map
	\begin{align}
		\DSuave_f(P) \tensor f^* M \isoto \iHom(P, f^! M)
	\end{align}
	is an isomorphism.

	\item \label{rslt:prim-obj-induce-twist-of-shriek-functors} If $P$ is $f$-prim then for all $M \in \D(X)$ the natural map
	\begin{align}
		f_!(\DPrim_f(P) \tensor M) \isoto f_* \iHom(P, M)
	\end{align}
	is an isomorphism.
\end{lemenum}
\end{lemma}
\begin{proof}
In the language of \cref{rslt:pointwise-criterion-for-adjunction} (cf.\ the proofs of \cref{rslt:criterion-for-suave-object,rslt:criterion-for-prim-object}) the maps in (i) and (ii) are exactly the natural maps $G_S(\id_S) \comp M \to G_X(M)$ and hence isomorphisms.
\end{proof}

A consequence of \cref{rslt:suave-and-prim-obj-induce-twist-of-shriek-functors} is that there is a close relation between suave and prim objects and compact objects:

\begin{lemma} \label{rslt:relation-between-suave-prim-and-compact}
Let $\D$ be a 3-functor formalism on some geometric setup $(\cat C, E)$ and let $f\colon X \to S$ be a map in $E$. Assume that $\D$ is compatible with small colimits.
\begin{lemenum}
\item If $\Delta_{f!} \one \in \D(X \times_S X)$ is compact, then every $f$-suave object in $\D(X)$ is compact.

\item If $\one \in \D(S)$ is compact, then every $f$-prim object in $\D(X)$ is compact.
\end{lemenum}
\end{lemma}
\begin{proof}
Part (i) follows directly from \cref{rslt:criterion-for-suave-object} and part (ii) follows directly from \cref{rslt:criterion-for-prim-object}. See \cite[Propositions~6.15,~6.16]{Scholze:Six-Functor-Formalism} for details.
\end{proof}

Another consequence of \cref{rslt:prim-obj-induce-twist-of-shriek-functors} is the following tight relation between suave duality and prim duality:

\begin{lemma} \label{rslt:relation-between-suave-and-prim-duality}
Let $\D$ be a 6-functor formalism on some geometric setup $(\cat C, E)$, let $f\colon X \to S$ be a map in $E$, let $P \in \D(X)$ be an $f$-suave object and $Q \in \D(X)$ an $f$-prim object. Then there is a natural isomorphism
\begin{align}
	f_* \iHom(Q, \DSuave_f(P)) = f_* \iHom(\DPrim_f(Q), P)^\vee
\end{align}
\end{lemma}
\begin{proof}
Using \cref{rslt:prim-obj-induce-twist-of-shriek-functors} we have
\begin{align}
	f_* \iHom(Q, \DSuave_f(P)) &= f_!(\DPrim_f(Q) \tensor \DSuave_f(P)) = \DPrim_f(Q) \comp \DSuave_f(P),\\
	f_* \iHom(\DPrim_f(Q), P)^\vee &= (f_!(Q \tensor P))^\vee = (Q \comp P)^\vee
\end{align}
as morphisms $S \to S$ in $\cat K_{\D,S}$, where we view $P$ and $Q$ as right adjoint morphisms $P\colon S \to X$ and $Q\colon X \to S$. Thus the observation follows from the fact that passing to adjoints is functorial with respect to composition (see \cref{rslt:passing-to-adjoints}) and that the adjoint of a right adjoint morphism $S \to S$ is given by its internal dual in $\D(S)$.
\end{proof}

There is also a very general base-change result for pushforward and pullback functors in the presence of suave and prim objects. It is rather technical to formulate in full generality, so we refer the reader to \cref{rslt:suave-and-prim-base-change} for the special case where the suave or prim object is the unit; see \cref{rslt:generalizations-of-suave-prim-base-change} for some generalizations.

\subsection{Suave and prim maps} \label{sec:kerncat.suave-prim-map}

In the previous subsection we defined and studied suave and prim objects in a 3-functor formalism. We now turn our focus to the very related notion of suave and prim maps and their behavior with respect to suave and prim objects. Most of the results in this subsection are just special cases of results from the previous subsection, but they are easier to read and often enough in practice, so we state them anyway.

Let us start with the definition of suave and prim maps. Recall that for every 3-functor formalism $\D$ on a geometric setup $(\cat C, E)$, there is a natural functor $\cat C \to \cat K_\D$ (see the remarks after \cref{rslt:functors-from-and-to-kercat}). In particular every map $f\colon Y \to X$ in $\cat C$ induces a map $Y \to X$ in $\cat K_{\D,X}$, which by abuse of notation we will often also denote by $f$. This map is represented by the object $\one \in \D(Y) = \Fun_{\cat K_{\D,X}}(Y, X)$.

\begin{definition} \label{def:suave-and-prim-map}
Let $\D$ be a 3-functor formalism on a geometric setup $(\cat C, E)$ and let $f\colon Y \to X$ be a morphism in $E$.
\begin{defenum}
	\item We say that $f$ is \emph{$\D$-suave} if it represents a left adjoint morphism in $\cat K_{\D,X}$ (equivalently, if $\one \in \D(Y)$ is $f$-suave). In this case we denote by
	\begin{align}
		\omega_f \coloneqq \DSuave_f(\one) \in \D(Y)
	\end{align}
	the associated right adjoint morphism and call it the \emph{dualizing complex of $f$}. We say that $f$ is \emph{$\D$-smooth} if it is $\D$-suave and $\omega_f$ is invertible.
	
	\item We say that $f$ is \emph{$\D$-prim} if it represents a right adjoint morphism in $\cat K_{\D,X}$ (equivalently, if $\one \in \D(Y)$ is $f$-prim). In this case we denote by
	\begin{align}
		\delta_f \coloneqq \DPrim_f(\one) \in \D(Y)
	\end{align}
	the associated right adjoint morphism and call it the \emph{codualizing complex of $f$}.
\end{defenum}
If $\D$ is clear from context we will sometimes drop it from the notation.
\end{definition}

\begin{remark}\label{rmk:cohomologically-smooth-maps}
Our notion of $\D$-smooth maps usually appears as \enquote{cohomologically smooth maps} in the literature. One common definition of a map $f\colon Y \to X$ to be cohomologically smooth is to require that $f^! \one$ is invertible, the map $f^* \tensor f^! \one \isoto f^!$ is an isomorphism of functors $\D(X) \to \D(Y)$, and the same is true after every base-change of $f$. By the results in this subsection, these properties are also satisfied by $\D$-smooth maps. Moreover, if $f$ is $\D$-smooth then the formation of $f^!$ commutes with pullback along $f$, which is an important condition for the notion of smoothness to be well-behaved (in the definitions in the literature, this condition is often automatic from the other conditions and hence omitted). In fact, the compatibility of $f^!$ with pullback along $f$ is already enough to prove $\D$-suaveness of $f$ (see \cref{rslt:criterion-for-suave-map}). To our knowledge, this was first observed in \cite[Proposition~IV.2.33]{Fargues-Scholze:Geometrization}.
\end{remark}

\begin{remark}
The notion of $\D$-prim maps seems not to be studied much in the classical literature. It roughly says that $f_!$ and $f_*$ are the same up to a twist by $\delta_f$ (see \cref{rslt:prim-maps-induce-twist-of-shriek-functors}). One occurrence known to us is in the case of local systems on anima, see \cref{ex:primness-of-anima}. $\D$-prim maps also occur frequently in the recent theory of analytic stacks by Clausen--Scholze, for example in the context of de Rham-stacks; see for example \cite{Carmago:Analytic-de-Rham-Stack} (where prim maps are called co-smooth maps) or \cite{Scholze:Real-Local-Langlands}. Of course the more restrictive notion of $\D$-proper maps, which is studied in the next subsection, has plenty of examples.
\end{remark}

The criteria for suave and prim objects immediately specialize to criteria for suave and prim maps, which read as follows:

\begin{lemma} \label{rslt:criterion-for-suave-map}
Let $\D$ be a 6-functor formalism on some geometric setup $(\cat C, E)$ and let $f\colon Y \to X$ be a map in $E$. Then $f$ is $\D$-suave if and only if the formation of $f^! \one$ commutes with base-change along $f$, i.e.\ the natural map $\pi_1^* f^! \one \isoto \pi_2^! \one$ is an isomorphism, where $\pi_1, \pi_2\colon Y \times_X Y \to Y$ are the projections; if this is the case then
\begin{align}
	\omega_f = f^! \one.
\end{align}
In fact, it is enough to show that the isomorphism $\pi_1^* f^! \one \isoto \pi_2^! \one$ holds after applying $\Hom(\one, \Delta^!(\blank))$, where $\Delta\colon Y \to Y \times_X Y$ is the diagonal.
\end{lemma}
\begin{proof}
Apply \cref{rslt:criterion-for-suave-object} to $P = \one$.
\end{proof}

\begin{lemma} \label{rslt:criterion-for-prim-map}
Let $\D$ be a 6-functor formalism on some geometric setup $(\cat C, E)$ and let $f\colon Y\to X$ be a map in $E$. Then $f$ is $\D$-prim if and only if the natural map $f_!\pi_{2*}\Delta_!\one \isoto f_*\one$ is an isomorphism, where $\pi_2\colon Y \times_X Y \to Y$ is the projection and $\Delta\colon Y \to Y \times_X Y$ is the diagonal; if this is the case, then
\begin{align}
	\delta_f = \pi_{2*}\Delta_!\one.
\end{align}
In fact, it is enough to show that the isomorphism $f_!\pi_{2*}\Delta_!\one \isoto f_*\one$ holds after applying $\Hom(\one,\blank)$.
\end{lemma}
\begin{proof}
Apply \cref{rslt:criterion-for-prim-object} to $P = \one$. By unraveling the definitions, we observe that the map $f_!\pi_{2*}\Delta_!\one \to f_*\pi_{1!}\Delta_!\one = f_*$ is induced by the natural map $f_!\pi_{2*} \to f_*\pi_{1!}$, which arises as the right mate of the base-change isomorphism $f^*f_! \isoto \pi_{1!}\pi_2^*$.
\end{proof}

By dualizing the above criteria, one can obtain additional criteria for suave and prim maps which make use of the \emph{left} adjoint of $f^*$ instead of its right adjoint (this also works for suave and prim objects instead of maps, of course). Let us spell this out explicitly in the case of suave maps:

\begin{lemma} \label{rslt:dual-criterion-for-suave-map}
Let $\D$ be a 3-functor formalism on some geometric setup $(\cat C, E)$ and let $f\colon Y \to X$ be a map in $E$. Then $f$ is $\D$-suave if and only if $f^*$ and $\pi_2^*$ have left adjoints $f_{\natural}$ and $\pi_{2\natural}$ and the natural map $f_{\natural}\one \isoto f_!\pi_{2\natural}\Delta_!\one$ is an isomorphism; here $\pi_1, \pi_2\colon Y \times_X Y \to Y$ are the projections and $\Delta\colon Y \to Y \times_X Y$ is the diagonal. If $f$ is $\D$-suave then
\begin{align}
	\omega_f = \pi_{2\natural} \Delta_! \one.
\end{align}
\end{lemma}
\begin{proof}
The proof is completely dual to the proof of \cref{rslt:criterion-for-prim-map}, using a dual version of \cref{rslt:criterion-for-prim-object}, which ultimately boils down to (the $\co$-dual version of) \cref{rslt:pointwise-criterion-for-adjunction}. As in the proof of \cref{rslt:criterion-for-prim-map}, we observe that the map $f_{\natural}\one = f_{\natural}\pi_{1!}\Delta_!\one \to f_!\pi_{2\natural}\Delta_!\one$ is induced by the natural map $f_{\natural}\pi_{1!} \to f_!\pi_{2\natural}$, which arises as the left mate of the base-change isomorphism $\pi_{1!}\pi_2^* \isoto f^*f_!$.
\end{proof}

Next we discuss locality of suave and prim maps on source and target. The following results are important special cases of similar results from the previous subsection:

\begin{lemma}\label{rslt:suave-and-prim-map-is-local-on-target}
Let $\D$ be a 3-functor formalism on some geometric setup $(\cat C, E)$ and let $f\colon Y \to X$ be a map in $E$. Then $f$ is $\D$-suave (resp.\ $\D$-prim) if it is so on a universal $\D^*$-cover of $X$.
\end{lemma}
\begin{proof}
This is a special case of \cref{rslt:suave-and-prim-obj-is-local-on-target} for $P = \one$.
\end{proof}

\begin{lemma}
Let $\D$ be a 3-functor formalism on some geometric setup $(\cat C, E)$ and let $f\colon Y \to X$ be a map in $E$.
\begin{lemenum}
\item \label{rslt:suave-map-is-local-on-source} Assume that $\D$ is compatible with small colimits. Then $f$ is $\D$-suave as soon as it is so on a universal $\D^*$-covering family of $Y$ by $\D$-suave maps.

\item \label{rslt:prim-map-is-local-on-source} The map $f$ is $\D$-prim if it is so on a universal $\D^*$-covering family of $Y$ by $\D$-suave maps such that for the generated sieve $\cat U$ (in $\cat C_E$) both $\D(Y)$ and $\D(X)$ admit $\cat U^\op$-indexed limits and $f_!$ preserves them.
\end{lemenum}
\end{lemma}
\begin{proof}
This is a special case of \cref{rslt:suave-prim-is-local-on-source} for $P = \one$.
\end{proof}

The previous two results show that suaveness and primness can be checked on a cover of source and target, but they do not say anything about the stability of suaveness and primness under passing to such a cover. This is solved by the next result:

\begin{lemma} \label{rslt:stability-of-suave-and-prim-maps}
Let $\D$ be a 3-functor formalism on some geometric setup $(\cat C, E)$.
\begin{lemenum}
	\item $\D$-suave maps in $E$ are stable under composition and every base-change.

	\item \label{rslt:cancellativeness-of-suave-and-prim-maps} Let $f\colon Y \to X$ and $g\colon Z \to Y$ be maps in $E$. If $\Delta_f\colon Y \to Y \times_X Y$ and $fg$ are $\D$-suave then $g$ is $\D$-suave.
\end{lemenum}
The same is true for prim maps.
\end{lemma}
\begin{proof}
Stability under base-change follows from \cref{rslt:suave-and-prim-obj-is-local-on-target} and stability under composition follows from \cref{rslt:pullback-of-suave-and-prim-objects}, proving (i). Part (ii) is a special case of \cref{rslt:suave-and-prim-along-g-for-suave-and-prim-Delta-f} below, but also follows from (i) by a standard argument: We factor $g$ as the composition $Z \to Z \times_X Y \to Y$. The second map is a base-change of $fg$ and hence $\D$-suave. The first map is a base-change of $\Delta_f$ along the map $(g,\id)\colon Z \times_X Y \to Y \times_X Y$ and hence also $\D$-suave.
\end{proof}

The next result allows one to detect suave and prim maps in the case where the exceptional pushforward is conservative. We wrote it in the version that we believe is the most common one in applications, but many different versions are true; see \cref{rslt:conservativity-criterion-for-suave-and-prim-objects} for the general result.

\begin{lemma} \label{rslt:conservative-prim-criterion-for-suave-map}
Let $\D$ be a 6-functor formalism on some geometric setup $(\cat C, E)$ and let $f\colon Y \to X$ be a $\D$-prim map such that $\pi_{2*}\colon \D(Y \times_X Y) \to \D(Y)$ is conservative. Then an object $P \in \D(Y)$ is $f$-suave if and only if $f_! P \in \D(X)$ is dualizable.
\end{lemma}
\begin{proof}
Apply \cref{rslt:conservativity-criterion-for-suave-and-prim-objects} to $Q_i = \one$ and use \cref{rslt:prim-obj-induce-twist-of-shriek-functors}.
\end{proof}

So far we have discussed various ways to detect suave and prim maps. Next up we investigate the nice properties that these maps enjoy. We start with the observation that suave and prim maps provide a close relation between $!$- and $*$-functors:

\begin{corollary} \label{rslt:suave-and-prim-maps-induce-twist-of-shriek-functors}
Let $\D$ be a 6-functor formalism on some geometric setup $(\cat C, E)$ and let $f\colon Y \to X$ be a map in $E$.
\begin{corenum}
	\item \label{rslt:suave-maps-induce-twist-of-shriek-functors} If $f$ is $\D$-suave then the natural maps
	\begin{align}
		\omega_f \tensor f^* \isoto f^!, \qquad f^* \isoto \iHom(\omega_f, f^!)
	\end{align}
	are isomorphisms of functors $\D(X) \to \D(Y)$.

	\item \label{rslt:prim-maps-induce-twist-of-shriek-functors} If $f$ is $\D$-prim then the natural maps
	\begin{align}
		f_!(\delta_f \tensor \blank) \isoto f_*, \qquad f_! \isoto f_* \iHom(\delta_f, \blank)
	\end{align}
	are isomorphisms of functors $\D(Y) \to \D(X)$.
\end{corenum}
\end{corollary}
\begin{proof}
Apply \cref{rslt:suave-and-prim-obj-induce-twist-of-shriek-functors} to $P = \one$ and $P = \omega_f$ or $P = \delta_f$, respectively.
\end{proof}

\begin{remark} \label{rslt:dualizable-dualizing-complex-is-invertible}
Suppose $f\colon Y \to X$ is a $\D$-suave map for some 6-functor formalism and assume that $\omega_f \in \D(Y)$ is dualizable. Then it is automatically invertible, so that $f$ is $\D$-smooth. Namely, by \cref{rslt:suave-maps-induce-twist-of-shriek-functors} we have $\one = \iHom(\omega_f, \omega_f)$ and by dualizability we get $\iHom(\omega_f, \omega_f) = \omega_f^\vee \tensor \omega_f$. This observation was pointed out to us by David Hansen.
\end{remark}

Another important property of suave and prim maps is that they allow general base-change results for push and pull functors (e.g.\ \enquote{suave base-change}). Namely, a simple principle in $\cat K_\D$ results in the following list of base-changes:

\begin{lemma} \label{rslt:suave-and-prim-base-change}
Let $\D$ be a 6-functor formalism on some geometric setup $(\cat C, E)$ and let
\begin{equation}\begin{tikzcd}
	Y' \arrow[r,"g'"] \arrow[d,"f'"'] & Y \arrow[d,"f"]\\
	X' \arrow[r,"g"'] & X 
\end{tikzcd}\end{equation}
be a cartesian diagram in $\cat C_E$. 
\begin{lemenum}
	\item\label{rslt:suave-base-change} Assume that $g$ is $\D$-suave. Then the natural morphisms
	\begin{align}
		g^*f_* &\isoto f'_*g'^*, & f'_!g'^! &\isoto g^!f_!, & f'^*g^! &\isoto g'^!f^*, &\text{and}&& g'^*f^! &\isoto f'^!g^*
	\end{align}
	are isomorphisms of functors.

	\item\label{rslt:prim-base-change} Assume that $g$ is $\D$-prim. Then the natural morphisms
	\begin{align}
		f^*g_* &\isoto g'_*f'^*, & g'_!f'^! &\isoto f^!g_!, & g_!f'_* &\isoto f_*g'_!, & \text{and} && f_!g'_* &\isoto g_*f'_!
	\end{align}
	are isomorphisms of functors.
\end{lemenum}
\end{lemma}
\begin{proof}
In each of the cases, it is not hard to construct \emph{some} isomorphism between the two functors. For example, if $g$ is $\D$-suave then the isomorphism $g^* f_* \isom f'_* g'^*$ follows by writing $g^* = \iHom(\omega_g, g^!)$ (and similarly for $g'$) and using the isomorphism $g^! f_* \isom f'_* g'^!$ encoded in the 6-functor formalism. However, with this strategy it is not clear why the constructed isomorphism is the natural one, especially since we made use of several non-trivial isomorphisms encoded by $\D$.

We now present a conceptual proof that the isomorphism is indeed the correct one, using the 2-category of kernels. We start with the proof of (i), so assume that $g$ is $\D$-suave. Note that the suaveness of $g$ implies that $g^*$ has a left adjoint $g_\natural\colon \D(X') \to \D(X)$ (this follows easily from the fact that $\Psi_{\D,X}\colon \cat K_{\D,X} \to \Cat$ preserves adjunctions, or alternatively by \cref{rslt:suave-maps-induce-twist-of-shriek-functors}). By \cref{rslt:stability-of-suave-and-prim-maps} also $g'$ is $\D$-suave and we get a similar functor $g'_\natural$. By passing to left adjoints in the first of the claimed isomorphisms of functors, we reduce to showing that the natural map $g'_\natural f'^* \isoto f^* g_\natural$ is an isomorphism of functors $\D(X') \to \D(Y)$. The advantage of this perspective is that $g_\natural$ and $g'_\natural$ can be represented in the category of kernels (i.e.\ lie in the image of $\Psi_{\D,X}$), so we reduce to showing a similar statement in $\cat K_{\D,X}$.

Recall from \cref{rmk:kerncat-symmetric-monoidal-2-cat} the 2-functor $\cat K_{\D,X} \times \cat K_{\D,X} \to \cat K_{\D,X}$ given by $(W,Z)\mapsto W\times_XZ$. Precomposing with the 2-functor $\cat K_{\D,X}\times [1] \to \cat K_{\D,X}\times \cat K_{\D,X}$ induced by $\one_Y\colon X\to Y$ yields a natural transformation $\sigma\colon \id_{\cat K_{\D,X}} \to \blank\times_XY$. Thus, for any morphism $P\colon W\to Z$ in $\cat K_{\D,X}$ we obtain a commutative diagram
\begin{equation}
\begin{tikzcd}
W \ar[d,"P"'] \ar[r,"\sigma_W"] & W\times_XY \ar[d,"f^*P"] \\
Z \ar[r,"\sigma_Z"'] \ar[ur,Rightarrow,shorten=4mm, sloped,"\sim", "\sigma_P"'] & Z\times_XY
\end{tikzcd}
\end{equation}
in $\cat K_{\D,X}$. Now, let $P\colon W\rightleftarrows Z \noloc Q$ be an adjunction in $\cat K_{\D,X}$. By 2-functoriality we obtain the adjunction $f^*P \colon W\times_XY \rightleftarrows Z\times_XY \noloc f^*Q$ and then an easy direct computation (the lazy reader is referred to \cite[Proposition~2.5]{Kelly-Street.1974}) shows that the left mate
\begin{align} \label{eq:suave-prim-base-change-left-mate}
	f^*P\comp \sigma_W \isoto \sigma_Z\comp P
\end{align}
of $\sigma_Q\colon \sigma_W\comp Q\isoto f^*Q\comp \sigma_Z$ is an isomorphism with inverse $\sigma_P$. Applying \eqref{eq:suave-prim-base-change-left-mate} to the adjunction $\omega_g\colon X'\rightleftarrows X\noloc \one_{X'}$ in $\cat K_{\D,X}$ and then applying the 2-functor $\Psi_{\D,X}\colon \cat K_{\D,X}\to \Cat$ from \cref{rslt:functors-from-and-to-kercat} shows that $g'_{\natural} f'^* \isoto f^*g_{\natural}$ is an isomorphism, as desired. Similarly, applying \eqref{eq:suave-prim-base-change-left-mate} to the adjunction $\one_{X'}\colon X'\rightleftarrows X \noloc \omega_g$ instead yields the isomorphism $f'^*g^!\isoto g'^!f^*$.

Running the argument in $\cat K_{\D,X}^\op$, $\cat K_{\D,X}^{\co}$ and $\cat K_{\D,X}^{\co,\op}$ yields the other six isomorphisms.
\end{proof}

\begin{remark} \label{rslt:generalizations-of-suave-prim-base-change}
Suppose we are in the context of \cref{rslt:suave-and-prim-base-change}.
\begin{remarksenum}
\item The proof of \cref{rslt:suave-and-prim-base-change} works slightly more generally for a \emph{3}-functor formalism $\D$ if one makes the appropriate modifications. For example, if $g$ is $\D$-suave then there is a natural isomorphism $g'_\natural f'^* \isoto f^* g_\natural$ of functors $\D(X') \to \D(Y)$, where $g_\natural\colon \D(X') \to \D(X)$ is the left adjoint of $g^*$ (which exists by suaveness). 
\item One can even generalize further and replace the $\D$-suaveness requirement by an arbitrary $g$-suave sheaf. For example, if $P \in \D(X')$ is $g$-suave then the natural map
\begin{align}
	P \tensor g^* f_* \isoto f'_* (f'^* P \tensor g'^*)
\end{align}
is an isomorphism of functors $\D(Y) \to \D(X')$. The proof is exactly the same as in \cref{rslt:suave-and-prim-base-change}.

\item If $g$ is $\D$-suave (resp.\ $\D$-prim), then for the isomorphisms $g^*f_*\isoto f'_*g'^*$ and $f'^*g^!\isoto g'^!f^*$ (resp.\ $f^*g_*\isoto g'_*f'^*$ and $g_!f'_* \isoto f_*g'_!$) to hold it is not necessary to require $f,f'\in E$---but then one has to slightly modify the proof by replacing $\sigma$ with the natural transformation $\sigma\colon \Psi_{\D,X} \to \Psi_{\D,Y} \comp f^*$ of 2-functors $\cat K_{\D,X} \to \Cat$ (and then the same argument goes through). 
\end{remarksenum}
\end{remark}

Next we discuss the interaction of suave and prim maps with suave and prim objects, specializing some results from the previous subsection. We start with the following preservation of suaveness and primness under suave and prim pullback and pushforward:

\begin{lemma} \label{rslt:suave-pullback-and-prim-pushforward-preservations}
Let $\D$ be a 3-functor formalism on some geometric setup $(\cat C, E)$ and let $f\colon X \to S$ and $g\colon Y \to X$ be maps in $E$.
\begin{lemenum}
	\item If $g$ is $\D$-suave then $g^*\colon \D(X) \to \D(Y)$ preserves suave objects over $S$.
	\item If $g$ is $\D$-prim then $g_!\colon \D(Y) \to \D(X)$ preserves suave objects over $S$.
\end{lemenum}
The same is true with \enquote{suave} and \enquote{prim} swapped.
\end{lemma}
\begin{proof}
This is a special case of \cref{rslt:composition-of-suave-and-prim-objects}.
\end{proof}

As discussed in \cref{rslt:suave-and-prim-for-identity-is-dualizable}, suave and prim objects for the identity on some object are the same as dualizable objects. This relation extends to some more general situations, as the following results show.

\begin{lemma} \label{rslt:suave-and-prim-over-suave-and-prim-map-and-diagonal}
Let $\D$ be a 3-functor formalism on some geometric setup $(\cat C, E)$ and let $f\colon X \to S$ and $g\colon Y \to X$ be maps in $E$.
\begin{lemenum}
	\item Suppose that $f$ is $\D$-suave. Then every $g$-suave object $P$ in $\D(Y)$ is $fg$-suave, and $\DSuave_{fg}(P) = g^*\omega_f\tensor \DSuave_g(P)$.

	\item \label{rslt:suave-and-prim-along-g-for-suave-and-prim-Delta-f} Suppose that $\Delta_f$ is $\D$-suave. Then every $fg$-suave object $P$ in $\D(Y)$ is $g$-suave and $\DSuave_g(P) = \DSuave_{fg}(P) \tensor \omega_a$, where $a$ is the $\D$-suave map $(g,\id_X)\colon Y \to X \times_S Y$.
\end{lemenum}
The same is true for primness instead of suaveness (with $\delta$ in place of $\omega$).
\end{lemma}
\begin{proof}
Part (i) is a special case of \cref{rslt:pullback-of-suave-and-prim-objects}. To prove (ii), suppose that $\Delta_f$ is $\D$-suave and let $P \in \D(Y)$ be an $fg$-suave object. We consider the diagram
\begin{equation}\begin{tikzcd}
	Y \arrow[r,"a"] & X \times_S Y \arrow[r,"\pi_2"] \arrow[d,"\pi_1"] & Y \arrow[d,"fg"]\\
	& X \arrow[r,"f"] & S.
\end{tikzcd}\end{equation}
Here $a = (g, \id_Y)$. By \cref{rslt:suave-and-prim-obj-is-local-on-target} we know that $\pi_2^* P$ is $\pi_1$-suave. Note that $a$ is the base-change of $\Delta_f$ along the map $g\times\id\colon Y \times_S X \to X \times_S X$, hence $a$ is $\D$-suave by \cref{rslt:stability-of-suave-and-prim-maps}. But then \cref{rslt:pullback-of-suave-and-prim-objects} implies that $P = a^* \pi_2^* P$ is $g = \pi_1 a$-suave and
\begin{align}
	\DSuave_g(P) = \DSuave_{fg}(P) \tensor \omega_a,
\end{align}
as desired.
\end{proof}

\begin{corollary} \label{rslt:relation-between-suave-prim-and-dualizable}
Let $\D$ be a 3-functor formalism on some geometric setup $(\cat C, E)$ and let $f\colon X \to S$ be a map in $E$.
\begin{corenum}
\item Suppose $f$ is $\D$-suave. Then every dualizable object $P \in \D(X)$ is $f$-suave and $\DSuave_f(P) = P^\vee \tensor \omega_f$.

\item Suppose $\Delta_f$ is $\D$-suave. Then every $f$-suave object $P \in \D(X)$ is dualizable and $\DSuave_f(P) = P^\vee \tensor \omega_{\Delta_f}^{-1}$.
\end{corenum}
The same is true for primness instead of suaveness (with $\delta$ in place of $\omega$).
\end{corollary}
\begin{proof}
This is the special case of \cref{rslt:suave-and-prim-over-suave-and-prim-map-and-diagonal} with $Y = X$ and $g = \id_X$.
\end{proof}

\subsection{Étale and proper maps} \label{sec:kerncat.etale-proper}

In the previous subsection we introduced suave and prim maps, which provide close relations between $!$- and $*$-functors by realizing one as the twist of the other. There is a special class of suave and prim maps where the twist is canonically trivialized, namely the étale and proper maps, respectively. We will introduce and study these maps now. Let us start with the definition following \cite[Definitions~6.10,~6.12]{Scholze:Six-Functor-Formalism} (recall the notion of truncated maps from \cref{defn:truncated}).

\begin{definition} \label{def:etale-and-proper-maps}
Let $\D$ be a 3-functor formalism on some geometric setup $(\cat C, E)$ and let $f\colon Y \to X$ be a truncated map in $E$.
\begin{defenum}
\item We say that $f$ is \emph{$\D$-étale} if it is $\D$-suave and $\Delta_f$ is $\D$-étale or an isomorphism.

\item We say that $f$ is \emph{$\D$-proper} if it is $\D$-prim and $\Delta_f$ is $\D$-proper or an isomorphism.
\end{defenum}
\end{definition}

\begin{remark}
By definition of truncated maps, if $f$ is $n$-truncated then $\Delta_f$ is $(n-1)$-truncated. Thus, by repeatedly applying the diagonal we eventually end up with an isomorphism, so that \cref{def:etale-and-proper-maps} is sensible.
\end{remark}

The above results on suave and prim maps restrict to similar results about étale and proper maps, which we state here to reference them later:

\begin{lemma}
Let $\D$ be a 3-functor formalism on some geometric setup $(\cat C, E)$.
\begin{lemenum}
	\item \label{rslt:stability-of-etale-and-proper-maps} $\D$-étale maps are stable under composition and base-change. Moreover, if $f$ and $g$ are composable maps in $E$ such that $fg$ and $\Delta_f$ are $\D$-étale, then $g$ is $\D$-étale.

	\item \label{rslt:localness-of-etale-and-proper-maps} $\D$-étaleness of a truncated morphism is $\D^*$-local on the target. If $\D$ is compatible with small colimits then $\D$-étaleness of truncated morphisms can be checked on a universal $\D$-étale $\D^*$-cover (in $\cat C_E$) of the source.

	\item \label{rslt:etale-means-suave=dualizable} If $f\colon Y\to X$ is $\D$-étale, then an object $P\in \D(Y)$ is $f$-suave if and only if it is dualizable, and in this case $\DSuave_f(P) = P^\vee \coloneqq \iHom(P,\one)$.
\end{lemenum}
The same properties hold for $\D$-proper maps, where in the second part of (ii) one has the following statement: 
\begin{lemenum}[label=(\roman*'), ref=\thelemma.(\roman*')]
	\stepcounter{lemenumi}
	\item \label{rslt:localness-of-proper-on-source} A morphism $f$ in $E$ is $\D$-proper if it is so locally on a universal $\D$-proper $\D^*$-cover $\cat U$ on the source such that $f_!$ commutes with $\cat U^{\op}$-indexed limits.
\end{lemenum}
\end{lemma}
\begin{proof}
Part (iii) is a special case of \cref{rslt:suave-and-prim-over-suave-and-prim-map-and-diagonal}, where the identification of the suave dual follows from the formula in the proof.

We prove stability under base-change in (i), so let $f\colon Y\to X$ be a $\D$-étale map, and denote $f'\colon Y'\to X'$ the pullback of $f$ along some map $X'\to X$. We prove the claim by induction on $n$ such that $f$ is $n$-truncated. By \cref{rslt:truncated-pullback} and \cref{rslt:stability-of-suave-and-prim-maps} it follows, respectively, that $f'$ is $n$-truncated and $\D$-suave. Moreover, $\Delta_{f'}$ is $(n-1)$-truncated and a pullback of $\Delta_f$. By the induction hypothesis, $\Delta_{f'}$ is $\D$-étale. Hence, $f'$ is $\D$-étale.

For stability under composition, let $g\colon Z\to Y$ and $f\colon Y\to X$ be $\D$-étale. Writing $\Delta_{fg} = \Delta'_f\comp \Delta_g$, where $\Delta'_f\colon Z\times_YZ\to Z\times_XZ$ is a pullback of $\Delta_f$. Arguing as before, but using \cref{rslt:truncated-composition} and stability under pullback, we see that $fg$ is $\D$-suave and $\Delta_{fg}$ is $\D$-étale, hence $fg$ is $\D$-étale.

Assume now that $fg$ and $\Delta_f$ are $\D$-étale. Then $g$ is $\D$-suave by \cref{rslt:stability-of-suave-and-prim-maps}. If $n\ge -2$ is such that $fg$ and $\Delta_f$ are $n$-truncated, then $g$ is $n$-truncated by \cref{rslt:truncated-composition}. Finally, we write $\Delta_{fg} = \Delta'_{f} \comp \Delta_g$ as above. We have shown that $\Delta_{fg}$ and $\Delta_{\Delta'_{f}}$ are $\D$-étale and $(n-1)$-truncated. It follows by induction on $n$ that $\Delta_g$ is $\D$-étale. Therefore, also $g$ is $\D$-étale, which proves (i).

We now prove (ii). Let $f\colon Y\to X$ be an $n$-truncated morphism in $E$. Assume that $f$ is locally on some universal $\D^*$-cover $\cat U$ of $X$ a $\D$-étale map. We prove by induction on $n$ that $f$ is $\D$-étale, the case $n=-2$ being trivial. \Cref{rslt:suave-and-prim-map-is-local-on-target} shows that $f$ is $\D$-suave. Denoting by $\pi\colon Y\times_XY\to X$ the natural map, we observe that $\Delta_f$ is locally on the universal $\D^*$-cover $\pi^*\cat U$ a $\D$-étale map. By induction on $n$ it follows that $\Delta_f$ is $\D$-étale. Hence $f$ is $\D$-étale.

We now prove the second statement in (ii). Let $f$ be as above. Assume that $f$ is locally on some universal $\D$-étale $\D^*$-cover $\cat U$ of $Y$ a $\D$-étale map. We prove by induction on $n$ that $f$ is $\D$-étale, the case $n=-2$ being trivial. \Cref{rslt:suave-map-is-local-on-source} shows that $f$ is $\D$-suave. For any $g\colon U\to Y$ in $\cat U$ we observe that $\Delta_f\comp g$ factors as the composite of $\Delta_{fg}\colon U\to U\times_XU$ and $g\times g\colon U\times_XU\to Y\times_XY$. $\Delta_{fg}$ is $\D$-étale by assumption and $g\times g$ is $\D$-étale by part (i). Hence $\Delta_f\comp g$ is $\D$-étale. By induction on $n$ we deduce that $\Delta_f$ is $\D$-étale. Therefore, also $f$ is $\D$-étale.
\end{proof}

We refer the reader to \cref{rmk:how-to-apply-locality-of-primness-on-source} for the standard way of how to apply \cref{rslt:localness-of-proper-on-source} in practice. The next result (taken from \cite[Propositions~6.11,~6.13]{Scholze:Six-Functor-Formalism}) shows that étale and proper maps do indeed guarantee an identification of $!$- and $*$-functors. In fact, this also provides an easy way of checking étaleness and properness.

\begin{lemma}
Let $\D$ be a 6-functor formalism on some geometric setup $(\cat C, E)$ and let $f\colon Y \to X$ be a map in $E$.
\begin{lemenum}
	\item \label{rslt:etale-equiv-iso-of-shriek-and-star} Suppose that $\Delta_f$ is $\D$-étale. Then there is a natural map $f^! \to f^*$ of functors $\D(X) \to \D(Y)$ and the following are equivalent:
	\begin{enumerate}[(a)]
		\item $f$ is $\D$-étale.
		\item The map $f^! \one \isoto f^* \one$ is an isomorphism in $\D(X)$.
		\item The map $f^! \isoto f^*$ is an isomorphism of functors $\D(Y) \to \D(X)$.
	\end{enumerate}

	\item \label{rslt:proper-equiv-iso-of-shriek-and-star} Suppose that $\Delta_f$ is $\D$-proper. Then there is a natural map $f_! \to f_*$ of functors $\D(Y) \to \D(X)$ and the following are equivalent:
	\begin{enumerate}[(a)]
		\item $f$ is $\D$-proper.
		\item The map $f_! \one \isoto f_* \one$ is an isomorphism in $\D(X)$.
		\item The map $f_! \isoto f_*$ is an isomorphism of functors $\D(Y) \to \D(X)$.
	\end{enumerate}
\end{lemenum}
\end{lemma}
\begin{proof}
This is proved in \cite[Propositions~6.11,~6.13]{Scholze:Six-Functor-Formalism}. We repeat the argument in a slightly more elaborate fashion.

We prove (i). By induction on the number $n$ such that $f$ is $n$-truncated we can assume that there is a natural isomorphism $\Delta_f^!= \Delta_f^*$. The map $f^!\to f^*$ is defined as the composite
\begin{align}
	f^! = \Delta_f^*\pi_2^*f^! = \Delta_f^!\pi_2^*f^! \to \Delta_f^!\pi_1^!f^* = f^*,
\end{align}
where the intermediate map is induced by the map $\pi_2^*f^!\to \pi_1^!f^*$, which arises as the right mate of the base-change isomorphism $\pi_{1!}\pi_2^* \isoto f^*f_!$. It is clear that (c) implies (b). In order to show that (b) implies (a), it suffices to show that $f$ is $\D$-suave. But note that (b) implies that the natural map $\pi_2^* f^! \one \isoto \pi_1^! \one$ becomes an isomorphism after applying $\Delta_f^!$, hence $\D$-suaveness follows directly from \cref{rslt:criterion-for-suave-map}. Finally, (a) implies (c) as $\pi_2^*f^!\isoto \pi_1^!f^*$ is an isomorphism by \cref{rslt:suave-base-change}.

Let us prove (ii). By induction on the number $n$ such that $f$ is $n$-truncated we can assume that there is a natural isomorphism $\Delta_{f!} = \Delta_{f*}$. The map $f_!\to f_*$ is defined as the composite
\begin{align}
	f_! = f_!\pi_{2*}\Delta_{f*} = f_!\pi_{2*}\Delta_{f!} \to f_*\pi_{1!}\Delta_{f!} = f_*,
\end{align}
where the map is induced by the map $f_!\pi_{2*} \to f_*\pi_{1!}$, which arises as the right mate of the base-change isomorphism $f^*f_! \isoto \pi_{1!}\pi_2^*$. It is clear that (c) implies (b). In order to show that (b) implies (a), it suffices to show that $f$ is $\D$-prim, which is the content of \cref{rslt:criterion-for-prim-map}. Finally, (a) implies (c) as $f_!\pi_{2*}\isoto f_*\pi_{1!}$ is an isomorphism by \cref{rslt:prim-base-change}.
\end{proof}

\begin{remark}
Given a 3-functor formalism $\D$ on some geometric setup $(\cat C, E)$, the classes $\mathit{et} \subseteq E$ and $\mathit{prop} \subseteq E$ of $\D$-étale and $\D$-proper maps are stable under composition and base-change and are right cancellative (cf. \cref{rslt:right-cancellative-fiber-products} for the definition). Hence both $(\cat C, \mathit{et})$ and $(\cat C, \mathit{prop})$ form geometric setups.
\end{remark}

\subsection{Exceptional descent} \label{sec:kerncat.excdescent}

In the extension results in \cref{sec:6ff.extend} we needed $\D^!$-descent in order to localize the construction of shriek functors on the source. With the results in the previous subsections we can now provide some basic criteria for proving $!$-descent in practice. In the following recall the definition of (universal) $\D^!$-covers from \cref{def:D*-and-D!-covers}. As a special case we say that a map $f\colon Y \to X$ is a (universal) $\D^!$-cover if the sieve generated by $f$ is so.

In the case that $f$ is suave, $!$-descent can be reduced to $*$-descent. This idea already appeared in \cite[Proposition~3.16]{Gulotta-Hansen-Weinstein:Enhanced-6FF-on-vStacks} and was used in several sources after. The cleanest and most general formulation is found in \cite[Proposition~6.18]{Scholze:Six-Functor-Formalism}, which we recall here:

\begin{lemma} \label{rslt:!-descent-along-suave-map}
Let $\D$ be a 6-functor formalism on some geometric setup $(\cat C, E)$ and let $f\colon Y \to X$ be a map in $E$ such that $\D(X)$ has all countable limits and colimits. If $f$ is $\D$-suave and $f^*\colon \D(X) \to \D(Y)$ is conservative, then $f$ is a universal $\D^*$-cover and a universal $\D^!$-cover.
\end{lemma}
\begin{proof}
Let us first discuss $\D^!$-descent, for which we need to show that the natural functor $\D^!(X) \isoto \varprojlim_{n\in\bbDelta} \D^!(Y_n)$ is an equivalence. We will employ Lurie's Beck--Chevalley condition \cite[Corollary~4.7.5.3]{HA}. Condition (2) of loc.\ cit.\ is implied by suave base-change (see \cref{rslt:suave-base-change}) and the remaining conditions of loc.\ cit.\ amount to showing that 
\begin{enumerate*}[(a)] 
\item $f^!\colon \D(X) \to \D(Y)$ is conservative and
\item $f^!$ preserves colimits over $f^!$-split simplicial objects in $\D(X)$.
\end{enumerate*}
Part (a) follows immediately from the fact that $f^*$ is conservative and $f^* = \iHom(\omega_f, f^!)$ by \cref{rslt:suave-maps-induce-twist-of-shriek-functors}. For (b) we note that $f^!$ preserves all colimits by the identity $f^! = f^* \tensor \omega_f$ (see \cref{rslt:suave-maps-induce-twist-of-shriek-functors}).

The $\D^*$-descent is proved very similarly, using the dual version of Lurie's Beck--Chevalley condition. One then reduces to showing that $f^*$ is conservative and preserves $f^*$-split totalizations of cosimplicial objects in $\D(X)$. But the former is true by assumption and the latter by $f^* = \iHom(\omega_f, f^!)$.
\end{proof}

The previous result provides a powerful criterion for checking $\D^*$- and $\D^!$-descent along a suave map. There is a similar result for \emph{prim} maps, but under much stronger assumptions. To our knowledge the following criterion first appeared in the proof of \cite[Lemma~3.11]{Hansen-Mann:Mod-p-Stacky-6FF} and is used there to construct the 6-functor formalism for classifying stacks of $p$-adic Lie groups (see \cref{rslt:BG-is-proper} below for a very similar argument). We state it in the abstract form as in \cite[Proposition~6.19]{Scholze:Six-Functor-Formalism}. The criterion makes use of Mathew's notion of descendable algebras (see \cite[Definition~3.18]{Mathew:2016}), which we recall here (in a slightly more general form):

\begin{definition} \label{def:descendability}
Let $\cat C$ be a stable monoidal category. For every object $A \in \cat C$ we denote by $\langle A \rangle \subseteq \cat C$ the smallest full subcategory containing $A$ which is stable under finite (co)limits, retracts and tensor products. We say that $A$ is \emph{descendable} if $\langle A \rangle$ contains the tensor unit.
\end{definition}

\begin{definition}[{cf. \cite[\S3.2]{Mathew:2016}}]
Let $\cat C$ be a category with finite colimits.
\begin{defenum}
	\item A filtered diagram $F\colon I \to \cat C$ is called \emph{ind-constant} if, when viewed as an element in $\Ind(\cat C)$, it lies in the essential image of the Yoneda embedding $\cat C \injto \Ind(\cat C)$. By dualizing, we define the notion of a \emph{pro-constant} cofiltered diagram.

	\item A simplicial object $M_\bullet\colon \bbDelta^\op \to \cat C$ is called \emph{ind-constant} if the associated diagram $\ZZ_{\ge0} \to \cat C$, $n \mapsto \varinjlim_{m \in \bbDelta_{\le n}^\op} M_n$ is ind-constant. We similarly define \emph{pro-constant} cosimplicial objects.
\end{defenum}
\end{definition}

\begin{lemma} \label{rslt:!-descent-along-prim-morphism}
Let $\D$ be a stable 6-functor formalism on some geometric setup $(\cat C, E)$ and let $f\colon Y \to X$ be a map in $E$. Suppose that $f$ is $\D$-prim and $f_*\one \in \D(X)$ is descendable. Then:
\begin{lemenum}
	\item $f$ is a universal $\D^*$-cover and a universal $\D^!$-cover.
	
	\item \label{rslt:!-descent-along-prim-morphism-implies-pro-const} Every $f^!$-split simplicial object in $\D(X)$ is ind-constant. Every $f^*$-split cosimplicial object in $\D(X)$ is pro-constant.

	\item Let $f_\bullet\colon Y_\bullet \to X$ be the Čech nerve of $f$. For every $M \in \D(X)$ the simplicial object $(f_{n!} f_n^! M)_n$ is ind-constant and the cosimplicial object $(f_{n*} f_n^* M)_n$ is pro-constant.
\end{lemenum}
\end{lemma}
\begin{proof}
Let us first prove universal $\D^!$-descent along $f$, for which we argue similarly as in the proof of \cref{rslt:!-descent-along-suave-map}. Using prim base-change in place of suave base-change (both are proved in \cref{rslt:suave-and-prim-base-change}), we see that the problem of showing $\D^!$-descent is reduced to proving the following claim:
\begin{itemize}
	\item[($*$)] $f^!$ is conservative and preserves colimits of $f^!$-split simplicial objects in $\D(X)$.
\end{itemize}
More concretely we need to show that $f^!$ preserves colimits of the form $\varinjlim_{n\in\bbDelta} f_{n!} N^n$ for $N^\bullet \in \varprojlim_{n\in\bbDelta} \D^!(Y_n)$, but these are colimits of $f^!$-split simplicial objects. Now observe that ($*$) follows from the following claim:
\begin{itemize}
	\item[($**$)] The object $f_! f^! \in \Fun(\D(X), \D(X))$ is descendable (for the composition monoidal structure).
\end{itemize}
Let us check that ($**$) implies ($*$). Firstly, it follows easily from ($**$) that the functor $f_! f^!$ is conservative, hence so is $f^!$. Moreover, if $M^\bullet$ is any $f^!$-split simplicial object in $\D(X)$ then it is also $f_! f^!$-split and in particular $f_! f^! M^\bullet$ is ind-constant (see \cite[Example~3.11]{Mathew:2016}). It follows that $\varphi(M^\bullet)$ is ind-constant for every functor $\varphi \in \langle f_! f^! \rangle$ (this follows from the criterion in \cite[Proposition~3.10]{Mathew:2016}). Hence $M^\bullet$ itself is ind-constant and thus its colimit commutes with $f^!$.

It remains to show that ($**$) is satisfied. But by \cref{rslt:prim-maps-induce-twist-of-shriek-functors,rslt:enriched-adjunction-of-shriek-functors} we have $f_! f^! = f_* \iHom(D_f, f^!) = \iHom(f_! D_f, \blank) = \iHom(f_* \one, \blank)$. This easily reduces ($**$) to the claim that $f_* \one$ is descendable, as desired. This finishes the proof of $\D^!$-descent along $f$. It follows from \cref{rslt:stability-of-suave-and-prim-maps} that the assumptions on $f$ are stable under any base-change, hence we get universal $\D^!$-descent.

The proof of $\D^*$-descent is similar and reduces to the observation that the functor $f_* f^* \in \Fun(\D(X), \D(X))$ is descendable, because by \cref{rslt:prim-maps-induce-twist-of-shriek-functors} we can write it as $f_* f^* \isom f_! (D_f \tensor f^*) = f_! D_f \tensor \blank \isom f_*\one \tensor \blank$. This finishes the proof of (i). Part (ii) follows from the proof (i), as indicated in the above explanation why ($**$) implies the claim. Part (iii) is a special case of (ii).
\end{proof}

In general it is hard to check primness locally on the source because it requires certain functors to commute with limits (see \cref{rslt:prim-map-is-local-on-source}). However, \cref{rslt:!-descent-along-prim-morphism} provides a useful exception:

\begin{corollary}
Let $\D$ be a stable 6-functor formalism on some geometric setup $(\cat C, E)$, let $f\colon Y \to X$ and $g\colon Z \to Y$ be maps in $E$ and assume that $g$ is $\D$-prim and $g_*\one$ is descendable.
\begin{corenum}
	\item\label{rslt:check-prim-object-on-!-descendable-cover} Let $P \in \D(Y)$ be such that $g^* P$ is $gf$-prim. Then $P$ is $f$-prim.
	\item If $gf$ is $\D$-prim then so is $f$.
	\item \label{rslt:check-properness-on-!-descendable-cover} If $f$ is truncated and $fg$ is $\D$-proper then $f$ is $\D$-proper.
\end{corenum}
\end{corollary}
\begin{proof}
Parts (ii) and (iii) follow formally from (i) as in the proof of \cref{rslt:localness-of-proper-on-source}. We now prove (i), so let $P \in \D(Y)$ be such that $g^* P$ is $gf$-prim. We want to apply \cref{rslt:suave-prim-is-local-on-source}. Unfortunately we cannot apply this result directly as we cannot guarantee a priori that $f_!(P \tensor \blank)$ preserves enough limits. However, looking at the proof of \cref{rslt:suave-prim-is-local-on-source} (and using \cref{rslt:descent-data-for-generated-sieve}) we observe that it is enough to show that certain functors commute with the limit $\varprojlim_{n\in\bbDelta} g_{n!} \DPrim_{fg_n}(g_n^* P)$, where $g_\bullet\colon Z_\bullet \to Y$ is the Čech nerve of $g$. Since all relevant functors are exact and hence commute with finite limits, it is enough to show that the cosimplicial object $(M^n)_n \coloneqq (g_{n!} \DPrim_{fg_n}(g_n^* P))_n$ in $\D(Y)$ is pro-constant. By \cref{rslt:!-descent-along-prim-morphism-implies-pro-const} this reduces to showing that $g^* M^\bullet$ is split.

Viewing $P$ as a morphism $Y \to X$ in $\cat K_{\D,X}$ we can describe $M^n$ as the composition of the left adjoint of $P \comp g_n$ with $g_n$. Let $Q \coloneqq \pi_{2*} \iHom(\pi_1^* P, \Delta_! \one) \in \D(Y)$ be the candidate left adjoint of $P$, where $\pi_i\colon Y \times_X Y \to Y$ are the two projections (see \cref{rslt:criterion-for-prim-object}) and let $g_n^\natural$ denote a left adjoint of $g_!$ (i.e.\ the image of the left adjoint of $g_n$ in $\cat K_{\D,Y}$ under $\Psi_{\D,Y}$). From the description of $\DSuave_{g_nf}(g_n^* P)$ in terms of the category of kernels it follows immediately that this sheaf identifies with $g_n^\natural Q$. We therefore obtain $M^\bullet = g_{\bullet!} g_\bullet^\natural Q$. We claim that every such cosimplicial object (with arbitrary $Q$) is $g^*$-split. By prim base-change (see \cref{rslt:prim-base-change}) we see that $g^* M^\bullet = g'_{\bullet!} g'^\natural_\bullet (g^* Q)$, where $g'_\bullet\colon Z'_\bullet \to Z$ is the base-change of $g_\bullet$ along $g$. But $g'_\bullet$ is a \emph{split} simplicial object and hence the same is true for $g^* M^\bullet$.
\end{proof}

\subsection{Examples} \label{sec:kerncat.examples}

We apply the results from the previous subsections to the 6-functor formalism of $\Lambda$-valued sheaves on condensed anima constructed in \cref{sec:6ff.example}. Our main application in that regard will be the study of classifying stacks of locally profinite groups, to which we devote the whole \cref{sec:reptheory}. In contrast, the present subsection is devoted to relating the 6-functor formalism on condensed anima to classical constructions in topology. We start with some basic definitions.

\begin{definition}
\begin{defenum}
	\item A map $U \to X$ of condensed anima is called an \emph{open immersion} (resp.\ \emph{closed immersion}) if after pullback to every profinite set $S$, it corresponds to an open (resp.\ closed) subset of $S$ via \cref{rslt:embed-Top-into-Cond-Ani}. By an \emph{open cover} (resp.\ a \emph{closed cover}) of $X$ we mean a family $(U_i \to X)_{i\in I}$ of open immersions (resp.\ closed immersions) which forms a cover, i.e.\ $\bigdunion_i U_i \surjto X$ is surjective.

	\item Let $\Lambda$ be an $\mathbb E_\infty$-ring. A map $Y \to X$ of condensed anima is called \emph{$\Lambda$-étale} if it is $\Lambda$-fine and $\D(\blank,\Lambda)$-étale in the sense of \cref{def:etale-and-proper-maps} (in particular truncated). We similarly define $\Lambda$-proper, $\Lambda$-suave and $\Lambda$-prim maps.
\end{defenum}
\end{definition}

Let us start with some easy basic computations, showing compatibility of open and closed immersions with the 6-functor formalism. In the following, recall the definition of quasiseparated condensed sets from \cref{def:quasiseparated-condensed-anima}.

\begin{lemma}
Let $\Lambda$ be an $\Einfty$-ring.
\begin{lemenum}
	\item \label{rslt:immersions-are-etale-resp-proper} Every open immersion of condensed anima is $\Lambda$-étale and every closed immersion of condensed anima is $\Lambda$-proper.

	\item \label{rslt:maps-from-profinite-sets-are-proper} Every map from a profinite set to a quasi-separated condensed set (e.g.\ a locally compact Hausdorff space) is $\Lambda$-proper.
\end{lemenum}
\end{lemma}
\begin{proof}
We start with the étaleness claim in (i), so fix an open immersion $j\colon U \injto X$. Since $\Lambda$-étaleness is local on the target by \cref{rslt:localness-of-etale-and-proper-maps} (note that $\Lambda$-fineness can also be checked locally on the target by construction, see \cref{rslt:extend-6ff-to-stacks-and-stacky-maps}), we can w.l.o.g.\ assume that $X$ is a profinite set and hence $U$ corresponds to an open subset of $X$. Suppose first that $U$ is compact, i.e.\ clopen. Then $U$ is profinite and hence $j$ is $\Lambda$-fine by construction. Moreover, $\Delta_j$ is an isomorphism and hence induces a map $j^! \to j^*$ of functors $\D(X, \Lambda) \to \D(U, \Lambda)$, which we need to show to be an isomorphism (see \cref{rslt:etale-equiv-iso-of-shriek-and-star}). Here we recall that $\D(U,\Lambda) = \Mod_{\Lambda(U)}$ and $\D(X,\Lambda) = \Mod_{\Lambda(X)}$, so that $j_! \isom j_*$ is the forgetful functor, $j^* = \blank \tensor_{\Lambda(X)} \Lambda(U)$ and $j^! = \iHom_{\Lambda(X)}(\Lambda(U), \blank)$. Note that $j_*$ is fully faithful, and hence the unit $\id \isoto j^! j_*$ and the counit $j^* j_* \isoto \id$ are isomorphisms. By going through the construction of the map $j^! \to j^*$ we see easily that it is composed of isomorphisms and hence is itself an isomorphism. This finishes the proof that $j$ is $\Lambda$-étale if $j$ is a clopen immersion.

Now let $j\colon U \injto X$ be a general open immersion (with $X$ still profinite). Pick an open cover $U = \bigunion_i U_i$ by clopen subsets $U_i \subseteq U$. Then by the previous paragraph each $U_i \injto U$ and each $U_i \injto X$ is $\Lambda$-étale. By \cref{rslt:Lambda-fineness-and-disjoint-unions} the map $\bigdunion_i U_i \to U$ is $\Lambda$-fine and by \cref{rslt:localness-of-etale-and-proper-maps} it is $\Lambda$-étale and in particular $\Lambda$-suave. It now follows from \cref{rslt:!-descent-along-suave-map} that this map has universal $!$-descent. Since also $\bigdunion_i U_i \to X$ is $\Lambda$-fine, we deduce by $!$-locality of $\Lambda$-fineness that $j$ is $\Lambda$-fine. Now \cref{rslt:localness-of-etale-and-proper-maps} immediately implies that it is $\Lambda$-étale.

It remains to prove (ii), which has the closed immersion claim of (i) as a special case. Thus let $f\colon S \to X$ be a map from a profinite set $S$ to a quasi-separated condensed set $X$. To show that $f$ is $\Lambda$-proper (and in particular $\Lambda$-fine), it is enough to do so after pullback along any map $T \to X$ from a profinite set $T$, i.e.\ we need to show that the map $S \times_X T \to T$ is $\Lambda$-proper. We claim that $S \times_X T$ is a profinite set. Indeed, by quasi-separatedness of $X$ it is a quasi-compact condensed set, and there is an injective map $S \times_X T \injto S \times T$. This is already enough to show that $S \times_X T$ is a profinite set (this is an easy exercise, see e.g.\ \cite[Lemma~2.7]{Clausen-Scholze:Complex-Geometry} for a reference). Changing notation slightly, we can now assume that $f\colon S \to T$ is a map of profinite sets. Then by construction $f$ is $\Lambda$-fine and $f_!$ is right adjoint to $f^*$ (in fact, $f_!$ is just the forgetful functor between module categories), and it is straightforward to check that this indeed implies $\Lambda$-properness.
\end{proof}

\begin{lemma} \label{rslt:finite-closed-cover-has-!-descent}
Let $(Z_i \injto X)_{i\in I}$ be a finite closed cover of a condensed anima $X$ and let $f\colon Z \coloneqq \bigdunion_i Z_i \to X$ be the induced map. Then $f_* \one \in \D(X, \Lambda)$ is descendable, for every $\Einfty$-ring $\Lambda$. In particular $f$ has universal $\D^!$-descent.
\end{lemma}
\begin{proof}
The \enquote{in particular} claim follows from \cref{rslt:!-descent-along-prim-morphism}, using that $f$ is $\Lambda$-proper (and in particular $\Lambda$-prim) by \cref{rslt:immersions-are-etale-resp-proper,rslt:localness-of-proper-on-source}. Let us now prove that $f_* \one$ is descendable. This is a standard argument, see e.g.\ \cite[Lemma~2.6.11]{Mann.2022a} for a related result. For every non-empty subset $J \subseteq I$ let us denote $Z_J \coloneqq \prod_{j\in J}^{/X} Z_j$ and let us write $f_J\colon Z_J \to X$ for the obvious map. Then $f_*\one = \prod_i f_{\{i\}*} \one$, so that each $f_{\{i\}*} \one$ is a retract of $f_* \one$ and hence lies in $\langle f_* \one \rangle$. Moreover, by the projection formula we have $f_{J*} \one \isom \bigtensor_{j\in J} f_{\{j\}*} \one$ and hence $f_{J*} \one \in \langle f_* \one \rangle$ for every non-empty $J \subset I$. But by \cref{rslt:descent-data-for-generated-sieve-by-monomorphisms} we have $\D(X,\Lambda) = \varprojlim_J \D(Z_J,\Lambda)$ and in particular the natural map $\one \isoto \varprojlim_J f_{J*} \one$ is an isomorphism. Since all terms in the limit on the right lie in $\langle f_* \one \rangle$ the same follows for the limit, i.e.\ $\one \in \langle f_* \one \rangle$. Thus $f_* \one$ is indeed descendable.
\end{proof}

\begin{corollary}
Let $f\colon Y \to X$ be a map of condensed anima and $\Lambda$ an $\Einfty$-ring.
\begin{corenum}
	\item \label{rslt:open-local-properties-for-condensed-anima} $f$ is $\Lambda$-fine (resp.\ $\Lambda$-suave, resp.\ $\Lambda$-étale) if it is so on an open cover of $X$ and $Y$.

	\item \label{rslt:closed-local-properties-for-condensed-anima} $f$ is $\Lambda$-fine (resp.\ $\Lambda$-suave, resp.\ $\Lambda$-prim, resp.\ $\Lambda$-proper) if it is so on a finite closed cover of $X$ and $Y$.
\end{corenum}
\end{corollary}
\begin{proof}
We first prove (i). Note that every cover of condensed anima is in particular a $\D^*$-cover by descent for $\D(\blank,\Lambda)$. Now the claim for $\Lambda$-fineness follows from $\D^!$-locality of this property together with the observation that an open cover is in particular a universal $\D^!$-cover. To prove the latter, we argue as in the proof of \cref{rslt:immersions-are-etale-resp-proper}: If $(U_i \to X)_i$ is an open cover of $X$ then $\bigdunion_i U_i \to X$ is $\Lambda$-fine and $\Lambda$-étale, hence a universal $\D^!$-cover by \cref{rslt:!-descent-along-suave-map}. This implies that $(U_i \to X)_i$ is a universal $\D^!$-cover. The other claims in (i) follow from locality of the respective properties, see \cref{rslt:suave-map-is-local-on-source,rslt:localness-of-etale-and-proper-maps}.

Part (ii) can proved very similarly to (i), using \cref{rslt:finite-closed-cover-has-!-descent} to get the necessary $!$-descent and using \cref{rslt:prim-map-is-local-on-source,rslt:localness-of-proper-on-source} for the localness of the respective properties (see also \cref{rmk:how-to-apply-locality-of-primness-on-source}).
\end{proof}

Next up we show two basic properties that one expects from a theory of \enquote{sheaves on spaces}: conservativity of stalks and the excision fiber sequences.

\begin{lemma} \label{rslt:cond-ani-stalks-are-conservative}
Let $X$ be a condensed anima. Then the collection of functors $p^*\colon \D(X,\Lambda) \to \D(*,\Lambda)$, for all maps $p\colon * \to X$, is conservative.
\end{lemma}
\begin{proof}
By covering $X$ by profinite sets (pullback to which form a conservative family of functors by construction), we can reduce to the case that $X$ is a profinite set, so that $\D(X,\Lambda) = \Mod_{\Lambda(X)}$. But then the claim follows from (the proof of) \cref{rslt:identify-sheaves-on-profinite-sets-with-D}.
\end{proof}

\begin{proposition} \label{rslt:excision-for-condensed-anima}
Let $X$ be a locally compact Hausdorff space, $i\colon Z \injto X$ a closed subset and $j\colon U \injto X$ the complement open subset. Then for every $\Einfty$-ring $\Lambda$ and every $M \in \D(X,\Lambda)$, we have the following fiber sequences in $\D(X,\Lambda)$:
\begin{align}
	j_! j^* M \to M \to i_* i^* M, \qquad i_* i^! M \to M \to j_* j^* M.
\end{align}
In particular, the two pairs of functors $(j^*, i^*)$ and $(j^*, i^!)$ form conservative families of functors on $\D(X,\Lambda)$.
\end{proposition}
\begin{proof}
By \cref{rslt:immersions-are-etale-resp-proper} we know that $j$ is $\D$-étale and $i$ is $\D$-proper, so that we get identifications $i_* = i_!$ and $j^* = j^!$. Using the projection formula and \cref{rslt:enriched-adjunction-of-shriek-functors} we can rewrite the claimed fiber sequences as
\begin{align}
	M \tensor j_! \one \to M \tensor \one \to M \tensor i_* \one, \qquad \iHom(i_* \one, M) \to \iHom(\one, M) \to \iHom(j_! \one, M),
\end{align}
hence the claim reduces to constructing a fiber sequence of the form $j_! \one \to \one \to i_* \one$. But note that by proper base-change we have $j^* i_* \one = 0$ (because $U \times_X Z = \emptyset$), hence $j^* \fib(\one \to i_* \one) = j^* \one = \one$. This produces a map $j_! \one \to \fib(\one \to i_* \one)$ via adjunctions, and we need to show that this map is an isomorphism. This is clearly true after applying $j^*$ and $i^*$, so it remains to show that the pair $(j^*, i^*)$ forms a conservative family of functors on $\D(X,\Lambda)$. This follows immediately from \cref{rslt:cond-ani-stalks-are-conservative}.
\end{proof}

We now come to the promised discussion of classical algebraic topology in the context of our 6-functor formalism. Here the most important computation is for the unit interval $I = [0, 1]$:

\begin{lemma} \label{rslt:unit-interval-computations}
Let $I = [0, 1] \subset \RR$ be the unit interval and $f\colon I \to *$ be the projection. Let $\Lambda$ be an $\Einfty$-ring. Then:
\begin{lemenum}
	\item $f$ is $\Lambda$-proper (and in particular $\Lambda$-fine).
	\item \label{rslt:unit-interval-cohomology} For any condensed anima $X$ with projection $f_X\colon X \times I \to X$, the pullback functor
	\begin{align}
		f_X^*\colon \D(X,\Lambda) \injto \D(X \times I,\Lambda)
	\end{align}
	is fully faithful. In particular $f_{X*} \one = \one$ and $\Gamma(I,\Lambda) = \Lambda$.
\end{lemenum}
\end{lemma}
\begin{proof}
Let $C \isom \prod_{\ZZ_{>0}} \{ 0, 1 \}$ be the Cantor set. This is a profinite set with a continuous surjective map $g\colon C \to I$ given by $(a_n)_{n>0} \mapsto \sum_{n=1}^\infty a_n 2^{-n}$. By \cref{rslt:maps-from-profinite-sets-are-proper} $g$ is $\Lambda$-proper. We claim that it has universal $\D^!$-descent, for which by \cref{rslt:!-descent-along-prim-morphism} it is enough to show that $A \coloneqq g_* \one \in \D(I,\Lambda)$ is descendable. The following argument is due to Clausen--Scholze (see lecture 20, minute 35 in \cite{Clausen-Scholze:Youtube-Analytic-Geometry}). For each $n > 0$ let $C_n$ be the disjoint union of $2^n$ consecutive compact subintervals of $I$ of equal length, e.g.\ $C_1 = [0,\frac12] \dunion [\frac12,1]$. Then there are natural maps $C \to C_n \to C_m \to I$ for $n \ge m > 0$ and one easily checks that the induced map $C \isoto \varprojlim_n C_n$ is an isomorphism of topological spaces (both sides are compact Hausdorff and the map is continuous, so one only needs to verify that it is a bijection). Let $g_n\colon C_n \to I$ be the projection and $A_n \coloneqq g_{n*} \one \in \D(I,\Lambda)$. Then the natural map $A \isoto \varinjlim_n A_n$ is a colimit of algebras. Indeed, this can be checked after pullback to any profinite set; then all the involved spaces become profinite sets and hence the claim reduces to a simple property of module categories.

Now recall the notion of the \emph{index} of descendability from \cite[Definition~11.18]{Bhatt-Scholze:2017} and also recall the fact that descendability is preserved by countable sequential colimits with bounded index (see \cite[Lemma~11.22]{Bhatt-Scholze:2017}). Therefore, to prove descendability of $A$ it is enough to show that each $A_n$ is descendable of index independent of $n$. By \cref{rslt:finite-closed-cover-has-!-descent} $A_n$ is descendable. To bound the index, we need to find a uniform integer $d \ge 0$ such that for all $n > 0$ the map $F_n^{\tensor d} \to \Lambda$ is zero, where $F_n \coloneqq \fib(\Lambda \to A_n)$ (by \cite[Lemma~11.20]{Bhatt-Scholze:2017} there is such a $d$ for each $n$, but we need to make sure that $d$ can be chosen independently of $n$). Note that $F_n$ is supported at finitely many points, namely those points where the intervals comprising $C_n$ intersect. Thus the question whether $F_n^{\tensor d} \to \Lambda$ is zero reduces to a local question around each of these points. But locally around a point, the situation is isomorphic to the case $n = 1$, hence the same $d$ that works for $n = 1$ works for all $n \ge 1$. This finishes the proof that $g$ has universal $!$-descent. By $!$-locality of $\Lambda$-fine maps it follows that $f$ is $\Lambda$-fine. For part (i) it only remains to prove that $f$ is $\Lambda$-proper. But this follows immediately from \cref{rslt:check-properness-on-!-descendable-cover} because $fg$ is $\Lambda$-proper by \cref{rslt:maps-from-profinite-sets-are-proper}.

We now prove part (ii), so let the condensed anima $X$ be given. We need to show that for every $M \in \D(X, \Lambda)$ the unit $M \isoto f_{X*} f_X^* M$ is an isomorphism. By (i) we know that $f_X$ is $\Lambda$-proper, so that $f_X^*$ and $f_{X*}$ are represented by adjoint morphisms in $\cat K_{\D(\blank,\Lambda)}$. Hence it is enough to show that the unit of their adjunction in $\cat K_{\D(\blank,\Lambda)}$ is an isomorphism, which immediately boils down to showing that the map $\Lambda \isoto f_{X*} f_X^* \Lambda$ is an isomorphism (secretly, the previous two sentences used the projection formula for $f_{X*}$ to reduce the case of general $M$ to that of $M = \Lambda$, but using the category of kernels it is immediately clear that the \enquote{natural maps} are the correct ones). Using conservativeness of stalks (see \cref{rslt:cond-ani-stalks-are-conservative}) and proper base-change (see \cref{rslt:suave-and-prim-base-change}) we can reduce to the case $X = *$. So all in all we are reduced to showing that the natural map $\Lambda \isoto \Gamma(I,\Lambda)$ is an isomorphism. This is of course a standard fact from algebraic topology and can therefore be deduced from \cref{rslt:identify-6ff-with-sheaves-on-top-spaces} and classical results. However, in the following we provide a self-contained argument using only the tools developed in this paper.	

Let $\CHaus$ be the category of compact Hausdorff spaces, equipped with the site of finite jointly surjective covers. For every $X \in \CHaus$ we get an associated profinite set $\pi_0(X)$ of connected components (see \cref{rslt:profinite-completion}) and we can hence construct the $\Lambda$-module $\Lambda(X)' \coloneqq \Lambda(\pi_0(X))$ using \cref{def:Einfty-ring-associated-to-profin-set}. This defines a functor $\Lambda(\blank)'\colon \CHaus^\op \to \Mod_\Lambda$ which agrees with $\Gamma(\blank,\Lambda)$ on profinite sets. Since profinite sets form a basis of $\CHaus$, the sheafiness of $\Gamma(\blank,\Lambda)$ (see \cref{rslt:cohomology-is-sheafy}) induces a natural transformation $\alpha\colon \Lambda(\blank)' \to \Gamma(\blank,\Lambda)$. Our goal is to show that $\alpha_I\colon \Lambda(I)' \isoto \Gamma(I,\Lambda)$ is an isomorphism. Since $\alpha_S$ is an isomorphism for all profinite sets $S$ and $\Gamma(\blank,\Lambda)$ satisfies descent, it is enough to show that $\Lambda(\blank)'$ descends along the map $C \to I$, i.e.\ we need to show that the natural map
\begin{align}
	\Lambda = \Lambda(I)' \isoto \varprojlim_{n\in\bbDelta} \Lambda(C^{\times_I (n+1)})' = \varprojlim_{n\in\bbDelta} \Gamma(C^{\times_I (n+1)},\Lambda)
\end{align}
is an isomorphism. Clearly the map $\Lambda \to \Gamma(C,\Lambda)$ is descendable (i.e.\ $\Gamma(C, \Lambda)$ is a descendable object in $\D(\Lambda)$), it is even split. This implies that the map $\Lambda \isoto \varprojlim_{n\in\bbDelta} \Lambda(C)^{\tensor (n+1)}$ is an isomorphism. There is a natural morphism $\Lambda(C)^{\tensor (\bullet+1)} \to \Gamma(C^{\times_I (\bullet+1)},\Lambda)$ of cosimplicial $\Lambda$-modules, which can e.g.\ be obtained by observing that $\Lambda(C)^{\tensor (\bullet+1)}$ is a left Kan extension of its restriction to $\{ [0] \} \subset \bbDelta$ when viewed as a cosimplicial $\Lambda$-algebra. All in all this reduces the problem to showing that the natural maps
\begin{align}
	\Lambda(C)' \tensor \dots \tensor \Lambda(C)' \isoto \Lambda(C \times_I \dots \times_I C)'
\end{align}
are isomorphisms. Both sides transform cofiltered limits in $C$ into colimits (using that $\pi_0(\blank)$ commutes with cofiltered limits by \cref{rslt:profinite-completion}), so we can replace $C$ by $C_n$. But then both sides of the claimed isomorphism are given by free $\Lambda$-modules of the same rank (e.g.\ $\Lambda(C_n)' = \Lambda^{2^n}$) and one easily checks that the induced map is the canonical isomorphism.
\end{proof}

\Cref{rslt:unit-interval-computations} provides an explicit description of the whole 6-functor formalism on the unit interval $I$ and the map $f\colon I \to *$. Namely, $\D(I, \Lambda)$ is the category of $\Lambda$-valued sheaves on $I$ (see \cref{rslt:identify-6ff-with-sheaves-on-top-spaces}), $f^*\colon \D(\Lambda) \to \D(I, \Lambda)$ is the constant sheaf functor, $f_*\colon \D(I, \Lambda) \to \D(\Lambda)$ is the cohomology functor, $f_!$ agrees with $f_*$ and $f^!$ is the right adjoint of $f_*$. We can now use \cref{rslt:unit-interval-computations} to deduce the following general results from classical algebraic topology.

\begin{definition}
\begin{defenum}
	\item \label{def:finite-dimensional-CHaus-space} A locally compact Hausdorff space is \emph{locally finite dimensional} if on an open cover it admits a locally closed immersion into $I^n$ for varying $n \ge 0$.

	\item A \emph{manifold} is a Hausdorff topological space which can be covered by open subsets that are homeomorphic to an open subset of $\RR^n$ for some $n \ge 0$.
\end{defenum}
\end{definition}

\begin{theorem} \label{rslt:6ff-on-manifolds}
Let $\Lambda$ be an $\Einfty$-ring and let $X$ and $Y$ locally finite dimensional locally compact Hausdorff spaces.
\begin{thmenum}
	\item Every map $Y \to X$ is $\Lambda$-fine.
	\item If $X$ and $Y$ are compact, then every map $Y \to X$ is $\Lambda$-proper.
	\item The map $X \to *$ is $\Lambda$-suave and the associated dualizing complex $\omega_X$ is invertible and locally isomorphic to $\Lambda[n]$, where $n$ is the local dimension.
\end{thmenum}
\end{theorem}
\begin{proof}
We first prove (i), which cancellativeness reduces to showing that $X \to *$ is $\Lambda$-fine. By \cref{rslt:open-local-properties-for-condensed-anima} this can be checked on an open cover of $X$, so we can reduce to the case that $X$ is a locally closed subset of $I^n$. Since open and closed immersions are $\Lambda$-fine (see \cref{rslt:immersions-are-etale-resp-proper}) we reduce to the case $X = I^n$ and then to the case $X = I$. But this case was handled in \cref{rslt:unit-interval-computations}.

We now prove (ii), so assume that $X$ and $Y$ are compact and fix a map $f\colon Y \to X$. By \cref{rslt:closed-local-properties-for-condensed-anima} it is enough to show that $f$ is $\Lambda$-proper on a finite closed cover of $X$ and $Y$. This allows us to assume that $X$ and $Y$ are closed subsets of $I^n$ and $I^m$ respectively, for some integers $n, m \ge 0$. Using cancellativeness of $\Lambda$-properness (which follows from \cref{rslt:stability-of-etale-and-proper-maps}) we easily reduce the claim to showing that $I^n \to *$ is $\Lambda$-proper, which follows immediately from \cref{rslt:unit-interval-computations}.

It remains to prove (iii). Since $\Lambda$-suaveness is local on $X$ (see \cref{rslt:open-local-properties-for-condensed-anima}) we immediately reduce to the case $X = \RR^n$ for some $n$. Since $\Lambda$-suaveness is also stable under composition (see \cref{rslt:stability-of-suave-and-prim-maps}) we can even assume $X = \RR$. Write $\RR = \bigunion_n I_n$, where $I_n = [-n, n] \subset \RR$. Then
\begin{align}
	\Gamma(\RR,\Lambda) = \varprojlim_n \Gamma(I_n,\Lambda) = \varprojlim_n \Lambda = \Lambda	
\end{align}
by \cref{rslt:unit-interval-cohomology}. We similarly compute $\Gamma(\RR^n,\Lambda) = \Lambda$. Now let us denote $p\colon I \to *$ the projections and recall the definition of $\Gamma_c(\RR,\Lambda)$ from \cref{def:Gamma-c}. By choosing an embedding $j\colon \RR \isoto (0,1) \injto [0,1]$ and using the excision sequence from \cref{rslt:excision-for-condensed-anima} as well as \cref{rslt:unit-interval-cohomology} we deduce that
\begin{align}
	\Gamma_c(\RR,\Lambda) = p_* j_! \one = \fib(\Gamma(I,\Lambda) \to \Gamma(\{ 0 \} \dunion \{ 1 \}, \Lambda)) \isom \Lambda[-1].
\end{align}
A similar computation works in the relative setting and in particular we deduce $\Gamma_c(\RR^n,\Lambda) \isom \Lambda[-n]$ inductively.

We now claim that the projection $f\colon \RR \to *$ is left adjoint to the map $g \coloneqq \one[1]\colon * \to \RR$ in $\K_{\D(\blank,\Lambda)}$; this proves the desired $\Lambda$-suaveness of $\RR$. To prove the adjunction, we follow the proof of \cite[Proposition~5.10]{Scholze:Six-Functor-Formalism} and explicitly construct the unit and counit maps.

For the counit map we observe that $f \comp g = f_! \one[1] = \Gamma_c(\RR,\Lambda)[1] \isom \one$ by the above computations. For the unit map we first observe that $g \comp f = \one[1]$ (viewed as a sheaf on $\RR^2$). Now let $j\colon U \injto \RR^2$ be the complement of the diagonal, so that \cref{rslt:excision-for-condensed-anima} gives us a fiber sequence $j_!\one \to \one \to \Delta_! \one$ in $\D(\RR^2, \Lambda)$. Using the fact that $U \isom \RR^2 \dunion \RR^2$ we compute
\begin{align}
	&\Hom(\Delta_! \one, \one[1]) = \Hom(\Delta_! \one, \one)[1] = \fib(\Hom(\one, \one) \to \Hom(j_! \one, \one))[1] \\&\qquad= \fib(\Gamma(\RR^2,\Lambda) \to \Gamma(U,\Lambda))[1] \isom \fib(\Lambda \to \Lambda^2)[1] = \Lambda
\end{align}
by the above computations; here the map $\Lambda \to \Lambda^2$ is the diagonal embedding. In the terminology of the category of kernels, the above computation reads $\Hom(\id_\RR, g \comp f) \isom \Lambda$, which provides the desired unit map.

It remains to check that the triangle identities hold for $f$ and $g$, i.e.\ that the compositions $f \to f g f \to f$ and $g \to g f g \to g$ are isomorphic to the identity. In fact, by \cref{rslt:adjoints-from-weak-triangle-identity} it is enough to show that these maps are isomorphisms. Since the counit $f g \isoto \one$ is an isomorphism by construction, this reduces the claim to showing that the maps $f \to fg f$ and $g \to g f g$ are isomorphisms. We only discuss the first map, the second being similar. Note that both $f$ and $f g f$ are represented by $\one \in \D(\RR,\Lambda)$. Moreover, the map $f \to fg f$ is obtained as $\pi_{2!}(\Delta_! \one \to \one[1])$, where $\pi_2\colon \RR^2 \to \RR$ is the projection to the second factor and $\Delta_! \one \to \one[1]$ is the unit map from above. In other words, we need to see that the map $\pi_{2!}\colon \Lambda \isom \Hom(\Delta_! \one_\RR, \one_{\RR^2}[1]) \to \Hom(\one_\RR, \one_\RR) \isom \Lambda$ is an isomorphism (then it certainly sends $1 \in \Lambda$ to a unit in $\Lambda$ and hence an isomorphism). Note first that $f_!\colon \Hom(\one_\RR, \one_\RR) \to \Hom(\one_*,\one_*)$ is an isomorphism because it is a map of $\Lambda$-algebras which are both isomorphic to $\Lambda$ by the above computations. It is therefore enough to show that the composed map $f_!\pi_{2!}\colon \Hom(\Delta_! \one_\RR, \one_{\RR^2}[1]) \to \Hom(\one_*, \one_*)$ is an isomorphism. The fiber sequence computing $\Hom(\Delta_! \one, \one[1])$ induces the following commuting diagram of fiber sequences in $\Mod_\Lambda$:
\begin{equation}\begin{tikzcd}
	\Hom(\Delta_! \one_\RR, \one_{\RR^2}[1]) \arrow[r] \arrow[d,"f_!\pi_{2!}"] & \Hom(\one_{\RR^2}, \one_{\RR^2}[1]) \arrow[d,"f_!\pi_{2!}"] \arrow[r] & \Hom(j_! \one_U, \one_{\RR^2}[1]) \arrow[d,"f_!\pi_{2!}"]\\
	\Hom(\one_*[-1], \one_*[-1]) \arrow[r] & \Hom(\one_*[-2], \one_*[-1]) \arrow[r] & \Hom(\one_*[-2]^2, \one_*[-1])
\end{tikzcd}\end{equation}
Here $j\colon U \injto \RR^2$ is the inclusion of the complement of the diagonal and we recall $U \isom \RR^2 \dunion \RR^2$. It is enough to show that the middle and right-hand vertical maps are isomorphisms. For the middle map this follows by the same argument as for $f_!$ (because both $\Hom$'s are endomorphism spaces up to a shift). For the right-hand vertical map we can split $j$ into the disjoint union of two open immersions $j_1, j_2\colon \RR^2 \injto \RR^2$ given by the part above the diagonal and below the diagonal, respectively. This reduces the claim to showing that the map $f_!\pi_{2!}\colon \Hom(j_{i!} \one_{\RR^2}, \one_{\RR^2}[1]) \isoto \Hom(\one_*[-2], \one_*[-1])$ is an isomorphism for $i = 1, 2$. Now we apply the excision sequence for $j_i$ and its closed complement to the target $\one_{\RR^2}[1]$ to split everything into simple pieces that can be handled by the above methods.
\end{proof}

\begin{remark}
From \cref{rslt:6ff-on-manifolds} one easily deduces classical Poincaré duality results on manifolds, the Künneth formula on manifolds and more. Of course it should be possible to extend these results to far more general spaces, e.g.\ CW complexes. However, a general map of compact Hausdorff spaces may not be $\Lambda$-proper (and probably not even $\Lambda$-fine): For the map $f\colon I^\NN \to *$ the functor $f_*\colon \D(I^\NN, \Lambda) \to \D(\Lambda)$ does not satisfy base-change (see \cite[Counterexample~6.5.4.2]{HTT}). This is the price we pay for assuming strong descent on $\D(\blank,\Lambda)$ (which results in working with hypercomplete sheaves rather than all sheaves on $X$). One can also start with compact Hausdorff spaces instead of profinite sets to avoid these issues, but at the cost of having less descent to work with (see \cite[\S13.3]{Scholze:Six-Functor-Formalism}).
\end{remark}

We see the above results as a proof of concept for our 6-functor formalism and do not attempt to develop the whole theory of sheaves on (nice) locally compact Hausdorff spaces here, as our focus lies on applications to representation theory of locally profinite groups. \Cref{rslt:6ff-on-manifolds} essentially contains Poincaré duality on manifolds; let us restate it in the following form:

\begin{corollary} \label{rslt:Poincare-duality-on-manifolds}
Let $M$ be a manifold and $\Lambda$ be an $\Einfty$-ring. Then there is a natural isomorphism
\begin{align}
	\Gamma(M,\omega_M) = \Gamma_c(M,\Lambda)^\vee
\end{align}
of $\Lambda$-modules, where $\omega_M$ is the sheaf from \cref{rslt:6ff-on-manifolds}.
\end{corollary}
\begin{proof}
Let $f\colon M \to *$ be the projection. Then
\begin{align}
	\Gamma(M,\omega_M) = f_* f^! \one = f_* \iHom(\one, f^! \one) = \iHom(f_! \one, \one) = \Gamma_c(M,\Lambda)^\vee,
\end{align}
as desired.
\end{proof}

Note that $\Gamma_c(M,\Lambda)$ plays the role of cohomology with compact support here. Indeed, it follows easily from the definitions that we can compute it as the colimit of $\Gamma_c(U,\Lambda)$ over open subsets $U \subseteq M$. Moreover, if $U$ embeds as an open subset into $[0,1]^n$ for some $n$, with complement $Z$, then $\Gamma_c(U,\Lambda) = \fib(\Gamma([0,1]^n,\Lambda) \to \Gamma(Z,\Lambda))$.

We have now discussed how cohomology and cohomology with compact support can be represented by our 6-functor formalism and we have used the six functors to prove that these two cohomologies are related by Poincaré duality on a manifold. What is still missing from the picture is homology, which we will explain next. Here it becomes useful to work with condensed \emph{anima} and not just condensed sets. We start with the following observation:

\begin{proposition} \label{rslt:maps-of-anima-are-etale}
Let $\Lambda$ be an $\Einfty$-ring and let $f\colon Y \to X$ be a truncated map of anima. Then $f$ is $\Lambda$-étale.
\end{proposition}
\begin{proof}
Let $n \ge 0$ be such that $f$ is $n$-truncated. We argue by induction on $n$. By \cref{rslt:localness-of-etale-and-proper-maps} the claim is local on $X$, so after passing to a covering of $X$ by copies of $*$ we can assume that $X = *$. In particular $Y$ is an $n$-truncated anima for some $n$. The case $n = 0$ is obvious, as then $Y$ is a disjoint union of points (use e.g.\ \cref{rslt:open-local-properties-for-condensed-anima}). In general we can cover $Y$ by maps $* \to Y$ (since $\Ani = \Shv(*)$), so it is enough to prove that every such map is étale (then $f$ is $\Lambda$-fine by $!$-localness of $\Lambda$-fine maps and it is étale by \cref{rslt:localness-of-etale-and-proper-maps}). But by \cref{rslt:truncated-composition} every map $* \to Y$ is $(n-1)$-truncated, so we finish by induction.
\end{proof}

\begin{remark} \label{rmk:require-truncated-anima}
The fact that we need to assume truncatedness in \cref{rslt:maps-of-anima-are-etale} is a deficiency of the current construction of the 6-functor formalism. We aim to fix this in future work.
\end{remark}

Note that for an anima $X$ we have $\D(X,\Lambda) = \Fun(X,\Mod_\Lambda)$ as both sides send colimits in $X$ to limits and agree on $X = *$. In particular $\D(X,\Lambda)$ is the usual category of \enquote{$\Lambda$-valued local systems on $X$} from topology. For a map $f\colon Y \to X$ of anima, the functor $f^*\colon \D(X,\Lambda) \to \D(Y,\Lambda)$ is given by precomposition with $f$, and by \cref{rslt:maps-of-anima-are-etale} (assuming $f$ is truncated) the functor $f_!$ is left adjoint to $f^*$ and hence computes \enquote{relative homology}. In particular, for a (truncated) anima $X$, $\Gamma_c(X,\Lambda)$ computes the $\Lambda$-valued homology of $X$.

\begin{example} \label{ex:primness-of-anima}
We can make the following fun construction. Let $X$ be a CW complex (viewed as a condensed set) and let $\abs X$ denote the associated anima. Assume that $\abs X$ is truncated (this is necessary for now, see \cref{rmk:require-truncated-anima}). By \cite[Lemma~11.9]{Clausen-Scholze:Analytic-Geometry} there is a map $\alpha_X\colon X \to \abs X$ of condensed anima and we expect $\alpha_X^*\colon \D(\abs X,\Lambda) \injto \D(X,\Lambda)$ to be the fully faithful embedding of $\Lambda$-valued local systems into all $\Lambda$-valued sheaves on $X$. Suppose that $X$ is a compact manifold and hence $\Lambda$-suave and $\Lambda$-proper by \cref{rslt:6ff-on-manifolds}. Since $\abs X$ is $\Lambda$-étale by \cref{rslt:maps-of-anima-are-etale}, we deduce from \cref{rslt:cancellativeness-of-suave-and-prim-maps} that $\alpha_X$ is $\Lambda$-suave and hence a left adjoint in $\cat K_{\D(\blank,\Lambda)}$. Denote its right adjoint by $\rho_X$. Then the fully faithfulness of $\alpha_X^*$ (which we use conjecturally here) implies that $\alpha_X \rho_X = \id_{\abs X}$ and hence the map $f\colon \abs{X} \to *$ can be written as the composition of $\rho_X\colon \abs{X} \to X$ and the projection $X \to *$. This implies that $f$ is right adjoint in $\cat K_{\D(\blank,\Lambda)}$ and hence $\Lambda$-prim! One checks that $\delta_f = \omega_X$, which results in the well-known identity
\begin{align}
	f_* L = f_!(L \tensor \omega_X)
\end{align}
for every local system $L$ on $X$, where $f_*$ is the cohomology and $f_!$ is the homology functor.
\end{example}

\section{Representation theory} \label{sec:reptheory}

In the following we explain how the theory of 6-functor formalisms developed above (and specifically the 6-functor formalism $\D(\blank,\Lambda)$ on condensed anima) can be used to gain a new perspective on classical constructions appearing in the theory of smooth representations. We believe this perspective to be especially useful in the setting of modular representations in natural characteristic, where many constructions are forced to be carried out on the derived level and hence harder to grasp using conventional methods.

\subsection{Derived categories} \label{sec:reptheory.descent}

We will make use of the 6-functor formalism $\D(\blank,\Lambda)$ on condensed anima constructed in \cref{sec:6ff.example}. If $\Lambda$ is static (i.e.\ an ordinary ring) and $X$ is a profinite set then $\D(X,\Lambda) = \Mod_{\Lambda(X)}$ is the derived category of its heart (see \cref{ex:Einfty-rings-and-modules}), where $\Lambda(X)$ is the $\Lambda$-algebra of locally constant functions $X\to \Lambda$ with pointwise addition and multiplication. In the following we investigate to what extent this property holds for more general condensed anima $X$. 

Recall from \cref{rslt:Lambda-sheaves-t-structure} that the categories $\D(X,\Lambda)$ admit a natural left-complete t-structure, for every condensed anima $X$ and every connective $\Einfty$-ring $\Lambda$, and all pullback functors in $\D(\blank,\Lambda)$ are t-exact. We can now ask when $\D(X,\Lambda)$ is the derived category of its heart and this heart is constructed via abelian descent; more concretely we make the following definition:

\begin{definition}
Let $\Lambda$ be a commutative ring. Consider the hypersheaf
\begin{align}
	\cat A(\blank,\Lambda)\colon \Cond(\Ani)^\op \to \Cat, \qquad X \mapsto \cat A(X,\Lambda) \coloneqq \D(X,\Lambda)^{\heartsuit}
\end{align}
of abelian categories. Note that $\cat A(S,\Lambda)$ is Grothendieck whenever $S$ is a profinite set; by descent it follows that $\cat A(X,\Lambda)$ is Grothendieck for every condensed anima $X$ (see also the proof of \cref{rslt:Grothendieck-descent}). We say that a condensed anima $X$ is \emph{$\Lambda$-naively derived} if the induced functor
\begin{align}
	\hat{\D}(\cat A(X,\Lambda)) \isoto \D(X,\Lambda)
\end{align}
is an isomorphism (here $\hat{\D}$ denotes the left-completion of $\D$).
\end{definition}

\begin{example}\label{ex:profinite-sets-are-naively-derived}
Every profinite set is $\Lambda$-naively derived for every commutative ring $\Lambda$.
\end{example}

\begin{lemma}\label{rslt:disjoint-union-is-naively-derived}
Let $\Lambda$ be a commutative ring and $(X_i)_i$ a family of $\Lambda$-naively derived condensed anima. Then $\bigdunion_{i} X_i$ is $\Lambda$-naively derived.
\end{lemma}
\begin{proof}
This is clear, since both functors $\hat{\D}(\cat A(\blank,\Lambda))$ and $\D(\blank, \Lambda)$ send disjoint unions to products.
\end{proof}

\begin{lemma} \label{rslt:product-with-profin-set-stays-naively-derived}
Let $\Lambda$ be a commutative ring and $X$ a $\Lambda$-naively derived condensed anima. Then for every profinite set $S$, also $X \times S$ is $\Lambda$-naively derived.
\end{lemma}
\begin{proof}
Fix a profinite set $S$. Denote by $p_X\colon X\times S\to X$ the projection, and let $f\colon X\to *$ and $g\colon S\to *$ be the obvious maps. We will proceed in several steps.
\medskip

\textit{Step~1:} The functor $p_{X*}$ is monadic, t-exact, and satisfies $p_{X*}(\one) \simeq \Lambda(S)$ (where we view $\Lambda(S)$ as an object of $\D(X,\Lambda)$ via pullback along $X\to *$); in particular, $p_{X*}$ induces a (t-exact) isomorphism
\begin{align}\label{eq:product-with-profin-set-monadic}
\D(X\times S,\Lambda) \isoto \Mod_{\Lambda(S)}(\D(X,\Lambda)) .
\end{align}
Note that $p_X$ is $\Lambda$-proper as the base-change of the $\Lambda$-proper map $S\to *$ (see \cref{rslt:maps-from-profinite-sets-are-proper}). In particular, $p_{X*}$ satisfies the projection formula, so that we have a natural isomorphism $p_{X*}p_X^* \simeq p_{X*}(\one)\tensor \blank$ of functors $\D(X,\Lambda) \to \D(X,\Lambda)$. Moreover, by base-change, we have an isomorphism
\begin{align}
p_{X*}(\one) = p_{X*}p_S^*(\one) \simeq f^*g_*(\one) = f^*\Lambda(S).
\end{align}
We show that $p_{X*}$ is t-exact and conservative. This is obvious if $X$ is profinite. For general $X$, both properties can be checked after pullback along $h\colon T\to X$, where $T$ is a profinite set: This follows from the base-change formula $h^*p_{X*} \simeq g_*(h\times \id_S)^*$ and the fact that the functors $(h\times \id_S)^*$ are jointly conservative by \cref{rslt:cond-ani-stalks-are-conservative}. As $p_{X*}$ commutes with colimits (it admits $p_X^!$ as a right adjoint), we deduce that $p_{X*}$ is monadic. The isomorphism \eqref{eq:product-with-profin-set-monadic} now follows from the Barr--Beck theorem \cite[Theorem~4.7.3.5]{HA}. 
\medskip

\textit{Step~2:} The natural map
\begin{align}
\hat{\D}\bigl(\Mod_{\Lambda(S)}(\cat A(X,\Lambda))\bigr) \isoto \Mod_{\Lambda(S)}(\D(X,\Lambda))
\end{align}
is an isomorphism. Let $U\colon \Mod_{\Lambda(S)}(\cat A(X,\Lambda)) \to \cat A(X,\Lambda)$ be the forgetful functor and denote $F$ its left adjoint. Note that $U$ and $F$ are exact functors between Grothendieck abelian categories (cf. \cref{rslt:Lambda-sheaves-have-hyperdescent-on-ProFin}). These extend to an adjunction 
\begin{align}
\R F\colon \hat{\D}(\cat A(X,\Lambda)) \rightleftarrows
\hat{\D}\bigl(\Mod_{\Lambda(S)}(\cat A(X,\Lambda))\bigr) \noloc \R U.
\end{align}
Observe that $\R U$ is conservative, because this can be checked on cohomology, where it follows from the fact that $U$ is conservative. Moreover, $\R U$ commutes with colimits: As an exact functor it commutes with finite colimits, and since direct sums in the Grothendieck abelian category $\Mod_{\Lambda(S)}(\cat A(X,\Lambda))$ are exact and $U$ commutes with them, the same applies to $\R U$. Thus, another application of Barr--Beck shows that $\R U$ is monadic. It remains to show that $\R U\comp \R F \simeq \Lambda(S)\tensor\blank$, but this follows from $UF\simeq \Lambda(S)\tensor\blank$ and the fact that $\Lambda(S)$ is a flat $\Lambda$-module. This proves the claim.
\medskip

\textit{Step~3:} Conclusion. Consider the commutative diagram
\begin{equation}
\begin{tikzcd}
\hat{\D}(\cat A(X\times S,\Lambda)) \ar[d] \ar[r,"\sim"] &
\hat{\D}\bigl(\Mod_{\Lambda(S)}(\cat A(X,\Lambda))\bigr) \ar[d,"\sim"] \\
\D(X\times S,\Lambda) \ar[r,"\sim"] & \Mod_{\Lambda(S)}(\D(X,\Lambda)).
\end{tikzcd}
\end{equation}
The bottom horizontal map is an isomorphism by Step~1. By passing to the hearts, we obtain an isomorphism $\cat A(X\times S,\Lambda) \isoto \Mod_{\Lambda(S)}(\cat A(X,\Lambda))$, which shows that the top horizontal map is an isomorphism. The right vertical map is an isomorphism by Step~2. We conclude that the left vertical map is an isomorphism. Hence, $X\times S$ is $\Lambda$-naively derived.
\end{proof}

\begin{proposition} \label{rslt:naively-derived-along-cover}
Let $\Lambda$ be a commutative ring and $f\colon Y \to X$ a cover of condensed anima with associated Čech nerve $f_\bullet\colon Y_\bullet \to X$. Assume that the following conditions are satisfied:
\begin{enumerate}[(a)]
	\item All $Y_n$ are $\Lambda$-naively derived.
	\item The functor $\pi_{1!}\colon \D(Y \times_X Y) \to \D(Y)$ is t-exact, where $\pi_1$ is the projection.
	\item $f$ is $\Lambda$-proper or $\Lambda$-étale.
\end{enumerate}
Then $X$ is $\Lambda$-naively derived.
\end{proposition}
\begin{proof}
Note that we have a commutative diagram
\begin{equation}
\begin{tikzcd}
\hat{\D}(\cat A(X,\Lambda)) \ar[d] \ar[r] & \varprojlim_{n\in\bbDelta} \hat{\D}(\cat A(Y_n,\Lambda)) \ar[d,"\isom"] \\
\D(X,\Lambda) \ar[r,"\isom"'] & \varprojlim_{n\in\bbDelta}\D(Y_n,\Lambda).
\end{tikzcd}
\end{equation}
Here, the right vertical map is an isomorphism by condition (a) and the bottom map is an isomorphism since $\D(\blank,\Lambda)$ is a sheaf. Hence, in order to prove that the left vertical map is an isomorphism, it suffices to show that the top horizontal map is an isomorphism. For the proof we employ \cref{rslt:Grothendieck-descent}. Since each $Y_n$ is $\Lambda$-naively derived, the pullback functor $\alpha^*\colon \hat{\D}(\cat A(Y_n,\Lambda)) \to \hat{\D}(\cat A(Y_m,\Lambda))$ is t-exact and admits a right adjoint $\alpha_*$, for each $\alpha \colon [n]\to [m]$ in $\bbDelta$. By condition (c), the diagram
\begin{equation}
\begin{tikzcd}
\hat{\D}(\cat A(Y_{n+1},\Lambda)) \ar[d,"\alpha'^*"'] \ar[r,"d^0_*"] & \hat{\D}(\cat A(Y_n, \Lambda)) \ar[d,"\alpha^*"] \\
\hat{\D}(\cat A(Y_{m+1}, \Lambda)) \ar[r,"d^0_*"'] & \hat{\D}(\cat A(Y_m,\Lambda))
\end{tikzcd}
\end{equation}
satisfies base-change: If $f$ is $\Lambda$-prim, then this follows from prim base-change (\cref{rslt:prim-base-change}); if on the other hand both $f$ and $\Delta_f$ are $\Lambda$-suave, then all the $\alpha\colon Y_m\to Y_n$ are $\Lambda$-suave and then the claim follows from suave base-change (\cref{rslt:suave-base-change}). It remains to check that $(d^{1*})^{\heartsuit} \colon \cat A(Y,\Lambda) \to \cat A(Y_1,\Lambda)$ sends injective objects to $(d^0_*)^{\heartsuit}$-acyclic objects. Note that we have $Y_1 = Y \times_X Y$ and $d^0 = \pi_1$, $d^1 = \pi_2$ are the two projections $Y \times_X Y \to Y$. If $f$ is $\Lambda$-proper then also $\pi_1$ is $\Lambda$-proper and hence $\pi_{1*} = \pi_{1!}$ is t-exact by (b), which immediately implies the claim. On the other hand, if $f$ is $\Lambda$-étale then $\pi_2$ is $\Lambda$-étale and hence $\pi_2^*$ is right adjoint to $\pi_{2!}$. By (b) $\pi_{2!}$ is t-exact and thus $\pi_2^*$ preserves injective objects, which again proves the claim.
\end{proof}

\begin{remark}
By looking at the proof we note that the assumptions in \cref{rslt:naively-derived-along-cover} can be weakened as follows. One certainly needs to assume (a). Then instead of (b) and (c) it is enough to assume one of the following: Either $f$ is $\Lambda$-prim and $f_*$ is t-exact or both $f$ and $\Delta_f$ are $\Lambda$-suave and $\pi_{1\natural}$ is t-exact (where $\pi_{1\natural}$ is the left adjoint of $\pi_1^*$).
\end{remark}

\begin{proposition}\label{rslt:naively-derived-quotient-by-group-action}
Let $\Lambda$ be a commutative ring, $G$ a locally profinite group and $X$ a $\Lambda$-naively derived condensed anima with an action of $G$. Then $X/G$ is $\Lambda$-naively derived.
\end{proposition}
\begin{proof}
First assume that $G$ is profinite. Then we apply \cref{rslt:naively-derived-along-cover} to the cover $X \to X/G$, where we note that condition (a) is satisfied by \cref{rslt:product-with-profin-set-stays-naively-derived} because all terms in the Čech nerve are of the form $G^n \times X$ (see \cref{rslt:cech-nerve-of-quotient-cover}), condition (b) is satisfied because $\pi_1$ is equivalent to the projection $G \times X \to X$ and hence pushforward along this map is t-exact (see e.g.\ the proof of \cref{rslt:product-with-profin-set-stays-naively-derived}) and condition (c) is satisfied because $f$ is $\Lambda$-proper (this can be checked on a cover of $X/G$ and thus follows from the fact that $\pi_1$ is $\Lambda$-proper).

Now let $G$ be general and let $K \subseteq G$ be a compact open subgroup. We apply \cref{rslt:naively-derived-along-cover} to the $\Lambda$-étale cover $X/K \to X/G$. Denote by $Y_{\bullet} \to X/G$ the associated Čech nerve, so that $Y_n = X/K \times_{X/G} \dotsb \times_{X/G}X/K$ ($n+1$ factors). We now show that the $Y_n$ are $\Lambda$-naively derived. Note that there is a natural map $Y_n \to (*/K)^{n+1} = */K^{n+1}$. By \cref{rslt:group-action-equals-map-to-*/G}, the object $X_n\coloneqq Y_n \times_{*/K^{n+1}} *$ naturally carries a $K^{n+1}$-action such that $X_n/K^{n+1} \simeq Y_n$. By the first part of the proof we are reduced to showing that $X_n$ is $\Lambda$-naively derived. Using $X \simeq X/K \times_{*/K} * \simeq X/G\times_{*/G} *$ and the fact that limits commute with fiber products, we compute
\begin{align}
X_n &= Y_n \times_{*/K^{n+1}} * \\
&\simeq (X/K\times_{*/K}*) \times_{X/G} \dotsb \times_{X/G} (X/K\times_{*/K}*) \\
&\simeq X\times_{X/G} \dotsb \times_{X/G}X \\
&\simeq G^n\times X.
\end{align}
The last isomorphism follows from \cref{rslt:cech-nerve-of-quotient-cover}. Choosing coset representatives for $K$ in $G$, we can write $G^n\times X$ as a disjoint union of copies of $K^n\times X$. Hence, \cref{rslt:product-with-profin-set-stays-naively-derived,rslt:disjoint-union-is-naively-derived} show that $X_n$ is $\Lambda$-naively derived.

In order to apply \cref{rslt:naively-derived-along-cover}, it only remains to verify condition (b) of that result, i.e.\ that $\pi_{1!}$ is t-exact. Now $\pi_1$ corresponds to a map $(G \times X)/K^2 \to X/K$. By base-change, the t-exactness of $\pi_{1!}$ can be checked after pullback along the projection $X \to X/K$, where it becomes the projection $X \times G/K \to X$. But $G/K$ is a discrete set, so the t-exactness of exceptional pushforward along this map reduces to the observation that direct sums are t-exact in $\D(X, \Lambda)$, which is clear.
\end{proof}

\begin{lemma}\label{rslt:modules-over-locally-profinite-set}
Let $\Lambda$ be a commutative ring and $X = \bigdunion_{i\in I}X_i$ a disjoint union of profinite sets. Consider the non-unital $\Lambda$-algebra $\Lambda_c(X) = \bigoplus_{i} \Lambda(X_i)$ and let $\Mod_{\Lambda_c(X)}^{\heartsuit}$ be the abelian category of $\Lambda_c(X)$-modules $M$ satisfying $\Lambda_c(X)M = M$. 

Then $\Lambda_c(X)$ is independent of the decomposition of $X$ and base-change along the projection maps $\Lambda_c(X) \to \Lambda(X_i)$ induces an equivalence of categories
\begin{align}
\D\bigl(\Mod_{\Lambda_c(X)}^{\heartsuit}\bigr) \simeq \prod_{i\in I} \Mod_{\Lambda(X_i)}
\end{align}
with quasi-inverse given by $(M_i)_i \mapsto \bigoplus_i M_i$.
\end{lemma}

\begin{proof}
If $X = \bigdunion_{j\in J} X'_j$ is another decomposition of $X$ into profinite sets, then $X = \bigdunion_{i,j} X_i\cap X'_j$ is a common refinement. For each $i\in I$ there is a finite subset $J_i\subseteq J$ such that $X_i\cap X'_j \neq \emptyset$ implies $j\in J_i$. Thus, $\Lambda(X_i) = \bigoplus_{j\in J} \Lambda(X_i\cap X'_j)$ for all $i\in I$, and hence $\Lambda_c(X) = \bigoplus_{i,j} \Lambda(X_i\cap X'_j)$. The identification $\bigoplus_{i,j} \Lambda(X_i\cap X'_j) = \bigoplus_{j} \Lambda(X'_j)$ follows symmetrically.

Note that the condition $\Lambda_c(X)M = M$ is equivalent to $M = \bigoplus_{i} e_iM$, where $e_i \in \Lambda_c(X)$ is the idempotent corresponding to the unit in $\Lambda(X_i)$. The last statement now follows from $\prod_{i\in I} \Mod_{\Lambda(X_i)} = \D(\prod_{i\in I} \Mod_{\Lambda(X_i)}^{\heartsuit})$.
\end{proof}

\begin{remark}
Let $\Lambda$ be a commutative ring and $X$ a disjoint union of profinite sets. 
\begin{enumerate}[(a)]
\item Consider the $\Lambda$-algebra $\Lambda(X)$ of locally constant functions $X\to \Lambda$ with pointwise addition and multiplication. Then $\Lambda_c(X) \subseteq \Lambda(X)$ is the ideal of locally constant functions $X\to \Lambda$ with compact support.

\item Let $Y$ be another disjoint union of profinite sets and $f\colon Y\to X$ a continuous map. The map $\Lambda(X)\to \Lambda(Y)$ given by pullback along $f$ turns $\Lambda_c(Y)$ into a $\Lambda_c(X)$-module, and the induced functor $f^*\colon \Mod_{\Lambda_c(X)}^{\heartsuit} \to \Mod_{\Lambda_c(Y)}^{\heartsuit}$ is given by $M\mapsto \Lambda_c(Y) \tensor_{\Lambda_c(X)} M$.
\end{enumerate}
\end{remark}

\begin{definition}
Let $\Lambda$ be a commutative ring and $G$ a locally profinite group. A $G$-representation on a $\Lambda$-module $V$ is called \emph{smooth} if the stabilizer of every vector $v\in V$ is open in $G$.

We denote by $\Rep_{\Lambda}(G)^{\heartsuit}$ the Grothendieck abelian category of smooth $G$-representations on $\Lambda$-modules, by $\Rep_{\Lambda}(G)$ its unbounded derived category, and by $\hatRep_{\Lambda}(G)$ the left-completion of $\Rep_{\Lambda}(G)$.

Note that both $\Rep_{\Lambda}(G)$ and $\hatRep_{\Lambda}(G)$ are endowed with a canonical t-structure whose heart is $\Rep_{\Lambda}(G)^{\heartsuit}$, so that the notation is consistent.
\end{definition}

\begin{proposition}\label{rslt:identification-of-sheaves-and-representations}
Let $\Lambda$ be a commutative ring and $G$ a locally profinite group. Then there is a t-exact equivalence $\hatRep_{\Lambda}(G) \isoto \D(*/G,\Lambda)$ of categories, which is natural with respect to continuous group homomorphisms.
\end{proposition}
\begin{proof}
By \cref{ex:profinite-sets-are-naively-derived,rslt:naively-derived-quotient-by-group-action}, we know that $*/G$ is $\Lambda$-naively derived. It therefore suffices to show that the functor $\Rep_{\Lambda}(G)^{\heartsuit} \isoto \cat A(*/G,\Lambda) = \D(*/G,\Lambda)^{\heartsuit}$ is an equivalence. By \cref{rslt:modules-over-locally-profinite-set} we have $\D(G^n, \Lambda) = \D(\Mod_{\Lambda_c(G^n)}^{\heartsuit})$ for all $n\ge0$. The category $\cat A(*/G,\Lambda)$ is the limit of the diagram 
\begin{equation}
\begin{tikzcd}
\Mod_{\Lambda}^{\heartsuit} \ar[r,shift left] \ar[r,shift right] & \Mod_{\Lambda_c(G)}^{\heartsuit} \ar[r,shift left=2,"\pi_2^*"] \ar[r, "m^*" {description, inner sep=0pt}] \ar[r,shift right=2,"\pi_1^*"'] & \Mod_{\Lambda_c(G\times G)}^{\heartsuit},
\end{tikzcd}
\end{equation}
where $\pi_1,\pi_2,m\colon G\times G\to G$ are the projections and multiplication map, respectively. According to \cref{rslt:Grothendieck-descent} the objects of $\cat A(*/G,\Lambda)$ are pairs $(V,\alpha)$ consisting of a $\Lambda$-module $V$ and a $\Lambda_c(G)$-linear isomorphism $\alpha\colon \Lambda_c(G)\tensor_{\Lambda} V \isoto \Lambda_c(G) \tensor_{\Lambda} V$ such that
\begin{align}\label{eq:cocycle-smooth-reps}
\pi_1^*\alpha \comp \pi_2^*\alpha = m^*\alpha
\end{align}
in $\Mod_{\Lambda_c(G\times G)}^{\heartsuit}$. We now construct a functor
\begin{align}\label{eq:identification-of-G-reps}
\Rep_{\Lambda}(G)^{\heartsuit} \to \cat A(*/G,\Lambda),
\end{align}
sending a smooth $G$-representation $(V,\rho)$ to the pair $(V,\alpha_{\rho})$, where $\alpha_{\rho}$ is defined as follows: Identifying $\Lambda_c(G)\tensor_{\Lambda} V$ with the $\Lambda$-module of locally constant functions $G\to V$ with compact support, we define a $\Lambda_c(G)$-linear map $\alpha_{\rho}\colon \Lambda_c(G) \tensor_{\Lambda} V \to \Lambda_c(G)\tensor_{\Lambda} V$ by
\begin{align}
\alpha_{\rho}(f\tensor v)(g) \coloneqq f(g)\cdot \rho(g)v
\end{align}
and $\Lambda$-linear extension (where $f\in \Lambda_c(G)$, $v\in V$ and $g\in G$). We need to prove that $\alpha_{\rho}$ is an isomorphism and that \eqref{eq:cocycle-smooth-reps} holds. For $g,h\in G$ we consider the maps $\varepsilon_g\colon *\to G$ and $\varepsilon_{g,h}\colon *\to G\times G$ given by $g$ and $(g,h)$, respectively. The fact that $\alpha_{\rho}$ is an isomorphism can be checked after applying $\varepsilon_g^*$ for all $g\in G$; but this is clear, since $\varepsilon_g^*\alpha_{\rho} = \rho(g)$ by construction. Similarly, the identity \eqref{eq:cocycle-smooth-reps} can be checked after applying $\varepsilon_{g,h}^*$ for all $(g,h) \in G\times G$, where it follows from the fact that $\rho$ is a group homomorphism. Let now $(V,\rho)$, $(W,\sigma)$ be two smooth $G$-representations and $f\colon V\to W$ a $\Lambda$-linear map. The equality $(\id\tensor f)\comp \alpha_{\rho} = \alpha_{\sigma} \comp (\id\tensor f)$ is equivalent (by applying the $\varepsilon_g^*$) to $f\comp \rho(g) = \sigma(g)\comp f$ for all $g\in G$. It follows that the assignment $(V,\rho)\mapsto (V,\alpha_{\rho})$ defines a fully faithful functor. It remains to prove that it is essentially surjective. So let $(V,\alpha) \in \cat A(*/G,\Lambda)$. Put $\rho(g)\coloneqq \varepsilon_g^*\alpha \colon V\isoto V$ for each $g\in G$. Applying $\varepsilon_{g,h}^*$ to \eqref{eq:cocycle-smooth-reps} yields the identity $\rho(g) \comp \rho(h) = \rho(gh)$. Hence, $(V,\rho)$ is a $G$-representation, and we need to verify that it is smooth. Let $v\in V$ be arbitrary and fix some compact open subgroup $K\subseteq G$. We write $\alpha(\one_K\tensor v) = \sum_{i} f_i\tensor v_i$, where $\one_K$ denotes the characteristic function on $K$. There is a compact open subgroup $H\subseteq K$ such that each $f_i$ is constant on the left cosets of $H$. Then also the orbit map $K\to V$, $k\mapsto \rho(k)v$ is constant on the left cosets of $H$, which implies that $(V,\rho)$ is smooth.

We now address the naturality of the isomorphism $\hatRep_{\Lambda}(G) \isoto \D(*/G,\Lambda)$. So let $f\colon G\to H$ be a homomorphism of locally profinite groups. Since $f^*$ is t-exact, the naturality reduces to checking that the diagram
\begin{equation}
\begin{tikzcd}
\Rep_{\Lambda}(H)^{\heartsuit} \ar[d,"\Res^H_G"'] \ar[r] & \cat A(*/H,\Lambda) \ar[d,"f^*"] \\
\Rep_{\Lambda}(G)^{\heartsuit} \ar[r] & \cat A(*/G,\Lambda)
\end{tikzcd}
\end{equation}
commutes. Given $(V,\rho)\in \Rep_{\Lambda}(H)^{\heartsuit}$, we need to show $\alpha_{\rho\comp f} = f^*\alpha_{\rho}$. As the equality can be checked after applying $\varepsilon_g^*$ for every $g\in G$, it suffices to show $\varepsilon_g^*f^* \simeq \varepsilon_{f(g)}^*$ as functors $\Mod_{\Lambda_c(H)}^{\heartsuit} \to \Mod_{\Lambda}^{\heartsuit}$. But this is clear, as the $\Lambda_c(H)$-module structure on $\Lambda_c(G)$ is given by restricting the natural $\Lambda(G)$-action on $\Lambda_c(G)$ along the composition $\Lambda_c(H) \subseteq \Lambda(H) \to \Lambda(G)$. 
\end{proof}

\begin{corollary}\label{rslt:*-functors-smooth-representations}
Let $\Lambda$ be a commutative ring and $\varphi\colon H\to G$ a homomorphism of locally profinite groups with associated morphism $f\colon */H \to */G$ of classifying stacks. Consider the adjunction $f^*\colon \D(*/G,\Lambda) \rightleftarrows \D(*/H,\Lambda) \noloc f_*$.
\begin{enumerate}[(i)]
\item If $\varphi$ is the inclusion of a closed subgroup, then $f^*$ identifies with restriction of representations, and $f_*$ identifies with the right derived functor $\RInd_H^G$ of smooth induction, where $\Ind_H^G\colon \Rep_{\Lambda}(H)^{\heartsuit} \to \Rep_{\Lambda}(G)^{\heartsuit}$ is explicitly given by
\begin{align}
\Ind_H^G(V,\rho) \coloneqq \set{\tau\colon G\to V}{\begin{array}{l}
\text{there exists an open subgroup $K\subseteq G$ such that} \\
\text{$\tau(gk) = \tau(g)$ for all $k\in K$, $g\in G$, and} \\
\text{$\tau(hg) = \rho(h)\tau(g)$ for all $h\in H$, $g\in G$}
\end{array}};
\end{align}
here the $G$-action is given by right translation: $(g\tau)(g') = \tau(g'g)$.

\item If $\varphi$ is a topological quotient map with kernel $U$, then $f^*$ identifies with inflation of representations, and $f_*$ identifies with the derived invariants functor $(\blank)^U$, that is, the right derived functor of the functor of $U$-invariants $\H^0(U,\blank) \colon \Rep_{\Lambda}(H)^{\heartsuit}\to \Rep_{\Lambda}(G)^{\heartsuit}$, which is given explicitly by
\begin{align}
\H^0(U,(V,\rho)) \coloneqq \set{v\in V}{\text{$\rho(u)v = v$ for all $u\in U$}}
\end{align}
and the induced action of $H/U = G$. For every $n\in\ZZ$ we write
\begin{align}
\H^n(U,V)\coloneqq \H^n(V^U)
\end{align}
for every $V \in \hatRep_{\Lambda}(G)$.
\end{enumerate}
\end{corollary}
\begin{proof}
The identification of $f^*$ was proved in \cref{rslt:identification-of-sheaves-and-representations}. The identification for $f_*$ follows from this and the fact that, by definition, $f_*$ is the right adjoint of $f^*$.
\end{proof}

\begin{remark}\label{rmk:tensor-product-for-representations}
Let $\Lambda$ be a commutative ring and $G$ a locally profinite group. Denote $f\colon *\to */G$ the obvious map. Consider the commutative diagram
\begin{equation}
\begin{tikzcd}
\D(*/G,\Lambda) \times \D(*/G,\Lambda) \ar[d,"f^*\times f^*"'] \ar[r,"\boxtimes"] & \D(*/G\times */G,\Lambda) \ar[d,"{(f\times f)^*}"] \ar[r,"\Delta^*"] & \D(*/G,\Lambda) \ar[d,"f^*"] \\
\Mod_{\Lambda} \times \Mod_{\Lambda} \ar[r,"\tensor_{\Lambda}"'] & \Mod_{\Lambda} \ar[r,equals] & \Mod_{\Lambda}.
\end{tikzcd}
\end{equation}
This implies that the symmetric monoidal structure on $\hatRep_{\Lambda}(G)$ is given under the identification of \cref{rslt:identification-of-sheaves-and-representations} by the underlying tensor product of $\Lambda$-modules with the diagonal $G$-action. One could also define a priori a symmetric monoidal structure on $\hatRep_{\Lambda}(G)$ extending the tensor product on $\Rep_{\Lambda}(G)^{\heartsuit}$ and show that the isomorphism in \cref{rslt:identification-of-sheaves-and-representations} is in fact symmetric monoidal, but we do not pursue this question further.
\end{remark}

\subsection{Cohomological dimension}

For a general locally profinite group $G$ and commutative ring $\Lambda$ the map $*/G \to *$ may not be $\Lambda$-fine, so that $\Rep_\Lambda(G)$ is out of reach of our 6-functor formalism. However, $\Lambda$-fineness is satisfied under mild assumptions on the $\Lambda$-cohomological dimension of $G$. In the present subsection we define the $\Lambda$-cohomological dimension and provide some basic properties of it, using classical representation theoretic methods (see \cref{rmk:how-to-avoid-reptheory-for-fineness} for a completely different approach using condensed mathematics instead).

\begin{definition}
Let $\Lambda$ be a connective $\Einfty$-ring and $G$ a profinite group. We call
\begin{equation}
	\cd_{\Lambda}G \coloneqq \sup\set{n}{\text{$\H^n(G,V) \neq 0$ for some $V\in \Rep_{\Lambda}(G)^{\heartsuit}$}} \in \ZZ_{\ge0} \cup \{\infty\}
\end{equation}
the \emph{$\Lambda$-cohomological dimension} of $G$. We say that a locally profinite group $G$ has \emph{locally finite $\Lambda$-cohomological dimension} if there exists an open profinite subgroup $K\subseteq G$ such that $\cd_{\Lambda}K$ is finite.
\end{definition}

\begin{example}
\begin{enumerate}[(i)]
\item Let $G$ be a $p$-torsionfree compact $p$-adic Lie group and $\Lambda$ an algebra over $\ZZ[n^{-1}]$, where $n$ is the finite part of the pro-order of $G$. Then $G$ has finite $\Lambda$-cohomological dimension: Note that $\cd_{\Lambda}(G) \le \cd_p(G)+1$ in the notation of \cite{Serre.1997} (in view of the fact that the cohomology groups are $p$-power torsion, this follows from \opcit{} Ch.~I, Proposition~13). Hence the claim follows from \cite[Corollaire (1)]{Serre.1965}.
\item Every profinite group whose pro-order is not divisible by $p$ has $\FF_p$-cohomological dimension $0$.
\item The open subgroups of $\GL_n(\FF_p\llbracket t\rrbracket)$ \emph{never} have finite $\FF_p$-cohomological dimension, because they all contain $p$-torsion.
\item Having locally finite $\Lambda$-cohomological dimension is closed under finite products and passage to closed subgroups.
\end{enumerate}
\end{example}

\begin{lemma}\label{rslt:cohomological-dimension}
Let $\Lambda$ be a commutative ring and $G$ a locally profinite group. Then
\begin{align}
\cd_{\Lambda}(G) -1 \le \sup_{\substack{\textnormal{$p$ prime,} \\
p\notin \Lambda^{\times}}} \cd_{\FF_p}(G) \le \cd_{\Lambda}(G).
\end{align}
Moreover, there exists a subring $R\subseteq \QQ$ such that $\Lambda$ is an $R$-algebra and
\begin{align}
\cd_{\Lambda}(G)\le \cd_{R}(G) \le \cd_{\Lambda}(G)+1.
\end{align}
\end{lemma}
\begin{proof}
Let $p$ be a prime with $p\notin \Lambda^{\times}$. Then $\Lambda/p$ is a $\Lambda$-algebra, so that $\cd_{\Lambda/p}(G) \le \cd_{\Lambda}(G)$. But note that we have
\begin{align}
\H^i(G, \Lambda/p \tensor_{\FF_p}V) = \Lambda/p\tensor_{\FF_p}\H^i(G, V)
\end{align}
for all $V\in \Rep_{\FF_p}(G)^{\heartsuit}$ and $i\ge0$, which follows easily from \cite[Ch.~I, Proposition~8]{Serre.1997} by writing $\Lambda/p$ as a filtered colimit of its finite-dimensional $\FF_p$-subvector spaces, which shows $\cd_{\Lambda/p}(G) = \cd_{\FF_p}(G)$. 

It remains to find a prime $p$ with $\cd_{\Lambda}(G) -1 \le \cd_{\Lambda/p}(G)$. Let $V \in \Rep_{\Lambda}(G)^{\heartsuit}$ be such that $\H^{d}(G,V) \neq 0$, where $d \le \cd_{\Lambda}(G)$. As the higher cohomology groups are torsion (see Corollary~3 in \loccit), we find a prime $p$ such that multiplication by $p$ induces a non-injective endomorphism on $\H^d(G,V)$. But note that we have short exact sequences $0\to V[p] \to V \xto{p\cdot} pV \to 0$ and $0\to pV\to V\to V/pV \to 0$ and hence exact sequences
\begin{align}
\H^d(G,V[p]) \to \H^d(G,V) \xto{\alpha} \H^d(G,pV) \qquad \text{and}\qquad \H^{d-1}(G,V/pV) \to \H^d(G,pV) \xto{\beta} \H^d(G,V).
\end{align}
As $\beta\circ\alpha$ is given by multiplication with $p$, either $\alpha$ or $\beta$ is not injective. But this means either $H^d(G,V[p]) \neq 0$ or $\H^{d-1}(G,V/pV) \neq 0$ or, in other words, $\cd_{\Lambda/p}(G) \ge d-1$.

The last assertion holds for the subring $R\subseteq \QQ$ which is obtained from $\ZZ$ by inverting all primes $p$ with $p\in \Lambda^\times$. Then $\Lambda$ is an $R$-algebra (whence $\cd_{\Lambda}(G) \le \cd_{R}(G)$) and the second estimate follows by applying the lemma to $R$ noting that $\set{p}{p\notin \Lambda^{\times}} = \set{p}{p\notin R^{\times}}$.
\end{proof}

\begin{lemma}\label{rslt:vanishing-Ext}
Let $\Lambda \subseteq \QQ$ be subring and $G$ a profinite group with $\cd_{\Lambda}(G) < \infty$. For all $V,W \in \Rep_{\Lambda}(G)^{\heartsuit}$ and $i>\cd_{\Lambda}(G)+1$ we have
\begin{align}
\Ext_G^i(V,W) = 0.
\end{align}
\end{lemma}
\begin{proof}
This follows from \cite[Corollary 6(iii)]{Martinez.1973}, which says that the injective dimension of each $W\in \Rep_{\Lambda}(G)^{\heartsuit}$ is bounded by $\cd_{\Lambda}(G)+1$; note that our $\cd_{\Lambda}(G)$ coincides with the strict cohomological dimension $\operatorname{scd}(G)$ of \loccit 

For the reader's convenience we sketch the outline of the argument, which is based on the following facts: 
\begin{fact*}
Recall that a representation $A$ of a finite group $\Gamma$ is called \emph{cohomologically trivial} if $\hat{\H}^i(\Gamma', A) = 0$ for all $i\in \ZZ$ and all subgroups $\Gamma'\subseteq \Gamma$, where $\hat{\H}^i$ denote the Tate cohomology groups.
\begin{enumerate}[(1)]
\item \cite[IX, \S6, Théorème 10]{Serre.1980} Let $A$ be a cohomologically trivial $\Gamma$-representation whose underlying $\ZZ$-module is divisible. Then $A$ is an injective $\Gamma$-representation.

\item \cite[Appendix~1 to Ch.~I, Lemma~1]{Serre.1997} Let $G$ be a profinite group and $V \in \Rep_{\ZZ}(G)^{\heartsuit}$. Put
\begin{align}
\cd(G,V) \coloneqq \sup\set{r\ge0}{\text{$\H^r(H,V) \neq 0$ for all closed subgroups $H\subseteq G$}}.
\end{align}
Then $\cd(G,V) = 0$ if and only if for all open normal subgroups $U\subseteq G$ the $G/U$-representation $\H^0(U,V)$ is cohomologically trivial.
\end{enumerate}
\end{fact*}

The category $\Rep_{\ZZ}(G)^{\heartsuit}$ (which contains $\Rep_{\Lambda}(G)^{\heartsuit}$ as a full subcategory) is Grothendieck abelian and hence admits injective hulls (this is well-known, but see \cite[Théorème 2, p.~362]{Gabriel.1962} for a proof). We now explain the argument in a series of steps.
\medskip

\textit{Step~1:} If $I\in \Rep_{\Lambda}(G)^{\heartsuit}$ is injective, then $\H^0(H,I)$ is an injective $G/H$-representation for every closed normal subgroup $H\subseteq G$. This follows from the observation that $\H^0(H,\blank)$ admits an exact left adjoint (namely inflation) and hence preserves injectivity.
\medskip

\textit{Step~2:} If $U\subseteq G$ is an open normal subgroup and $A \in \Rep_{\Lambda}(G/U)^{\heartsuit}$ is injective, then the underlying $\Lambda$-module of $A$ is divisible. As $\Res^{G/U}_1$ admits $\cInd_1^{G/U}$ as an exact left adjoint, it preserves injectivity. Hence the underlying $\Lambda$-module of $A$ is injective or, equivalently, divisible.
\medskip

\textit{Step~3:} Let $W\in \Rep_{\Lambda}(G)^{\heartsuit}$ and consider a resolution
\begin{align}
0\to W \to I_0 \xto{d_0} I_1 \to \dotsb \to I_{d-1}\xto{d_{d-1}} I_d \xto{\varepsilon} J \to 0
\end{align}
in $\Rep_{\Lambda}(G)^{\heartsuit}$, where each $I_i$ is injective and $d\coloneqq \cd_{\Lambda}(G) \ge \cd(G,W)$. We need to show that $J$ is injective. 

Let $U\subseteq G$ be an arbitrary but fixed open normal subgroup. We first show that $\H^0(U,J) \in \Rep_{\ZZ}(G/U)^{\heartsuit}$ injective. By Facts (1) and (2) it suffices to show:
\begin{enumerate}[(i)]
\item $\cd(G,J) = 0$, and
\item $\H^0(U,J)$ is divisible.
\end{enumerate}
We first prove (i). A dimension shifting argument shows
\begin{align}
\H^i(H,J) = \H^{i+1}(H,\Ker(\varepsilon)) = \H^{i+2}(H,\Ker(d_{d-1})) = \dotsb = \H^{i+d}(H, \Ker(d_1)) = 0
\end{align}
for all $i>0$ and all closed subgroups $H\subseteq G$, because $\cd(G,\Ker(d_1)) \le d$. This shows (i).

For (ii) one argues as follows: Note that the same argument as above shows $\cd(G,\Ker(\varepsilon)) = 0$. Applying cohomology to the short exact sequence $0\to \Ker(\varepsilon) \to I_d \to J\to 0$, we obtain an exact sequence
\begin{align}
\H^0(U, I_d) \to \H^0(U, J) \to \H^1(U, \Ker(\varepsilon)) = 0.
\end{align}
Now, $\H^0(U, I_d)$ is injective by Step~1 and hence divisible by Step~2. Consequently $\H^0(U,J)$ is divisible as a quotient of a divisible abelian group.

We deduce that $\H^0(U,J)$ injective in $\Rep_{\ZZ}(G/U)^{\heartsuit}$. Consider now an exact sequence $0\to J\to I \to Q \to 0$ in $\Rep_{\ZZ}(G)^{\heartsuit}$, where $I$ is an injective hull of $J$. Since $\cd(G,J) = 0$ by (i), we obtain an exact sequence
\begin{align}
0\to \H^0(U,J) \to \H^0(U,I) \to \H^0(U,Q) \to 0.
\end{align}
By Step~1, $\H^0(U,I)$ is an injective $G/U$-representation and is clearly an essential extension of $\H^0(U,J)$. But as $\H^0(U,J)$ is injective, the short exact sequence splits, from which we deduce $\H^0(U,Q) = 0$. As $U$ was arbitrary, it follows that $Q = \bigcup_U \H^0(U,Q) = 0$. Therefore, $J = I$ is injective in $\Rep_{\ZZ}(G)^{\heartsuit}$, hence also in $\Rep_{\Lambda}(G)^{\heartsuit}$.
\end{proof}

The above results accumulate in the following result, which is the key technical input we need in order to construct the shriek functors on $*/G$. Recall the definition of descendable objects from \cref{def:descendability}.

\begin{proposition} \label{rslt:descendability-of-regular-representation}
Let $\Lambda$ be a commutative ring and $G$ a profinite group such that $\cd_\Lambda G < \infty$. Then the $G$-representation $\Cont(G, \Lambda)$ is descendable in $\hatRep_\Lambda(G)$.
\end{proposition}
\begin{proof}
We first assume that $\Lambda$ is a subring of $\QQ$. Note that under the equivalence from \cref{rslt:identification-of-sheaves-and-representations}, $\Cont(G, \Lambda)$ identifies with $g_* \one$, where $g\colon * \to */G$ is the canonical map; note that $g$ is $\Lambda$-fine and even $\Lambda$-proper, because its base-change to $*$ is given by $G \to *$ (see \cref{rslt:suave-and-prim-map-is-local-on-target}). Let $g_{\bullet} \colon G^{\bullet}\to */G$ be the Čech nerve of $g\colon *\to */G$. From the equivalence $\D(*/G,\Lambda) \isoto \varprojlim_{n\in\bbDelta} \D(G^n,\Lambda)$ we obtain
\begin{align}
	\one = \varprojlim_{n\in\bbDelta}g_{n*}\one = \Tot(\Cont(G, \Lambda)^{\tensor(\bullet+1)}),
\end{align}
where $g_{n*}\one = \Cont(G,\Lambda)^{\tensor(n+1)}$ follows by induction on $n$ from base-change and the projection formula. Put $d\coloneqq \cd_\Lambda G$ and denote $\Tot_d$ the limit over $\bbDelta_{\le d}$. Denote $X$ the cofiber of the natural map $\one = \Tot(\Cont(G,\Lambda)^{\tensor(\bullet+1)}) \to \Tot_d(\Cont(G,\Lambda)^{\tensor(\bullet+1)})$. For any $i\ge0$ we get an exact sequence
\begin{align}
\H^i(\Tot(\Cont(G,\Lambda)^{\tensor(\bullet+1)})) \xto{\alpha^i} \H^i(\Tot_d(\Cont(G,\Lambda)^{\tensor(\bullet+1)})) \to \H^i(X) \to \H^{i+1}(\one) = 0.
\end{align}
By (the dual of) \cite[Proposition~1.2.4.5]{HA} the map $\alpha^i$ is an isomorphism for $i<d$. We conclude that $\H^i(X) = 0$ for $i<d$, and hence $X \in \D^{\ge d}(*/G,\Lambda)$. By \cref{rslt:vanishing-Ext} and a spectral sequence argument, we compute $\Hom(X,\one[1]) = 0$. But this means that the (co)fiber sequence
\begin{align}
\one \to \Tot_d(\Lambda(G)^{\tensor(\bullet+1)}) \to X
\end{align}
splits, which shows $\Tot_d(\Cont(G,\Lambda)^{\tensor(\bullet+1)}) = \one \oplus X$. Hence $\one$ is a retract of a finite limit built out of elements of the form $\Cont(G,\Lambda)^{\tensor n}$, proving that $\Cont(G,\Lambda)$ is indeed descendable.

Now let $\Lambda$ be general. By \cref{rslt:cohomological-dimension} there is a subring $R \subseteq \QQ$ such that $\Lambda$ is an $R$-algebra and $\cd_R G < \infty$. By the above argument we deduce that $\Cont(G,R)$ is descendable in $\hatRep_R(G)$. Moreover, by \cref{rslt:functoriality-of-6ff-on-Cond-Ani-in-Lambda} the ring map $R \to \Lambda$ induces a symmetric monoidal exact functor $\blank \tensor_R \Lambda\colon \hatRep_R(G) \to \hatRep_\Lambda(G)$ sending $\Cont(G, R)$ to $\Cont(G, \Lambda)$. This shows that $\Cont(G,\Lambda)$ is descendable, as desired.
\end{proof}

\begin{remark} \label{rmk:how-to-avoid-reptheory-for-fineness}
\Cref{rslt:descendability-of-regular-representation} is the main ingredient for our proof that $*/G \to *$ is $\Lambda$-fine in the next subsection. This result relies on \cref{rslt:vanishing-Ext}, which in turn relies on non-trivial results in classical representation theory. A natural question is if this reliance can be avoided, which is indeed the case in many cases:
\begin{enumerate}[(a)]
	\item Let $G$ be a pro-$p$ group and $\Lambda$ a $\ZZ[1/p]$-algebra. The fact that there exists a  $\Lambda$-valued right Haar measure $\mu$ on $G$ with $\mu(G) = 1$ means that $\one$ is a direct summand of $\Cont(G,\Lambda)$, which immediately implies that $\Cont(G,\Lambda)$ is descendable. Concretely, we can write $\Cont(G,\Lambda) = \varinjlim_{K} \Cont(G/K,\Lambda)$, where $K\subseteq G$ runs through the open subgroups. The inclusion $\one \to \Cont(G/K,\Lambda)$ into the constant functions admits a retract given explicitly by the assignment $f\mapsto [G:K]^{-1} \sum_{g\in G/K} f(g)$. These retracts assemble into a retract for $\one \injto \Cont(G,\Lambda)$.

	\item Let $G$ be a profinite group with $\cd_{\FF_p}(G) < \infty$. In this case one can show \cref{rslt:vanishing-Ext} using a completely different method involving condensed mathematics, see \cite[Lemma~3.11]{Hansen-Mann:Mod-p-Stacky-6FF}. The basic idea is as follows: Fix static $G$-representations $V$ and $W$ and consider $\RHom_G(V, W) \in \Mod_\Lambda$. Our goal is to show that $\RHom_G(V, W)$ is right bounded with a bound depending only on $\cd_{\FF_p}(G)$. We now observe that
	\begin{align}
		\RHom_G(V, W) = \iHom_G(V, W)^G,
	\end{align}
	where $\iHom(V, W)$ denotes the $G$-representation on the $\FF_p$-vector space $\iHom_{\FF_p}(V, W)$ and $(\blank)^G$ denotes $G$-cohomology. Clearly $\iHom_G(V, W)$ is static, so one is tempted to conclude using the finite $\FF_p$-cohomological dimension of $G$. However, $\iHom_G(V, W)$ is in general not smooth anymore, so we have no a priori control over its $G$-cohomology. This can be fixed by equipping $\iHom_G(V, W)$ with its natural topology (more precisely, its natural condensed structure) and observing that the finite $G$-cohomological dimension extends from smooth (i.e.\ discrete) representations to nice enough topological representations; more precisely it extends to all \emph{solid} $G$-representations.
\end{enumerate}
\end{remark}

\subsection{Six functors on the classifying stack}

In \cref{rslt:*-functors-smooth-representations} we have seen that the $*$-functors identify with the classical functors (restriction/inflation and smooth induction/cohomology) in representation theory. The goal of this section is to show that for locally profinite groups with bounded cohomological dimension the $!$-functors exist and to determine the suave and prim sheaves on $*/G$. Finally, we will show that for a $p$-adic Lie group $G$ the structure map $*/G\to *$ is $\Lambda$-suave.

\begin{remark}\label{rmk:etale-proper-maps-classifying-stacks}
Let $\Lambda$ be an $\Einfty$-ring. In the following, we will often make use of the following simple facts: Let $G$ be a locally profinite group and let $H\subseteq K\subseteq G$ be compact subgroups with $K$ open and $H$ closed.
\begin{enumerate}[(a)]
\item The map $*/K\to */G$ is $\Lambda$-étale. Indeed, by \cref{rslt:stability-of-etale-and-proper-maps}, this can be checked after pullback along the cover $*\to */G$, where the claim becomes that $G/K \to *$ is $\Lambda$-étale. But $G/K$ is a discrete set, and hence the claim follows from \cref{rslt:open-local-properties-for-condensed-anima}.

\item The map $*/H \to */K$ is $\Lambda$-proper: As before, this can be checked after pullback along the cover $*\to */K$, where the claim becomes that $K/H \to *$ is $\Lambda$-proper. As $K/H$ is a profinite set, the claim follows from \cref{rslt:maps-from-profinite-sets-are-proper}.
\end{enumerate}
\end{remark}

\begin{proposition}\label{rslt:BG-is-proper}
Let $\Lambda$ be a commutative ring and $G$ a profinite group with $\cd_{\Lambda}(G) < \infty$. Then $*/G\to *$ is $\Lambda$-proper (and in particular $\Lambda$-fine).
\end{proposition}
\begin{proof}
Following the strategy in \cite[Lemma~3.11]{Hansen-Mann:Mod-p-Stacky-6FF}, we want to apply \cref{rslt:check-properness-on-!-descendable-cover} to $f\colon */G\to *$ and $g\colon *\to */G$. From \cref{rmk:etale-proper-maps-classifying-stacks} we know that $g$ is $\Lambda$-prim (even $\Lambda$-proper); moreover, it is clear that $f$ is truncated and $fg = \id_*$ is $\Lambda$-proper. To conclude, it remains to show that $g_*\one = C(G, \Lambda)$ is descendable. But this was shown in \cref{rslt:descendability-of-regular-representation}.
\end{proof}

\begin{corollary}\label{rslt:maps-of-class-stacks-are-fine}
Let $\Lambda$ be a commutative ring and $H\to G$ a morphism of locally profinite groups which have locally finite $\Lambda$-cohomological dimension. Then the induced map $*/H\to */G$ is $\Lambda$-fine.
\end{corollary}
\begin{proof}
Since the property of being $\Lambda$-fine is cancellative, it suffices to show that the map $*/G\to *$ is $\Lambda$-fine. This can be checked after restriction to $*/K$ for some compact open subgroup $K \subseteq G$ with $\cd_{\Lambda}K < \infty$ (combine \cref{rmk:etale-proper-maps-classifying-stacks} and \cref{rslt:!-descent-along-suave-map}). Now the claim follows from \cref{rslt:BG-is-proper}.
\end{proof}

\begin{corollary}\label{rslt:finite-coh-dim-implies-unit-compact}
Let $\Lambda$ be a commutative ring and $G$ a profinite group with $\cd_{\Lambda}(G)< \infty$. Then $\one \in \D(*/G,\Lambda)$ is compact.
\end{corollary}
\begin{proof}
The functor
\begin{equation}
\RHom_{\D(*/G,\Lambda)}(\one, \blank) = f_*\iHom(\one,\blank) = f_* = f_!
\end{equation}
commutes with colimits by \cref{rslt:BG-is-proper}. 
\end{proof}

Without the assumption on cohomological dimension, the pushforward along $*/G \to *$ still commutes with filtered colimits of uniformly left-bounded representations, which we believe to be of independent interest:

\begin{lemma}\label{rslt:coh-of-profin-group-commutes-with-filtered-colim}
Let $\Lambda$ be a connective $\Einfty$-ring and $G$ a profinite group. Denote $f\colon */G\to *$ the tautological map. Then $f_*\colon \D(*/G,\Lambda)^{\ge0} \to \Mod_{\Lambda}$ commutes with filtered colimits.
\end{lemma}
\begin{proof}
Denote $q_{\bullet} \colon G^\bullet \to */G$ the Čech nerve of $*\to */G$. The functor $f_*$ is then computed by
\begin{align}
\D(*/G,\Lambda) \simeq \varprojlim_{n\in\bbDelta} \D(G^n,\Lambda) &\to \Mod_{\Lambda}, \\
M &\mapsto \Tot\bigl( (fq_\bullet)_*q_\bullet^*M\bigr).
\end{align}
Note that the functors $(fq_n)_*$ and $q_n^*$ commute with colimits and are t-exact (for the first functor we use that $G^n$ is profinite).

Let now $(M_i)_{i\in I}$ be a filtered diagram in $\D(*/G,\Lambda)^{\ge0}$. For any $d$, denoting $\Tot_d$ the limit over $\bbDelta_{\le d}$, we have a commutative diagram
\begin{equation}
\begin{tikzcd}
\varinjlim_i \Tot\bigl((fq_{\bullet})_*q_{\bullet}^*M_i\bigr) \ar[d] \ar[r] & \Tot\bigl( \varinjlim_i (fq_{\bullet})_*q_{\bullet}^*M_i\bigr) \ar[r,"\sim"] \ar[d] & \Tot\bigl( (fq_{\bullet})_*q_{\bullet}^* \varinjlim_i M_i\bigr) \ar[d] \\
\varinjlim_i \Tot_{d}\bigl((fq_{\bullet})_*q_{\bullet}^*M_i\bigr) \ar[r,"\sim"] & \Tot_d\bigl( \varinjlim_i (fq_{\bullet})_*q_{\bullet}^*M_i\bigr) \ar[r,"\sim"] & \Tot_{d}\bigl( (fq_{\bullet})_*q_{\bullet}^* \varinjlim_i M_i\bigr).
\end{tikzcd}
\end{equation}
We need to show that the upper left horizontal map is an isomorphism, which can be checked on cohomology. The lower left horizontal map is an isomorphism, because $\Tot_d$ is a finite limit and hence commutes with filtered colimits. By (the dual of) \cite[Proposition~1.2.4.5]{HA}, the vertical maps become isomorphisms after applying $\H^n$, for all $n<d$. It follows that the top left horizontal map becomes an isomorphism after applying $\H^n$, for all $n<d$. Varying $d$, this proves the claim.
\end{proof}

\begin{notation}
Let $\Lambda$ be a commutative ring and $G$ a locally profinite group with structure map $f\colon */G\to *$. For any $V \in \D(*/G,\Lambda)$ we write
\begin{align}
V^G \coloneqq \Gamma(*/G,V) = f_*V
\end{align}
for the derived invariants of $V$. Beware that, if the pro-order of $G$ is not invertible in $\Lambda$, then the functor $V\mapsto V^G$ is generally not t-exact and hence should not be confused with the underived functor of invariants.
\end{notation}

A classical result in smooth representation theory is that a smooth $G$-representation can be recovered as the filtered colimit over its invariants under compact open subgroups $K \subseteq G$. We recover this result (over general coefficient rings and with derived invariants) from basic properties of the 6-functor formalism:

\begin{lemma}\label{rslt:colimit-of-invariants}
Let $\Lambda$ be a commutative ring. Let $G$ be a locally profinite group which has locally finite $\Lambda$-cohomological dimension. Then for every $V \in \D(*/G,\Lambda)$, the natural map
\begin{align}
\varinjlim_{K} V^K \isoto V
\end{align}
is an isomorphism in $\Mod_{\Lambda}$, where on the right-hand side we implicitly pass to the underlying $\Lambda$-module via pullback along $q\colon *\to */G$. Here the colimit runs over the compact open subgroups $K\subseteq G$ with $\cd_{\Lambda}(K) < \infty$.
\end{lemma}
\begin{proof}
For each compact open subgroup $K\subseteq G$ we write $i_K\colon */K\to */G$ and $f_K\colon */K\to *$ for the obvious maps, so that $(\blank)^{K} = f_{K*}i_{K}^*$. 

Restricting to a compact open subgroup if necessary, we may assume from the start that $G$ is profinite and $\cd_{\Lambda}G <\infty$. Then $\cd_{\Lambda}(K) = \cd_{\Lambda}G < \infty$ for all open subgroups $K\subseteq G$. We will construct a natural (in $V$) isomorphism
\begin{align}\label{eq:colimit-of-inductions}
\varinjlim_{K} i_{K*}i_K^*V \isoto q_*q^*V
\end{align}
in $\D(*/G,\Lambda)$. Applying $f_*$, which by \cref{rslt:BG-is-proper} commutes with colimits, to \eqref{eq:colimit-of-inductions} then implies the claim.

Similarly as in the proof of \cref{rslt:cohomology-is-sheafy} (but working in $\bigl(\Cond(\Ani)_{/(*/G)}\big)^{\op}$ and replacing $\Mod_{\Lambda}$ with $\D(*/G,\Lambda)$), we obtain a functor 
\begin{align}\label{eq:functor-of-inductions}
\bigl(\Cond(\Ani)_{/(*/G)}\bigr)^{\op} &\to \Fun(\D(*/G,\Lambda), \D(*/G,\Lambda))
\end{align}
sending a condensed anima $X\xrightarrow{\alpha_X} */G$ to $\alpha_{X*}\alpha_X^*$. Let $\Omega$ be the poset of compact open subgroups $K\subseteq G$, ordered by inclusion. Then $\Omega^{\triangleleft}$ embeds into $\Cond(\Ani)_{/(*/G)}$ via $K\mapsto */K$ and $\emptyset \mapsto *$. Restricting \eqref{eq:functor-of-inductions} along $(\Omega^{\triangleleft})^{\op}$ then gives the map \eqref{eq:colimit-of-inductions}. We will check that this map is an isomorphism, which we can do after applying the conservative functor $q^*$. We contemplate the diagram
\begin{equation}
\begin{tikzcd}
G \ar[d] \ar[r] & G/K \ar[d] \ar[r] & * \ar[d,"q"] \\
* \ar[dr,"q"'] \ar[r] & */K \ar[d,"i_K"] \ar[r,"i_K"] & */G \\
& */G
\end{tikzcd}
\end{equation}
in which the squares are cartesian. Note that all maps are $\Lambda$-proper by \cref{rmk:etale-proper-maps-classifying-stacks} and hence satisfy base-change and the projection formula. The pullback of \eqref{eq:colimit-of-inductions} along $q^*$ thus becomes
\begin{align}
\varinjlim_{K} \Lambda(G/K)\tensor_{\Lambda} q^*V \isoto \Lambda(G)\tensor_{\Lambda}q^*V,
\end{align}
which is an isomorphism since $G = \varprojlim_K G/K$ (see also \cref{def:Einfty-ring-associated-to-profin-set}). 
\end{proof}

As a corollary of \cref{rslt:colimit-of-invariants} we can show that if $G$ has locally finite $\Lambda$-cohomological dimension then the distinction between $\Rep_\Lambda(G)$ and $\hatRep_\Lambda(G)$ vanishes:

\begin{proposition} \label{rslt:Rep-G-left-complete}
Let $G$ be a locally profinite group and $\Lambda$ a commutative ring such that $G$ has locally finite $\Lambda$-cohomological dimension. Then $\Rep_\Lambda(G)$ is left-complete, i.e.\ we have
\begin{align}
	\Rep_\Lambda(G) = \hatRep_\Lambda(G) = \D(*/G,\Lambda).
\end{align}
\end{proposition}
\begin{proof}
It is enough to show that countable products in $\Rep_\Lambda(G)$ have finite cohomological dimension, then the claim follows from \cite[Proposition~1.2.1.19]{HA} (more precisely, \loccit{} requires that countable products have cohomological dimension $0$, but the proof works the same with a positive finite cohomological dimension); note that condition (2) in \loccit{} is satisfied, because $\Rep_{\Lambda}(G)$ is a derived category.

Now let $(V_n)_n$ be a countable family of objects in $\Rep_\Lambda(G)^\heartsuit$; we need to show that $\prod_n V_n \in \Rep_\Lambda(G)$ is right-bounded, with a bound independent of $(V_n)_n$. The boundedness can be checked on the underlying $\Lambda$-modules, so by \cref{rslt:colimit-of-invariants} we have
\begin{align}
	\prod_n V_n = \varinjlim_K \Bigl(\prod_n V_n\Bigr)^K = \varinjlim_K \prod_n V_n^K.
\end{align}
Here $K$ ranges over compact open subgroups of $G$ with finite $\Lambda$-cohomological dimension. Moreover, fixing one such $K$, say $K_0$, we have $\cd_\Lambda K \le \cd_\Lambda K_0$ for all $K \subseteq K_0$. This shows that $\prod_n V_n$ is bounded by $\cd_\Lambda K_0$, as desired.
\end{proof}

We have now established the 6-functor formalism on $*/G$, i.e.\ for smooth $G$-representations. In the following we will apply the results from \cref{sec:kerncat} to study the suave and prim objects in $*/G$ and in particular describe them in more classical terms. We first determine the dualizable objects in $\Rep_\Lambda(G)$ (recall the notion of being dualizable from \cref{def:dualizable-object}):

\begin{lemma}\label{rslt:pullback-preserves-dualizability}
Let $\Lambda$ be a commutative ring, $G$ a locally profinite group and $K\subseteq G$ a compact open subgroup with $\cd_{\Lambda}(K) < \infty$. For an object $V \in \D(*/G,\Lambda)$ the following are equivalent:
\begin{enumerate}[(a)]
\item $V$ is dualizable; 
\item $i_K^*V$ is dualizable in $\D(*/K,\Lambda)$, where $i_K\colon */K\to */G$ is the obvious map;
\item the underlying $\Lambda$-module of $V$ is perfect (i.e.\ dualizable in $\Mod_{\Lambda}$).
\end{enumerate}
\end{lemma}
\begin{proof} 
The implication \enquote{(a)$\implies$(c)} follows from the fact that the forgetful functor $\D(*/G,\Lambda) \to \Mod_{\Lambda}$ is symmetric monoidal.

We now prove that (b) implies (a). By \cref{rslt:properties-of-dualizable-objects} it suffices to show that the canonical map $V\tensor V^{\vee} \to \iHom(V,V)$ is a $G$-equivariant isomorphism, where $V^{\vee} = \iHom(V,\one)$ denotes the smooth dual of $V$. As $i_K^*$ is conservative, it suffices to check that the top horizontal map in the commutative diagram
\begin{equation}
\begin{tikzcd}
i_K^*(V\tensor V^\vee) \ar[d,equals] \ar[rr] && i_K^*\iHom(V,V) \ar[d,"\sim"] \\
i_K^*V \tensor i_K^*(V^\vee) \ar[r,"\sim"] & i_K^*V \tensor (i_K^*V)^\vee \ar[r,"\sim"] & \iHom(i_K^*V, i_K^*V),
\end{tikzcd}
\end{equation} 
is an isomorphism. Observe that, since $i_{K!}$ is left adjoint to $i_K^*$, it follows that the projection formula $i_{K!}\comp (i_K^*V\tensor \blank) \isoto V\tensor i_{K!}(\blank)$ gives rise to an isomorphism
\begin{align}\label{eq:pullback-preserves-dualizability}
i_K^*\iHom(V,W) \isoto \iHom(i_K^*V, i_K^*W),
\end{align}
for all $W\in \D(*/G,\Lambda)$, by passing to the right adjoints. This implies that the lower left horizontal and the right vertical maps are isomorphisms. The lower right horizontal map is an isomorphism by the assumption (b) (together with \cref{rslt:properties-of-dualizable-objects}). It follows that the top horizontal map is an isomorphism, whence (a).

Finally, we prove that (c) implies (b). Since $f_K\colon */K\to *$ is $\Lambda$-proper, \cref{rslt:etale-means-suave=dualizable} shows that the $f_K$-prim and the dualizable objects in $\D(*/K,\Lambda)$ agree. We are thus reduced to prove the following: An object $V \in \D(*/K,\Lambda)$ is $f_K$-prim as soon as $q^*V$ is prim (i.e.\ dualizable), where $q\colon *\to */K$ is the obvious map. Now we note that $q$ is $\Lambda$-prim by \cref{rmk:etale-proper-maps-classifying-stacks} and $q_*\one$ is descendable by the proof of \cref{rslt:BG-is-proper}. Hence, the claim follows from \cref{rslt:check-prim-object-on-!-descendable-cover}.
\end{proof}

We can now determine the suave and prim objects in $\D(*/G,\Lambda)$. The suave objects will turn out to be exactly the \emph{admissible} representations, for the following (derived) notion of admissibility:

\begin{definition} \label{def:admissible-representation}
Let $\Lambda$ be commutative ring and $G$ a locally profinite group which has locally finite $\Lambda$-cohomological dimension. An object $V \in \D(*/G,\Lambda)$ is called \emph{admissible} if $V^K \in \Mod_{\Lambda}$ is dualizable for all compact open subgroups $K\subseteq G$ with $\cd_{\Lambda}K < \infty$.
\end{definition}

\begin{remark}
Let $\Lambda$ be a commutative ring and $G$ a locally profinite group. In each of the following cases $G$ has locally finite $\Lambda$-cohomological dimension, hence the definition of admissibility applies.
\begin{remarksenum}
\item Assume that $\Lambda$ is regular and that the pro-order of $G$ is invertible in $\Lambda$. Then every compact open subgroup $K\subseteq G$ satisfies $\cd_{\Lambda}K = 0$ so that the functor $V\mapsto V^K$ is t-exact. Hence, in view of \cite[\href{https://stacks.math.columbia.edu/tag/066Z}{Tag 066Z}]{stacks-project}, a complex of representations $V \in \Rep_{\Lambda}(G)$ is admissible if and only if it satisfies the following conditions:
\begin{itemize}
\item all cohomologies $\H^i(V)$ are admissible in the sense that $\H^i(V)^K = \H^i(V^K)$ is a perfect $\Lambda$-module for all compact open subgroups $K\subseteq G$; 
\item for every compact open subgroup $K\subseteq G$ we have $\H^i(V)^K \neq 0$ for only finitely many $i\in \ZZ$.
\end{itemize}

\item Assume that $G$ is a $p$-adic Lie group and $\Lambda$ is a field of characteristic $p$. For $V\in \Rep_{\Lambda}(G)$ the following are equivalent:
\begin{enumerate}[(a)]
\item $V$ is admissible;
\item $V^I$ is dualizable for \emph{some} torsionfree open pro-$p$ subgroup $I\subseteq G$;
\item $V$ is bounded and all cohomologies $\H^i(V)$ are admissible smooth $G$-representations in the classical sense, i.e.\ for some (hence every) torsionfree open pro-$p$ subgroup $I\subseteq G$ the space $\H^0(I,\H^i(V))$ is finite-dimensional. 
\end{enumerate}
Clearly, (b) is a special case of (a). Conversely, note that for any torsionfree open pro-$p$ subgroups $U\subseteq I\subseteq G$, it follows from \cite[Proposition~6]{Schneider.2015} that $\one$ and $\cInd_U^I\one$ generate the same full subcategory of $\D(*/I,\Lambda)$ under (co)fibers and retracts. Thus, $V^U = \RHom_I(\cInd_U^I\one, V)$ is dualizable if and only if $V^I = \RHom_I(\one, V)$ is dualizable. We deduce that (b) implies (a). For the equivalence of (b) and (c) we refer to \cite[Corollary~4.12]{Schneider-Sorensen.2023a}.
\end{remarksenum}
\end{remark}

\begin{proposition} \label{rslt:identify-suave-and-prim-representations}
Let $\Lambda$ be commutative ring and $G$ a locally profinite group which has locally finite $\Lambda$-cohomological dimension and denote $f\colon */G\to *$ the structure map. An object of $\D(*/G,\Lambda)$ is:
\begin{propenum}
\item\label{rslt:prim=compact} $f$-prim if and only if it is compact.
\item\label{rslt:suave=admissible} $f$-suave if and only if it is admissible.
\end{propenum}
\end{proposition}
\begin{proof}
We first prove (i). Observe that $\Lambda$ is compact in $\Mod_{\Lambda}$. Then \cref{rslt:relation-between-suave-prim-and-compact} shows that every $f$-prim object is compact in $\D(*/G,\Lambda)$. For the converse, we make the following claim:
\begin{claim*}
Let $\cat G$ be the set of objects of the form $i_{K!}\one$, where $i_K\colon */K\to */G$ and $K\subseteq G$ runs through the set of compact open subgroups with $\cd_{\Lambda}(K) < \infty$. Then $\cat G$ is a generating set of $\D(*/G,\Lambda)$ consisting of compact and $f$-prim objects.
\end{claim*}
To prove the claim, note that for $K$ as above we know by \cref{rslt:BG-is-proper} and  \cref{rslt:relation-between-suave-prim-and-dualizable} that the dualizable object $\one \in \D(*/K,\Lambda)$ is $f_K$-prim and by \cref{rslt:finite-coh-dim-implies-unit-compact} that it is compact, where $f_K\colon */K\to *$ is the structure map. As $i_K$ is $\Lambda$-suave by \cref{rmk:etale-proper-maps-classifying-stacks}, it follows from \cref{rslt:pushforward-of-suave-and-prim-objects} that $i_{K!}\one$ is $f$-prim. Since $i_{K!}$ preserves compact objects (its right adjoint $i_{K}^*$ admits a further right adjoint $i_{K*}$ and hence commutes with colimits), we conclude that $i_{K_!}\one$ is compact and $f$-prim. The fact that $\cat G$ is a generating set follows from \cref{rslt:colimit-of-invariants} and the observation 
\begin{align}
P^K = f_{K*}i_K^*P = f_{*}i_{K*}\iHom(\one, i_K^*P) = f_{*}\iHom(i_{K!}\one, P) = \RHom_{\D(*/G,\Lambda)}(i_{K!}\one, P);
\end{align}
here we denote $f_K\colon */K\to *$ and $f\colon */G\to *$, and the third equality holds by \cref{rslt:enriched-adjunction-of-shriek-functors}. This finishes the proof of the claim.

Let now $\cat G$ be as in the claim and denote $\langle \cat G\rangle \subseteq \D(*/G,\Lambda)$ the full subcategory generated by $\cat G$ under (co)fibers and retracts. Then \cref{rslt:suave-prim-stable-under-retracts} implies $\langle \cat G\rangle \subseteq \Prim(*/G)$. Since $\cat G$ consists of a set of compact generators of $\D(*/G,\Lambda)$, the functor $\Ind(\langle \cat G\rangle) \isoto \D(*/G,\Lambda)$ is an equivalence (see \cite[Proposition~5.3.5.11]{HTT}). By passing to compact objects and using \cite[Lemma~5.4.2.4]{HTT}, we obtain $\langle \cat G\rangle = \D(*/G,\Lambda)^{\omega}$ finishing the proof of (i).

Part (ii) follows from \cref{rslt:prim-generation-criterion-for-suave-objects} applied to the family $(Q_i)_{i\in I} = \cat G$ from the claim. One only needs to observe that the $\pi_1^*i_{K!}\one \tensor \pi_2^*i_{K'!}\one = i_{(K\times K')!}\one$ generate $\D(*/G\times G,\Lambda)$ by \cref{rslt:colimit-of-invariants}, where the $K, K'$ run through the compact open subgroups of $G$ with finite $\Lambda$-cohomological dimension.
\end{proof}

\begin{corollary}\label{rslt:rigidity-of-proper-BG}
Let $\Lambda$ be a commutative ring and $G$ a profinite group with $\cd_{\Lambda}(G) < \infty$. Then $\D(*/G,\Lambda)$ is \emph{rigid}. More precisely, $\D(*/G,\Lambda)$ is compactly generated and an object is compact if and only if it is dualizable.
\end{corollary}
\begin{proof}
By \cref{rslt:prim=compact} and the claim in its proof, the category $\D(*/G,\Lambda)$ is compactly generated, and the $f$-prim and the compact objects agree, where $f\colon */G\to *$ is the structure map. It remains to show that an object in $\D(*/G,\Lambda)$ is $f$-prim if and only if it is dualizable. Since $*/G\to *$ is $\Lambda$-proper by \cref{rslt:BG-is-proper}, this follows from \cref{rslt:etale-means-suave=dualizable}.
\end{proof}

\begin{remark}
It should be true that $\D(X, \Lambda)$ is rigid (in the sense of presentable stable symmetric monoidal categories) for every condensed anima $X$ such that $f\colon X \to *$ is $\Lambda$-fine and both $f$ and $\Delta_f$ are $\Lambda$-prim. In fact, this follows if one can show $\D(X \times X, \Lambda) = \D(X,\Lambda) \tensor_{\D(\Lambda)} \D(X,\Lambda)$ in $\PrL$, which by $!$-descent should reduce to the case that $X$ is a profinite set (where it is clear). We may pursue this line of ideas in future work.
\end{remark}

We now turn to the question of when $*/G \to *$ is $\Lambda$-suave. To this end, we make the following definition:

\begin{definition}
Let $\Lambda$ be a commutative ring. Let $G$ be profinite group with $d\coloneqq \cd_{\Lambda}(G) < \infty$. We say that $G$ is \emph{$\Lambda$-Poincaré (of dimension $d$)} if $f_* \colon \D(*/G,\Lambda) \to \Mod_{\Lambda}$ preserves dualizable objects, where $f\colon */G\to *$ denotes the structure map.

A locally profinite group is called \emph{locally $\Lambda$-Poincaré (of dimension $d$)} if it admits an open profinite subgroup which is $\Lambda$-Poincaré (of dimension $d$).
\end{definition}

The following example highlights the difference between being $\Lambda$-Poincaré and having finite $\Lambda$-cohomological dimension.

\begin{example}\label{example:non-Poincare}
Let $G = F(I)$ be the free pro-$p$ group on a set $I$. Then $G$ has $\FF_p$-cohomological dimension $\le 1$, and there is an isomorphism $\H^1(G,\FF_p) \cong \FF_p^{\lvert I\rvert}$, see \cite[\S4.2]{Serre.1997}. Therefore, if $I$ is infinite, then $G$ has finite $\FF_p$-cohomological dimension but is not $\FF_p$-Poincaré.
\end{example}

\begin{proposition}\label{rslt:suave-equiv-to-locally-Poincare}
Let $\Lambda$ be a commutative ring and $G$ a locally profinite group which has locally finite $\Lambda$-cohomological dimension. Then $*/G$ is $\Lambda$-suave if and only if $G$ is locally $\Lambda$-Poincaré.
\end{proposition}
\begin{proof}
For any compact open subgroup $K\subseteq G$ we denote $f_K\colon */K\to *$ and $i_K\colon */K\to */G$ the obvious maps.
Suppose that $G$ is locally $\Lambda$-Poincaré and fix a compact open subgroup $H\subseteq G$ which is $\Lambda$-Poincaré. Since $*/H\to */G$ is a $\Lambda$-suave cover, it suffices to show by \cref{rslt:suave-map-is-local-on-source} that $*/H\to *$ is $\Lambda$-suave. We may thus assume from the beginning that $G$ is $\Lambda$-Poincaré. Let $K\subseteq G$ be an open subgroup (necessarily with $\cd_{\Lambda}(K) < \infty$). Then $i_{K*}\one = i_{K!}\one$ is dualizable (combine \cref{rslt:finite-coh-dim-implies-unit-compact,rslt:rigidity-of-proper-BG}) and hence so is $f_{K*}\one = f_*i_{K*}\one$, as $G$ is $\Lambda$-Poincaré. It follows that $\one \in \D(*/G,\Lambda)$ is admissible, hence $*/G$ is $\Lambda$-suave by \cref{rslt:suave=admissible}.

Conversely, suppose that $*/G$ is $\Lambda$-suave. Replacing $G$ by a compact open subgroup if necessary, we may assume that $G$ is profinite with $\cd_{\Lambda}G <\infty$. As $\one \in \D(*/G,\Lambda)$ is admissible by \cref{rslt:suave=admissible}, we have that $f_*i_{K*}\one = f_{K*}\one$ is dualizable in $\Mod_{\Lambda}$, for every compact open subgroup $K\subseteq G$. But as the $i_{K*}\one$ generate the category of dualizable objects in $\D(*/G,\Lambda)$ under (co)fibers and retractions (cf.~\cref{rslt:rigidity-of-proper-BG}), we deduce that $f_*$ preserves dualizable objects. In other words, $G$ is $\Lambda$-Poincaré.
\end{proof}

\begin{lemma}\label{rslt:functoriality-of-Poincare-in-base-ring}
Let $\varphi\colon \Lambda \to \Lambda'$ be a morphism of commutative rings and $G$ a locally profinite group. If $G$ is locally $\Lambda$-Poincaré (resp.\ with invertible dualizing complex), then $G$ is locally $\Lambda'$-Poincaré (resp.\ with invertible dualizing complex).
\end{lemma}
\begin{proof}
Assume that $G$ is locally $\Lambda$-Poincaré. In particular, $G$ has locally finite $\Lambda'$-cohomological dimension. In view of \cref{rslt:suave-equiv-to-locally-Poincare}, it remains to show that $\Lambda$-suaveness implies $\Lambda'$-suaveness and that the invertibility of the dualizing complex over $\Lambda$ implies invertibility over $\Lambda'$. But both statements follow immediately from the fact that the natural transformation $\D(\blank, \Lambda) \to \D(\blank,\Lambda')$ (which is pointwise symmetric monoidal) of 3-functor formalisms induced by $\varphi$ (see \cref{rslt:functoriality-of-6ff-on-Cond-Ani-in-Lambda}) gives rise to a 2-functor $\cat K_{\D(\blank,\Lambda)} \to \cat K_{\D(\blank,\Lambda')}$ by \cref{rslt:functoriality-of-kerncat-in-D}. 
\end{proof}

There are essentially two situations in which we can prove that a locally profinite group is locally $\Lambda$-Poincaré and determine its dualizing complex. This is discussed in the following two examples.

\begin{example}\label{example:dualizing-complex-away-from-p}
Let $\Lambda$ be a $\ZZ[1/p]$-algebra and $G$ a locally pro-$p$ group. We claim that $G$ is locally $\Lambda$-Poincaré and compute its dualizing complex.

By \cref{rslt:functoriality-of-Poincare-in-base-ring} we may assume $\Lambda = \ZZ[1/p]$. By passing to an open subgroup if necessary, we may further assume that $G$ is pro-$p$. In this case, the category $\Rep_{\ZZ[1/p]}(G)^{\heartsuit}$ is semisimple so that $V^G$ is a retract of $V$, for every $V\in \Rep_{\ZZ[1/p]}(G)^{\heartsuit}$. It follows that $G$ is $\ZZ[1/p]$-Poincaré with $\cd_{\ZZ[1/p]}G = 0$. Moreover, the computation \eqref{eq:invertibility-of-dualizing-complex} shows that $\omega_f = \one$ is trivial.

For a general locally pro-$p$ group $G$ we can compute the dualizing complex $\omega_f$ explicitly as follows: Let $\Delta\colon */G\to */G\times */G$ be the diagonal map and make the identification $\Delta_!\one = \cInd_{\Delta(G)}^{G\times G}\one \cong \Cont_c(G)$ via restriction along $G \cong \{1\}\times G \injto G\times G$; note that the $G\times G$-action on $\Cont_c(G)$ is then given by $((a,b)\hat{f})(g) = \hat{f}(a^{-1}gb)$ for all $a,b,g\in G$. Now, the formula for $\omega_f$ in \cref{rslt:dual-criterion-for-suave-map} implies
\begin{align}
\Hom_{\Lambda}(\omega_f,\one) = \Hom_G(\Cont_c(G), \Lambda),
\end{align}
where $G$ acts on $\Cont_c(G)$ via right translation. The right hand side is by definition the space of $\Lambda$-valued right Haar measures on $G$, and $G$ acts through its left translation action on $\Cont_c(G)$. Thus, $G$ acts on $\omega_f$ via the modulus character $\delta\colon G\to \Lambda^{\times}$; concretely, it is given by the generalized index $\delta(g) = [gKg^{-1}:K]$, where $K\subseteq G$ is \emph{any} fixed open pro-$p$ subgroup, cf.\ \cite[I, 2.6]{Vigneras.1996}.
\end{example}

\begin{example}\label{example:dualizing-complex-mod-p}
Let $\Lambda$ be a $\ZZ/p^n$-algebra (for some $n\ge1$) and $G$ a $p$-adic Lie group. We claim that $G$ is locally $\Lambda$-Poincaré and has an invertible dualizing complex, which we will compute explicitly.

As in the previous example, we may assume $\Lambda = \ZZ/p^n$ and that $G$ is torsionfree and pro-$p$. Now, \cite[Corollaire (1)]{Serre.1965} shows $\cd_{\ZZ/p^n}G < \infty$ and hence $G$ is $\FF_p$-Poincaré by \cite[V, (2.5.8)]{Lazard.1965}. From the finite resolution in \cite[Proposition~4.1.1 (4)]{Symonds-Weigel.2000} we deduce by base-change to $\ZZ/p^n$ that $G$ is also $\ZZ/p^n$-Poincaré. Using \cref{rslt:colimit-of-invariants} one computes
\begin{align}\label{eq:invertibility-of-dualizing-complex}
q^*f^!\one &= \varinjlim_{K} \RHom_{\ZZ/p^n}(\one^K, \one),
\end{align}
where $q\colon *\to */G$ is the obvious map and $K\subseteq G$ runs through the open subgroups. By \cite[Steps (5) and (4) in Proposition~30]{Serre.1997} (and since each $\one^K$ is dualizable) we deduce $\omega_f = f^!\one = \one[\cd_{\ZZ/p^n}G]$. (See also \cite[Proposition~3.1.10]{Heyer.2023} for a similar computation.)

If $G$ is a general $p$-adic Lie group of dimension $d$, then $G$ acts on the dualizing complex $\omega_f$ via the character $\chi\colon G\to \Lambda^{\times}$, which is defined as follows: Let $\lie{g}$ be the $\QQ_p$-Lie algebra of $G$ in the sense of \cite[Definition on p.~100]{Schneider:Lie-Groups}. The determinant of the adjoint action of $G$ on $\lie g$ yields a character $\det(\lie g)\colon G\to \QQ_p^{\times}$. If $\lvert \cdot\rvert_p\colon \QQ_p^{\times} \to p^{\ZZ}$ denotes the $p$-adic absolute value, normalized by $\lvert p\rvert_p = p^{-1}$, then $\delta = \lvert \det(\lie g)\rvert_p \colon G\to p^{\ZZ}$ is the usual modulus character of $G$. We obtain a character
\begin{align}
\mathfrak{d}_G\coloneqq \det(\lie{g})\cdot \lvert\det(\lie g)\rvert_p \colon G\to \ZZ_p^{\times}.
\end{align}
Then $\chi$ is given by the composition of $\mathfrak{d}_G$ with the natural projection $\ZZ_p^{\times} \to \Lambda^{\times}$. For $\Lambda = \FF_p$, this is in \cite[Proposition~2.3.4]{Heyer.2023} formally deduced from \cite[Theorem~5.1]{Kohlhaase.2017}. But note that the proof of \loccit{} actually explains that this also works for $\Lambda= \ZZ/p^n$, $n>1$.
\end{example}

\subsection{Identifying the exceptional pushforward}

We fix a commutative ring $\Lambda$ and a locally profinite group $G$ which has locally finite $\Lambda$-cohomological dimension. 

If $H\subseteq G$ is a closed subgroup, the structure map $f\colon */H\to */G$ is $\Lambda$-fine by \cref{rslt:maps-of-class-stacks-are-fine}. One goal of this section is to identify the functor
\begin{align}
f_!\colon \D(*/H,\Lambda) \to \D(*/G,\Lambda)
\end{align}
on the level of smooth representations with the compact induction functor under the equivalence in \cref{rslt:identification-of-sheaves-and-representations}.

\begin{definition}
Let $H\subseteq G$ be a closed subgroup of a locally profinite group. For $V\in \Rep_{\Lambda}(H)^{\heartsuit}$ we define 
\begin{align}
\cInd_H^G(V) &\coloneqq \set{f\colon G\to V}{\begin{array}{l} \text{$f$ is locally constant,} \\
\text{$f(hg) = hf(g)$ for all $h\in H$, $g\in G$, and} \\
\text{the image of $\Supp(f)$ in $H\backslash G$ is compact} \end{array}}.
\end{align}
Letting $G$ act via right translation, $\cInd_H^G(V)$ becomes a smooth $G$-representation, which is usually referred to as \emph{compact induction}. The functor $\cInd_H^G$ is exact and hence we denote its induced derived functor again by
\begin{align}
\cInd_H^G \colon \hatRep_{\Lambda}(H) \to \hatRep_{\Lambda}(G).
\end{align}
\end{definition}

The following lemma reduces us to comparing the functors $f_!$ and $\cInd_H^G$ on the hearts.

\begin{lemma}\label{rslt:lower-shriek-is-t-exact}
Let $\Lambda$ be a commutative ring and $H\subseteq G$ a closed subgroup in a locally profinite group with locally finite $\Lambda$-cohomological dimension. Then the functor $f_!\colon \D(*/H,\Lambda) \to \D(*/G,\Lambda)$ is t-exact.
\end{lemma}
\begin{proof}
Note that we have the following cartesian square:
\begin{equation}\label{eq:lower-shriek-is-t-exact}
\begin{tikzcd}
H\backslash G \ar[d,"q'"'] \ar[r,"f'"] & * \ar[d,"q"]\\
*/H \ar[r,"f"'] & */G,
\end{tikzcd}
\end{equation}
where we remark that $H\backslash G$ is a condensed \emph{set}, since the action of $H$ on $G$ by left multiplication is free. By base-change, we thus have $q^*f_! \simeq f'_!q'^*$. Write $H\backslash G = \bigdunion_{i\in I} X_i$ as a disjoint union of profinite sets $X_i$, and denote by $f'_i\colon X_i\to *$ and $q_i\colon X_i \injto H\backslash G$ the obvious maps. For every $V\in \D(*/H,\Lambda)$ we then have
\begin{align}
q^*f_!V = f'_!q'^*V = \bigoplus_{i\in I}f'_{i*}q_i^*q'^*V.
\end{align}
Since all functors on the right hand side are t-exact, we deduce that $f_!$ is t-exact as well.
\end{proof}

\begin{proposition}\label{rslt:compact-induction-lower-shriek}
Let $\Lambda$ be a commutative ring and $H\subseteq G$ a closed subgroup in a locally profinite group with locally finite $\Lambda$-cohomological dimension. Then the diagram
\begin{equation}
\begin{tikzcd}
\Rep_{\Lambda}(H) \ar[d,"\simeq"'] \ar[r,"\cInd_H^G"] & \Rep_{\Lambda}(G) \ar[d,"\simeq"] \\
\D(*/H,\Lambda) \ar[r,"f_!"'] & \D(*/G,\Lambda)
\end{tikzcd}
\end{equation}
commutes, where the vertical maps are the isomorphisms from \cref{rslt:identification-of-sheaves-and-representations}.
\end{proposition}
\begin{proof}
By \cref{rslt:lower-shriek-is-t-exact} it suffices to prove the assertion on the hearts. We start by identifying the functor $f_!\colon \cat A(*/H,\Lambda) \to \cat A(*/G,\Lambda)$. By \cite[Lemma~2.3]{Abe-Henniart-Vigneras.2019}, the projection $\pi\colon G\to H\backslash G$ admits a continuous section, say $s$. Since $\pi$ is a morphism over $*/H$, so is $s$. Likewise, the inclusion of the unit $*\to G$ induces a morphism $* \to H\backslash G$ over $*/H$. The composite $H\backslash G \to G\to *$ identifies with $f'$ in \eqref{eq:lower-shriek-is-t-exact} (because $*$ is final) and hence the diagram
\begin{equation}
\begin{tikzcd}
* \ar[dr,"p"'] \ar[rr,shift left,"e"] & & \ar[ll,shift left, "f'"] H\backslash G \ar[dl,"q'"] \\
& */H
\end{tikzcd}
\end{equation}
commutes in both ways (where the isomorphism $pf' \simeq q'$ depends on the section $s$). Passing to the Čech nerves of $p$ and $q'$, we obtain morphisms
\begin{align}
e^{\bullet} \colon H^\bullet \rightleftarrows G^{\bullet}\times H\backslash G \noloc f'^{\bullet},
\end{align}
of simplicial diagrams $\bbDelta^{\op} \to \Cond(\Ani)$, which induce the identity on $*/H$ after passing to the colimit. The maps $e^n$ are given by $(h_1,\dotsc,h_n)\mapsto (h_1,\dotsc,h_n,H)$, whereas the maps $f'^n$ can be explicitly described as
\begin{align}
f'^1(g,\gamma) &=\; \text{unique element $h\in H$ with $hs(\gamma)g^{-1} = s(\gamma g^{-1})$}, \\
f'^2(g_1,g_2,\gamma) &= \begin{array}[t]{l} \text{unique element $(h_1,h_2) \in H^2$ such that} \\
\text{$h_2s(\gamma)g_2^{-1} = s(\gamma g_2^{-1})$ and $h_1s(\gamma g_2^{-1})g_1^{-1} = s(\gamma g_2^{-1}g_1^{-1})$,}
\end{array} \\
&\;\;\vdots \\
f'^n(g_1,\dotsc,g_n,\gamma) &= \begin{array}[t]{l} \text{unique element $(h_1,\dotsc,h_n) \in H^n$ such that} \\
\text{$h_is(\gamma g_n^{-1}\dotsm g_{i+1}^{-1})g_i^{-1} = s(\gamma g_n^{-1}\dotsm g_i^{-1})$ for all $i$,}
\end{array}
\end{align}
where $\gamma \in H\backslash G$ and $g, g_i \in G$. The functor $f_!\colon \cat A(*/H,\Lambda) \to \cat A(*/G,\Lambda)$ can now be computed as
\begin{equation}
\begin{tikzcd}
\cat A(*/H,\Lambda) \ar[d,equals] \ar[r] & \Mod^{\heartsuit}_{\Lambda} \ar[d,"f'^*"] \ar[r,shift left] \ar[r,shift right] & \Mod^{\heartsuit}_{\Lambda_c(H)} \ar[r,shift left=2] \ar[r] \ar[r,shift right=2] \ar[d,"f'^{1*}"] & \Mod^{\heartsuit}_{\Lambda_c(H^2)} \ar[d,"f'^{2*}"] \\
\cat A(*/H,\Lambda) \ar[r] \ar[d,"f_!"'] & \Mod^{\heartsuit}_{\Lambda_c(H\backslash G)} \ar[r,shift left] \ar[r,shift right] \ar[d,"f'_!"] & \Mod^{\heartsuit}_{\Lambda_c(H\backslash G\times G)} \ar[r,shift left=2] \ar[r] \ar[r,shift right=2] \ar[d] & \Mod^{\heartsuit}_{\Lambda_c(H\backslash G\times G^2)} \ar[d] \\
\cat A(*/G, \Lambda) \ar[r] & \Mod^{\heartsuit}_{\Lambda} \ar[r,shift left] \ar[r,shift right] & \Mod^{\heartsuit}_{\Lambda_c(G)} \ar[r,shift left=2] \ar[r] \ar[r,shift right=2] & \Mod^{\heartsuit}_{\Lambda_c(G^2)},
\end{tikzcd}
\end{equation}
where the rows are limit diagrams and the vertical maps in the second row (except $f_!$) are the forgetful maps. It follows that for every $(V, \alpha) \in \cat A(*/H)$, with $V\in \Mod^{\heartsuit}_{\Lambda}$ and $\alpha\colon d^{0*}V \isoto d^{1*}V$ an isomorphism in $\Mod^{\heartsuit}_{\Lambda_c(H)}$ satisfying the cocycle condition (cf. \cref{rslt:Grothendieck-descent}), we obtain an isomorphism
\begin{align}
f_!(V,\alpha) \isoto \bigl(\Lambda_c(H\backslash G) \tensor_{\Lambda} V, \beta\bigr) = \bigl(\Cont_c(H\backslash G,V), \beta\bigr)
\end{align}
in $\Mod^{\heartsuit}_{\Lambda_c(G)}$, where $\Cont_c(H\backslash G,V)$ is the $\Lambda$-module of locally constant maps $H\backslash G\to V$ with compact support and $\beta$ is the image of 
\begin{align}
f'^{1*}\alpha \colon \Lambda_c(G\times H\backslash G) \dotimes_{\Lambda_c(H\backslash G)} \Cont_c(H\backslash G,V) \isoto \Lambda_c(G\times H\backslash G) \dotimes_{\Lambda_c(H\backslash G)} \Cont_c(H\backslash G,V),
\end{align}
under the forgetful functor $\Mod^{\heartsuit}_{\Lambda_c(G\times H\backslash G)} \to \Mod^{\heartsuit}_{\Lambda_c(G)}$; beware that on the left hand side the $\Lambda_c(H\backslash G)$-module structure on $\Lambda_c(H\backslash G\times G)$ is induced by the projection $G\times H\backslash G\to H\backslash G$ whereas on the right it is induced by $(g,\gamma) \mapsto \gamma g^{-1}$. If $(V,\rho_{\alpha})$ denotes the smooth $H$-representation associated with $(V,\alpha)$ under the equivalence in \cref{rslt:identification-of-sheaves-and-representations}, then $(\Cont_c(H\backslash G, V), \beta)$ corresponds to the smooth $G$-representation $(\Cont_c(H\backslash G,V),\rho_{\beta})$, where
\begin{align}
(\rho_{\beta}(g)\varphi)(\gamma) &= \rho_{\alpha}\bigl(f'^1(g,\gamma g)\bigr) \varphi(\gamma g)
\end{align}
for all $g\in G$, $\varphi\in \Cont_c(H\backslash G,V)$ and $\gamma\in H\backslash G$.

Finally, it remains to show that the map
\begin{align}
\Phi\colon \cInd_H^G(V,\rho_{\alpha}) &\isoto \bigl(\Cont_c(H\backslash G,V), \rho_{\beta}\bigr),\\
\varphi &\mapsto \varphi\comp s
\end{align}
is a $G$-equivariant isomorphism. It is clearly a $\Lambda$-linear isomorphism with inverse
\begin{align}
\Cont_c(H\backslash G, V) &\to \cInd_H^G(V), \\
\psi &\mapsto [g\mapsto h_g\psi(Hg)], 
\end{align}
where $h_g\in H$ is the unique element satisfying $h_gs(Hg) = g$. In order to show that $\Phi$ is $G$-equivariant we compute, for $g\in G$, $\varphi\in \cInd_H^G(V)$ and $\gamma\in H\backslash G$,
\begin{align}
\bigl[\rho_{\beta}(g)\Phi(\varphi)\bigr](\gamma) &= \rho_{\alpha}\bigl(f'^1(g,\gamma g)\bigr) \Phi(\varphi)(\gamma g) = \rho_{\alpha}\bigl(f'^1(g,\gamma g)\bigr) \varphi(s(\gamma g)) \\
&= \varphi\bigl(f'^1(g,\gamma g)s(\gamma g)\bigr) = \varphi\bigl(s(\gamma) g\bigr) = (g\cdot \varphi)(s(\gamma)) = \Phi\bigl(g\cdot \varphi\bigr)(\gamma).
\end{align}
\end{proof}

\begin{remark}
The analog of \cref{rslt:compact-induction-lower-shriek} in the étale 6-functor formalism on v-stacks had been observed already in \cite[Example~4.2.3]{Hansen-Kaletha-Weinstein.2022}, where a different argument is sketched; the details have been worked out in Shram's Master's thesis \cite[Proposition~7.6]{Shram}.
\end{remark}

\begin{lemma}\label{rslt:derived-coinvariants}
Let $\Lambda$ be a commutative ring. Let $f\colon G\epito H$ be a quotient map between locally profinite groups such that $\Ker(f)$ is locally $\Lambda$-Poincaré. Then the inflation functor $\Inf^H_G$ admits a left adjoint, called the functor of \emph{derived coinvariants}, 
\[
\L^H_G\colon \hatRep_{\Lambda}(G) \to \hatRep_{\Lambda}(H).
\]
If, moreover, the dualizing complex $\omega_f \in \hatRep_{\Lambda}(G)$ is dualizable, then $f_!\colon \D(*/G,\Lambda) \to \D(*/H,\Lambda)$ identifies with $\L^H_G\iHom_H(\omega_f,\blank)$ under the equivalence in \cref{rslt:identification-of-sheaves-and-representations}.
\end{lemma}
\begin{proof}
The induced map $f\colon */G\to */H$ is $\Lambda$-suave, since the pullback along $*\to */H$ is the $\Lambda$-suave map $*/\Ker(f) \to *$ (cf. \cref{rslt:suave-equiv-to-locally-Poincare} and \cref{rslt:suave-and-prim-map-is-local-on-target}). By \cref{rslt:suave-maps-induce-twist-of-shriek-functors} the functor $f^*\colon \D(*/H,\Lambda) \to \D(*/G,\Lambda)$ admits a left adjoint $f_{\natural} \simeq f_!(\omega_f\tensor\blank)$, which under the equivalence of \cref{rslt:identification-of-sheaves-and-representations} corresponds to the left adjoint of the inflation functor $\Inf^H_G$. If $\omega_f$ is dualizable, then $\iHom_H(\omega_f,\blank)\colon \D(*/H,\Lambda) \to \D(*/H,\Lambda)$ is  left (and right) adjoint to $\omega_f\tensor\blank$, which implies the final claim.
\end{proof}

\begin{remark} \label{rmk:Heyer-coinvariants-functor}
Let $\Lambda$ be a field of characteristic $p>0$ and $f\colon G\epito H$ a quotient map between $p$-adic Lie groups. The derived coinvariants functor $\L^H_G \colon \Rep_{\Lambda}(G) \to \Rep_{\Lambda}(H)$ has been studied in \cite{Heyer.2023} using representation theoretic methods to prove the existence. The existence had been independently observed in \cite{Hansen-Mann:Mod-p-Stacky-6FF}, whose proof strategy we follow in this paper.
\end{remark}

\begin{remark}
Let $\Lambda$ be a commutative ring and $f\colon G\to H$ a morphism of locally $\Lambda$-Poincaré groups. Since $f$ factors as a composition $G\injto G\times H \epito H$, where the first map is the graph of $f$ and the second map is the projection, one can combine \cref{rslt:compact-induction-lower-shriek,rslt:derived-coinvariants} to give a representation theoretic description of $f_!\colon \D(*/G,\Lambda)\to \D(*/H,\Lambda)$ (assuming that the dualizing complex of $*/G\to *$ is dualizable).
\end{remark}

\subsection{The canonical anti-involution on derived Hecke algebras}
In this section we will apply our 6-functor formalism to the construction of a certain involutive anti-automorphism on derived Hecke algebras. In the literature such an anti-involution has been considered in various degrees of generality, and we explain its historical development before giving our own construction in \cref{rslt:anti-involution-Hecke}.

\begin{remark}
Let $\Lambda$ be a field of characteristic $p>0$. We fix a locally profinite group $G$.
\begin{remarksenum}
\item Let $K\subseteq G$ be a compact open subgroup and $V\in \Rep_{\Lambda}(K)^{\heartsuit}$ a finite-dimensional representation. The algebra $\Hecke(G,K,V) \coloneqq \End_G(\cInd_K^GV)$ is called the \emph{Hecke algebra} associated with $(G,K,V)$. The map
\begin{align}
\Hecke(G,K,V) &\isoto \set{f\colon G\to \End_{\Lambda}(V)}{\begin{array}{l} \text{$f(kgk') = kf(g)k'$ for all $k,k'\in K$ and $g\in G$} \\
\text{and $f$ has compact support}
\end{array}}, \\
T &\mapsto \bigl[g\mapsto [v\mapsto T([1,v])(g)]\bigr],
\end{align}
is an isomorphism, where for all $g\in G$ and $v\in V$ we define $[g,v]$ as the unique map in $\cInd_K^GV$ with support $Kg^{-1}$ and taking the value $v$ at $g^{-1}$. It was observed in \cite[\S2.3]{Herzig.2011}, but see also \cite[above Proposition~2.4]{Henniart-Vigneras.2012}, that under this identification the map 
\begin{align}
\iota\colon \Hecke(G,K,V) \isoto \Hecke(G,K,V^*), \qquad \iota(T)(g) \coloneqq (T(g^{-1}))^*
\end{align}
is an involutive anti-isomorphism of algebras; here $V^* = \Hom_{\Lambda}(V,\Lambda)$ denotes the $\Lambda$-linear dual and $(T(g^{-1}))^*$ the dual map of $T(g^{-1})$.

\item Assume that $G$ is a $p$-adic Lie group. Let $I\subseteq G$ be the pro-$p$ Iwahori subgroup and assume that $I$ is $p$-torsionfree. In \cite{Ollivier-Schneider.2019}, Ollivier--Schneider introduce the pro-$p$ Iwahori--Hecke $\Ext$-algebra
\begin{align}
E^* \coloneqq \Ext^*_G\bigl(\cInd_I^G\Lambda, \cInd_I^G\Lambda\bigr),
\end{align}
which becomes a graded algebra with the Yoneda product.\footnote{Strictly speaking, \cite{Ollivier-Schneider.2019} consider the opposite algebra of $E^*$ but for the sake of exposition we ignore this difference.} The subalgebra $E^0$ identifies with the Hecke algebra $\Hecke(G,I,\Lambda)$ from (i). In \S6 of \opcit{} Ollivier--Schneider construct an anti-involution on $E^*$ as follows: Put $I_g \coloneqq I\cap gIg^{-1}$ for $g\in G$ and fix a set $W\subseteq G$ of representatives for $I\backslash G/I$. Conjugation by $w$ induces on cohomology a $\Lambda$-linear isomorphism
\begin{align}
\conj_w\colon \H^i(I_{w^{-1}}, \Lambda) \isoto \H^i(I_w,\Lambda)
\end{align}
for all $i\ge0$ (\cite[\S6]{Ollivier-Schneider.2019} also explains how to deal with the ambiguity coming from the choice of $W$). Denoting by $\cInd_I^{IwI}\Lambda \subseteq \Res^G_I\cInd_I^{G}\Lambda$ the subrepresentation of those functions with support in $IwI$, the map $f\mapsto [x\mapsto f(wx)]$ constitutes an $I$-equivariant isomorphism $\cInd_I^{IwI}\Lambda \isoto \cInd_{I_{w^{-1}}}^I\Lambda$. By Shapiro we obtain a canonical isomorphism $\H^i(I,\cInd_I^{IwI}\Lambda) \isom \H^i(I_{w^{-1}}, \Lambda)$. Now define a $\Lambda$-linear isomorphism $\Inv^i_w$ by the commutativity of the diagram
\begin{equation}
\begin{tikzcd}
\H^i(I, \cInd_I^{IwI}\Lambda) \ar[d,"\sim"'] \ar[r,"\Inv^i_w", "\sim"'] & \H^i(I, \cInd_I^{Iw^{-1}I}\Lambda) \ar[d,"\sim"] \\
\H^i(I_{w^{-1}}, \Lambda) \ar[r,"\conj_w"', "\sim"] & \H^i(I_w, \Lambda).
\end{tikzcd}
\end{equation}
Under the isomorphisms $E^* \isom \H^*(I, \cInd_I^G\Lambda) = \bigoplus_{w\in W} \H^*(I, \cInd_I^{IwI}\Lambda)$ the $\Inv_w^*$ assemble into a graded linear isomorphism
\begin{align}
\Inv^*_{\mathrm{OS}}\colon E^* \isoto E^*,
\end{align}
which by \cite[Proposition~6.1]{Ollivier-Schneider.2019} is a graded involutive anti-isomorphism of algebras and restricts to the anti-involution $\iota \colon E^0\isoto E^0$ from (i).

\item Assume that $G$ is a $p$-adic Lie group. Let $K\subseteq G$ be a compact open subgroup. In \cite{Schneider-Sorensen.2023b} the authors make the crucial observation that the anti-involution on the Hecke $\Ext$-algebra $E^* = \Ext^*_G(\cInd_K^G\Lambda, \cInd_K^G\Lambda)$ should be induced by the map
\begin{align}\label{eq:swapping-factors}
\varsigma \colon \cInd_K^G\Lambda \tensor \cInd_K^G\Lambda \isoto \cInd_K^G\Lambda \tensor \cInd_K^G\Lambda
\end{align}
given by swapping the factors. 
To see this, Schneider--Sorensen explicitly define in \cite[\S2.2.2]{Schneider-Sorensen.2023b} a natural automorphism $J$ on $\Hom_{K}(V, \Ind_K^GW)$ via
\begin{align}
J(\alpha)(v)(g) \coloneqq g\cdot \bigl(\alpha(gv)(g^{-1})\bigr),
\end{align}
where $V, W \in \Rep_{\Lambda}(G)^{\heartsuit}$, $\alpha\in \Hom_K(V, \Ind_K^GW)$, $v\in V$ and $g\in G$.\footnote{Beware that Schneider--Sorensen follow a different convention for the smooth induction $\Ind_K^GW$, and hence the explicit description of $J$ in \loccit{} looks slightly different.} 
Note that we omit the restriction functor $\Rep_{\Lambda}(G)^{\heartsuit}\to \Rep_{\Lambda}(K)^{\heartsuit}$ from the notation to increase legibility. From the definition it is immediate that $J$ is involutive and maps $\Hom_K(V, \cInd_K^{KgK}W)$ bijectively onto $\Hom_K(V, \cInd_K^{Kg^{-1}K}W)$ for all $g\in G$. Specializing to $V = W = \Lambda$ and taking an injective resolution $\Lambda \to \mathcal{I}^{\bullet}$ in $\Rep_{\Lambda}(G)^{\heartsuit}$ yields an involutive automorphism $J$ on
\begin{align}
\Hom^{\bullet}_{G}(\cInd_K^G\Lambda, \cInd_K^G\mathcal{I}^{\bullet}) \isom \Hom^{\bullet}_{K}(\Lambda, \cInd_K^G\mathcal{I}^{\bullet}) \isom \bigoplus_{g\in K\backslash G/K} \Hom^{\bullet}_{K}(\Lambda, \cInd_K^{KgK}\mathcal{I}^{\bullet}).
\end{align}
Passing to cohomology and using that $\cInd_K^G\Lambda \to \cInd_K^G\mathcal{I}^{\bullet}$ is an $\H^0(K,\blank)$-acyclic resolution yields a graded anti-involution
\begin{align}\label{eq:InvSS}
\Inv^*_{\mathrm{SS}} \colon E^* \isoto E^*;
\end{align}
this is proved in \cite[Proposition~2.6]{Schneider-Sorensen.2023b} by showing $\Inv^*_{\mathrm{SS}} = \Inv^*_{\mathrm{OS}}$. The relation between $J$ and $\varsigma$ is explained by the existence of the commutative diagram
\begin{equation}
\begin{tikzcd}
\Hom_{G}(\cInd_K^G\Lambda, \Ind_K^GV) \ar[d,"\sim"'] \ar[r,"J"] & \Hom_{G}(\cInd_K^G\Lambda, \Ind_K^GV) \ar[d,"\sim"] \\
\Hom_{G}(\cInd_K^G\cInd_K^G\Lambda, V) \ar[d,"\sim"'] & \Hom_{G}(\cInd_K^G\cInd_K^G\Lambda, V) \ar[d,"\sim"] \\
\Hom_{G}(\cInd_K^G\Lambda \tensor \cInd_K^G\Lambda, V) \ar[r,"\varsigma^*"'] & \Hom_{G}(\cInd_K^G\Lambda \tensor \cInd_K^G\Lambda, V)
\end{tikzcd}
\end{equation}
for all $V \in \Rep_{\Lambda}(G)^{\heartsuit}$: see \cite[Lemma~2.7]{Schneider-Sorensen.2023b} and the preceding discussion. Here, the top vertical maps are given by applying Frobenius reciprocity twice and the lower vertical maps are induced by (the inverse of) the $G$-equivariant isomorphism
\begin{align}\label{eq:compact-induction-projection-formulas}
\cInd_K^G\cInd_K^G W & \isoto \cInd_K^G\Lambda \tensor \cInd_K^G\Lambda \tensor W, \\
\bigl[g, [g',w]\bigr] &\mapsto [g] \tensor [gg'] \tensor gg'w,
\end{align}
for $W \in \Rep_{\Lambda}(G)^{\heartsuit}$ (applied with $W= \Lambda$). Here, the symbols $[g',w]$ were defined in (i) and we abbreviate $[g] \coloneqq [g,1]$.
\end{remarksenum}
\end{remark}

We now present our construction of the anti-involution, which works for a general commutative ring $\Lambda$ and for the derived Hecke algebra (i.e.\ not just its associated Ext-algebra). We start by recollecting the necessary ingredients. In the following, fix a commutative ring $\Lambda$ and a locally profinite group $G$ such that $G$ has locally finite $\Lambda$-cohomological dimension.
\begin{itemize}
\item Our 6-functor formalism $\D(\blank,\Lambda)$ factors through $\LMod_{\Mod_{\Lambda}}(\PrL)^{\tensor}$ by \cref{rslt:6ff-with-presentable-cat}; in particular, the categories $\D(X,\Lambda)$ and all six functors are automatically enriched over $\Mod_{\Lambda}$.

\item There is a canonical 2-isomorphism $\theta\colon \cat K_{\D(\blank,\Lambda)}^\op \isoto \cat K_{\D(\blank,\Lambda)}$ (\cref{rslt:kern-cat-is-self-dual}).

\item For every condensed anima $X$ there is a canonical self-inverse equivalence $\DPrim\colon \Prim(X)^{\op} \isoto \Prim(X)$ (\cref{rslt:SD-PD-are-self-inverse-equivalences}), which is enriched over $\Mod_{\Lambda}$ in view of the explicit formula in \cref{rslt:criterion-for-prim-object}.

\item There is a canonical isomorphism $\Rep_{\Lambda}(G) \isom \D(*/G,\Lambda)$ (\cref{rslt:identification-of-sheaves-and-representations}).

\item The full subcategory $\Prim(*/G)$ of $f$-prim objects corresponds to the full subcategory $\Rep_{\Lambda}(G)^{\omega}$ of compact objects (\cref{rslt:prim=compact}), where $f\colon */G\to *$ is the structure map.

\item Suppose that $G$ is compact and $\cd_\Lambda G < \infty$. Then the full subcategory $\Prim(*/G)$ of $f$-prim objects corresponds to the full subcategory $\Rep_{\Lambda}(G)^{\dbl}$ of dualizable objects in $\Rep_{\Lambda}(G)$ (\cref{rslt:rigidity-of-proper-BG}), where $f\colon */G\to *$ is the structure map.

\item Let $K \subseteq G$ be an open subgroup. Recall that the map $i_K\colon */K\to */G$ is $\Lambda$-étale (\cref{rmk:etale-proper-maps-classifying-stacks}). The functor $\cInd_K^G$, which is computed by $i_{K!}$ (see \cref{rslt:compact-induction-lower-shriek}), restricts to a ($\Mod_{\Lambda}$-enriched) functor
\begin{align}
\cInd_K^G\colon \Rep_{\Lambda}(K)^{\dbl} \to \Rep_{\Lambda}(G)^{\omega}
\end{align}
by \cref{rslt:pushforward-of-suave-and-prim-objects}.
\end{itemize}

\begin{notation}
Let $\Lambda$ be a commutative ring, $G$ a locally profinite group and $K\subseteq G$ a compact open subgroup with $\cd_{\Lambda}(K) < \infty$. Let $\Hecke_K$ be the $\Mod_{\Lambda}$-enriched category whose objects are the dualizable objects in $\Rep_{\Lambda}(K)$ and with morphisms
\begin{align}
\Hecke_K(V,W) \coloneqq \RHom_G(\cInd_K^GV, \cInd_K^GW) \in \Mod_{\Lambda}.
\end{align}
More formally, $\Hecke_K$ is the full $\Mod_{\Lambda}$-enriched subcategory of $\Rep_{\Lambda}(G)^{\omega}$ spanned by the essential image of $\cInd_K^G\big|_{\Rep_{\Lambda}(K)^{\dbl}}$.

We call $\Hecke_K(V,V)$ the \emph{derived Hecke algebra of weight $V$}; it is in fact an algebra in $\Mod_{\Lambda}$ (combine \cref{ex:full-subcategory-of-enriched-category,ex:algebra-is-enriched-category}). In the special case $V = \one$, we further abbreviate
\begin{align}
\Hecke_K^{\bullet} \coloneqq \Hecke_K(\one,\one).
\end{align}
\end{notation}

\begin{proposition}\label{rslt:anti-involution-Hecke}
Let $\Lambda$ be a commutative ring. Let $G$ be a locally profinite group and $K\subseteq G$ a compact open subgroup with $\cd_{\Lambda}K <\infty$. Then prim duality $\DPrim$ on $\Prim(*/G)$ induces an involutive equivalence
\begin{align}
\Hecke_K^{\op} \isoto \Hecke_K
\end{align}
of $\Mod_{\Lambda}$-enriched categories, which is given on objects by $V\mapsto V^{\vee} = \RHom_{\Lambda}(V,\Lambda)$. 
\end{proposition}
\begin{proof}
Note that $\Hecke_K$ is the essential image of the $\Mod_{\Lambda}$-enriched functor $\cInd_K^G\colon \Rep_{\Lambda}(K)^{\dbl} \to \Rep_{\Lambda}(G)^{\omega}$. Denoting $f\colon */G\to *$ and $f_K\colon */K \to *$ the obvious maps, we need to show that the diagram
\begin{equation}
\begin{tikzcd}[column sep=4em]
\Rep_{\Lambda}(K)^{\dbl,\op} \ar[d,"\DPrim_{f_K}"'] \ar[r,"(\cInd_K^G)^{\op}"] & \Rep_{\Lambda}(G)^{\omega,\op} \ar[d,"\DPrim_{f}"] \\
\Rep_{\Lambda}(K)^{\dbl} \ar[r,"\cInd_K^G"'] & \Rep_{\Lambda}(G)^{\omega}
\end{tikzcd}
\end{equation}
is commutative and that $\DPrim_{f_K} = [V\mapsto V^{\vee}]$. 

We will prove the claim by translating it into the 2-category $\cat K\coloneqq \cat K_{\D(\blank,\Lambda)}$. To simplify the notation, we write $\cat K_{G} \coloneqq \cat K_{\D(\blank, \Lambda), */G}$. Note that $f$ induces a 2-functor $f_!\colon \cat K_{G} \to \cat K$ by \cref{rslt:functoriality-of-kerncat-pushforward}. Since $i_K$ is $\Lambda$-étale, we have an adjunction $\one \colon */K \rightleftarrows */G \noloc \one = i_K^!(\one)$ in $\cat K_{G}$ by \cref{rslt:etale-equiv-iso-of-shriek-and-star}. Since every 2-functor preserves adjunctions, by applying $f_!$ we thus obtain an adjunction $j'_!(\one) \colon */K \rightleftarrows */G \noloc j_!(\one)$ in $\cat K$, where $j = (\id, i_K) \colon */K\to */K\times */G$ and $j' = (i_K, \id) \colon */K\to */G\times*/K$. Now, for every $P \colon */K\to *$ in $\cat K$ one easily computes $i_{K!}(P) = P\comp j_!(\one)$; similarly, for every $Q \colon *\to */K$ in $\cat K$ one has $i_{K!}(Q) = j'_!(\one) \comp Q$. Since adjoints compose, we obtain an isomorphism 
\begin{align}
\DPrim_{f} \comp (\blank \comp j_!(\one)) \isom (j'_!(\one) \comp \blank)\comp \DPrim_{f_K}
\end{align}
of functors $\Prim(*/K)^{\op} \to \Prim(*/G)$. This translates to the commutativity of the displayed diagram.

It remains to prove $\DPrim_{f_K}(V) = V^{\vee}$ for every $V \in \Rep_{\Lambda}(V)^{\dbl}$. But note that $f_K$ is $\Lambda$-proper by \cref{rslt:BG-is-proper}, hence this follows from \cref{rslt:etale-means-suave=dualizable}. Here, we observe that $V^{\vee} = \iHom(V,\one)$ identifies with $\RHom_{\Lambda}(V,\Lambda)$ under the forgetful functor $\Rep_{\Lambda}(K) \to \Mod_{\Lambda}$ (because it is symmetric monoidal and hence preserves duals).
\end{proof}

In the rest of this subsection we compare the above construction with the construction by Schneider--Sorensen. We first need to discuss a technical lemma on double cosets. Recall that the inclusion $\Cond(\Set) \injto \Cond(\Ani)$ has a left adjoint $X\mapsto \pi_0(X)$, which informally assigns to a condensed anima its condensed set of isomorphism classes (where we view a condensed anima as a \enquote{topologically enriched $\infty$-groupoid}).

\begin{lemma}\label{rslt:decomposition-of-double-cosets}
Let $G$ be a locally profinite group and $H, K\subseteq G$ closed subgroups with $K$ open. For the action of $H\times K$ on $G$ via $(h,k)\cdot g = hgk^{-1}$ we consider the quotient stack $H\backslash G/K \coloneqq G/(H\times K)$ in $\Cond(\Ani)$. Then: 
\begin{enumerate}[(i)]
\item The canonical map
\begin{align}
\bigdunion_{g\in \pi_0(H\backslash G/K)} */H\cap gKg^{-1} \isoto */H\bigtimes_{*/G} */K \isom H\backslash G/K 
\end{align}
is an isomorphism.

\item Under the identification in (i), the projections are induced by the inclusion $H\cap gKg^{-1} \injto H$ and the $g^{-1}$-conjugation $H\cap gKg^{-1} \injto gKg^{-1} \xto{\conj_{g^{-1}}} K$.

\item Under the identification in (i), the swap isomorphism $s\colon */H\times_{*/G} */K \isoto */K \times_{*/G} */H$ is induced by conjugation:
\begin{align}
\bigdunion_{g\in \pi_0(H\backslash G/K)} */H\cap gKg^{-1} \xto{\dunion_g \conj_{g^{-1}}} \bigdunion_{g\in \pi_0(H\backslash G/K)} */K\cap g^{-1}Hg,
\end{align}
where on the right hand side we identified $\pi_0(K\backslash G/H)$ with $\pi_0(H\backslash G/K)$ via $g\mapsto g^{-1}$.
\end{enumerate}
\end{lemma}
\begin{proof}
In order to prove (i), we first argue that $H\backslash G/K \isom */H \times_{*/G} */K \eqqcolon X$. Consider the commutative diagram
\begin{equation}
\begin{tikzcd}
G \ar[d] \ar[r] & H\backslash G \ar[d] \ar[r] & * \ar[d] \\
G/K \ar[d] \ar[r] & X \ar[d] \ar[r] & */K \ar[d] \\
* \ar[r] & */H \ar[r] & */G,
\end{tikzcd}
\end{equation}
where every square is cartesian.
Observe that the actions of $H$ and $K$ on $G$ commute with each other, and we obtain a functor $\bbDelta^{\op} \times \bbDelta^{\op} \to \Cond_{\kappa}(\Ani)$ (which in view of \cref{rslt:embed-Top-into-Cond-Ani} comes from a functor of ordinary categories) given on objects by $([m], [n]) \mapsto H^m\times K^n\times G$ such that precomposition with the diagonal $\bbDelta^{\op} \to \bbDelta^{\op}\times \bbDelta^{\op}$ yields the Čech nerve $H^\bullet\times K^\bullet \times G \to G/(H\times K)$ of $G\to G/(H\times K)$. As $\bbDelta^{\op}$ is sifted (see \cite[\href{https://kerodon.net/tag/02QP}{Tag 02QP}]{kerodon}), the claim follows from the computation
\begin{align}
H\backslash G/K &= \varinjlim_{[n]\in \bbDelta^{\op}} H^n\times K^n\times G 
= \varinjlim_{[m]} \varinjlim_{[n]} H^m\times K^n\times G \\
&= \varinjlim_{[m]} H^m \times \Bigl(\varinjlim_{[n]} K^n\times G\Bigr) 
= \varinjlim_{[m]} H^m \times G/K \\
&= \varinjlim_{[m]} H^m \times \bigl(*\times_{*/H}X\bigr) 
= \varinjlim_{[m]} H^m\times_{*/H}X 
= */H \times_{*/H} X \\
&= X
\end{align}
for the obvious action of $H$ on $G/K$; here, we have used that $\Cond_{\kappa}(\Ani)$ is a topos, so that in particular $\blank \times_{*/H}X$ commutes with colimits.

We next show that $\pi_0(X)$ is discrete, where now $\pi_0$ denotes the left adjoint to the inclusion $\Cond_{\kappa}(\Set) \injto \Cond_{\kappa}(\Ani)$. Since $K\subseteq G$ is open, the condensed anima $G/K$ is discrete. As $* \epito */H$ is an effective epimorphism, so is its pullback $G/K \epito X$. We obtain a (necessarily effective) epimorphism $f\colon \pi_0(G/K) \epito \pi_0(X)$ in $\Cond_{\kappa}(\Set)$, where the source is discrete. Note that the full subcategory of discrete condensed sets is closed under small colimits and subobjects. Thus, if $Y_{\bullet} \to \pi_0(X)$ denotes the Čech nerve of $f$, then each $Y_n \subseteq \pi_0(G/K)^{n+1}$ is discrete and hence $\pi_0(X) = \varinjlim_{[n]} Y_n$ is discrete.

For every element $g\in \pi_0(X)$ we fix a lift $g\colon *\to X$ (denoted with the same symbol). So far we have a decomposition $X = \bigdunion_{g\in \pi_0(X)} X_g$, where $X_g = \{g\} \times_{\pi_0(X)}X$ is a condensed anima with $\pi_0(X_g) = *$. We now show that $X_g = */H\cap gKg^{-1}$. Since limits commute with limits and $\{g\} \times_{\pi_0(X)} \{g\} = *$, we have
\begin{align}
\{g\} \times_{X_g} \{g\} = * \bigtimes_{\{g\}\bigtimes\limits_{\pi_0(X)}\{g\}} \bigl(\{g\} \times_X\{g\}\bigr) = \{g\} \times_{X} \{g\} \eqqcolon \Omega(X,g).
\end{align}
Now, $\{g\} \epito X_g$ is an effective epimorphism, and hence taking the colimit of its Čech nerve, the above computation shows $X_g = */\Omega(X,g)$. It remains to prove $\Omega(X,g) = H\cap gKg^{-1}$ as groups. Note that the automorphism on $G$ given by $\gamma \mapsto \gamma g^{-1}$ induces an isomorphism
\begin{align}
\rho_g\colon H^{\bullet} \times K^{\bullet} \times G \isoto H^{\bullet} \times (gKg^{-1})^{\bullet} \times G
\end{align}
of simplicial diagrams in $\Cond(\Ani)$. Denoting by $X'$ the geometric realization of the latter, we obtain a natural isomorphism $\Omega(X,g) \isom \{1\} \times_{X'} \{1\}$. Now consider the following commutative diagram
\begin{equation}
\begin{tikzcd}
H\cap gKg^{-1} \ar[dd] \ar[rr] &[3em] &  gKg^{-1} \ar[d,"gKg^{-1}\times\{1\}"] \ar[r] & \{1\} \ar[d] \\
& & gKg^{-1}\times G \ar[d,"\act"] \ar[r] & G \ar[d] \\
H \ar[r,"H\times \{1\}"] \ar[d] & H\times G \ar[d] \ar[r,"\act"] & G \ar[d] \ar[r] & H\backslash G \ar[d] \\
\{1\} \ar[r] & G \ar[r] & G/gKg^{-1} \ar[r] & X',
\end{tikzcd}
\end{equation}
where every square is cartesian. We deduce canonical isomorphisms $H\cap gKg^{-1} \isom \{1\} \times_{X'} \{1\} \isom \Omega(X,g)$ as desired. This finishes the proof that
\begin{align}
\bigdunion_{g\in \pi_0(H\backslash G/K)} */H\cap gKg^{-1} \isoto H\backslash G/K.
\end{align}
This finishes the proof of (i).

We now identify the projections from $X = */H \times_{*/G} */K$. It is clear from the construction that the restriction of the projections to $*/H\cap K$ is induced by the inclusions. As the projections are compatible with the isomorphisms $X\isoto X'$ induced by $\rho_g$, the general claim follows by contemplating the commutative diagram
\begin{equation}
\begin{tikzcd}
X \ar[d] \ar[r] & G/K \ar[d,"\cdot g^{-1}"] \ar[r] & */K \ar[d,"\conj_{g}"] \\
X' \ar[r] & G/gKg^{-1} \ar[r] & */gKg^{-1},
\end{tikzcd}
\end{equation}
in which all squares are cartesian. The case for the projection $X\to */H$ is similar. This proves (ii).

For part (iii), it is immediate from (i) and (ii) that the displayed map coincides with $s$.
\end{proof}

\begin{proposition}\label{rslt:anti-involution-comparison-with-Schneider-Sorensen}
Let $\Lambda$ be a field of characteristic $p>0$, let $G$ be a $p$-adic Lie group and $I\subseteq G$ a $p$-torsionfree compact open subgroup. Then the anti-involution
\begin{align}
\Inv \colon (\Hecke_I^\bullet)^{\op} \isoto \Hecke_I^\bullet
\end{align}
from \cref{rslt:anti-involution-Hecke} coincides with $\Inv_{\mathrm{SS}}^*$ \eqref{eq:InvSS} after passage to the cohomology algebras.
\end{proposition}
\begin{proof}
We first give a different description of $\Inv$ and then show how it relates to $\Inv_{\mathrm{SS}}^*$. Consider the correspondence $[*/G \xfrom{i_I} */I \xto{f_I} *]$ and denote by $A = i_{I!}\one \colon */G \to *$ its image in $\cat K \coloneqq \cat K_{\D(\blank,\Lambda)}$ (cf.\ \cref{rslt:functors-from-and-to-kercat}), and write $B = \DPrim(A) = i_{J!}\one \colon * \to */G$, where we put $J\coloneqq I$ in order to keep track of the different roles $I$ plays. We contemplate the following commutative diagram of isomorphisms, where we further write $A' = A$ and $B' = B$ to increase legibility:
\begin{equation}
\begin{tikzcd}
\RHom_G(A,A') \ar[r,"\rho"] & \RHom_{\Lambda}(\id_*, A' \comp B) \ar[d,"\theta"'] \ar[r,"\lambda"] & \RHom_G(B', B) \ar[d,"\theta"] \\
& \RHom_{\Lambda}(\id_*, \theta(B) \comp \theta(A')) \ar[r,"\lambda"'] & \RHom_G(\theta(B'), \theta(B)).
\end{tikzcd}
\end{equation}
Here, $\rho$ and $\lambda$ denote the passage to right and left mates, respectively; they are reviewed in \cref{rslt:mate-correspondence}. The top-right circuit in the diagram defines the prim duality $(\Hecke_I^{\bullet})^{\op} \isoto \Hecke_I^{\bullet}$. We now analyze the bottom-left circuit. One easily computes
\begin{align}
\RHom_{\Lambda}(\id_*, A'\comp B) &= A'\comp B = f_!(i_{I!}\one \tensor i_{J!}\one) \\
&= f_!i_{I!}i_I^*i_{J!}\one = f_{I!}i_I^*i_{J!}\one ,
\end{align}
where $f\colon */G \to *$ denotes the structure map; the third equality uses the projection formula. We similarly have $\RHom_{\Lambda}(\id_*, \theta(B)\comp \theta(A')) \isom f_{J!}i_J^*i_{I!}\one$. Next we identify the isomorphism $f_{I!}i_I^*i_{J!}\one \isoto f_{J!}i_J^*i_{I!}\one$ induced by $\theta$. The isomorphism $\theta$ above factors into isomorphisms $A' \comp B \isoto \theta(A' \comp B) \isoto \theta(B) \comp \theta(A')$. These are induced by the following 2-isomorphisms in the correspondence category:
\begin{equation}\label{eq:Involution-swap}
\begin{tikzcd}
& & X \ar[d, equals,"\id"] \ar[dl, bend right, "i_J\pi_1"'] \ar[dr, bend left, "i_I\pi_2"] \\
* & */G \ar[l,"f"'] & X \ar[l, "i_I\pi_2"'] \ar[r, "i_J\pi_1"] \ar[d, Rightarrow, "s"] & */G \ar[r,"f"] & * \\
& & \theta(X) \ar[ul, bend left, "i_I\pi_1"] \ar[ur, bend right, "i_J\pi_2"'],
\end{tikzcd}
\end{equation}
where we denote $X \coloneqq */J \times_{*/G} */I$ and $\theta(X) = */I \times_{*/G} */J$, and where $s\colon X\to \theta(X)$ is the swap isomorphism from \cref{rslt:decomposition-of-double-cosets}.(iii); keep in mind that the correspondence $\alpha = [*/G \from X \to */G]$ is to be understood as the composite
\begin{equation}
\alpha \colon 
\begin{tikzcd}
& & X \ar[dl,"\pi_1"'] \ar[dr, "\pi_2"] \\
& */J \ar[dl,"i_J"'] \ar[dr,"i_J"] & & */I \ar[dl,"i_I"'] \ar[dr,"i_I"] \\
*/G & & */G & & */G
\end{tikzcd}
\end{equation}
in the correspondence category. Note that the image of $*/G \xfrom{i_I} */I \xto{i_I} */G$ in $\cat K_G \coloneqq \cat K_{\D(\blank,\Lambda), */G}$ is given by $i_{I!}\one$, and hence the image of $\alpha$ is $i_{I!}\one \tensor i_{J!}\one$. Under these identifications the isomorphism $s$ visibly induces the swap isomorphism
\begin{align}
\varsigma\colon i_{I!}\one \tensor i_{J!}\one \isom i_{J!}\one \tensor i_{I!}\one.
\end{align}
Taking the image under the 2-functor $f_!\colon \cat K_G \to \cat K$, and then pre- and postcomposing with $\one \colon *\to */G$ and $\one \colon */G \to *$ in $\cat K$, respectively, shows that $s$ induces the swap isomorphism $f_!\varsigma \colon f_!(i_{I!}\one \tensor i_{J!}\one) \isoto f_!(i_{J!}\one \tensor i_{I!}\one)$ in $\cat K$.
To summarize, we deduce that the diagram
\begin{equation}
\begin{tikzcd}
\RHom_{\Lambda}(\id_*, A' \comp B) \ar[r,"\sim"] \ar[d,"\theta"'] & f_!(i_{I!}\one \tensor i_{J!}\one) \ar[d,"f_!\varsigma"] \\
\RHom_{\Lambda}(\id_*, \theta(B) \comp \theta(A')) \ar[r,"\sim"'] & f_!(i_{J!}\one \tensor i_{I!}\one)
\end{tikzcd}
\end{equation}
is commutative. 

Let $\iota_{\gamma}\colon */J\cap \gamma I\gamma^{-1} \injto \bigsqcup_{\gamma\in \pi_0(J\backslash G/I)} */J\cap \gamma I\gamma^{-1} \cong X$ be the inclusion (see \cref{rslt:decomposition-of-double-cosets} for the isomorphism). By the description of $s$ in \cref{rslt:decomposition-of-double-cosets}.(iii), the isomorphisms $i_J\pi_1 \isom i_I\pi_1s$ and $i_I\pi_2 \isom i_J\pi_2s$ in \eqref{eq:Involution-swap} are equivalently described by the commutative diagrams
\begin{equation}
\begin{tikzcd}[column sep=small]
*/J\cap \gamma I\gamma^{-1} \ar[dr,"i_{J\cap \gamma I\gamma^{-1}}"{description, name=inc}] \ar[ddr, bend right, "f_{J\cap \gamma I\gamma^{-1}}"'] \ar[rr,"\conj_{\gamma^{-1}}"] 
& & 
*/\gamma^{-1} J\gamma \cap I \ar[dl, "i_{\gamma^{-1}J\gamma \cap I}" description] \ar[ddl, bend left, "f_{\gamma^{-1}J\gamma\cap I}"] \ar[to=inc, Rightarrow, shorten <=6mm, shorten >=2mm, "c_{\gamma}"']
\\[3em]
& */G \ar[d,"f"] 
\\
& *
\end{tikzcd}
\end{equation}
for all $\gamma \in \pi_0(J\backslash G/I)$, where the 2-isomorphism $c_{\gamma}$ is induced by conjugation with $\gamma$. It follows that under the identification
\begin{align}
f_!(i_{I!}\one \tensor i_{J!}\one) &= f_!i_{I!}i_I^*i_{J!}\one = f_!i_{I!}\pi_{2!}\pi_1^*\one \\
&= \bigoplus_{\gamma} f_!i_{I!}\pi_{2!}\iota_{\gamma!} \iota_{\gamma}^* \pi_1^*\one = \bigoplus_{\gamma} f_{\gamma^{-1}J\gamma \cap I !} \one,
\end{align}
the swap isomorphism $f_!\varsigma$ corresponds to the isomorphism
\begin{align}
\bigoplus_{\gamma\in \pi_0(J\backslash G/I)} f_{\gamma^{-1}J\gamma \cap I!}\one \xto{(c_{\gamma})_{\gamma}}
\bigoplus_{\gamma\in \pi_0(J\backslash G/I)} f_{J\cap \gamma I\gamma^{-1}!}\one.
\end{align}
To conclude, we have shown that the anti-involution $\Inv\colon (\Hecke_I^{\bullet})^{\op} \isoto \Hecke_I^{\bullet}$ is given as the composite
\begin{align}
\RHom_{G}(i_{I!}\one, i_{I!}\one) &\isom \bigoplus_{\gamma} f_{\gamma^{-1}J\gamma \cap I!} \one
\xto{(c_{\gamma})_{\gamma}} \bigoplus_{\gamma} f_{J\cap \gamma I\gamma^{-1}!}\one = \RHom_{G}(i_{J!}\one, i_{J!}\one).
\end{align}
Hence, after passing to cohomology, we obtain the involution $\Inv^*_{\mathrm{SS}}$ as asserted.
\end{proof}

\appendix
\section{\texorpdfstring{$\infty$}{∞}-Categories} \label{sec:cat}

The past two decades have seen a surge in the development of the theory of $\infty$-categories, i.e.\ $(\infty,1)$-categories. 
Intuitively an $\infty$-category is a structure consisting of objects, $1$-morphisms between objects and invertible $i$-morphisms between $(i-1)$-morphisms for all $i\ge2$; these $i$-morphisms can be composed in a meaningful way. 
It is a non-trivial task to keep track of all the data involved, and over the years there have been established several equivalent models of $\infty$-categories based on simplicial sets, e.g.\ simplicially enriched categories, complete Segal spaces and quasicategories. 
A vast portion of the theory was developed using the language of quasicategories starting with foundational work by Joyal \cite{Joyal.2002, Joyal.2008} and Lurie \cite{HTT, HA}. 

This paper is deeply rooted in the theory of $\infty$-categories, and since their application to representation theory is relatively new, it seems necessary to write a few words to orient the reader. 
Of course, this is not the place to give a thorough introduction into this complex topic for which there are already several excellent resources like \cite{kerodon, Haugseng.2017, Rezk.2017, Land.2021} to name a few. Instead, this appendix serves the purpose to set up notation, introduce the tools that are used in the paper and to point to the relevant literature.

\subsection{Quasicategories and Kan complexes}\label{sec:quasicategories}
Let $\bbDelta$ be the category of non-empty totally ordered finite sets together with order-preserving maps. Up to isomorphism the objects of $\bbDelta$ are given by the sets $[n]\coloneqq \{0,1,\dotsc,n\}$, for $n\ge0$, with their natural ordering. Recall that a \emph{simplicial set} is a functor $X\colon \bbDelta^\op \to \Set$, $[n] \mapsto X([n])\eqqcolon X_n$, and that simplicial sets together with natural transformations form a category $\sSet$. For each $n\ge0$ we let $\Delta^n\in \sSet$ be the functor represented by $[n]$, and for $0\le i\le n$ we denote $\Horn^n_i \subseteq \Delta^n$ the \emph{$i$-th horn} obtained by deleting the interior and the face opposite the $i$-th vertex in $\Delta^n$.

\begin{definition} 
A simplicial set $X \in \sSet$ is called:
\begin{defenum}
\item \label{def:quasicategory} a \emph{quasicategory} if $X$ admits $i$-th horn fillers for all $0<i<n$, i.e.\ every map $\Horn^n_i\to X$ extends along the inclusion $\Horn^n_i \injto \Delta^n$ to an $n$-simplex $\Delta^n\to X$.

\item a \emph{Kan complex} if $X$ admits $i$-th horn fillers for all $0\le i\le n$.
\end{defenum}
We denote $\Kan \subseteq \qCat \subseteq \sSet$ the full subcategories of Kan complexes and quasicategories, respectively.
\end{definition} 

If $\cat C$ is a quasicategory, we call $\cat C_0$ the set of \emph{objects} and $\cat C_1$ the set of \emph{morphisms} of $\cat C$. A map $\Horn^2_1 \to \cat C$ is then the datum of a sequence $X\xto f Y \xto g Z$ of morphisms in $\cat C$, and the horn-filling condition precisely says that this can be extended to a $2$-simplex $\sigma\colon \Delta^2\to \cat C$ depicted as
\begin{equation}
\begin{tikzcd}[column sep=small]
	& Y \ar[dr,"g"]\\
	X \ar[ur,"f"] \ar[rr,"h"'] & & Z.
\end{tikzcd}
\end{equation}
We say that $\sigma$ exhibits $h$ as a composition of $f$ and $g$. By abuse of notation one usually writes $h = g\comp f$; beware that a composite is rarely unique, but that the \enquote{space of composites is contractible} \cite[\href{https://kerodon.net/tag/007A}{Tag 007A}]{kerodon}.
\medskip

In the literature Kan complexes are also called $\infty$-groupoids, spaces, or anima. The inclusion $\Kan \subseteq \qCat$ admits a right adjoint sending a quasicategory $\cat C$ to its maximal subsimplicial set $\cat C^{\simeq}$ which is a Kan complex.
We call $\cat C^{\simeq}$ the \emph{underlying anima} of $\cat{C}$, cf.~\cite[\href{https://kerodon.net/tag/01CZ}{Tag 01CZ}]{kerodon}.

Given any two simplicial sets $X,Y$, we can form the simplicial set $\Fun(X,Y)$ whose set of $n$-simplices is given by $\Hom_{\sSet}(\Delta^n\times X,Y)$. 
This exhibits $\sSet$ as a cartesian closed category, and in particular $\sSet$ is simplicially enriched. 
We remark that, if $Y$ is a quasicategory (resp.\ a Kan complex), then so is $\Fun(X,Y)$ by \cite[\href{https://kerodon.net/tag/0066}{Tag 0066}]{kerodon} (resp.\ \cite[\href{https://kerodon.net/tag/00TN}{Tag 00TN}]{kerodon}).

The following construction provides many examples of quasicategories by associating them with simpler objects that appear throughout mathematics.

\begin{principle}\label{principle:adjunctions}
Let $\cat C$ be a cocomplete ordinary category and $C^\bullet\colon \bbDelta \to \cat C$ a cosimplicial object. By the density of the Yoneda embedding $\bbDelta \to \sSet$ we obtain an adjunction
\begin{equation}
\begin{tikzcd}
	\lvert\,\cdot\,\rvert_{C^\bullet} \colon \sSet \ar[r,shift left] &
	\ar[l,shift left] \cat C \noloc \Sing_{C^\bullet},
\end{tikzcd}
\end{equation}
where $\lvert\,\cdot\,\rvert_{C^\bullet}$ is the unique colimit preserving functor satisfying $\lvert \Delta^n\rvert_{C^\bullet} = C^n$ for all $n\ge0$, and for each $X\in \cat C$, the simplicial set $\Sing_{C^\bullet}(X)$ is given by
\begin{align}
\Sing_{C^\bullet}(X)_n \coloneqq \Hom_{\cat C}(C^n, X).
\end{align}
\end{principle}

\begin{example}
Let $\cat C = \Cat_{\ordinary}$ be the category of ordinary categories and let $C^\bullet\colon \bbDelta \injto \Cat_{\ordinary}$ be given by the inclusion. For any category $\cat A\in \Cat_{\ordinary}$, the simplicial set $\Nerve(\cat A) \coloneqq \Sing_{C^\bullet}(\cat A)$ is called the \emph{nerve} of $\cat A$. Then $\Nerve(\cat A)$ is always a quasicategory \cite[\href{https://kerodon.net/tag/0032}{Tag 0032}]{kerodon}, and $\Nerve(\cat A)$ is a Kan complex if and only if $\cat A$ is a groupoid. For example, $\Delta^n = \Nerve([n])$ is a quasicategory but a Kan complex only if $n=0$.

Since the nerve functor $\Nerve$ is fully faithful \cite[\href{https://kerodon.net/tag/002Z}{Tag 002Z}]{kerodon}, we will usually identify $\cat A$ with its associated quasicategory $\Nerve(\cat A)$. The special property of $\Nerve(\cat A)$ compared to general quasicategories is that (higher) composition is unique, i.e.\ the horn fillers in \cref{def:quasicategory} are unique.

For any simplicial set $X \in \sSet$, the category $\homotopy X \coloneqq \lvert X\rvert_{C^\bullet}$ is called the \emph{homotopy category} associated with $X$. If $X$ is a quasicategory then $\homotopy X$ is the ordinary category whose objects are those of $X$ but whose morphisms are given by isomorphism classes of morphisms in $X$. Note that the counit $\homotopy\Nerve(\cat A) \isoto \cat A$ is an equivalence of categories.
\end{example}

\begin{example}
Let $\cat C = \Top$ be the category of topological spaces. Let $\lvert\Delta^n\rvert$ be the topological $n$-simplex consisting of all $(t_0,\dotsc,t_n)\in [0,1]^{n+1}$ such that $\sum_{i=0}^nt_i = 1$ and consider the cosimplicial object $C^\bullet \coloneqq \lvert \Delta^\bullet\rvert$ of $\Top$. For any topological space $X$, the \emph{singular simplicial set} $\Sing(X) \coloneqq \Sing_{\lvert \Delta^\bullet\rvert}(X)$ is a Kan complex \cite[\href{https://kerodon.net/tag/002K}{Tag 002K}]{kerodon}.
\end{example}

\begin{example}
Let $\cat C = \Cat_{\Delta}$ be the category of simplicially enriched ordinary categories, and let $C^\bullet = \Path[\bullet]$ be the cosimplicial object of $\Cat_{\Delta}$ constructed in \cite[\href{https://kerodon.net/tag/00KN}{Tag 00KN}]{kerodon}. For any simplicially enriched category $\cat A$, the simplicial set $\Nerve^{\mathrm{hc}}(\cat A) \coloneqq \Sing_{\Path[\bullet]}(\cat A)$ is called the \emph{homotopy coherent nerve} of $\cat A$. Concretely, the objects of $\Nerve^{\mathrm{hc}}(\cat A)$ are the objects of $\cat A$ and a morphism $f\colon X\to Y$ in $\Nerve^{\mathrm{hc}}(\cat A)$ is a $0$-simplex of the simplicial set $\Hom_{\cat A}(X,Y)$. A $2$-simplex consists of morphisms $f\colon X\to Y$, $g\colon Y\to Z$ and $h\colon X\to Z$, together with a $1$-simplex $\sigma\colon g\comp f \To h$ in $\Hom_{\cat A}(X,Z)_1$.

By a Theorem of Cordier--Porter \cite[\href{https://kerodon.net/tag/00LJ}{Tag 00LJ}]{kerodon}, $\Nerve^{\mathrm{hc}}(\cat A)$ is a quasicategory provided $\cat A$ is \emph{locally Kan}, i.e.\ for all $X,Y \in \cat A$ the simplicial sets $\Hom_{\cat{A}}(X,Y)$ are Kan complexes. The following examples are fundamentally important:
\begin{exampleenum}
\item\label{ex:category-of-anima} The full subcategory $\Kan \subseteq \sSet$ of Kan complexes is a locally Kan simplicially enriched category. We denote
\begin{align}
\Ani \coloneqq \Nerve^{\mathrm{hc}}(\Kan)
\end{align}
the \emph{$\infty$-category of (small) anima}. In the literature, $\Ani$ is also denoted $\cat{S}$ and referred to as the $\infty$-category of spaces. Its objects are those quasicategories where all morphisms are invertible (i.e.\ \enquote{$\infty$-groupoids}).

There is a fully faithful embedding $\Set \injto \Ani$. It admits a left adjoint, which assigns to an anima $X$ its set of \emph{connected components} $\pi_0 X$.

\item\label{ex:category-of-categories} The full subcategory $\qCat \subseteq \sSet$ of quasicategories is enriched in Kan complexes via the underlying anima $\Fun(\cat A, \cat B)^{\simeq}$ of the functor category, where $\cat A,\cat B \in \qCat$. We denote
\begin{align}
\Cat \coloneqq \Nerve^{\mathrm{hc}}(\qCat)
\end{align}
the \emph{quasicategory of (small) $\infty$-categories}.

There is a unique (up to contractible choice) non-trivial automorphism on $\Cat$, which assigns to an $\infty$-category $\cat C$ its \emph{opposite category} $\cat C^{\op}$ \cite[Theorem~4.4.1]{Lurie.2009}.
\end{exampleenum}
\end{example}

\begin{remark}\label{remark:model-independence}
The classical way to deal with higher categorical structures is to use the theory of model categories, which provides a notion of \emph{fibrations} and \emph{cofibrations} on an underlying ordinary category. For example, one can model the $\infty$-category $\Cat$ of $\infty$-categories using the Joyal model structure on $\sSet$ (see \cite[\S2.2.5]{HTT}) and the $\infty$-category $\Ani$ using the Kan--Quillen model structure on $\sSet$ (see \cite[\S A.2.7]{HTT}). Both are uniquely determined by their cofibrations and fibrant objects, see \cite[Proposition~E.1.10]{Joyal.2008}.


Our approach in this paper is to work in the categories $\Ani$ and $\Cat$ directly rather than in $\sSet$ endowed with the Kan--Quillen and Joyal model structures, respectively. This model-independent approach allows us to make conceptual constructions, which hopefully improves the clarity of the exposition.
\end{remark}

\begin{example}
Let $\cat C = \Ch(\ZZ)$ be the category of chain complexes of abelian groups. For $n\ge0$, let $C^n = C(\Delta^n)$ be the chain complex which in degree $m$ is the free abelian group on the non-degenerate $m$-simplices of $\Delta^n$, and where the differential $d\colon C^n_m\to C^n_{m-1}$ is given by $d(\sigma) = \sum_{i=0}^m (-1)^i d_i(\sigma)$. (Here, $d_i$ denotes the $i$-th face map of the $m$-simplex $\sigma$.) Then $C^\bullet$ is a cosimplicial object in $\Ch(\ZZ)$. The Dold--Kan correspondence \cite[\href{https://kerodon.net/tag/00QQ}{Tag 00QQ}]{kerodon} says that the adjunction
\begin{equation}
\begin{tikzcd}
	\lvert\,\cdot\,\rvert_{C^\bullet}\colon \Fun(\bbDelta^\op,\Ab) \ar[r,shift left,"\sim"] & \ar[l,shift left] \Ch(\ZZ)_{\ge0} \noloc \Sing_{C^\bullet}.
\end{tikzcd}
\end{equation}
is an equivalence of categories. Endowing $\Fun(\bbDelta^\op,\Ab)$ with the levelwise tensor product and $\Ch(\ZZ)_{\ge0}$ with the tensor product of chain complexes, there is a colax (and also a lax) symmetric monoidal structure on $\lvert\,\cdot\,\rvert_{C^\bullet}$ \cite[\href{https://kerodon.net/tag/00S0}{Tag 00S0}]{kerodon}, and hence $\Sing_{C^\bullet}$ is lax symmetric monoidal.

Let now $\cat A$ be a differential graded category, i.e.\ a category enriched in $\Ch(\ZZ)$. Transferring the enrichment along $\Sing_{C^\bullet}$, we obtain a locally Kan, simplicially enriched category $\cat A^\Delta$. The quasicategory $\Nerve^{\mathrm{hc}}(\cat A^{\Delta})$ may be identified with the \emph{differential graded nerve} $\Nerve^{\mathrm{dg}}(\cat A)$ in view of \cite[\href{https://kerodon.net/tag/00SV}{Tags 00SV}, \href{https://kerodon.net/tag/00PK}{00PK}]{kerodon}. We refer to \loccit{} for an explicit description of the $n$-simplices.
\end{example}

Let $\cat C$ be a quasicategory. For any two objects $X,Y \in \cat C$, we define the \emph{mapping anima} $\Hom_{\cat C}(X, Y)$ as the fiber product
\begin{equation}\begin{tikzcd}
	\Hom_{\cat C}(X,Y) \ar[d] \ar[r] \ar[dr,phantom,very near start, "\lrcorner"] & \Fun(\Delta^1,\cat C) \ar[d,"{(\ev_0,\ev_1)}"]\\
	\Delta^0 \ar[r,"{(X,Y)}"'] & \cat C\times \cat C.
\end{tikzcd}\end{equation}
Note that $\Hom_{\cat C}(X,Y)$ is indeed a Kan complex \cite[\href{https://kerodon.net/tag/01JC}{Tag 01JC}]{kerodon}. By \cite[\href{https://kerodon.net/tag/01PF}{Tag 01PF}]{kerodon} there is a (unique up to contractible choice) composition law $\comp\colon \Hom_{\cat C}(Y,Z) \times \Hom_{\cat C}(X,Y) \to \Hom_{\cat C}(X,Z)$, which enhances $\homotopy\cat C$ to a $\homotopy\Kan$-enriched category.
\medskip

Many constructions for ordinary categories admit a straightforward generalization to the $\infty$-categorical setting.

\begin{definition}\label{def:adjoint}
Let $G\colon \cat D \to \cat C$ be a functor between quasicategories. We say that \emph{$G$ admits a left adjoint at $X \in \cat C$} if there exists $Y\in \cat D$ and a morphism $u\colon X\to G(Y)$ in $\cat C$ such that for all $Z\in \cat D$ the composite
\begin{align}
\Hom_{\cat D}(Y,Z) \xto{G} \Hom_{\cat C}\bigl(G(Y), G(Z)\bigr) \xto{\comp u} \Hom_{\cat C}\bigl(X, G(Z)\bigr)
\end{align}
becomes an isomorphism in $\homotopy\Kan$. In this case, $Y$ is called the \emph{left adjoint} of $G$ at $X$.
By \cite[\href{https://kerodon.net/tag/02FV}{Tag 02FV}]{kerodon}, $G$ admits a left adjoint at every $X$ if and only if there exists a functor $F\colon \cat C \to \cat D$ and natural transformations $\eta \colon \id_{\cat C} \to G F$ and $\varepsilon \colon F G\to \id_{\cat D}$ such that the composites 
\begin{align}
F \xto{F\eta} F G F \xto{\varepsilon F} F\quad\text{and}\quad
G \xto{\eta G} G F G \xto{G \varepsilon} G
\end{align}
are (homotopic to) $\id_F$ and $\id_G$, respectively. Then $F$ is called the \emph{left adjoint of $G$}, and the natural transformations $\eta$ and $\varepsilon$ are called the \emph{unit} and \emph{counit}, respectively. The dual notion of a \emph{right adjoint} is defined in the obvious way. 
\end{definition} 

We remark that adjunctions really are a $2$-categorical concept \cite[\href{https://kerodon.net/tag/02CA}{Tag 02CA}]{kerodon}. We refer to \cite[\href{https://kerodon.net/tag/02EH}{Tag 02EH}]{kerodon} for an in-depth study of adjunctions between quasicategories, see also \cref{sec:2cat.adj} for the general concept of adjoint morphisms in an $(\infty,2)$-category.

\begin{example} \label{def:lim-and-colim-in-category}
Let $\cat I$, $\cat C$ be quasicategories and consider the functor
\begin{align}
\delta \colon \cat C = \Fun(\Delta^0,\cat C) \to \Fun(\cat I, \cat C).
\end{align}
The \emph{colimit} of a functor $F\colon \cat I \to \cat C$ is a left adjoint object $c \in \cat C$ of $\delta$ at $F$. Concretely, this means that there is a natural transformation $\eta\colon F \to \delta(c)$ of functors $\cat I \to \cat C$ such that for every $d \in \cat C$ the natural map
\begin{align}
\Hom_{\cat C}(c,d) \isoto \Hom_{\Fun(\cat I, \cat C)}(F, \delta(d))
\end{align}
is an isomorphism in $\homotopy\Kan$. The colimit $c$ is usually denoted by $\colim F$, $\colim_{i\in \cat I}F(i)$ or $\varinjlim_{i\in \cat I} F(i)$. If $\cat I = \emptyset$ is empty, then $\emptyset \coloneqq \colim F \in \cat C$ is called an \emph{initial object} of $\cat C$.

Similarly, a \emph{limit} of $F$ is a right adjoint object of $\delta$ at $F$, which we usually denote by $\lim F$, $\lim_{i\in \cat I} F(i)$ or $\varprojlim_{i\in \cat I} F(i)$. If $\cat I = \emptyset$ is empty, then $*\coloneqq \lim F \in \cat C$ is called a \emph{terminal object}.
\end{example}

The above definition of limits and colimits generalizes to the concept of Kan extensions, which is one of the most important constructions in category theory and is used abundantly in this paper. Let us briefly recall the main definitions and properties from \cite[\href{https://kerodon.net/tag/02Y1}{Tag 02Y1}]{kerodon}.

\begin{definition}\label{def:Kan-extension}
Let $\varphi\colon \cat I \to \cat J$ and $F\colon \cat I \to \cat C$ be functors of quasicategories.
\begin{defenum}
	\item A \emph{(pointwise) left Kan extension} of $F$ along $\varphi$ is a functor $\ol F\colon \cat J \to \cat C$ together with a natural transformation $\eta\colon F \to \ol F \comp \varphi$ such that for all $j\in \cat J$ the natural map
	\begin{align}
		\varinjlim_{i \in \cat I_{/j}} F(i) \isoto \ol F(j)
	\end{align}
	induced by $\eta$ is an isomorphism, where $\cat I_{/j} \coloneqq \cat I \times_{\cat J} \cat J_{/j}$ with $\cat J_{/j}$ being the over-category from \cref{def:over-category} below. 

	\item A \emph{(pointwise) right Kan extension} of $F$ along $\varphi$ is a functor $\ol F\colon \cat J \to \cat C$ together with a natural transformation $\varepsilon\colon \ol F \comp \varphi \to F$such that for all $j \in \cat J$ the natural map
	\begin{align}
		\ol{F}(j) \isoto \varprojlim_{i \in \cat I_{j/}} F(i)
	\end{align}
	induced by $\varepsilon$ is an isomorphism, where $\cat I_{j/} \coloneqq \cat I \times_{\cat J} \cat J_{j/}$ with $\cat J_{j/}$ being the under-category from \cref{def:over-category}.
\end{defenum}
\end{definition}

As an example, the left Kan extension of a functor $F\colon \cat I \to \cat C$ along the projection $\cat I \to *$ is the same as a colimit of $F$, while the right Kan extension is the same as a limit. In general, we can depict the Kan extension in the following diagram:
\begin{equation}\begin{tikzcd}
	\cat I \ar[r,"F"] \ar[d,"\varphi"'] & \cat C \\
	\cat J \ar[ur,"\ol{F}"',dashed]
\end{tikzcd}\end{equation}
One should see $\ol F$ as a universal way of extending $F$ to a functor $\cat J \to \cat C$. The power of Kan extensions comes from the fact that its existence can be checked pointwise and that its computation is functorial in $F$. This allows many constructions of functors which would otherwise be hard to perform (as in $\infty$-land, functors can almost never be constructed \enquote{by hand}). For example, given a sheaf theory $\D\colon \cat C^\op \to \Cat$ on some site $\cat C$, one can extend $\D$ from $\cat C$ to a functor $\D\colon \Shv(\cat C)^\op \to \Cat$ on all sheaves/stacks on $\cat C$ by performing a right Kan extension along the functor $\cat C^\op \to \Shv(\cat C)^\op$. In the following result we collect some properties of Kan extensions.

\begin{proposition}\label{rslt:properties-of-Kan-extensions}
Let 
\begin{equation}
\begin{tikzcd}
\cat I \ar[r,"F"{name=F}] \ar[d,"\varphi"'] & \cat C \\
\cat J \ar[ur,"\ol{F}"'] \ar[from=F,Rightarrow,shorten=3mm, "\eta"']
\end{tikzcd}
\end{equation}
be a (non-commutative) diagram of functors and $\eta\colon F \to \ol{F}\comp \varphi$ a natural transformation of functors $\cat I \to \cat C$.
\begin{enumerate}[(i)]
\item If for every $j\in \cat J$ the colimit $\varinjlim_{i \in \cat I_{/j}} F(i)$ exists in $\cat C$, then the left Kan extension of $F$ along $\varphi$ exists.

\item \label{rslt:Kan-extension-is-adjoint-to-restriction} Suppose that $\eta\colon F \to \ol{F}\comp \varphi$ exhibits $\ol{F}$ as the left Kan extension of $F$ along $\varphi$. Then for all functors $G\colon \cat J \to \cat C$ the natural map
\begin{align}
\Hom_{\Fun(\cat J, \cat C)}(\ol{F},G) \isoto \Hom_{\Fun(\cat I, \cat C)}(F, G\comp\varphi)
\end{align}
is an isomorphism in $\homotopy\Kan$. The converse holds if $\cat C$ admits $\cat I_{/j}$-indexed colimits for all $j\in \cat J$.

\item Suppose that $\varphi\colon \cat I \injto \cat J$ is fully faithful. If $\eta$ exhibits $\ol{F}$ as the left Kan extension of $F$ along $\varphi$, then $\eta$ is an isomorphism.

\item Suppose that $\varphi$ admits a right adjoint $\psi\colon \cat J \to \cat I$. Then the natural map $\eta\colon F \to F\psi\comp \varphi$ (which is induced by the unit $\id \to \psi\varphi$) exhibits $F\psi$ as the left Kan extension of $F$ along $\varphi$. In particular, the left Kan extension exists.

\item (Transitivity) Let $\alpha\colon \cat J \to \cat K$ and $\ol{\ol{F}}\colon \cat K \to \cat C$ be functors and $\ol{\eta} \colon \ol{\ol{F}} \to \ol{F}\alpha$ a natural transformation of functors $\cat J \to \cat C$. If $\eta \colon F \to \ol{F}\varphi$ exhibits $\ol{F}$ as a left Kan extension of $F$ along $\varphi$ and $\ol{\eta}\colon \ol{F} \to \ol{\ol{F}}\alpha$ exhibits $\ol{\ol{F}}$ as a left Kan extension of $\ol{F}$ along $\alpha$, then $\ol{\eta}\varphi\comp \eta \colon F \to \ol{\ol{F}}\alpha\varphi$ exhibits $\ol{\ol{F}}$ as a left Kan extension of $F$ along $\alpha\varphi$.
\end{enumerate}
The analogous results for right Kan extensions and limits also hold (where in (iv) we require $\varphi$ to admit a left adjoint).
\end{proposition}
\begin{proof}
Part (i) is \cite[\href{https://kerodon.net/tag/0300}{Tag 0300}]{kerodon}, (ii) follows from \cite[\href{https://kerodon.net/tag/030D}{Tag 030D}]{kerodon} and (i), (iii) is \cite[\href{https://kerodon.net/tag/02YN}{Tag 02YN}]{kerodon} and (v) follows from \cite[\href{https://kerodon.net/tag/031M}{Tag 031M}]{kerodon}. It remains to prove (iv). The adjunction $\varphi \colon \cat I \rightleftarrows \cat J \noloc \psi$ induces an adjunction $\psi^* \colon \Fun(\cat J,\cat C) \rightleftarrows \Fun(\cat I, \cat C) \noloc \varphi^*$. By the converse direction in (ii) it remains to prove that $\cat C$ admits $\cat I\times_{\cat J} \cat J_{/j}$-indexed limits for all $j\in \cat J$. But this is obvious from the fact that, since $\varphi$ is a left adjoint, $\cat I \times_{\cat J} \cat J_{/j}$ has a terminal object by \cite[\href{https://kerodon.net/tag/02J9}{Tag 02J9}]{kerodon}.
\end{proof}

Note that by \cref{rslt:Kan-extension-is-adjoint-to-restriction}, under the assumptions of enough colimits in $\cat C$ the left Kan extension $\varphi_!\colon F \mapsto \ol F$ defines a left adjoint to the restriction $\varphi^*\colon G \mapsto G \comp \varphi$. This in particular shows that Kan extension is functorial (see \cref{def:adjoint}).

\begin{definition}
A quasicategory $\cat C$ is called \emph{complete} (resp.\ \emph{cocomplete}) if every small diagram $\cat I\to \cat C$ admits a limit (resp.\ colimit) in $\cat C$.
\end{definition}

\begin{theorem}\label{rslt:Anima-Cat-adjunctions}
The quasicategories $\Ani$ and $\Cat$ are both complete and cocomplete, and one has adjunctions 
\begin{equation}
\begin{tikzcd}
\Ani \ar[r,hook] & \Cat \ar[l,shift right=2, "\lvert\blank\rvert"'] \ar[l,shift left=2, "(\blank)^\simeq"],
\end{tikzcd}
\end{equation}
where the \emph{geometric realization} $\lvert\blank\rvert$ is left adjoint and the \emph{underlying anima} (or \emph{core}) $\cat C\mapsto \cat C^\simeq$ is right adjoint to the inclusion. In particular, $\Ani \injto \Cat$ preserves limits and colimits.
\end{theorem} 
\begin{proof} 
The (co)completeness of $\Ani$ is proved in  \cite[\href{https://kerodon.net/tag/02VA}{Tag 02VA}]{kerodon}. The completeness and cocompleteness of $\Cat$ is proved in \cite[\href{https://kerodon.net/tag/02T0}{Tag 02T0}]{kerodon} and \cite[\href{https://kerodon.net/tag/02UN}{Tag 02UN}]{kerodon}, respectively. The fact that $(\blank)^\simeq$ is right adjoint to the inclusion follows from \cite[\href{https://kerodon.net/tag/01DA}{Tag 01DA}]{kerodon}. The geometric realization $\lvert \cat C\rvert$ is by \cite[\href{https://kerodon.net/tag/01MY}{Tag 01MY}]{kerodon} the localization of $\cat C$ at the collection of all edges.
\end{proof}

The final object of $\Cat$ and $\Ani$ is the ordinary category $*$ with a single object and a single endomorphism (the identity).

\begin{definition} \label{def:weakly-contractible-category}
A quasicategory $\cat C$ is called \emph{weakly contractible} if $\abs{\cat C} = *$ (cf.\ \cite[\href{https://kerodon.net/tag/04GW}{Tag 04GW}]{kerodon}).
\end{definition}

\begin{example}
If a category has an initial or a final object then it is weakly contractible (see \cite[\href{https://kerodon.net/tag/02P2}{Tag 02P2}]{kerodon}).
\end{example}

\begin{definition}
Let $f\colon \cat C \to \cat D$ be a functor between quasicategories.
\begin{defenum}
\item $f$ is called \emph{fully faithful} if for all $c,c'\in \cat C_0$ the map $\Hom_{\cat C}(c,c') \isoto \Hom_{\cat D}(fc, fc')$ is an isomorphism in $\Ani$. In this case we call $f$ an \emph{embedding} and $\cat C$ a full subcategory of $\cat D$, which we denote by $\cat C \injto \cat D$ or even $\cat C \subseteq \cat D$.

\item $\cat C$ is called a \emph{subcategory} of $\cat D$, and we write $\cat C \subset \cat D$, if for all $c,c'\in \cat C_0$ the functors $\Hom_{\cat C}(c,c') \injto \Hom_{\cat D}(fc,fc')$ and $\cat C^\simeq \injto \cat D^\simeq$ are fully faithful.\footnote{By \cite[\href{https://kerodon.net/tag/01K1}{Tag 01K1}]{kerodon} the former is equivalent to $\Hom_{\cat C}(c,c')$ being isomorphic to a union of path components of $\Hom_{\cat D}(fc,fc')$.}
\end{defenum}
\end{definition}

\begin{lemma}\label{rslt:limit-of-subcategories}
The full subcategory of $\Fun([1],\Cat)$ spanned by the subcategories (resp.\ embeddings) $\cat C\to \cat D$ is closed under limits.
\end{lemma}
\begin{proof}
We first show that the category of full subcategories is closed under limits in $\Fun([1],\Cat)$. Let $\cat I$ be a small quasicategory and $\cat I\to \Fun([1],\Cat)$ a small diagram whose components $f_i\colon \cat C_i\to \cat D_i$ exhibit $\cat C_i$ as a full subcategory of $\cat D_i$. Put $\cat C = \lim_i\cat C_i$ and $\cat D = \lim_i\cat D_i$ and $f\coloneqq \lim_if_i \colon \cat C\to \cat D$. We claim that $f$ exhibits $\cat C$ as a full subcategory of $\cat D$. Let $c,c' \in \cat C$ with respective images $c_i$ and $c'_i$ in $\cat C_i$. Then $\Hom_{\cat C}(c,c') = \lim_i \Hom_{\cat C_i}(c_i,c'_i)$ (cf.~\cite[\href{https://kerodon.net/tag/01J4}{Tag 01J4}]{kerodon}), and similarly for $\Hom_{\cat D}(fc,fc')$. Since each $\Hom_{\cat C_i}(c_i,c'_i) \isoto \Hom_{\cat D_i}(f_ic_i, f_ic'_i)$ is an isomorphism in $\Ani$, it follows that $\Hom_{\cat C}(c,c') \isoto \Hom_{\cat D}(fc,fc')$ is an isomorphism. Hence $\cat C$ is a full subcategory of $\cat D$.

The argument for (not necessarily full) subcategories is similar. With the same notation as above, one has to check that the map $\lim_i \Hom_{\cat C_i}(c_i,c_i) \to \lim_i \Hom_{\cat D_i}(f_ic_i, f_ic'_i)$ is fully faithful, which follows from the fully faithfulness of each $\Hom_{\cat C_i}(c_i, c'_i) \injto \Hom_{\cat D_i}(f_ic_i, f_ic'_i)$ and the discussion above.
\end{proof}

\begin{proposition}\label{rslt:subcategories-are-monomorphisms}
For a functor $f\colon \cat C\to \cat D$ the following are equivalent:
\begin{propenum}
\item $f$ exhibits $\cat C$ as a subcategory of $\cat D$.
\item $f$ is a \emph{monomorphism}, that is, the induced functor $\Hom_{\Cat}(\cat B, \cat C) \injto \Hom_{\Cat}(\cat B,\cat D)$ is fully faithful for any quasicategory $\cat B$.
\item The diagonal map $\cat C \isoto \cat C\times_{\cat D}\cat C$ is an isomorphism.
\end{propenum}
\end{proposition}
\begin{proof}
Combine \cite[\href{https://kerodon.net/tag/04W5}{Tag 04W5}]{kerodon} (where our notion of subcategory coincides with (2)) and \cite[\href{https://kerodon.net/tag/04VW}{Tag 04VW}]{kerodon}. 
\end{proof}

\begin{corollary}\label{rslt:fiber-product-over-subcategory}
Let $\cat C \subset \cat D$ be a subcategory and consider functors $\cat A\to \cat C$ and $\cat B\to \cat C$. Then the natural map $\cat A\times_{\cat C}\cat B \isoto \cat A \times_{\cat D}\cat B$ is an isomorphism of categories.
\end{corollary}
\begin{proof}
The diagram
\begin{equation}
\begin{tikzcd}
\cat A\times_{\cat C}\cat B \ar[d] \ar[r] \ar[dr,phantom, very near start, "\lrcorner"] & \cat A\times_{\cat D}\cat B \ar[d]\\
\cat C \ar[r,"\simeq", "\Delta"'] & \cat C\times_{\cat D}\cat C
\end{tikzcd}
\end{equation}
is cartesian. Now, $\Delta$ is an isomorphism by \cref{rslt:subcategories-are-monomorphisms}, hence so is the top map.
\end{proof}

The notion of subcategories considered above captures monomorphisms in $\Cat$, as shown in \cref{rslt:subcategories-are-monomorphisms}. There is a more general notion of $n$-truncated morphisms in a quasicategory, which we introduce now. Monomorphisms are then the special case $(-1)$-truncated morphisms.

\begin{definition}\label{defn:truncated}
Let $\cat C$ be a quasicategory which admits fiber products. An edge $f\colon Y\to X$ is called \emph{$(-2)$-truncated} if it is an isomorphism. For $n\ge -2$ we say that $f$ is \emph{$n$-truncated} if the diagonal $\Delta_f\colon Y\to Y\times_XY$ is $(n-1)$-truncated. We say that $f$ is \emph{truncated} if it is $n$-truncated for some $n\ge-2$.
\end{definition}

\begin{remark}
By \cite[\href{https://kerodon.net/tag/05FT}{Tag 05FT}]{kerodon}, our definition of truncated morphisms is equivalent to the rather technical definition in \cite[\href{https://kerodon.net/tag/05F9}{Tag 05F9}]{kerodon} (but which does not require $\cat C$ to admit fiber products).
\end{remark}

\begin{proposition}\label{rslt:truncated-properties}
Let $\cat C$ be a quasicategory.
\begin{propenum}
\item\label{rslt:truncated-monotonicity}
If an edge $f$ is $n$-truncated, then it is $m$-truncated for all $m\ge n$.

\item\label{rslt:truncated-composition}
Let $f$, $g$ be composable edges in $\cat C$.
\begin{enumerate}[label=(\alph*)]
\item If $f$ and $g$ are $n$-truncated, then $fg$ is $n$-truncated. 
\item If $fg$ is $n$-truncated and $f$ is $(n+1)$-truncated, then $g$ is $n$-truncated.
\end{enumerate}

\item\label{rslt:truncated-pullback}
The collection of $n$-truncated morphisms is closed under pullbacks along arbitrary morphisms in $\cat C$.
\end{propenum}
\end{proposition}
\begin{proof}
See, respectively, \cite[\href{https://kerodon.net/tag/05FC}{Tag 05FC}]{kerodon}, \cite[\href{https://kerodon.net/tag/05FQ}{Tag 05FQ}]{kerodon}, and \cite[\href{https://kerodon.net/tag/05FR}{Tag 05FR}]{kerodon}.
\end{proof}

\subsection{Cocartesian fibrations, Straightening and the Yoneda lemma} \label{sec:cat.str-yoneda}

In the present subsection we recall two core constructions in category theory: the straightening/unstraightening equivalence (classically also called the Grothendieck construction) and the Yoneda lemma. We will discuss these results on the level of $\infty$-categories, which makes the straightening construction more elegant (as it avoids talking about pseudofunctors).

From now on, we will fully adopt the $\infty$-categorical language and in particular call an $\infty$-category (i.e.\ a quasicategory) simply a category. Up to a few exceptions in this subsection, we will not rely on a particular model for $\infty$-categories (like quasicategories or complete Segal spaces), but instead work internally in $\Cat$. This is enabled by the foundational work in \cite{HTT,kerodon}. For example, (co)limits in $\Cat$ are understood as their $\infty$-categorical version (as defined in \cref{def:lim-and-colim-in-category}) rather than in some ambient ordinary category like $\sSet$.

We start our exhibition with a few basic definitions about slice and join categories, which are used throughout the paper.

\begin{definition} 
Let $\cat C$ be a category.
\begin{defenum}
\item\label{def:over-category} For an object $X$ in $\cat C$, we define the \emph{over-category} as the pullback 
\begin{align}
	\cat C_{/X} \coloneqq \Fun([1],\cat C) \times_{\cat C} \{X\},
\end{align}
where the map $\Fun([1], \cat C) \to \cat C$ is given by evaluation at $1$. Similarly, the \emph{under-category} $\cat C_{X/}$ is defined using evaluation at $0$.

\item More generally, if $F\colon \cat I\to \cat C$ is a functor, we put $\cat C_{/F} \coloneqq \cat C \times_{\Fun(\cat I,\cat C)} \Fun(\cat I,\cat C)_{/F}$ and define $\cat C_{F/}$ similarly. Note that there are canonical forgetful functors $\cat C_{/F} \to \cat C$ and $\cat C_{F/} \to \cat C$.

\item \label{def:join-of-categories} If $\cat D$ is another category, we define the \emph{join} as the pushout
\begin{align}
	\cat C \star \cat D \coloneqq \cat C \bigdunion_{\cat C\times\{0\}\times \cat D} \cat C\times[1] \times\cat D \bigdunion_{\cat C\times\{1\} \times \cat D} \cat D.
\end{align}
Note that the join comes with canonical maps $\cat C \injto \cat C\star \cat D \injfrom \cat D$. 
Specifically, we put $\cat C^{\triangleright} \coloneqq \cat C \star [0]$, which freely adjoins a terminal object to $\cat C$. Similarly, $\cat C^{\triangleleft} \coloneqq [0] \star \cat C$ freely adjoins an initial object to $\cat C$.
\end{defenum}
\end{definition} 

\begin{remarks} 
\begin{remarksenum}
	\item The fiber of the forgetful map $\cat C_{X/} \to \cat C$ at an object $Y$ (as well as the fiber of $\cat C_{/Y} \to \cat C$ at $X$) identifies with the mapping anima $\Hom_{\cat C}(X,Y)$.

	\item \label{rslt:universal-property-of-slice-category} The category $\cat C_{/F}$ has the following universal property: For any category $\cat D$ one has
	\begin{align}
		\Fun(\cat D, \cat C_{/F}) = \Fun(\cat D \star \cat I, \cat C) \times_{\Fun(\cat I,\cat C)}\{F\},
	\end{align}
	where the map $\Fun(\cat D\star \cat I, \cat C) \to \Fun(\cat I,\cat C)$ is induced by the inclusion $\cat I \injto \cat D\star \cat I$. Hence, an object of $\cat C_{/F}$ can be identified with an extension of $F$ to a functor $\cat I^{\triangleleft} \to \cat C$. Observe that under the above isomorphism the identity on $\cat C_{/F}$ corresponds to a canonical functor $\cat C_{/F}\star \cat I \to \cat C$ whose restriction to $\cat I$ is given by $F$.
	A similar property holds for $\cat C_{F/}$.

	\item \label{rslt:transitivity-of-slices} Let $F\colon \cat I \to \cat C$ and $G\colon \cat J \to \cat C_{/F}$ be diagrams. Then by the previous remark $G$ corresponds to a diagram $G'\colon \cat J \star \cat I \to \cat C$ and there is a natural equivalence
	\begin{align}
		(\cat C_{/F})_{/G} = \cat C_{/G'}.
	\end{align}
	Indeed, by the previous remark we compute, for every category $\cat D$,
	\begin{align}
		\Fun(\cat D, (\cat C_{/F})_{/G}) &= \Fun(\cat D \star \cat J, \cat C_{/F}) \times_{\Fun(\cat J, \cat C_{/F})} \{ G \}\\
		&= (\Fun(\cat D \star \cat J \star \cat I, \cat C) \times_{\Fun(\cat I, \cat C)} \{ F \}) \times_{(\Fun(\cat J \star \cat I, \cat C) \times_{\Fun(\cat I, \cat C)} \{ F \})} \{ G \}\\
		&= \Fun(\cat D \star \cat J \star \cat I, \cat C) \times_{\Fun(\cat J \star \cat I, \cat C)} \{ G' \}\\
		&= \Fun(\cat D, \cat C_{/G'}),
	\end{align}
	where in the third step we used that $G'$ is exactly the image of the map $\{ G \} \to \Fun(\cat J \star \cat I, \cat C)$.
\end{remarksenum}
\end{remarks}

\begin{lemma}
Let $F\colon \cat I \to \cat C$ be a functor and $G\colon \cat J \to \cat C_{/F}$ a diagram.
\begin{lemenum}
	\item \label{rslt:colimits-in-overcategory} If the projection of $G$ to $\cat C$ has a colimit then so does $G$ and the projection $\cat C_{/F} \to \cat C$ preserves that colimit.

	\item \label{rslt:limits-in-overcategory} If the induced diagram $G'\colon \cat J \star \cat I \to \cat C$ has a limit in $\cat C$ then $G$ has a limit and the projection $\cat C_{/F} \to \cat C$ sends the limit of $G$ to the limit of $G'$.
\end{lemenum}
\end{lemma}
\begin{proof}
Part (i) is \cite[\href{https://kerodon.net/tag/02KB}{Tag 02KB}]{kerodon}. For part (ii) we observe that a limit of $G$ corresponds to a final object of $(\cat C_{/F})_{/G}$, while a limit of $G'$ corresponds to a final object of $\cat C_{/G'}$. Thus the claim follows from \cref{rslt:transitivity-of-slices}.
\end{proof}

We now come to the discussion of the straightening/unstraightening correspondence. The idea is to encode a functor $F\colon \cat C \to \Cat$ from some category $\cat C$ via a \emph{fibration} $\cat E \to \cat C$, i.e.\ via a certain functor from another category $\cat E$ to $\cat C$. For every object $S \in \cat C$ the fiber $\cat E_S \coloneqq \cat E \times_{\cat C} \{ S \}$ identifies with the category $F(S)$. Moreover, fixing a morphism $\alpha\colon T \to S$ in $\cat C$ and objects $X \in \cat E_S = F(S)$ and $Y \in \cat E_T = F(T)$, then a morphism $f\colon Y \to X$ in $\cat E$ lying over $\alpha$ is the same as a morphism $F(\alpha)(Y) \to X$ in $F(S)$, i.e.
\begin{align}
	\Hom_\alpha(Y, X) = \Hom_{F(X)}(F(\alpha)(Y), X),
\end{align}
where $\Hom_\alpha(Y, X) \coloneqq \Hom(Y, X) \times_{\Hom(T, S)} \{ \alpha \}$. One can recover the functor $F$ from the fibration $\cat E \to \cat C$ and we will see below (see \cref{rslt:Straightening-Unstraightening}) that this provides an equivalence between the category of functors $\cat C \to \Cat$ and the category of \enquote{fibrations} $\cat E \to \cat C$. Before we come to this equivalence, we need to introduce the correct notion of \enquote{fibration} (i.e.\ we need to specify which maps $\cat E \to \cat C$ we allow):

\begin{definition} \label{def:cocartesian-morphism-and-fibration}
Let $F\colon \cat E \to \cat C$ be a functor.
\begin{defenum}
\item \label{def:cocartesian-morphism}
An edge $f\colon X \to Y$ in $\cat E$ is called \emph{$F$-cocartesian} if for all $Z\in \cat E$ the diagram
\begin{equation}\label{eq:cocartesian-definition}
\begin{tikzcd}
	\Hom_{\cat E}(Y,Z) \ar[d,"F"'] \ar[r,"\comp f"] \ar[dr,phantom,very near start, "\lrcorner"]
	&
	\Hom_{\cat E}(X,Z) \ar[d,"F"]
	\\
	\Hom_{\cat C}\bigl(F(Y), F(Z)\bigr) \ar[r,"\comp F(f)"']
	&
	\Hom_{\cat C}\bigl(F(X), F(Z)\bigr)
\end{tikzcd}
\end{equation}
is a pullback in $\Ani$. (See \cref{rmk:cocartesian-precise-formulation} for a more precise formulation.) Dually, \emph{$F$-cartesian} maps are defined.

\item \label{def:cocartesian-fibration}
We call $F$ a \emph{cocartesian fibration} if the diagram
\begin{equation}\begin{tikzcd}
	\Fun([1],\cat E)^{\cocart} \ar[r,"F\comp"] \ar[d,"\ev_0"'] \ar[dr,phantom, very near start, "\lrcorner"] & \Fun([1], \cat C) \ar[d,"\ev_0"]\\
	\cat E \ar[r,"F"'] & \cat C
\end{tikzcd}\end{equation}
is a pullback in $\Cat$, where $\Fun([1],\cat E)^{\cocart} \subseteq \Fun([1], \cat E)$ is the full subcategory consisting of $F$-cocartesian maps. Dually, one defines \emph{$F$-cartesian fibrations}.
\end{defenum}
\end{definition}

\begin{remarks} 
\begin{remarksenum}
\item\label{rmk:cocartesian-precise-formulation} Our definitions are essentially taken from \cite[\S2.1]{Haugseng.2023}, which also takes a model-independent approach. 
It is \emph{a priori} not obvious how to construct the commutative diagram \eqref{eq:cocartesian-definition}. We give here a precise (but rather technical and model dependent) description: Let $\cat E$ be a quasicategory, $f\colon X\to Y$ an edge and $Z$ an object in $\cat E$. We define the simplicial set $\Hom_{\cat E}(f,Z)$ by the pullback
\begin{equation}
\begin{tikzcd}
\Hom_{\cat E}(f,Z) \ar[d] \ar[r] \ar[dr,phantom, very near start, "\lrcorner"] & \Fun(\Delta^2, \cat E) \ar[d,"{(d_2,\ev_2)}"]\\
\Delta^0 \ar[r,"{(f,Z)}"'] & \Fun(\Delta^1,\cat E) \times \cat E
\end{tikzcd}
\end{equation}
in $\sSet$. Since $(d_2,\ev_0)$ is an isofibration (e.g.\ by \cite[\href{https://kerodon.net/tag/01G3}{Tag 01G3}]{kerodon} applied with $A = \Delta^{\{0,1\}} \dunion \Delta^{\{2\}}$ and $B = \Delta^2$), this is also a pullback in $\Cat$. Since $\Horn^2_1 \subseteq \Delta^2$ is inner anodyne, the map $d_0\colon \Hom_{\cat E}(f,Z) \to \Hom_{\cat E}(Y,Z)$ is a trivial fibration (\cite[\href{https://kerodon.net/tag/01BW}{Tag 01BW}]{kerodon}); in particular $\Hom_{\cat E}(f,Z)$ is an anima. If moreover $F\colon \cat E\to \cat C$ is a functor, then we obtain a commutative diagram
\begin{equation}
\begin{tikzcd}
\Hom_{\cat E}(Y,Z) \ar[d,"F"'] & \Hom_{\cat E}(f,Z) \ar[l,"\simeq"'] \ar[d,"F"] \ar[r,"d_1"] & \Hom_{\cat E}(X,Z) \ar[d,"F"]\\
\Hom_{\cat C}\bigl(F(Y),F(Z)\bigr) & \Hom_{\cat E}\bigl(F(f),F(Z)\bigr) \ar[l,"\simeq"] \ar[r,"d_1"'] & \Hom_{\cat E}\bigl(F(X),F(Z)\bigr)
\end{tikzcd}
\end{equation}
in $\sSet$. As the left horizontal arrows become invertible in $\Ani$, this yields the commutativity of the diagram \eqref{eq:cocartesian-definition}.

\item Lurie gives a slightly different definition of $F$-(co)cartesian maps and (co)cartesian fibrations, see \cite[\href{https://kerodon.net/tag/01TF}{Tag 01TF}]{kerodon} and \cite[\href{https://kerodon.net/tag/01UA}{Tag 01UA}]{kerodon}. The main difference is that our notion of (co)cartesian fibration is invariant under isomorphisms of categories: If $F\colon \cat E\to \cat C$ is a functor between quasicategories and $\cat E \isoto \cat E' \xto{F'} \cat C$ is any factorization of $F$ into a categorical equivalence and an isofibration $F'$, then $F$ (and hence $F'$) is a (co)cartesian fibration in our sense if and only if $F'$ is a (co)cartesian fibration in the sense of Lurie \cite[Corollary~3.4]{Mazel-Gee.2015}.

\item \label{rslt:cocartesian-fib-equiv-existence-of-lifts} One can show that the natural map $\Fun([1],\cat E)^{\cocart} \to \cat E \times_{\cat C} \Fun([1],\cat C)$ is always fully faithful. Hence, $F$ is a cocartesian fibration if and only if this map is essentially surjective. Concretely, this means that for all $X\in\cat E$ and $\ol{f}\colon \ol X \to \ol Y$ in $\cat C$ with $F(X) \simeq \ol X$, there exists an $F$-cocartesian edge $f \colon X\to Y$ in $\cat E$ such that $F(f) \simeq \ol{f}$ in $\Fun([1],\cat C)$.
\end{remarksenum}
\end{remarks} 

\begin{example}\label{ex:cocartesian-fibration} 
Let $\cat C$ be a category.
The functor $\ev_1\colon \Fun([1],\cat C) \to \cat C$ is a cocartesian fibration; it is a cartesian fibration if and only if $\cat C$ admits fiber products. 

Similarly, $\ev_0\colon \Fun([1], \cat C)\to \cat C$ is a cartesian fibration; it is a cocartesian fibration if and only if $\cat C$ admits pushouts.
\end{example} 

\begin{proposition}
Let $F\colon \cat E \to \cat D$ and $G\colon \cat D\to \cat C$ be cocartesian fibrations. Then:
\begin{propenum}
\item\label{rslt:composition-of-cocartesian-fibrations} 
The composite $G\comp F \colon \cat E\to \cat C$ is a cocartesian fibration.
\item\label{rslt:morphism-of-cocartesian-fibrations} 
$F$ is a morphism of cocartesian fibrations, i.e.\ $F$ sends $(G\comp F)$-cocartesian morphisms to $G$-cocartesian morphisms.
\end{propenum}
\end{proposition}
\begin{proof}
This is the dual of \cite[\href{https://kerodon.net/tag/01UL}{Tag 01UL}]{kerodon}.
\end{proof}

\begin{proposition}
\label{rslt:pullback-of-cocartesian-fibration}
Consider a pullback diagram
\begin{equation}
\begin{tikzcd}
\cat E' \ar[d,"F'"'] \ar[r,"G"] \ar[dr,phantom, very near start, "\lrcorner"] & \cat E \ar[d,"F"]\\
\cat C' \ar[r] & \cat C
\end{tikzcd}
\end{equation}
of categories, where $F$ is a cocartesian fibration. Then $F'$ is a cocartesian fibration, and a map $f'$ in $\cat E'$ is $F'$-cocartesian if and only if $G(f')$ is $F$-cocartesian.
\end{proposition}
\begin{proof}
See \cite[\href{https://kerodon.net/tag/01UF}{Tag 01UF}]{kerodon}.
\end{proof}

\begin{theorem}\label{rslt:exponentiation-of-cocartesian-fibrations}
Let $F\colon \cat E\to \cat C$ be a cocartesian fibration, and let $\cat A$ be a category. Then the induced functor $F'\colon \Fun(\cat A, \cat E) \to \Fun(\cat A, \cat C)$ is a cocartesian fibration. Moreover, an edge $f$ in $\Fun(\cat A, \cat E)$ is $F'$-cocartesian if and only if $\ev_A(f)$ is $F$-cocartesian for all $A\in \cat A$, where $\ev_A\colon \Fun(\cat A, \cat E) \to \cat E$ is the evaluation at $A$. 
\end{theorem}
\begin{proof}
See \cite[\href{https://kerodon.net/tag/01VG}{Tag 01VG}]{kerodon}. 
\end{proof}

\begin{proposition}\label{rslt:lifting-property-of-cocartesian-fibrations}
Let $F\colon \cat E \to \cat C$ be a cocartesian fibration.
Let $\Fun'(\Horn^2_0, \cat E) \subseteq \Fun(\Horn^2_0,\cat E)$ be the full subcategory of functors whose restriction to $[1] \simeq \{0, 1\} \subseteq \Horn^2_0$ classifies an $F$-cocartesian edge. Define $\Fun'([2], \cat E) \subseteq \Fun([2], \cat E)$ analogously. Then the functor $\Horn^2_0 \to [2]$ induces an isomorphism
\[
\Fun'([2], \cat E) \isoto \Fun'(\Horn^2_0, \cat E) \times_{\Fun(\Horn^2_0, \cat C)} \Fun([2], \cat C).
\]
\end{proposition}
\begin{proof}
We prove a more general (but model dependent) version: Let $F\colon \cat E \to \cat C$ be a cocartesian fibration of simplicial sets \cite[Definition~2.4.2.1]{HTT} and let $n\ge2$. We denote by $\Fun'(\Horn^n_0, \cat E)$ and $\Fun'(\Delta^n, \cat E)$ the full simplicial subsets of $\Fun(\Horn^n_0, \cat E)$ and $\Fun(\Delta^n, \cat E)$, respectively, generated by those vertices whose restriction to $\Delta^1 \simeq \{0, 1\} \subseteq \Horn^n_0 \subseteq \Delta^n$ classifies an $F$-cocartesian edge. Then the obvious map
\[
\Fun'(\Delta^n, \cat E) \to \Fun'(\Horn^n_0, \cat E) \times_{\Fun(\Horn^n_0, \cat C)} \Fun(\Delta^n, \cat C)
\]
is a trivial fibration of simplicial sets.

For the proof we will employ the theory of marked simplicial sets, see \cite[\S 3.1]{HTT}. Recall that a marked simplicial set is a pair $(X, E)$ consisting of a simplicial set $X$ and a collection $E\subseteq X_1$ of edges which contains the degenerate edges. The marked simplicial sets assemble into an ordinary category $\sSet^+$, where the morphisms are morphisms of simplicial sets which send marked edges to marked edges. The forgetful functor $\sSet^+\to \sSet$ admits a left adjoint $X\mapsto X^\flat$, where only the degenerate edges are marked, and a right adjoint $X\mapsto X^\sharp$, where all edges are marked. If $p\colon X\to S$ is a cocartesian fibration of simplicial sets, we denote $X^\natural$ the marked simplicial set $X$ with all $p$-cocartesian edges marked. For marked simplicial sets $(X,E_X), (Y, E_Y)$, we denote by $\Fun((X,E_X), (Y,E_Y))$ the full simplicial subset of $\Fun(X,Y)$ spanned by those vertices which send $E_X$ into $E_Y$.

Let now $E\subseteq (\Horn^n_0)_1$ be the subset containing $0\to 1$ and all degenerate edges. Translating the claim into the language of marked simplicial sets, we need to show that
\begin{equation}
\Fun\bigl((\Delta^n, E), \cat E^\natural\bigr) \to \Fun\bigl((\Lambda_0^n, E), \cat E^\natural\bigr) \times_{\Fun((\Horn^n_0, E), \cat C^\sharp)} \Fun\bigl((\Delta^n, E), \cat C^\sharp\bigr)
\end{equation}
is a trivial fibration of simplicial sets. We thus need to check the right lifting property against the inclusions $\partial\Delta^r\to \Delta^r$ for all $r\ge0$. Unraveling, we need to solve the following lifting problem:
\begin{equation}
\begin{tikzcd}
(\partial\Delta^r)^\flat \times (\Delta^n,E) \coprod\limits_{(\partial\Delta^r)^\flat \times (\Horn^n_0,E)} (\Delta^r)^\flat\times (\Horn^n_0,E)
\ar[d,hook, "i"'] \ar[r] &
\cat E^\natural \ar[d,"F"] \\
(\Delta^r)^\flat \times (\Delta^n,E) 
\ar[r] \ar[ur,dotted] &
\cat C^\sharp.
\end{tikzcd}
\end{equation}
By (the dual of) \cite[Proposition~3.1.2.3]{HTT}, the map $i$ is left marked anodyne. Since $F\colon \cat E\to \cat C$ is a cocartesian fibration, (the dual of) \cite[Proposition~3.1.1.6]{HTT} shows that the dotted arrow exists.
\end{proof}

\begin{definition}\label{def:left-fibration}
Let $F\colon \cat E\to \cat C$ be a functor. We call $F$ a \emph{left fibration} if $F$ is a cocartesian fibration and the following equivalent conditions are satisfied:\footnote{The equivalences only hold under the assumption that $F$ is cocartesian.}
\begin{defenum}
\item $F$ is \emph{conservative}, i.e.\ a morphism $f$ in $\cat E$ is an isomorphism as soon as $F(f)$ is an isomorphism in $\cat C$;
\item For each $X\in \cat C$, the fiber $\{X\} \times_{\cat C}\cat E$ is an anima;
\item Every morphism in $\cat E$ is cocartesian.
\end{defenum}
Similarly, $F$ is called a \emph{right fibration} if $F$ is a cartesian fibration satisfying the equivalent conditions (a), (b) and
\begin{defenum}[label=(c')]
\item Every morphism in $\cat E$ is cartesian.
\end{defenum}
\end{definition} 

\begin{example} 
Let $\cat C$ be a category containing an object $X$.
It follows from \cref{ex:cocartesian-fibration} that the forgetful functor $\cat C_{/X} \to \cat C$ is a right fibration.
Similarly, the forgetful functor $\cat C_{X/} \to \cat C$ is a left fibration.
\end{example} 

Cartesian and cocartesian fibrations over a category $\cat C$ can be organized into a subcategory of the slice category $\Cat_{/\cat C}$:

\begin{definition}
Let $\cat C$ be a category. We denote by $\Cat_{/\cat C}^{\cocart} \subset \Cat_{/\cat C}$ the (non-full) subcategory consisting of the cocartesian fibrations together with those functors which preserve cocartesian morphisms. Let $\Cat_{/\cat C}^{\lfib}$ be the full subcategory of $\Cat_{/\cat C}$ (and of $\Cat_{/\cat C}^{\cocart}$) of left fibrations. The categories $\Cat_{/\cat C}^{\cart}$ and $\Cat_{/\cat C}^{\rfib}$ of cartesian fibrations and right fibrations, respectively, are defined similarly. 
\end{definition}

\begin{lemma} \label{rslt:cocartesian-morphism-in-limit-of-fibrations}
Let $\cat C$ be a category and $(\cat E_i)_i$ a diagram in $\Cat_{/\cat C}$ with limit $\cat E \in \Cat_{/\cat C}$. Let $\ol{f}\colon \ol X \to \ol Y$ be a map in $\cat C$ and $X = (X_i)_{i\in I}$ an object in the fiber $\cat E_{\ol X}$ such that $\ol f$ admits a cocartesian lift $f_i\colon X_i \to Y_i$ in $\cat E_i$ for all $i$. Then the $f_i$'s assemble to a cocartesian lift $f\colon X \to Y$ of $\ol f$ in $\cat E$.
\end{lemma}
\begin{proof}
By \cref{rslt:limits-in-overcategory} the limit $\cat E$ exists and is computed as the limit in $\Cat$ of the diagram $I^{\triangleright} \to \Cat$ which sends $i$ to $\cat E_i$ and $*$ to $\cat C$. By \cref{rslt:cocartesian-fib-equiv-existence-of-lifts}, evaluation at the source defines a fully faithful functor
\begin{align}
	\Fun_{\cat C}(\Delta^1, \cat E_i) \supseteq \Fun_{\cat C}^\cocart(\Delta^1, \cat E_i) \injto (\cat E_i)_{\ol X}
\end{align}
and similarly for $\cat E$ in place of $\cat E_i$; here $\Fun_{\cat C}$ denotes functors over $\cat C$, where we view $\Delta^1$ as a category over $\cat C$ via $\ol f$. By assumption $X_i$ lies in the essential image and is represented by $f_i$ on the left. Passing to the limit over $I$ we obtain a fully faithful embedding
\begin{align}
	\Fun_{\cat C}(\Delta^1, \cat E) \supseteq \varprojlim_i \Fun_{\cat C}^\cocart(\Delta^1, \cat E_i) \injto \varprojlim_i (\cat E_i)_{\ol X} = \cat E_{\ol X}.
\end{align}
Now $X = (X_i)_i$ lies in this full subcategory (because this is true termwise for each $i$) and hence corresponds to a map $f\colon X \to Y$ which is assembled from the maps $f_i$. It remains to show that $f$ is cocartesian over $\cat C$. To this end we note that for every $Z = (Z_i)_i \in \cat E$ lying over some $\ol Z \in \cat C$ we have $\Hom_{\cat E}(Y, Z) = \varprojlim_i \Hom_{\cat E_i}(Y_i, Z_i)$ in $\Ani_{/\Hom_{\cat C}(\ol Y, \ol Z)}$ by \cref{rslt:limits-in-overcategory}, so the claim follows immediately from the definitions.
\end{proof}

\begin{corollary} \label{rslt:cocartesian-fibrations-stable-under-limits}
Let $\cat C$ be a category. The category $\Cat_{/\cat C}^{\cocart}$ has all small limits and the inclusion into $\Cat_{/\cat C}$ preserves them.
\end{corollary}
\begin{proof}
Let $\cat E = \varprojlim_i \cat E_i$ be a limit in $\Cat_{/\cat C}$, where each $\cat E_i \to \cat C$ is a cocartesian fibration. Then by \cref{rslt:cocartesian-morphism-in-limit-of-fibrations} $\cat E \to \cat C$ is a cocartesian fibration. It remains to show that for another cocartesian fibration $\cat E' \to \cat C$ and a functor $\alpha\colon \cat E' \to \cat E$ over $\cat C$, $\alpha$ preserves cocartesian edges if and only if all the compositions $\cat E' \to \cat E_i$ preserve cocartesian edges. The \enquote{only if} part follows from \cref{rslt:composition-of-cocartesian-fibrations} and the \enquote{if} part can be proved as in the last part of \cref{rslt:cocartesian-morphism-in-limit-of-fibrations}.
\end{proof}

We can now state the straightening/unstraightening correspondence, as promised above. It is due to Lurie \cite[Theorem~3.2.0.1]{HTT} and further simplified by Hebestreit--Heuts--Ruit \cite{Hebestreit-Heuts-Ruit.2021}.

\begin{theorem}[Straightening/Unstraightening] 
\label{rslt:Straightening-Unstraightening}
There are inverse isomorphisms of categories
\begin{equation}
\begin{tikzcd}[column sep=large]
\Cat_{/\cat C}^{\cocart} \ar[r,phantom,"\sim", yshift=-.7pt] \ar[r,shift left, "\Str^{\cocart}"] & \ar[l,shift left, "\Un^{\cocart}"] \Fun(\cat C, \Cat)
\end{tikzcd}
\quad\text{and}\quad
\begin{tikzcd}[column sep=large]
\Cat_{/\cat C}^{\cart} \ar[r,phantom,"\sim", yshift=-.7pt] \ar[r,shift left, "\Str^{\cart}"] & \ar[l,shift left, "\Un^{\cart}"] \Fun(\cat C^\op, \Cat),
\end{tikzcd}
\end{equation}
which respectively restrict to inverse isomorphisms
\begin{equation}
\begin{tikzcd}[column sep=large]
\Cat_{/\cat C}^{\lfib} \ar[r,phantom,"\sim", yshift=-.7pt] \ar[r,shift left, "\Str^\lfib"] & \ar[l,shift left, "\Un^\lfib"] \Fun(\cat C, \Ani)
\end{tikzcd}
\quad\text{and}\quad
\begin{tikzcd}[column sep=large]
\Cat_{/\cat C}^{\rfib} \ar[r,phantom,"\sim", yshift=-.7pt] \ar[r,shift left, "\Str^{\rfib}"] & \ar[l,shift left, "\Un^{\rfib}"] \Fun(\cat C^\op, \Ani).
\end{tikzcd}
\end{equation}
\end{theorem}

We have already explained above how the straightening/unstraightening correspondence encodes a functor $\cat C \to \Cat$ via a cocartesian fibration $\cat E \to \cat C$ (see the discussion before \cref{def:cocartesian-morphism-and-fibration}). Let us explain why this result is powerful and where the terminology comes from. Giving a functor $F\colon \cat C\to \Cat$ amounts to the following infinite data:
\begin{itemize}
	\item for each $X\in \cat C$ a category $F(X)$;
	\item for all morphisms $f\colon X\to Y$ in $\cat C$ a functor $F(f)\colon F(X) \to F(Y)$;
	\item for all morphisms $X\xto f Y \xto g Z$ in $\cat C$ a natural isomorphism $\alpha_{f,g}\colon F(g) \comp F(f) \isoto F(g\comp f)$;
	\item for all morphisms $W \xto f X \xto g Y \xto h Z$ a $3$-simplex of the form
	\begin{equation}\begin{tikzcd}
		& F(X) \ar[r,"F(g)"] & F(Y) \ar[dr,"F(h)"] \\
		F(W) \ar[ur,"F(f)"] \ar[rrr,"F(h\comp g\comp h)"'] \ar[urr,dashed, "F(g\comp f)" description] & & & F(Z) \ar[from=ull,crossing over, "F(h\comp g)" description],
	\end{tikzcd}\end{equation}
	where the faces are given by the natural transformations of the previous step; 
	\item $\dotsc$
\end{itemize}
It is generally impossible to write down these data explicitly. It is often easier to construct a functor $\cat E \to \cat C$ whose fiber at $X \in \cat C$ is the category $F(X)$. The straightening/unstraightening principle then says that the infinite \emph{data} needed to make $X\mapsto F(X)$ into a functor correspond to the \emph{condition} that the fibration $\cat E\to \cat C$ is cocartesian. The fibration $\cat E \to \cat C$ is an \enquote{unstraightening} of the functor $F$ in the sense that none of the homotopies making the functor $F$ (i.e.\ the transition functors $F(f)$, the isomorphisms $\alpha_{f,g}$ etc.) are chosen explicitly in the fibration. Instead, the straightening/unstraightening correspondence says that all of these data \emph{can be} chosen and such a choice is essentially unique.

\begin{remarks} 
\begin{remarksenum}
\item The straightening/unstraightening constructions in Theorem~\ref{rslt:Straightening-Unstraightening} can be made functorial in $\cat C$ \cite[Appendix~A]{Gepner-Haugseng-Nikolaus.2017}.

\item The left fibration $\Ani_{*/} \to \Ani$ is universal in the following sense: If $F\colon \cat C\to \Ani$ is a functor, then the corresponding left fibration is given by the base-change $\cat C\times_{\Ani}{\Ani_{*/}} \to \cat C$.

\item There exists a universal cocartesian fibration $\Cat_{\mathrm{Obj}} \to \Cat$ \cite[\href{https://kerodon.net/tag/0212}{Tag 0212}]{kerodon}. The objects of $\Cat_{\mathrm{Obj}}$ are pairs $(\cat C, X)$ consisting of a category $\cat C$ and an object $X\in \cat C$. A morphism from $(\cat C, X)$ to $(\cat D,Y)$ is a pair $(F,f)$ consisting of a functor $F\colon \cat C\to \cat D$ and an edge $f\colon F(X) \to Y$ in $\cat D$. 

Note that $\Cat_{\mathrm{Obj}}$ is not the same as $\Cat_{[0]/}$: The latter has the same objects as $\Cat_{\mathrm{Obj}}$, but the morphisms are pairs $(F,f)$ as above such that $f$ is an isomorphism.

\item The autoequivalence $\op\colon \Cat \isoto \Cat$ induces a commutative diagram
\begin{equation}\begin{tikzcd}
	\Cat_{/\cat C}^{\cocart} \ar[r,"\op"] \ar[d,shift right,"\Str^{\cocart}"']
	&
	\Cat_{/\cat C^\op}^{\cart} \ar[d,shift right, "\Str^{\cart}"']
	\\
	\Fun(\cat C,\Cat) \ar[r,"\op_*"'] \ar[u,shift right, "\Un^{\cocart}"']
	&
	\Fun(\cat C, \Cat) \ar[u, shift right, "\Un^{\cart}"'].
\end{tikzcd}\end{equation}
\end{remarksenum}
\end{remarks}

\begin{remark} \label{rmk:lax-trafos-via-unstraightening}
Apart from the fact that the straightening/unstraightening correspondence is a powerful tool to construct functors into $\Cat$, it is also useful because it provides access to some 2-categorical constructions internally in $\Cat$. Namely, given two functors $F, F'\colon \cat C \to \Cat$ with associated cocartesian fibrations $\cat E, \cat E' \to \cat C$, a natural transformation $\alpha\colon F \to F'$ is the same as a functor $A\colon \cat E \to \cat E'$ over $\cat C$ which preserves cocartesian morphisms. If we drop the requirement that $A$ preserves cocartesian edges, then we obtain the category $\Fun_{\cat C}(\cat E, \cat E')$ of \emph{lax natural transformations} $F \to F'$. Concretely, a lax natural transformation $\alpha\colon F \to F'$ consists of the following data:
\begin{itemize}
	\item For every $X \in \cat C$ a functor $\alpha_X\colon F(X) \to F'(X)$.
	\item For every morphism $f\colon Y \to X$ in $\cat C$ a 2-morphism
	\begin{equation}\begin{tikzcd}
		F(Y) \arrow[r,"\alpha_Y"] \arrow[d,"F(f)",swap] & F'(Y) \arrow[d,"F'(f)"] \arrow[dl,Rightarrow,shorten >=1.5ex]\\
		F(X) \arrow[r,"\alpha_X"] & F'(X)
	\end{tikzcd}\end{equation}
	i.e.\ a natural transformation of functors $F'(f) \comp \alpha_Y \to \alpha_X \comp F(f)$.
	\item Higher homotopies for the natural transformations, exhibiting their compatibility with composition of morphisms in $\cat C$.
\end{itemize}
One checks that these data are indeed encoded by a functor $\cat E \to \cat E'$ over $\cat C$ and that this functor preserves cocartesian morphisms if and only if the 2-morphisms above are isomorphisms.
\end{remark}

One can deduce from \cref{rslt:Straightening-Unstraightening} the following relative version, where we replace $\Cat$ with the slice $\Cat_{/\cat S}$ for a fixed base category $\cat S$:

\begin{corollary}\label{rslt:relative-Straightening}
Let $\cat C, \cat S \in \Cat$. Then straightening/unstraightening induces an equivalence
\begin{align}
	\Cat^{\cocart_{\cat C}}_{/\cat C \times \cat S} \isom \Fun(\cat C, \Cat_{/\cat S}),
\end{align}
where $\Cat^{\cocart_{\cat C}}_{/\cat C \times \cat S} \subset \Cat_{/\cat C \times \cat S}$ denotes the (non-full) subcategory consisting of those functors $\cat E \to \cat C \times \cat S$ such that $\cat E \to \cat C$ is a cocartesian fibration and $\cat E\to \cat S$ sends cocartesian edges in $\cat E$ to isomorphisms in $\cat S$, and where the morphisms are given by those functors (over $\cat C\times \cat S$) that preserve cocartesian edges over $\cat C$.

The analogous statements with \enquote{cocartesian} replaced with \enquote{cartesian} (resp.\ with $\Cat$ replaced with $\Ani$) also hold.
\end{corollary}
\begin{proof}
Recall that, by definition, we have $\Cat_{/\cat S} = \Fun([1], \Cat) \times_{\Fun(\{1\},\Cat)}\{ \cat S \}$. Therefore, we obtain natural isomorphisms
\begin{align}
\Fun(\cat C, \Cat_{/\cat S}) \cong \Fun(\cat C\times [1], \Cat) \times_{\Fun(\cat C\times\{1\}, \Cat)} \{\const_{\cat S}\} \cong \Fun(\cat C, \Cat)_{/\const_{\cat S}},
\end{align}
where $\const_{\cat S} \colon \cat C \to \Cat$ is the constant functor on $\cat S$. Note also that the unstraightening of $\const_{\cat S}$ is given by the projection $\pr\colon \cat C\times \cat S \to \cat C$. Hence, we obtain equivalences
\begin{align}
\Fun(\cat C, \Cat_{/\cat S}) \cong \Fun(\cat C, \Cat)_{/\const_{\cat S}} \cong \bigl(\Cat_{/\cat C}^{\cocart}\bigr)_{/\pr}.
\end{align}
Using \cref{rslt:transitivity-of-slices} we can identify $\bigl(\Cat_{/\cat C}^{\cocart}\bigr)_{/\pr} = \Cat^{\cocart_{\cat C}}_{/\cat C \times \cat S}$, as desired. The other cases are similar.
\end{proof}

The straightening/unstraightening equivalence gives an explicit way of computing limits and colimits in $\Ani$ and $\Cat$.

\begin{corollary}\label{rslt:(co)limits-in-Ani-Cat}
Let $\cat C$ be a small category.
\begin{corenum}
\item Given a diagram $F\colon \cat C \to \Ani$, we have
\begin{align}
\colim_{\cat C}F &= \lvert \Un^{\lfib}(F)\rvert, &\text{and} 
&&
\lim_{\cat C}F &= \Hom_{/\cat C}(\cat C, \Un^{\lfib}(F)).
\end{align}
\item Given a diagram $F\colon \cat C \to \Cat$, we have
\begin{align}
\colim_{\cat C}F &= \Un^{\cocart}(F)[\cocart^{-1}], &\text{and} 
&&
\lim_{\cat C}F &= \Fun^{\cocart}_{/\cat C}(\cat C, \Un^{\cocart}(F)),
\end{align}
where $\Fun^{\cocart}_{/\cat C}(\cat C, \Un^{\cocart}(F))$ denotes the category of sections of $p\colon \Un^{\cocart}(F) \to \cat C$ sending edges to $p$-cocartesian edges.
\end{corenum}
\end{corollary}
\begin{proof}
All formulas can be directly verified by checking the universal properties. Alternatively, see \cite[\href{https://kerodon.net/tag/02VF}{Tag 02VF}]{kerodon} and \cite[\href{https://kerodon.net/tag/02VD}{Tag 02VD}]{kerodon} for part (i) and \cite[\href{https://kerodon.net/tag/02V0}{Tag 02V0}]{kerodon} and \cite[\href{https://kerodon.net/tag/02TK}{Tag 02TK}]{kerodon} for part (ii).
\end{proof}

\begin{example}\label{example:(co)limits-in-Ani-Cat}
Let $\cat C$ be a small category
\begin{exampleenum}
\item Let $c_X\colon \cat C\to \Ani$ be the constant diagram on some anima $X$. The unstraightening of $c_X$ is then given by the projection $\cat C\times X \to \cat C$ and hence 
\begin{align}
\colim_{\cat C}c_X &= \lvert \cat C\rvert\times X, & \text{and} &&
\lim_{\cat C}c_X &= \Hom(\lvert\cat C\rvert, X).
\end{align}

\item Let $c_{\cat E} \colon \cat C\to \Cat$ be the constant diagram on some category $\cat E$. The unstraightening of $c_{\cat E}$ is then given by the projection $\cat C\times \cat E \to \cat C$ and hence 
\begin{align}
\colim_{\cat C}c_{\cat E} &= \lvert \cat C\rvert\times \cat E, & \text{and} &&
\lim_{\cat C} c_{\cat E} &= \Fun(\lvert \cat C\rvert, \cat E).
\end{align}
\end{exampleenum}
\end{example}

We now come to the second core result in this subsection: the Yoneda lemma. It allows one to embed every (small) category into a category of presheaves. In fact, the construction of this embedding is a nice application of the straightening/unstraightening correspondence:

\begin{example}
Fix a category $\cat C$. We start with the following morphism of cocartesian fibrations
\begin{equation}\begin{tikzcd}
	\Fun([1], \cat C) \ar[rr,"{(\ev_0,\ev_1)}"] \ar[dr,"\ev_1"'] & & \cat C \times \cat C \ar[dl,"\pr_2"]\\
	& \cat C.
\end{tikzcd}\end{equation}
Under cocartesian straightening we obtain a natural transformation $\alpha\colon \Str^{\cocart}(\ev_1) \to \Str^{\cocart}(\pr_2)$ of functors $\cat C\to \Cat$. 
Under the identification $\Fun([1], \Fun(\cat C, \Cat)) = \Fun(\cat C, \Fun([1], \Cat))$ the map $\alpha$ is given by sending $X \in \cat C$ to the right fibration $\cat C_{/X} \to \cat C$, and hence $\alpha$ factors through the subcategory $\Cat_{/\cat C}^{\rfib} \subseteq \Fun([1],\Cat)$. Composing $\alpha \colon \cat C\to \Cat_{/\cat C}^{\rfib}$ with $\Str^{\rfib}$ yields the desired Yoneda embedding
\begin{align}
\Yo\colon \cat C \to \Fun(\cat C^\op, \Ani) \eqqcolon \PSh(\cat C).
\end{align}
Unraveling of the definitions shows that $\Yo(X)(Y) = \Hom_{\cat C}(Y,X)$ for all $X,Y\in \cat C$. By currying, we obtain a functor
\begin{equation}\label{eq:Hom-functor}
\Hom_{\cat C}(\blank,\blank)\colon \cat C\times \cat C^\op \to \Ani.
\end{equation}
We call $\Tw(\cat C) \coloneqq \Un^{\lfib}(\Hom_{\cat C}(\blank,\blank)) \to \cat C\times \cat C^\op$ the \emph{twisted arrow category} of $\cat C$.
\end{example}

\begin{theorem}[Yoneda Lemma]\label{rslt:Yoneda-Lemma}
Let $\cat C$ be a category. The Yoneda functor $\Yo\colon \cat C \injto \PSh(\cat C)$ constructed above is fully faithful. Moreover, for every object $X \in \cat C$ the evaluation functor $\ev_X \colon \PSh(\cat C) \to \Ani$, $F \mapsto F(X)$ is corepresented by the pair $(\Yo(X), \id_X)$, i.e.\ the composite
\begin{align}
\Hom_{\PSh(\cat C)}\bigl(\Yo(X), F\bigr) \xto{\ev_X} \Hom_{\Ani}\bigl(\Hom(X,X), F(X)\bigr) \xrightarrow{\ev_{\id_X}} F(X)
\end{align}
is an isomorphism in $\Ani$.
\end{theorem}
\begin{proof}
See \cite[\href{https://kerodon.net/tag/03M5}{Tags 03M5}, \href{https://kerodon.net/tag/03NJ}{03NJ}]{kerodon}.
\end{proof}

There is another useful universal property of the Yoneda embedding $\Yo$ which exhibits $\PSh(\cat C)$ as the free cocompletion of $\cat C$. 

\begin{theorem}\label{rslt:Yoneda-cocompletion}
Let $\cat C$ be a small category. For any cocomplete category $\cat D$, restriction along $\Yo\colon \cat C \injto \PSh(\cat C)$ induces an isomorphism
\[
\Fun^L(\PSh(\cat C), \cat D) \isoto \Fun(\cat C, \cat D),
\]
where $\Fun^L(\PSh(\cat C), \cat D) \subseteq \Fun(\PSh(\cat C), \cat D)$ denotes the full subcategory of colimit-preserving functors. The inverse is given by left Kan extension along $\Yo$.
\end{theorem}
\begin{proof}
\cite[\href{https://kerodon.net/tag/04BE}{Tag 04BE}]{kerodon}.
\end{proof}

In order to demonstrate the utility of the Yoneda lemma, we deduce some functoriality properties of adjoints.

\begin{proposition}\label{rslt:right-adjoints-closed-under-limits}
Let $\cat C, \cat D$ be categories. 
\begin{enumerate}[(i)]
\item The full subcategory $\Fun^L(\cat C, \cat D)$ of $\Fun(\cat C, \cat D)$ spanned by the left adjoint functors is closed under all colimits which exist in $\Fun(\cat C, \cat D)$.

\item The full subcategory $\Fun^R(\cat D, \cat C)$ of $\Fun(\cat D, \cat C)$ spanned by the right adjoint functors is closed under all limits which exist in $\Fun(\cat D, \cat C)$.

\item There is a canonical equivalence of categories $\Fun^R(\cat D, \cat C) \isom \Fun^L(\cat C, \cat D)^\op$, which is given by passing to the adjoint.
\end{enumerate}
\end{proposition}
\begin{proof}
Part (iii) is well-known, see \cite[Proposition~5.2.6.2]{HTT}. We recall the proof in order to deduce the statements (i) and (ii) along the way. 

Note that (i) follows from (ii) because $\Fun^L(\cat C, \cat D)^\op \isom \Fun^R(\cat C^\op, \cat D^\op)$. We will now prove (ii). The Yoneda embeddings $\cat C \injto \Fun(\cat C^\op, \Ani)$ and $\cat D^\op \injto \Fun(\cat D, \Ani)$ preserve limits, and hence (e.g.\ by \cite[\href{https://kerodon.net/tag/02X9}{Tag 02X9}]{kerodon}) the induced embeddings
\begin{equation}
\Fun(\cat D, \cat C) \xinjto{\rho} \Fun(\cat D \times \cat C^\op, \Ani) \xinjfrom{\lambda} \Fun(\cat C^\op, \cat D^\op)
\end{equation}
preserve limits. Denote by $\cat R,\cat L \subseteq \Fun(\cat D\times \cat C^\op, \Ani)$ the essential images of $\rho$ and $\lambda$, respectively. Since $\cat R$ and $\cat L$ are closed under limits, so is their intersection $\cat R\cap \cat L$. It remains to prove that the map $\Fun^R(\cat D, \cat C) \injto \cat R\cap \cat L$ induced by $\rho$ is essentially surjective and hence an equivalence of categories. But note that a functor $H\colon \cat D\times \cat C^\op \to \Ani$ lies in $\cat R$ if and only if $H|_{\{d\}\times \cat C^\op}$ is representable for all $d\in \cat D$. Similarly, $H$ lies in $\cat L$ if and only if $H|_{\cat D\times \{c\}}$ is corepresentable for all $c\in \cat C^\op$. Thus, if $g \colon \cat D\to \cat C$ is a functor such that $\rho(g)$ lies in $\cat R\cap \cat L$, then for all $c\in \cat C^\op$ there exists $d\in \cat D$ and a natural isomorphism $\Hom_{\cat D}(d, \blank) = \Hom_{\cat C}(c, g(\blank))$. But it was observed in \cref{def:adjoint} that this is equivalent to $g$ being a right adjoint. This finishes the proof of (ii). Finally, the same argument as above shows that $\lambda$ induces an equivalence $\Fun^R(\cat C^\op, \cat D^\op) \isoto \cat R\cap \cat L$; the equivalence in (iii) is thus induced by composing $\rho$ with an inverse of $\lambda$.
\end{proof}

\subsection{Complete Segal anima}\label{sec:Segal-anima}

Complete Segal anima, introduced by Rezk \cite{Rezk.2001}, provide another way of looking at a category and can sometimes be useful to construct a category. The idea is to encode a category as a simplicial anima
\begin{equation}
\begin{tikzcd}
\dotsb 
\ar[r,shift left=3] \ar[r,shift left] \ar[r,shift right] \ar[r,shift right=3] 
& 
\cat C_2 
\ar[r,shift left=2] \ar[r] \ar[r,shift right=2] 
\ar[l,dashed,shorten=2mm, shift right=2]
\ar[l,dashed,shorten=2mm]
\ar[l,dashed,shorten=2mm, shift left=2]
& 
\cat C_1 
\ar[r,shift left] \ar[r,shift right] 
\ar[l,dashed,shorten=2mm,shift right]
\ar[l,dashed,shorten=2mm,shift left]
& 
\cat C_0,
\ar[l,dashed,shorten=2mm]
\end{tikzcd}
\end{equation}
where $\cat C_0$ encodes the anima of objects, $\cat C_1$ encodes the anima of morphisms together with source and target maps $\cat C_1 \rightrightarrows \cat C_0$ and a map $\cat C_0\to \cat C_1$ sending $X\in \cat C_0$ to $\id_{X}$, $\cat C_2$ encodes the anima of composites of two morphisms, etc.

We will first discuss the Segal condition on a simplicial anima, which expresses the idea that $n$-fold compositions in a category should be unique up to a contractible anima of choices.

\begin{definition}
A \emph{Segal anima} is a simplicial anima $\cat C\colon \bbDelta^\op \to \Ani$ satisfying the \emph{Segal condition}: For all $n\ge0$ the canonical map
\[
\cat C_n \xrightarrow[\sim]{(\rho^1,\dotsc,\rho^n)} \cat C_1 \times_{\cat C_0} \dotsb \times_{\cat C_0}\cat C_1
\]
is an isomorphism in $\Ani$, where $\rho^i$ is induced by the interval inclusion $[1] \cong \{i-1,i\} \injto [n]$.
\end{definition}

Thus, for a Segal anima $\cat C\colon \bbDelta^\op \to \Ani$ we can define a composition by choosing an inverse $\alpha$ of the Segal map $\cat C_2 \to \cat C_1\times_{\cat C_0}\cat C_1$ and putting
\[
\cat C_1 \times_{\cat C_0} \cat C_1 \xrightarrow{\alpha} \cat C_2 \xrightarrow{d_1} \cat C_1,
\]
where $d_1$ is induced by $[1] \cong \{0,2\} \subseteq [2]$. The spaces $\cat C_n$, for $n\ge3$, then witness the higher associativity of composition.

Observe that there is an $\infty$-categorical analog of \cref{principle:adjunctions}, which is useful for relating Segal anima with categories. 

\begin{example}
Consider the embedding $\Delta^\bullet\colon \bbDelta \injto \Cat_{\ordinary} \injto \Cat$. By \cref{rslt:Yoneda-cocompletion} it extends uniquely to a colimit preserving functor $\lvert\blank\rvert_\Delta\colon \Fun(\bbDelta^\op, \Ani) \to \Cat$. The adjoint functor theorem \cite[Corollary~5.5.2.9]{HTT} provides a right adjoint
\begin{equation}\label{eq:Segal-anima-categories-adjunction}
\begin{tikzcd}
\lvert\blank\rvert_{\Delta}\colon \Fun(\bbDelta^\op, \Ani) \ar[r,shift left] & \ar[l,shift left] \Cat \noloc \RezkNerve.
\end{tikzcd}
\end{equation}
The simplicial anima $\RezkNerve(\cat C)$ is called the \emph{Rezk nerve} of $\cat C$; it is concretely given by
\begin{align}
\RezkNerve(\cat C)_n = \Hom_{\Cat}(\Delta^n, \cat C).
\end{align}
Unsurprisingly, one can check that $\RezkNerve(\cat C)$ satisfies the Segal condition. But it satisfies another condition, which is not automatic: A functor $F\colon \Delta^1\to \cat C$ represents an isomorphism in $\cat C$ if and only if $F$ is naturally isomorphic to a composite $\Delta^1 \to \Delta^0 \to \cat C$. 
An abstraction of this idea leads to the notion of complete Segal anima, which we will discuss below.
\end{example}

We rephrase the completeness property of $\RezkNerve(\cat C)$ for a category $\cat C$ as follows: Let $J$ be the ordinary category with two objects and a unique isomorphism between them. Then $J$ is a contractible anima and $\Hom(J,\cat C) = \Hom(J, \cat C^\simeq) \isoto \Hom(\Delta^1,\cat C^\simeq)$ via any of the two inclusions $\Delta^1 \to J$. In other words, $\Hom(J,\cat C)$ identifies with the connected components of $\Hom(\Delta^1,\cat C)$ consisting of isomorphisms. Thus, $\Delta^1\to \cat C$ classifies an isomorphism if and only if it lifts uniquely (up to contractible choice) to a functor $J\to \cat C$, and so the completeness condition observed above says that the natural map $\Hom(\Delta^0,\cat C) \isoto \Hom(J,\cat C)$ induced by $J\to \Delta^0$ is an isomorphism.

\begin{definition}
Let $\cat C \colon \bbDelta^\op\to \Ani$ be a Segal anima. Then $\cat C$ is called \emph{complete} if the map
\[
\cat C_0 = \Hom(\Delta^0, \cat C) \isoto \Hom(\Nerve(J), \cat C)
\]
induced by $\Nerve(J)\to \Delta^0$ is an isomorphism in $\Ani$. (Here, the nerve $\Nerve(J)$ is viewed as a simplicial anima via the inclusion $\Set \injto \Ani$.)

We denote by $\CSegAni$ the full subcategory of $\Fun(\bbDelta^\op,\Ani)$ consisting of complete Segal anima.
\end{definition}

The following theorem is due to Joyal--Tierney \cite[Theorem~4.12]{Joyal-Tierney.2007}, but see also \cite{Hebestreit-Steinebrunner.2023} for a modern proof.

\begin{theorem}\label{rslt:Segal-anima-categories}
The Rezk nerve $\RezkNerve \colon \Cat \to \Fun(\bbDelta^\op,\Ani)$ is fully faithful with essential image $\CSegAni$. In other words, the adjunction \eqref{eq:Segal-anima-categories-adjunction} restricts to an isomorphism $\CSegAni \isoto \Cat$ of categories.
\end{theorem}

\subsection{Sheaves} \label{sec:sheaves}

Given an ordinary site $\cat C$ (e.g.\ the category of open subsets of a topological space), we have the following well-known notions associated with it:
\begin{enumerate}[(i)]
	\item A \emph{sheaf of sets} on $\cat C$ is a functor $\shv F\colon \cat C^\op \to \Set$ satisfying the following property: For every covering $(U_i \to U)_{i\in I}$ in $\cat C$ and every collection of elements $f_i \in \shv F(U_i)$ such that for all $i, j \in I$ we have $f_i|_{U_i \times_U U_j} = f_j|_{U_i \times_U U_j}$, there is a unique element $f \in \shv F(U)$ such that $f|_{U_i} = f_i$ for all $i$.

	\item Let $\shv F\colon \cat C^\op \to \Cat_\ordinary$ be a functor (where $\Cat_\ordinary$ denotes the category of ordinary categories). Given a covering $\cat U = (U_i \to U)_{i \in I}$ in $\cat C$ one introduces the category $\Desc(\cat U, \shv F)$ of \emph{descent data} of $\shv F$ along $\cat U$, which is roughly defined as follows: An object in $\Desc(\cat U, \shv F)$ consists of a family of objects $X_i \in \shv F(U_i)$ and for every pair $i, j \in I$ an isomorphism $\varphi_{ij}\colon X_i|_{U_i \times_U U_j} \isoto X_j|_{U_i \times_U U_j}$ such that for all triples $i, j, k \in I$ the isomorphisms $\varphi_{ij}, \varphi_{jk}$ and $\varphi_{ik}$ are compatible on $U_i \times_U U_j \times_U U_k$ (the latter condition is usually called the \emph{cocycle condition}). We say that $\shv F$ \emph{satisfies descent} on $\cat C$ if for every $\cat U$ as above, the natural functor
	\begin{align}
		\shv F(U) \isoto \Desc(\cat U, \shv F), \qquad X \mapsto ((X|_{U_i})_i, (\varphi_{ij})_{ij})
	\end{align}
	is an equivalence of categories.
\end{enumerate}
In the following we discuss how the above two notions can be adapted to the setting of higher categories. Note that both of the above examples generalize to the following: Let $\cat V$ be a category which has all small limits (e.g.\ $\cat V = \Set$ or $\cat V = \Cat_\ordinary$); a \emph{$\cat V$-valued sheaf} on $\cat C$ is a functor $\shv F\colon \cat C^\op \to \cat V$ such that for all coverings $(U_i \to U)_{i \in I}$ in $\cat C$ the natural morphism
\begin{align}
	\shv F(U) \isoto \varprojlim \left( \begin{tikzcd}[column sep=small,ampersand replacement=\&]
		\prod_i \shv F(U_i) \arrow[r,shift left] \arrow[r,shift right] \&
		\prod_{i,j} \shv F(U_i \times_U U_j) \arrow[r,shift left=2] \arrow[r] \arrow[r,shift right=2] \&
		\prod_{i,j,k} \shv F(U_i \times_U U_j \times_U U_k) \arrow[r,shift left=3] \arrow[r,shift left] \arrow[r,shift right] \arrow[r,shift right=3] \&
		\dots
	\end{tikzcd} \right)
\end{align}
is an isomorphism in $\cat V$. Here the limit on the right is a limit over $\bbDelta$. If $\cat V = \Set$ then this limit depends only on the first two steps, i.e.\ the above condition boils down to a bijection
\begin{align}
	\shv F(U) \isoto \varprojlim \left( \begin{tikzcd}[column sep=small,ampersand replacement=\&]
		\prod_i \shv F(U_i) \arrow[r,shift left] \arrow[r,shift right] \&
		\prod_{i,j} \shv F(U_i \times_U U_j)
	\end{tikzcd}
	\right).
\end{align}
One checks easily that this bijection is exactly the sheafiness requirement in example (i) above. Similarly, if $\cat V = \Cat_\ordinary$ then the limit on the right depends only on the first three steps and produces exactly the category $\Desc(\cat U, \shv F)$ of descent data discussed in example (ii). In general we need to take the whole limit over $\bbDelta$, which we will call \emph{descent data} of $\shv F$ along $\cat U$ (see \cref{def:descent-data} below). Note that in general the cocycle condition from example (ii) becomes an additional datum (a homotopy witnessing the commutativity of the usual triangle), which itself has to satisfy a higher cocycle condition on quadruple fiber products $U_i \times_U U_j \times_U U_k \times_U U_\ell$. The latter condition is in fact an even higher homotopy, satisfying even higher cocycle conditions, etc. Writing these cocycle data down explicitly is quite cumbersome, but fortunately the above limit over $\bbDelta$ is a very concise way of organizing these data.

With the above preparations, let us now come to formal definitions. In the following we roughly follow the standard literature (see e.g.\ \cite[\S1.3]{SAG}), but focus on specific properties of coverings that we did not find there, but which are needed in our paper. We first provide a precise notion of coverings and descent data. The notion of coverings used above (i.e.\ a family of maps $\cat U = (U_i \to U)_{i\in I}$) is too restrictive in general, so we replace it by the following more general notion of sieves:

\begin{definition} \label{def:sieves}
Let $\cat C$ be a category and $U \in \cat C$ an object.
\begin{defenum}
	\item A \emph{sieve on $U$} is a full subcategory $\cat U \subseteq \cat C_{/U}$ such that for every morphism $V' \to V$ in $\cat C_{/U}$, if $V$ belongs to $\cat U$ then so does $V'$.

	\item Let $f\colon U' \to U$ be a morphism in $\cat C$ and $\cat U$ a sieve on $U$. The \emph{pullback} $f^* \cat U$ is the sieve on $U'$ given by the full subcategory of $\cat C_{/U'}$ consisting of those maps $g\colon V \to U'$ such that $f \comp g \in \cat U$.
\end{defenum}
\end{definition}

\begin{definition} \label{def:descent-data}
Let $\cat C$ and $\cat V$ be categories, let $\cat U \subseteq \cat C_{/U}$ be a sieve on some object $U \in \cat C$, and let $\shv F\colon \cat C^\op \to \cat V$.
\begin{defenum}
	\item The \emph{descent data} of $\shv F$ along $\cat U$ are defined as (assuming it exists)
	\begin{align}
		\Desc(\cat U, \shv F) \coloneqq \varprojlim_{V \in \cat U^\op} \shv F(V).
	\end{align}

	\item We say that $\shv F$ \emph{descends along $\cat U$} if the natural map
	\begin{align}
		\shv F(U) \isoto \Desc(\cat U, \shv F)
	\end{align}
	is an isomorphism (in particular $\Desc(\cat U, \shv F)$ exists).

	\item We say that $\shv F$ \emph{descends universally along $\cat U$} if it descends along every pullback of $\cat U$.
\end{defenum}
\end{definition}

Let us first check how \cref{def:sieves,def:descent-data} relate to the previously introduced versions of descent data along covering families.

\begin{definition}
Let $\cat C$ be a category.
\begin{defenum}
 	\item Let $(U_i \to U)_{i\in I}$ be a family of maps in $\cat C$. The \emph{sieve generated by $(U_i \to U)_i$} is the sieve $\cat U \subseteq \cat C_{/U}$ consisting of all the maps $V \to U$ which factor over some $U_i$.

 	\item We say that a sieve $\cat U \subset \cat C_{/U}$ is \emph{small} if it admits a cofinal map from some small category. We say that $\cat U$ is \emph{universally small} if all pullbacks of $\cat U$ are small. We say that $\cat U$ is \emph{small generated} if it is generated by a small family of maps $(U_i \to U)_{i\in I}$.
\end{defenum} 
\end{definition}

\begin{remark} \label{rmk:universal-smallness-from-small-generation}
Every sieve is generated by the family of all objects in the sieve (viewed as maps to the common target). Smallness of the sieve ensures that descent data along the sieve can be computed by a small limit, which can sometimes be necessary for technical reasons. By \cref{rslt:descent-data-for-generated-sieve} below, if $\cat C$ admits fiber products then a small generated sieve is automatically universally small.
\end{remark}

\begin{definition} \label{def:bbDelta-I}
For a set $I$ we denote by $\bbDelta_I \to \bbDelta$ the right fibration associated with the functor $\bbDelta^\op \to \Ani$, $[n] \mapsto I^{n+1}$. Concretely, $\bbDelta_I$ is the following ordinary category: The objects of $\bbDelta_I$ are pairs $([n] \in \bbDelta, i_\bullet \in I^{n+1})$ and a morphism $([n], i_\bullet) \to ([m], j_\bullet)$ is a morphism $\alpha\colon [n] \to [m]$ in $\bbDelta$ such that $i_k = j_{\alpha(k)}$ for all $k \in [n]$.
\end{definition}

\begin{lemma} \label{rslt:descent-data-for-generated-sieve}
Let $\cat C$ be a category and $(U_i \to U)_{i\in I}$ a family of maps in $\cat C$ generating the sieve $\cat U \subseteq \cat C_{/U}$. Assume that all fiber products of the form $U_{i_0} \times_U \dots \times_U U_{i_n}$ exist in $\cat C$. Then:
\begin{lemenum}
	\item \label{rslt:descent-data-for-generated-sieve-functor-on-Delta-I} There is a functor $\bbDelta_I^\op \to \cat U$ sending $([n], i_\bullet) \mapsto U_{i_0} \times_U \dots \times_U U_{i_n}$.

	\item Let $\cat V$ be a category which has limits over $\bbDelta_I$. Then for every functor $\shv F\colon \cat C^\op \to \cat V$ we have
	\begin{align}
		\Desc(\cat U, \shv F) = \varprojlim_{([n], i_\bullet) \in \bbDelta_I} \shv F(U_{i_0} \times_U \dots \times_U U_{i_n}).
	\end{align}
\end{lemenum}
\end{lemma}
\begin{proof}
We first prove (i). Note that we have a full embedding of categories $I \subset \bbDelta_I^\op$, where we view $I$ as a category with only identity morphisms. The family $(U_i \to U)_i$ provides a functor $U_\bullet\colon I \to \cat U$ sending $i \mapsto U_i$. In order to obtain the desired functor $\bbDelta_I^\op \to \cat U$ we perform a right Kan extension along $I \subset \bbDelta_I^\op$. For every $x = ([n], i_\bullet) \in \bbDelta_I$ the value of this right Kan extension is the limit of $U_\bullet$ over $(\bbDelta_I^\op)_{x/} \times_{\bbDelta_I^\op} I = \bigdunion_{k\in[n]} \{ i_k \}$, which is exactly the product of $U_{i_0}, \dots, U_{i_n}$ in $\cat U$.

Part (ii) follows in the same way as in the proof of \cite[Lemma~A.5.13]{Mann.2022a}; let us sketch the argument. We need to show that the functor $U_\bullet\colon \bbDelta_I^\op \to \cat U$ is cofinal, i.e.\ for every $V \in \cat U$ the category $\cat X \coloneqq \bbDelta_I^\op \times_{\cat U} \cat U_{V/}$ is weakly contractible (see \cref{def:weakly-contractible-category}). By \cref{rslt:(co)limits-in-Ani-Cat} we can equivalently show that $\varinjlim \chi = *$, where $\chi\colon \bbDelta_I^\op \to \Ani$ is the functor associated with the left fibration $\cat X \to \bbDelta_I^\op$. Note that $\chi$ is explicitly computed as $\chi([n], i_\bullet) = \Hom_{\cat U}(V, U_{i_0} \times_U \dots \times_U U_{i_n})$. Let $\chi'\colon \bbDelta^\op \to \Ani$ be the left Kan extension of $\chi$ along the projection $\bbDelta^\op_I \to \bbDelta^\op$. Thus, for every $[n] \in \bbDelta^\op$, $\chi'([n])$ is the colimit of $\chi$ over $\bbDelta^\op_I \times_{\bbDelta^\op} (\bbDelta^\op)_{/[n]}$. It follows easily that
\begin{align}
	\chi'([n]) = \bigdunion_{i_\bullet \in I^{n+1}} \chi([n], i_\bullet) = \bigdunion_{i_\bullet \in I^{n+1}} \prod_{k=0}^n \Hom_{\cat U}(V, U_{i_k}).
\end{align}
Since $\chi'$ is obtained from $\chi$ via left Kan extension, we have $\varinjlim \chi = \varinjlim \chi'$, so it is enough to show that $\varinjlim \chi' = *$. On the other hand, one sees from the above formula that $\chi'$ is the Čech nerve of the map $\bigdunion_i \Hom(V, U_i) \to *$ in $\Ani$. Since $\Ani$ is a topos, the colimit over that Čech nerve is $*$ as soon as $\bigdunion_i \Hom(V, U_i)$ is non-empty, which is clearly the case.
\end{proof}

If a sieve is generated by \emph{injective} maps $U_i \injto U$ then the fiber products $U_{i_0} \times_U \dots \times_U U_{i_n}$ depend only on the \emph{set} of indices $\{ i_0, \dots, i_n \}$, i.e.\ multiple occurrences of the same index can be ignored. This is reflected by the following result, showing that descent data along an injective sieve can be computed as limit over the finite subsets of $I$:

\begin{lemma} \label{rslt:descent-data-for-generated-sieve-by-monomorphisms}
Let $\cat C$ be a category, $(U_i \to U)_{i\in I}$ a family of monomorphisms in $\cat C$ generating the sieve $\cat U \subseteq \cat C_{/U}$. Assume that all fiber products of the form $U_{i_0} \times_U \dots \times_U U_{i_n}$ exist in $\cat C$. Then:
\begin{lemenum}
	\item Let $P_I$ be the category of finite non-empty subsets of $I$ (where morphisms are given by inclusion). Then there is a functor $P_I^\op \to \cat U$ sending $J \mapsto U_J \coloneqq \prod_{j\in J}^{/U} U_j$.

	\item Let $\cat V$ be a category which has limits over $P_I$. Then for every functor $\shv F\colon \cat C^\op \to \cat V$ we have
	\begin{align}
		\Desc(\cat U, \shv F) = \varprojlim_{J \in P_I} \shv F(U_J).
	\end{align}
\end{lemenum}
\end{lemma}
\begin{proof}
Part (i) follows from right Kan extension along the inclusion $I \subseteq P_I$ as in the proof of \cref{rslt:descent-data-for-generated-sieve-functor-on-Delta-I}. We now prove (ii), for which we employ the same strategy as in the proof of \cref{rslt:descent-data-for-generated-sieve}. Namely, given $V \in \cat U$ we need to see that the category $\cat X \coloneqq P_I^\op \times_{\cat U} \cat U_{V/}$ is weakly contractible, for which it is enough to show that $\varinjlim \chi = *$, where $\chi\colon P_I^\op \to \Ani$ is the functor $J \mapsto \Hom_{\cat U}(V, U_J)$. Now since each $U_J \injto U$ is a monomorphism and $\Hom_{\cat U}(V, U) = *$, we have $\chi(J) \in \{ *, \emptyset \}$ for all $J$. Let $I' \subset I$ be the subset of those $i$ such that $\chi(\{ i \}) = *$. Then $\chi$ is the left Kan extension of its restriction to $P_{I'}^\op \subseteq P_I^\op$ and in particular the colimit coincides with the one over $P_{I'}$. Replacing $I$ by $I'$ we can now assume that $\chi(J) = *$ for all $J \in P_I^\op$, i.e.\ $\chi\colon P_I^\op \to \Ani$ is the constant functor with value $*$. Since $P_I$ is filtered, we can combine \cref{example:(co)limits-in-Ani-Cat} and \cite[\href{https://kerodon.net/tag/02PH}{Tag 02PH}]{kerodon} to deduce $\varinjlim\chi = \lvert P_I^\op\rvert = *$.
\end{proof}

In general there is no easy description of the descent data of a functor $\shv F$ along a sieve $\cat U$ on $U$. However, a good technique for proving that $\shv F$ descends along $\cat U$ is to decompose $\cat U$ into sieves of simpler shape, i.e.\ ones that are generated by a nice family of maps (in which case the descent of $\shv F$ can be checked via \cref{rslt:descent-data-for-generated-sieve} and standard descent results). This proof technique is enabled by the following results.

\begin{lemma} \label{rslt:descent-along-pullback-of-covers}
Let $\cat C$ and $\cat V$ be categories and $\shv F\colon \cat C^\op \to \cat V$ a functor. Let $\cat U, \cat U' \subseteq \cat C_{/U}$ be sieves on some object $U \in \cat C$ and assume that the following conditions are satisfied:
\begin{enumerate}[(a)]
	\item $\shv F$ descends along $\cat U'$ and along $f^* \cat U'$ for every $f\colon U' \to U$ in $\cat U$.

	\item $\shv F$ descends along $f^* \cat U$ for every $f\colon U' \to U$ in $\cat U'$.
\end{enumerate}
Then $\shv F$ descends along $\cat U$.
\end{lemma}
\begin{proof}
Let $\cat U'' \coloneqq \cat U \isect \cat U' \subset \cat C_{/U}$. We claim that $\shv F$ descends along $\cat U''$. In order to prove this claim, we first show that $\shv F|_{\cat U'^\op}$ is the right Kan extension of $\shv F|_{\cat U''^\op}$ along the inclusion $\cat U''^\op \injto \cat U'^\op$. Namely, given any $f\colon V \to U$ in $\cat U'$ we need to see that
\begin{align}
	\shv F(V) = \varprojlim_{W \in (\cat U''^\op)_{V/}} \shv F(W),
\end{align}
where we abbreviate $(\cat U''^\op)_{V/} = \cat U''^\op \times_{\cat U'^\op} (\cat U'^\op)_{V/}$. But note that this category identifies with $(f^* \cat U)^\op$, so the desired limit property follows from descent of $\shv F$ along $f^* \cat U$. From the Kan extension property we deduce that the limits of $\shv F$ over $\cat U''^\op$ and over $\cat U'^\op$ agree. Therefore the claimed descent of $\shv F$ along $\cat U''$ follows from the descent of $\shv F$ along $\cat U'$.

In order to prove that $\shv F$ descends along $\cat U$, we now claim that $\shv F|_{\cat U^\op}$ is the right Kan extension of $\shv F|_{\cat U''^\op}$ along the inclusion $\cat U''^\op \injto \cat U^\op$. By the same argument as above (swapping $\cat U$ and $\cat U'$) this boils down to the descent of $\shv F$ along $f^* \cat U'$ for every $f \in \cat U$, which is true by assumption. Finally, by the just proved Kan extension property, the limits of $\shv F$ over $\cat U^\op$ and over $\cat U''^\op$ agree, so the desired descent of $\shv F$ along $\cat U$ follows from the descent along $\cat U''$.
\end{proof}

\begin{corollary} \label{rslt:sieve-descent-from-other-sieves}
Let $\cat C$ and $\cat V$ be categories, $\shv F\colon \cat C^\op \to \cat V$ a functor and $\cat U \subseteq \cat C_{/U}$ a sieve on some object $U \in \cat C$.
\begin{corenum}
	\item \label{rslt:descent-along-sieve-from-subsieve} Assume that there is a subsieve $\cat U' \subseteq \cat U$ on $U$ such that $\shv F$ descends universally along $\cat U'$. Then $\shv F$ descends universally along $\cat U$.

	\item \label{rslt:descent-along-composition-of-covers} Assume that there is a sieve $\cat U' \supseteq \cat U$ on $U$ such that $\shv F$ descends along $\cat U'$ and descends universally along $f^* \cat U$ for all $f \in \cat U'$. Then $\shv F$ descends along $\cat U$.
\end{corenum}
\end{corollary}
\begin{proof}
Both claims follow by applying \cref{rslt:descent-along-pullback-of-covers} to $\cat U$ and $\cat U'$. In (i), condition (a) follows from universal descent of $\shv F$ along $\cat U'$ and condition (b) follows from the observation that $f^* \cat U = f^* \cat U'$ for $f \in \cat U'$. In (ii), condition (a) follows from the observation that $f^* \cat U' = f^* \cat U$ for $f \in \cat U$ and condition (b) is obvious.
\end{proof}

\begin{example}
In the setting of \cref{rslt:descent-along-sieve-from-subsieve}, assume that $\cat C$ has fiber products and assume that $\cat U'$ is generated by a family $(U_i \to U)_{i\in I}$. Then in order to prove descent of $\cat F$ along $\cat U$, it is enough to show that $\cat F$ descends along the family $(U_i \to U)_i$ and along any pullback of that family. This follows from \cref{rslt:descent-data-for-generated-sieve} together with the observation that for any map $f\colon U' \to U$ in $\cat C$, the sieve $f^* \cat U'$ is generated by the family $(U'_i \to U')_i$, where $U'_i = U_i \times_U U'$.
\end{example}

\begin{remark}
In the case that $\cat C$ has fiber products and all sieves in \cref{rslt:descent-along-pullback-of-covers,rslt:sieve-descent-from-other-sieves} are generated by a single map, we arrive at the results in \cite[Lemma~3.1.2]{Liu-Zheng.2012}.
\end{remark}

Having a good understanding of sieves and descent along them, we now come to the definition of sites and sheaves. With the above preparations, this is easy to define:

\begin{definition} \label{def:site}
A \emph{site} is a small category $\cat C$ equipped with a distinguished collection of sieves, called \emph{coverings}, satisfying the following conditions:
\begin{enumerate}[(i)]
	\item For every $U \in \cat C$, the sieve $\cat C_{/U}$ is a covering.
	\item Coverings are stable under pullback.
	\item Let $U \in \cat C$ and let $\cat U, \cat U' \subseteq \cat C_{/U}$ be two sieves on $U$. Suppose that $\cat U$ is a covering and that for every $f\colon V \to U$ in $\cat U$, the sieve $f^* \cat U'$ is a covering. Then $\cat U'$ is a covering.
\end{enumerate}
\end{definition}

\begin{definition} \label{def:sheaf}
Let $\cat C$ be a site and $\cat V$ a category. A \emph{$\cat V$-valued sheaf on $\cat C$} is a functor $\shv F\colon \cat C^\op \to \cat V$ which descends along every covering sieve of $\cat C$. We denote by
\begin{align}
	\Shv(\cat C, \cat V) \subseteq \Fun(\cat C^\op, \cat V)
\end{align}
the full subcategory of $\cat V$-valued sheaves. We also abbreviate $\Shv(\cat C) \coloneqq \Shv(\cat C, \Ani)$.
\end{definition}

\begin{example}
Suppose $\cat C$ is an ordinary site, e.g.\ as defined in  \cite[\href{https://stacks.math.columbia.edu/tag/03NH}{Definition~03NH}]{stacks-project}. Then one can equip $\cat C$ with the structure of a site as in \cref{def:site} by designating precisely those sieves $\cat U \subset \cat C_{/U}$ as coverings for which there is a covering family $(U_i \to U)_{i\in I}$ with all $U_i \in \cat U$. With this definition, \cref{rslt:descent-data-for-generated-sieve,rslt:descent-along-sieve-from-subsieve} imply that a $\Set$-valued or $\Cat_\ordinary$-valued sheaf on $\cat C$ is the same as was discussed at the beginning of this subsection.
\end{example}

One of the most important constructions in sheaf theory is the sheafification functor. It exists in the abstract setup above, at least under mild assumptions on $\cat V$:

\begin{lemma} \label{rslt:sheafification-exists}
Let $\cat C$ be a site and let $\cat V$ be a presentable category. Then $\Shv(\cat C, \cat V)$ is presentable and the inclusion $\Shv(\cat C, \cat V) \subseteq \Fun(\cat C^\op, \cat V)$ has a left adjoint
\begin{align}
	(\blank)^\sharp\colon \Fun(\cat C^\op, \cat V) \to \Shv(\cat C, \cat V), \qquad \shv F \mapsto \shv F^\sharp,
\end{align}
called the \emph{sheafification}.
\end{lemma}
\begin{proof}
Note that $\Fun(\cat C^\op, \cat V)$ is presentable because $\cat C^\op$ is small (see \cite[Proposition~5.5.3.6]{HTT}). The claim that $\Shv(\cat C, \cat V)$ is presentable and that $(\blank)^\sharp$ exists is exactly saying that $\Shv(\cat C, \cat V)$ is a strongly reflective subcategory of $\Fun(\cat C^\op, \cat V)$ (see the paragraph before \cite[Lemma~5.5.4.17]{HTT}). That this is indeed satisfied follows from the fact that $\Shv(\cat C, \cat V)$ is cut out by a small amount of limit conditions, so that we can apply \cite[Lemmas~5.5.4.17,~5.5.4.18,~5.5.4.19]{HTT}. See the proof of \cite[Lemma~A.3.4]{Mann.2022a} for a very similar argument.
\end{proof}

\begin{definition}
Let $\cat C$ be a site. By composing the Yoneda embedding with the sheafification functor from \cref{rslt:sheafification-exists} we obtain the functor
\begin{align}
	\cat C \to \Shv(\cat C), \qquad U \mapsto [V \mapsto \Hom(V, U)]^\sharp.
\end{align}
\begin{defenum}
	\item A sheaf $\shv F \in \Shv(\cat C)$ is called \emph{representable} if it lies in the essential image of the functor $\cat C \to \Shv(\cat C)$.

	\item The site $\cat C$ is called \emph{subcanonical} if the functor $\cat C \to \Shv(\cat C)$ is fully faithful. Equivalently, for every $U \in \cat C$ the presheaf $\Hom(\blank, U)$ is already a sheaf.
\end{defenum}
\end{definition}

The category $\cat X = \Shv(\cat C)$ of $\Ani$-valued sheaves on a site $\cat C$ is of particularly nice shape. It is an example of a so-called \emph{topos}, meaning that it is presentable, colimits are universal, coproducts are disjoint and every groupoid object is effective. We refer the reader to \cite[\S6.1]{HTT} for a precise definition of topoi (which are called \enquote{$\infty$-topoi} in the reference). For every topos $\cat X$ there is an associated notion of $\cat V$-valued sheaves on $\cat X$:

\begin{definition}
Let $\cat X$ be a topos and $\cat V$ a category which has all small limits. A \emph{$\cat V$-valued sheaf on $\cat X$} is a functor $\shv F\colon \cat X^\op \to \cat V$ which preserves all small limits. We denote by
\begin{align}
	\Shv(\cat X, \cat V) \subset \Fun(\cat X^\op, \cat V)
\end{align}
the full subcategory of $\cat V$-valued sheaves.
\end{definition}

We now have two notions of sheaves: one associated with sites and one associated with topoi. Fortunately, there is little room for ambiguity, as the following result shows.

\begin{lemma} \label{rslt:sheaves-on-topos-equiv-sheaves-on-site}
Let $\cat C$ be a site with associated topos $\cat X = \Shv(\cat C)$ and let $\cat V$ be a category that has all small limits. Then precomposition with the functor $\cat C \to \cat X$ induces an equivalence of categories
\begin{align}
	\Shv(\cat X, \cat V) \isoto \Shv(\cat C, \cat V).
\end{align}
\end{lemma}
\begin{proof}
See \cite[Proposition~1.3.1.7]{SAG}. The idea is that every sheaf $\shv F \in \cat X$ can be written as a colimit of representable sheaves (as this is already true for presheaves), hence a $\cat V$-valued sheaf on $\cat X$ is determined by its restriction to $\cat C$.
\end{proof}

The above notion of sheaves requires descent along sieves, which by \cref{rslt:descent-data-for-generated-sieve} amounts to Čech descent. In the context of higher categories there is a stronger descent property one can ask for, namely descent along all \emph{hypercovers}. Hyperdescent is automatic for truncated sheaves (see \cite[Lemma~6.5.2.9]{HTT}), but is an extra condition in the unbounded case. Unfortunately, on a general site the notion of hyperdescent is hard to grasp, so we are forced to make the following unsatisfying definition:

\begin{definition}
Let $\cat V$ be a category which has all small limits.
\begin{defenum}
	\item Let $\cat X$ be a topos. A $\cat V$-valued sheaf $\shv F\colon \cat X^\op \to \cat V$ is called \emph{hypercomplete} if for all hypercovers $U_\bullet \to X$ in $\cat X$ (see \cite[Definition~6.5.3.2]{HTT}) the associated map 
	\begin{align}
		\shv F(X) \isoto \varprojlim_{[n]\in\bbDelta} \shv F(U_n)
	\end{align}
	is an isomorphism in $\cat V$. We denote by
	\begin{align}
		\HypShv(\cat X, \cat V) \subseteq \Shv(\cat X, \cat V)
	\end{align}
	the full subcategory spanned by the hypercomplete $\cat V$-valued sheaves on $\cat X$. In the special case $\cat V = \Ani$ we write $\Hyp(\cat X) \coloneqq \HypShv(\cat X, \Ani) \subseteq \cat X$ and call it the \emph{hypercompletion} of $\cat X$. Note that $\Hyp(\cat X)$ is again a topos (see the remarks before \cite[Lemma~6.5.2.9]{HTT}).

	\item Let $\cat C$ be a site. A $\cat V$-valued sheaf $\shv F\colon \cat C^\op \to \cat V$ on $\cat C$ is called \emph{hypercomplete} if it is hypercomplete as a sheaf on the associated topos via \cref{rslt:sheaves-on-topos-equiv-sheaves-on-site}. We denote
	\begin{align}
		\HypShv(\cat C, \cat V) \subseteq \Shv(\cat C, \cat V)
	\end{align}
	the full subcategory of hypercomplete sheaves. In case $\cat V = \Ani$ we simply write $\HypShv(\cat C)$.

	\item We say that a site $\cat C$ is \emph{hyper-subcanonical} if $\cat C$ is subcanonical and all representable sheaves are hypercomplete.
\end{defenum}
\end{definition}

\begin{remark}
In general it is hard to verify hypercompleteness of a sheaf on a site, as it is a priori a requirement on hypercovers in the associated topos of $\Ani$-valued sheaves. However, under some additional assumptions on $\cat C$, it is enough to check descent along hypercovers in $\cat C$. For this to work one roughly needs to assume that fiber products and finite coproducts exist in $\cat C$ and interact nicely, and that every covering in $\cat C$ can be refined to a covering generated by a single map. We refer to \cite[Definition~A.3.12]{Mann.2022a} for the precise list of conditions and to \cite[Proposition~A.3.16]{Mann.2022a} for the promised characterization of hypersheaves. These results are slight variations of Lurie, see \cite[Proposition~A.5.7.2]{SAG}.
\end{remark}

\begin{remark} \label{rmk:subcanonical-ordinary-site-is-hyper-subcanonical}
Every subcanonical ordinary site is automatically hyper-subcanonical because all representable sheaves are then $0$-truncated (they are sheaves of sets), so we can apply \cite[Lemma~6.5.2.9]{HTT}.
\end{remark}

It follows easily from the definitions that we have the following version of \cref{rslt:sheaves-on-topos-equiv-sheaves-on-site} for hypercomplete sheaves. For the following result, note that by similar arguments as in \cref{rslt:sheafification-exists} there is a hypercompletion functor $\Shv(\cat C, \cat V) \to \HypShv(\cat C, \cat V)$, so in particular there is a natural functor $\cat C \to \HypShv(\cat C)$.

\begin{lemma} \label{rslt:sheaves-on-hypertopos-equiv-hypersheaves-on-site}
Let $\cat C$ be a site with associated hypercomplete topos $\cat X = \HypShv(\cat C)$ and let $\cat V$ be a category that has all small limits. Then precomposition with the functor $\cat C \to \cat X$ induces an equivalence of categories
\begin{align}
	\Shv(\cat X, \cat V) \isoto \HypShv(\cat C, \cat V).
\end{align}
\end{lemma}
\begin{proof}
This follows from the definitions, see also \cite[Lemma~A.3.6]{Mann.2022a}.
\end{proof}

We finish this section with the following result relating descent of abelian categories and their associated derived categories. A similar result appeared in \cite[Proposition~A.1.2]{Mann.2022a}, but the following result is slightly more general and removes some unnecessary hypotheses.

\begin{proposition}\label{rslt:Grothendieck-descent}
Let $\cat A^\bullet$ be a cosimplicial object of $\Cat$ with limit $\cat A \coloneqq \varprojlim_{n\in\bbDelta} \cat A^n$ in $\Cat$. Assume that the following conditions hold:
\begin{enumerate}[(a)]
\item All $\cat A^n$ are Grothendieck abelian categories. 

\item For every $\alpha \colon [n]\to [m]$ in $\bbDelta$, the transition functor $\alpha^* \colon \cat A^n \to \cat A^m$ is exact and admits a right adjoint $\alpha_*$.

\item Letting $d^i\colon [n]\to [n+1]$ denote the $i$-th coface map, then for all $\alpha\colon [n] \to [m]$ in $\bbDelta$ the following diagram satisfies base-change:
\begin{equation}
\begin{tikzcd}
\D^+(\cat A^{n+1}) \ar[d,"\alpha'^*"'] \ar[r,"\R d^0_*"] & \D^+(\cat A^n) \ar[d,"\alpha^*"] \\
\D^+(\cat A^{m+1}) \ar[r,"\R d^0_*"'] & \D^+(\cat A^m).
\end{tikzcd}
\end{equation}
In other words, the natural morphism of functors $\alpha^*\comp \R d^0_* \isoto \R d^0_* \comp \alpha'^*$ is an equivalence; here, $\alpha'$ is defined by $\id_{[0]}\star \alpha$ under the identifications $[n+1] = [0]\star [n]$ and $[m+1] = [0]\star [m]$.
\end{enumerate}
Then:
\begin{enumerate}[(i)]
\item The category $\cat A$ is Grothendieck abelian and admits the following description: An object of $\cat A$ is a pair $(X,\alpha)$, where $X$ is an object of $\cat A^0$ and $\alpha\colon d^{0*}X \isoto d^{1*}X$ is an isomorphism in $\cat A^1$ such that the following diagram in $\cat A^2$ commutes:
\begin{equation}
\begin{tikzcd}[row sep=tiny]
& d^{0*}d^{1*} X \ar[r,equals] & d^{2*}d^{0*}X \ar[dr,"d^{2*}\alpha"] \\
d^{0*}d^{0*}X \ar[ur, "d^{0*}\alpha"] \ar[dr,equals] & & & d^{2*}d^{1*}X \\
& d^{1*}d^{0*}X \ar[r,"d^{1*}\alpha"'] & d^{1*}d^{1*}X \ar[ur,equals].
\end{tikzcd}
\end{equation}
A morphism $(X_1,\alpha_1) \to (X_2,\alpha_2)$ in $\cat A$ is a morphism $f\colon X_1\to X_2$ in $\cat A^0$ such that $d^{1*}f \comp \alpha_1 = \alpha_2 \comp d^{0*}f$.

\item Suppose that $d^{1*}\colon \cat A^0\to \cat A^1$ sends injective objects to $d^0_*$-acyclic objects. Then there is a natural equivalence of categories
\begin{align}
\D^+(\cat A) \isoto \varprojlim_{n\in\bbDelta} \D^+(\cat A^n),
\end{align}
where the limit on the right is computed in $\Cat$. Moreover, if the categories $\D(\cat A)$ and $\D(\cat A^n)$ are left-complete, then the above equivalence extends to an equivalence
\begin{align}
\D(\cat A) \isoto \varprojlim_{n\in\bbDelta} \D(\cat A^n).
\end{align}
\end{enumerate}
\end{proposition}
\begin{proof}
Note first that $\Cat_{\ordinary}$ is equivalent to a $(2,1)$-category by \cite[Proposition~2.3.4.18]{HTT}, hence it follows from the proof of \cite[Lemma~1.3.3.10]{HA} that we have an equivalence $\cat A = \varprojlim_{n\in\bbDelta^{\le 2}} \cat A^n$. This implies the explicit description of $\cat A$. Given this description, it is clear that $\cat A$ is an abelian category and the forgetful functor $\cat A \to \cat A^0$ is conservative and exact. Since the transition maps are exact and commute with colimits, it follows that $\cat A$ satisfies the Grothendieck axiom AB5 (colimits exist and filtered colimits are exact). As Grothendieck abelian categories are precisely the presentable abelian AB5 categories (see \cite[Proposition~3.10]{Beke.2000}) and the limit of presentable categories is presentable (see \cite[Proposition~5.5.3.13]{HTT}), it follows that $\cat A$ is Grothendieck. This proves (i).

We now prove (ii). By \cite[Lemma~A.1.1]{Mann.2022a} we obtain an augmented cosimplicial object $\D^+(\cat A^\bullet)\colon \bbDelta_+ = \bbDelta^{\triangleleft} \to \Cat$, where we write $\cat A^{-1} \coloneqq \cat A$. 
Applying \cite[Theorem~4.7.5.2]{HA} to $(\cat A^\bullet)^\op$, we see that the forgetful functor $d^{0*}\colon \cat A\to \cat A^0$ admits a right adjoint $d^0_*$ and that the natural map $d^{0*} d^0_* \isoto d^0_* d^{1*}$ is an isomorphism of functors $\cat A^0\to \cat A^0$. Since $d^{1*}\colon \cat A^0\to \cat A^1$ sends injective objects to $d^0_*$-acyclic objects, the natural map $d^{0*}\comp \R d^0_*\isoto \R d^0_*\comp d^{1*}$ is an isomorphism of functors $\D^+(\cat A^0)\to \D^+(\cat A^0)$. 
In order to prove the claim, we want to apply (the dual of) Lurie's Beck--Chevalley descent result (see \cite[Corollary~4.7.5.3]{HA}), so we need to verify the following conditions:
\begin{enumerate}[(1)]
\item The forgetful functor $F\coloneqq d^{0*}\colon \D^+(\cat A) \to \D^+(\cat A^0)$ is conservative.

\item The category $\D^+(\cat A)$ admits totalizations of $F$-split cosimplicial objects, and these are preserved by $F$.

\item For every morphism $\alpha\colon [n]\to [m]$ in $\bbDelta_+$, the natural map $\alpha^*\comp \R d^0_* \isoto \R d^0_*\comp \alpha'^*$ is an isomorphism of functors $\D^+(\cat A^{n+1}) \to \D^+(\cat A^m)$.
\end{enumerate}
The discussion above shows that (3) is satisfied: For $n\ge0$ this is condition (c) and for $n=-1$ we use that the map $[-1] \to [m]$ factors through $[-1] \to [0]$. Condition (1) can be checked on cohomology groups, which reduces the claim to $d^{0*}\colon \cat A\to \cat A^0$ being conservative; but this is clear from the explicit description of $\cat A$. It therefore remains to prove (2). Let $X^\bullet$ be an $F$-split cosimplicial object in $\D^+(\cat A)$. In view of \cite[\href{https://kerodon.net/tag/04RE}{Tag 04RE}]{kerodon} we may instead work with the cosemisimplicial object underlying $X^\bullet$. Assuming for the moment that the totalization $\Tot(X^\bullet)$ exists, we need to verify that the natural map $F\Tot(X^\bullet) \to \Tot(FX^\bullet)$ is an isomorphism in $\D^+(\cat A^0)$, which can be checked after applying the cohomology functors $\H^n$. We will apply \cite[Corollary~1.2.4.12]{HA}, remarking that the conclusion of the result holds true (with essentially the same proof) if the condition \enquote{right complete} in \loccit{} is replaced with \enquote{right bounded}:\footnote{We point out the following typo in \cite[Corollary~1.2.4.12]{HA}: The conditions in the statement need to be verified for all integers $q$, not just for $q\ge0$.} Indeed, the right completeness is invoked twice to show that certain colimits exists---but if $\cat C$ is right bounded, then these colimit diagrams are eventually constant, so the proof also works in this case.

Our argument follows the reasoning in \cite[Proposition~D.6.4.6]{SAG}. For any $n\in\ZZ$ we consider the unnormalized chain complex
\begin{align}\label{eq:Grothendieck-descent}
0\to K(n) \to \H^nX^0 \xrightarrow{\partial^1} \H^nX^1 \to \H^nX^2 \to \dotsb
\end{align}
associated with $\H^nX^\bullet$, where we put $K(n)\coloneqq \Ker(\partial^1)$. Since  $F$ is exact and $FX^\bullet$ is split by assumption, it follows that applying $F$ to \eqref{eq:Grothendieck-descent} yields a split acyclic complex. (The exactness is used to identify $F(K(n))$ with the kernel of $\H^nFX^0\to \H^nFX^1$.) As $F$ is conservative, we deduce that \eqref{eq:Grothendieck-descent} is acyclic, for all $n\in\ZZ$. By \cite[Corollary~1.2.4.12]{HA} applied to $\cat C = \D^+(\cat A)^\op$, it follows that the totalization $X \coloneqq \Tot(X^\bullet)$ exists and that the map $X\to X^0$ induces an isomorphism $\H^nX \isoto K(n)$, for all $n$. Applying \loccit{} to $\cat C = \D^+(\cat A^0)^\op$, we deduce that also $FX^\bullet$ admits a totalization and the map $\H^n\Tot(FX^\bullet) \isoto F(K(n))$ is an isomorphism for all $n$. This shows that the natural map $\H^nFX \isoto \H^n\Tot(FX^\bullet)$ is an isomorphism for all $n$, as desired. Consequently, the natural map $\D^+(\cat A) \isoto \varprojlim \D^+(\cat A^\bullet)$ is an isomorphism. The last assertion follows from this by passing to the left-completions (and using \cite[Remark~1.2.1.18]{HA}).
\end{proof}

%

\subsection{Presentable categories} \label{sec:cat.pres}

A very important class of categories is formed by the presentable categories, which are large categories that can be generated by a small set of generators. As we will be using presentable categories in many places in this paper, we briefly recall the main definitions and constructions, see \cite[\S5]{HTT} for details.

For the following definition we use the notion of $\kappa$-filtered colimits, where $\kappa$ is a regular cardinal. We refer the reader to \cite[\S5.3.1]{HTT} for the basic definitions of filtered categories. Most importantly, \cite[Proposition~5.3.1.16]{HTT} shows that filtered colimits can always be computed via colimits over directed sets (in the classical sense).

\begin{definition}
Let $\cat C$ be a category and $\kappa$ a regular cardinal.
\begin{defenum}
	\item An object $X \in \cat C$ is called \emph{$\kappa$-compact} if the functor $\Hom(X, \blank)\colon \cat C \to \Ani$ preserves $\kappa$-filtered colimits. We denote by $\cat C^\kappa \subseteq \cat C$ the full subcategory spanned by the $\kappa$-compact objects.

	\item \label{def:kappa-compactly-generated-category} $\cat C$ is called \emph{$\kappa$-compactly generated} if $\cat C$ admits all small colimits and there is a (small) set $S$ of $\kappa$-compact objects in $\cat C$ such that every object of $\cat C$ is a (small) colimit of objects in $S$. We denote by $\PrL[\kappa]$ the category of $\kappa$-compactly generated categories and colimit preserving functors which preserve $\kappa$-compact objects.

	\item $\cat C$ is called \emph{presentable} if it is $\kappa$-compactly generated for some $\kappa$. We denote by $\PrL$ the category of presentable categories and colimit preserving functors.
\end{defenum}
\end{definition}

The most important case of \cref{def:kappa-compactly-generated-category} is $\kappa = \omega$, in which case we simply speak of \emph{compactly generated categories}. A $\kappa$-compactly generated category is essentially given by its full subcategory of $\kappa$-compact objects. This will be shown using the following general construction which freely attaches $\kappa$-filtered colimits to a small category:

\begin{definition}[{cf. \cite[Corollary~5.3.5.4]{HTT}}]
Let $\cat C$ be a small category and $\kappa$ a regular cardinal. We denote by
\begin{align}
	\Ind_\kappa(\cat C) \subseteq \PSh(\cat C) = \Fun(\cat C^\op, \Ani)
\end{align}
the full subcategory generated under $\kappa$-filtered colimits by the image of the Yoneda embedding $\cat C \injto \PSh(\cat C)$. In the case $\kappa = \omega$ we simply write $\Ind(\cat C)$. We also denote
\begin{align}
	\Pro_\kappa(\cat C) \coloneqq \Ind_\kappa(\cat C^\op)^\op
\end{align}
and $\Pro(\cat C) = \Pro_\omega(\cat C)$.
\end{definition}

\begin{lemma} \label{rslt:universal-property-of-Ind-kappa}
Let $\cat C$ be a small category, $\kappa$ a regular cardinal and $\cat D$ a category which admits $\kappa$-filtered colimits. Then left Kan extension along the inclusion $\cat C \injto \Ind_\kappa(\cat C)$ induces a fully faithful embedding
\begin{align}
	\Fun(\cat C, \cat D) \injto \Fun(\Ind_\kappa(\cat C), \cat D),
\end{align}
whose essential image are precisely those functors $\Ind_\kappa(\cat C) \to \cat D$ that preserve all $\kappa$-filtered colimits.
\end{lemma}
\begin{proof}
This is found in \cite[Proposition~5.3.5.10]{HTT}.
\end{proof}

It follows easily from the definitions that every object in $\cat C \subset \Ind_\kappa(\cat C)$ is $\kappa$-compact, hence $\Ind_\kappa(\cat C)$ is $\kappa$-compactly generated. In general the full subcategory $\Ind_\kappa(\cat C)^\kappa \subseteq \Ind_\kappa(\cat C)$ might be bigger than $\cat C$. Namely, it is the full subcategory of $\Ind_\kappa(\cat C)$ generated under $\kappa$-small colimits and retracts by $\cat C$. The next result shows that every $\kappa$-compactly generated category is obtained via $\Ind_\kappa$:

\begin{proposition} \label{rslt:characterizations-of-kappa-compactly-generated-categories}
For a category $\cat C$ and a regular cardinal $\kappa$, the following are equivalent:
\begin{propenum}
	\item $\cat C$ is $\kappa$-compactly generated.
	\item $\cat C \isom \Ind_\kappa(\cat C_0)$ for some small category $\cat C_0$ which has all $\kappa$-small colimits.
\end{propenum}
\end{proposition}
\begin{proof}
This is part of \cite[Theorem~5.5.1.1]{HTT}. We have already observed that (ii) implies (i). For the other direction, the key observation is that if $\cat C$ is $\kappa$-compactly generated then $\cat C = \Ind_\kappa(\cat C^\kappa)$.
\end{proof}

\begin{remark} \label{rslt:PrL-kappa-equiv-Cat-kappa}
By combining the above results one can show that $\PrL[\kappa]$ is equivalent to the category $\Cat_\kappa$ of small categories which have all $\kappa$-small colimits and retracts, and functors preserving $\kappa$-small colimits. The equivalence is given by
\begin{align}
	\Ind_\kappa\colon \Cat_\kappa \rightleftarrows \PrL[\kappa] \noloc (\blank)^\kappa.
\end{align}
See \cite[Proposition~5.5.7.10]{HTT} for details.
\end{remark}

Before we continue with the theory, let us provide some examples. Intuitively, every big \enquote{reasonable} category is presentable, and in practice even $\omega_1$-compactly generated.

\begin{examples}
\begin{exampleenum}
	\item $\Ani$ and $\Cat$ are compactly generated (for $\Cat$ one can e.g.\ use the description in \cref{rslt:Segal-anima-categories}, see also \cite[Proposition~A.2.10]{Mann.2022a}).

	\item If $\Lambda$ is an ordinary ring then the category of ordinary $\Lambda$-modules is compactly generated. Also the derived category $\D(\Lambda)$ is compactly generated (see \cite[\S1.3.5]{HA} for the construction of $\D(\Lambda)$ in the $\infty$-setting and see \cite[Proposition~1.3.5.21]{HA} for the claim about presentability).
\end{exampleenum}
\end{examples}

A very powerful tool in category theory is the adjoint functor theorem for presentable categories. It allows one to construct the adjoint functor to a given functor \enquote{out of nowhere}. In this paper this is often used to construct $f_*$ and $f^!$ as the right adjoints of $f^*$ and $f_!$, respectively.

\begin{theorem} \label{rslt:adjoint-functor-theorem}
Let $F\colon \cat C \to \cat D$ be a functor and assume that $\cat C$ is presentable.
\begin{thmenum}
	\item $F$ is left adjoint if and only if it preserves all small colimits (and in particular $\cat D$ has the relevant colimits).

	\item Suppose that $\cat D$ is presentable. Then $F$ is right adjoint if and only if it preserves all small limits and all $\kappa$-filtered colimits for some regular cardinal $\kappa$.
\end{thmenum}
\end{theorem}
\begin{proof}
See \cite[Corollary~5.5.2.9]{HTT} and \cite[Remark~5.5.2.10]{HTT}.
\end{proof}

The category $\PrL$ has all small limits and colimits und Lurie provides an explicit description of them. Namely, they can be computed as follows:
\begin{itemize}
	\item Limits in $\PrL$ commute with the forgetful functor $\PrL \to \Cat$, i.e.\ they can be computed in $\Cat$ (see \cite[Proposition~5.5.3.13]{HTT}).

	\item Let $\PrR$ be the category of presentable categories and right adjoint functors. Then $\PrL = (\PrR)^\op$ (see \cite[Corollary~5.5.3.4]{HTT} or \cref{rslt:passing-to-adjoints} below), hence colimits in $\PrL$ can be computed as limits in $\PrR$. Moreover, limits in $\PrR$ commute with the forgetful functor to $\Cat$, i.e.\ are computed on underlying categories.
\end{itemize}
As an example, we have the following computation of colimits:

\begin{lemma} \label{rslt:PrL-kappa-stable-under-colim-in-PrL}
For every regular cardinal $\kappa$ the forgetful functor $\PrL[\kappa] \to \PrL$ preserves all small colimits.
\end{lemma}
\begin{proof}
This is found in \cite[Lemma~2.5.7]{Krause-Nikolaus:Sheaves-on-Manifolds}, which we repeat in slightly more verbose language. Suppose we are given a colimit $\cat C = \varinjlim_{i \in I} \cat C_i$ in $\PrL$ such that the diagram $(\cat C_i)_i$ factors over $\PrL[\kappa]$. Let us first check that $\cat C$ is $\kappa$-compactly generated. By the above description of colimits in $\PrL$ we can compute $\cat C = \varprojlim_{i \in I^\op} \cat C_i$, where the limit is now taken in $\Cat$ and where the transition functors are the right adjoints of the transition functors in the original diagram. One checks easily that these right adjoints preserve $\kappa$-filtered colimits, which means that $\kappa$-filtered colimits in $\cat C$ are computed pointwise in each $\cat C_i$. In particular each projection $\cat C \to \cat C_i$ preserves $\kappa$-filtered colimits, hence their left adjoints $\cat C_i \to \cat C$ preserve $\kappa$-compact objects. Since the collection of functors $(\cat C \to \cat C_i)_i$ is conservative, the subcategory of $\cat C$ spanned by the images of $\cat C_i^\kappa \to \cat C$ forms a generating set of $\kappa$-compact objects in $\cat C$, hence $\cat C$ is $\kappa$-compactly generated (cf. \cite[Lemma~A.2.1(i)]{Mann.2022a}).

To finish the proof, it remains to show that a colimit preserving functor $\cat C \to \cat D$ to a $\kappa$-compactly generated category $\cat D$ preserves $\kappa$-compact objects if and only if each of the induced functors $\cat C_i \to \cat D$ does. Similar to the previous paragraph this translates to the claim that a right adjoint functor $\cat D \to \cat C$ preserves $\kappa$-filtered colimits if and only each composition $\cat D \to \cat C_i$ does, which follows from the fact that these $\kappa$-filtered colimits are computed pointwise.
\end{proof}

\begin{corollary} \label{rslt:Ind-preserves-filtered-colim}
Let $\kappa$ be a regular cardinal and $(\cat C_i)_i$ be a $\kappa$-filtered diagram of small categories such that all $\cat C_i$ have all $\kappa$-small colimits and all transition functors preserve $\kappa$-small colimits. Then
\begin{align}
	\varinjlim_i \Ind_\kappa(\cat C_i) = \Ind_\kappa(\varinjlim_i \cat C_i),
\end{align}
where the colimit on the left is formed in $\PrL$ and the colimit on the right is formed in $\Cat$.
\end{corollary}
\begin{proof}
By \cref{rslt:PrL-kappa-stable-under-colim-in-PrL} the colimit on the left-hand side can be computed in $\PrL[\kappa]$. By \cite[Proposition~5.5.7.10]{HTT} we deduce the claimed identity if the colimit $\varinjlim_i \cat C_i$ is computed in the category of small categories with $\kappa$-small colimits and functors preserving them. But by \cite[Proposition~5.5.7.11]{HTT}, $\kappa$-filtered colimits in this category can be computed in $\Cat$.
\end{proof}

The category $\PrL$ has the binary operation $\times$, which coincides with the one in $\Cat$. However, in \cite[Proposition~4.8.1.15]{HA} Lurie constructs a better behaved tensor product $\tensor$ on $\PrL$ such that $\cat C \tensor \cat D$ is the universal presentable category that admits a functor $\cat C \times \cat D \to \cat C \tensor \cat D$ which preserves colimits separately in each argument. We make use of this tensor product in the definition of presentable 6-functor formalisms.

\section{Algebra} 

In this paper we make heavy use of algebra in the context of higher categories. The present section collects the main definitions and results from the literature and provides some basic results which we did not find elsewhere.

\subsection{Operads} \label{sec:algebra.operads}

The main purpose of operads is to provide a solid framework to talk about algebraic structures, like (symmetric) monoidal categories and algebras and modules therein. We refer the reader to the beginning of \cite[\S2]{HA} for a good introduction to operads. We briefly recall the main definitions and constructions. A core definition is that of $\Fin_*$:

\begin{definition} \label{def:category-of-finite-pointed-sets}
We denote by $\Fin_*$ the (ordinary) category of finite pointed sets. For $n \ge 0$ we write $\langle n\rangle^\circ \coloneqq \{1,2,\dotsc,n\}$ and $\langle n \rangle \coloneqq \{*\}\dunion \langle n\rangle^\circ \in \Fin_*$. We say that a morphism $\alpha\colon \langle m \rangle \to \langle n \rangle$ is \emph{inert} if for all $i = 1, \dots, n$ the inverse image $\alpha^{-1}(i)$ has exactly one element.
\end{definition}

With the category $\Fin_*$ at hand, we can now introduce operads. To motivate the definition, let us briefly recall the definition of a symmetric monoidal category. This is a category $\cat C$ together with an operation $\tensor\colon \cat C \times \cat C \to \cat C$ which is associative and symmetric and has a unit object $\one \in \cat C$. Of course, associativity and symmetry are themselves given by certain isomorphisms and higher homotopies and thus must be added to the data defining the symmetric monoidal structure on $\cat C$. In the case that $\cat C$ is an ordinary category, one can make these data explicit, albeit somewhat cumbersome. For higher categories, one needs a more elegant way of encoding the data, which can be done as follows: A symmetric monoidal category is the same as a functor $F\colon \Fin_* \to \Cat$ such that for all $n \ge 0$ the map $F([n]) \to F([1])^n$ induced by the $n$ inert maps $\langle n \rangle \to \langle 1 \rangle$, sending all but one index to $*$, is an equivalence. Given $F$, we denote $\cat C \coloneqq F(\langle 1 \rangle)$. Then the map $\langle 2 \rangle \to \langle 1 \rangle$ sending $1 \mapsto 1$, $2 \mapsto 1$ induces a functor
\begin{align}
	\tensor\colon F(\langle 2 \rangle) = \cat C \times \cat C \to F(\langle 1 \rangle) = \cat C
\end{align}
and the unique map $\langle 0 \rangle \to \langle 1 \rangle$ induces a functor $* \to \cat C$, i.e.\ an object $\one \in \cat C$. One checks that the coherences in the functor $F$ encode the associativity and symmetric constraints on $\cat C$.

We have seen how one can define symmetric monoidal categories via functors $\Fin_* \to \Cat$. Then symmetric monoidal functors are simply natural transformations of such functors. However, in practice it is very useful to have weaker versions of these notions, e.g.\ \emph{lax} symmetric monoidal functors. In the spirit of \cref{rmk:lax-trafos-via-unstraightening} it is therefore convenient to pass to fibrations $\cat O^\tensor \to \Fin_*$ in place of functors $\Fin_* \to \Cat$, as they give us much more flexibility. This motivates the following definition:

\begin{definition}[{see \cite[Definition 2.1.1.10]{HA}}] \label{def:operad}
An \emph{operad} is a category $\cat O^\tensor$ together with a functor $p\colon \cat O^\tensor \to \Fin_*$ satisfying the following conditions:
\begin{enumerate}[(i)]
	\item For every inert morphism $\alpha\colon \langle m \rangle \to \langle n \rangle$ in $\Fin_*$ and every $X \in \cat O^\tensor$ there exists a $p$-cocartesian lift $\colon X \to X'$ of $\alpha$. We call $p$-cocartesian lifts of inert morphisms in $\Fin_*$ the \emph{inert} morphisms in $\cat O^\tensor$.

	\item Let $\cat O^\tensor_{\langle n \rangle}$ denote the fiber of $\cat O^\tensor$ over $\langle n \rangle \in \Fin_*$ and let $\cat O \coloneqq \cat O^\tensor_{\langle 1 \rangle}$. For $i = 1, \dots, n$ let $\rho_{i!}\colon \cat O^\tensor_{\langle n \rangle} \to \cat O$ denote the functor induced by straightening the $i$-th inert map $\rho_i\colon \langle n \rangle \to \langle 1 \rangle$ using (i). Then the collection of maps $(\rho_{i!})_i$ induces an equivalence
	\begin{align}
		\cat O^\tensor_{\langle n \rangle} \isoto \cat O^n.
	\end{align}

	\item Let $\alpha\colon \langle m \rangle \to \langle n \rangle$ be a map in $\Fin_*$ and for $X \in \cat O_{\langle n \rangle}$ and $Y \in \cat O_{\langle m \rangle}$ we denote by $\Hom_{\cat O^\tensor}^\alpha(Y, X) \subseteq \Hom_{\cat O^\tensor}(Y, X)$ the full subanima spanned by the morphisms over $\alpha$. Then
	\begin{align}
		\Hom^\alpha_{\cat O^\tensor}(Y, X) = \prod_{i=1}^n \Hom^{\rho_i \alpha}_{\cat O^\tensor}(Y, \rho_{i!}(X)).
	\end{align}
\end{enumerate}
We call $\cat O = \cat O^\tensor_{\langle 1 \rangle}$ the \emph{underlying category} of $\cat O^\tensor$. A \emph{map of operads} $\cat O^\tensor \to \cat O'^\tensor$ is a functor over $\Fin_*$ which preserves inert morphisms; we denote by
\begin{align}
	\Alg(\cat O, \cat O') \subseteq \Fun_{\Fin_*}(\cat O^\tensor, \cat O'^\tensor)
\end{align}
the full subcategory spanned by the operad maps. Moreover, we let
\begin{align}
	\Op \subset \Cat_{/\Fin_*}
\end{align}
denote the non-full subcategory of operads and operad maps.
\end{definition}

\begin{remark}
Our notation $\Alg(\cat O, \cat O')$ is non-standard: in \cite[Definition~2.1.3.1]{HA} the notation $\Alg_{\cat O}(\cat O')$ is used instead. We chose to deviate from that notation to make it more clear that $\Alg(\cat O, \cat O')$ is a morphism category (cf. \cref{rslt:2-categorical-enhancement-of-Op}). More generally, we will use $\Alg_{\cat O}$ to denote the morphism categories in $\Op_{/\cat O^\tensor}$ for any fixed operad $\cat O^\tensor$.
\end{remark}

By property (ii) in \cref{def:operad} we can identify objects in $\cat O^\tensor$ by tuples $(X_1, \dots, X_n)$ of objects in $\cat O$; we will usually denote the corresponding object in $\cat O^\tensor$ by $X_1 \dsum \dots \dsum X_n$. The homomorphism anima in $\cat O^\tensor$ encode \enquote{operations} in $\cat O$, see e.g.\ the discussion after \cref{def:monoidal-categories} below.

\begin{definition} \label{def:basic-operads}
\begin{defenum}
	\item $\Comm^\tensor$ is the operad $\Comm^\tensor \coloneqq \Fin_*$. We call it the \emph{commutative operad}. In the literature this operad is often also denoted by $\mathbb E_\infty$.

	\item $\Ass^\tensor$ is the category with underlying objects the same as $\Fin_*$, but where a morphism $\langle n \rangle \to \langle m \rangle$ is a pair $(\alpha, (\preceq_i)_{1 \le i \le m})$ such that $\alpha\colon \langle n \rangle \to \langle m \rangle$ is a map in $\Fin_*$ and for each $i$, $\preceq_i$ is a linear ordering on $\alpha^{-1}(i)$ (see \cite[Remark 4.1.1.4]{HA}). The projection $\Ass^\tensor \to \Fin_*$ makes $\Ass^\tensor$ into an operad, which we call the \emph{associative operad}. In the literature this operad is often also denoted by $\mathbb E_1$.

	\item $\LM^\tensor$ is the operad from \cite[Notation~4.2.1.6]{HA}. Its objects are pairs $(\langle n \rangle, S)$, where $\langle n \rangle$ is an object of $\Fin_*$ and $S$ is a subset of $\langle n \rangle^\circ$. A morphism $(\langle n \rangle, S) \to (\langle n' \rangle, S')$ is a morphism $\alpha\colon \langle n \rangle \to \langle n' \rangle$ in $\Ass^\tensor$ such that $\alpha(S \union \{ * \}) \subset S' \union \{ * \}$ and for every $s' \in S'$, $\alpha^{-1}(s')$ contains exactly one element of $S$ and that element is maximal with respect to the linear ordering on $\alpha^{-1}(s')$ induced by $\alpha$.

	We denote $\mathfrak a \coloneqq (\langle 1 \rangle, \emptyset)$ and $\mathfrak m \coloneqq (\langle 1 \rangle, \langle 1 \rangle^\circ)$, which are the two objects of the underlying category $\LM$ of $\LM^\tensor$. Note that there is a natural embedding $\Ass^\tensor \injto \LM^\tensor$ sending $\langle 1 \rangle$ to $\mathfrak a$ and a natural projection $\LM^\tensor \to \Ass^\tensor$ forgetting the set $S$.
\end{defenum}
\end{definition}

\begin{remark} \label{rmk:non-symmetric-vs-symmetric-operads}
In the literature there is a second approach to operads, often called \emph{non-symmetric operads}, see e.g.\ \cite[\S3]{Gepner-Haugseng.2015}. These non-symmetric operads are also introduced in \cite[\S4.1.3]{HA}, where they are called \emph{planar operads}. The definition of non-symmetric operads is similar to the definition of operads introduced above, except that the role of $\Fin_*$ is now played by $\bbDelta^\op$. By \cite[Theorem~4.1.3.14]{HA} the category of non-symmetric operads is isomorphic to the category $\Op_{/\Ass^\tensor}$ of operads over $\Ass^\tensor$, so in this paper we will usually ignore the difference.
\end{remark}

The operad $\Comm^\tensor$ encodes a commutative and associative binary operation, the operad $\Ass^\tensor$ encodes an associative binary operation and the operad $\LM^\tensor$ encodes a left action of an associative algebra on a module. This intuition becomes more concrete by the following definitions.

\begin{definition} \label{def:monoidal-categories}
\begin{defenum}
	\item Given an operad $\cat O^\tensor$, an \emph{$\cat O$-monoidal category} is a cocartesian fibration of operads $\cat C^\tensor \to \cat O^\tensor$ (see \cite[Definition~2.1.2.13]{HA}). Given $\cat O$-monoidal categories $\cat C_1^\tensor$ and $\cat C_2^\tensor$, an \emph{$\cat O$-monoidal functor $\cat C_1 \to \cat C_2$} is a morphism $F\colon \cat C_1^\tensor \to \cat C_2^\tensor$ of operads over $\cat O^\tensor$ which respects cocartesian morphisms over $\cat O$ (see \cite[Definition~2.1.3.7]{HA}).

	\item A \emph{symmetric monoidal category} is a $\Comm$-monoidal category. A \emph{lax symmetric monoidal functor} of symmetric monoidal categories is a morphism of operads (over $\Comm^\tensor$), and a \emph{symmetric monoidal functor} of symmetric monoidal categories is a $\Comm$-monoidal functor. We denote by $\CMon \subset \Op$ the non-full subcategory of symmetric monoidal categories, where the morphisms are the symmetric monoidal functors.

	\item A \emph{monoidal category} is an $\Ass$-monoidal category. A \emph{lax monoidal functor} of monoidal categories is a morphism of operads over $\Ass^\tensor$, and a \emph{monoidal functor} of monoidal categories is an $\Ass$-monoidal functor. We denote by $\Mon \subset \Op_{/\Ass^\tensor}$ the non-full subcategory of monoidal categories, where the morphisms are the monoidal functors.

	\item We denote by $\LMod \subset \Op_{/\LM^\tensor}$ the non-full subcategory of $\LM$-monoidal categories and $\LM$-monoidal functors. Pullback along the map $\Ass^\tensor \to \LM^\tensor$ provides a functor $\LMod \to \Mon$. For every monoidal category $\cat V^\tensor \in \Mon$ we denote by $\LMod_{\cat V}$ the fiber over $\cat V$. The objects of $\LMod_{\cat V}$ are called the \emph{$\cat V$-linear categories} and the morphisms are called the \emph{$\cat V$-linear functors}. We can similarly define \emph{lax $\cat V$-linear functors} by allowing all operad maps over $\LM^\tensor$ instead of only the $\LM$-monoidal ones.
\end{defenum}
By abuse of notation we will often denote a (symmetric) monoidal category $\cat C^\tensor$ simply by $\cat C$, with the operad $\cat C^\tensor$ being implicit. Similarly, we will often denote a $\cat V$-linear category $\cat O^\tensor$ simply by $\cat C$, where $\cat C$ is the fiber $\cat O^\tensor_{\mathfrak m}$.
\end{definition}

Given a symmetric monoidal category $\cat C^\tensor$, we can straighten the cocartesian fibration $\cat C^\tensor \to \Fin_*$ to obtain a functor $F\colon \Fin_* \to \Cat$ such that $F(\langle n \rangle) = \cat C^n$ for all $n$. This recovers the discussion at the beginning of this subsection and thus equips $\cat C$ with a binary operation $\tensor\colon \cat C \times \cat C \to \cat C$ together with a unit $\one \in \cat C$. Furthermore, we see that for $X_1, X_2, Y \in \cat C$, the maps $X_1 \dsum X_2 \to Y$ in $\cat C^\tensor$ lying over the non-inert map $\alpha\colon \langle 2 \rangle \to \langle 1 \rangle$ are equivalent to the maps $X_1 \tensor X_2 \to Y$ in $\cat C$.

Given two symmetric monoidal categories $\cat C_1^\tensor$ and $\cat C_2^\tensor$, a lax symmetric monoidal functor $\cat C_1 \to \cat C_2$ is a functor $f\colon \cat C_1 \to \cat C_2$ together with maps
\begin{align}
	\one_{\cat C_2} \to f(\one_{\cat C_1}), \qquad f(X) \tensor_{\cat C_2} f(Y) \to f(X \tensor_{\cat C_1} Y)
\end{align}
for all $X, Y \in \cat C_1$, plus higher homotopies exhibiting the \enquote{naturality} of these maps. A symmetric monoidal functor is a lax symmetric monoidal functor where all of the above maps are isomorphisms.

The above discussion of symmetric monoidal categories applies analogously also to monoidal categories, with the main difference being that $\tensor$ is not symmetric. Note that every symmetric monoidal category can be made into a monoidal category by applying $\blank \times_{\Comm^\tensor} \Ass^\tensor$; this amounts to forgetting the commutator (and its higher coherences).

By a similar reasoning one can make explicit the definition of linear categories. Namely, given a monoidal category $\cat V$, a $\cat V$-linear category is a category $\cat C$ with an operation
\begin{align}
	\tensor\colon \cat V \times \cat C \to \cat C, \qquad (V, X) \mapsto V \tensor X
\end{align}
together with higher homotopies exhibiting the associativity and unit of that action with respect to $\cat V$. A $\cat V$-linear functor $\cat C_1 \to \cat C_2$ of $\cat V$-linear categories is a functor $f\colon \cat C_1 \to \cat C_2$ together with isomorphisms
\begin{align}
	V \tensor f(X) \isoto f(V \tensor X)
\end{align}
for all $V \in \cat V$ and $X \in \cat C$, plus higher homotopies exhibiting the \enquote{naturality} of these isomorphisms.

\begin{example} \label{ex:cocartesian-monoidal-structure}
For every category $\cat C$ there is an associated operad $\cat C^\amalg \to \Comm^\tensor$, called the \emph{cocartesian operad associated with $\cat C$}, see \cite[Construction~2.4.3.1]{HA}. The underlying category of $\cat C^\amalg$ is $\cat C$. If $\cat C$ admits finite coproducts then $\cat C^\amalg$ is a symmetric monoidal category such that the tensor operation is given by the coproduct. 
\end{example}

\begin{remark} \label{rslt:functor-from-C-amalg-to-C}
Let $\cat C$ be a category with finite coproducts. Then there is a natural functor $\cat C^\amalg \to \cat C$ sending $X_1 \dsum \dots \dsum X_n \mapsto X_1 \dunion \dots \dunion X_n$. Namely, by the explicit construction of $\cat C^\amalg$, it is enough to construct, for every category $\cat K$ over $\Fin_*$, a natural functor
\begin{align}
	\Fun(\Gamma^* \times_{\Fin_*} \cat K, \cat C) \to \Fun(\cat K, \cat C),
\end{align}
where $\Gamma^*$ is the category from \cite[Construction~2.4.3.1]{HA}. There is a natural functor from right to left induced by the projection $\Gamma^* \to \Fin_*$. We claim that this functor has a left adjoint, producing the desired functor. In other words, we need to see that the left Kan extension of a functor $F\colon \Gamma^* \times_{\Fin_*} \cat K \to \cat C$ along $\Gamma^* \times_{\Fin_*} \cat K \to \cat K$ exists. Given $X \in \cat K$, we need to see that the colimit of $F$ over
\begin{align}
	\cat I_X \coloneqq (\Gamma^* \times_{\Fin_*} \cat K) \times_{\cat K} \cat K_{/X} = \Gamma^* \times_{\Fin_*} \cat K_{/X}
\end{align}
exists. Let $\langle n \rangle$ be the image of $X$ in $\Fin_*$ and consider the full subcategory of $\cat I_X$ spanned by the objects $((\langle n \rangle, i), X)$ for $i = 1, \dots, n$. This subcategory is cofinal in $\cat I_X$ and it is isomorphic to a disjoint union of $n$ copies of $*$. Thus the colimit over $\cat I_X$ is computed by an $n$-fold coproduct in $\cat C$, which exists by assumption.
\end{remark}

\begin{example} \label{ex:cartesian-monoidal-structure}
If $\cat C$ is a category which admits finite products, then there is an associated symmetric monoidal category $\cat C^\times \to \Comm^\tensor$, see \cite[Proposition~2.4.1.5.(4)]{HA}. The underlying category of $\cat C^\times$ is $\cat C$ and the tensor operation is the product $\times$. The operad $\cat C^\times$ can be characterized by the following universal property (see \cite[Proposition~2.4.1.7]{HA}): Given any operad $\cat O^\tensor$, there is a natural isomorphism
\begin{align}
	\Hom(\cat O^\tensor, \cat C^\times) = \Fun^\lax(\cat O^\tensor, \cat C)^\simeq,
\end{align}
where $\Fun^\lax(\cat O^\tensor, \cat C)$ denotes the full subcategory of $\Fun(\cat O^\tensor, \cat C)$ consisting of the \emph{lax cartesian structures}, i.e.\ those functors $F\colon \cat O^\tensor \to \cat C$ such that for all $X_1, \dots, X_n \in \cat O$ the natural map $F(X_1 \dsum \dots \dsum X_n) \isoto F(X_1) \times \dots \times F(X_n)$ is an isomorphism (see \cite[Definition~2.4.1.1]{HA}).
\end{example}

\begin{example} \label{ex:lax-monoidal-functor-corep-by-1}
By \cref{ex:cartesian-monoidal-structure} there is a symmetric monoidal category $\Ani^\times$ whose tensor operation is the product of anima. Given an operad $\cat C^\tensor$, let $\emptyset_{\cat C} \in \cat C^\tensor_{\langle 0 \rangle}$ be the unique object in the fiber over $\langle 0 \rangle \in \Fin_*$. Then the functor $\Hom(\emptyset_{\cat C}, \blank)\colon \cat C^\tensor \to \Ani$ is a lax cartesian structure and hence defines a map of operads $\cat C^\tensor \to \Ani^\times$.

Suppose now that $\cat C^\tensor$ is a (symmetric) monoidal category. Then the map of operads we just constructed is a lax (symmetric) monoidal functor
\begin{align}
	\cat C \to \Ani, \qquad X \mapsto \Hom(\one_{\cat C}, X).
\end{align}
For any two objects $X, Y \in \cat C$, the lax monoidal structure on the functor provides a natural map
\begin{align}
	\Hom(\one_{\cat C}, X) \times \Hom(\one_{\cat C}, Y) \to \Hom(\one_{\cat C}, X \tensor Y),
\end{align}
which is obtained by sending two maps $f\colon \one_{\cat C} \to X$, $g\colon \one_{\cat C} \to Y$ to the map $\one_{\cat C} \isoto \one_{\cat C} \tensor \one_{\cat C} \xto{f \tensor g} X \tensor Y$.
\end{example}

To a monoidal structure on a category we can attach internal hom objects, assuming they exist. They will play an important role in the theory of enriched categories in \cref{sec:enr}.

\begin{definition} \label{def:internal-hom}
Let $\cat C$ be a monoidal category. Given objects $X, Y \in \cat C$, an \emph{internal hom object} for $X$ and $Y$ is an object $\iHom_{\cat C}(Y, X) \in \cat C$ together with a map $\iHom_{\cat C}(Y, X) \tensor Y \to X$ such that for all $Z$ the induced map
\begin{align}
	\Hom(Z, \iHom_{\cat C}(Y, X)) \isoto \Hom(Z \tensor_{\cat C} Y, X)
\end{align}
is an isomorphism. We say that $\cat C$ is a \emph{closed} monoidal category if all internal hom objects exist.
\end{definition}

Next we discuss the behavior of monoidal functors with respect to adjunctions. For example, if a symmetric monoidal functor $f^*\colon \cat C \to \cat C'$ has a right adjoint $f_*$ then for all $X, Y \in \cat C'$ there is a natural map $f_* X \tensor f_* Y \to f_*(X \tensor Y)$, obtained via adjunction from $f^* f_* X \tensor f^* f_* Y \to X \tensor Y$ (which in turn comes from the counit $f^* f_* \to \id$). This should equip $f_*$ with a lax symmetric monoidal structure. The following results show that this intuition is indeed correct:

\begin{lemma} \label{rslt:monoidal-right-adjoint}
\begin{lemenum}
	\item Let $f^*\colon \cat C \to \cat C'$ be a (symmetric) monoidal functor of (symmetric) mo\-noi\-dal categories. Assume that the underlying functor admits a right adjoint $f_*$. Then $f_*$ lifts to a lax (symmetric) monoidal functor $f_*\colon \cat C' \to \cat C$ which is right adjoint to $f^*$ as a functor of operads.

	\item \label{rslt:monoidal-linear-right-adjoints} Let $\cat V$ be a monoidal category and let $f^*\colon \cat C \to \cat C'$ be a $\cat V$-linear functor of $\cat V$-linear categories. Assume that the underlying functor admits a right adjoint $f_*$. Then $f_*$ lifts to a lax $\cat V$-linear functor $\cat C' \to \cat C$ which is right adjoint to $f^*$ as a functor of $\LM^\tensor$-operads.
\end{lemenum}
\end{lemma}
\begin{proof}
Both statements are a special case of \cite[Corollary~7.3.2.7]{HA}.
\end{proof}

\begin{lemma} \label{rslt:monoidal-left-adjoints}
\begin{lemenum}
	\item Let $f^*\colon \cat C \to \cat C'$ be a lax (symmetric) monoidal functor of (symmetric) monoidal categories. Assume that the underlying functor admits a left adjoint $f_!$ such that for all $X, Y \in \cat C'$ the natural map
	\begin{align}
		f_!(X \tensor_{\cat C'} Y) \isoto f_! X \tensor_{\cat C} f_! Y
	\end{align}
	is an isomorphism. Then $f_!$ lifts to a (symmetric) monoidal functor $f_!\colon \cat C' \to \cat C$ which is left adjoint to $f^*$ as a functor of operads.

	\item \label{rslt:monoidal-linear-left-adjoints} Let $\cat V$ be a monoidal category and let $f^*\colon \cat C \to \cat C'$ be a $\cat V$-linear functor of $\cat V$-linear categories. Assume that the underlying functor admits a left adjoint $f_!$ such that for all $V \in \cat V$ and $M \in \cat C'$ the natural map
	\begin{align}
		f_!(V \tensor M) \isoto V \tensor f_! M
	\end{align}
	is an isomorphism. Then $f_!$ lifts to a $\cat V$-linear functor $f_!\colon \cat C' \to \cat C$ which is left adjoint to $f^*$ as a functor of $\LM^\tensor$-operads.
\end{lemenum}
\end{lemma}
\begin{proof}
Both statements are a special case of \cite[Corollary~7.3.2.12]{HA}.
\end{proof}

\begin{remark}
There is an apparent asymmetry between \cref{rslt:monoidal-left-adjoints} and \cref{rslt:monoidal-right-adjoint}. One can make the statements symmetric by introducing \emph{oplax} monoidal functors. We refer the reader to \cite[Corollary~3.4.8]{Haugseng.2023} for a precise statement.
\end{remark}

An important concept in the theory of monoidal categories is that of dualizable objects, which form a particularly nice subcategory of a monoidal category. In the following we recall the definition and basic properties of dualizable objects. For simplicity, we will stick to the case of \emph{symmetric} monoidal categories.

\begin{definition} \label{def:dualizable-object}
Let $\cat C$ be a symmetric monoidal category. An object $P \in \cat C$ is called \emph{dualizable} if there are an object $P^\vee \in \cat C$, called the \emph{dual} of $P$, and morphisms $\ev_P\colon P^\vee \tensor P \to \one_{\cat C}$ and $i_P\colon \one_{\cat C} \to P \tensor P^\vee$, called the \emph{evaluation map} and \emph{coevaluation map} respectively, such that there are commutative diagrams
\begin{equation}
	\begin{tikzcd}[column sep=3em]
		P \arrow[dr,"\id",swap] \arrow[r,"i_P \tensor \id"] & P \tensor P^\vee \tensor P \arrow[d,"\id \tensor \ev_P"] \\
		& P
	\end{tikzcd}
	\qquad
	\begin{tikzcd}[column sep=3em]
		P^\vee \arrow[dr,"\id"] \arrow[r,"\id \tensor i_P"] & P^\vee \tensor P \tensor P^\vee \arrow[d,"\ev_P \tensor \id"] \\
		& P^\vee
	\end{tikzcd}
\end{equation}
in $\cat C$.
\end{definition}

\begin{lemma} \label{rslt:properties-of-dualizable-objects}
Let $\cat C$ be a symmetric monoidal category and $P \in \cat C$ an object. We denote $P^\vee \coloneqq \iHom_{\cat C}(P, \one_{\cat C})$ if it exists. Then the following are equivalent:
\begin{lemenum}
	\item $P$ is dualizable.
	\item The objects $P^\vee$ and $\iHom_{\cat C}(P, P)$ exist and the natural map $P \tensor_{\cat C} P^\vee \isoto \iHom_{\cat C}(P, P)$ is an isomorphism.
	\item For every $M \in \cat C$ the object $\iHom_{\cat C}(P, M)$ exists and the natural map
	\begin{align}
		P^\vee \tensor_{\cat C} M \isoto \iHom_{\cat C}(P, M)
	\end{align}
	is an isomorphism.
\end{lemenum}
If this is the case then $P^\vee$ is the dual of $P$. Moreover, $P^\vee$ is itself dualizable and the natural map $P \isoto (P^\vee)^\vee$ is an isomorphism.
\end{lemma}
\begin{proof}
Note that if $P$ is dualizable with dual $P^*$ then for all $M, N \in \cat C$ the natural map
\begin{align}
	\Hom(N, M \tensor P^*) \isoto \Hom(N \tensor P, M),
\end{align}
given by sending $f\colon N \to M \tensor P^*$ to the composition $N \tensor P \xto{f \tensor \id} M \tensor P^* \tensor P \xto{\id \tensor \ev_P} M$, is an isomorphism (one can explicitly define the inverse of that map). This implies that $P^* = P^\vee$ (in particular $P^\vee$ exists) and that (iii) holds. Part (ii) is a special case of (iii), so it only remains to show that (ii) implies (i). This is an easy exercise, see \cite[Lemma~6.2]{Mann.2022b}.
\end{proof}

\begin{example}
If $A$ is a ring and $\cat C = \D(A)$ is the (derived) category of $A$-modules then for an object $P \in \D(A)$ the following are equivalent:
\begin{enumerate}[(i)]
	\item $P$ is dualizable.
	\item $P$ is compact, i.e.\ $\Hom(P, \blank)\colon \D(A) \to \Ani$ preserves all filtered colimits.
	\item $P$ is perfect, i.e.\ can be obtained from $A$ in finitely many steps using finite (co)limits (i.e.\ fiber sequences) and retracts.
\end{enumerate}
This is a classical result, albeit formulated in a modern language. We sketch the proof. The core observation is that $A$ is compact in $\D(A)$ and a generator (i.e.\ $\Hom(A, \blank)$ is conservative). This immediately implies that compact objects of $\D(A)$ are exactly those that can be obtained from $A$ using finite colimits and retracts (see e.g.\ \cite[Lemma~A.2.1]{Mann.2022a}), proving the equivalence of (ii) and (iii). Now if $P$ is dualizable then $\Hom(P, \blank) = \Hom(A, P^\vee \tensor \blank)$, proving that $P$ is compact. On the other hand, it follows easily from \cref{rslt:properties-of-dualizable-objects} that dualizable objects in $\D(A)$ are stable under finite colimits and retracts, hence the dualizability of $A$ implies that every perfect object is dualizable.
\end{example}

There are many abstract constructions one can perform with operads, some of which are covered in \cite[\S2.2]{HA}. In this paper we additionally need to work with limits of operads, for which we did not find a good reference. Therefore we provide the corresponding result ourselves:

\begin{lemma} \label{rslt:limits-of-operads}
Let $\cat O^\tensor$ be an operad. Then the category $\Op_{/\cat O^\tensor}$ has all small limits and the forgetful functor $\Op_{/\cat O^\tensor} \to \Cat_{/\cat O^\tensor}$ preserves them.
\end{lemma}
\begin{proof}
By the computation of limits in slice categories in terms of limits in the ambient category (see \cref{rslt:limits-in-overcategory}), we immediately reduce to the case that $\cat O^\tensor = \Fin_*$, i.e.\ we need to show that the category $\Op$ has all small limits and that the forgetful functor $\Op \to \Cat_{/\Fin_*}$ preserves them. To this end, let $(\cat V^\tensor_i)_{i \in I}$ be a diagram in $\Op$ and let $\cat V^\tensor \coloneqq \varprojlim_i \cat V^\tensor_i$ be their limit in $\Cat_{/\Fin_*}$. We need to show the following:
\begin{enumerate}[(a)]
 	\item The functor $\cat V^\tensor \to \Fin_*$ makes $\cat V^\tensor$ into an operad.
 	\item For every operad $\cat V'^\tensor$ with compatible functors $f_i\colon \cat V'^\tensor \to \cat V_i^\tensor$ over $\Fin_*$, the induced functor $f\colon \cat V'^\tensor \to \cat V^\tensor$ is an operad map if and only if all $f_i$ are operad maps.
\end{enumerate}
The diagram $I \to \Cat_{/\Fin_*}$, $i \mapsto \cat V^\tensor_i$, induces a diagram $I^\triangleright \to \Cat$ sending $* \mapsto \Fin_*$. The limit $\cat V^\tensor$ can be computed as the limit over $I^\triangleright$ in $\Cat$ (with the canonical induced functor to $\Fin_*$). In particular the objects of $\cat V^\tensor$ can be identified with compatible families of objects $(X_i)_{i \in I}$ in $\cat V^\tensor_i$ lying over a fixed object in $\Fin_*$, and the morphism anima in $\cat V^\tensor$ are obtained as certain limits over $I$. From these observations it is easy to see that $\cat V^\tensor$ satisfies conditions (i), (ii) and (iii) in \cref{def:operad} (for (i) we use that cocartesian fibrations are stable under limits over the base, see \cref{rslt:cocartesian-morphism-in-limit-of-fibrations}), proving (a). Moreover, we note that a morphism in $\cat V^\tensor$ is inert if and only if its projection to all $\cat V^\tensor_i$ is inert (for the \enquote{only if} part use the uniqueness of cocartesian lifts); this easily implies (b).
\end{proof}

\subsection{Algebras and modules}

A fundamental notion associated with operads is that of (commutative) algebras, i.e.\ objects equipped with a homotopy associative (and commutative) multiplication, as well as modules over such algebras. They are defined as follows:

\begin{definition} \label{def:algebras-and-modules}
Let $\cat C^\tensor$ be an operad.
\begin{defenum}
	\item We denote by $\CAlg(\cat C) \coloneqq \Alg(\Comm, \cat C)$ the functor category of operad maps $\Comm^\tensor \to \cat C^\tensor$ (see \cite[Definition~2.1.3.1]{HA}). An object in $\CAlg(\cat C)$ is called a \emph{commutative algebra in $\cat C$}.

	\item If $\cat C^\tensor$ is equipped with a map to $\Ass^\tensor$ then we denote by $\Alg(\cat C) \coloneqq \Alg_\Ass(\Ass, \cat C)$ the functor category of operad maps $\Ass^\tensor \to \cat C^\tensor$ over $\Ass^\tensor$ (see \cite[Definition~4.1.1.6]{HA}). An object in $\Alg(\cat C)$ is called an \emph{algebra in $\cat C$}.

	\item If $\cat C^\tensor$ is equipped with a map to $\LM^\tensor$ then we denote by $\LMod(\cat C) \coloneqq \Alg_\LM(\LM, \cat C)$ the functor category of operad maps $\LM^\tensor \to \cat C^\tensor$ over $\LM^\tensor$. Restriction along the map $\Ass^\tensor \to \LM^\tensor$ provides a functor $\LMod(\cat C) \to \Alg(\cat C)$. For every algebra $A \in \Alg(\cat C)$ we denote by $\LMod_A(\cat C)$ the fiber over $A$ and call it the \emph{category of left $A$-modules in $\cat C$}.
\end{defenum}
\end{definition}

Spelling out the definition, we see that an algebra object in a monoidal category $\cat C$ is given by an object $A \in \cat C$ together with a unit map $\one_{\cat C} \to A$ and a multiplication $A \tensor A \to A$ as well as higher coherences witnessing the associativity of the multiplication. A commutative algebra object contains additional coherences witnessing the commutativity of the operation. An object in $\LMod(\cat C)$ consists of an algebra $A \in \Alg(\cat C)$ and an object $M \in \cat C$ (the image of $\mathfrak m \in \LM^\tensor$) together with a map $A \tensor M \to M$ as well as higher coherences witnessing the associativity of this $A$-multiplication.

\begin{example} \label{ex:trivial-algebra}
If $\cat C$ is a symmetric monoidal category then $\CAlg(\cat C)$ has an initial object, which is given by the canonical commutative algebra structure on $\one_{\cat C}$ coming from the equivalence $\one_{\cat C} \tensor_{\cat C} \one_{\cat C} \isoto \one_{\cat C}$. Similarly, if $\cat C$ is a monoidal category then $\Alg(\cat C)$ has an initial object given by $\one_{\cat C}$. See \cite[Lemma~3.2.1.10]{HA} and \cite[Example~3.2.1.6]{HA} for a proof.
\end{example}

\begin{example} \label{ex:monoidal-categories-as-algebras}
By \cite[Remark~2.4.2.6]{HA} we have $\CMon = \CAlg(\Cat^\times)$, i.e.\ a symmetric monoidal category is the same as a commutative algebra in the cartesian symmetric monoidal structure on $\Cat$ (defined in \cref{ex:cartesian-monoidal-structure}). Similarly we have $\Mon = \Alg(\Cat^\times)$ and $\LMod = \LMod(\Cat^\times)$.
\end{example}

In the following we recall some basic results about algebras and modules in a (symmetric) monoidal category, which are generalizations of well-known results about classical rings and modules.

\begin{lemma} \label{rslt:lim-and-colim-of-algebras}
Let $\cat C$ be a symmetric monoidal category.
\begin{lemenum}
	\item \label{rslt:lim-of-algebras} The forgetful functor $\CAlg(\cat C) \to \cat C$ is conservative and preserves all limits. Moreover, the existence of limits in $\CAlg(\cat C)$ can be checked on underlying objects in $\cat C$.

	\item Suppose that $\cat C$ has all small colimits and that the tensor product on $\cat C$ preserves small colimits in each argument. Then $\CAlg(\cat C)$ has all small colimits and the forgetful functor $\CAlg(\cat C) \to \cat C$ preserves all sifted colimits.

	\item The coproduct of two commutative algebras $A, B \in \CAlg(\cat C)$ exists and is computed as $A \tensor_{\cat C} B$ on underlying objects.
\end{lemenum}
Part (i) and (ii) also hold if $\cat C$ is a monoidal category and $\CAlg$ is replaced by $\Alg$.
\end{lemma}
\begin{proof}
Part (i) follows from \cite[Lemma~3.2.2.6]{HA} and \cite[Corollary~3.2.2.5]{HA} and part (ii) follows from \cite[Corollary~3.2.3.3]{HA}. These references also work for $\Alg$ instead of $\CAlg$. Part (iii) follows from \cite[Proposition~3.2.4.7]{HA}.
\end{proof}

\begin{lemma} \label{rslt:modules-over-commutative-algebras}
Let $\cat C$ be a symmetric monoidal category and assume that $\cat C$ admits colimits over $\bbDelta^\op$ and that the tensor product on $\cat C$ preserves these colimits in each argument. Then:
\begin{lemenum}
	\item For every $A \in \CAlg(\cat C)$ the category $\LMod_A(\cat C)$ of (left) $A$-modules inherits a symmetric monoidal structure $\tensor_A$ with tensor unit $A$.

	\item For every map $f\colon A \to B$ in $\CAlg(\cat C)$ there is a forgetful functor $f_*\colon \LMod_B(\cat C) \to \LMod_A(\cat C)$ which is the identity on the underlying object of $\cat C$. It admits a symmetric monoidal left adjoint $f^*\colon \LMod_A(\cat C) \to \LMod_B(\cat C)$ which acts as $f^*(M) = M \tensor_A B$ on underlying $A$-modules.

	\item \label{rslt:functoriality-of-modules-over-algebras} The functors $f^*$ from (ii) are functorial in $f$, i.e.\ there is a functor $\CAlg(\cat C) \to \CMon$ sending $A \mapsto \LMod_A(\cat C)^\tensor$ and $f \mapsto f^*$.

	\item For every $A \in \CAlg(\cat C)$ there is a natural equivalence of categories $\CAlg(\LMod_A(\cat C)) = \CAlg(\cat C)_{A/}$, i.e.\ an algebra in the category of $A$-modules is the same as an algebra in $\cat C$ with a map from $A$.
\end{lemenum}
\end{lemma}
\begin{proof}
Part (i) is found in \cite[Theorem~4.5.2.1]{HA}. The functor $f_*$ in (ii) is constructed in \cite[Corollary~3.4.3.3]{HA} (using \cite[Corollary~4.5.1.6]{HA}) and the functor $f^*$ in (ii) is constructed in \cite[Theorem~4.5.3.1]{HA}, which also shows (iii). Part (iv) is shown in \cite[Corollary~3.4.1.7]{HA}.
\end{proof}

\Cref{rslt:lim-and-colim-of-algebras} is self-explanatory (the reader who is unfamiliar with sifted colimits may note that an important special case of sifted colimits are filtered colimits): It tells us that under mild assumptions on $\cat C$, limits and colimits of (commutative) algebras in $\cat C$ exist and are computed as expected. Regarding \cref{rslt:modules-over-commutative-algebras}, let us explain how the relative tensor product $\tensor_A$ is constructed. Given $A \in \CAlg(\cat C)$ and $M, N \in \LMod_A(\cat C)$ one can construct the following diagram $\bbDelta^\op \to \cat C$ (where we omit the degeneracy maps):
\begin{equation}\begin{tikzcd}
	\dots \arrow[r,shift right=3] \arrow[r,shift right=1] \arrow[r,shift left=1] \arrow[r,shift left=3] & M \tensor A \tensor A \tensor N \arrow[r,shift right=2] \arrow[r] \arrow[r,shift left=2] & M \tensor A \tensor N \arrow[r,shift right] \arrow[r,shift left] & M \tensor N
\end{tikzcd}\end{equation}
Here the two maps $M \tensor A \tensor N \to M \tensor N$ are given by $m \tensor a \tensor n \mapsto am \tensor n$ and $m \tensor a \tensor n \mapsto m \tensor an$, respectively. The relative tensor product $M \tensor_A N$ is defined to be the colimit over the above diagram. In case $\cat C$ is an ordinary category, the colimit over this diagram only depends on its truncation to the last two objects and computes the quotient of the two maps $M \tensor A \tensor N \to M \tensor N$, resulting in the usual explicit construction of the tensor product.

\begin{example} \label{ex:Einfty-rings-and-modules}
Let $\Sp$ denote the category of spectra (see \cite[\S1.4]{HA}). This is a symmetric monoidal presentable stable category (see \cite[\S4.8.2]{HA}), so in particular we can consider algebras and modules inside of it. We denote $\CAlg \coloneqq \CAlg(\Sp)$ and call the objects $\Lambda \in \CAlg$ the \emph{$\Einfty$-rings}. For every $\Einfty$-ring $\Lambda$ we denote $\Mod_\Lambda \coloneqq \LMod_\Lambda(\Sp)$ the associated symmetric monoidal category of \emph{$\Lambda$-modules}. Recall that $\Sp$ comes equipped with a t-structure whose heart is the category of abelian groups, inducing a lax symmetric monoidal embedding $\Ab \injto \Sp$. In particular every ordinary ring $\Lambda$ can be seen as an $\Einfty$-ring (see also \cite[Proposition~7.1.3.18]{HA}). Then $\Mod_\Lambda = \D(\Lambda)$ identifies with the derived category of $\Lambda$-modules (see \cite[Remark~7.1.1.16]{HA}).

For every $\Einfty$-ring $\Lambda$ the category $\Mod_\Lambda$ is compactly generated (see \cref{def:kappa-compactly-generated-category}) because the shifts of $\Lambda$ form compact generators (in the sense of \cite[Lemma~A.2.1]{Mann.2022a}). Now let $\Lambda = \varinjlim_i \Lambda_i$ be a filtered colimit of $\Einfty$-rings. Then
\begin{align}
	\Mod_\Lambda = \varinjlim_i \Mod_{\Lambda_i}
\end{align}
where the colimit on the right is taken in the category $\PrL$ of presentable categories. Indeed, by \cref{rslt:Ind-preserves-filtered-colim} the claim reduces to showing that $\Mod_\Lambda^\omega = \varinjlim_i \Mod_{\Lambda_i}^\omega$, where the colimit is now taken in $\Cat$. There is a natural functor from right to left and it is an easy exercise to show that this functor is fully faithful and essentially surjective (see \cite[Lemma~2.7.4]{Mann.2022a} for a similar argument).
\end{example}

\begin{example} \label{ex:algebra-and-modules-in-Cat}
The category $\Cat$ has all small limits and colimits and the product $\times$ preserves all small colimits in each argument (indeed, for a category $\cat C$ the functor $\cat C \times \blank$ has a right adjoint given by $\Fun(\cat C, \blank)$). Therefore, from \cref{ex:monoidal-categories-as-algebras,rslt:lim-and-colim-of-algebras,rslt:modules-over-commutative-algebras} we deduce that for every symmetric monoidal category $\cat V^\tensor$ the category $\LMod_{\cat V}$ admits a symmetric monoidal structure $\tensor_{\cat V}$ and for every symmetric monoidal functor $\cat V^\tensor \to \cat W^\tensor$ of symmetric monoidal categories there is an adjunction
\begin{align}
	- \tensor_{\cat V} \cat W\colon \LMod_{\cat V} \rightleftarrows \LMod_{\cat W} \noloc \operatorname{forget},
\end{align}
where $- \tensor_{\cat V} \cat W$ is additionally equipped with a symmetric monoidal structure. Moreover, an algebra in $\LMod_{\cat V}$ is the same as a symmetric monoidal map $\cat V^\tensor \to \cat V'^\tensor$ of symmetric monoidal categories and the coproduct of two such algebras $\cat V^\tensor \to \cat V'^\tensor, \cat V''^\tensor$ is computed as $\cat V' \tensor_{\cat V} \cat V''$. In particular this also computes the pushout.
\end{example}

\begin{example}\label{ex:modules-over-trivial-algebra}
Let $\cat C$ be a symmetric monoidal category and recall the trivial commutative algebra $\one_{\cat C}$ from \cref{ex:trivial-algebra}. Then the forgetful map $\LMod_{\one_{\cat C}}(\cat C)\to \cat C$ is an equivalence of symmetric monoidal categories by \cite[Proposition~3.4.2.1]{HA}.
\end{example}

\begin{lemma}\label{rslt:extension-of-monoidal-functor-to-modules}
Let $F\colon \cat C^{\tensor} \to \cat D^{\tensor}$ be a lax symmetric monoidal functor of symmetric monoidal categories and $A \in \CAlg(\cat C)$ a commutative algebra. Then $F$ enhances to an operad map
\begin{align}
\LMod_A(\cat C)^{\tensor} \to \LMod_{F(A)}(\cat D)^{\tensor}.
\end{align}
\end{lemma}
\begin{proof}
Let $\LMod(\cat C)^{\tensor} \to \CAlg(\cat C) \times \Fin_*$ be the map of generalized operads from \cite[Definition~3.3.3.8]{HA} (applied with $\cat O = \Fin_*$). It satisfies $\LMod_A(\cat C)^{\tensor} = \LMod(\cat C)^{\tensor} \times_{\CAlg(\cat C)} \{A\}$ and we observe that this definition is functorial in a map of operads, i.e.\ there is a commutative diagram
\begin{equation}
\begin{tikzcd}
\LMod(\cat C)^{\tensor} \ar[d] \ar[r] & \LMod(\cat D)^{\tensor} \ar[d] \\
\CAlg(\cat C) \times \Fin_* \ar[r] & \CAlg(\cat D) \times \Fin_*.
\end{tikzcd}
\end{equation}
By restricting the left-hand side to the fiber over $A$, we obtain the desired operad map.
\end{proof}

In the setting of monoidal categories one can reverse the tensor operation, which will play an important role when dealing with enrichment:

\begin{definition}
There is a unique non-trivial map of operads $\rev\colon \Ass^\tensor \to \Ass^\tensor$ (denoted $\sigma$ in \cite[Remark~4.1.1.7]{HA}). It is the identity on objects and reverses total orderings on morphisms. Note that $\rev$ is involutive, i.e.\ $\rev^2 = \id_{\Ass^\tensor}$. Given an operad $\cat C^\tensor$ over $\Ass^\tensor$, we denote by $\cat C^{\tensor,\rev}$ the operad over $\Ass^\tensor$ which is given by the composition
\begin{align}
	\cat C^\tensor \to \Ass^\tensor \xto{\rev} \Ass^\tensor.
\end{align}
We call $\cat C^{\tensor,\rev}$ the \emph{reverse operad} over $\Ass^\tensor$. Note that if $\cat C^\tensor$ is a monoidal category then so is $\cat C^{\tensor,\rev}$ and we call it the \emph{reverse monoidal category}.
\end{definition}

If $\cat C^\tensor$ is a monoidal category, given by a tensor operation $\tensor_{\cat C}$ on $\cat C$, then $\cat C^{\tensor,\rev}$ encodes the monoidal structure on $\cat C$ with tensor operation
\begin{align}
	X \tensor_{\cat C^\rev} Y = Y \tensor_{\cat C} X.
\end{align}
Note that if $\cat C^\tensor$ is a \emph{symmetric} monoidal category then the associated monoidal category is naturally isomorphic to its own reverse category.

\begin{remark} \label{rmk:rev-algebras-same-as-algebras}
Precomposing with $\rev\colon \Ass^\tensor \isoto \Ass^\tensor$ induces an equivalence
\begin{align}
	\Alg(\cat C^\rev) = \Alg(\cat C)
\end{align}
for every operad $\cat C^\tensor \to \Ass^\tensor$. The image of an algebra $A$ under this equivalence is usually denoted $A^{\op}$ and referred to as the \emph{opposite algebra} of $A$.
\end{remark}

\subsection{Groups, actions and quotients} \label{sec:alg.grp}

In this subsection we recall basic notions about group actions and their quotients in an ambient category $\cat C$. Usually $\cat C$ will be a category of sheaves, e.g.\ the category $\Cond(\Ani)$ of condensed anima.

\begin{definition} \label{def:monoid-group-action}
Let $\cat C$ be a category with finite products.
\begin{defenum}
	\item A \emph{monoid} in $\cat C$ is an algebra for the cartesian symmetric monoidal structure on $\cat C$ (see \cref{ex:cartesian-monoidal-structure}). We denote $\Mon(\cat C) \coloneqq \Alg(\cat C^\times)$ the category of monoids in $\cat C$.

	\item A \emph{group} in $\cat C$ is a monoid $G \in \Mon(\cat C)$ such that the map $G^2 \to G^2$ representing $(g, g') \mapsto (g, gg')$ is an isomorphism. We denote by $\Grp(\cat C) \subseteq \Mon(\cat C)$ the full subcategory spanned by the groups in $\cat C$.

	\item A \emph{(left) action} of a monoid $M \in \Mon(\cat C)$ on an object $X \in \cat C$ is a left module for the cartesian symmetric monoidal structure on $\cat C$. We denote by $\LAct_M(\cat C) \coloneqq \LMod_M(\cat C)$ the category of left actions of $M$ and by $\LAct(\cat C) \coloneqq \LMod(\cat C)$ the category of left actions.
\end{defenum}
\end{definition}

Let us unpack \cref{def:monoid-group-action}. A monoid $G$ in $\cat C$ consists of an object $G \in \cat C$ together with a multiplication map $G \times G \to G$ and a unit $* \to G$ plus higher homotopies exhibiting the associativity. A $G$-action on an object $X \in \cat C$ consists of a map $G \times X \to X$ together with higher homotopies exhibiting the associativity of this action. If $\cat C$ is an ordinary category, all higher homotopies are fixed and hence reduce to the usual commutativity requirements of certain diagrams. Note that by \cref{rslt:lim-of-algebras} the category of monoids has a final object $*$ whose underlying object in $\cat C$ is the final object $*$. There is a different description of monoids and actions:

\begin{lemma} \label{rslt:monoids-and-actions-via-Delta}
Let $\cat C$ be a category with finite products.
\begin{lemenum}
	\item There is a fully faithful embedding $\Mon(\cat C) \injto \Fun(\bbDelta^\op, \cat C)$ whose essential image consists of those functors $M\colon \bbDelta^\op \to \cat C$ such that for all $n \ge 0$ we have
	\begin{align}
		M([n]) = M(\{0,1\}) \times \dots \times M(\{n-1, n\})
	\end{align}
	via the obvious face maps. In particular we require $M([0]) = *$.

	\item \label{rslt:actions-via-Delta} There is a fully faithful embedding $\LAct(\cat C) \injto \Fun(\bbDelta^\op \times \Delta^1, \cat C)$ whose essential image consists of those natural transformations $X \to M$ of simplicial objects $X, M\colon \bbDelta^\op \to \cat C$ such that $M$ satisfies the condition in (i) and for all $n \ge 0$ we have
	\begin{align}
		X([n]) = M([n]) \times X([0]),
	\end{align}
	where the map $X([n]) \to X([0])$ is the one induced by the inclusion $\{ n \} \injto [n]$.
\end{lemenum}
\end{lemma}
\begin{proof}
Both claims follow from the fact that algebra and module objects in an operad are defined as certain functors from $\Ass^\tensor$ and $\LM^\tensor$. In the language of non-symmetric operads, $\Ass^\tensor$ corresponds to $\bbDelta^\op$ and $\LM^\tensor$ corresponds to $\bbDelta^\op \times \Delta^1$, so we can conclude by \cref{rmk:non-symmetric-vs-symmetric-operads}. For more details, see \cite[Proposition~2.4.2.5]{HA} together with \cite[Proposition~4.1.2.10]{HA} for (i) and \cite[Proposition~4.2.2.9]{HA} for (ii).
\end{proof}

Unpacking \cref{rslt:monoids-and-actions-via-Delta}, we see that a monoid $M \in \Mon(\cat C)$ can be equivalently described by a simplicial diagram of the form
\begin{equation}\begin{tikzcd}
	\dots \quad M \times M \times M \arrow[r,shift right=6] \arrow[r,shift right=2] \arrow[r,shift left=2] \arrow[r,shift left=6] & M \times M \arrow[r,shift right=4] \arrow[r] \arrow[r,shift left=4] \arrow[l] \arrow[l,shift right=4] \arrow[l,shift left=4] & M \arrow[r,shift right=2] \arrow[r,shift left=2] \arrow[l,shift right=2] \arrow[l,shift left=2] & * \arrow[l].
\end{tikzcd}\end{equation}
Here the two maps $M \to *$ are the projections and the map $* \to M$ is the unit map. Moreover, the three maps $M \times M \to M$ are given by $(m, m') \mapsto m', mm', m$ and the maps $M \to M \times M$ are given by $m \mapsto (1, m), (m, 1)$. Similarly, an action of $M$ on an object $X \in \cat C$ induces a simplicial object
\begin{equation}\begin{tikzcd}
	\dots \quad M^3 \times X \arrow[r,shift right=6] \arrow[r,shift right=2] \arrow[r,shift left=2] \arrow[r,shift left=6] & M^2 \times X \arrow[r,shift right=4] \arrow[r] \arrow[r,shift left=4] \arrow[l] \arrow[l,shift right=4] \arrow[l,shift left=4] & M \times X \arrow[r,shift right=2] \arrow[r,shift left=2] \arrow[l,shift right=2] \arrow[l,shift left=2] & X \arrow[l].
\end{tikzcd}\end{equation}
Here the maps $M^n \times X \to M^{n-1} \times X$ are given by
\begin{align}
	(m_1, \dots, m_n, x) \mapsto\: & (m_2, \dots, m_n, x),\\
		& (m_1, \dots, m_i m_{i+1}, \dots, m_n, x),\\
		& (m_1, \dots, m_{n-1}, m_n x),
\end{align}
for $i = 1, \dots, n-1$. Moreover, the maps $M^{n-1} \times X \to M^n \times X$ are given by inserting the unit of $M$ at the $i$-th place, for $i = 0, \dots, n-1$.

Having defined monoids and their actions, we now come to the definition of the quotient of such an action. With the above description in terms of simplicial objects, the definition of quotients is straightforward:

\begin{definition}
Let $\cat C$ be a category with finite products and let $M$ be a monoid in $\cat C$ acting on some object $X \in \cat C$. We denote $X_\bullet \to M_\bullet$ the associated natural transformation of simplicial objects in $\cat C$ via \cref{rslt:actions-via-Delta}. We then define the \emph{quotient of $X$ by $M$} as
\begin{align}
	X/M \coloneqq \varinjlim_{[n]\in\bbDelta^\op} X_n = \varinjlim_{[n]\in\bbDelta^\op} (M^n \times X),
\end{align}
assuming it exists in $\cat C$.
\end{definition}

By construction the quotient $X/M$ receives a map $\pi\colon X \to X/M$ whose precompositions with the two maps $M \times X \to X$, $(m, x) \mapsto x, mx$ agree. In general quotients may not behave nicely (in particular they may not even exist), but the situation becomes better in the case that $\cat C$ is a category of sheaves (i.e.\ \enquote{stacks}, as everything is defined in terms of higher categories). More generally, we have the following results (recall the notion of groupoids and effective groupoids from \cite[Definitions~6.1.2.14,~6.1.2.7]{HTT}).

\begin{lemma} \label{rslt:quotients-are-groupoids}
Let $\cat C$ be a category with finite products, let $G$ be a group in $\cat C$ acting on some object $X \in \cat C$ and let $X_\bullet \to G_\bullet$ be the associated map of simplicial objects in $\cat C$ (via \cref{rslt:actions-via-Delta}). Then $X_\bullet$ is a groupoid object in $\cat C$.
\end{lemma}
\begin{proof}
Fix some $n \ge 0$ and a partition $[n] = S \union S'$ with $S \isect S' = \{i \}$ for some $0 \le i \le n$. By \cite[Proposition~6.1.2.6.(4$''$)]{HTT} we need to show that the diagram
\begin{equation}\begin{tikzcd}
	X([n]) \arrow[r] \arrow[d] & X(S) \arrow[d]\\
	X(S') \arrow[r] & X(\{i\})
\end{tikzcd}\end{equation}
is a pullback diagram in $\cat C$. This diagram has the shape
\begin{equation}\begin{tikzcd}
	G^n \times X \arrow[r] \arrow[d] & G^{\abs{S}-1} \times X \arrow[d]\\
	G^{\abs{S'}-1} \times X \arrow[r] & X
\end{tikzcd}\end{equation}
One can now inductively reduce to the case $n = 1$ using the above explicit description of the transition maps in $X_\bullet$ and the group property of $G$. We leave the details to the reader.
\end{proof}

\begin{corollary} \label{rslt:cech-nerve-of-quotient-cover}
Let $\cat C$ be a category with finite limits and assume that all groupoids in $\cat C$ are effective (e.g.\ $\cat C$ is a topos). Then for every group $G$ in $\cat C$ acting on an object $X \in \cat C$, the Čech nerve of the map $X \to X/G$ is given by
\begin{equation}\begin{tikzcd}
	\dots \arrow[r,shift right=3] \arrow[r,shift right=1] \arrow[r,shift left=1] \arrow[r,shift left=3] & G^2 \times X \arrow[r,shift right=2] \arrow[r] \arrow[r,shift left=2] & G \times X \arrow[r,shift right] \arrow[r,shift left] & X \arrow[r] & X/G
\end{tikzcd}\end{equation}
More precisely, if the $G$-action on $X$ corresponds to the map $X_\bullet \to G_\bullet$ of simplicial objects in $\cat C$, then $X_\bullet$ is the Čech nerve of $X \to X/G$.
\end{corollary}
\begin{proof}
This follows immediately from \cref{rslt:quotients-are-groupoids} and the definition of effective groupoids (see \cite[Definition~6.1.2.14]{HTT}).
\end{proof}

\begin{remark}
By applying \cref{rslt:quotients-are-groupoids} to the trivial action of a group $G$ on the final object $*$, we deduce that a group in $\cat C$ is equivalently a groupoid $G_\bullet$ in $\cat C$ with $G_0 = *$.
\end{remark}

\begin{proposition}\label{rslt:group-action-equals-map-to-*/G}
Suppose that $\cat C$ is a topos. Let $G$ be a group in $\cat C$ and let $X\to */G$ in $\cat C$. Then there exists a unique (up to contractible choice) $G$-action on $P\coloneqq *\times_{*/G}X$ such that the projection $P\to X$ induces an isomorphism $P/G \simeq X$.
\end{proposition}
\begin{proof}
See \cite[Propositions~3.8 and~3.16]{Nikolaus-Schreiber-Stevenson.2015} and Theorem~3.17 in \loccit{} for a stronger result. We give the proof for the convenience of the reader.

Let $G_\bullet \to */G$ be the Čech nerve of $*\to */G$ and define $P_\bullet \coloneqq G_\bullet \times_{*/G}X$. For every $n\ge0$ the map $P_n \isoto G_n\times P$ is an isomorphism, which shows that $P_\bullet\to G_\bullet$ exhibits an action of $G$ on $P = P_0$. Since $\blank\times_{*/G}X$ commutes with colimits (here it is used that $\cat C$ is a topos), we deduce that $P/G = \varinjlim_{[n]\in\bbDelta^{\op}} P_n \simeq X$, which finishes the proof of the existence claim. The uniqueness follows from \cref{rslt:cech-nerve-of-quotient-cover}. 
\end{proof}



\subsection{Cartesian left actions} \label{sec:alg.cartact}

Let $\cat V$ be a category with finite products and equip it with the cartesian symmetric monoidal structure from \cref{ex:cartesian-monoidal-structure}. Then one way to construct a $\cat V$-linear structure on some category $\cat C$ is to specify a functor $\cat V \to \cat C$ that preserves finite products. Indeed, this functor then upgrades to a symmetric monoidal functor of symmetric monoidal cartesian categories, i.e.\ to a map in $\CMon = \CAlg(\Cat^\times)$. By applying the forgetful functor $\LMod_{\cat C} \to \LMod_{\cat V}$ to $\cat C$ we arrive at the desired $\cat V$-action of $\cat C$.

The purpose of this subsection is to provide an axiomatic way of characterizing $\cat V$-linear actions that arise from the construction described in the previous paragraph. The main consequence of this is to show that such an action is uniquely determined by the map $\cat V \to \cat C$, which helps to show that different constructions of enrichment yield the same answer (see e.g.\ \cref{rslt:2-yoneda} below). The results in this section should be seen as an analog of the results in \cite[\S2.4.1]{HA} for $\Ass^\tensor$ and $\LM^\tensor$ in place of $\Comm^\tensor$. Let us start by discussing cartesian (non-symmetric) monoidal categories:

\begin{definition} \label{def:cartesian-monoidal-category}
A monoidal category $\cat C^\tensor \to \Ass^\tensor$ is called \emph{cartesian} if the tensor unit $\one \in \cat C$ is final and for all $X, Y \in \cat C$ the canonical maps $X = X \tensor \one \from X \tensor Y \to \one \tensor Y = Y$ induce an isomorphism $X \tensor Y = X \times Y$.
\end{definition}

\begin{definition} \label{def:weak-cartesian-structure}
Let $\cat O^\tensor$ be an operad, $\cat C^\tensor \to \cat O^\tensor$ an $\cat O$-monoidal category and $\cat D$ a category with finite products. A functor $F\colon \cat C^\tensor \to \cat D$ is called an \emph{$\cat O$-weak cartesian structure} if it is a lax cartesian structure (in the sense of \cref{ex:cartesian-monoidal-structure}) and sends cocartesian active morphisms over $\cat O^\tensor$ to isomorphisms.
\end{definition}

The following results should be seen as an analog of \cite[Proposition~2.4.1.6]{HA}. We will deduce afterwards that cartesian monoidal structures are uniquely determined (see \cref{rslt:cartesian-monoidal-structure-unique} below)

\begin{proposition} \label{rslt:functors-out-of-cartesian-monoidal-category}
Let $\cat C^\tensor \in \Op_{/\Ass^\tensor}$ be a cartesian monoidal category and let $\cat D$ be a category which has all finite products. Then restriction induces an equivalence of categories
\begin{align}
	\Fun^\times(\cat C^\tensor, \cat D) \isoto \Fun^\times(\cat C, \cat D).
\end{align}
Here the left-hand side denotes the category of weak cartesian structures and the right-hand side denotes the full subcategory of $\Fun(\cat C, \cat D)$ spanned by those functors that preserve finite products.
\end{proposition}
\begin{proof}
The argument is very similar to the proof of \cite[Proposition~2.4.1.6]{HA}, but some additional care is required to handle the differences between $\Comm^\tensor$ and $\Ass^\tensor$. In the following, we say that a map in $\alpha\colon \langle m \rangle \to \langle n \rangle$ in $\Ass^\tensor$ is \emph{order-preserving} if the underlying map $\langle m \rangle \to \langle n \rangle$ is monotonous (where we ignore all $j \in \langle m \rangle^\circ$ that get sent to $*$) and for all $i = 1, \dots, n$ the induced ordering on $\alpha^{-1}(i)$ is the canonical one. We say that a map in $\cat C^\tensor$ is order-preserving if it lies over an order-preserving map in $\Ass^\tensor$. With this terminology at hand, we now introduce the following subcategories of $\cat C^\tensor$:
\begin{itemize}
	\item $\cat C_1^\tensor \subseteq \cat C^\tensor$ is the subcategory with all objects, but where we only allow order-preserving active morphisms.

	\item $\cat C_2^\tensor \subseteq \cat C^\tensor$ is the subcategory with all objects, but where we only allow order-preserving morphisms.
\end{itemize}
The restriction functor $\Fun^\times(\cat C^\tensor, \cat D) \to \Fun^\times(\cat C, \cat D)$ factors as
\begin{align}
	\Fun^\times(\cat C^\tensor, \cat D) \xto{R_3} \Fun^\times(\cat C_2^\tensor, \cat D) \xto{R_2} \Fun^\times(\cat C_1^\tensor, \cat D) \xto{R_1} \Fun^\times(\cat C, \cat D),
\end{align}
where $\Fun^\times(\cat C_2^\tensor, \cat D)$ is defined similarly as in \cref{def:weak-cartesian-structure} and $\Fun^\times(\cat C_1^\tensor, \cat D)$ is the full subcategory of $\Fun(\cat C_1^\tensor, \cat D)$ consisting of those functors $F$ which restrict to a finite-product preserving functor $\cat C \to \cat D$ and send cocartesian active morphisms over $\Ass^\tensor$ to isomorphisms. In the following we will show that $R_1$, $R_2$ and $R_3$ are equivalences by constructing their inverses via Kan extensions.

Let us first show that $R_1$ is an equivalence. This is very similar to the first step in the proof of \cite[Proposition~2.4.1.6]{HA}, so we will be brief. Namely, the inverse of $R_1$ is given by right Kan extension, which follows easily by observing that for every $X_\bullet \in \cat C_1^\tensor$, the category $\cat C_{X_\bullet/}$ has an initial object given by the cocartesian map $X_\bullet \to \prod_i X_i$ lying over $\alpha_n$.

We now show that $R_2$ is an equivalence. We claim that its inverse is given by left Kan extension along the inclusion $\cat C_1^\tensor \to \cat C_2^\tensor$. Fix a functor $F\in \Fun^\times(\cat C_1^\tensor, \cat D)$ and an object $X_\bullet \in \cat C^\tensor$ lying over $\langle n \rangle \in \Ass^\tensor$. Let us denote by $F'\colon \cat C_2^\tensor \to \cat D$ the left Kan extension of $F$ (assuming it exists) and let $\cat I \coloneqq \cat C_1^\tensor \times_{\cat C_2^\tensor} (\cat C_2^\tensor)_{/X_\bullet}$. Then $F'(X_\bullet) = \varinjlim_{Y_\bullet \in \cat I} F(Y_\bullet)$, so in order to prove that the left Kan extension of $F$ exists and defines an inverse of $R_2$, we need to verify that the natural map
\begin{align}
	F(X_\bullet) \isoto \varinjlim_{Y_\bullet \in \cat I} F(Y_\bullet)
\end{align}
is an isomorphism. Let $\cat I' \subseteq \cat I$ be the full subcategory spanned by those maps $(Y_j)_{1\le j \le m} \to (X_i)_{1\le i \le n}$ lying over $\alpha\colon \langle m \rangle \to \langle n \rangle$ such that for each $j = 1, \dots, m$ for which $\alpha(j) = *$, the object $Y_j \in \cat C$ is final. We claim that the inclusion $\cat I' \injto \cat I$ is right cofinal. Fix an object $f\colon Y_\bullet \to X_\bullet$ in $\cat I$ lying over some map $\alpha\colon \langle m \rangle \to \langle n \rangle$ in $\Ass^\tensor$; we need to see that the category $\cat I'_{f/}$ is weakly contractible (see \cite[\href{https://kerodon.net/tag/02NY}{Tag 02NY}]{kerodon}). But this category has an initial element given by the map $Y_\bullet \to Y'_\bullet$ lying over the identity $\langle m \rangle \to \langle m \rangle$, where $Y'_j = Y_j$ if $\alpha(j) \ne *$ and $Y'_j = \one$ otherwise. This proves that $\cat I' \injto \cat I$ is right cofinal, so we are reduced to showing that the natural map
\begin{align}
	F(X_\bullet) \isoto \varinjlim_{Y_\bullet \in \cat I'} F(Y_\bullet)
\end{align}
is an isomorphism. Now let $\cat I'' \subset \cat I'$ be the full subcategory spanned by the maps $Y_\bullet \to X_\bullet$ that lie in $\cat C_1^\tensor$ (in other words, $\cat I'' = (\cat C_1^\tensor)_{/X_\bullet}$). We claim that $F\colon \cat I' \to \cat D$ is the left Kan extension of its restriction to $\cat I''$. To see this, denote this Kan extension by $F''$ and fix some $f\colon Y_\bullet \to X_\bullet$ in $\cat I'$. Suppose $f$ lies over $\alpha\colon \langle m \rangle \to \langle n \rangle$ in $\Ass^\tensor$ and let $S \subseteq \langle m \rangle^\circ$ be the subset of those indices $j$ for which $\alpha(j) = *$. Then $\cat I''_{/f}$ has a final element $Y'_\bullet$ lying over $\langle m \rangle \setminus S$, where $Y'_j = Y_j$ for all $j \in \langle m \rangle^\circ \setminus S$. Therefore $F''(Y_\bullet) = F(Y'_\bullet)$ and we need to see that the natural map $F(Y'_\bullet) = F''(Y_\bullet) \isoto F(Y_\bullet)$ is an isomorphism. But this is clear, as $F(Y'_\bullet)$ and $F(Y_\bullet)$ differ only by a product of copies of $\one$. By transitivity of Kan extensions, we deduce that
\begin{align}
	\varinjlim_{Y_\bullet \in \cat I'} F(Y_\bullet) = \varinjlim_{Y_\bullet \in \cat I''} F(Y_\bullet) = F(X_\bullet),
\end{align}
where the second identity is implied by the fact that $\cat I''$ has a final element given by the identity on $X_\bullet$. This finishes the proof that $R_2$ is an equivalence.

It remains to show that $R_3$ is an equivalence. We claim that its inverse is given by right Kan extension. Let $G\in \Fun^\times(\cat C_2^\tensor, \cat D)$ and let $G'\colon \cat C^\tensor \to \cat D$ be its right Kan extension (assuming it exists). Fix an object $X_\bullet \in \cat C^\tensor$ lying over $\langle n \rangle \in \Ass^\tensor$ and let $\cat J \coloneqq \cat C_2^\tensor \times_{\cat C^\tensor} (\cat C^\tensor)_{/X_\bullet}$. Similarly to the previous step, the proof that $R_3$ is an equivalence reduces to showing that the natural map
\begin{align}
	\varprojlim_{Y_\bullet \in \cat J} G(Y_\bullet) = G'(X_\bullet) \isoto G(X_\bullet)
\end{align}
is an isomorphism. Let $\cat J' \subseteq \cat J$ be the full subcategory spanned by the inert morphisms $X_\bullet \to Y_\bullet$ in $\cat J$; in other words, $Y_\bullet$ is obtained from $X_\bullet$ by potentially dropping some $X_i$'s and reordering the rest. The map $\cat J' \injto \cat J$ is left cofinal: By \cite[\href{https://kerodon.net/tag/02NY}{Tag 02NY}]{kerodon} this follows from the observation that for every object $f\colon \cat X_\bullet \to Y_\bullet$ in $\cat J$ the category $\cat J'_{/f}$ has a final object, dictated by the underlying map $\alpha\colon \langle n \rangle \to \langle m \rangle$ in $\Ass^\tensor$. By the cofinality we are reduced to showing that the natural map
\begin{align}
	\varprojlim_{Y_\bullet \in \cat J'} G(Y_\bullet) \isoto G(X_\bullet)
\end{align}
is an isomorphism. Let $\cat J'' \subseteq \cat J'$ be the full subcategory spanned by those (inert) maps $X_\bullet \to Y_\bullet$ where $Y_\bullet$ lies over $\langle 1 \rangle$. We claim that $G\colon \cat J' \to \cat D$ is the right Kan extension of its restriction to $\cat J''$. Indeed, this right Kan extension maps $Y_\bullet \mapsto \prod_i G(Y_i)$ (because $J''_{Y_\bullet/} = \bigdunion_i \{ Y_i \}$), which easily implies the claim. By a similar computation we have
\begin{align}
	\varprojlim_{Y_\bullet \in \cat J'} G(Y_\bullet) = \varprojlim_{Y \in \cat J''} G(Y) = \prod_i G(X_i) = G(X_\bullet),
\end{align}
as desired. This proves that $R_3$ is an equivalence.
\end{proof}

\begin{corollary} \label{rslt:cartesian-monoidal-structure-unique}
Let $\cat C^\tensor \to \Ass^\tensor$ be a cartesian monoidal category.
\begin{corenum}
	\item $\cat C^\tensor$ extends uniquely to a cartesian \emph{symmetric} monoidal category. More precisely, there is an isomorphism $\cat C^\tensor = \cat C^\times \times_{\Comm^\tensor} \Ass^\tensor$ of monoidal categories which is uniquely determined by the requirement that it induces the identity on $\cat C$.

	\item If $\cat C'^\tensor \to \Ass^\tensor$ is another cartesian monoidal category then restriction induces an equivalence
	\begin{align}
		\Fun^\tensor(\cat C^\tensor, \cat C'^\tensor) \isoto \Fun^\times(\cat C, \cat C'),
	\end{align}
	where $\Fun^\tensor \subseteq \Fun_{\Ass^\tensor}$ denotes the full subcategory of monoidal functors.
\end{corenum}
\end{corollary}
\begin{proof}
By \cref{rslt:functors-out-of-cartesian-monoidal-category} and \cite[Proposition~2.4.1.7]{HA} we deduce
\begin{align}
	\Fun^\times(\cat C, \cat C') = \Fun^\times(\cat C^\tensor, \cat C') \subseteq \Fun^\lax(\cat C^\tensor, \cat C') = \Alg(\cat C^\tensor, \cat C'^\times) = \Alg_{\Ass^\tensor}(\cat C^\tensor, \cat C'^\times \times_{\Comm^\tensor} \Ass^\tensor),
\end{align}
where $\Fun^\lax$ denotes the category of lax cartesian structures. By analyzing the essential image of the above fully faithful embedding, we obtain an equivalence
\begin{align}
	\Fun^\times(\cat C, \cat C') = \Fun^\tensor(\cat C^\tensor, \cat C'^\times \times_{\Comm^\tensor} \Ass^\tensor).
\end{align}
Claims (i) and (ii) follow easily from this. Taking $\cat C' = \cat C$ we see that the identity $\cat C \to \cat C$ induces a monoidal functor $\cat C^\tensor \to \cat C'^\times \times_{\Comm^\tensor} \Ass^\tensor$ which acts as the identity on underlying categories and is therefore an isomorphism; this proves (i). Now (ii) follows from the above identity using (i) applied to $\cat C'$.
\end{proof}

To summarize, we have seen that if a monoidal structure on some category $\cat C$ is given by the product in $\cat C$, then this property uniquely determines all the data that come with the monoidal structure, i.e.\ all (higher) associators and unit relations. In fact, one can even uniquely add (higher) commutators in order to make the monoidal structure into a \emph{symmetric} monoidal structure.

We now apply similar ideas to the setting of cartesian left actions. We start with the definition and then continue with the analogs of \cref{rslt:functors-out-of-cartesian-monoidal-category,rslt:cartesian-monoidal-structure-unique}.

\begin{definition}
A \emph{cartesian left action} is an $\LM$-monoidal category $\cat O^\tensor \to \LM^\tensor$ exhibiting the left action of the monoidal category $\cat V^\tensor \coloneqq \cat O^\tensor \times_{\LM^\tensor} \Ass^\tensor$ on the category $\cat C \coloneqq \cat O^\tensor_{\mathfrak m}$ satisfying the following properties:
\begin{enumerate}[(i)]
	\item $\cat V^\tensor$ is a cartesian monoidal category.
	\item $\cat C$ has a final object $* \in \cat C$. Moreover, for all $V \in \cat V$ and $X \in \cat C$ the maps $V \tensor * \from V \tensor X \to \one_{\cat V} \tensor X = X$ exhibit an isomorphism $V \tensor X = (V \tensor *) \times X$.
\end{enumerate}
\end{definition}

\begin{proposition} \label{rslt:functors-out-of-cartesian-left-action}
Let $\cat O^\tensor \to \LM^\tensor$ be a cartesian left action of $\cat V^\tensor$ on $\cat C$ and let $\cat D$ be a category with finite products. Then restriction induces an equivalence of categories
\begin{align}
	\Fun^\times(\cat O^\tensor, \cat D) \isoto \Fun'(\cat C, \cat D),
\end{align}
where $\Fun'(\cat C, \cat D) \subseteq \Fun(\cat C, \cat D)$ denotes the full subcategory spanned by those functors $F\colon \cat C \to \cat D$ which satisfy $F(V \tensor X) = F(V \tensor *) \times F(X)$ for all $V \in \cat V$ and $X \in \cat C$.
\end{proposition}
\begin{proof}
We argue similarly as in \cref{rslt:functors-out-of-cartesian-monoidal-category}. Recall the notion of order-preserving maps from the beginning of the proof of \cref{rslt:functors-out-of-cartesian-monoidal-category}, which we can apply to $\cat O^\tensor$ via the projection $\LM^\tensor \to \Ass^\tensor$. Note that every object in $\LM^\tensor$ is of the form $\mathfrak x_1 \dsum \dots \dsum \mathfrak x_n$, where $\mathfrak x_i \in \{ \mathfrak a, \mathfrak m \}$ for all $i$. Consequently, each object in $\cat O^\tensor$ is of the form $X_1 \dsum \dots \dsum X_n$, where each $X_i$ lies either in $\cat V = \cat O^\tensor_{\mathfrak a}$ or in $\cat C = \cat O^\tensor_{\mathfrak m}$. Fix an order-preserving map $\alpha\colon (\mathfrak y_j)_{1 \le j \le m} \to (\mathfrak x_i)_{1 \le i \le n}$ in $\LM^\tensor$. We introduce the following terminology:
\begin{itemize}
	\item We say that $\alpha$ \emph{respects blocks} if for all $j \in \{ 1, \dots, m \}$ such that $\mathfrak y_j = \mathfrak m$ and $\alpha(j) = *$, either $\alpha(j') = *$ for all $j' \ge j$ or $\alpha(j') = *$ for all $j'$ with $j_0 < j' \le j$, where $j_0 < j$ is the largest index for which $\mathfrak y_{j_0} = \mathfrak m$.

	\item Suppose that $\alpha$ respects blocks. We say that $\alpha$ is \emph{active in blocks} if for every $j \in \{ 1, \dots, m \}$ such that $\mathfrak y_j = \mathfrak a$ and $\alpha(j) = *$, we have $\alpha(j_1) = *$, where $j_1 > j$ is the smallest index for which $\mathfrak y_{j_1} = \mathfrak m$.
\end{itemize}
We say that a map in $\cat O^\tensor$ respects blocks or is active in blocks if this is true for the underlying map in $\LM^\tensor$. Note that these conditions are stable under composition. We can now define the following subcategories of $\cat O^\tensor$:
\begin{itemize}
 	\item $\cat O_1^\tensor \subset \cat O^\tensor$ is the subcategory consisting of the objects $(X_i)_{1 \le i \le n}$ such that $X_n \in \cat C$ (in particular we require $n > 0$) and the order-preserving morphisms that respect blocks and are active in blocks.

 	\item $\cat O_2^\tensor \subset \cat O^\tensor$ is the subcategory with all objects, where the morphisms are the order-preserving maps that respect blocks.
\end{itemize}
Then the restriction $\Fun^\times(\cat O^\tensor, \cat D) \to \Fun'(\cat C, \cat D)$ factors as
\begin{align}
	\Fun^\times(\cat O^\tensor, \cat D) \xto{R_3} \Fun^\times(\cat O_2^\tensor, \cat D) \xto{R_2} \Fun'(\cat O_1^\tensor, \cat D) \xto{R_1} \Fun'(\cat C, \cat D),
\end{align}
where $\Fun^\times(\cat O_2^\tensor, \cat D) \subseteq \Fun(\cat O_2^\tensor, \cat D)$ is a full subcategory defined analogous to the definition of weak cartesian structures and $\Fun'(\cat O_1^\tensor, \cat D) \subseteq \Fun(\cat O_1^\tensor, \cat D)$ is the full subcategory consisting of the functors that restrict to a functor in $\Fun'(\cat C, \cat D)$ and that send active cocartesian morphisms to isomorphisms. In the following we show that $R_1$, $R_2$ and $R_3$ are equivalences by constructing their inverses via Kan extensions. The arguments are very similar to the proof of \cref{rslt:functors-out-of-cartesian-monoidal-category}, so we will focus on the steps that are different and only sketch the rest.

The inverse of $R_1$ is given by right Kan extension. Indeed, fix some $(X_i)_{1\le i \le n} \in \cat O_1^\tensor$ and consider $\cat K \coloneqq \cat C_{X_\bullet/}$. Let $S_1, \dots, S_N \subset \{ 1, \dots, n \}$ be the \emph{blocks} of $X_\bullet$, i.e.\ they form a disjoint partition of $\{ 1, \dots, n \}$ into subsets of consecutive indices such that each subset has exactly one index $i$ for which $X_i \in \cat C$ and this index is maximal in the subset. For each $I \in \{ 1, \dots, N \}$ there is an object $X_{S_I} \in \cat C$ and a cocartesian edge $X_\bullet \to X_{S_I}$ in $\cat O_1^\tensor$ such that the underlying map in $\Fin_*$ is the projection $\langle n \rangle \surjto S_I \dunion \{ * \}$. Explicitly, if $S_I = \{ i_1, \dots, i_k \}$ then $X_{S_I} = (X_{i_1} \tensor \dots \tensor X_{i_{k-1}}) \tensor X_{i_k}$. Now $\{ X_{S_1}, \dots, X_{S_N} \} \injto \cat K$ is left cofinal, hence the limit of a functor over $\cat K$ is the same as the product of the values of that functor on the $X_{S_I}$; this easily implies that right Kan extension is indeed the inverse of $R_1$.

The inverse of $R_2$ is given by left Kan extension. To show this, fix $(X_i)_{1\le i \le n}$ in $\cat O^\tensor$ and define $\cat I \coloneqq \cat O_1^\tensor \times_{\cat O_2^\tensor} (\cat O_2^\tensor)_{/X_\bullet}$. We need to understand colimits of functors $\cat I \to \cat D$ (see the proof of \cref{rslt:functors-out-of-cartesian-monoidal-category}). As in the proof of \cref{rslt:functors-out-of-cartesian-monoidal-category} we consider the full subcategory $\cat I' \subseteq \cat I$ spanned by the maps $(Y_j)_{1\le j \le m} \to X_\bullet$ over $\alpha\colon \langle m \rangle \to \langle n \rangle$ in $\Fin_*$ such that for all $j \in \{ 1, \dots, m \}$ for which $\alpha(j) = *$, the object $Y_j$ is final. The inclusion $\cat I' \injto \cat I$ is right cofinal, so we can restrict to computing colimits over $\cat I'$. By a left Kan extension argument as in the proof of \cref{rslt:functors-out-of-cartesian-monoidal-category} we can further reduce to the full subcategory $\cat I'' \subset \cat I'$ consisting of the maps $Y_\bullet \to X_\bullet$ which are active in blocks. If $n > 0$ and $X_n \in \cat C$ then $\cat I''$ has a final object given by the identity on $X_\bullet$, which provides the desired colimit computation. Otherwise $\cat I''$ has a final object given by $(X'_i)_{1 \le i \le n+1}$ with $X'_i = X_i$ for $1 \le i \le n$ and $X'_{n+1} = * \in \cat C$, where the map $X'_\bullet \to X_\bullet$ is the unique inert map. This forces the colimit over $\cat I$ to be the value on $X'_\bullet$, which is exactly what is desired (in the easiest example $X_\bullet = V \in \cat V$, the colimit sends this to the value on $V \tensor * \in \cat C$).

The inverse of $R_3$ is given by right Kan extension. The proof is more or less the same as in \cref{rslt:functors-out-of-cartesian-monoidal-category}, so we do not repeat it here. The crucial observation is that for every $X_\bullet \in \cat O^\tensor$, all the inert maps from $X_\bullet$ to some object in $\cat O^\tensor_{\langle 1 \rangle}$ lie in $\cat O_2^\tensor$.
\end{proof}

With \cref{rslt:functors-out-of-cartesian-left-action} we can now show that cartesian left actions are uniquely determined by the underlying finite-product preserving functor $\cat V \to \cat C$ (see \cref{rslt:uniqueness-of-cartesian-left-actions} below). In fact, the following definition provides a canonical (and a posteriori unique) cartesian left action of this form.

\begin{definition} \label{def:induced-cartesian-left-action}
Let $\cat V$ and $\cat C$ be categories which have all finite products and let $F\colon \cat V \to \cat C$ be a functor that preserves finite products. We denote by
\begin{align}
	(\cat V, \cat C, F)^\times \in \LMod_{\cat V}
\end{align}
the cartesian left action of $\cat V$ on $\cat C$ \emph{induced by $F$}. It is constructed as follows: $F$ induces a map $\cat V^\times \to \cat C^\times$ of cartesian symmetric monoidal categories and hence of algebras in $\Cat^\times$; then $(\cat V, \cat C, F)^\times$ is obtained by applying the forgetful functor $\LMod_{\cat C} \to \LMod_{\cat V}$ to $\cat C$.
\end{definition}

\begin{corollary} \label{rslt:uniqueness-of-cartesian-left-actions}
Let $\cat O^\tensor \to \LM^\tensor$ be a cartesian left action of $\cat V^\tensor$ on $\cat C$ and assume that $\cat C$ has finite products. Then there is a unique isomorphism
\begin{align}
	\cat O^\tensor = (\cat V, \cat C, F)^\times
\end{align}
of $\LM$-monoidal categories that acts as the identity on $\cat C$. Here $F\colon \cat V \to \cat C$ is the functor $V \mapsto V \tensor *$, where $* \in \cat C$ is a final object.
\end{corollary}
\begin{proof}
From \cref{rslt:functors-out-of-cartesian-left-action} and \cite[Proposition~2.4.1.7]{HA} we deduce that the identity on $\cat C$ induces a map of operads $\cat O^\tensor \to \cat C^\times$ and hence a map of $\LM$-monoidal categories $\cat O^\tensor \to \cat C'^\tensor \coloneqq \cat C^\times \times_{\Comm^\tensor} \LM^\tensor$. Recall that the functor $\Op_{/\LM^\tensor} \to \Op_{\Ass^\tensor}$ is a cartesian fibration (this follows e.g.\ from the description in \cite[Proposition~3.33]{Heine.2023}), so we can pull $\cat C'^\tensor$ along $\cat V^\times \to \cat C^\times$ to obtain the $\LM$-monoidal category $\cat C''^\tensor$ exhibiting the $\cat V$-action on $\cat C$ induced by the map $\cat V \to \cat C$. Then the map $\cat O^\tensor \to \cat C'^\tensor$ induces a unique map $\cat O^\tensor \to \cat C''^\tensor$ of $\cat V$-linear categories, which acts as the identity on the underlying category (which is $\cat C$ in both cases) and thus is an isomorphism of operads. This shows that any two cartesian left actions of $\cat V$ on $\cat C$ inducing the same functor $F$ are uniquely isomorphic, which in particular implies the claim.
\end{proof}

\begin{remark}
It should be possible to refine \cref{rslt:uniqueness-of-cartesian-left-actions} in order to obtain an equivalence of the category of cartesian left actions with a certain category of product-preserving functors (in the spirit of \cite[Corollary~2.4.1.9]{HA}. We do not pursue this idea further, as we do not need it.
\end{remark}

\section{Enriched categories} \label{sec:enr}

In the context of the category of kernels, one of the main objects of study in this paper, we will make heavy use of the language of enriched categories and in particular 2-categories. The basic language of enriched categories has been established in the literature by many authors, but we find it a bit cumbersome to digest for non-experts. In the following we provide a (hopefully) easy-to-follow introduction to enriched category theory by summarizing the results from the literature and adding some of our own results that are needed in this paper.

\subsection{Basic definitions}

Fix a monoidal category $\cat V$. A $\cat V$-enriched category $\cat C$ consists of the following data:
\begin{enumerate}[(a)]
	\item a set of objects $S$,
	\item for every pair of objects $X, Y \in S$ an object $\enrHom[\cat V](X, Y) \in \cat V$,
	\item for every object $X \in S$ a map $\one_{\cat V} \to \enrHom[\cat V](X, X)$, called the \emph{identity map},
	\item for every triple of objects $X, Y, Z \in S$ a map $\enrHom[\cat V](Y, Z) \tensor \enrHom[\cat V](X, Y) \to \enrHom[\cat V](X, Z)$, called the \emph{composition},
\end{enumerate}
together with higher compatibilities. It is a non-trivial question how to encode these data appropriately, and several models for enriched categories have appeared in the literature. Up to easy equivalence, there are two main approaches to enriched categories: the one using weak actions by Lurie (see \cite[Definition~4.2.1.28]{HA}) and the one based on quivers by Gepner--Haugseng and Hinich (see \cite{Gepner-Haugseng.2015,Hinich.2021}). These two approaches have been shown to be equivalent in a strong sense by \cite{Heine.2023}, so we are usually free to choose any model we like for a task at hand.

In the following we recall both approaches and how they are related. We start with the quiver approach, as this is probably the more intuitive one. Suppose we are given the monoidal category $\cat V$ and an anima $S$ of objects (working with anima is more natural than working with sets in the higher categorical framework, see also the discussion below \cref{def:quiver-fibration}). Then part (b) of the above data can be encoded by a functor
\begin{align}
	F\colon S \times S \to \cat V, \qquad (X, Y) \mapsto F(X, Y) \eqqcolon \enrHom[\cat V](X, Y). \label{eq:functor-F-defining-an-enriched-category}
\end{align}
It turns out that parts (c) and (d) of the above data can be encoded by some algebra structure on $F$ for a certain monoidal structure on $\Fun(S \times S, \cat V)$:

\begin{definition} \label{def:quiver-category}
Let $\cat V$ be an operad over $\Ass^\tensor$ and $S$ an anima. We denote by
\begin{align}
	\Quiv_S(\cat V)^\tensor \coloneqq \Fun(\Ass_S, \cat V)^\tensor \to \Ass^\tensor
\end{align}
the operad constructed in \cite[Notation~4.17]{Heine.2023} (their $X$ is our $S$ and their $S$ is $*$ in our case). It is defined to be the Day convolution operad of functors on a certain generalized operad\footnote{See \cite[Definition~2.3.2.1]{HA} for the definition of \emph{generalized} operads, which uses a slightly weaker condition on the map $\cat O^\tensor \to \Fin_*$. By \cite[Proposition~2.3.2.5]{HA} the main difference to operads is that a generalized operad $\cat O^\tensor$ may have a non-contractible fiber $\cat O^\tensor_{\langle 0 \rangle}$.} $\Ass_S^\tensor \to \Ass^\tensor$ with fibers $(\Ass_S^\tensor)_{\langle n\rangle} = S^{\times(n+1)}$ for all $n \ge 0$ (see \cite[Notation~4.1]{Heine.2023} for the definition of $\Ass_S^\tensor$ as a non-symmetric operad).
\end{definition}

Let us explain $\Quiv_S(\cat V)^\tensor$ in more concrete terms. Day convolution is a general way of equipping the functor category of two operads with an operad structure itself, see \cite[\S2.2.6]{HA} for more explanations. In the setting of \cref{def:quiver-category} we need a version of Day convolution for \emph{generalized} operads, which is developed in \cite[\S11]{Heine.2023}. In general, $\Quiv_S(\cat V)^\tensor$ is not a monoidal category, but it is so if $\cat V$ is a monoidal category and admits enough colimits (depending on $S$). More concretely, assume for simplicity that $S$ is a set (i.e.\ $\Hom_S(X, Y) \in \{ \emptyset, * \}$ for all $X, Y \in S$); then by \cite[Remark~4.28]{Heine.2023} we can describe $\Quiv_S(\cat V)^\tensor$ as follows:
\begin{itemize}
	\item The underlying category of $\Quiv_S(\cat V)^\tensor$ is
	\begin{align}
		\Quiv_S(\cat V) = \Fun(S \times S, \cat V).
	\end{align}
	This follows from the fact that the underlying category of $\Ass_S^\tensor$ is $S \times S$.

	\item The tensor unit of $\Quiv_S(\cat V)$ is the functor
	\begin{align}
		\one\colon (X, Y) \mapsto \begin{cases}
			\one_{\cat V}, & \text{if $\Hom(X, Y) = *$},\\
			\emptyset, & \text{if $\Hom(X, Y) = \emptyset$}.
		\end{cases}
	\end{align}
	Here $\emptyset \in \cat V$ denotes an initial object (assuming it exists).

	\item Given two functors $F, G\colon S \times S \to \cat V$, their tensor product is given by
	\begin{align}
		(F \circledast G)(X, Y) = \bigdunion_{Z \in S} F(Z, Y) \tensor G(X, Z),
	\end{align}
	where the right-hand side denotes the coproduct in $\cat V$ (assuming it exists).
\end{itemize}
From this description it follows that for the functor $F$ from \eqref{eq:functor-F-defining-an-enriched-category} the data (c) and (d) from the beginning of the section can be encoded as natural transformations $\one \to F$ and $F \circledast F \to F$, respectively. This motivates the following definition:

\begin{definition} \label{def:enriched-precategory}
Let $\cat V$ be an operad over $\Ass^\tensor$ and $S$ an anima. A \emph{$\cat V$-enriched precategory on $S$} is an algebra in $\Quiv_S(\cat V)$ (see \cite[Notation~4.40]{Heine.2023}).
\end{definition}

\begin{remarks} \label{rmk:comparison-of-precategories-to-Hinich-and-GH}
We compare \cref{def:enriched-precategory} to the definitions of Gepner--Haugseng and Hinich:
\begin{remarksenum}
	\item \label{rmk:comparison-of-precategories-to-GH} It follows easily from the universal property of Day convolution that an algebra object in $\Quiv_S(\cat V)$ is the same as a map of generalized operads $\Ass_S^\tensor \to \cat V^\tensor$ over $\Ass^\tensor$ (cf. \cite[Example~2.2.6.3]{HA}, the same argument applies to generalized operads). Therefore enriched precategories as defined above coincide with categorical algebras as defined by Gepner--Haugseng in \cite[Remark 2.4.7]{Gepner-Haugseng.2015} (but see \cite[Remark 4.5]{Heine.2023}).

	\item Hinich \cite{Hinich.2021} constructs a quiver operad and defines enriched precategories as algebras in that operad. It was shown by \cite{MacPherson.2020} that Hinich's quiver operad is equivalent to the one defined above (see also \cite[Remark 4.20]{Heine.2023}), so that Hinich's definition of enriched precategories coincides with \cref{def:enriched-precategory}.
\end{remarksenum}
\end{remarks}

Let us make the definition of enriched precategories functorial in the anima $S$ and in the enriching operad $\cat V$:

\begin{definition}
\begin{defenum}
	\item \label{def:quiver-fibration} Let
	\begin{align}
		\Quiv^\tensor \to \Ani \times \Op_{/\Ass^\tensor} \times \Ass^\tensor
	\end{align}
	be the map of generalized operads constructed in \cite[Notation~4.37]{Heine.2023}. The fiber of $\Quiv^\tensor$ over an anima $S \in \Ani$ and a non-symmetric operad $\cat V \in \Op_{/\Ass^\tensor}$ is the operad $\Quiv_S(\cat V)^\tensor$.

	\item \label{def:category-of-enriched-precategories} Let
	\begin{align}
		\PreEnr \coloneqq \Alg^{\Ani \times \Op_{/\Ass^\tensor}}(\Quiv^\tensor) \to \Ani \times \Op_{/\Ass^\tensor}
	\end{align}
	be the functor of categories constructed in \cite[Notation~4.40]{Heine.2023}. The notation $\Alg^{(\blank)}$ indicates that we take $\Alg$ fiberwise (see \cite[Notation~2.56]{Heine.2023}). In particular, the fiber of $\PreEnr$ over an anima $S \in \Ani$ and a non-symmetric operad $\cat V \in \Op_{/\Ass^\tensor}$ is the category of $\cat V$-enriched precategories on $S$. We call $\PreEnr$ the \emph{category of enriched precategories}.
\end{defenum}
\end{definition}

The name \emph{pre}-category indicates that something is still missing. Indeed, the definition of $\cat V$-enriched precategory has some redundancy in its data: The anima $S$ of objects should coincide with the space of equivalences in the enriched category, but in the above definition of enriched precategory this is not necessarily the case. For example, $S$ could be a set, while an enriched precategory on $S$ may introduce many new isomorphisms between objects in $S$. This behavior results in undesired functoriality behavior (e.g.\ fully faithful and essentially surjective functors of enriched precategories need not be equivalences). This can be fixed by defining the anima of equivalences $\iota \cat C$ of a $\cat V$-enriched precategory $\cat C$ on $S$ and then defining $\cat C$ to be a \emph{$\cat V$-enriched category} if the natural map $S \to \iota \cat C$ is an isomorphism of anima (cf. \cite[Definition~5.2.2]{Gepner-Haugseng.2015} and \cite[Definition~6.4]{Heine.2023}). We will now introduce a second model for enriched categories, due to Lurie, where this redundancy does not occur and we will use this model to pass from enriched precategories to enriched categories.

\begin{definition}
Let $\cat O^\tensor \to \Ass^\tensor$ be a map of operads. For all $X_1, \dots, X_n, Y \in \cat O$ we denote
\begin{align}
	\Mul_{\cat O}(X_1, \dots, X_n; Y) \subseteq \Hom_{\cat O^\tensor}(X_1 \dsum \dots \dsum X_n, Y)
\end{align}
the full subanima spanned by those maps which lie over the morphism $\alpha\colon \langle n \rangle \to \langle 1 \rangle$ in $\Ass^\tensor$ which is given by the natural ordering on $\langle n \rangle^\circ$.
\end{definition}

\begin{definition}
Let $\cat O^\tensor \to \LM^\tensor$ be a map of operads and denote $\cat V^\tensor \coloneqq \cat O^{\tensor} \times_{\LM^\tensor} \Ass^\tensor$ and $\cat C := \cat O_{\mathfrak m}$, so that $\cat V^\tensor$ is an operad over $\Ass^\tensor$ and $\cat C$ is a category.
\begin{defenum}
	\item \label{def:morphism-objects-for-enrichment} A \emph{morphism object} of $X, Y \in \cat C$ for $\cat O^\tensor$ is an object $\enrHom[\cat V](X, Y) \in \cat V$ together with a map in $\Mul_{\cat O}(\enrHom[\cat V](X, Y), X; Y)$ such that for all $V_1, \dots, V_n \in \cat V$ the induced map
	\begin{align}
		\Mul_{\cat V}(V_1, \dots, V_n; \enrHom[\cat V](X, Y)) \isoto \Mul_{\cat O}(V_1, \dots, V_n, X; Y)
	\end{align}
	is an isomorphism.\footnote{We require this isomorphism also for $n = 0$, where it reads $\Mul_{\cat V}(\one_{\cat V}; \enrHom[\cat V](X, Y)) \isoto \Mul_{\cat O}(X, Y)$; here $\one_{\cat V}$ denotes the unique object in $\cat V^\tensor_{\langle 0 \rangle}$.}

	\item \label{def:Lurie-enriched-category} We say that \emph{$\cat O^\tensor$ exhibits $\cat C$ as enriched over $\cat V$} if for all $X, Y \in \cat C$ there is a morphism object $\enrHom[\cat V](X, Y)$.
\end{defenum}
\end{definition}

\begin{remark} \label{rmk:comparison-of-definitions-for-Lurie-enriched-cat}
One checks that \cref{def:Lurie-enriched-category} is equivalent to \cite[Definition~3.122]{Heine.2023}. Indeed, by \cite[Proposition~3.33]{Heine.2023} and \cref{rmk:non-symmetric-vs-symmetric-operads} the notion of \enquote{weakly left tensored $\infty$-category} used in \cite{Heine.2023} coincides with a map of operads $\cat O^\tensor \to \LM^\tensor$ in our language, which immediately implies the equivalence of the definitions of enriched categories. In case the operad $\cat V^\tensor$ is a monoidal category, this definition of enriched categories also coincides with the one given in \cite[Definition~4.2.1.28]{HA} (see \cite[Lemma~3.123]{Heine.2023}).
\end{remark}

To get a better intuition for \cref{def:Lurie-enriched-category}, suppose that the operad $\cat O^\tensor \to \LM^\tensor$ exhibits $\cat C$ as a $\cat V$-linear category. Then for all $X, Y \in \cat C$ the morphism object $\enrHom[\cat V](X, Y) \in \cat V$ must satisfy
\begin{align}
	\Hom_{\cat V}(V, \enrHom[\cat V](X, Y)) = \Hom_{\cat C}(V \tensor X, Y)
\end{align}
for all $V \in \cat V$. Hence we see that $\enrHom[\cat V](X, \blank)\colon \cat C \to \cat V$ is a right adjoint to the functor $\blank \tensor X\colon \cat V \to \cat C$. As a special case, if $\cat C = \cat V$ with the natural $\cat V$-linear structure, then $\cat V$ is enriched over itself precisely if it is closed, i.e.\ if it admits an internal hom (cf. \cref{ex:monoidal-category-enriched-in-itself} below).

In Lurie's model of enriched categories it is easy to make everything functorial:

\begin{definition}
\begin{defenum}
	\item We denote by
	\begin{align}
		\Enr \subseteq \Op_{/\LM^\tensor}
	\end{align}
	the full subcategory spanned by the operad maps $\cat O^\tensor \to \LM^\tensor$ for which $\cat O^\tensor$ exhibits $\cat C \coloneqq \cat O_{\mathfrak m}$ as enriched over $\cat V^\tensor \coloneqq \cat O^\tensor \times_{\LM^\tensor} \Ass^\tensor$.

	\item \label{def:category-of-Lurie-enriched-categories} There is a natural functor
	\begin{align}
		\Enr \to \Ani \times \Op_{/\Ass^\tensor}, \qquad \cat O^\tensor \mapsto ((\cat O_{\mathfrak m})^\simeq, \cat O^\tensor \times_{\LM^\tensor} \Ass^\tensor)
	\end{align}
	Given an operad map $\cat V^\tensor \to \Ass^\tensor$, we call the fiber $\Enr_{\cat V}$ the \emph{category of $\cat V$-enriched categories}.
\end{defenum}
\end{definition}

We have introduced two models for enriched categories, see \cref{def:category-of-enriched-precategories,def:category-of-Lurie-enriched-categories}). It is one of the main results of \cite{Heine.2023} that these two models are equivalent:

\begin{theorem} \label{rslt:Heine-equivalence-of-models-for-enrichment}
\begin{thmenum}
	\item There is a fully faithful embedding
	\begin{align}
		\chi\colon \Enr \injto \PreEnr
	\end{align}
	over $\Ani \times \Op_{/\Ass^\tensor}$. Its essential image consists of exactly those enriched precategories $\cat C$ on anima $S$ for which the anima of equivalences of $\cat C$ coincides with $S$ (see the discussion below \cref{def:quiver-fibration}).

	\item \label{rslt:categorification-of-enriched-precategories} The embedding $\chi$ has a left adjoint $L\colon \PreEnr \to \Enr$ which restricts to fiberwise left adjoint functors
	\begin{align}
		L_{\cat V}\colon \PreEnr_{\cat V} \to \Enr_{\cat V}
	\end{align}
	for all operads $\cat V \in \Op_{/\Ass^\tensor}$. Moreover, for every $\cat V \in \Op_{/\Ass^\tensor}$ and every $\cat V$-enriched precategory $\cat C$, the unit of the adjunction $\cat C \to L_{\cat V} \cat C$ induces an essentially surjective map on the underlying anima of objects.
\end{thmenum}
\end{theorem}
\begin{proof}
Part (i) is \cite[Corollary~6.13]{Heine.2023} (see \cref{rmk:comparison-of-definitions-for-Lurie-enriched-cat} for an argument why our definitions coincide with the ones in \loccit). Part (ii) follows from the proof of (i), see \cite[Theorem~6.22]{Heine.2023}. For the convenience of the reader we sketch the construction of the functor $\chi$ as outlined in \cite[\S4.2]{Heine.2023}.

Suppose we are given an anima $S$, an operad $\cat V \in \Op_{/\Ass^\tensor}$ and a $\cat V$-enriched category $\cat C$ exhibited by some operad $\cat O^\tensor \to \LM^\tensor$. Using Day convolution for $\LM^\tensor$-operads one constructs an operad $\cat O'^\tensor \to \LM^\tensor$ exhibiting the category $\Fun(S, \cat C)$ as enriched over $\Quiv_S(\cat V)^\tensor$ (see \cite[Notations~4.17]{Heine.2023} and \cite[Proposition~4.30]{Heine.2023}). Explicitly, for any two functors $F, G \in \Fun(S, \cat C)$ the enriched hom from $F$ to $G$ is given by
\begin{align}
	\enrHom[\Quiv_S(\cat V)](F, G) = [S \times S \isom S^\op \times S \xto{F^\op \times G} \cat C^\op \times \cat C \xto{\enrHom[\cat V](\blank,\blank)} \cat V].
\end{align}
Now take $S = \cat C^\simeq$ and consider the canonical functor $F \in \Fun(\cat C^\simeq, \cat C)$. Its enriched endomorphisms are given by
\begin{align}
	\enrHom[\Quiv_{\cat C^\simeq}(\cat V)](F, F) = [\cat C^\simeq \times \cat C^\simeq \xto{\enrHom[\cat V](\blank,\blank)} \cat V].
\end{align}
They come naturally equipped with an algebra structure, as does every enriched endomorphism object. This algebra is a $\cat V$-enriched precategory and hence the desired $\chi(\cat C)$. One can perform the above construction functorially over $\Ani \times \Op_{/\Ass^\tensor}$ to produce the desired functor $\chi$ from (i).
\end{proof}

From now on we will use the functor $\chi$ to implicitly identify $\Enr$ as a full subcategory of $\PreEnr$ with left adjoint $L\colon \PreEnr \to \Enr$. If a specific model of enrichment is needed for a construction, we refer to the above two models as the \enquote{quiver model} and the \enquote{Lurie model}. It is high time for some examples:

\begin{example} \label{ex:enrichment-of-linear-category}
Let $\cat V$ be a monoidal category and let $\cat C$ be a $\cat V$-linear category. Suppose that for every $X \in \cat C$ the functor $\blank \tensor X\colon \cat V \to \cat C$ has a right adjoint $\enrHom[\cat V](X,\blank)$. Then $\cat C$ is naturally a $\cat V$-enriched category with $\cat V$-enriched homomorphisms given by $\enrHom[\cat V]$. Indeed, if $\cat O^\tensor \to \LM^\tensor$ denotes the operad exhibiting the $\cat V$-linear structure on $\cat C$ then it follows immediately from the definition of the Lurie model that $\cat O^\tensor$ exhibits the $\cat V$-enrichment of $\cat C$.
\end{example}

\begin{example} \label{ex:monoidal-category-enriched-in-itself}
Suppose that $\cat V$ is a \emph{closed} monoidal category, i.e.\ for all $V \in \cat V$ the functor $\blank \tensor V\colon \cat V \to \cat V$ has a right adjoint $\iHom(V, \blank)$. Then $\cat V$ is naturally enriched over itself with $\enrHom[\cat V] = \iHom$. Indeed, if $\cat V^\tensor \to \Ass^\tensor$ denotes the operad exhibiting the monoidal structure on $\cat V$ then $\cat V^\tensor \times_{\Ass^\tensor} \LM^\tensor$ exhibits $\cat V$ as a $\cat V$-linear category, so we reduce to \cref{ex:enrichment-of-linear-category}.
\end{example}

\begin{example} \label{ex:full-subcategory-of-enriched-category}
Let $\cat C$ be a $\cat V$-enriched category for some operad $\cat V^\tensor \to \Ass^\tensor$. Then for any collection $S$ of objects in $\cat C$ we can define a $\cat V$-enriched category $\cat C_S$ whose objects are given by $S$ and whose morphisms are the same as the ones in $\cat C$; this is easily seen in the Lurie model of enrichment, where we can take the appropriate full subcategory of the operad $\cat O^\tensor \to \LM^\tensor$ exhibiting the enrichment of $\cat C$. We call $\cat C_S$ the \emph{full subcategory of $\cat C$} spanned by $S$.
\end{example}

\begin{example} \label{ex:algebra-is-enriched-category}
Let $\cat V^\tensor \to \Ass^\tensor$ be an operad and let $A \in \Alg(\cat V)$ be an algebra. Let $S = *$ denote the final anima. Then $A$ is a map of operads $\Ass^\tensor = \Ass^\tensor_S \to \cat V^\tensor$, i.e.\ an algebra in $\Quiv_S(\cat V)$ (see \cref{rmk:comparison-of-precategories-to-GH}). We can thus canonically view $A$ as a $\cat V$-enriched precategory, and we denote the associated $\cat V$-enriched category by $*_A$. It has a single object $*$, whose $\cat V$-enriched endomorphisms are given by $A$.
\end{example}

\begin{example} \label{ex:trivial-enriched-category}
Let $\cat V$ be a monoidal category. Then by \cref{ex:trivial-algebra} there is a trivial algebra $\one_{\cat V} \in \Alg(\cat V)$, which by \cref{ex:algebra-is-enriched-category} induces a $\cat V$-enriched category $*_{\one_{\cat V}}$. We call this category the \emph{trivial $\cat V$-enriched category} and denote it by $*_{\cat V}$.
\end{example}

One important construction of categories is that of taking the opposite category. In the enriched setting, this construction has a simple implementation using the quiver model:

\begin{definition} \label{def:opposite-enriched-category}
Let $\cat C$ be a $\cat V$-enriched precategory for some operad $\cat V^\tensor \to \Ass^\tensor$ and suppose that $\cat C$ is represented by a quiver algebra $F \in \Alg(\Quiv_S(\cat V))$ for some anima $S$. By using the natural equivalence of operads
\begin{align}
	\Quiv_S(\cat V)^\rev = \Quiv_S(\cat V^\rev)
\end{align}
from \cite[Remark~4.23]{Heine.2023} together with \cref{rmk:rev-algebras-same-as-algebras} we can interpret $F$ as an algebra in $\Quiv_S(\cat V^\rev)$. This $\cat V^\rev$-enriched precategory is called the \emph{opposite category of $\cat C$} and denoted $\cat C^\op$.
\end{definition}

\begin{remark}
There is an explicit construction of the opposite enriched category in terms of the Lurie model, see \cite[Notation~10.2]{Heine.2023} and \cite[Corollary~10.5]{Heine.2023}. In particular, the opposite enriched precategory of an enriched category is automatically an enriched category.
\end{remark}

Explicitly the equivalence $\Quiv_S(\cat V)^\rev = \Quiv_S(\cat V^\rev)$ maps a quiver $F\colon S \times S \to \cat V$ to the quiver $F^\rev\colon S \times S \to \cat V$ with $F^\rev(X, Y) = F(Y, X)$. In particular if $\cat C$ is a $\cat V$-enriched category with $\cat V$-enriched homomorphisms $\enrHom[\cat V]_{\cat C}(X, Y)$ for any two $X, Y \in \cat C$ then $\cat C^\op$ is the $\cat V^\rev$-enriched category with the same anima of objects and with $\cat V^\rev$-enriched homomorphisms
\begin{align}
	\enrHom[\cat V^\rev]_{\cat C^\op}(X, Y) = \enrHom[\cat V]_{\cat C}(Y, X).
\end{align}
for all $X, Y \in \cat C$. It is also immediately clear that $(\cat C^\op)^\op = \cat C$.

By construction every enriched category has an underlying anima of objects. In fact, the Lurie model allows us to even define an underlying \emph{category}:

\begin{definition} \label{def:underlying-category-of-enriched-category}
Let $\cat V^\tensor \to \Ass^\tensor$ be an operad and let $\cat C$ be a $\cat V$-enriched category. If $\cat O^\tensor \to \LM^\tensor$ exhibits the $\cat V$-enrichment of $\cat C$ in the Lurie model then we denote
\begin{align}
	\ul{\cat C} \coloneqq \cat O^\tensor_{\mathfrak m}
\end{align}
and call it the \emph{underlying category of $\cat C$}. This construction provides a functor $\Enr \to \Cat$. See \cref{rslt:underlying-category-via-transfer-of-enrichment} for a more conceptual and model independent definition of the underlying category.
\end{definition}

\begin{example} \label{ex:Ani-enriched-category-same-as-category}
Equip $\Ani$ with the cartesian monoidal structure (see \cref{ex:cartesian-monoidal-structure}). Then the assignment $\cat C \mapsto \ul{\cat C}$ induces an equivalence
\begin{align}
	\Enr_\Ani = \Cat,
\end{align}
i.e.\ an $\Ani$-enriched category is the same as a category. See \cite[Lemma~6.18]{Heine.2023} for a proof.
\end{example}

Under mild assumptions on a monoidal category $\cat V$, the category of $\cat V$-enriched categories has all small limits and these limits have a very explicit description:

\begin{lemma} \label{rslt:limits-of-enriched-categories}
Let $I$ be a small category and let $\cat V^\tensor$ be an operad over $\Ass^\tensor$ such that the underlying category $\cat V$ admits limits indexed by $I$. Then $\Enr_{\cat V}$ admits limits indexed by $I$ and they are computed as follows. For a diagram $(\cat C_i)_{i\in I}$ in $\Enr_{\cat V}$ with limit $\cat C = \varprojlim_i \cat C_i$ we have
\begin{lemenum}
	\item $\ul{\cat C} = \varprojlim_i \ul{\cat C_i}$.
	\item For all objects $X = (X_i)_i$ and $Y = (Y_i)_i$ in $\cat C$ we have
	\begin{align}
		\enrHom[\cat V]_{\cat C}(Y, X) = \varprojlim_i \enrHom[\cat V]_{\cat C_i}(Y_i, X_i).
	\end{align}
\end{lemenum}
\end{lemma}
\begin{proof}
We denote
\begin{align}
	\Op^{\cat V}_{/\LM^\tensor} \coloneqq \Op_{/\LM^\tensor} \times_{\Op_{/\Ass^\tensor}} \{ \cat V^\tensor \},
\end{align}
where the map $\Op_{/\LM^\tensor} \to \Op_{/\Ass^\tensor}$ is given by base-change along the map $\Ass^\tensor \to \LM^\tensor$. In other words, $\Op^{\cat V}_{/\LM^\tensor}$ is the category of operads over $\LM^\tensor$ whose algebra part is identified with $\cat V^\tensor$. By \cref{rmk:non-symmetric-vs-symmetric-operads} and \cite[Proposition~3.33]{Heine.2023} we can identify $\Op^{\cat V}_{/\LM^\tensor}$ with a non-full subcategory of $\Cat_{/\cat V^\tensor}$, whose objects we denote by $\cat M^\circledast \to \cat V^\tensor$ as in \cite[Definition~3.11]{Heine.2023}. Analogous to the proof of \cref{rslt:limits-of-operads} one shows that $\Op^{\cat V}_{/\LM^\tensor}$ has all small limits and the forgetful functor $\Op^{\cat V}_{/\LM^\tensor} \to \Cat_{/\cat V^\tensor}$ preserves them. Note that by definition $\Enr_{\cat V} \subseteq \Op^{\cat V}_{/\LM^\tensor}$ is a full subcategory, so we need to show that this subcategory is stable under limits indexed by $I$.

Let now $(\cat C_i)_i$ be given as in the claim and let $\cat C \coloneqq \varprojlim_i \cat C_i$ be the limit in $\Op^{\cat V}_{/\LM^\tensor}$. By using the fact that this limit can be computed in $\Cat_{/\cat V^\tensor}$ (which gives us an explicit description of all the $\Hom$ anima) and then passing through the above identifications we deduce that $\cat C$ is an $\LM^\tensor$-operad whose fiber over $\mathfrak a$ is identified with $\cat V$, whose fiber over $\mathfrak m$ is identified with $\ul{\cat C} = \varprojlim_i \ul{\cat C_i}$, and whose $\Hom$'s satisfy the following property: For all $V_1, \dots, V_n \in \cat V$ and $(X_i)_i, (Y_i)_i \in \ul{\cat C}$ we have
\begin{align}
	\Mul_{\cat C}(V_1, \dots, V_n, (X_i)_i; (Y_i)_i) = \varprojlim_i \Mul_{\cat C_i}(V_1, \dots, V_n, X_i; Y_i).
\end{align}
We see immediately that
\begin{align}
	\enrHom[\cat V]((X_i)_i, (Y_i)_i) \coloneqq \varprojlim_i \enrHom[\cat V](X_i, Y_i)
\end{align}
defines a morphism object from $(X_i)_i$ to $(Y_i)_i$ in $\cat C$, hence $\cat C$ is indeed $\cat V$-enriched in the sense of \cref{def:Lurie-enriched-category}. From the proof we immediately deduce the description of $\cat C$ in (i) and (ii).
\end{proof}

\subsection{Enriched functors and Yoneda}

In the previous subsection we introduced the category $\Enr_{\cat V}$ of $\cat V$-enriched categories, for a fixed monoidal category $\cat V$. In the present subsection we will have a closer look at the homomorphisms in that category.

\begin{definition} \label{def:enriched-functors}
Let $\cat C_1$ and $\cat C_2$ be categories enriched in some operad $\cat V^\tensor \to \Ass^\tensor$ and let $\cat O_1^\tensor, \cat O_2^\tensor \to \LM^\tensor$ be operads exhibiting the $\cat V$-enrichment in the Lurie model.
\begin{defenum}
	\item We denote by
	\begin{align}
		\enrFun[\cat V](\cat C_1, \cat C_2) = \Alg_\LM(\cat O_1, \cat O_2) \times_{\Alg_\Ass(\cat V, \cat V)} \{ \id \}
	\end{align}
	the \emph{category of $\cat V$-enriched functors} from $\cat C_1$ to $\cat C_2$.

	\item \label{def:forgetful-functor-on-enriched-functor-categories} By restricting to the fiber over $\mathfrak m \in \LM$ we get a conservative forgetful functor
	\begin{align}
		\enrFun[\cat V](\cat C_1, \cat C_2) \to \Fun(\ul{\cat C_1}, \ul{\cat C_2}), \qquad f \mapsto \ul{f}
	\end{align}
	(see \cite[Lemma~A.4.5.(a)]{Mann.2022a} and note that its proof works without assuming that $\cat V$ is a monoidal category). We call the image $\ul{F}\colon \ul{\cat C_1} \to \ul{\cat C_2}$ of a $\cat V$-enriched functor $F\colon \cat C_1 \to \cat C_2$ the \emph{underlying functor of $F$}.
\end{defenum}
\end{definition}

By unravelling the definitions, we see that a $\cat V$-enriched functor $\cat C_1 \to \cat C_2$ consists of a map $F\colon \cat C_1^\simeq \to \cat C_2^\simeq$ of the underlying anima and for any two objects $X, Y \in \cat C_1$ a map
\begin{align}
	\enrHom[\cat V]_{\cat C_1}(X, Y) \to \enrHom[\cat V]_{\cat C_2}(F(X), F(Y))
\end{align}
in $\cat V$, together with coherences exhibiting the compatibilities of the above maps with composition and the identity maps (cf. the discussion following \cite[Definition~4.4]{Mann.2022a}).

\begin{example}\label{ex:lax-linear-equals-enriched}
Let $\cat V$ be a (symmetric) monoidal category and let $\cat C$, $\cat C'$ be $\cat V$-linear categories. Then, by definition, the category of lax $\cat V$-linear functors agrees with the category of $\cat V$-enriched functors.
\end{example}

\begin{remark} \label{rmk:compatibility-of-enriched-functors-with-Heine}
It is clear from the definition that the underlying anima of $\enrFun[\cat V](\cat C_1, \cat C_2)$ is $\Hom_{\Enr_{\cat V}}(\cat C_1, \cat C_2)$. One checks that under the identifications of \cref{rmk:comparison-of-definitions-for-Lurie-enriched-cat}, $\enrFun[\cat V]$ coincides with $\operatorname{LaxLinFun}_{\cat V}$ from \cite[Notation~3.72]{Heine.2023}. It is shown in \cite[Lemma~3.64]{Heine.2023} and \cite[Remark~3.90]{Heine.2023} that there is a left action of $\Cat$ on $\Enr_{\cat V}$ such that the associated morphism objects are given by $\enrFun[\cat V](\blank,\blank)$. By \cref{ex:enrichment-of-linear-category} we can thus view $\Enr_{\cat V}$ as a $\Cat$-enriched category, i.e.\ a 2-category (see \cref{rslt:2-categorical-enhancement-of-Enr} for a different construction of this 2-category). 
\end{remark}

\begin{remark}
There is in general no natural $\cat V$-enrichment on $\enrFun[\cat V](\cat C_1, \cat C_2)$. This is similar to the fact that for a non-commutative ring $A$ and $A$-modules $M$ and $N$, there is no natural $A$-module structure on $\Hom_A(M, N)$. In the setting of $A$-modules this can be fixed by putting additional right $A$-module structures on $M$ or $N$. A similar phenomenon happens for enrichment, cf. \cite[Proposition~3.86]{Heine.2023}.
\end{remark}

\begin{example} \label{rslt:universal-property-of-trivial-enriched-category}
Let $\cat V$ be a monoidal category and let $*_{\cat V}$ be the trivial $\cat V$-enriched category from \cref{ex:trivial-enriched-category}. Then for every $\cat V$-enriched category $\cat C$ there is a natural isomorphism
\begin{align}
	\enrFun[\cat V](*_{\cat V}, \cat C) = \ul{\cat C}.
\end{align}
Namely, by an explicit computation using operadic Kan extensions this can be shown to be true for the category $*'_{\cat V}$ from \cite[Definition~A.4.6]{Mann.2022a} (which is an explicit definition of $*_{\cat V}$ in Lurie's model). In particular $*'_{\cat V}$ is characterized by the universal property that $\Hom_{\Enr_{\cat V}}(*'_{\cat V}, \cat C) = \ul{\cat C}^\simeq$. The same universal property holds for $*_{\cat V}$ by \cite[Lemma~5.1.2]{Gepner-Haugseng.2015}, hence there is a natural equivalence $*'_{\cat V} = *_{\cat V}$ and we conclude.
\end{example}

Let us single out two important special cases of enriched functors: the fully faithful functors and the essentially surjective functors. As we will see shortly, a fully faithful and essentially surjective functor is an equivalence, as expected.

\begin{definition}
Let $\cat V^\tensor \to \Ass^\tensor$ be an operad and let $F\colon \cat C_1 \to \cat C_2$ be a $\cat V$-enriched functor of $\cat V$-enriched categories.
\begin{defenum}
	\item \label{def:fully-faithful-enriched-functor} $F$ is called \emph{fully faithful} if for all objects $X, Y \in \cat C_1$ the induced map $\enrHom[\cat V]_{\cat C_1}(X, Y) \isoto \enrHom[\cat V]_{\cat C_2}(F(X), F(Y))$ is an isomorphism in $\cat V$.

	\item \label{def:essentially-surjective-enriched-functor} $F$ is called \emph{essentially surjective} if the underlying functor $\ul{F}$ is essentially surjective.
\end{defenum}
\end{definition}

\begin{lemma} \label{rslt:enriched-surjective-plus-fully-faithful-implies-isom}
Let $\cat V^\tensor \to \Ass^\tensor$ be an operad and let $F\colon \cat C_1 \to \cat C_2$ be a $\cat V$-enriched functor of $\cat V$-enriched categories. Then $F$ is an isomorphism if and only if it is fully faithful and essentially surjective.
\end{lemma}
\begin{proof}
Let $\cat O_1^\tensor, \cat O_2^\tensor \to \LM^\tensor$ be the operads exhibiting the $\cat V$-enrichments of $\cat C_1$ and $\cat C_2$ in the Lurie model. Then $F$ is essentially surjective in the sense of \cref{def:essentially-surjective-enriched-functor} if and only if it is essentially surjective as a functor $\cat O_1^\tensor \to \cat O_2^\tensor$. Similarly, one checks easily from the definition of morphism objects (see \cref{def:morphism-objects-for-enrichment}) that $F$ is fully faithful in the sense of \cref{def:fully-faithful-enriched-functor} if and only if it is fully faithful as a functor $\cat O_1^\tensor \to \cat O_2^\tensor$. This implies the claim.
\end{proof}

We now come to the enriched version of the Yoneda lemma, which is a helpful tool for constructing enriched functors. We first need to introduce the enriched presheaf category $\PSh[\cat V](\cat C)$ of a $\cat V$-enriched category $\cat C$:

\begin{definition} \label{def:enriched-yoneda-category}
Let $\cat V$ be a monoidal category such that $\cat V^\rev$ is closed. Then $\cat V$ is enriched in $\cat V^\rev$ by \cref{ex:monoidal-category-enriched-in-itself}. For a $\cat V$-enriched category $\cat C$ we denote
\begin{align}
	\PSh[\cat V](\cat C) \coloneqq \enrFun[\cat V^\rev](\cat C^\op, \cat V)
\end{align}
and call it the \emph{category of $\cat V$-enriched presheaves on $\cat C$}.
\end{definition}

\begin{remark} \label{rmk:generalized-enriched-presheaves}
The category of $\cat V$-enriched presheaves can be constructed for any operad $\cat V^\tensor \to \Ass^\tensor$, without assuming it to be a monoidal category with $\cat V^\rev$ closed. Then $\cat V$ is not necessarily enriched in $\cat V^\rev$, but it is still \emph{weakly enriched} in the sense that one can still form the $\LM$-operad $\cat V'^\tensor := \cat V^{\rev,\tensor} \times_{\Ass^\tensor} \LM^\tensor$. This is enough to define the category of $\cat V^\rev$-enriched functors $\cat C^\op \to \cat V$ as in \cref{def:enriched-functors}. One can then even upgrade $\PSh[\cat V](\cat C)$ to a weakly $\cat V$-enriched category. The reader is invited to have a look at \cite[Proposition~3.86]{Heine.2023} for these definitions; we decided to keep things as simple as possible in this paper.
\end{remark}

Having defined the category of enriched presheaves, we can now construct the Yoneda embedding. We stick to its simplest formulation, as this is all we need in this paper.

\begin{theorem} \label{rslt:enriched-Yoneda-embedding}
Let $\cat V$ be a monoidal category such that $\cat V^\rev$ is closed. Then for every $\cat V$-enriched category $\cat C$ there is a natural fully faithful embedding
\begin{align}
	\Yo\colon \ul{\cat C} \injto \PSh[\cat V](\cat C)
\end{align}
such that for every object $X \in \cat C$ the underlying functor of $\Yo(X)$ is the functor $\enrHom[\cat V](\blank, X)\colon \cat C^\op \to \cat V$. Moreover, for every $f \in \PSh[\cat V](\cat C)$ and every $X \in \cat C$ there is a natural isomorphism of anima
\begin{align}
	\Hom(\enrHom[\cat V](\blank,X), f) = \Hom(\one_{\cat V}, f(X)).
\end{align}
\end{theorem}
\begin{proof}
See \cite[Corollary~A.4.9]{Mann.2022a} or \cite[Lemma 10.1]{Heine.2023}, where the latter reference also works in the more general setting discussed in \cref{rmk:generalized-enriched-presheaves}. We sketch the proof of the first reference. For every $X \in \cat C$ we consider the functor $\ev_X\colon \PSh[\cat V](\cat C) \to \cat V$ sending an enriched functor $f\colon \cat C^\op \to \cat V$ to $f(X) \in \cat V$. Using operadic Kan extensions one shows that this functor admits a left adjoint $\ell_X\colon \cat V \to \PSh[\cat V](\cat C)$, which can be explicitly computed as $\ell_X(V) = \enrHom[\cat V](\blank, X) \tensor V$. We have thus constructed the $\cat V^\rev$-enriched functor $\enrHom[\cat V](\blank, X)\colon \cat C^\op \to \cat V$ for fixed $X$. By a standard routine we can make this construction functorial in $X$: Consider the functor
\begin{align}
	\alpha\colon \ul{\cat C}^\op \times \cat V^\op \times \PSh[\cat V](\cat C) \to \Ani, \qquad (X, V, f) \mapsto \Hom(V, f(X)).
\end{align}
We can view this functor as a functor $\alpha\colon \ul{\cat C}^\op \times \cat V^\op \to \PSh(\PSh[\cat V](\cat C)^\op)$, where $\PSh(\blank) = \Fun((\blank)^\op, \Ani)$ is the usual presheaf category. By what we have shown above, for each $(X, V) \in \ul{\cat C}^\op \times \cat V^\op$ the object $\alpha(X, Y) \in \PSh(\PSh[\cat V](\cat C)^\op)$ is representable by $\enrHom[\cat V](\blank, X) \tensor V$, i.e.\ lies in the image of the Yoneda embedding $\PSh[\cat V](\cat C)^\op \injto \PSh(\PSh[\cat V](\cat C)^\op)$. Composing $\alpha$ with an inverse of the Yoneda embedding and applying $(\blank)^\op$ yields a functor $\ul{\cat C} \times \cat V \to \PSh[\cat V](\cat C)$. Now restrict this functor to the subcategory $\ul{\cat C} = \ul{\cat C} \times \{ \one_{\cat V} \} \subset \ul{\cat C} \times \cat V$ to obtain the desired Yoneda functor $\ul{\cat C} \to \PSh[\cat V](\cat C)$. The full faithfulness of this functor follows from the above description of $\enrHom[\cat V](\blank, X)$ via the left adjoint functor $\ell_X$, which also shows the claimed $\Hom$-identity.
\end{proof}

\begin{remark}
See \cite[Theorem~5.1]{Heine.2023} for a much stronger version of the enriched Yoneda embedding.
\end{remark}

In this paper we mainly need the enriched Yoneda embedding to construct the functor $\enrHom[\cat V](\blank, X)$ as a $\cat V^\rev$-enriched functor. In fact, we easily obtain the following:

\begin{corollary}\label{rslt:symmetric-monoidal-self-duality}
Let $\cat V$ be a symmetric monoidal category in which every object is dualizable (see \cref{def:dualizable-object}). Then $\cat V$ is enriched in itself and the assignment $X \mapsto X^\vee$ defines a $\cat V$-enriched equivalence $\cat V^\op \isoto \cat V$.
\end{corollary}
\begin{proof}
By \cref{rslt:properties-of-dualizable-objects} one sees immediately that $\cat V$ is closed. Hence $\cat V$ is enriched in itself by \cref{ex:monoidal-category-enriched-in-itself}. Using the isomorphism $\cat V = \cat V^\rev$ induced by the symmetry of the monoidal structure, we obtain from \cref{rslt:enriched-Yoneda-embedding} a functor $\cat V \injto \enrFun[\cat V](\cat V^\op, \cat V)$ sending $Y \in \cat V$ to $\enrHom[\cat V](\blank, Y) = \iHom_{\cat V}(\blank, Y)$. Taking $Y = \one_{\cat V}$ we obtain the desired $\cat V$-enriched functor $\cat V^\op \to \cat V$ sending $X \mapsto X^\vee$. To check that this functor is an equivalence, it is enough to show that it is fully faithful and essentially surjective (see \cref{rslt:enriched-surjective-plus-fully-faithful-implies-isom}). The essential surjectivity follows immediately from the fact that $P = (P^\vee)^\vee$ for all $P \in \cat V$. The full faithfulness requires us to show that for all $P, Q \in \cat V$ the natural map
\begin{align}
	\iHom_{\cat V}(Q, P) \isoto \iHom_{\cat V}(P^\vee, Q^\vee)
\end{align}
is an isomorphism, which follows immediately from \cref{rslt:properties-of-dualizable-objects}.
\end{proof}

Under mild assumptions the category $\enrFun[\cat V](\cat C_1, \cat C_2)$ has the same limits and colimits as those in $\cat D$ and they are computed pointwise:

\begin{lemma} \label{rslt:lim-and-colim-in-enrFun}
Let $\cat V^\tensor$ be a monoidal category and let $\cat C_1$ and $\cat C_2$ be $\cat V$-enriched categories such that $\cat C_2$ is $\cat V$-linear. Let $\cat I$ be some category.
\begin{lemenum}
	\item \label{rslt:limits-in-enrFun} Suppose $\cat C_2$ admits $\cat I$-indexed limits. Then $\enrFun[\cat V](\cat C_1, \cat C_2)$ admits $\cat I$-indexed limits and the forgetful functor $\enrFun[\cat V](\cat C_1, \cat C_2) \to \Fun(\cat C_1, \cat C_2)$ preserves them.

	\item \label{rslt:colimits-in-enrFun} Suppose that $\cat C_2$ admits $\cat I$-indexed colimits and the operation $\tensor\colon \cat V \times \cat C_2 \to \cat C_2$ preserves $\cat I$-indexed colimits in the second argument. Then $\enrFun[\cat V](\cat C_1, \cat C_2)$ admits $\cat I$-indexed colimits and the forgetful functor $\enrFun[\cat V](\cat C_1, \cat C_2) \to \Fun(\cat C_1, \cat C_2)$ preserves them.
\end{lemenum}
\end{lemma}
\begin{proof}
This is proved in \cite[Lemma~A.4.5]{Mann.2022a} but the proof is flawed, as it does not fix the enriching category $\cat V$ properly. However, the same argument works in the setting of non-symmetric operads, with the definition of enriched functors from \cite[Notation~3.72]{Heine.2023} (and \cref{rmk:compatibility-of-enriched-functors-with-Heine}), see also \cite[Lemma~3.47(1)]{Heine.2023}.
\end{proof}

\subsection{Transfer of enrichment} \label{sec:enr.transfer}

Given monoidal categories $\cat V$ and $\cat W$, a lax monoidal functor $\alpha\colon \cat V \to \cat W$, and a $\cat V$-enriched category $\cat C$, one can transfer the enrichment of $\cat C$ along $\alpha$ to construct a $\cat W$-enriched category with the same objects, where the $\cat W$-enriched homs are obtained from the $\cat V$-enriched homs of $\cat C$ by applying $\alpha$. This construction plays a crucial role in this paper, so in the following we provide a proper definition and some basic properties of it.

Note that by \cite[Remark~4.41]{Heine.2023} the functor $\PreEnr \to \Op_{/\Ass^\tensor}$ is a cocartesian fibration. By \cref{rslt:categorification-of-enriched-precategories} the fully faithful embedding $\Enr \subset \PreEnr$ has a fiberwise left adjoint over $\Op_{/\Ass^\tensor}$, which implies that $\Enr \to \Op_{/\Ass^\tensor}$ is still a cocartesian fibration. This allows us to make the following definition:

\begin{definition} \label{def:transfer-of-enrichment}
The cocartesian fibration $\Enr \to \Op_{/\Ass^\tensor}$ induces a functor
\begin{align}
	\Op_{/\Ass^\tensor} \to \Cat, \qquad \cat V^\tensor \mapsto \Enr_{\cat V}.
\end{align}
For every map $\alpha\colon \cat V^\tensor \to \cat W^\tensor$ in $\Op_{/\Ass^\tensor}$ we denote by
\begin{align}
	\tau_\alpha\colon \Enr_{\cat V} \to \Enr_{\cat W}
\end{align}
the induced functor and call it the \emph{transfer of enrichment along $\alpha$}.
\end{definition}

Let us unpack the definition. Given a map $\alpha\colon \cat V^\tensor \to \cat W^\tensor$ in $\Op_{/\Ass^\tensor}$ and an anima $S$, composition with $\alpha$ induces a map of operads
\begin{align}
	\tau'_\alpha\colon \Quiv_S(\cat V)^\tensor \to \Quiv_S(\cat W)^\tensor.
\end{align}
This induces a functor
\begin{align}
	\tau''_\alpha\colon \Alg(\Quiv_S(\cat V)) \to \Alg(\Quiv_S(\cat W)),
\end{align}
i.e.\ a functor from $\cat V$-enriched precategories on $S$ to $\cat W$-enriched precategories on $S$. This is the functor that is encoded in the cocartesian fibration $\PreEnr \to \Op_{/\Ass^\tensor}$. If $\cat C$ is a $\cat V$-enriched category then $\tau''_\alpha(\cat C)$ may still be a \emph{pre}category, because the transfer of enrichment along $\alpha$ may very well introduce new isomorphisms (so that $S$ is not the correct anima of objects for the resulting $\cat W$-enriched precategory). We therefore also need to apply $L_{\cat W}\colon \PreEnr_{\cat W} \to \Enr_{\cat W}$ to arrive at the functor $\tau_\alpha = L_{\cat W} \comp \tau''_\alpha$.

More concretely, for any $\cat V$-enriched category $\cat C$, the $\cat W$-enriched category $\tau_\alpha(\cat C)$ has the same objects as $\cat C$ and for all $X, Y \in \cat C$ we have
\begin{align}
	\enrHom[\cat W]_{\tau_\alpha(\cat C)}(X, Y) = \alpha(\enrHom[\cat V]_{\cat C}(X, Y)).
\end{align}

\begin{remark}
The transfer of enrichment makes crucial use of the quiver model for enriched categories. It is hard to perform a similar construction directly in the Lurie model.
\end{remark}

The transfer of enrichment $\tau_\alpha$ along a morphism $\alpha\colon \cat V^\tensor \to \cat W^\tensor$ should satisfy a certain 2-categorical functoriality in $\alpha$. In this paper we will content ourselves with the following:

\begin{lemma} \label{rslt:2-functoriality-of-transfer-of-enrichment}
Let $\cat V^\tensor$ and $\cat W^\tensor$ be operads over $\Ass^\tensor$. Then the assignment $\alpha \mapsto \tau_\alpha$ defines a functor
\begin{align}
	\tau\colon \Alg_{\Ass}(\cat V, \cat W) \to \Fun(\Enr_{\cat V}, \Enr_{\cat W}),
\end{align}
where the left-hand side denotes the functor category of operad maps $\cat V^\tensor \to \cat W^\tensor$ over $\Ass^\tensor$ (see \cite[Definition 2.1.3.1]{HA}). Given a natural transformation $\nu\colon \alpha \to \alpha'$ of operad maps $\cat V^\tensor \to \cat W^\tensor$, the associated natural transformation $\tau_\alpha$ sends a $\cat V$-enriched category $\cat C$ to the $\cat W$-enriched functor $\tau_\alpha(\cat C) \to \tau_{\alpha'}(\cat C)$ which acts on objects as the identity and on morphisms as $\nu$.
\end{lemma}
\begin{proof}
We start with a preliminary remark. Note that on underlying anima the desired functor reduces to the map
\begin{align}
	\Hom_{\Op_{/\Ass^\tensor}}(\cat V^\tensor, \cat W^\tensor) \to \Hom_{\Cat}(\Enr_{\cat V}, \Enr_{\cat W}) 
\end{align}
which comes from the straightening $\Op_{/\Ass^\tensor} \to \Cat$, $\cat V^\tensor \mapsto \Enr_{\cat V}$ of the cocartesian fibration $\Enr \to \Op_{/\Ass^\tensor}$. This straightening should upgrade to a 2-functor, but we lack the tools to show this. However, we can at least show that this desired 2-functor exists on functor categories, which is the content of this proof.

Let us now come to the proof. Let us denote by $\Quiv(\cat V)^\tensor \to \Ani \times \Ass^\tensor$ the fiber of $\Quiv^\tensor$ (see \cref{def:quiver-fibration}) over $\cat V^\tensor$; we similarly denote $\Quiv(\cat W)^\tensor$. Then by definition we have $\PreEnr_{\cat V} = \Alg^\Ani(\Quiv(\cat V)^\tensor)$, where the notation $\Alg^{(\blank)}$ is taken from \cite[Notation~2.56]{Heine.2023} (here we take $\cat O = \Ass$ and omit it from the notation). From the explicit construction of $\Alg^{(\blank)}$ in \cite[Remark~2.57]{Heine.2023} it follows that $\Alg^\Ani(\blank)$ is 2-functorial (it is essentially a functor category) and hence induces a functor
\begin{align}
	\Alg_{\Ani \times \Ass^\tensor}(\Quiv(\cat V)^\tensor, \Quiv(\cat W)^\tensor) \to \Fun(\PreEnr_{\cat V}, \PreEnr_{\cat W}),
\end{align}
where the left-hand side denotes the functor category of generalized operads, as defined in \cite[Notation~2.34]{Heine.2023}. By using the embedding $\Enr_{\cat V} \injto \PreEnr_{\cat V}$ and the left adjoint functor $\PreEnr_{\cat W} \to \Enr_{\cat W}$ from \cref{rslt:categorification-of-enriched-precategories} we obtain a functor
\begin{align}
	\Fun(\PreEnr_{\cat V}, \PreEnr_{\cat W}) \to \Fun(\Enr_{\cat V}, \Enr_{\cat W}).
\end{align}
Combining this with the above 2-functoriality of $\Alg^\Ani(\blank)$, we are reduced to constructing a functor
\begin{align}
	\tau'\colon \Alg_\Ass(\cat V, \cat W) \to \Alg_{\Ani \times \Ass^\tensor}(\Quiv(\cat V)^\tensor, \Quiv(\cat W)^\tensor).
\end{align}
In other words, we need to see that the functor $\Quiv(\blank)^\tensor$ is 2-functorial in the appropriate sense. By the explicit construction of $\Quiv(\blank)^\tensor$ (see \cite[Proposition~11.4]{Heine.2023}) this reduces to showing that the construction $\Fun_S^T(\cat C, \blank)$ from \cite[Notation~3.80]{Heine.2023} is 2-functorial. But this construction is right adjoint to a $\Cat$-linear functor and hence lax $\Cat$-linear (equivalently, $\Cat$-enriched, see \cref{ex:lax-linear-equals-enriched}) by \cref{rslt:monoidal-right-adjoint}. This finishes the construction of the claimed functor. Note that for fixed $\alpha\colon \cat V^\tensor \to \cat W^\tensor$ we recover the explicit construction of $\tau_\alpha$ and $\tau'_\alpha$ above. One easily verifies that $\tau$ has the claimed description on natural transformations.
\end{proof}

From the explicit description of the transfer of enrichment it seems clear that it commutes with taking opposite categories. Here is a formal proof of this fact:

\begin{lemma} \label{rslt:transfer-of-enrichment-commutes-with-op}
Let $\alpha\colon \cat V^\tensor \to \cat W^\tensor$ be a map in $\Op_{/\Ass^\tensor}$. Then for every $\cat V$-enriched category $\cat C$ there is a natural equivalence
\begin{align}
	(\tau_\alpha \cat C)^\op = \tau_\alpha(\cat C^\op)
\end{align}
of $\cat W^\rev$-enriched categories.
\end{lemma}
\begin{proof}
The claim is clear on enriched precategories (i.e.\ for the functor $\tau''_\alpha$ above), as there the transfer of enrichment is induced by the postcomposition of a quiver with the map $\alpha$. It remains to check that taking opposite categories commutes with completion of enriched precategories, i.e.\ that the following diagram commutes:
\begin{equation}\begin{tikzcd}
	\PreEnr_{\cat W} \arrow[r,"L_{\cat W}"] \arrow[d,"(\blank)^\op",swap] & \Enr_{\cat W} \arrow[d,"(\blank)^\op"] \\
	\PreEnr_{\cat W^\rev} \arrow[r,"L_{\cat W^\rev}"] & \Enr_{\cat W^\rev}
\end{tikzcd}\end{equation}
To prove this, we note that the embedding $\Enr_{\cat W} \injto \PreEnr_{\cat W}$ commutes with taking opposite categories, hence there is a natural transformation of functors $L_{\cat W^\rev} \comp (\blank)^\op \to (\blank)^\op \comp L_{\cat W}$. To see that this natural transformation is an equivalence, by \cite[Theorem~6.7]{Heine.2023} it is enough to show that $(\blank)^\op$ preserves essentially surjective and fully faithful functors of enriched precategories, which is clear. Here \emph{essentially surjective} means that the induced map on the anima of objects is essentially surjective, and \emph{fully faithful} means that the induced maps on enriched hom-objects are isomorphisms.
\end{proof}

Using the transfer of enrichment, we get a new interpretation of the underlying category of an enriched category:

\begin{lemma} \label{rslt:underlying-category-via-transfer-of-enrichment}
Let $\cat V^\tensor \to \Ass^\tensor$ be an operad and let $\pi\colon \cat V^\tensor \to \Ani^\times$ be the map of operads from \cref{ex:lax-monoidal-functor-corep-by-1}. Then under the identification $\Enr_\Ani = \Cat$ from \cref{ex:Ani-enriched-category-same-as-category} there is a natural isomorphism
\begin{align}
	\ul{\cat C} = \tau_\pi(\cat C)
\end{align}
for every $\cat V$-enriched category $\cat C$. If $\cat V$ is a monoidal category then $\pi$ is the functor $\Hom_{\cat V}(\one_{\cat V}, \blank)$, so in particular $\ul{\cat C}$ is the category with the same anima of objects as $\cat C$ and with hom's given by
\begin{align}
	\Hom_{\ul{\cat C}}(X, Y) = \Hom_{\cat V}(\one_{\cat V}, \enrHom[\cat V]_{\cat C}(X, Y))
\end{align}
for $X, Y \in \cat C$.
\end{lemma}
\begin{proof}
We will show that there is a natural isomorphism $\ul{\cat C} = \tau_\pi(\cat C)$ as functors $\Enr_{\cat V} \to \Cat$. The hardest part about this statement is to construct the natural transformation in the first place, for which we perform the following more general construction: Given a map $\alpha\colon \cat V^\tensor \to \cat W^\tensor$ of operads over $\Ass^\tensor$, we will construct a natural transformation $\ul{\cat C} \to \ul{\tau_\alpha(\cat C)}$ of functors $\Enr_{\cat V} \to \Cat$.

By \cref{rslt:exponentiation-of-cocartesian-fibrations} the cocartesian fibration $\Enr \to \Op_{/\Ass^\tensor}$ induces a cocartesian fibration
\begin{align}
\Fun(\Enr_{\cat V}, \Enr) \to \Fun(\Enr_{\cat V}, \Op_{/\Ass^\tensor}).
\end{align}
Let $\iota\colon \Enr_{\cat V} \injto \Enr$ be the inclusion and denote $\varphi\colon \iota \to \alpha_*(\iota)$ the cocartesian lift of $\delta(\alpha)$, where $\delta \colon \Op_{/\Ass^\tensor} \to \Fun(\Enr_{\cat V}, \Op_{/\Ass^\tensor})$ is the diagonal functor defined by precomposing with the unique map $\Enr_{\cat V} \to *$. Observe (again by \cref{rslt:exponentiation-of-cocartesian-fibrations}) that the evaluation $\varphi(\cat C)$ identifies with the cocartesian edge $\cat C\to \tau_\alpha(\cat C)$ for all $\cat C \in \Enr_{\cat V}$. Whiskering $\varphi$ with the functor $\Enr \to \Cat$, $\cat C \mapsto \ul{\cat C}$ yields the natural transformation $\ul{\cat C} \to \ul{\tau_{\alpha}(\cat C)}$ of functors $\Enr_{\cat V} \to \Cat$.

We now apply the above general construction to the case $\cat W = \Ani$ and $\alpha = \pi$. We obtain a natural transformation $\ul{\cat C} \to \ul{\tau_\pi(\cat C)}$ of functors $\Enr_{\cat V} \to \Cat$. We need to show that this natural transformation is an isomorphism, i.e.\ for a fixed $\cat V$-enriched category $\cat C$ we need to see that the functor $\ul{\cat C} \isoto \ul{\tau_\pi(\cat C)}$ is an equivalence. Let $\cat O^\tensor \to \LM^\tensor$ be the operad exhibiting the $\cat V$-enrichment of $\cat C$, so that $\ul{\cat C} = \cat O_{\mathfrak m}$. From the definition of morphism objects in $\cat O^\tensor$ one sees immediately that for all $X, Y \in \cat C$ there is a natural isomorphism
\begin{align}
	\Hom_{\ul{\cat C}}(X, Y) = \Hom_{\cat V^\tensor}(\emptyset_{\cat V}, \enrHom[\cat V]_{\cat C}(X, Y)) = \pi(\enrHom[\cat V]_{\cat C}(X, Y)) = \Hom_{\ul{\tau_\pi(\cat C)}}(X, Y),
\end{align}
where $\emptyset_{\cat V}$ denotes the unique object in $\cat V^\tensor_{\langle 0 \rangle}$. By going through the equivalence in \cref{rslt:Heine-equivalence-of-models-for-enrichment} one checks that this natural isomorphism is the one induced by the functor $\ul{\cat C} \to \ul{\tau_\pi(\cat C)}$, so that this functor is fully faithful. It follows from the explicit description of the transfer of enrichment that the functor is also essentially surjective.
\end{proof}

As part of the proof of \cref{rslt:underlying-category-via-transfer-of-enrichment} we showed that the transfer of enrichment induces a functor on underlying categories. Let us state this more explicitly:

\begin{lemma} \label{rslt:functor-on-underlying-categories-of-transfer-of-enrichment}
Let $\alpha\colon \cat V^\tensor \to \cat W^\tensor$ be a map in $\Op_{/\Ass^\tensor}$. Then for every $\cat V$-enriched category $\cat C$ there is a natural essentially surjective functor
\begin{align}
	F_\alpha\colon \ul{\cat C} \to \ul{\tau_\alpha(\cat C)}
\end{align}
which is the identity on objects. If $\cat V^\tensor$ and $\cat W^\tensor$ are monoidal categories then the above functor induces the natural map
\begin{align}
	\Hom_{\ul{\cat C}}(X, Y) &= \Hom_{\cat V}(\one_{\cat V}, \enrHom[\cat V]_{\cat C}(X, Y))\\
	&\to \Hom_{\cat W}(\alpha(\one_{\cat V}), \alpha(\enrHom[\cat V]_{\cat C}(X, Y)))\\
	&\to \Hom_{\cat W}(\one_{\cat W}, \alpha(\enrHom[\cat V]_{\cat C}(X, Y)))\\
	&= \Hom_{\ul{\tau_\alpha(\cat C)}}(X, Y)
\end{align}
on hom anima for all $X, Y \in \cat C$. Here the second map is induced by the natural map $\one_{\cat W} \to \alpha(\one_{\cat V})$.
\end{lemma}
\begin{proof}
The functor $F_\alpha$ is obtained from the cocartesian edge $\cat C \to \tau_\alpha(\cat C)$ in $\Enr$ by applying the functor $\Enr \to \Cat$, $\cat C \mapsto \ul{\cat C}$. Naturality of this construction was proved in \cref{rslt:underlying-category-via-transfer-of-enrichment}. The explicit description of $F_\alpha$ on enriched hom's follows from the construction of the cocartesian fibration $\PreEnr \to \Op_{/\Ass^\tensor}$, we leave the details to the reader.
\end{proof}

In this paper we also need the following somewhat exotic construction of an enriched functor associated with the transfer of enrichment:

\begin{lemma} \label{rslt:functor-from-transfer-of-self-enrichment}
Let $\alpha\colon \cat V \to \cat W$ be a lax monoidal functor of closed monoidal categories. Then $\cat V$ and $\cat W$ are enriched in themselves and there is a natural $\cat W$-enriched functor
\begin{align}
	G_\alpha\colon \tau_\alpha(\cat V) \to \cat W, \qquad X \mapsto \alpha(X).
\end{align}
On $\cat W$-enriched hom's it acts as the natural map
\begin{align}
	\enrHom[\cat W]_{\tau_\alpha(\cat V)}(X, Y) = \alpha(\iHom_{\cat V}(X, Y)) \to \iHom_{\cat W}(\alpha(X), \alpha(Y)) = \enrHom[\cat W]_{\cat W}(\alpha(X), \alpha(Y))
\end{align}
for all $X, Y \in \cat V$. It can be described as the functor $G_\alpha = \enrHom[\cat W]_{\tau_\alpha(\cat V)}(\one_{\cat V}, \blank)$.
\end{lemma}
\begin{proof}
One can use the last part of the claim to construct $G_\alpha$, but we will perform a different construction in the following and then show that it is equivalent to $\enrHom[\cat W]_{\tau_\alpha(\cat V)}(\one_{\cat V}, \blank)$. This alternative construction will be crucial in \cref{rslt:factor-map-of-monoidal-categories-into-enriched-maps} below.

The self-enrichment of $\cat V$ and $\cat W$ is discussed in \cref{ex:monoidal-category-enriched-in-itself}. Note that in the Lurie model, these enrichments are exhibited by the operads $\cat V'^\tensor \coloneqq \cat V^\tensor \times_{\Ass^\tensor} \LM^\tensor$ and $\cat W'^\tensor \coloneqq \cat W^\tensor \times_{\Ass^\tensor} \LM^\tensor$. The map $\alpha$ thus induces a map $\alpha'\colon \cat V'^\tensor \to \cat W'^\tensor$ of operads over $\LM^\tensor$, i.e.\ a map in $\Enr$ lifting $\alpha$. From now on we will write $\cat V$ for $\cat V'^\tensor$ and $\cat W$ for $\cat W'^\tensor$. Since $\Enr \to \Op_{/\Ass^\tensor}$ is a cocartesian fibration, there is a cocartesian edge over $\alpha$ with source $\cat V'^\tensor$, whose target is necessarily $\tau_\alpha(\cat V)$ (by definition of the transfer of enrichment via precisely these cocartesian edges). Let us denote this edge by $\gamma\colon \cat V \to \tau_\alpha(\cat V)$. By the definition of cocartesian edges, $\alpha'\colon \cat V \to \cat W$ factors uniquely over $\gamma$, which induces the desired $\cat W$-enriched functor $G_\alpha\colon \tau_\alpha(\cat V) \to \cat W$. By going through the construction of the functor $\chi$ in \cref{rslt:Heine-equivalence-of-models-for-enrichment} one checks that $\alpha'$ induces the expected map on enriched hom's, hence the same is true for $G_\alpha$.

It remains to prove that $G_\alpha = \enrHom[\cat W]_{\tau_\alpha(\cat V)}(\one_{\cat V}, \blank)$. By the enriched Yoneda lemma (see \cref{rslt:enriched-Yoneda-embedding}), a map from right to left is equivalently given by a map $\one_{\cat W} \to G_\alpha(\one_{\cat V}) = \alpha(\one_{\cat V})$. Such a map exists by the lax monoidal structure on $\alpha$. We can now prove the equivalence of the induced map of functors $\tau_\alpha(\cat V) \to \cat W$ after evaluating both sides on each object $X \in \tau_\alpha(\cat V)$; but it is clear that both functors give $\alpha(X)$, as desired.
\end{proof}

The construction of $G_\alpha$ in the previous result is functorial in $\alpha$ in a lax sense. We do not have the 2-categorical machinery available to state this result in the most general way, but we can formulate the following consequence:

\begin{lemma} \label{rslt:functoriality-of-functor-from-transfer-of-self-enrichment}
Let
\begin{equation}\begin{tikzcd}
	\cat V \arrow[d,"f",swap] \arrow[dr,"\alpha",""{name=alpha,below}]\\
	\cat V' \arrow[r,"\alpha'",swap] \arrow[Rightarrow,from=alpha,"\sigma"] & \cat W
\end{tikzcd}\end{equation}
be a diagram of closed monoidal categories $\cat V$, $\cat V'$ and $\cat W$, lax monoidal functors $f$, $\alpha$ and $\alpha'$, and a natural transformation $\sigma\colon \alpha \to \alpha' f$ of lax monoidal functors. Then there is an induced diagram
\begin{equation}\begin{tikzcd}
	\tau_\alpha(\cat V) \arrow[d,"\tilde f",swap] \arrow[dr,"G_\alpha",""{name=alpha,below}]\\
	\tau_{\alpha'}(\cat V') \arrow[r,"G_{\alpha'}",swap] \arrow[Rightarrow,from=alpha,"\tilde\sigma"] & \cat W
\end{tikzcd}\end{equation}
of $\cat W$-enriched categories, $\cat W$-enriched functors and a $\cat W$-enriched natural transformation $\tilde\sigma\colon G_\alpha \to G_{\alpha'} \tilde f$. More explicitly:
\begin{lemenum}
	\item $\tilde f$ is the composition of $\tau_\alpha \cat V \to \tau_{\alpha' f} \cat V$ (induced by $\sigma$) and $\tau_{\alpha'}(G_f)\colon \tau_{\alpha'}(\tau_f \cat V) \to \tau_{\alpha'} \cat V'$. It acts on objects in the same way as $f$ does.

	\item For every $X \in \tau_\alpha \cat V$ the morphism $\tilde\alpha(X)\colon G_\alpha(X) \to G_{\alpha'}(\tilde f(X))$ is the same as the one induced by $\alpha$ after identifying the objects in the categories of the upper and lower diagram above.
\end{lemenum}
\end{lemma}
\begin{proof}
This should follow from a 2-categorical version of cocartesian fibrations. As we do not have this machinery available here, we provide an argument via direct computation. By \cref{rslt:functor-from-transfer-of-self-enrichment} we can write $G_\alpha = \enrHom[\cat W]_{\tau_\alpha(\cat V)}(\one_{\cat V}, \blank)$, hence by the enriched Yoneda lemma (see \cref{rslt:enriched-Yoneda-embedding}), a morphism $G_\alpha \to G_{\alpha'} \tilde f$ is equivalently given by a morphism $\one_{\cat W} \to G_{\alpha'}(\tilde f(\one_{\cat V}))$. We compute
\begin{align}
	G_{\alpha'}(\tilde f(\one_{\cat V})) = G_{\alpha'}(f(\one_{\cat V})) = \enrHom[\cat W]_{\tau_{\alpha'}}(\one_{\cat V'}, f(\one_{\cat V})) = \alpha'(\enrHom[\cat V']_{\cat V'}(\one_{\cat V'}, f(\one_{\cat V}))) = \alpha' f(\one_{\cat V}).
\end{align}
Now for the desired map $\one_{\cat V} \to \alpha' f(\one_{\cat V})$ we can simply choose the one provided by the lax monoidal structure on $\alpha' f$. By unravelling the construction, we deduce (ii). Here we note that while the map $\one_{\cat V} \to \alpha' f(\one_{\cat V})$ does not depend on $\sigma$, the identification with a map $G_\alpha \to G_{\alpha'} \tilde f$ \emph{does} depend on it. 
\end{proof}

Combining the above constructions of $F_\alpha$ and $G_\alpha$, we end up with the following canonical decomposition of a lax monoidal functor into certain enriched functors. This construction will be very useful for studying the category of kernels in the main part of the paper.

\begin{proposition} \label{rslt:factor-map-of-monoidal-categories-into-enriched-maps}
Let $\alpha\colon \cat V \to \cat W$ be a lax monoidal functor of closed monoidal categories, which we view as enriched over themselves. Then $\alpha$ decomposes as
\begin{align}
	\cat V = \ul{\cat V} \xlongto{F_\alpha} \ul{\tau_\alpha(\cat V)} \xlongto{\ul{G_\alpha}} \ul{\cat W} = \cat W,
\end{align}
where $F_\alpha$ is the functor from \cref{rslt:functor-on-underlying-categories-of-transfer-of-enrichment} and $G_\alpha$ is the functor from \cref{rslt:functor-from-transfer-of-self-enrichment}.
\end{proposition}
\begin{proof}
The map $\alpha$ upgrades to a map $\alpha'$ in $\Enr$ lying over $\alpha$, as in the proof of \cref{rslt:functor-from-transfer-of-self-enrichment}, and it is clear that by applying the functor $\Enr \to \Cat$, $\cat C \mapsto \ul{\cat C}$ we get the functor $\alpha$ back. Now the functors $F_\alpha$ and $G_\alpha$ are obtained from a decomposition of $\alpha'$ into two maps $\cat V \to \tau_\alpha(\cat V)$ and $\tau_\alpha(\cat V) \to \cat W$ in $\Enr$ (the first one being a cocartesian map). This immediately implies the claim.
\end{proof}

So far we have considered the transfer of enrichment along an operad map $\cat V^\tensor \to \cat W^\tensor$ for \emph{fixed} $\cat V^\tensor$ and $\cat W^\tensor$. Our next result deals with a generalization of this construction to families. The result makes use of the notion of $\cat S$-families of non-symmetric operads, for which we refer the reader to \cite[Definition~2.54]{Heine.2023}. In the following result we slightly abuse notation and implicitly work with non-symmetric operads throughout (we do not know if \cref{rmk:non-symmetric-vs-symmetric-operads} generalizes to generalized operads).

\begin{lemma} \label{rslt:transfer-of-enrichment-in-families}
Let $\cat S$ be a category, $F\colon \cat S \to \Enr$ a functor and $\cat V^\tensor \to \Ass^\tensor \times \cat S$ the cocartesian $\cat S$-family of operads associated with the composition of $F$ with the projection $\Enr \to \Op_{/\Ass^\tensor}$. Let $\cat W^\tensor \to \Ass^\tensor \times \cat S$ be another cocartesian $\cat S$-family of operads and let $\alpha\colon \cat V^\tensor \to \cat W^\tensor$ be a map of $\cat S$-families of operads. Then there is a functor $\tau_\alpha(F)\colon \cat S \to \Enr$ such that for each $s \in \cat S$ we have
\begin{align}
	\tau_\alpha(F)(s) = \tau_{\alpha_s}(F(s)),
\end{align}
where $\alpha_s\colon \cat V^\tensor_s \to \cat W^\tensor_s$ denotes the fiber of $\alpha$ over $s$. Moreover, given a morphism $f\colon s \to t$ in $\cat S$, the associated morphism $\tau_{\alpha_s}(F(s)) \to \tau_{\alpha_t}(F(t))$ acts on objects in the same way as $F(f)$ does and acts on enriched morphisms via the natural transformation $w_f \alpha_s \to \alpha_t v_f$, where $v_f\colon \cat V^\tensor_s \to \cat V^\tensor_t$ and $w_f\colon \cat W^\tensor_s \to \cat W^\tensor_t$ are the induced functors.
\end{lemma}

\begin{remark}
The $\cat S$-families $\cat V$ and $\cat W$ of operads in \cref{rslt:transfer-of-enrichment-in-families} are required to be cocartesian and hence correspond to functors $\cat S \to \Op_{/\Ass^\tensor}$. However, it is crucial for our applications that the map $\alpha$ is only a map of $\cat S$-families of operads (i.e.\ a map of generalized operads over $\Ass^\tensor \times \cat S$), so it does \emph{not} need to preserve cocartesian edges. In particular, $\alpha$ does not necessarily correspond to a natural transformation of functors $\cat S \to \Op_{/\Ass^\tensor}$.
\end{remark}

\begin{proof}[Proof of \cref{rslt:transfer-of-enrichment-in-families}]
We can replace $\Enr$ with $\PreEnr$ for this proof, as we can always pass back to $\Enr$ using the left adjoint $L$ to the inclusion from \cref{rslt:categorification-of-enriched-precategories}. We thus view $F$ as a functor $\cat S \to \PreEnr = \Alg^{\Ani \times \Op_{/\Ass^\tensor}}(\Quiv^\tensor)$, which by definition of $\Alg^{(\blank)}$ (see \cite[Notation~2.56]{Heine.2023}) corresponds to a map of generalized operads
\begin{align}
	F'\colon \Ass^\tensor \times \cat S \to \Quiv^\tensor
\end{align}
over $\Ass^\tensor \times \Ani \times \Op_{/\Ass^\tensor}$. Here the implicit map $\cat S \to \Ani \times \Op_{/\Ass^\tensor}$ sends $s \mapsto (F(s)^\simeq, \cat V_s^\tensor)$, i.e.\ it is the composition of $F$ with the projection from $\PreEnr$. By \cite[Remark~4.18]{Heine.2023} we have $\Quiv^\tensor \times_{\Ani \times \Op_{/\Ass^\tensor}} \cat S = \Quiv^{\cat S}_{\cat X}(\cat V)^\tensor$, where we use the notation from \cite[Notation~4.17]{Heine.2023} and we write $\cat X \to \cat S$ for the unstraightening of the functor $\cat S \to \Ani$, $s \mapsto F(s)^\simeq$. Therefore $F'$ corresponds to a map of generalized operads
\begin{align}
	F''\colon \Ass^\tensor \times \cat S \to \Quiv^{\cat S}_{\cat X}(\cat V)^\tensor
\end{align}
over $\cat S$. But $\Quiv^{\cat S}_{\cat X}(\blank)^\tensor$ is functorial in generalized operads over $\cat S$, hence $\alpha$ induces a map $\Quiv^{\cat S}_{\cat X}(\cat V)^\tensor \to \Quiv^{\cat S}_{\cat X}(\cat W)^\tensor$. By composing $F''$ with this map and then going the above equivalent descriptions of functors backwards, we arrive at the desired functor $\tau_\alpha(F)$.
\end{proof}

\section{2-categories} \label{sec:2cat}

The main object of study in this paper is the 2-category of kernels $\cat K_\D$ associated with a 6-functor formalism $\D$. In the following we provide some general tools for studying 2-categories, many of which we have not found in the literature.

\subsection{Definition and examples} \label{sec:2cat.def}

As was hinted at before, we define 2-categories simply as $\Cat$-enriched categories. Let us introduce the following notation and terminology:

\begin{definition}
\begin{defenum}
	\item\label{def:2-category-of-2-categories} A \emph{2-category} is a $\Cat$-enriched category and a \emph{2-functor} between 2-categories is a $\Cat$-enriched functor, where we equip $\Cat$ with the cartesian monoidal structure from \cref{ex:cartesian-monoidal-structure}. We denote by
	\begin{align}
		\TwoCat \coloneqq \Enr_\Cat
	\end{align}
	the category of 2-categories and 2-functors.

	\item\label{def:Hom-category-in-2-Cat} For a 2-category $\cat C \in \TwoCat$ and objects $X, Y \in \cat C$ we denote
	\begin{align}
		\Fun_{\cat C}(X, Y) \coloneqq \enrHom[\Cat]_{\cat C}(X, Y)
	\end{align}
	the morphism category from $X$ to $Y$ in $\cat C$.
\end{defenum}
\end{definition}

\begin{remark}
In the literature, 2-categories are often referred to as $(\infty,2)$-categories, but we stick with the simpler naming scheme introduced above. There are many different models for 2-categories, see e.g.\ \cite[\S2]{Haugseng.2021}, and they are known to be equivalent. We are therefore free to choose whichever model we prefer, so we stick to the one based on enriched categories.
\end{remark}

We can apply the constructions from \cref{sec:enr} to the specific setting of $\cat V = \Cat$, providing us with the following useful definitions:

\begin{definition}
\begin{defenum}
	\item We denote $\ul{(\blank)}\colon \TwoCat \to \Cat$ the functor sending a 2-category $\cat C$ to its underlying category $\ul{\cat C}$. Then $\ul{\cat C}$ has the same objects as $\cat C$ and for $X, Y \in \cat C$ we have $\Hom_{\ul{\cat C}}(X, Y) = \Fun_{\cat C}(X, Y)^\simeq$, cf. \cref{rslt:underlying-category-via-transfer-of-enrichment}.

	\item A 2-functor $F\colon \cat C_1 \to \cat C_2$ of 2-categories is called \emph{fully faithful} if it is so as a $\Cat$-enriched functor, i.e.\ if for all $X, Y \in \cat C_1$ the induced functor $\Fun_{\cat C_1}(X, Y) \isoto \Fun_{\cat C_2}(F(X), F(Y))$ is an equivalence. We say that $F$ is \emph{2-fully faithful} if the induced functors $\Fun_{\cat C_1}(X, Y) \injto \Fun_{\cat C_2}(F(X), F(Y))$ are fully faithful. We say that $F$ is \emph{essentially surjective} if it is so on underlying categories. Note that $F$ is an equivalence of 2-categories if and only if it is fully faithful and essentially surjective, see \cref{rslt:enriched-surjective-plus-fully-faithful-implies-isom}.

	\item For a 2-category $\cat C$ we denote by $\cat C^\op \in \TwoCat$ the \emph{(1-)opposite 2-category}, cf. \cref{def:opposite-enriched-category}. It has the same objects as $\cat C$, but for $X, Y \in \cat C$ we have $\Fun_{\cat C^\op}(X, Y) = \Fun_{\cat C}(Y, X)$.

	\item For a 2-category $\cat C$ we denote
	\begin{align}
		\cat C^\co \coloneqq \tau_{(\blank)^\op}(\cat C) \in \TwoCat
	\end{align}
	and call it the \emph{2-opposite 2-category} of $\cat C$. Here the right-hand side denotes the transfer of enrichment along the symmetric monoidal equivalence $(\blank)^\op\colon \Cat \to \Cat$. Explicitly, $\cat C^\co$ has the same objects as $\cat C$, but for $X, Y \in \cat C$ we have $\Fun_{\cat C^\co}(X, Y) = \Fun_{\cat C}(X, Y)^\op$. Note that by \cref{rslt:transfer-of-enrichment-commutes-with-op}, the endofunctors $(\blank)^\op$ and $(\blank)^\co$ of $\TwoCat$ commute.

	\item \label{def:inclusion-of-Cat-into-TwoCat} The underlying category functor $\ul{(\blank)}$ has a fully faithful left adjoint $\Cat \injto \TwoCat$ which sends a category $\cat C$ to the 2-category with the same objects, where the hom categories are simply the hom anima of $\cat C$. Indeed, this left adjoint can be defined as the transfer of enrichment along $\Ani \injto \Cat$ (cf. \cref{ex:Ani-enriched-category-same-as-category}). To show that this is a left adjoint, we recall that $\ul{(\blank)}$ is the transfer of enrichment along the functor $\Cat \to \Ani$, $\cat C \mapsto \cat C^\simeq$ (see \cref{rslt:underlying-category-via-transfer-of-enrichment}), hence the adjunction of the two functors $\Ani \rightleftarrows \Cat$ provides the necessary unit and counit for the desired adjunction.

	\item \label{def:2-subcategory-generated-by-objs-and-homs} Let $\cat C$ be a 2-category and suppose we have chosen a collection $\cat S$ of isomorpism classes of objects in $\cat C$ and for all $X, Y \in \cat S$ a collection $\cat S_{X,Y}$ of isomorphism classes in $\Fun_{\cat C}(X, Y)$ such that $\cat S_{X,X}$ contains the identity and the $\cat S_{X,Y}$'s are stable under composition. Then there is a 2-full subcategory $\cat C_0 \subseteq \cat C$ with objects $\cat S$ and morphism categories the full subcategories spanned by $\cat S_{X,Y}$. We call $\cat C_0$ the 2-category \emph{generated} by the above data. Its construction can easily be performed in the Lurie model for enrichment.

	\item \label{def:final-2-category} The category $\TwoCat$ has all limits (by \cref{rslt:limits-of-enriched-categories}) and they have an explicit description on objects and hom categories. In particular $\TwoCat$ has a final object, which we denote by $* \in \TwoCat$. It coincides with the trivial $\Cat$-enriched category from \cref{ex:trivial-enriched-category} and with the image of $* \in \Cat$ under the embedding $\Cat \injto \TwoCat$ (to see these identifications, use e.g.\ the explicit description of limits in \cref{rslt:limits-of-enriched-categories}). By \cref{rslt:universal-property-of-trivial-enriched-category}, for every 2-category $\cat C$, a 2-functor $* \to \cat C$ is the same as an object in $\cat C$.

	\item For a 2-category $\cat C$ we denote by
	\begin{align}
		\cat P_2(\cat C) \coloneqq \PSh[\Cat](\cat C) = \TwoFun(\cat C^\op, \Cat)
	\end{align}
	the $\Cat$-enriched Yoneda category (see \cref{def:enriched-yoneda-category}). We will see in \cref{rslt:2-yoneda} that $\cat P_2(\cat C)$ can naturally be enhanced to a 2-category.
\end{defenum}
\end{definition}

By \cref{def:inclusion-of-Cat-into-TwoCat} we can view every category implicitly as a 2-category, which we will often do in this paper. Moreover, a 2-functor from a category $\cat C_1$ to a 2-category $\cat C_2$ is the same as a functor from $\cat C_1$ to $\ul{\cat C_2}$.

Many categories that appear in practice come equipped with a natural 2-categorical enhancement. In the following we will provide examples of such enhancements for many categories that appear in this paper.

\begin{example}\label{ex:2-categorical-enhancement-of-Cat}
The category $\Cat$ of categories, equipped with the cartesian monoidal structure, is a closed symmetric monoidal category and hence enriched in itself by \cref{ex:monoidal-category-enriched-in-itself}. Thus $\Cat$ is naturally a 2-category with $\Fun_{\Cat}(\cat C_1, \cat C_2) = \Fun(\cat C_1, \cat C_2)$ for all $\cat C_1, \cat C_2 \in \Cat$.
\end{example}

\begin{example} \label{rslt:2-categorical-enhancement-of-Cat-over-S}
Let $\cat S$ be a fixed category. Then $\blank \times \cat S\colon \Cat \to \Cat_{/\cat S}$ defines a symmetric monoidal functor for the cartesian monoidal structures and in particular equips $\Cat_{/\cat S}$ with a $\Cat$-linear structure. One checks easily that this $\Cat$-linear structure has all morphism objects and hence equips $\Cat_{/\cat S}$ with a 2-categorical enhancement such that for $\cat C_1, \cat C_2 \in \Cat_{/\cat S}$,
\begin{align}
	\Fun_{\Cat_{/\cat S}}(\cat C_1, \cat C_2) = \Fun_{\cat S}(\cat C_1, \cat C_2)
\end{align}
is the category of functors $\cat C_1 \to \cat C_2$ over $\cat S$. For every functor $\cat T \to \cat S$, the induced base-change functor $\blank \times_{\cat S} \cat T$ is symmetric monoidal for the cartesian structures and in particular $\Cat$-linear and therefore induces a 2-functor $\Cat_{/\cat S} \to \Cat_{/\cat T}$. Its left adjoint is the forgetful functor sending $\cat C \to \cat T$ to $\cat C \to \cat S$ and using \cref{rslt:monoidal-linear-left-adjoints} one sees easily that this functor is $\Cat$-linear and thus a 2-functor.
\end{example}

Using similar techniques one can construct many more 2-categorical enhancements and upgrade functors to 2-functors. In the following we will focus on the category of operads and of enriched categories. We will always implicitly equip $\Cat_{/\cat S}$ with the 2-categorical enhancement from \cref{rslt:2-categorical-enhancement-of-Cat-over-S}.

\begin{proposition} \label{rslt:2-categorical-enhancement-of-Op}
Let $\cat O^\tensor$ be an operad. Then $\Op_{/\cat O^\tensor}$ is naturally equipped with a 2-categorical enhancement such that for all $\cat V^\tensor, \cat W^\tensor \in \Op_{/\cat O^\tensor}$ we have
\begin{align}
	\Fun_{\Op_{/\cat O^\tensor}}(\cat V^\tensor, \cat W^\tensor) = \Alg_{\cat O}(\cat V, \cat W).
\end{align}
The forgetful functor $\Op_{/\cat O^\tensor} \to \Cat_{/\cat O^\tensor}$ upgrades to a 2-fully faithful 2-functor. Moreover, for every map of operads $\cat O'^\tensor \to \cat O^\tensor$ the base-change functor $\blank \times_{\cat O^\tensor} \cat O'^\tensor\colon \Op_{/\cat O^\tensor} \to \Op_{/\cat O'^\tensor}$ and the forgetful functor $\Op_{/\cat O'^\tensor} \to \Op_{/\cat O^\tensor}$ upgrade to 2-functors.
\end{proposition}
\begin{proof}
By \cref{def:2-subcategory-generated-by-objs-and-homs} we can construct the 2-categorical enhancement on $\Op_{/\cat O^\tensor}$ by the requirement that the forgetful functor $\Op_{/\cat O^\tensor} \to \Cat_{/\cat O^\tensor}$ is a 2-fully faithful 2-functor. However, we prefer the following more natural construction. Let $\GOp \subset \Cat_{/\Fin_*}$ denote the non-full subcategory of \emph{generalized} operads, as defined in \cite[Definition~2.3.2.1]{HA}. Then $\Op \subseteq \GOp$ is a full subcategory (see \cite[Notation~2.3.2.7]{HA}), so it is enough to construct a 2-categorical enhancement on $\GOp$. On the other hand, one checks in the same way as in \cref{rslt:limits-of-operads} that $\GOp$ has all small limits and the forgetful functor $\GOp \to \Cat_{/\Fin_*}$ preserves these limits. In particular the functor $\cat C \mapsto \cat C \times \Fin_*$ defines a product preserving functor $\Cat \to \GOp$ (see \cite[Proposition~2.3.2.9]{HA}) and hence equips $\GOp$ with a $\Cat$-linear structure. One checks that this $\Cat$-linear structure has enough morphism objects given by $\Alg$ (by \cref{rslt:2-categorical-enhancement-of-Cat-over-S} this reduces to basic properties of inert morphisms) and hence defines the desired 2-categorical enhancement of $\GOp$. By construction the functor $\GOp \to \Cat_{/\Fin_*}$ is $\Cat$-linear and in particular a 2-functor.

It is now formal to pass to slice categories, i.e.\ categories of the form $\GOp_{/\cat O^\tensor}$, using the same technique as in \cref{rslt:2-categorical-enhancement-of-Cat-over-S}.
\end{proof}

\begin{proposition} \label{rslt:2-categorical-enhancement-of-Enr}
Let $\cat V$ be an operad over $\Ass^\tensor$. Then the category $\Enr_{\cat V}$ of $\cat V$-enriched categories admits a 2-categorical enhancement satisfying the following properties:
\begin{propenum}
	\item For $\cat V$-enriched categories $\cat C_1, \cat C_2 \in \Enr_{\cat V}$ we have
	\begin{align}
		\Fun_{\Enr_{\cat V}}(\cat C_1, \cat C_2) = \enrFun[\cat V](\cat C_1, \cat C_2),
	\end{align}
	where $\enrFun[\cat V]$ denotes the category of $\cat V$-enriched functors from \cref{def:enriched-functors}.

	\item The functor $\ul{(\blank)}\colon \Enr_{\cat V} \to \Cat$ upgrades to a 2-functor, which on morphism categories induces the forgetful functors from \cref{def:forgetful-functor-on-enriched-functor-categories}.
\end{propenum}
\end{proposition}
\begin{proof}
By construction, $\Enr_{\cat V}$ is a full subcategory of $\Op_{/\LM^\tensor} \times_{\Op_{/\Ass^\tensor}} \{ \cat V^\tensor \}$. But by \cref{rslt:2-categorical-enhancement-of-Op,def:final-2-category} all of the categories on the right admit a 2-categorical enhancement and the functors to $\Op_{/\Ass^\tensor}$ upgrade to 2-functors. We can thus form the same fiber product in $\TwoCat$ to obtain the desired 2-category, which by the explicit description of limits in \cref{rslt:limits-of-enriched-categories} satisfies (i). For (ii) we note that the functor $\ul{(\blank)}$ can be constructed as the composition of the forgetful functor $\Enr_{\cat V} \to \Cat_{/\LM^\tensor}$ and the base-change functor $\blank \times_{\LM^\tensor} \{ \mathfrak m \}$, both of which are 2-functors by \cref{rslt:2-categorical-enhancement-of-Op} and \cref{rslt:2-categorical-enhancement-of-Cat-over-S}, respectively.
\end{proof}

As a special case of \cref{rslt:2-categorical-enhancement-of-Enr} we can naturally view $\TwoCat$ as a 2-category, where the hom-categories are given by $\TwoFun$. In the following we will always implicitly equip $\Enr_{\cat V}$ with its 2-categorical enhancement.

All of the 2-categorical enhancements constructed so far ultimately relied on the following construction: Suppose we have a category $\cat C$ which has finite products, and a functor $\Cat \to \cat C$ that preserves finite products; then this functor induces a $\Cat$-linear structure on $\cat C$ which (under suitable assumptions on $\cat C$) defines the desired 2-categorical enhancement of $\cat C$. By \cref{rslt:uniqueness-of-cartesian-left-actions} such a 2-categorical enhancement is uniquely characterized by the property that $\Cat$ acts via products. As an application of this result, we obtain a nice way of describing the 2-categorical Yoneda lemma:

\begin{proposition} \label{rslt:2-yoneda}
For every small 2-category $\cat C$ there is a natural 2-categorical enhancement of $\cat P_2(\cat C)$ such that the Yoneda embedding upgrades to a fully faithful 2-functor
\begin{align}
	\cat C \injto \cat P_2(\cat C), \qquad X \mapsto \Fun(\blank, X).
\end{align}
Explicitly, the $\Cat$-enrichment on $\cat P_2(\cat C)$ is $\Cat$-linear and comes via \cref{def:induced-cartesian-left-action} from the functor $\Cat \to \cat P_2(\cat C)$ sending each category $V \in \Cat$ to the constant functor $[X \mapsto V] \in \cat P_2(\cat C)$.
\end{proposition}
\begin{proof}
In the setting of an arbitrary operad $\cat V^\tensor$ over $\Ass^\tensor$ and a $\cat V$-enriched category $\cat C$, \cite[Notation~3.84]{Heine.2023} and \cite[Proposition~3.86(5)]{Heine.2023} introduce a left action of $\cat V^\tensor$ on $\cat P_{\cat V}(\cat C) = \enrFun[\cat V](\cat C^\op, \cat V)$ (here we view $\cat C^\op$ as a weak right $\cat V$-module and $\cat V$ as a $\cat V$-bimodule). By \cite[Theorem~10.11]{Heine.2023} the Yoneda embedding from \cref{rslt:enriched-Yoneda-embedding} upgrades to a fully faithful embedding $\cat C \injto \cat P_{\cat V}(\cat C)$ of categories with (weak) left $\cat V^\tensor$-action (\cite[Theorem~10.11]{Heine.2023} is a bit difficult to parse, but it does indeed contain the claimed Yoneda embedding by combining the previous results in \cite[\S10]{Heine.2023} and in particular using the definition of opposite enriched category from \cite[Notation~10.2]{Heine.2023}).

We apply the results from the previous paragraph to the case $\cat V^\tensor = \Cat^\times$. We thus obtain a left $\Cat$-action on $\cat P_2(\cat C)$. We claim that this $\Cat$-action is cartesian. First note that by \cref{rslt:limits-in-enrFun} the category $\cat P_2(\cat C)$ has all small limits and in particular all finite products. In the case $\cat C = *$ we have $\cat P_2(*) = \Cat$ with the obvious $\Cat$-linear structure (this follows directly from the definitions or from the universal property in \cite[Theorem~5.1]{Heine.2023}). Moreover, given an object $X \in \cat C$, i.e.\ a functor $* \to \cat C$ then precomposition with this functor induces a $\Cat$-linear functor $\cat P_2(\cat C) \to \Cat$, $f \mapsto f(X)$, by functoriality of the left action from \cite[Notation~3.84]{Heine.2023}. This means that for $X \in \cat C$, $V \in \Cat$ and $f\in \cat P_2(\cat C)$ there is a natural isomorphism $(V \tensor f)(X) = V \times f(X)$. From this we deduce that $V \tensor f = (V \tensor *) \times f$, because this can be checked element-wise (recall that the forgetful functor $\TwoFun(\cat C^\op, \Cat) \to \Fun(\ul{\cat C}^\op, \Cat)$ is conservative, see \cref{def:enriched-functors}). In other words, the $\Cat$-action on $\cat P_2(\cat C)$ is indeed cartesian. From \cref{rslt:uniqueness-of-cartesian-left-actions} we deduce that the $\Cat$-linear structure on $\cat P_2(\cat C)$ is induced by the functor $\Cat \to \cat P_2(\cat C)$, $V \mapsto V \tensor *$, and $V \tensor *$ is the constant functor $X \mapsto V$. Finally, we notice that for fixed $f \in \cat P_2(\cat C)$ the functor $\Cat \to \cat P_2(\cat C)$, $V \mapsto V \tensor f$, commutes with colimits (by \cref{rslt:lim-and-colim-in-enrFun} this reduces to the underlying functor category and then to the claim that $\times$ preserves colimits in each argument in $\Cat$), so by the adjoint functor theorem (see \cite[Remark~5.5.2.10]{HTT}) this functor has a right adjoint, providing enough morphism objects for $\cat P_2(\cat C)$. Thus $\cat P_2(\cat C)$ is $\Cat$-enriched, as desired.
\end{proof}

\begin{remark}
A similar argument as in \cref{rslt:2-yoneda} gives the same result for every cartesian symmetric monoidal presentable category $\cat V$ in place of $\Cat$ (assuming that $\times\colon \cat V \times \cat V \to \cat V$ preserves colimits in each argument).
\end{remark}



\subsection{Adjunctions} \label{sec:2cat.adj}

A fundamental concept of 2-categories is that of adjoint morphisms, which also play a major role in the study of 6-functor formalisms through the associated 2-category of kernels. In the following we introduce the notion of adjoint morphisms and provide some basic properties. In the next subsection we will show that passing to the adjoint functor is a \enquote{natural} operation.

\begin{definition} \label{def:adjoint-morphisms}
Let $\cat C$ be a 2-category with objects $X, Y \in \cat C$. We say that a morphism $f\colon Y \to X$ is a \emph{left adjoint} in $\cat C$ if there is a morphism $g\colon X \to Y$, called the \emph{right adjoint of $f$}, and 2-morphisms $\eta\colon \id_Y \to gf$ and $\varepsilon\colon fg \to \id_X$, called the \emph{unit} and \emph{counit} respectively, such that there are commutative diagrams
\begin{equation}
	\begin{tikzcd}
		f \arrow[r,"f \eta"] \arrow[dr,"\id",swap] & fgf \arrow[d,"\varepsilon f"]\\
		& f
	\end{tikzcd}
	\qquad
	\begin{tikzcd}
		g \arrow[r,"\eta g"] \arrow[dr,"\id",swap] & gfg \arrow[d,"g \varepsilon"]\\
		& g
	\end{tikzcd}
\end{equation}
in $\Fun_{\cat C}(Y, X)$ and $\Fun_{\cat C}(X, Y)$, respectively. We similarly define the notion of being a \emph{right adjoint}. We denote by
\begin{align}
	\Fun^L_{\cat C}(Y, X), \Fun^R_{\cat C}(Y, X) \subseteq \Fun_{\cat C}(Y, X)
\end{align}
the full subcategories of left adjoint and right adjoint morphisms, respectively.
\end{definition}

\begin{example}
A functor $F\colon \cat C \to \cat D$ of categories is a left adjoint functor if and only if it is a left adjoint morphism in the 2-category $\Cat$.
\end{example}

\begin{example} \label{rslt:adj-in-*-V-same-as-dualizable}
Let $\cat V^\tensor$ be a symmetric monoidal category and let $*_{\cat V}$ be the 2-category with a single object $*$ and $\End(*) = \cat V$ (see \cref{ex:algebra-is-enriched-category}). By comparing \cref{def:dualizable-object,def:adjoint-morphisms} we see that a morphism in $*_{\cat V}$ is left adjoint if and only if it is dualizable as an object in $\cat V$, and if this is the case then the associated right adjoint corresponds to the dual. We can thus see adjoint morphisms as a generalization of dualizable objects.
\end{example}

Note that the datum of an adjunction in a 2-category $\cat C$ only depends on the underlying ordinary 2-category $\cat C$, i.e.\ it does not depend on the higher categorical structure of $\cat C$. This allows one to prove basic properties of adjunctions \enquote{by hand}. It is sometimes useful to note that the triangle identities in the definition of adjoint morphisms can be relaxed:

\begin{lemma} \label{rslt:adjoints-from-weak-triangle-identity}
Let $f\colon Y \rightleftarrows X\noloc g$ be morphisms in a 2-category $\cat C$ and let $\eta\colon \id_Y \to gf$ and $\varepsilon\colon fg \to \id_X$ be 2-morphisms such that the composed 2-morphisms
\begin{align}
	f \to fgf \to f, \qquad g \to gfg \to g
\end{align}
are isomorphisms. Then there is a 2-morphism $\eta'\colon \id_Y \to gf$ such that $f$ is left adjoint to $g$ with unit $\eta'$ and counit $\varepsilon$.
\end{lemma}
\begin{proof}
This is \cite[Lemma~5.11]{Scholze:Six-Functor-Formalism}; for the convenience of the reader we repeat the argument. Let $\eta'$ be the composition of $\eta$ with the inverse of the isomorphism $gf \isoto gf$ induced by the second isomorphism $g \isoto g$ in the claim. We similarly construct $\eta''$ by using the first isomorphism in the claim. Then with $\eta'$ in place of $\eta$, the composition $g \to gfg \to g$ is the identity, and with $\eta''$ in place of $\eta$ the composition $f \to fgf \to f$ is the identity. It is therefore enough to show that $\eta' \isom \eta''$. But
\begin{align}
	[1 \xto{\eta'} gf] &\isom [1 \xto{\eta'} gf \xto{gf \eta''} gfgf \xto{g\varepsilon f} gf]\\
	&\isom [1 \xto{\eta''} gf \xto{\eta'gf} gfgf \xto{g\varepsilon f} gf]\\
	&\isom [1 \xto{\eta''} gf],
\end{align}
as desired.
\end{proof}

The next criterion often allows to reduce questions about adjoint morphisms in $\cat C$ to questions about adjoint functors in $\Cat$ (see \cref{rslt:pointwise-criterion-for-adjunction} below for a more powerful result of similar flavor):

\begin{lemma} \label{rslt:adjunctions-in-2-category-via-Hom-categories}
Let $f\colon Y \rightleftarrows X\noloc g$ be morphisms in a 2-category $\cat C$ and let $\eta\colon \id_Y \to gf$ be a 2-morphism. Then the following are equivalent:
\begin{lemenum}
	\item $f$ is left adjoint to $g$ with unit $\eta$.
	
	\item \label{rslt:adjunctions-in-2-category-induce-adjunction-on-Hom-categories} For every $Z \in \cat C$, the composition with $f$ and $g$ induces an adjunction
	\begin{align}
		f_*\colon \Fun_{\cat C}(Z, Y) \rightleftarrows \Fun_{\cat C}(Z, X)\noloc g_*
	\end{align}
	with unit induced by $\eta$.

	\item \label{rslt:adjunctions-in-2-category-via-Hom-category-mates} For every $Z \in \cat C$ and morphisms $h\colon Z \to Y$ and $h'\colon Z \to X$, the natural map
	\begin{align}
		\Hom(fh, h') \isoto \Hom(h, gh'),
	\end{align}
	induced by $\eta$, is an isomorphism of anima.
\end{lemenum}
Moreover, in (ii) and (iii) it is enough to require the statement for $Z \in \{ X, Y \}$.
\end{lemma}
\begin{proof}
By applying the Yoneda embedding $\cat C^\op \injto \cat P_2(\cat C^\op)$ from \cref{rslt:2-yoneda} to an object $Z \in \cat C^\op$ we see that the functor $\Fun(Z, \blank)\colon \cat C \to \Cat$ is naturally a 2-functor. As such, it preserves adjunctions (as it preserves all the data and conditions in \cref{def:adjoint-morphisms}, see also \cref{rslt:2-functors-preserve-adjunctions} below). Therefore, (i) implies (ii). Moreover, (ii) and (iii) are equivalent by the usual criterion for adjoint functors (see e.g.\ \cite[\href{https://kerodon.net/tag/02FX}{Tag 02FX}]{kerodon}).

It remains to prove that (iii) implies (i), so assume that (iii) is satisfied. Let $\varepsilon\colon fg \to \id_X$ be a preimage of $\id_g$ under the natural isomorphism $\Hom(fg,\id) \isoto \Hom(g,g)$ from (iii). We claim that the tuple $(f,g,\eta, \varepsilon)$ provides the necessary data for an adjunction between $f$ and $g$. To see this, we need to verify the commutativity of the two triangles in \cref{def:adjoint-morphisms}.

The right triangle commutes by the definition of $\varepsilon$. To verify that the left triangle commutes, we write $\varphi \coloneqq \varepsilon f \circ f\eta \colon f\to f$ and have to show that $\varphi = \id_f$. Equivalently, it suffices to check that the images of $\varphi$ and $\id_f$ under the isomorphism $\Hom(f,f) \isoto \Hom(\id_Y, gf)$ from (iii) coincide. Unraveling the definitions, we have to check that the diagram
\begin{equation}
\begin{tikzcd}
\id_Y \ar[r,"\eta"] \ar[d,equals] & gf \ar[r,"gf\eta"] & gfgf \ar[r,"g\varepsilon f"] & gf\\
\id_Y \ar[rr,"\eta"'] & & gf \ar[u,"\eta gf"] \ar[ur,equals]
\end{tikzcd}
\end{equation}
commutes. But the square obviously commutes and the triangle commutes by the definition of $\varepsilon$ as argued above. We conclude that $(f,g,\eta,\varepsilon)$ is an adjunction as desired.
\end{proof}

The next result contains a list of easy-to-prove properties of adjunctions in a 2-category $\cat C$ which we will use throughout the paper.

\begin{lemma}
Let $\cat C$ be a 2-category.
\begin{lemenum}
	\item A morphism $f$ in $\cat C$ is left adjoint if and only if it is right adjoint as a morphism in $\cat C^\co$.

	\item A morphism $f\colon Y \to X$ in $\cat C$ is left adjoint if and only if the induced morphism $f^\op\colon X \to Y$ in $\cat C^\op$ is right adjoint.

	\item \label{rslt:adjoints-stable-under-composition} If $f$ and $g$ are composable left adjoint morphisms in a 2-category $\cat C$ then also $fg$ is left adjoint.

	\item \label{rslt:isomorphisms-are-adjoint} Every isomorphism in $\cat C$ is both a left and a right adjoint morphism.

	\item \label{rslt:adjoint-morphisms-are-unique} Let $f$ be a left adjoint with two right adjoints $g_1, g_2$. Then the composite $g_1 \to g_2fg_1 \to g_2$ is an isomorphism. In other words, any two right adjoint morphisms of a fixed left adjoint morphism are canonically isomorphic.
\end{lemenum}
\end{lemma}
\begin{proof}
For (i) we note that the unit and counit of an adjunction in $\cat C$ provide the counit and unit of an adjunction in $\cat C^\co$. A similar observation proves (ii).

For (iii) let $f\colon Y \to X$ and $f'\colon Z \to Y$ be left adjoint morphisms in $\cat C$ with right adjoints $g\colon X \to Y$ and $g'\colon Y \to Z$. We claim that $ff'$ is left adjoint with right adjoint $g'g$. Indeed, if $\eta\colon \id_Y \to gf$ and $\eta'\colon \id_Z \to g'f'$ are the adjunction units then the unit of the adjunction for $ff'$ and $g'g$ is the composite
\begin{align}
	\eta''\colon \id_Z \xto{\eta'} g'f' \xto{g'\eta f'} g'gf'f.
\end{align}
This follows at once from the criterion (iii) in \cref{rslt:adjunctions-in-2-category-via-Hom-categories}.

For (iv) we note that if $f\colon Y \to X$ is an isomorphism with inverse $g\colon X \to Y$ then there are isomorphisms $\eta\colon \id_X \isoto gf$ and $\varepsilon\colon fg \isoto \id_Y$. It follows immediately from \cref{rslt:adjoints-from-weak-triangle-identity} that these data provide the desired adjunction.

For (v), let $f\colon Y\rightleftarrows X\noloc g_i$, $i=1,2$, be adjunctions with unit $\eta_i\colon \id_Y \to g_if$ and counit $\varepsilon_i \colon fg_i \to \id_X$. Then an easy computation shows that, for all $h\colon Z\to Y$, $h'\colon Z\to X$, the composite of the isomorphisms (see \cref{rslt:adjunctions-in-2-category-via-Hom-category-mates})
\begin{align}
\Hom(h,g_1h') \isoto \Hom(fh,h') \isoto \Hom(h,g_2h')
\end{align}
is given by composition with $g_1 \xto{\eta_2g_1} g_2fg_1 \xto{g_2\varepsilon_1} g_1$, which by the Yoneda lemma is an isomorphism.
\end{proof}

\Cref{rslt:adjunctions-in-2-category-via-Hom-category-mates} is a special case of the classical \enquote{mate correspondence} (see \cite[Proposition~2.1]{Kelly-Street.1974}), which we now describe in full generality:

\begin{lemma}\label{rslt:mate-correspondence}
Let $\cat C$ be a 2-category. Let $f\colon A\rightleftarrows B \noloc u$ and $f'\colon A'\rightleftarrows B' \noloc u'$ be adjunctions in $\cat C$ with unit-counit pairs $(\eta, \varepsilon)$ and $(\eta',\varepsilon')$, respectively; we abbreviate the adjunctions by $(f,u,\eta,\varepsilon)$ and $(f',u',\eta',\varepsilon')$.
Let $a\colon A\to A'$ and $b\colon B\to B'$ be two morphisms. 
\begin{enumerate}[(i)]
\item There is a natural bijection
\begin{align}
	\rho\colon \Hom_{\Fun_{\cat C}(A,B')} (f'a, bf) & \quad \isotofrom \quad \Hom_{\Fun_{\cat C}(B,A')} (au, u'b) \noloc \lambda\\
	\varphi & \quad \mapsto \quad (u'b\varepsilon) \comp (u'\varphi u) \comp (\eta' au)\\
	(\varepsilon' bf) \comp (f'\psi f) \comp (f'a\eta)  & \quad \mapsfrom \quad \psi
\end{align}
\item Passage from $\cat C$ to $\cat C^{\op}$ interchanges $\rho$ and $\lambda$. In other words, $(u^{\op}, f^{\op}, \varepsilon^{\op}, \eta^{\op})$ is an adjunction in $\cat C^{\op}$, and the diagram
\begin{equation}
\begin{tikzcd}
\Hom_{\Fun_{\cat C}(A,B')} (f'a,bf) \ar[r,shift left,"\rho"] \ar[d,equals] & \Hom_{\Fun_{\cat C}(B,A')} (au, u'b) \ar[l,shift left, "\lambda"] \ar[d,equals] \\
\Hom_{\Fun_{\cat C^{\op}}(B',A)}(a^{\op} f'^{\op}, f^{\op} b^{\op}) \ar[r,shift left,"\lambda^{\op}"] & \Hom_{\Fun_{\cat C^{\op}}(A',B)} (u^{\op} a^{\op}, b^{\op} u'^{\op}) \ar[l,shift left, "\rho^{\op}"]
\end{tikzcd}
\end{equation}
commutes.

\item Every 2-functor $\theta\colon \cat C \to \cat K$ preserves the mate correspondence. In other words, one has adjunctions $(\theta(f),\theta(u), \theta(\eta), \theta(\varepsilon))$ and $(\theta(f'), \theta(u'), \theta(\eta'), \theta(\varepsilon'))$ in $\cat K$, and the diagram
\begin{equation}
\begin{tikzcd}
\Hom_{\Fun_{\cat C}(A,B')} (f'a, bf) \ar[r,shift left,"\rho"] \ar[d,"\theta"] & \Hom_{\Fun_{\cat C}(B,A')}(au, u'b) \ar[d,"\theta"] \ar[l,shift left,"\lambda"] \\
\Hom_{\Fun_{\cat K}(\theta A, \theta B')} \bigl(\theta(f')\theta(a), \theta(b)\theta(f)\bigr) \ar[r,shift left,"\rho"] & \Hom_{\Fun_{\cat K}(\theta B, \theta A')} \bigl(\theta(a) \theta(u), \theta(u') \theta(b)\bigr) \ar[l,shift left, "\lambda"] 
\end{tikzcd}
\end{equation}
commutes.
\end{enumerate}
\end{lemma}
\begin{proof}
The statements are straightforward to check.
\end{proof}

A general strategy for proving statements about adjunctions in a 2-category is to reduce the statements to adjoint functors of categories. The criterion in \cref{rslt:adjunctions-in-2-category-via-Hom-categories} is a useful tool in that regard, but it still requires the a priori construction of a candidate adjoint together with the candidate unit of the adjunction. In the following we will provide a more powerful criterion for adjunctions which allows to construct the adjoint morphism and the adjunction unit \emph{implicitly}:

\begin{proposition} \label{rslt:pointwise-criterion-for-adjunction}
Let $f\colon Y \to X$ be a morphism in a 2-category $\cat C$. Then the following are equivalent:
\begin{propenum}
	\item $f$ is a left adjoint morphism.

	\item For $Z \in \{ X, Y \}$ the right adjoint $G_Z$ of $f_*\colon \Fun_{\cat C}(Z, Y) \to \Fun_{\cat C}(Z, X)$ exists and the natural map $G_X(\id_X) \comp f \to G_Y(f)$ in $\Fun_{\cat C}(Y, Y)$ becomes an isomorphism upon applying $\Hom(\id_Y,\blank)$.

	\item For all $Z\in \cat C$ the right adjoint $G_Z$ of $f_*\colon \Fun_{\cat C}(Z,Y) \to \Fun_{\cat C}(Z,X)$ exists and for all $h\colon Z \to Z'$ in $\cat C$ the diagram
	\begin{equation}\begin{tikzcd}
		\Fun_{\cat C}(Z', Y) \arrow[d,"h^*"'] & \Fun_{\cat C}(Z', X) \arrow[l,"G_{Z'}"'] \arrow[d,"h^*"] \\
		\Fun_{\cat C}(Z, Y) & \Fun_{\cat C}(Z, X) \arrow[l,"G_Z"]
	\end{tikzcd}\end{equation}
	commutes via the natural 2-morphism, i.e.\ the natural map $h^* G_{Z'} \isoto G_Z h^*$ is an isomorphism of functors.
\end{propenum}
In this case, the right adjoint of $f$ is given by $G_X(\id_X)\colon X \to Y$.
\end{proposition}
\begin{proof}
Clearly (iii) implies (ii). We next show that (i) implies (iii), so assume that (i) is satisfied and denote by $g\colon X \to Y$ a right adjoint morphism of $f$ together with unit $\eta$ and counit $\varepsilon$. Then for every $Z \in \cat C$ the right adjoint $G_Z$ exists by \cref{rslt:adjunctions-in-2-category-via-Hom-categories} and is given by $g_*\colon \Fun_{\cat C}(Z, X) \to \Fun_{\cat C}(Z, Y)$. Moreover, for fixed $h\colon Z \to Z'$ we have natural isomorphisms $f_* h^* = h^* f_*$ and $g_* h^* = h^* g_*$ because post-composition and pre-composition commute (this follows easily from the associativity constraint in the quiver model structure for enriched categories). Now the map $h^* G_{Z'} \to G_Z h^*$ is given by the composition
\begin{align}
	h^* G_{Z'} = h^* g_* \xto{\eta h^* g_*} g_* f_* h^* g_* = g_* h^* f_* g_* \xto{g_* h^* \varepsilon} g_* h^* = G_Z h^*
\end{align}
and we need to show that it is an isomorphism. But this follows from the fact that we can commute $h^*$ with $g_*$ and $f_*$ and this is compatible with $\eta$ and $\varepsilon$, hence the claim reduces to the triangle identity for $\eta$ and $\varepsilon$.

It remains to prove that (ii) implies (i), so assume that (ii) is satisfied. In the following we write $g\coloneqq G_X(\id_X)$ to shorten the notation. The counit of the adjunction between $f_*$ and $G_X$ produces a 2-morphism $\varepsilon\colon fg = f_* G_X(\id_X) \to \id_X$. Moreover, by (ii) there is a unique 2-morphism $\eta\colon \id_Y \to gf$ such that the composite with the natural map $gf\to G_Y(f) = G_Yf_*(\id_Y)$ is induced by the unit of the adjunction between $f_*$ and $G_Y$. We claim that $(f, g, \eta, \varepsilon)$ provide the data of an adjunction, which will prove (i). To see this, we verify the commutativity of the two triangles in \cref{def:adjoint-morphisms}.

We start with the left-hand triangle. By the adjunction between $f_*$ and $G_Y$ the composition of maps $f_* \to f_* G_Y f_* \to f_*$ is the identity. By plugging in $\id_Y$ and using the map in (ii) we obtain a certain commutative diagram
\begin{equation}
\begin{tikzcd}
f \ar[d,equals] \ar[r,"f\eta"] & fgf \ar[r] \ar[d] & f \ar[d,equals] \\
f_*(\id_Y) \ar[r] & f_*G_Y(f) \ar[r] & f_*(\id_Y).
\end{tikzcd}
\end{equation}
It remains to verify that the map $fgf \to f$ is given by $\varepsilon f$. But note that the map $gf \to G_Y(f)$ in (ii) is given as the composite 
\begin{align}
	gf \xto{u(gf)} G_Yf_*(gf) = G_Y(fgf) \xto{G_Y(\varepsilon f)} G_Y(f),	
\end{align}
where $u\colon \id \to G_Yf_*$ denotes the unit map. Hence, by definition, the map $fgf \to f$ is given by the top-right circuit in the diagram
\begin{equation}\begin{tikzcd}[column sep=large]
	fgf = f_*(gf) \ar[r,"f_*u(gf)"] \ar[dr,equals] & 
	f_*G_Yf_*(gf) \ar[d,"c f_*(gf)"] \ar[r,equals] & 
	f_*G_Y(fgf) \ar[d,"c(fgf)"] \ar[r,"f_*G_Y(\varepsilon f)"] & f_*G_Y(f) \ar[d,"c(f)"] \\
	& f_*(gf) \ar[r,equals] & fgf \ar[r,"\varepsilon f"'] & f,
\end{tikzcd}\end{equation}
where $c\colon f_*G_Y\to \id$ denotes the counit map. Since this diagram 
commutes, we deduce that the map $fgf\to f$ is indeed given by $\varepsilon f$.

We now check that the right-hand triangle in \cref{def:adjoint-morphisms} commutes. We have the following diagram in $\Fun_{\cat C}(X, Y)$:
\begin{equation}\begin{tikzcd}
	g \arrow[r,"\eta g"] \arrow[dr] & gfg \arrow[r,"g\varepsilon"] \arrow[d] & g\\
	& G_Y(f) \comp g \arrow[d]\\
	& G_X(fg) \arrow[uur, bend right=10]
\end{tikzcd}\end{equation}
Here the upper vertical map is induced by the map from the map in (ii), and the right-hand diagonal map is obtained by applying $G_X$ to $\varepsilon\colon fg \to \id_X$. By definition of $\eta$ the left-hand triangle commutes. The right-hand triangle also commutes, as is observed by the compatibility of the natural transformation $g_* \to G_X$ of functors $\Fun_{\cat C}(X, X) \to \Fun_{\cat C}(X, Y)$ with the counit $\varepsilon\colon fg \to \id_X$. Altogether this implies that the composition of the upper two horizontal maps agrees with the composition $g \to G_X(fg) \to g$ along the bottom part of the diagram. But note that $G_X(fg) = (G_X f_* G_X)(\id_X)$ and hence the adjunction of $G_X$ and $f_*$ implies that this composition is the identity.
\end{proof}

The next result provides another criterion for adjunctions that comes in especially handy in the context of 6-functor formalisms:

\begin{proposition} \label{rslt:conservativity-criterion-for-adjunction}
Let $\cat C$ be a 2-category and let $(\ell_i\colon W_i \rightleftarrows Y\noloc r_i)_{i\in I}$ be a family of adjunctions in $\cat C$ such that the family of functors $r_{i*}\colon \Fun(Y, Y) \to \Fun(Y, W_i)$ is conservative. Then a morphism $f\colon Y \to X$ is left adjoint if and only if the following is true:
\begin{enumerate}[(a)]
	\item For $Z \in \{ X, Y \}$ the functor $f_*\colon \Fun_{\cat C}(Z, Y) \to \Fun_{\cat C}(Z, X)$ admits a right adjoint $G_Z$.

	\item All the compositions $f \ell_i\colon W_i \to X$ are left adjoint.
\end{enumerate}
In this case the right adjoint of $f$ is given by $G_X(\id_X)$.
\end{proposition}
\begin{proof}
The \enquote{only if} part is clear. To prove the \enquote{if} part, assume that (a) and (b) are satisfied. For each $i$ we consider the diagram
\begin{equation}\begin{tikzcd}
	\Fun_{\cat C}(X, W_i) \arrow[d,"f^*"] & \Fun_{\cat C}(X, Y) \arrow[d,"f^*"] \arrow[l,"r_{i*}"'] & \Fun_{\cat C}(X, X) \arrow[d,"f^*"] \arrow[l,"G_X"']\\
	\Fun_{\cat C}(Y, W_i) & \Fun_{\cat C}(Y, Y) \arrow[l,"r_{i*}"] & \Fun_{\cat C}(Y, X) \arrow[l,"G_Y"]
\end{tikzcd}\end{equation}
By \cref{rslt:pointwise-criterion-for-adjunction} it is enough to show that the right-hand square commutes via the natural map. By assumption the bottom left horizontal map $r_{i*}$ is conservative (if we apply it for all $i$), hence the desired commutativity can be checked after applying $r_{i*}$. But then it follows from the fact that both the left-hand square and the outer rectangle commute by \cref{rslt:pointwise-criterion-for-adjunction} and (b).
\end{proof}

So far we have developed several basic properties and criteria for adjunctions in a 2-category. Next we provide some examples by giving explicit descriptions of adjoint morphisms in the 2-categories that appeared in \cref{sec:2cat.def}:

\begin{proposition} \label{rslt:adjunctions-in-Cat-over-S}
Let
\begin{equation}\begin{tikzcd}
	\cat C \arrow[dr,"p"'] && \cat C' \arrow[dl,"q"] \arrow[ll,"G"'] \\
	& \cat S
\end{tikzcd}\end{equation}
be a commuting diagram of categories. Then the following are equivalent:
\begin{propenum}
	\item $G$ is a right adjoint morphism in $\Cat_{/\cat S}$.
	\item $G$ is \emph{relatively right adjoint over $\cat S$} in the sense of \cite[Definition~7.3.2.2]{HA}, i.e.\ $G$ admits a left adjoint $F\colon \cat C \to \cat C'$ such that for every $X \in \cat C$ the unit $X \to GF(X)$ is sent to an isomorphism under $p$.
\end{propenum}
\end{proposition}
Beware that in (ii) we require the unit to be sent under $p$ only to an isomorphism rather than the identity. Hence, the subtlety in the implication \enquote{(ii)$\implies$(i)} is to construct a new unit which lifts the identity.
\begin{proof}
It is easy to see that (i) implies (ii) by using that the forgetful functor $\Cat_{/\cat S} \to \Cat$ is a 2-functor (see \cref{rslt:2-categorical-enhancement-of-Cat-over-S}). For the other implication, the hardest part is to show that $F$ and the unit $\id \to GF$ lift from $\Cat$ to $\Cat_{/S}$. This is essentially the content of \cite[Proposition~7.3.2.1]{HA}, which we rephrase here in a less model dependent language. By the relative version of straightening (see \cref{rslt:relative-Straightening}), a functor $\cat I^{\op} \to \Cat_{/\cat S}$ from some fixed category $\cat I$ is the same as a functor $\cat E \to \cat I \times \cat S$ such that the functor $\cat E \to \cat I$ is a cartesian fibration and the functor $\cat E \to \cat S$ sends cartesian edges in $\cat E$ to isomorphisms in $\cat S$. In particular $G$ induces a functor $\cat E_G \to \Delta^1 \times \cat S$ such that $\cat C \isom \{0\}\times_{\Delta^1}\cat E_G $ and $\cat C' \isom \{1\}\times_{\Delta^1}\cat E_G$. The existence of the left adjoint $F$ is equivalent to the fact that the cartesian fibration $\cat E_G \to \Delta^1$ is also a cocartesian fibration (cf. \cite[\href{https://kerodon.net/tag/02FN}{Tag 02FN}]{kerodon}). The condition on the unit for the adjunction between $F$ and $G$ implies that the functor $\cat E_G \to \cat S$ sends cartesian edges to isomorphisms and hence by relative straightening $F$ upgrades to a functor over $\cat S$.

In order to lift the unit to $\Cat_{/\cat S}$, we consider the induced functor $\Fun(\cat C, \cat E_G) \to \Fun(\cat C, \Delta^1) \times \Fun(\cat C, \cat S)$ whose projection to $\Fun(\cat C, \Delta^1)$ is a cartesian and cocartesian fibration by \cref{rslt:exponentiation-of-cocartesian-fibrations}. By \cite[\href{https://kerodon.net/tag/02VV}{Tag 02VV}]{kerodon} the pullback along $\{ p \} \to \Fun(\cat C, \cat S)$, i.e.\ the functor
\begin{align}
\Fun_{\cat S}(\cat C, \cat E_G) \to \Fun(\cat C, \Delta^1),
\end{align}
is still a cartesian and cocartesian fibration. We now construct the unit as a morphism in $\Fun_{\cat S}(\cat C, \cat C)$. Denote $\alpha$ the morphism in $\Fun(\cat C,\Delta^1)$ corresponding the projection $\Delta^1\times \cat C \to \Delta^1$. Let $\iota_0 \colon \cat C \injto \cat E_G$ be the embedding. Then $\alpha$ admits a cocartesian lift $\iota_0 \to \alpha_!(\iota_0)$ and a cartesian lift $\alpha^*\alpha_!(\iota_0) \to \alpha_!(\iota_0)$. By the cartesian property, there exists a unique map $\eta\colon \iota_0\to \alpha^*\alpha_!(\iota_0)$ making the triangle
\begin{equation}
\begin{tikzcd}
\iota_0 \ar[d,dashed,"\exists! \eta"'] \ar[r] & \alpha_!(\iota_0) \\
\alpha^*\alpha_!(\iota_0) \ar[ur]
\end{tikzcd}
\end{equation}
commute in $\Fun_{\cat S}(\cat C, \cat E_G)$. By construction, the image of $\eta$ in $\Fun(\cat C,\Delta^1)$ is the identity of the constant functor $\const_{0}\colon \cat C\to \{0\} \injto \Delta^1$. Hence, $\eta$ corresponds to the desired unit $\id_{\cat C} \to GF$ under the identification $\Fun_{\cat S}(\cat C, \cat C) \isom \Fun_{\cat S}(\cat C, \cat E_G) \times_{\Fun(\cat C,\Delta^1)} \{\const_0\}$.

By dualizing the above argument one can also lift the adjunction counit of $F$ and $G$ to $\Cat_{/\cat S}$. In order to verify the adjunction of $F$ and $G$ in $\Cat_{/\cat S}$, by \cref{rslt:adjoints-from-weak-triangle-identity} it is enough to verify that the induced maps $F \to FGF \to F$ and $G \to GFG \to G$ are isomorphisms in $\Fun_{\cat S}(\cat C, \cat C')$ and $\Fun_{\cat S}(\cat C', \cat C)$, respectively. But the forgetful functor $\Fun_{\cat S} \to \Fun$ is conservative, hence the claim follows from the adjointness of $F$ and $G$ as morphisms in $\Cat$.
\end{proof}

\begin{proposition} \label{rslt:adjunctions-in-Enr-V}
Let $\cat V$ be an operad over $\Ass^\tensor$, let $F\colon \cat C \rightleftarrows \cat C' \noloc G$ be a diagram of $\cat V$-enriched categories and let $\eta\colon \id_{\cat C} \to GF$ be a morphism of $\cat V$-enriched functors. Then the following are equivalent:
\begin{propenum}
	\item $F$ is left adjoint to $G$ with unit $\eta$.

	\item For all $X \in \cat C$ and $Y \in \cat C'$ the induced morphism
	\begin{align}
		\enrHom[\cat V]_{\cat C'}(F(X), Y) \isoto \enrHom[\cat V]_{\cat C}(X, G(Y))
	\end{align}
	is an isomorphism in $\cat V$.
\end{propenum}
\end{proposition}
\begin{proof}
By the construction in \cref{rslt:2-categorical-enhancement-of-Enr} there is a 2-functor $\Enr_{\cat V} \to \Cat_{/\LM^\tensor}$ which is conservative on $\Hom$-categories, hence by \cref{rslt:adjoints-from-weak-triangle-identity} the adjointness of $F$ and $G$ can equivalently be checked in $\Cat_{/\LM^\tensor}$. The unit $\id \to GF$ is clearly sent to the identity under the projection $\cat C^\tensor \to \LM^\tensor$, hence by \cref{rslt:adjunctions-in-Cat-over-S} we can finally deduce that $F$ and $G$ are adjoint in $\Enr_{\cat V}$ if and only if they form an adjoint pair of functors $\cat C^\tensor \rightleftarrows \cat C'^\tensor$, where $\cat C^\tensor$ and $\cat C'^\tensor$ are the $\LM^\tensor$-operads exhibiting the $\cat V$-enrichment of $\cat C$ and $\cat C'$.

Now $F$ and $G$ are adjoint as functors on $\cat C^\tensor$ and $\cat C'^\tensor$ if and only if for all $X_\bullet \in \cat C^\tensor$ and $Y_\bullet \in \cat C'^\tensor$ the induced map $\Hom(F(X_\bullet), Y_\bullet) \isoto \Hom(X_\bullet, G(Y_\bullet))$ is an isomorphism. This can be checked in the fibers over maps in $\LM^\tensor$ and then reduces to the fibers over \emph{active} maps. We can then further reduce to the case that $Y_\bullet = Y$ lives over $\mathfrak m \in \LM^\tensor$ and $X_\bullet = V_1 \dsum \dots \dsum V_n \dsum X$ lives over $\mathfrak a^{\dsum n} \dsum \mathfrak m$ and that we consider maps over the unique active map $\alpha\colon \mathfrak a^{\dsum n} \dsum \mathfrak m \to \mathfrak m$. In other words, the adjointness of $F$ and $G$ is equivalent to showing that the induced map
\begin{align}
	\Hom_\alpha(F(V_1 \dsum \dots \dsum V_n \dsum X), Y) \isoto \Hom_\alpha(V_1 \dsum \dots \dsum V_n \dsum X, G(Y))
\end{align}
is an isomorphism. But $F(V_1 \dsum \dots \dsum V_n \dsum X) = V_1 \dsum \dots \dsum V_n \dsum F(X)$, hence by the definition if $\enrHom[\cat V]$ (see \cref{def:morphism-objects-for-enrichment}) the above isomorphism is equivalent to the isomorphism
\begin{align}
	\Mul_{\cat V}(V_1, \dots, V_n; \enrHom[\cat V]_{\cat C'}(F(X), Y)) \isoto \Mul_{\cat V}(V_1, \dots, V_n; \enrHom[\cat V]_{\cat C}(X, G(Y))).
\end{align}
This shows that (ii) implies (i). By taking $n = 1$ we deduce that (i) implies (ii).
\end{proof}

\begin{proposition} \label{rslt:adjoint-morphisms-in-functor-category}
Let $\cat C$ be a small (1-)category, $F, G\colon \cat C \to \Cat$ functors and $\alpha\colon F \to G$ a natural transformation, i.e.\ a morphism in $\cat P_2(\cat C^\op) = \TwoFun(\cat C, \Cat)$. Then $\alpha$ is a left adjoint morphism if and only if the following is true:
\begin{enumerate}[(a)]
	\item For every $X \in \cat C$ the functor $\alpha(X)\colon F(X) \to G(X)$ has a right adjoint $\beta_X$.

	\item For every morphism $f\colon Y \to X$ in $\cat C$ the diagram
	\begin{equation}\begin{tikzcd}
		F(Y) \arrow[d,"F(f)",swap] & G(Y) \arrow[l,"\beta_Y"'] \arrow[d,"G(f)"] \\
		F(X) & G(X) \arrow[l,"\beta_X"]
	\end{tikzcd}\end{equation}
	commutes via the natural 2-morphism, i.e.\ the natural map $F(f) \comp \beta_Y \isoto \beta_X \comp G(f)$ is an isomorphism of functors $G(Y) \to F(X)$.
\end{enumerate}
A similar criterion holds for right-adjointness of $\alpha$.
\end{proposition}
\begin{proof}
Note that $\TwoFun(\cat C, \Cat) = \Fun(\cat C, \Cat)$ is the category of functors $\cat C \to \Cat$ (see \cref{def:inclusion-of-Cat-into-TwoCat}) and the 2-categorical enhancement of this functor category is induced by the cartesian left action of $\Cat$ coming from the functor $A\colon \Cat \to \Fun(\cat C, \Cat)$ sending a category $\cat V$ to the constant functor with value $\cat V$ (see \cref{rslt:2-yoneda}). By straightening we have an equivalence of categories $\Fun(\cat C, \Cat) = \Cat^\cocart_{/\cat C}$. We claim that the 2-categorical enhancement on $\Fun(\cat C, \Cat)$ corresponds to the natural 2-categorical enhancement of $\Cat^\cocart_{/\cat C}$ as a 2-full subcategory of the 2-category $\Cat_{/\cat C}$ from \cref{rslt:2-categorical-enhancement-of-Cat-over-S}. Namely, the latter is induced by the cartesian left action of $\Cat$ coming from the functor $A'\colon \Cat \to \Cat^\cocart_{/\cat C}$, $\cat V \mapsto \cat C \times \cat V$. By \cref{rslt:uniqueness-of-cartesian-left-actions} it is enough to verify that the functors $A$ and $A'$ coincide under the straightening equivalence. But this follows immediately from functoriality of straightening in $\cat C$ (along the unique map $\cat C \to *$).

The functors $F$ and $G$ correspond to cocartesian fibrations $\cat E_F, \cat E_G \to \cat C$ and the natural transformation $\alpha$ corresponds to a functor $\alpha'\colon \cat E_F \to \cat E_G$ over $\cat C$ which preserves cocartesian edges. By \cref{rslt:adjunctions-in-Cat-over-S} $\alpha'$ is a left adjoint morphism in $\Cat_{/\cat C}$ if and only if it is a relative left adjoint functor over $\cat C$ in the sense of \cite[Definition~7.3.2.2]{HA}. By (the dual version of) \cite[Proposition~7.3.2.6]{HA} this is true if and only if for all $X \in \cat C$ the induced map of fibers $\alpha'_X\colon (\cat E_F)_X \to (\cat E_G)_X$ is a left adjoint functor; but this functor is exactly the functor $\alpha(X)\colon F(X) \to G(X)$. We conclude that $\alpha'$ is a left adjoint morphism in $\Cat_{/\cat C}$ if and only if condition (a) is satisfied.

Since $\Cat^\cocart_{/\cat C} \subset \Cat_{/\cat C}$ is a 2-full subcategory, it follows easily from \cref{rslt:adjoint-morphisms-are-unique} that a morphism in $\Cat^\cocart_{/\cat C}$ is a left adjoint if and only if it is a left adjoint in $\Cat_{/\cat C}$ and the associated right adjoint lies in $\Cat^\cocart_{/\cat C}$. We deduce that $\alpha'$ is a left adjoint morphism if and only if condition (a) is satisfied and the induced functor $\beta'\colon \cat E_G \to \cat E_F$ preserves cocartesian edges over $\cat C$. The latter part of the condition translates exactly to (b), as desired.
\end{proof}

\begin{remark}
\Cref{rslt:adjoint-morphisms-in-functor-category} holds more generally for a \emph{2-}category $\cat C$, but this seems to be surprisingly subtle to prove with the tools in this appendix. See \cite[Corollary~5.15]{Heine:Monadicity} for a reference.
\end{remark}

As an application of \cref{rslt:adjoint-morphisms-in-functor-category} we extract the following observation, which is of independent interest:

\begin{lemma} \label{rslt:termwise-adjoint-of-limit}
Let $\alpha\colon (\cat C_i)_i \to (\cat D_i)_i$ be a map of $I$-indexed diagrams of categories. Assume that every $\alpha_i\colon \cat C_i \to \cat D_i$ admits a right adjoint $\beta_i$ which commutes with the transition maps in the diagram via the natural transformations. Then there is an adjunction
\begin{align}
	\varprojlim \alpha\colon \varprojlim_i \cat C_i \rightleftarrows \varprojlim_i \cat D_i \noloc \varprojlim \beta,
\end{align}
where $\beta$ is assembled from the $\beta_i$'s.
\end{lemma}
\begin{proof}
We view $\alpha$ as a morphism in $\Fun(I, \Cat)$. Then by \cref{rslt:adjoint-morphisms-in-functor-category} $\alpha$ admits a right adjoint $\beta$ which is termwise given by $\beta_i$. The claim follows by applying the 2-functor $\varprojlim\colon \Fun(I, \Cat) \to \Cat$. To see that this functor is indeed a 2-functor (i.e.\ a lax $\Cat$-linear functor, see \cref{ex:lax-linear-equals-enriched}), we observe that it is right adjoint to the $\Cat$-linear functor given by precomposition with $I \to *$, so the lax linearity follows from \cref{rslt:monoidal-linear-right-adjoints}.
\end{proof}

We have studied basic properties of adjoint morphisms and discussed their shape in several 2-categories. Next up we introduce and study the subcategory of adjoint morphisms. Note that the identity morphism of every object in a 2-category is left and right adjoint to itself (e.g.\ by \cref{rslt:isomorphisms-are-adjoint}), so using \cref{rslt:adjoints-stable-under-composition} we can form the following subcategories: 

\begin{definition}
Let $\cat C$ be a 2-category. We denote by
\begin{align}
	\cat C^L, \cat C^R \subseteq \cat C
\end{align}
the 2-full subcategories spanned by all objects, but where we only allow left adjoint and right adjoint morphisms, respectively.
\end{definition}

The main goal of the next subsection is to show that for every 2-category $\cat C$ we have $\cat C^R \isom (\cat C^L)^{\co,\op}$, i.e.\ one can pass from left adjoints to right adjoints in a functorial way. Before we can come to that, let us start with some basic properties of the assignments $\cat C \mapsto \cat C^L$ and $\cat C \mapsto \cat C^R$. In the following result we implicitly talk about limits in the 2-category $\TwoCat$, which have not been introduced yet; in this case however, they agree with limits in the underlying category of $\TwoCat$ (see \cref{rslt:limits-in-linear-2-cat-coincide-with-limits-in-ul} below).

\begin{proposition} \label{rslt:functoriality-of-C-L-R}
The assignments $\cat C \mapsto \cat C^L$ and $\cat C \mapsto \cat C^R$ upgrade to 2-functors
\begin{align}
	(\blank)^L&\colon \TwoCat \to \TwoCat, \qquad \cat C \mapsto \cat C^L,\\
	(\blank)^R&\colon \TwoCat \to \TwoCat, \qquad \cat C \mapsto \cat C^R,
\end{align}
together with natural transformations $(\blank)^L \to \id$ and $(\blank)^R \to \id$ of 2-functors. Both $(\blank)^L$ and $(\blank)^R$ preserve limits and fully faithful 2-functors. More concretely, the following is true:
\begin{propenum}
	\item \label{rslt:2-functors-preserve-adjunctions} Every 2-functor $F\colon \cat C_1 \to \cat C_2$ preserves left adjoints and hence restricts to a 2-functor $F^L\colon \cat C_1^L \to \cat C_2^L$.

	\item \label{rslt:adjoint-morphism-in-limit-of-2-categories} Let $\cat C = \varprojlim_i \cat C_i$ be a limit in $\TwoCat$. Then a morphism in $\cat C$ is left adjoint if and only if all of the induced morphisms in $\cat C_i$ are left adjoint. In particular we have $\cat C^L = \varprojlim_i \cat C_i^L$.

	\item If $\cat C_1 \subseteq \cat C_2$ is a fully faithful inclusion of 2-categories, then a morphism in $\cat C_1$ is left adjoint if and only if it is left adjoint as a morphism in $\cat C_2$. In particular we have $\cat C_1^L = \cat C_1 \isect \cat C_2^L$ inside $\cat C_2$.
\end{propenum}
The same statements are true for right adjoints in place of left adjoints.
\end{proposition}
\begin{proof}
Part (i) is straightforward to check using the definition of adjoint morphisms, because all of the occurring data in that definition are preserved by a 2-functor. Part (iii) is obvious from the definitions. To make the construction $\cat C \mapsto \cat C^L$ functorial, we need to employ a version of straightening for 2-categories. In the following we use an ad-hoc version of 2-straightening that can be implemented without requiring more theory. Suppose we have some 2-category $\cat I$ and want to construct a 2-functor $\cat I \to \TwoCat$. If $\cat I^\tensor$ and $\TwoCat^\tensor$ denote the $\LM$-operads exhibiting the $\Cat$-enrichment, then the desired 2-functor is an $\LM$-operad map $\cat I^\tensor \to \TwoCat^\tensor$ together with an isomorphism of the fiber of this map over $\Ass^\tensor$ to the identity. By construction we have $\TwoCat^\tensor = (\Cat, \TwoCat, \iota)^\times$ (see \cref{def:induced-cartesian-left-action}), where $\iota\colon \Cat \injto \TwoCat$ is the inclusion. This implies
\begin{align}
	\Fun_2(\cat I, \TwoCat) = \Alg(\cat I^\tensor, \TwoCat^\tensor) \times_{\Alg(\cat I^\tensor_{\mathfrak a}, \TwoCat^\times)} \{ \iota_{\cat I} \},
\end{align}
where $\iota_{\cat I}$ is the map $\cat I^\tensor_{\mathfrak a} \isom \Cat^\times \injto \TwoCat^\times$. We can further identify operad maps into $\TwoCat^\times$ with lax cartesian structures (see \cref{ex:cartesian-monoidal-structure}) to see that there is a fully faithful embedding
\begin{align}
	\Fun_2(\cat I, \TwoCat) \subseteq \Fun(\cat I^\tensor, \TwoCat) \times_{\Fun(\cat I^\tensor_{\mathfrak a}, \TwoCat)} \{ \iota'_{\cat I} \},
\end{align}
whose essential image consists of the lax cartesian structures. Now $\TwoCat$ is a full subcategory of $\Op_{/\LM^\tensor} \times_{\Op_{/\Ass^\tensor}} \{ \Cat^\times \}$. Now $\Op_{/\Ass^\tensor} \subset \Cat_{/\Ass^\tensor}$ is a (non-full) subcategory that contains all isomorphisms, hence $\Op_{/\LM^\tensor} \times_{\Op_{/\Ass^\tensor}} \{ \Cat^\times \} = \Op_{/\LM^\tensor} \times_{\Cat_{/\Ass^\tensor}} \{ \Cat^\times \}$ (see \cref{rslt:fiber-product-over-subcategory}). By applying the same reasoning several times we deduce that
\begin{align}
	\Fun_2(\cat I, \TwoCat) &\subset (\Fun(\cat I^\tensor, \Cat_{/\LM^\tensor}) \times_{\Fun(\cat I^\tensor, \Cat_{/\Ass^\tensor})} \{ \Cat^\times \}) \\&\qquad\qquad\times_{(\Fun(\cat I_{\mathfrak a}^\tensor, \Cat_{/\LM^\tensor}) \times_{\Fun(\cat I_{\mathfrak a}^\tensor, \Cat_{/\Ass^\tensor})} \{ \Cat^\times \})} \{ \iota'_{\cat I} \}
\end{align}
is a non-full subcategory, where we identify $\Cat^\times$ with the constant functor $\cat I^\tensor \to \Op_{/\Ass^\tensor}$. By the straightening equivalence we have $\Fun(\cat I^\tensor, \Cat_{/\LM^\tensor}) = (\Cat^\cocart_{/\cat I^\tensor})_{/\cat I^\tensor \times \LM^\tensor}$, i.e.\ an object of $\Fun(\cat I^\tensor, \Cat_{/\LM^\tensor})$ can be identified with a pair $(\cat E, \rho)$, where $\cat E \to \cat I^\tensor$ is a cocartesian fibration and $\rho\colon \cat E \to \cat I^\tensor \times \LM^\tensor$ is a map of cocartesian fibrations over $\cat I^\tensor$. Under this correspondence, the functor $\Fun(\cat I^\tensor, \Cat_{/\LM^\tensor}) \to \Fun(\cat I^\tensor, \Cat_{/\Ass^\tensor})$ corresponds to taking the fiber product $\blank \times_{\LM^\tensor} \Ass^\tensor$ and the functor $\Fun(\cat I^\tensor, \Cat_{/\LM^\tensor}) \to \Fun(\cat I^\tensor_{\mathfrak a}, \Cat_{/\LM^\tensor})$ corresponds to taking the fiber product $\blank \times_{\cat I^\tensor} \cat I_{\mathfrak a}^\tensor$. Altogether we see that we can identify a 2-functor $\cat I \to \TwoCat$ with a map of cocartesian fibrations $\cat E \to \cat I^\tensor \times \LM^\tensor$ together with certain identifications of $\cat E \times_{\LM^\tensor} \Ass^\tensor$ and $\cat E \times_{\cat I^\tensor} \cat I_{\mathfrak a}^\tensor$.

We now apply the above \enquote{ad-hoc 2-straightening} to the case $\cat I = \TwoCat$. The identity 2-functor $\TwoCat \to \TwoCat$ (strictly speaking this should be replaced by the inclusion $\TwoCat \subset \wh{\TwoCat}$, where the right-hand side denotes the 2-category of large 2-categories) thus corresponds to a map of cocartesian fibrations $\cat E \to \TwoCat^\tensor \times \LM^\tensor$ over $\TwoCat^\tensor$. By going through the identifications, we see that for a 2-category $\cat C \in \TwoCat = (\TwoCat^\tensor)_{\mathfrak m}$, the fiber of $\cat E$ over $\{ \cat C \} \times \{ \mathfrak m \}$ is given by $\ul{\cat C}$, whereas the fiber over $\{ \cat C \} \times \{ \mathfrak a \}$ is given by $\Cat$. We now define $\cat E' \subset \cat E$ to be the non-full subcategory where we only allow left adjoint morphisms in $\cat C$ (in the obvious sense). Note that this requirement does not change the fibers over $\Ass^\tensor$ and $(\TwoCat^\tensor)_{\mathfrak a}$, so we can keep the induced identifications from $\cat E$. It follows from (i) that $\cat E' \to \TwoCat$ is still a cocartesian fibration and thus defines the desired 2-functor $(\blank)^L\colon \TwoCat \to \TwoCat$. In fact, from the construction we also obtain a natural transformation $(\blank)^L \to \id$ of endo-2-functors of $\TwoCat$.

It remains to prove (ii), so let $\cat C = \varprojlim_i \cat C_i$ be given as in the claim. Let us first explain how the first part of the claim implies the second part. By \cref{rslt:isomorphisms-are-adjoint} every isomorphism is a left adjoint morphism, showing that we have an identification $(\cat C)^\simeq = (\cat C^L)^\simeq$ of underlying anima. By \cref{rslt:limits-of-enriched-categories} we deduce $(\cat C^L)^\simeq = (\varprojlim_i \cat C_i^L)^\simeq$, hence to prove that $\cat C^L = \varprojlim_i \cat C_i^L$ it is enough to show that the functor $\cat C^L \to \varprojlim_i \cat C_i^L$ is fully faithful. By \cref{rslt:limits-of-enriched-categories} this amounts to showing that $\Fun^L_{\cat C}(Y, X) = \varprojlim_i \Fun^L_{\cat C_i}(Y_i, X_i)$ for all $X = (X_i)_i$ and $Y = (Y_i)_i$ in $\cat C$. Since we know that the same is true for $\Fun$ in place of $\Fun^L$, the claimed limit property reduces to showing that a morphism $f\colon Y \to X$ is left adjoint in $\cat C$ if and only if all the induced morphisms $f_i\colon Y_i \to X_i$ are left adjoint.

We now prove the first part of (ii), so let $f\colon Y \to X$ and $f_i\colon Y_i \to X_i$ be given as above. It follows from (i) that if $f$ is left adjoint then so are all $f_i$. We now prove the converse, so assume that all $f_i$ are left adjoint. Fix an object $Z = (Z_i)_i \in \cat C$. Using the quiver model of enriched categories we see that the $f_i$'s induce a compatible system of functors $(f_{i*}\colon \Fun_{\cat C_i}(Z_i, Y_i) \to \Fun_{\cat C_i}(Z_i, X_i))_i$ given by composition with $f_i$, whose limit is the functor $f_*\colon \Fun_{\cat C}(Z, Y) \to \Fun_{\cat C}(Z, X)$. By \cref{rslt:adjunctions-in-2-category-via-Hom-categories} each $f_{i*}$ has a right adjoint $g_{i*}\colon \Fun_{\cat C_i}(Z_i, X_i) \to \Fun_{\cat C_i}(Z_i, Y_i)$ given by composition with a right adjoint $g_i\colon X_i \to Y_i$ of $f_i$. From the fact that 2-functors preserve adjunctions we deduce that for every transition 2-functor $\varphi\colon \cat C_i \to \cat C_j$ in the diagram $(\cat C_i)_i$, the diagram
\begin{equation}\begin{tikzcd}
	\Fun_{\cat C_i}(Z_i, Y_i) \arrow[d,"\varphi"'] & \Fun_{\cat C_i}(Z_i, X_i) \arrow[l,"g_{i*}"'] \arrow[d,"\varphi"] \\
	\Fun_{\cat C_j}(Z_j, Y_j) & \Fun_{\cat C_j}(Z_j, Y_j) \arrow[l,"g_{j*}"]
\end{tikzcd}\end{equation}
commutes via the isomorphism $\varphi g_{i*} \isoto g_{j*} \varphi$ coming from the adjunctions. By \cref{rslt:termwise-adjoint-of-limit} the $g_{i*}$'s assemble into a right adjoint $g_*\colon \Fun_{\cat C}(Z, X) \to \Fun_{\cat C}(Z, Y)$ of $f_*$. Now take $Z = X$ and let $g \coloneqq g_*(\id_X)\colon X \to Y$. We claim that $g$ is a right adjoint of $f$, which will prove (ii).

For general $Z$, there is a natural transformation $g \comp (\blank) \to g_*$ of functors $\Fun_{\cat C}(Z, X) \to \Fun_{\cat C}(Z, Y)$, coming via adjunction from the map $fg = (f_* g_*)(\id_X) \to \id_X$ (the latter map uses $Z = X$). We claim that this natural transformation is an equivalence: This can be checked after passage to each $\cat C_i$, because the family of functors $(\Fun_{\cat C}(Z, X) \to \Fun_{\cat C_i}(Z_i, X_i))_i$ is conservative; but in $\cat C_i$ the desired isomorphism is true by construction of $g_{i*}$. Therefore we see that $g_*$ is the functor induced by composition with $g$. Taking $Z = Y$, the unit $\id \to g_* f_*$ induces a map $\eta\colon \id_Y \to (g_* f_*)(\id_Y) = gf$, which we claim to be the unit of the desired adjunction of $f$ and $g$. But this follows from \cref{rslt:adjunctions-in-2-category-via-Hom-categories}.
\end{proof}

\subsection{Passing to the adjoint morphism}

In the previous subsection we defined the notion of adjoint morphisms in a 2-category $\cat C$. By definition, if $f\colon Y \to X$ is left adjoint, then there is an associated right adjoint morphism $g\colon X \to Y$, which by \cref{rslt:adjoint-morphisms-are-unique} is unique up to isomorphism. In the following we will show that the association $f \mapsto g$ extends to a 2-functorial equivalence $\cat C^L \to (\cat C^R)^{\co,\op}$, i.e.\ passing to right adjoints can be made functorial.

Let us first explain the general idea by looking at the case where we fix two objects $X, Y \in \cat C$ and only consider morphisms $Y \to X$ and $X \to Y$. In this case we can use the following trick (cf. \cite[Construction~5.2.1.9]{HA}):

\begin{lemma} \label{rslt:pairing-of-categories}
Let $\cat C$ and $\cat D$ be categories and let $\mu\colon \cat D^\op \times \cat C \to \Ani$ be a functor such that for all $X \in \cat C$ and $Y \in \cat D^\op$, the functors
\begin{align}
	\mu(\blank, X)\colon \cat D^\op \to \Ani, \qquad \mu(Y, \blank) \colon \cat C \to \Ani
\end{align}
are representable and corepresentable, respectively. Then there is an induced pair of adjoint functors
\begin{align}
	f\colon \cat D \rightleftarrows \cat C\noloc g
\end{align}
such that
\begin{align}
	\Hom_{\cat D}(Y, g(X)) = \mu(Y, X) = \Hom_{\cat C}(f(Y), X)
\end{align}
for all $X \in \cat C$ and $Y \in \cat D$.
\end{lemma}
\begin{proof}
The functor $\mu$ induces a functor $\cat C \to \Fun(\cat D^\op, \Ani) = \PSh(\cat D)$ and by the representability assumption this functor factors over the Yoneda embedding $\cat D \injto \PSh(\cat D)$. This defines the functor $g\colon \cat C \to \cat D$. We can similarly define the functor $f\colon \cat D \to \cat C$, where we implicitly need to take opposite categories at the end. Then by construction we have an equivalence $\Hom(f(\blank), \blank) = \mu(\blank, \blank) = \Hom(\blank, g(\blank))$ of functors $\cat D^\op \times \cat C \to \Ani$. This implies that $f$ and $g$ are indeed adjoint (cf.\ \cref{def:adjoint}).
\end{proof}

\begin{example} \label{ex:internal-dual-via-pairing}
Let $\cat C$ be a closed symmetric monoidal category (symmetry is only assumed for simplicity) and consider the functor
\begin{align}
	\mu\colon \cat C^\op \times \cat C^\op \to \Ani, \qquad (P, Q) \mapsto \Hom(P \tensor Q, \one).
\end{align}
Then $\mu$ satisfies the assumptions of \cref{rslt:pairing-of-categories} and hence produces an adjunction $f\colon \cat C \rightleftarrows \cat C^{\op}\noloc g$. Here both $f$ and $g$ are given by $(\blank)^\vee \coloneqq \iHom(\blank, \one)$, i.e.\ they compute the internal dual. The unit and counit of the adjunction are given by the natural map $P \to P^{\vee\vee}$.
\end{example}

The construction in \cref{ex:internal-dual-via-pairing} can be seen as a way of passing between adjoints in the 2-category $*_{\cat C}$ with only one object $*$ and endomorphisms $\cat C$. It generalizes immediately to an arbitrary 2-category, if we fix the objects for our morphisms:

\begin{proposition} \label{rslt:passing-to-adjoints-for-fixed-objects}
Let $\cat C$ be a 2-category and $X, Y \in \cat C$ two objects. Then there is a natural equivalence
\begin{align}
	\Fun^L_{\cat C}(Y, X) = \Fun^R_{\cat C}(X, Y)^\op.
\end{align}
\end{proposition}
\begin{proof}
We consider the functor
\begin{align}
	\mu\colon \Fun^R_{\cat C}(X, Y) \times \Fun^L_{\cat C}(Y, X) \to \Ani, \qquad (g, f) \mapsto \Hom(\id_Y, gf).
\end{align}
By \cref{rslt:adjunctions-in-2-category-via-Hom-category-mates}, for every fixed left adjoint $f\colon Y \to X$, the functor $g \mapsto \Hom(\id_Y, gf)$ is corepresented by the right adjoint $f^r$ of $f$, i.e.\ $\Hom(\id_Y, gf) = \Hom(f^r, g)$. Similarly, for a fixed right adjoint morphism $g\colon X \to Y$ the functor $f \mapsto \Hom(\id_Y, gf)$ is represented by a left adjoint $g^{\ell}$ of $g$. Thus the assumptions of \cref{rslt:pairing-of-categories} are satisfied and produce an adjunction
\begin{align}
\ell\colon \Fun^R_{\cat C}(X, Y)^\op \rightleftarrows \Fun^L_{\cat C}(Y, X) \noloc r.
\end{align}
It remains to check that the unit and counit are isomorphisms. The adjunction is given by the composite of the canonical isomorphisms
\begin{align}
\Hom(g^{\ell}, f) = \Hom(\id_Y, gf) = \Hom(f^r,g).
\end{align}
Tracing the image of $\id_{g^{\ell}}$, we observe that the counit is given by the composite $g^{\ell r} \to gg^{\ell}g^{\ell r} \to g$, which by \cref{rslt:adjoint-morphisms-are-unique} is an isomorphism. Similary, the unit is an isomorphism being given by the composite $f^{r\ell} \to f^{r\ell}f^rf \to f$. 
\end{proof}

In the following we want to lift \cref{rslt:passing-to-adjoints-for-fixed-objects} to an equivalence $\cat C^L = (\cat C^R)^{\co,\op}$, i.e.\ we want the equivalence $\Fun^L_{\cat C}(Y, X) = \Fun^R_{\cat C}(X, Y)^\op$ to be functorial with respect to compositions in $\cat C$. The general strategy will be the same as before, i.e.\ we want to construct an appropriate pairing $\mu$. However, this is already non-trivial in the case that $\cat C$ has only a single object, i.e.\ corresponds to a monoidal category $\cat D$. In this case we need to construct a \emph{lax monoidal} functor $\mu\colon \cat D^\rev \times \cat D \to \Ani$, $(P, Q) \mapsto \Hom(\one, P \tensor Q)$. The challenge is that the functor $(P, Q) \mapsto P \tensor Q$ is not lax monoidal unless $\cat D$ is \emph{symmetric} monoidal. In the following we provide a solution to circumvent this problem, by explicitly constructing the unstraightening corresponding to the desired functor.

In order to apply our construction to an arbitrary 2-category, we first need to introduce yet another model for 2-categories. In the following we will implicitly work with generalized non-symmetric operads (as defined in \cite[Definition~2.9]{Heine.2023} or \cite[Definition~3.1.13]{Gepner-Haugseng.2015}), as this is the framework in which most of the relevant literature is written. It should be possible to develop a similar theory with generalized operads (in Lurie's sense), but we did not attempt this. By abuse of notation we denote by $\Ass^\tensor \coloneqq \bbDelta^\op$ the non-symmetric base operad.

We start things off with the following analog and generalization of \cite[Proposition~2.4.1.7]{HA}. In the following we say that a non-symmetric monoidal category is \emph{cartesian} if under the equivalence of non-symmetric operads and operads (see \cref{rmk:non-symmetric-vs-symmetric-operads}) it corresponds to a cartesian monoidal category in the sense of \cref{def:cartesian-monoidal-category}.

\begin{definition}
Let $\cat O^\tensor \to \Ass^\tensor$ be a generalized non-symmetric operad and $\cat C$ a category. A functor $F\colon \cat O^\tensor \to \cat C$ is called a \emph{lax cartesian structure} if it satisfies the following condition: For every $n \ge 0$ and every
\begin{align}
	X = (X_1, \dots, X_n) \in \cat O^\tensor_{[n]} = \cat O^\tensor_{[1]} \times_{\cat O^\tensor_{[0]}} \dots \times_{\cat O^\tensor_{[0]}} \cat O^\tensor_{[1]}
\end{align}
the natural maps induce an isomorphism $F(X) \isoto \prod_{i=1}^n F(X_i)$ in $\cat C$; in particular $F$ sends all of $\cat O^\tensor_{[0]}$ to a final object in $\cat C$. We denote by
\begin{align}
	\Fun^{\times,\lax}(\cat O^\tensor, \cat C) \subseteq \Fun(\cat O^\tensor, \cat C)
\end{align}
the full subcategory spanned by the lax cartesian structures.
\end{definition}

\begin{proposition} \label{rslt:generalized-ns-algebras-in-monoidal-category-equiv-lax-cartesian-str}
Let $\cat C^\tensor$ be a cartesian (non-symmetric) monoidal category with underlying category $\cat C$ and let $\cat O^\tensor \to \Ass^\tensor$ a generalized non-symmetric operad. Then composition with the natural projection $\cat C^\tensor \to \cat C$ induces an equivalence
\begin{align}
	\Alg(\cat O, \cat C) \isoto \Fun^{\times,\lax}(\cat O^\tensor, \cat C).
\end{align}
\end{proposition}
\begin{proof}
By \cref{rslt:cartesian-monoidal-structure-unique} we have $\cat C^\tensor = \cat C^\times \times_{\Comm^\tensor} \Ass^\tensor$ where $\cat C^\times$ is the cartesian symmetric monoidal category from \cref{ex:cartesian-monoidal-structure}. Here we implicitly use the map $\Ass^\tensor = \bbDelta^\op \to \Comm^\tensor = \Fin_*$ from \cite[Construction~4.1.2.9]{HA} to translate between non-symmetric operads and operads. We now follow the argument of \cite[Proposition~2.4.1.7]{HA}. First note that we have the following inclusions of full subcategories:
\begin{align}
	&\Alg(\cat O, \cat C) \subseteq \Fun_{\Ass^\tensor}(\cat O^\tensor, \cat C^\tensor) = \Fun_{\Fin_*}(\cat O^\tensor, \cat C^\times) \subseteq \Fun_{\Fin_*}(\cat O^\tensor, \widetilde{\cat C}^\times) =\\&\qquad= \Fun(\cat O^\tensor \times_{\Fin_*} \Gamma^\times, \cat C),
\end{align}
where we recall the definition of $\Gamma^\times$ and $\widetilde{\cat C}^\times$ from \cite[Notation~2.4.1.2]{HA} and \cite[Construction~2.4.1.4]{HA}, respectively. Concretely, $\Gamma^\times$ is the category of pairs $(\langle n \rangle, S)$ for $\langle n \rangle \in \Fin_*$ and $S \subseteq \langle n \rangle^\circ$ a subset.

Unwinding the definitions, we see that $\Alg(\cat O, \cat C)$ is the category of functors $F\colon \cat O^\tensor \times_{\Fin_*} \Gamma^\times \to \cat C$ satisfying the following conditions:
\begin{enumerate}[(a)]
	\item For every inert morphism $Y \to X$ in $\cat O^\tensor$ lying over a map $\alpha\colon \langle m \rangle \to \langle n \rangle$ in $\Fin_*$ and every subset $S \subseteq \langle n \rangle^\circ$, the induced map $F(Y, \alpha^{-1}(S)) \to F(X, S)$ is an isomorphism in $\cat C$.

	\item For every object $X \in \cat O^\tensor_{[n]}$ and every subset $S \subseteq \langle n \rangle^\circ$ the functor $F$ induces an isomorphism $F(X, S) \isoto \prod_{i\in S} F(X, \{ i \})$.
\end{enumerate}
Namely, condition (a) comes from the first of the above inclusions and characterizes the property that $F\colon \cat O^\tensor \to \cat C^\tensor$ preserves inert morphisms (see \cite[Proposition~2.4.1.5.(2)]{HA} for the description of inert morphisms in $\cat C^\tensor$). Condition (b) comes from the second of the above inclusions and ensures that $F$ factors over $\cat C^\times \subseteq \widetilde{\cat C}^\times$.

Now note that the projection $\Gamma^\times \to \Fin_*$ has a section given by $\langle n \rangle \mapsto (\langle n \rangle, \langle n \rangle^\circ)$, inducing a fully faithful embedding $\cat O^\tensor \injto \cat O^\tensor \times_{\Fin_*} \Gamma^\times$. We claim that right Kan extension along this embedding induces an embedding
\begin{align}
	\Fun^{\times,\lax}(\cat O^\tensor, \cat C) \injto \Fun(\cat O^\tensor \times_{\Fin_*} \Gamma^\times, \cat C)
\end{align}
To show the existence of the Kan extension, fix some $n \ge 0$ and an object $(X, S)$ with $X \in \cat O^\tensor_{[n]}$ and $S \subset \langle n \rangle^\circ$. Let $S_1, \dots, S_k \subset S$ be the collection of disjoint maximal subintervals of $S$ inside $\langle n \rangle^\circ$. For each $i$ denote $n_i \coloneqq \abs{S_i}$ and denote by $\alpha_i\colon [n] \to [n_i]$ the inert map in $\Ass^\tensor = \bbDelta^\op$ corresponding to the subinterval $S_i$. Then there is a cocartesian lift $X \to X_i$ of $\alpha_i$ in $\cat O^\tensor$ and we can consider the inclusion
\begin{align}
	\{ (X_1, S_1), \dots, (X_k, S_k) \} \subset (\cat O^\tensor)_{(X, S)/},
\end{align}
where we view $(X_i, S_i)$ as an object over $S_{i*} \in \Fin_*$. We observe that this inclusion is final. Indeed, the image of $\Ass^\tensor = \bbDelta^\op \to \Fin_*$ contains only maps $\alpha\colon \langle n \rangle \to \langle m \rangle$ such that $\alpha^{-1}(\langle m \rangle^\circ)$ is a subinterval of $\langle n \rangle^\circ$. This implies that every map $(\langle n \rangle, S) \to (\langle m \rangle, \langle m \rangle^\circ)$ in $\Gamma^\times$ whose projection to $\Fin_*$ has a lift to $\cat O^\tensor$ factors uniquely over $(S_{i*}, S_i)$ for some $i$; this proves the claimed finality.

We now see that for the right Kan extension $F'$ of a lax cartesian structure $F\colon \cat O^\tensor \to \cat C$ along the inclusion $\cat O^\tensor \injto \cat O^\tensor \times_{\Fin_*} \Gamma^\times$ we have
\begin{align}
	F'(X, S) = \prod_{i=1}^k F(X_i).
\end{align}
This shows that the right Kan extension exists (and is necessarily fully faithful). It is straightforward to check that the essential image of the Kan extension functor are exactly the functors satisfying (a) and (b) above, as desired.
\end{proof}

\begin{remark}
Our \cref{rslt:generalized-ns-algebras-in-monoidal-category-equiv-lax-cartesian-str} is also claimed in \cite[Proposition~2.42]{Heine.2023}, but the proof of \loccit{} simply refers to \cite[Proposition~2.4.1.7]{HA}. While the essential idea of the proof is indeed given by \cite[Proposition~2.4.1.7]{HA}, one has to be a bit careful with the right Kan extension at the end of the argument, which is why we present a full proof above.
\end{remark}

With \cref{rslt:generalized-ns-algebras-in-monoidal-category-equiv-lax-cartesian-str} at hand, we can now provide the promised model for 2-categories, at least with fixed anima of objects:

\begin{definition}
Let $S$ be an anima.
\begin{defenum}
 	\item Let $\Ass_S^\tensor$ be the generalized non-symmetric operad from \cite[Notation~4.1]{Heine.2023} or \cite[Definition~4.1.1]{Gepner-Haugseng.2015}. It comes equipped with a left fibration $\Ass_S^\tensor \to \Ass^\tensor = \bbDelta^\op$ whose fiber over $[n]$ is $S^{n+1}$, with transition maps given by projections and diagonals. There is an involution $\rev\colon \Ass_S^\tensor \isoto \Ass_S^\tensor$ which induces the reversion map $S^{n+1} \isoto S^{n+1}$, $(X_0, \dots, X_n) \mapsto (X_n, \dots, X_0)$ for all $n \ge 0$.

 	\item An \emph{$\Ass_S$-monoidal category} is a functor $F\colon \Ass_S^\tensor \to \Cat$ such that for all $n \ge 0$ and $(X_0, \dots, X_n) \in (\Ass_S^\tensor)_{[n]}$ the induced map
 	\begin{align}
 		F(X_0, \dots, X_n) \isoto F(X_0, X_1) \times F(X_1,X_2) \times \dots \times F(X_{n-1}, X_n) \label{eq:Segal-property-for-AssS-monoidal-categories}
 	\end{align}
 	is an isomorphism of categories; in particular $F$ sends all objects in $(\Ass_S^\tensor)_{[0]}$ to the final category $*$. We denote by
 	\begin{align}
 		\Mon_{\Ass_S}(\Cat) \subseteq \Fun(\Ass_S^\tensor, \Cat)
 	\end{align}
 	the full subcategory spanned by the $\Ass_S$-monoidal categories.

 	\item Given an $\Ass_S$-monoidal category $F\colon \Ass_S^\tensor \to \Cat$, we denote by
 	\begin{align}
 		F^\rev\colon \Ass_S^\tensor \xto{\rev} \Ass_S^\tensor \xto{F} \Cat
 	\end{align}
 	the reversed $\Ass_S$-monoidal category.
\end{defenum} 
\end{definition}

\begin{remark} \label{rmk:fibration-model-for-AssS-monoidal-categories}
By the straightening correspondence, an $\Ass_S$-monoidal category is equivalently given by a cocartesian fibration $\cat C^\tensor \to \Ass_S^\tensor$ such that for all $n \ge 0$ and $(X_0, \dots, X_n) \in S^{n+1}$ the natural map
\begin{align}
	\cat C^\tensor_{(X_0, \dots, X_n)} \isoto \cat C^\tensor_{(X_0, X_1)} \times \dots \times \cat C^\tensor_{(X_{n-1}, X_n)}
\end{align}
is an equivalence. In the following we will use both descriptions of $\Ass_S$-monoidal categories interchangeably.
\end{remark}

\begin{proposition} \label{rslt:2-categories-via-Ass-S-monoidal-categories}
Fix an anima $S$ and denote by $(\TwoCat)_S$ the category of 2-categories with underlying object anima $S$. Then there is a natural fully faithful embedding
\begin{align}
	(\TwoCat)_S \injto \Mon_{\Ass_S}(\Cat)
\end{align}
and the right-hand side identifies with $\Cat$-enriched precategories with object anima $S$. Under this equivalence, taking the opposite 2-category corresponds to taking the reversed $\Ass_S$-monoidal category.
\end{proposition}
\begin{proof}
For the first part of the claim we show that there is a natural equivalence $(\PreEnr_\Cat)_S = \Mon_{\Ass_S}(\Cat)$. Note that
\begin{align}
	(\PreEnr_\Cat)_S = \Alg(\Quiv_S(\Cat)) = \Alg(\Ass_S, \Cat),
\end{align}
where the first identification comes from the definition (see \cref{def:enriched-precategory}) and the second identification comes from the universal property of the Day convolution monoidal structure on $\Quiv_S(\Cat)$ (cf. the remark following \cite[Notation~4.42]{Heine.2023}). Since $\Cat$ is equipped with the cartesian monoidal structure, we deduce from \cref{rslt:generalized-ns-algebras-in-monoidal-category-equiv-lax-cartesian-str} that
\begin{align}
	\Alg(\Ass_S, \Cat) = \Fun^{\times,\lax}(\Ass_S^\tensor, \Cat) = \Mon_{\Ass_S}(\Cat),
\end{align}
as desired. The second part of the claim, i.e.\ the identification of opposite and reverse categories, follows easily from the definitions of both constructions (cf. the construction in \cite[Remark~4.43]{Heine.2023}).
\end{proof}

Explicitly, the above model for a 2-category $\cat C$ can be described as follows. Let $F\colon \Ass_S^\tensor \to \Cat$ be the functor corresponding to $\cat C$. Then the underlying anima of objects of $\cat C$ is $S$ and for two objects $X, Y \in S$ the category of maps from $X$ to $Y$ is given by
\begin{align}
	\Fun_{\cat C}(X, Y) = F(Y, X) \in \Cat,
\end{align}
where by $F(Y, X)$ we denote the image of $F$ on the element $(Y, X) \in S^2 = (\Ass_S^\tensor)_{[1]}$. For every object $X \in S$ there is a diagonal inclusion $X \to (X, X)$ in $\Ass_S^\tensor$, resulting in a functor $* = F(X) \to F(X, X) = \Fun_{\cat C}(X, X)$ which exhibits the identity morphism $X \to X$ in $\cat C$. Moreover, for objects $X, Y, Z \in S$ the projection $(Z, Y, X) \to (Z, X)$ in $\Ass_S^\tensor$ produces a functor
\begin{align}
F(Z, Y) \times F(Y, X) = F(Z, Y, X) \to F(Z, X),
\end{align}
the composition. 

With our new model for 2-categories at hand, we can now come to the promised construction of the isomorphism $\cat C^L = (\cat C^R)^{\co,\op}$. As explained above, the hard part of the construction is to produce the appropriate pairing $\mu$ extending the one appearing in \cref{rslt:passing-to-adjoints-for-fixed-objects}. In order to obtain this pairing, we need one auxiliary construction: the pyramid category. We take the following definition from the main text, cf. \cref{def:Sigma-n-and-Lambda-n} for a more detailed definition and explanation.

\begin{definition}
\begin{defenum}
	\item For $n \ge 0$ we denote by $\bbSigma^n \subseteq [n]\times [n]^{\op}$ the poset of pairs $(i, j)$ with $0 \le i \le j \le n$. This defines a cosimplicial category $\bbSigma^\bullet\colon \bbDelta \to \Cat$.

	\item Let $\cat C$ be a category. We define the category $\cat C^\Lambda$ as the cocartesian unstraightening of the functor
	\begin{align}
		\bbDelta^\op \to \Cat, \qquad [n] \mapsto \Fun((\bbSigma^n)^\op, \cat C).
	\end{align}
	In particular, $\cat C^\Lambda$ comes equipped with a map to $\bbDelta^\op$. We call $\cat C^\Lambda$ the \emph{category of pyramids in $\cat C$}.
\end{defenum}
\end{definition}

Let us make the pyramid category more explicit. By construction, an object in $\cat C^\Lambda$ is a pair $(n, X)$, where $n \in \ZZ_{\ge0}$ and $X$ is a functor $X\colon (\bbSigma^n)^\op \to \Cat$. The map $\cat C^\Lambda \to \bbDelta^\op$ sends this pair to $[n]$. We can picture the object $(n, X)$ as follows, in the special case $n = 2$:
\begin{equation}\begin{tikzcd}[column sep=0pt,row sep=tiny]
	&& X_{0,2}\\
	& X_{1,2} \arrow[ur] && X_{0,1} \arrow[ul]\\
	X_{2,2} \arrow[ur] && X_{1,1} \arrow[ur] \arrow[ul] && X_{0,0} \arrow[ul]
\end{tikzcd}\end{equation}
Here the above diagram is a diagram in $\cat C$. Morphisms in $\cat C^\Lambda$ are hard to draw in general, so we explain them in an example. First note that morphisms over the identity $[n] \to [n]$ in $\bbDelta^\op$ are just natural transformations of functors $(\bbSigma^n)^\op \to \cat C$, i.e.\ translate to a transformation of the above diagram to another diagram of similar shape. To explain morphisms in $\cat C^\Lambda$ which lie over more general morphisms in $\bbDelta^\op$, we look at the following example (which captures the general pattern). Suppose we are given objects $X\colon (\bbSigma^2)^\op \to \cat C$ and $Y\colon (\bbSigma^1)^\op \to \cat C$ in $\cat C^\Lambda$. A morphism $X \to Y$ in $\cat C^\Lambda$ lies over a map $\alpha\colon [2] \to [1]$ in $\bbDelta^\op$ and there are three choices of $\alpha$. We depict the corresponding diagrams for a morphism $X \to Y$ below. The outer two diagrams correspond to the inert maps $\alpha$ (i.e.\ those $\alpha$ which correspond to a subinterval embedding of $[1]$ in $[2]$) and the middle diagram corresponds to the active morphism $\alpha$. In all three cases the morphism is encoded by the dashed arrows such that the diagram commutes.
\begin{equation}
	\begin{tikzcd}[column sep=0pt,row sep=tiny]
		& X_{0,1} \arrow[dd,dashed]\\
		X_{1,1} \arrow[ur] \arrow[dd,dashed] && X_{0,0} \arrow[ul] \arrow[dd,dashed]\\
		& Y_{0,1}\\
		Y_{1,1} \arrow[ur] && Y_{0,0} \arrow[ul]\\
	\end{tikzcd}
	\qquad
	\begin{tikzcd}[column sep=0pt,row sep=tiny]
		& X_{0,2} \arrow[dd,dashed]\\
		X_{2,2} \arrow[ur] \arrow[dd,dashed] && X_{0,0} \arrow[ul] \arrow[dd,dashed]\\
		& Y_{0,1}\\
		Y_{1,1} \arrow[ur] && Y_{0,0} \arrow[ul]\\
	\end{tikzcd}
	\qquad
	\begin{tikzcd}[column sep=0pt,row sep=tiny]
		& X_{1,2} \arrow[dd,dashed]\\
		X_{2,2} \arrow[ur] \arrow[dd,dashed] && X_{1,1} \arrow[ul] \arrow[dd,dashed]\\
		& Y_{0,1}\\
		Y_{1,1} \arrow[ur] && Y_{0,0} \arrow[ul]\\
	\end{tikzcd}
\end{equation}

We will from now on use the notation $\Ass^\tensor$ for $\bbDelta^\op$. In particular, for every category $\cat C$ we have the cocartesian fibration $\cat C^\Lambda \to \Ass^\tensor$. In the special case that $\cat C = \Ass^\tensor$, this fibration has a natural section:

\begin{lemma} \label{rslt:section-of-Ass-Lambda-to-Ass}
The functor $\Ass^{\tensor,\Lambda} \to \Ass^\tensor$ has a section
\begin{align}
	s\colon \Ass^\tensor \to \Ass^{\tensor,\Lambda}, \qquad [n] \mapsto s([n]) = \big[ (i, j) \mapsto [i] \big].
\end{align}
Here the relevant transition maps are given by the obvious interval inclusions.
\end{lemma}
\begin{proof}
Both $\Ass^\tensor$ and $\Ass^{\tensor,\Lambda}$ are ordinary categories, so we can construct $s$ by hand. The description on objects is provided above. Now suppose we are given a morphism $\alpha\colon [n] \to [m]$ in $\Ass^\tensor$ (i.e.\ $\alpha$ is given by a map $[m] \to [n]$ of ordered sets which we will also denote $\alpha$ in the following). Then $\alpha$ induces a functor $(\bbSigma^m)^\op \to (\bbSigma^n)^\op$ via $(i, j) \mapsto (\alpha(i), \alpha(j))$. Now $s(\alpha)$ needs to be a natural transformation $\alpha^* s([n]) \to s([m])$, i.e.\ we need to specify natural maps $[\alpha(i)] \to [i]$ in $\Ass^\tensor$ for $i = 0, \dots, m$. But there is a canonical such map, given by the obvious map $[i] \to [\alpha(i)]$ of ordered sets (the restriction of $\alpha$).

One can check by hand that this construction is natural. Alternatively, we now provide a more conceptual definition. Note that the composite $\bbSigma^n \subseteq [n]\times [n]^{\op} \to [n] \to \bbDelta$ is natural in $[n]$, where the first map is given by projection and the second sends $i\mapsto [i]$ and $(i\le j)$ to the natural inclusion $[i]\subseteq [j]$. In other words, we have a map $(\bbSigma^{\bullet})^{\op} \to \const_{\bbDelta^\op}$ of cosimplicial categories. Consider now the composite map
\begin{align}
\const_* \xto{\id_{\bbDelta^{\op}}} \const_{\Fun(\bbDelta^{\op},\Ass^{\tensor})} \to \Fun\bigl(\const_{\bbDelta^{\op}}, \Ass^{\tensor}\bigr) \to \Fun\bigl((\bbSigma^{\bullet})^{\op}, \Ass^{\tensor}\bigr)
\end{align}
of simplicial categories $\Ass^{\tensor} = \bbDelta^{\op} \to \Cat$. Via unstraightening we obtain the desired section $s\colon \Ass^{\tensor} \to \Ass^{\tensor,\Lambda}$.
\end{proof}

For the desired construction of the relative pairing $\mu$ we need the following observation about the section $s$ from \cref{rslt:section-of-Ass-Lambda-to-Ass}:

\begin{lemma} \label{rslt:pyramid-subcategory-equiv-original}
Let $\cat C^\tensor$ be a category with a functor $p\colon \cat C^\tensor \to \Ass^\tensor$ and assume that $p$ admits all cocartesian lifts over inert morphisms in $\Ass^\tensor$. Consider the diagram
\begin{equation}\begin{tikzcd}
	\cat C'^\tensor \arrow[r,hook] \arrow[dr,dashed,"p'",swap] & \cat C^\tensor_s \arrow[r] \arrow[d] &\cat C^{\tensor,\Lambda} \arrow[d,"p^\Lambda"]\\
	& \Ass^\tensor \arrow[r,"s"] & \Ass^{\tensor,\Lambda}
\end{tikzcd}\end{equation}
where $s$ is the map from \cref{rslt:section-of-Ass-Lambda-to-Ass}, the square is cartesian and $\cat C'^\tensor$ is the full subcategory of the fiber product spanned by those pyramids where all transition maps are cocartesian over $\Ass^\tensor$. Then there is a natural equivalence $\cat C'^\tensor = \cat C^\tensor$ which identifies $p'$ with $p$.
\end{lemma}
\begin{proof}
Let us denote by $\iota\colon \Ass \injto \Ass^{\tensor,\Lambda}$ the fiber over $\{ [0] \} \injto \Ass^\tensor$, so that we have pullback squares
\begin{equation}\begin{tikzcd}
	\cat C^\tensor \arrow[r,"p"] \arrow[d,hook] & \Ass^\tensor \arrow[r] \arrow[d,hook] & {\{[0]\}} \arrow[d,hook]\\
	\cat C^{\tensor,\Lambda} \arrow[r,"p^\Lambda"] & \Ass^{\tensor,\Lambda} \arrow[r] & \Ass^\tensor
\end{tikzcd}\end{equation}
Here we note that for every category $\cat D$ the fiber of $\cat D^\Lambda$ over $\{ [0] \}$ is $\Fun((\bbSigma^0)^\op, \cat D) = \cat D$. There is a natural transformation $\iota \to s$ of functors $\Ass^\tensor \to \Ass^{\tensor,\Lambda}$ which is given by the natural maps $[n] \to s([n])$ that are induced by the obvious maps $[n] \to [i]$ for $i = 0, \dots, n$ (induced by the interval inclusions).  Let us denote by $\sigma\colon \Ass^\tensor \times [1] \to \Ass^{\tensor,\Lambda}$ the functor capturing this natural transformation and let $\cat C^\tensor_\sigma$ be the fiber product of $\sigma$ and $p^\Lambda$, i.e.\ we have a natural map $\cat C^\tensor_\sigma \to \Ass^\tensor \times [1]$ whose fiber over $\Ass^\tensor \times \{ 0 \}$ is $\cat C^\tensor$ and whose fiber over $\Ass^\tensor \times \{ 1 \}$ is $\cat C^\tensor_s$.

We denote by $\cat C'^\tensor_\sigma \subseteq \cat C^\tensor_\sigma$ the full subcategory where we take all objects in the fiber over $0 \in [1]$, but only allow the ones in $\cat C'^\tensor$ over $1 \in [1]$. We claim that $\cat C'^\tensor_\sigma$ is a cocartesian fibration over $[1]$, i.e.\ it has all cocartesian lifts of edges of the form $\{ [n] \} \times [0 \to 1]$ in $\Ass^\tensor \times [1]$. Fix $n \ge 0$ and an object $X$ in the fiber of $\cat C^{\tensor,\Lambda}$ over $[n] \in \Ass^{\tensor,\Lambda}$, viewed as a pyramid of height $0$ (i.e.\ $X \in \cat C^{\tensor,\Lambda}_{[n]} = \cat C^\tensor_{[n]}$ by the above pullback squares). We construct an object $X'$ over $s([n])$ together with a map $X \to X'$ in $\cat C^{\tensor,\Lambda}$ which lies over the map $[n] \to s([n])$. The pyramid $X'$ needs to be a functor $(\bbSigma^n)^\op \to \cat C^\tensor$ and we construct it from $X$ as follows. Denote by $\cat I^n \subseteq (\bbSigma^n)^\op$ the full subcategory of those pairs $(i,j)$ with $j = n$. We view $X$ as a functor $* \to \cat C^{\tensor,\Lambda}$ over $\Ass^{\tensor,\Lambda}$ and perform a relative left Kan extension of $X$ along the inclusion $* \injto \cat I^n$ mapping $*$ to $(n, n)$. We then perform a relative right Kan extension along the inclusion $\cat I^n \injto (\bbSigma^n)^\op$. Here, we implicitly use the functor $(\bbSigma^n)^\op \to \Ass^\tensor$ given by $s([n])$ in order to perform the relative Kan extensions. The existence of the Kan extensions follows from our assumption that $\cat C^\tensor$ has cocartesian lifts of inert morphisms (specifically, we need to lift morphisms of the form $[k] \to [k-1]$), see \cite[\href{https://kerodon.net/tag/02KL}{Tag 02KL}]{kerodon}. The resulting pyramid $X'$ looks as follows (exemplatory for $n = 2$):
\begin{equation}\begin{tikzcd}[column sep=0pt,row sep=small]
	&& X_0\\
	& X_1 \arrow[ur] && X_0 \arrow[ul,equal]\\
	X_2 \arrow[ur] && X_1 \arrow[ur] \arrow[ul,equal] && X_0 \arrow[ul,equal]
\end{tikzcd}\end{equation}
Here $X_n = X$ and the maps $X_k \to X_{k-1}$ are cocartesian lifts of the maps $[k] \to [k-1]$. We have constructed a map $X \to X'$ in $\cat C'^\tensor_\sigma$ which lies over the map $[n] \times [0 \to 1]$ in $\Ass^\tensor \times [1]$. It remains to show that the map $X \to X'$ is cocartesian.

Note that $\Ass^\tensor \times [1]$ is an ordinary category and in particular all $\Hom$-anima are static (i.e.\ sets). Therefore, the property of $X \to X'$ being cocartesian (as defined in \cref{def:cocartesian-morphism}) translates to the following claim:
\begin{itemize}
	\item[($*$)] For every $m \ge 0$, every map $\alpha\colon [n] \to [m]$ in $\Ass^\tensor$ and every $Y \in \cat C'^\tensor_\sigma$ lying over $[m] \times 1$, precomposition with $X \to X'$ induces an isomorphism of anima
	\begin{align}
		\Hom_{\alpha \times \id_1}(X', Y) \isoto \Hom_{\alpha \times [0 \to 1]}(X, Y).
	\end{align}
	Here $\Hom_{\alpha \times \id_1}$ denotes the connected component of $\Hom$ consisting of those maps that lie over $\alpha \times \id_1$ in $\Ass^\tensor \times [1]$ (and similarly on the right-hand side of the desired isomorphism).
\end{itemize}
To prove the claim, we compute both sides of the desired isomorphism, starting with the left-hand side. We first observe that the left-hand side identifies with $\Hom_{s(\alpha)}(X', Y)$, where we implicitly view $X'$ and $Y$ as objects in $\cat C^{\tensor,\Lambda}$. By definition of the pyramid category, this in turn identifies with the connected component
\begin{align}
	H(X', Y, \alpha) \subseteq \Hom_{\Fun((\bbSigma^m)^\op, \cat C^\tensor)}(\alpha^* X', Y)
\end{align}
of those natural transformations $\alpha^* X' \to Y$ which lie over the natural transformation $\alpha^* s([n]) \to s([m])$ when composed with $p\colon \cat C^\tensor \to \Ass^\tensor$. Note that $Y$ is a relative right Kan extension along $\cat I^m \injto (\bbSigma^m)^\op$. Hence by the universal property of relative Kan extensions (see \cite[\href{https://kerodon.net/tag/030F}{Tag 030F}]{kerodon}) we can equivalently replace $(\bbSigma^m)^\op$ by $\cat I^m$ and all the occuring functors by their restriction to $\cat I^m$ without changing the relevant subspace $H(X', Y, \alpha)$ above. The restriction of $\alpha^* X'$ to $\cat I^m$ is represented by the diagram $X_{\alpha(m)} \to X_{\alpha(m-1)} \to \dots \to X_{\alpha(0)}$, hence $H(X', Y, \alpha)$ is the anima of commuting diagrams of the form
\begin{equation}\begin{tikzcd}
	X_{\alpha(m)} \arrow[r,dashed] \arrow[d] & X_{\alpha(m-1)} \arrow[r] \arrow[d,dashed] & \dots \arrow[r] & X_{\alpha(0)} \arrow[d,dashed]\\
	Y_m \arrow[r] & Y_{m-1} \arrow[r] & \dots \arrow[r] & Y_0, 
\end{tikzcd}\end{equation}
lying over the fixed diagram in $\Ass^\tensor$ induced from $s(\alpha)$; here the lower row depicts $Y$ restricted to $\cat I^m$. Note that all the maps in the upper (and lower) row are cocartesian over $\Ass^\tensor$, which shows that the upper row is a relative left Kan extension from its restriction to the left-most object. By employing the universal property of relative Kan extensions again (see \cite[\href{https://kerodon.net/tag/030F}{Tag 030F}]{kerodon}), we deduce that
\begin{align}
	\Hom_{\alpha \times \id_1}(X', Y) = H(X', Y, \alpha) = \Hom_\alpha(X_n, Y_m). \label{eq:pyramid-computation-of-left-hand-Hom}
\end{align}
This finishes the computation of the left-hand side in claim ($*$). The right-hand side is very similar (but easier, because the upper row in the above diagram becomes constant $X$), proving claim ($*$).

From \cref{rslt:relative-Straightening} we deduce that $\cat C'^\tensor_\sigma$ defines a map $\cat C^\tensor \to \cat C'^\tensor$ in $\Cat_{/\Ass^\tensor}$ and it remains to see that this map is an equivalence. It is essentially surjective because clearly every object in $\cat C'^\tensor$ is constructed via relative right and left Kan extension (like $X'$). The fully faithfulness follows from the computation \eqref{eq:pyramid-computation-of-left-hand-Hom}.
\end{proof}

In order to construct the promised relative pairing $\mu$, we require a second section of $\Ass^{\tensor,\Lambda} \to \Ass^\tensor$, which has the following description:

\begin{lemma} \label{rslt:second-section-from-Ass-to-Ass-Lambda}
The functor $\Ass^{\tensor,\Lambda} \to \Ass^\tensor$ has a second section
\begin{align}
	t\colon \Ass^\tensor \to \Ass^{\tensor,\Lambda}, \qquad [n] \mapsto t([n]) = \big[ (i,j) \mapsto [2i + 1 - \delta_{ij}] \big]
\end{align}
Here all transition maps in $t([n])$ are of the following shape, for varying $i$:
\begin{enumerate}[(a)]
	\item $[2i] \to [2i]$: Identity.
	\item $[2i+1] \to [2i-1]$: The map corresponding to the inclusion $[2i-1] \subset [2i+1]$ of ordered sets that misses exactly the elements $i$ and $i+1$.
	\item $[2i] \to [2i-1]$: The map corresponding to the inclusion $[2i-1] \subset [2i]$ of ordered sets that misses exactly the element $i$.
	\item $[2i] \to [2i+1]$: The map corresponding to the map $[2i+1] \to [2i]$ of ordered sets which sends both $i$ and $i+1$ to $i$ and is one-to-one otherwise.
\end{enumerate}
Moreover, we have the following:
\begin{lemenum}
	\item There is a natural transformation $t \to s$ induced by the inclusion of the first half-intervals.
	\item \label{rslt:second-section-Ass-to-Ass-Lambda-is-symmetric} There is a natural isomorphism $t \isom \rev^\Lambda \comp t$, where $\rev^\Lambda\colon \Ass^{\tensor,\Lambda} \isoto \Ass^{\tensor,\Lambda}$ is the isomorphism induced by $\rev\colon \Ass^\tensor \isoto \Ass^\tensor$.
\end{lemenum}
\end{lemma}
\begin{proof}
As in the proof of \cref{rslt:section-of-Ass-Lambda-to-Ass} we can construct the section $t$ by hand and we are already given its description on objects. Given a morphism $\alpha\colon [n] \to [m]$ in $\Ass^\tensor$ we need to construct the morphism $t(\alpha)\colon t([n]) \to t([m])$, which amounts to constructing natural maps
\begin{align}
	[2\alpha(i) + 1 - \delta_{\alpha(i)\alpha(j)}] \to [2i + 1 - \delta_{ij}].
\end{align}
These maps amount to maps of ordered sets in the other direction. We choose them as follows: On the first $i + 1$ elements they act as $\alpha$ and on the last $i + 1$ elements they act as the reverse of $\alpha$ (in the sense of \cite[Remark~2.29]{Heine.2023}). One checks that this is indeed a natural transformation by separately looking at the first half and the second half of elements of the ordered sets appearing in the required commuting diagrams. This finishes the construction of $t$. From the construction it is easy to see that the claimed natural transformation $t \to s$ exists. Moreover, since everything in the construction in $t$ is symmetric with respect to the center of the ordered sets, one checks easily that $t \isom \rev^\Lambda \comp t$.
\end{proof}

We have finally obtained enough understanding of the pyramid category in order to produce the promised relative pairing $\mu$. Here it is:

\begin{proposition} \label{rslt:relative-adjunction-pairing-mu}
Let $S$ be an anima and $\cat C^\tensor \to \Ass_S^\tensor$ be an $\Ass_S$-monoidal category. Then there is a natural functor
\begin{align}
	\mu\colon \cat C^\tensor \times_{\Ass_S^\tensor} \cat C^{\tensor,\rev} \to \Ani
\end{align}
which has the following explicit description. Fix $X = (X_0, \dots, X_n) \in S^{n+1}$ and let $g_\bullet \in \cat C^\tensor_X$ and $f_{n-\bullet+1} \in \cat C^{\tensor,\rev}_X$; concretely, $g_i \in \cat C^\tensor_{(X_{i-1}, X_{i})}$ and $f_i \in \cat C^\tensor_{(X_{i}, X_{i-1})}$ for $i = 1, \dots, n$. Then
\begin{align}
	\mu(g_\bullet, f_\bullet) = \prod_{i=1}^n \Hom(\id_{X_{i-1}}, g_i f_i).
\end{align}
\end{proposition}
\begin{proof}
We will explicitly construct the left fibration corresponding to $\mu$. Let $\cat C^\tensor_t \coloneqq \cat C^{\tensor,\Lambda} \times_{\Ass^{\tensor,\Lambda}} \Ass^\tensor$, where the map $\Ass^\tensor \to \Ass^{\tensor,\Lambda}$ is the section $t$ from \cref{rslt:second-section-from-Ass-to-Ass-Lambda}. In particular, an object of $\cat C^\tensor_t$ is a pyramid in $\cat C^\tensor$ whose underlying pyramid in $\Ass^\tensor$ is of the form $t([n])$ for some $n \ge 0$. We let $\cat M \subseteq \cat C^\tensor_t$ be the full subcategory spanned by the pyramids $X\colon (\bbSigma^n)^\op \to \cat C^\tensor$ with the following properties:
\begin{enumerate}[(i)]
	\item For all $(i, j) \in (\bbSigma^n)^\op$ with $i > 0$, the morphism $X_{i,j} \to X_{i-1,j}$ is cocartesian over $\Ass^\tensor$.
	\item Fix $i \in \{ 0, \dots, n - 1 \}$ and let $X_{i,i}^l, X_{i,i+1}^l \in \cat C^\tensor$ denote the endpoints of the cocartesian lifts of the inert maps $[2i] \to [i]$ and $[2i+1] \to [i]$ (induced by the inclusion of the first half) starting at $X_{i,i}$ and $X_{i,i+1}$, respectively. Then the induced map $X_{i,i}^l \to X_{i,i+1}^l$ is an isomorphism. Symmetrically, we require the same for the cocartesian lifts over the inclusions of the second halves.
\end{enumerate}
We can similarly define the full subcategory $\Ass^\tensor_{\cat M} \subseteq \Ass^\tensor_{S,t} = \Ass^{\tensor,\Lambda}_S \times_{\Ass^{\tensor,\Lambda}} \Ass^\tensor$. The following two diagrams depict a general object of $\cat M$ and $\Ass^\tensor_{\cat M}$ of size 2:
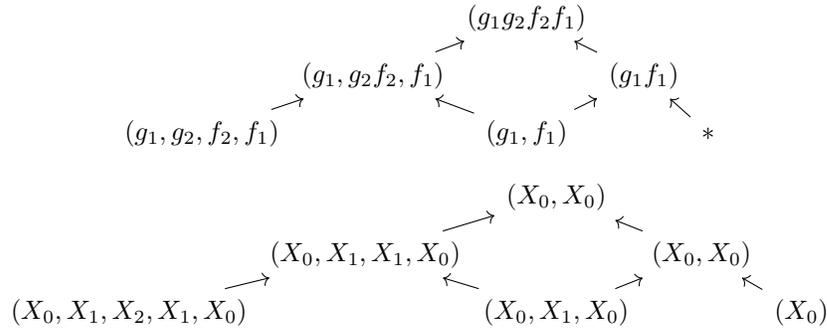
\begin{equation}\begin{tikzcd}[column sep=0pt,row sep=tiny]
	&& (g_1 g_2 f_2 f_1)\\
	& (g_1, g_2 f_2, f_1) \arrow[ur] && (g_1 f_1) \arrow[ul]\\
	(g_1, g_2, f_2, f_1) \arrow[ur] && (g_1, f_1) \arrow[ur] \arrow[ul] && * \arrow[ul]
\end{tikzcd}\end{equation}
\begin{equation}\begin{tikzcd}[column sep=0pt,row sep=tiny]
	&& (X_0, X_0)\\
	& (X_0, X_1, X_1, X_0) \arrow[ur] && (X_0, X_0) \arrow[ul]\\
	(X_0, X_1, X_2, X_1, X_0) \arrow[ur] && (X_0, X_1, X_0) \arrow[ur] \arrow[ul] && (X_0) \arrow[ul]
\end{tikzcd}\end{equation}
In the bottom diagram, all maps are the obvious ones (given by projections or diagonals). In the top diagram, maps going to the right are cocartesian and hence contain no information. There are three maps going to the left: The map $(g_1, f_1) \to (g_1, g_2f_2, f_1)$ corresponds to a map $\id_{X_1} \to g_2 f_2$, the map $* \to (g_1f_1)$ corresponds to a map $\id_{X_0} \to g_1 f_1$ and the map $(g_1f_1) \to (g_1g_2f_2f_1)$ is uniquely determined by the other maps. In summary, the top pyramid exactly encodes the choice of composable maps $f_1, f_2, g_2, g_1$ together with maps $\id_{X_1} \to g_2 f_2$ and $\id_{X_0} \to g_1 f_1$.

As in the proof of \cref{rslt:pyramid-subcategory-equiv-original}, the natural transformation $t \to s$ induces a functor $\Ass^\tensor_{\cat M} \to \Ass^\tensor_{S,s}$ whose image are exactly those pyramids where all maps are cocartesian over $\Ass^\tensor$. By \cref{rslt:pyramid-subcategory-equiv-original} this subcategory is equivalent to $\Ass^\tensor_S$, so we get a functor $\Ass^\tensor_{\cat M} \to \Ass^\tensor_S$. For example, this functor sends the above pyramid to $(X_0, X_1, X_2)$. Note that $\Ass^\tensor_{\cat M}$ is exactly the full subcategory of $\Ass^\tensor_{S,t}$ where all maps in the pyramid are cocartesian, hence by the same relative Kan extension techniques as in the proof of \cref{rslt:pyramid-subcategory-equiv-original} we see that the functor is an isomorphism, i.e.
\begin{align}
	\Ass^\tensor_{\cat M} = \Ass^\tensor_S.
\end{align}
By the same arguments as in the first half of the previous paragraph, the natural transformation $t \to s$ also induces a functor $\cat M \to \cat C^\tensor$ sending the above pyramid to $(g_1, g_2)$. By construction this functor lies over the isomorphism $\Ass^\tensor_{\cat M} \isoto \Ass^\tensor_S$ from above. Furthermore, the isomorphism $t = \rev^\Lambda \comp t$ from \cref{rslt:second-section-Ass-to-Ass-Lambda-is-symmetric} induces a natural transformation $t \to \rev^\Lambda \comp s$, which we can use as above to construct a functor $\cat M \to \cat C^{\tensor,\rev}$ sending the above pyramid to $(f_2, f_1)$ and lying over the isomorphism $\Ass^\tensor_{\cat M} \isoto \Ass^\tensor_S$. Combining the above constructions yields a functor
\begin{align}
	p\colon \cat M \to \cat C^\tensor \times_{\Ass^\tensor_S} \cat C^{\tensor,\rev}
\end{align}
which sends the above pyramid to the pair $((g_1, g_2), (f_2, f_1))$. Our goal is to show that $p$ is the left fibration which straightens to the desired functor $\mu$.

Let us first understand the fibers of $p$, so fix an object $(g_\bullet, f_\bullet) \in \cat C^\tensor \times_{\Ass^\tensor_S} \cat C^{\tensor,\rev}$, say over $[n] \in \Ass^\tensor$. Let $X_\bullet \in \Ass^\tensor_S$ be the element over which both $g_\bullet$ and $f_\bullet$ live. To compute the fiber of $p$ over $(g_\bullet, f_\bullet)$, we may first restrict all relevant categories to the fibers over $X_\bullet$, i.e.\ we look at the map $\cat M_{X_\bullet} \to \cat C^\tensor_{X_\bullet} \times \cat C^{\tensor,\rev}_{X_\bullet}$. By going through the above construction, we see that this map equivalent to the map
\begin{align}
	p_{X_\bullet}\colon \Fun'_{t,X_\bullet}((\bbSigma^n)^\op, \cat C^\tensor) \to \Fun'_{s,X_\bullet}((\bbSigma^n)^\op, \cat C^\tensor) \times \Fun'_{s,X_\bullet}((\bbSigma^n)^\op, \cat C^{\tensor,\rev}).
\end{align}
Here $\Fun_{t,X_\bullet}((\bbSigma^n)^\op, \cat C^\tensor)$ denotes the category of functors $(\bbSigma^n)^\op \to \cat C^\tensor$ over $\Ass_S^\tensor$, where the map $(\bbSigma^n)^{\op} \to \Ass_S^{\tensor}$ is the element of $\Ass_{\cat M}^{\tensor}$ corresponding to $X_{\bullet}$ under the identification $\Ass_{\cat M}^{\tensor} = \Ass_S^{\tensor}$. Moreover, $\Fun'_{t,X_\bullet} \subseteq \Fun_{t,X_\bullet}$ is the full subcategory of those functors which satisfy (i) and (ii) above; we similarly define $\Fun_{s,X_\bullet}$ and $\Fun'_{s,X_\bullet}$ (where in place of (i) and (ii) we require that all maps are cocartesian over $\Ass^\tensor$). The object $(f_\bullet, g_\bullet)$ corresponds to an object $(G_{\bullet,\bullet}, F_{\bullet,\bullet})$ on the right-hand side above. Here $F_{\bullet,\bullet}$ and $G_{\bullet,\bullet}$ are pyramids in $\cat C^\tensor$ of the following form (depicted for $n = 2$):
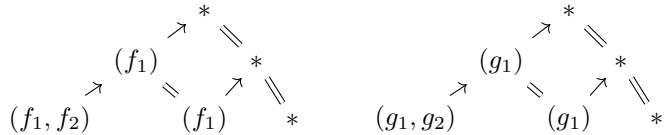
\begin{equation}
	\begin{tikzcd}[column sep=0pt,row sep=tiny]
		&& *\\
		& (f_1) \arrow[ur] && * \arrow[ul,equal]\\
		(f_1, f_2) \arrow[ur] && (f_1) \arrow[ur] \arrow[ul,equal] && * \arrow[ul,equal]
	\end{tikzcd}
	\qquad
	\begin{tikzcd}[column sep=0pt,row sep=tiny]
		&& *\\
		& (g_1) \arrow[ur] && * \arrow[ul,equal]\\
		(g_1, g_2) \arrow[ur] && (g_1) \arrow[ur] \arrow[ul,equal] && * \arrow[ul,equal]
	\end{tikzcd}
\end{equation}
In these pyramids, all maps are inert. We need to understand the fiber of $p_{X_\bullet}$ over $(G_{\bullet,\bullet}, F_{\bullet,\bullet})$. To this end, let $\cat J^n \subset \cat K^n \subset (\bbSigma^n)^\op$ be the following subcategories: $\cat K^n$ is the subcategory with all objects but only forward morphisms (i.e.\ only morphisms of the form $(i,j) \to (i',j)$ for $i' \le i \le j$) and $\cat J^n$ is the full subcategory of objects of the form $(i,i)$. By restricting functors we obtain maps
\begin{align}
	\Fun'_{t,X_\bullet}((\Sigma^n)^\op, \cat C^\tensor) \xto{q_{X_\bullet}} \Fun'_{t,X_\bullet}(\cat K^n, \cat C^\tensor) \isoto \Fun'_{t,X_\bullet}(\cat J^n, \cat C^\tensor) = \prod_{i=0}^n \cat C^\tensor_{(X_0, \dots, X_i, \dots, X_0)},
\end{align}
where the second map is an isomorphism by (i) (which says exactly that the inverse of the restriction is given by left Kan extension). A similar computation holds for the target of $p_{X_\bullet}$. We obtain a commuting diagram
\begin{equation}\begin{tikzcd}
	\Fun'_{t,X_\bullet}((\Sigma^n)^\op, \cat C^\tensor) \arrow[r,"p_{X_\bullet}"] \arrow[d,"q_{X_\bullet}"] & \Fun'_{s,X_\bullet}((\bbSigma^n)^\op, \cat C^\tensor) \times \Fun'_{s,X_\bullet}((\bbSigma^n)^\op, \cat C^{\tensor,\rev}) \arrow[d]\\
	\Fun'_{t,X_\bullet}(\cat K^n, \cat C^\tensor) \arrow[r,"\sim"] \arrow[d,equal] & \Fun'_{s,X_\bullet}(\cat K^n, \cat C^\tensor) \times \Fun'_{s,X_\bullet}(\cat K^n, \cat C^{\tensor,\rev}) \arrow[d,equal]\\
	\prod_{i=0}^n \cat C^\tensor_{(X_0, \dots, X_i, \dots, X_0)} \arrow[r,"\sim"] & \left( \prod_{i=0}^n \cat C^\tensor_{(X_0, \dots, X_i)} \right) \times \left( \prod_{i=0}^n \cat C^\tensor_{(X_i, \dots, X_0)} \right)
\end{tikzcd}\end{equation}
Here the bottom horizontal map is an isomorphism by the definition of $\Ass_S$-monoidal categories, and hence the middle horizontal map is also an isomorphism. Therefore, in order to compute the fiber of $p_{X_\bullet}$ over $(G_{\bullet,\bullet}, F_{\bullet,\bullet})$, we can first pass to the fiber over the associated object $H_{\bullet,\bullet}$ in the middle row above. When viewed as an object in $\Fun'_{t,X_\bullet}(\cat K^n, \cat C^\tensor)$ the object $H_{\bullet,\bullet}$ corresponds to the following diagram of solid arrows:
\begin{equation}\begin{tikzcd}[column sep=0pt,row sep=tiny]
	&& (g_1 g_2 f_2 f_1)\\
	& (g_1, g_2 f_2, f_1) \arrow[ur] && (g_1 f_1) \arrow[ul,dashed]\\
	(g_1, g_2, f_2, f_1) \arrow[ur] && (g_1, f_1) \arrow[ur] \arrow[ul,dashed] && * \arrow[ul,dashed]
\end{tikzcd}\end{equation}
The dashed arrows indicate the additional arrows in $(\bbSigma^n)^\op$ compared to $\cat K^n$. We observe that we can equivalently view $H_{\bullet,\bullet}$ as a collection of functors $(h_j)_{0\le j \le n}$, where $h_j$ is a functor $[j] \to \cat C^\tensor$; here we identify the ordered set $[j]$ with the associated category, which is equivalent to the full subcategory of objects $(i, j)$ (with fixed $j$ and varying $i$) in $(\bbSigma^n)^\op$. We observe
\begin{align}
	\big( \Fun((\bbSigma^n)^\op, \cat C^\tensor) \big)_{H_{\bullet,\bullet}} = \prod_{j=0}^{n-1} \Hom_{\Fun([j], \cat C^\tensor)}(h_j, h_{j+1}|_{[j]}).
\end{align}
Reintroducing the decorations $\Fun'_{t,X_\bullet}$ and using property (i) above (which says that each $h_j$ is a relative left Kan extension of its restriction to $\{ 0 \}$) together with the universal property of relative Kan extensions (see \cite[\href{https://kerodon.net/tag/030F}{Tag 030F}]{kerodon}), we deduce
\begin{align}
	&\big( \Fun'_{t,X_\bullet}((\bbSigma^n)^\op, \cat C^\tensor) \big)_{H_{\bullet,\bullet}} = \prod_{j=0}^{n-1} \Hom_{\alpha^n_j}(h_j(0), h_{j+1}(1))\\
	&\qquad\qquad= \prod_{j=0}^{n-1} \Hom_{\alpha^n_j}((g_1, \dots, g_j, f_j, \dots, f_1), (g_1, \dots, g_j, g_{j+1} f_{j+1}, f_j, \dots, f_1)),
\end{align}
where $\alpha^n_j\colon (X_0, \dots, X_j, \dots, X_0) \to (X_0, \dots, X_j, X_j, \dots, X_0)$ is the diagonal map in $\Ass^\tensor_S$. A similar description holds for the fiber of $\Fun'_{s,X_\bullet}(\cat K^n, \cat C^\tensor) \times \Fun'_{s,X_\bullet}(\cat K^n, \cat C^{\tensor,\rev})$ over $H_{\bullet,\bullet}$. Taking $(g_\bullet,f_\bullet)$ in that latter fiber and then passing to the fiber of $p_{X_\bullet}$ over this element, we finally deduce
\begin{align}
	\cat M_{(g_\bullet,f_\bullet)} &= \big(\big( \Fun'_{t,X_\bullet}((\bbSigma^n)^\op, \cat C^\tensor) \big)_{H_{\bullet,\bullet}} \big)_{(g_\bullet,f_\bullet)}\\
	&= \prod_{j=0}^{n-1} \Hom_{\beta^n_j}((g_1, \dots, g_j, f_j, \dots, f_1), (g_1, \dots, g_j, g_{j+1} f_{j+1}, f_j, \dots, f_1))\\
	&= \prod_{j=0}^{n-1} \Hom(\id_{X_j}, g_{j+1} f_{j+1}),
\end{align}
where $\beta^n_j$ is the identity map on $((g_1, \dots, g_j), (f_1, \dots, f_j))$ in $\cat C^\tensor \times \cat C^{\tensor,\rev}$. This finally shows that the fiber of $\cat M$ over $(g_\bullet, f_\bullet)$ is exactly what it needs to be.

It remains to show that $p$ is a cocartesian fibration. By the above computation the fibers of $p$ are anima, hence by \cite[\href{https://kerodon.net/tag/046V}{Tag 046V}]{kerodon} it is enough to show that $p$ is a \emph{locally} cocartesian fibration. This reduces to the following two claims:
\begin{itemize}
\item[($*$)] The projection $\cat M \to \Ass^\tensor$ is a cocartesian fibration.
\item[($**$)] The map $p\colon \cat M \to \cat C^{\tensor} \times_{\Ass_S^{\tensor}} \cat C^{\tensor, \rev}$ preserves cocartesian maps over $\Ass_S^{\tensor}$.
\item[($*{*}*$)] For every $X_\bullet \in \Ass^\tensor_S$, the fiber of $p$ over $X_\bullet$ is a locally cocartesian fibration. 
\end{itemize}
Indeed, suppose ($*$), ($**$) and ($*{*}*$) are satisfied. Then also $\cat M \to \Ass^\tensor_S$ is a cocartesian fibration: Every map in $\Ass^\tensor_S$ is cocartesian over $\Ass^{\tensor}$ (see \cref{def:left-fibration}), so that in view of ($*$) this follows from \cite[\href{https://kerodon.net/tag/01TU}{Tag 01TU}]{kerodon} and \cite[\href{https://kerodon.net/tag/01U5}{Tag 01U5}]{kerodon}. By ($**$), the same argument shows that every edge in $\cat C^\tensor \times_{\Ass^\tensor_S} \cat C^{\tensor,\rev}$ which is cocartesian over $\Ass^\tensor_S$ admits cocartesian lifts in $\cat M$. Using \cite[\href{https://kerodon.net/tag/01U4}{Tag 01U4}]{kerodon} and ($*{*}*$) we easily deduce that $p$ is a locally cocartesian fibration.


It now remains to prove claims ($*$), ($**$) and ($*{*}*$). We start with ($*$), so fix a map $\alpha\colon [n] \to [m]$ in $\Ass^\tensor$ and a pyramid $H_{\bullet,\bullet}\colon (\bbSigma^n)^\op \to \cat C^\tensor$ in $\cat M_{[n]}$. We need to construct a cocartesian map $H_{\bullet,\bullet} \to H'_{\bullet,\bullet}$ in $\cat M$ which lives over $\alpha$. This amounts to constructing a pyramid $H'_{\bullet,\bullet}$ of size $m$ and a natural transformation $\alpha^* H_{\bullet,\bullet} \to H'_{\bullet,\bullet}$ of functors $(\bbSigma^m)^\op \to \cat C^\tensor$. Note that the map $\Fun((\bbSigma^m)^{\op}, \cat C^{\tensor}) \to \Fun((\bbSigma^m)^{\op}, \Ass^{\tensor})$ is a cocartesian fibration by \cref{rslt:exponentiation-of-cocartesian-fibrations}. Then the natural transformation $h\colon \alpha^*H_{\bullet,\bullet} \to H'_{\bullet,\bullet}$ is defined as the cocartesian lift of the map $t(\alpha) \colon \alpha^*t([n]) \to t([m])$ in $\Fun((\bbSigma^m)^{\op}, \Ass^{\tensor})$. It is straightforward to check that $H'_{\bullet,\bullet}$ satisfies conditions (i) and (ii) above and hence defines an element of $\cat M$. Moreover, $h$ provides a map $H_{\bullet,\bullet} \to H'_{\bullet,\bullet}$ in $\cat M$, and similarly as in the proof of \cref{rslt:pyramid-subcategory-equiv-original} one checks that this map is cocartesian. This proves ($*$). By the construction of the cocartesian map $H_{\bullet,\bullet} \to H'_{\bullet,\bullet}$ it is clear that $p$ preserves cocartesian edges over $\Ass^{\tensor}$, which proves ($**$).


It remains to prove ($*{*}*$), so fix $X_\bullet \in \Ass^\tensor_S$. We want to show that the map $p_{X_\bullet}\colon \cat M_{X_\bullet} \to \cat C^\tensor_{X_\bullet} \times \cat C^{\tensor,\rev}_{X_\bullet}$ is a locally cocartesian fibration. A map $(\beta, \alpha)\colon (g_\bullet, f_\bullet) \to (g'_\bullet, f'_\bullet)$ is the same as a collection of maps $\beta_i\colon g_i \to g'_i$ and $\alpha_i\colon f_i \to f'_i$. Suppose $H_{\bullet,\bullet}$ and $H'_{\bullet,\bullet}$ are lifts of $(g_\bullet, f_\bullet)$ and $(g'_\bullet, f'_\bullet)$, respectively. Then similarly to the above computations one shows that
\begin{align}
	\Hom_{(\beta,\alpha)}(H_{\bullet,\bullet}, H'_{\bullet,\bullet}) = \prod_{j=0}^{n-1} R_j,
\end{align}
where $R_j$ is the anima of triangles in $\cat C^\tensor_{(X_j, X_j)}$ of the form
\begin{equation}\begin{tikzcd}
	\id_{X_j} \arrow[r] \arrow[dr] & g_{j+1}f_{j+1} \arrow[d,"\beta_{j+1}\alpha_{j+1}"]\\
	& g'_{j+1} f'_{j+1}
\end{tikzcd}\end{equation}
where all three edges are fixed (the ones from left to right are encoded by $H_{\bullet,\bullet}$ and $H'_{\bullet,\bullet}$). From this explicit description it is not hard to deduce that $p_{X_\bullet}$ is indeed a locally cocartesian fibration (using \cite[\href{https://kerodon.net/tag/01U8}{Tag 01U8}]{kerodon}).
\end{proof}

\begin{example}
Let $\cat C$ be a monoidal category. Then \cref{rslt:relative-adjunction-pairing-mu} with $S = *$ produces a lax monoidal functor
\begin{align}
	\mu\colon \cat C \times \cat C^\rev \to \Ani, \qquad (P, Q) \mapsto \Hom(\one, P \tensor Q),
\end{align}
where $\Ani$ is equipped with the cartesian monoidal structure. The lax monoidal structure on $\mu$ is encoded by natural maps
\begin{align}
	\Hom(\one, P \tensor Q) \times \Hom(\one, P' \tensor Q') \to \Hom(\one, P \tensor P' \tensor Q' \tensor Q)
\end{align}
for all $P, P', Q, Q' \in \cat C$. This natural map is obtained by sending a pair of maps $\eta\colon \one \to P \tensor Q$ and $\eta'\colon \one \to P' \tensor Q'$ to the map
\begin{align}
	\one \xto{\eta} P \tensor Q = P \tensor \one \tensor Q \xto{P \tensor \eta' \tensor Q} P \tensor P' \tensor Q' \tensor Q.
\end{align}
If $\cat C$ is \emph{symmetric} monoidal then $\mu$ can be constructed as the composition of lax monoidal functors $\cat C \times \cat C \xto{\tensor} \cat C \xto{\Hom(\one, \blank)} \Ani$. However, if $\cat C$ is not symmetric, then such a construction does not work because $\tensor$ is not a lax monoidal functor anymore. In the non-symmetric case, a more involved construction like the one above seems necessary (the assumption $S = *$ does not simplify much in that strategy). In the case $S = *$ one may also try to argue via an enriched monoidal Yoneda lemma.
\end{example}

With \cref{rslt:relative-adjunction-pairing-mu} at hand, we can finally come to the construction of \enquote{passing to the adjoint morphism}, by generalizing \cref{rslt:passing-to-adjoints-for-fixed-objects}:

\begin{theorem} \label{rslt:passing-to-adjoints}
Let $\cat C$ be a 2-category with anima of objects $\cat C^\simeq$. Then there is a natural equivalence
\begin{align}
	\cat C^R \isom \cat C^{L,\co,\op}
\end{align}
in $(\TwoCat)_{\cat C^\simeq}$ which acts on morphisms as in \cref{rslt:passing-to-adjoints-for-fixed-objects}.
\end{theorem}
\begin{proof}
We denote $S \coloneqq \cat C^\simeq$, we let $\cat C^\tensor \to \Ass^\tensor_S$ be the $\Ass_S$-monoidal category corresponding to $\cat C$ via \cref{rslt:2-categories-via-Ass-S-monoidal-categories} and we let $\mu\colon \cat C^\tensor \times_{\Ass_S^\tensor} \cat C^{\tensor,\rev} \to \Ani$ be the functor from \cref{rslt:relative-adjunction-pairing-mu}. We now perform the construction from \cref{rslt:pairing-of-categories} in the relative setting over $\Ass_S^\tensor$ using Heine's generalized Day convolution. First note that by \cref{rslt:generalized-ns-algebras-in-monoidal-category-equiv-lax-cartesian-str} we can view $\mu$ as a map of generalized non-symmetric operads $\cat C^\tensor \times_{\Ass^\tensor_S} \cat C^{\tensor,\rev} \to \Ani^\tensor$, where $\Ani^\tensor$ is the cartesian monoidal structure on $\Ani$. We now apply \cite[Theorem~11.23]{Heine.2023} with $\cat O = \Ass_S^\tensor$, $\cat O' = \Ass^\tensor$, $\cat C = \cat C^{\tensor,\rev}$ and $\cat D = \Ani^\tensor$ in order to obtain a generalized $\Ass_S$-operad
\begin{align}
	\Fun^\tensor(\cat C^\rev, \Ani) \coloneqq \Fun_{\Ass^\tensor}^{\Ass_S^\tensor}(\cat C^{\tensor,\rev}, \Ani^\tensor) \to \Ass_S^\tensor.
\end{align}
Using \cite[Remark~11.7]{Heine.2023} we can compute the fibers of $\Fun^\tensor(\cat C^\rev, \Ani)$ over objects $(X_0), (X_0, X_1) \in \Ass_S^\tensor$:
\begin{align}
	\Fun^\tensor(\cat C^\rev, \Ani)_{(X_0)} &= \Fun_{\Ass^\tensor}(\{ (X_0) \} \times_{\Ass_S^\tensor} \cat C^{\tensor,\rev}, \Ani^\tensor) = *,\\
	\Fun^\tensor(\cat C^\rev, \Ani)_{(X_0, X_1)} &= \Fun_{\Ass^\tensor}^\inert((\Ass^\tensor_{S,(X_0,X_1)/})^\inert \times_{\Ass^\tensor} \cat C^{\tensor,\rev}, \Ani^\tensor) = \Fun(\cat C^{\tensor,\rev}_{(X_0, X_1)}, \Ani).
\end{align}
Here $\Fun^\inert$ denotes the full subcategory of those functors which preserve cocartesian lifts of inert maps in $\Ass^\tensor$ and $(\Ass^\tensor_{S,(X_0,X_1)/})^\inert \subset \Ass^\tensor_{S,(X_0,X_1)/}$ is the full subcategory spanned by the inert maps. The first of the above computations follows immediately from the fact that $\Ani^\tensor_{[0]} = *$. For the second computation, we note that $(\Ass^\tensor_{S,(X_0,X_1)/})^\inert$ contains only three objects: $(X_0, X_1)$, $(X_0)$ and $(X_1)$. Moreover, the requirement on $\Fun^\inert$ forces all edges in $(\Ass^\tensor_{S,(X_0,X_1)/})^\inert$ to be sent to cocartesian edges over $\Ass^\tensor$, which says exactly that the functors are the left Kan extension of their restriction to $\{ (X_0, X_1) \}$. This implies the computation. From the definition of generalized operads (see \cite[Definition~2.9]{Heine.2023}) we deduce for all $(X_0, \dots, X_n) \in \Ass_S^\tensor$:
\begin{align}
	\Fun^\tensor(\cat C^\rev, \Ani)_{(X_0, \dots, X_n)} = \prod_{i=1}^n \Fun(\cat C^{\tensor,\rev}_{(X_{i-1}, X_i)}, \Ani) = \prod_{i=1}^n \PSh(\Fun_{\cat C}(X_{i-1}, X_{i})^{\op}).
\end{align}
By definition of $\Fun^\tensor(\cat C^\rev, \Ani)$, the map $\mu$ of generalized operads induces a map of generalized $\Ass_S$-operads
\begin{align}
	\ell\colon \cat C^\tensor \to \Fun^\tensor(\cat C^\rev, \Ani),
\end{align}
i.e.\ $\ell$ is a functor over $\Ass_S^\tensor$ which preserves inert morphisms. Concretely, in the fiber over $(X_0, \dots, X_n)$ the functor $\ell$ induces the functor
\begin{align}
	\prod_{i=1}^n \Fun_{\cat C}(X_{i}, X_{i-1}) \to \prod_{i=1}^n \PSh(\Fun_{\cat C}(X_{i-1}, X_{i})^{\op}), \qquad g_\bullet \mapsto (f_i \mapsto \Hom(\id_{X_{i-1}}, g_i f_i))_i
\end{align}
In particular, by the proof of \cref{rslt:passing-to-adjoints-for-fixed-objects}, if we restrict $\ell$ to the full subcategory $\cat C^{R,\tensor}$ of right adjoint morphisms, then the image of $\ell$ will land in the full subcategory $\cat C'^\tensor \subseteq \Fun^\tensor(\cat C^\rev, \Ani)$ whose objects are the representable presheaves. Let $\cat C^{\tensor,\rev,\vee}$ denote the dual cocartesian fibration of $\cat C^{\tensor,\rev,\vee}$, i.e.\ for all $X_\bullet \in \Ass_S^\tensor$ we have $\cat C^{\tensor,\rev,\vee}_{X_\bullet} = (\cat C^{\tensor,\rev}_{X_\bullet})^\op$. Then by \cite[\S5]{Barwick-Glasman-Nardin.2018} there is a pairing
\begin{align}
	\cat C^{\tensor,\rev,\vee} \times_{\Ass_S^\tensor} \cat C^{\tensor,\rev} \to \Ani, \qquad (g_\bullet, f_\bullet) \mapsto \prod_i \Hom(g_i, f_i).
\end{align}
This induces a fully faithful embedding
\begin{align}
	\cat C^{\tensor,\rev,\vee} \injto \Fun^\tensor(\cat C^\rev, \Ani)
\end{align}
whose image consists exactly of the representable presheaves. Altogether we see that $\mu$ induces a map
\begin{align}
	\ell\colon \cat C^{R,\tensor} \to \cat C^{L,\tensor,\rev,\vee}
\end{align}
of $\Ass_S$-monoidal categories which acts on objects by sending a right adjoint morphism to its corresponding left adjoint. By \cref{rslt:passing-to-adjoints-for-fixed-objects} $\ell$ is fiberwise an equivalence. One checks easily that $\ell$ preserves all cocartesian edges: This amounts to showing that passing to left adjoint morphisms is compatible with composition via the natural map induced by $\ell$, which was shown in \cref{rslt:adjoints-stable-under-composition}. We can altogether deduce that $\ell$ is an equivalence of $\Ass_S$-monoidal categories. But observe that $\cat C^{L,\tensor,\rev,\vee}$ is the $\Ass_S$-monoidal category corresponding to the 2-category $\cat C^{L,\op,\co}$, hence by \cref{rslt:2-categories-via-Ass-S-monoidal-categories} we obtain the desired equivalence of 2-categories.
\end{proof}

\begin{examples}
\begin{exampleenum}
	\item We have $\Cat^L \isom (\Cat^R)^{\co,\op}$. This also follows immediately from \cite[Theorem~B]{Haugseng.2023} for $B = *$.

	\item Let $\Pres \subseteq \Cat$ denote the full sub-2-category spanned by the presentable categories. Then $\PrL \isom (\PrR)^{\co,\op}$. The underlying 1-categorical equivalence is well-known and was proved in \cite[Corollary~5.5.3.4]{HTT}.
\end{exampleenum}
\end{examples}

\subsection{Limits and adjunctions}

In \cref{rslt:functoriality-of-C-L-R} we showed that adjoint morphisms in a limit of 2-categories can be detected componentwise. In the following we study the related question of how adjoint morphisms behave with respect to limits inside a fixed 2-category. To keep things simple, we work with a very basic notion of limits and colimits in a 2-category:

\begin{definition} \label{def:limits-in-a-2-category}
Let $\cat C$ be a 2-category and $I$ a category. We say that a 2-functor $F\colon I^\triangleleft \to \cat C$ is a \emph{limit diagram} if for all objects $Z \in \cat C$ the induced diagram
\begin{align}
	\Fun_{\cat C}(Z, F(\blank))\colon I^\triangleleft \to \Cat
\end{align}
is a limit diagram. In this case we denote $\varprojlim_{i\in I} F(i) \coloneqq F(\emptyset)$. We define colimit diagrams in $\cat C$ dually.
\end{definition}

More concretely, given a diagram $(X_i)_{i\in I}$ of objects in the 2-category $\cat C$ then its limit $X = \varprojlim_i X_i$ is an object $X \in \cat C$ such that for all objects $Z \in \cat C$ we have
\begin{align}
	\Fun_{\cat C}(Z, X) = \varprojlim_i \Fun_{\cat C}(Z, X_i).
\end{align}
By passing to the underlying anima we deduce that $X = \varprojlim_i X_i$ also holds in $\ul{\cat C}$ and hence $X$ is uniquely determined by that property. In other words, being a limit in $\cat C$ is a \emph{property} on a limit in $\ul{\cat C}$. The same discussion applies to colimits.

\begin{remark}
The reader who is familiar with weighted limits and colimits in the context of enriched categories may view \cref{def:limits-in-a-2-category} as the special case of a weighted limit where the weight $I \to \Cat$ is given by the constant functor $i \mapsto *$.
\end{remark}

\begin{remark}
If $\cat C$ is a 2-category, then limits in $\cat C$ and $\cat C^\co$ agree, while limits in $\cat C^\op$ correspond to colimits in $\cat C$.
\end{remark}

The above discussion indicates that limits and colimits in a 2-category $\cat C$ are closely related to those in the underlying category $\ul{\cat C}$. In fact, in many cases these two notions coincide:

\begin{lemma}
Let $\cat C$ be a 2-category such that the 2-categorical structure is induced by a $\Cat$-action on $\cat C$, i.e.\ the associated $\LM$-operad is cocartesian over $\LM^\tensor$. Then:
\begin{lemenum}
	\item \label{rslt:limits-in-linear-2-cat-coincide-with-limits-in-ul} Limits in $\cat C$ and $\ul{\cat C}$ coincide.
	\item \label{rslt:colimits-in-linear-2-cat-coincide-with-colimits-in-ul} Suppose that the action map $\tensor\colon \Cat \times \ul{\cat C} \to \ul{\cat C}$ preserves colimits in the second argument. Then colimits in $\cat C$ and $\ul{\cat C}$ coincide.
\end{lemenum}
\end{lemma}
\begin{proof}
Suppose $X = \varprojlim_i X_i$ is a limit in $\ul{\cat C}$; we need to show that it is also a limit in $\cat C$. Fix any $Z \in \cat C$ and some test category $V \in \Cat$. Then we need to verify that the natural map
\begin{align}
	\Hom(V, \Fun_{\cat C}(Z, X)) \isoto \varprojlim_i \Hom(V, \Fun_{\cat C}(Z, X_i))
\end{align}
is an isomorphism. But we can rewrite this map as
\begin{align}
	\Hom_{\ul{\cat C}}(V \tensor Z, X) \to \varprojlim_i \Hom_{\ul{\cat C}}(V \tensor Z, X_i),
\end{align}
which is clearly an isomorphism. This proves (i). To prove (ii), let now $Y = \varinjlim_j Y_j$ be a colimit in $\ul{\cat C}$. For $Z$ and $V$ as above we need to verify that the natural map
\begin{align}
	\Hom(V, \Fun_{\cat C}(Y, Z)) \isoto \varprojlim_j \Hom(V, \Fun_{\cat C}(Y_j, Z))
\end{align}
is an isomorphism. This map can be rewritten as
\begin{align}
	\Hom_{\ul{\cat C}}(V \tensor Y, Z) \to \varprojlim_j \Hom_{\ul{\cat C}}(V \tensor Y_j, Z),
\end{align}
hence the claim follows immediately from $V \tensor Y = \varinjlim_j (V \tensor Y_j)$, which is true by assumption.
\end{proof}

\begin{example}
By \cref{rslt:limits-in-linear-2-cat-coincide-with-limits-in-ul} limits in the 2-categories $\Cat$, $\Enr_{\cat V}$, $\Cat_{/S}$ and $\TwoCat$ are the same as the ones in the underlying 1-categories. A similar statement is probably true for colimits, but we have not studied these colimits well enough in order to verify the criterion in \cref{rslt:colimits-in-linear-2-cat-coincide-with-colimits-in-ul} (but for $\Cat$ it is true).
\end{example}

In the following we will frequently employ the pointwise criterion for adjunctions (see \cref{rslt:pointwise-criterion-for-adjunction}), for which the following property of a 2-category $\cat C$ comes in handy:

\begin{definition}
Let $\cat C$ be a 2-category and $I$ a category. We say that $\cat C$ is \emph{compatible with $I$-indexed limits} if the following conditions are satisfied:
\begin{enumerate}[(i)]
	\item For all objects $X, Z \in \cat C$, the category $\Fun_{\cat C}(Z, X)$ has all $I$-indexed limits.

	\item For all objects $X, Z, Z' \in \cat C$ and morphisms $Z \to Z'$ the induced functor $\Fun_{\cat C}(Z', X) \to \Fun_{\cat C}(Z, X)$ preserves $I$-indexed limits.
\end{enumerate}
We similarly say that $\cat C$ is \emph{compatible with $I$-indexed colimits} if the above properties are satisfied for $I$-indexed colimits in place of $I$-indexed limits.
\end{definition}

We now study how adjoint morphisms behave with respect to limits in a 2-category. We start with the most basic 2-category: $\Cat$. Here we have the following result:

\begin{lemma}
Let $(\cat D_i)_{i\in I}$ be a diagram of categories, $\cat C$ another category and $F\colon \cat C \to \varprojlim_i \cat D_i$ a functor with induced functors $F_i\colon \cat C \to \cat D_i$.
\begin{lemenum}
	\item Suppose that each $F_i$ has a left adjoint $L_i\colon \cat D_i \to \cat C$ and for all objects $(X_i)_i \in \varprojlim_i \cat D_i$ the colimit $\varinjlim_{i\in I^\op} L_i(X_i)$ exists in $\cat C$. Then $F$ has the left adjoint
	\begin{align}
		L\colon \varprojlim_{i\in I} \cat D_i \to \cat C, \qquad (X_i)_i \mapsto \varinjlim_{i\in I^\op} L_i(X_i).
	\end{align}

	\item \label{rslt:limits-and-adjoints-in-Cat} Suppose that each $F_i$ has a right adjoint $R_i\colon \cat D_i \to \cat C$ and for all objects $(X_i)_i \in \varprojlim_i \cat D_i$ the limit $\varprojlim_{i\in I} R_i(X_i)$ exists in $\cat C$. Then $F$ has the right adjoint
	\begin{align}
		R\colon \varprojlim_{i\in I} \cat D_i \to \cat C, \qquad (X_i)_i \mapsto \varprojlim_{i\in I} R_i(X_i).
	\end{align} 
\end{lemenum}
\end{lemma}
\begin{proof}
It is enough to prove (ii), then (i) follows by passing to opposite categories. So assume now that the $R_i$'s and the limits $\varprojlim_i R_i(X_i)$ exist. Let $p\colon \cat E \to I$ be the cocartesian fibration corresponding to the diagram $I \to \Cat$, $i \mapsto \cat D_i$. By \cite[\href{https://kerodon.net/tag/02TK}{Tag 02TK}]{kerodon} we have $\varprojlim_i \cat D_i = \Fun_I^\cocart(I, \cat E)$, i.e.\ $\varprojlim_i \cat D_i$ identifies with the category of cocartesian sections $I \to \cat E$ of $p$. Now $F$ provides a natural transformation from the constant functor $I \to \Cat$, $i \mapsto \cat C$ to the diagram $(\cat D_i)_i$, which translates into a map of cocartesian fibrations $F'\colon \cat C \times I \to \cat E$. We can write $F$ as the composition
\begin{align}
	F\colon \cat C \xto{F_1} \Fun(I, \cat C) = \Fun_I(I, \cat C \times I) \xto{F_2} \Fun_I(I, \cat E),
\end{align}
where $F_1$ sends an object $X \in \cat C$ to the constant diagram $X$, and $F_2$ is the given by composition with $F'$; we note that $F_2 \comp F_1$ lands in the full subcategory $\Fun_I^\cocart(I, \cat E) \subseteq \Fun_I(I, \cat E)$.

By (the dual version of) \cite[Proposition~7.3.2.6]{HA} the right adjoints $R_i$ assemble to a right adjoint $R'\colon \cat E \to \cat C \times I$ of $F'$ relative to $I$. Composition with $R'$ yields the functor
\begin{align}
	R_2\colon \Fun_I(I, \cat E) \to \Fun_I(I, \cat C \times I) = \Fun(I, \cat C), \qquad (X_i)_i \mapsto (R_i(X_i))_i.
\end{align}
This functor is right adjoint to $F_2$. Moreover, on the image of $\varprojlim_i \cat D_i \subseteq \Fun_I(I, \cat E)$ under $R_2$ the functor $F_1$ has a local right adjoint given by $\varprojlim_i$ (here we use that $\varprojlim_i R_i(X_i)$ exists by assumption). By the local existence criterion for adjoint functors (see \cite[\href{https://kerodon.net/tag/02FV}{Tag 02FV}]{kerodon}) we obtain the desired right adjoint $R$ of $F$.
\end{proof}

Having understood adjunctions and limits in $\Cat$, we now come to general 2-categories. By the usual Yoneda shenanigans we obtain the following result:

\begin{proposition} \label{rslt:limits-and-adjunctions}
Let $X = \varprojlim_i X_i$ be a limit in a 2-category $\cat C$, let $Y \in \cat C$ be another object and $f\colon Y \to \varprojlim_i X_i$ a morphism with induced morphisms $f_i\colon Y \to X_i$. Assume that:
\begin{enumerate}[(a)]
\item Every $f_i$ has a right adjoint $g_i\colon X_i \to Y$.
\item For $Z \in \{ X, Y \}$ and every morphism $h = (h_i)_i\colon Z \to \varprojlim_i X_i$ the limit $\varprojlim_i g_i h_i$ exists in $\Fun_{\cat C}(Z, Y)$. Moreover, for $Z=X$ this limit is preserved by precomposition along $f\colon Y \to X$.
\end{enumerate}
Then $f$ has a right adjoint $g\colon X \to Y$ given by $g = \varprojlim_i g_i \pi_i$, where $\pi_i\colon X \to X_i$ are the projections.
\end{proposition}
\begin{proof}
We employ the pointwise criterion for adjunctions from \cref{rslt:pointwise-criterion-for-adjunction}. To check condition (a) of \loccit{}, let $Z \in \{ X, Y \}$ be given. We need to see that the functor
\begin{align}
	f_*\colon \Fun_{\cat C}(Z, Y) \to \Fun_{\cat C}(Z, \varprojlim_i X_i) = \varprojlim_i \Fun_{\cat C}(Z, X_i)
\end{align}
has a right adjoint. By \cref{rslt:limits-and-adjoints-in-Cat} this follows because each $f_{i*}$ has a right adjoint and the relevant limits exist in $\Fun_{\cat C}(Z, Y)$ by assumption (b). Condition (b) of \cref{rslt:pointwise-criterion-for-adjunction} follows easily from the fact that precomposition along $f$ preserves the relevant limits by (b).
\end{proof}

\begin{corollary} \label{rslt:limits-in-C-L}
Let $\cat C$ be a 2-category and $I$ a category and assume that $\cat C$ is compatible with $I$-indexed limits. Let $X = \varprojlim_{i\in I} X_i$ be a limit in $\cat C$ and $Y \in \cat C$ an object. Then a morphism $f\colon Y \to \varprojlim_i X_i$ is left adjoint as soon as all the induced morphisms $f_i\colon Y \to X_i$ are left adjoint.
\end{corollary}
\begin{proof}
This follows immediately from \cref{rslt:limits-and-adjunctions}.
\end{proof}

\begin{remark} \label{rmk:symmetry-of-limits-and-adjunctions}
By applying \cref{rslt:limits-and-adjunctions,rslt:limits-in-C-L} to $\cat C^\op$, $\cat C^\co$ and $\cat C^{\co,\op}$ one obtains variants of these results for colimits and right adjoint morphisms.
\end{remark}

The above results provide a good understanding of how adjunctions interact with limits of objects in a 2-category. In the following we provide a related result that is concerned with limits of \emph{morphisms}. Again we start with the version in $\Cat$:

\begin{lemma} \label{rslt:adjunction-for-limit-of-morphisms-in-Cat}
Let $\cat C$, $\cat D$ and $I$ be categories and let $(f_i)_{i\in I}$ be a diagram in $\Fun(\cat C, \cat D)$. Suppose that every $f_i$ has a right adjoint $g_i\colon \cat D \to \cat C$. Assume furthermore that one of the following conditions is satisfied:
\begin{enumerate}[(a)]
\item the colimit $\varinjlim_{i\in I} f_i$ exists in $\Fun(\cat C,\cat D)$, or
\item the limit $\varprojlim_{i\in I^\op} g_i$ exists in $\Fun(\cat D,\cat C)$.
\end{enumerate}
Then there is an adjunction
\begin{align}
\varinjlim_{i\in I} f_i \colon \cat C \rightleftarrows \cat D \noloc \varprojlim_{i\in I^\op} g_i.
\end{align}
\end{lemma}
\begin{proof}
This is immediate from \cref{rslt:right-adjoints-closed-under-limits}.
\end{proof}

By combining \cref{rslt:adjunction-for-limit-of-morphisms-in-Cat} with the pointwise criterion for adjunctions, we obtain the following version for a general 2-category:

\begin{proposition} \label{rslt:adjunction-for-limit-of-morphisms}
Let $X$ and $Y$ be objects in a 2-category $\cat C$ and let $(f_i)_{i\in I}$ be a diagram in $\Fun_{\cat C}(Y, X)$. Assume that the following conditions are satisfied:
\begin{enumerate}[(a)]
\item Every $f_i$ has a right adjoint $g_i\colon X\to Y$.

\item $\varinjlim_{i\in I} f_i \in \Fun_{\cat C}(Y,X)$ exists, and for all $Z\in \{X,Y\}$ and all morphisms $h\colon Z\to Y$ the natural map $\varinjlim_{i} (f_{i}\comp h) \isoto (\varinjlim_{i}f_{i})\comp h$ is an isomorphism.

\item For all $Z\in \{X,Y\}$ and all $h\colon Z\to X$ the natural maps $(\varprojlim_{i\in I^{\op}}g_{i*})(h) \isoto \varprojlim_{i}g_{i*}(h)$ and $(\varprojlim_ig_i)\comp h \isoto \varprojlim_i (g_i\comp h)$ are isomorphisms (provided the left-hand limits exist).
\end{enumerate}
%
%
Then there is an adjunction
\begin{align}
\varinjlim_{i\in I} f_i\colon Y \rightleftarrows X \noloc \varprojlim_{i\in I^\op} g_i
\end{align}
in $\cat C$.
\end{proposition}
\begin{proof}
Fix some $Z \in \{X,Y\}$ and write $f\coloneqq \varinjlim_i f_i$. Then the morphisms $f_i,f\colon Y\to X$ induce functors
\begin{align}
f_{i*}, f_*\colon \Fun_{\cat C}(Z,Y) \to \Fun_{\cat C}(Z,X).
\end{align}
From condition (b) we easily deduce that the natural map $\varinjlim_i f_{i*} \isoto f_*$ is an isomorphism (since $\varinjlim_if_i$ exists, this can be checked pointwise). Moreover, by condition (a) each $f_{i*}$ has a right adjoint given by $g_{i*}\colon \Fun_{\cat C}(Z,X) \to \Fun_{\cat C}(Z,Y)$. By \cref{rslt:adjunction-for-limit-of-morphisms-in-Cat} the limit $G_Z\coloneqq \varprojlim_{i\in I^\op} g_{i*}$ exists in $\Fun(\Fun_{\cat C}(Z,X), \Fun_{\cat C}(Z,Y))$ and is a right adjoint of $f_*$. By condition (c) also $g\coloneqq \varprojlim_{i\in I^{\op}}g_i = G_X(\id_X)$ exists in $\Fun_{\cat C}(X,Y)$, and the natural map $g\comp f\isoto G_Y(f)$ is an isomorphism. Now apply \cref{rslt:pointwise-criterion-for-adjunction}.
\end{proof}
%

\begin{corollary} \label{rslt:Fun-L-stable-under-colim-if-compatible-2-cat}
Let $\cat C$ be a 2-category and $I$ a category such that $\cat C$ is compatible with colimits over $I$ and limits over $I^\op$. Then for all $X, Y \in \cat C$ the category $\Fun^L(Y, X)$ has all colimits over $I$ and the inclusion $\Fun^L(Y, X) \subseteq \Fun(Y, X)$ preserves them.
\end{corollary}
\begin{proof}
This reduces easily to \cref{rslt:adjunction-for-limit-of-morphisms}.
\end{proof}

\begin{corollary} \label{rslt:Fun-L-stable-under-retracts}
Let $\cat C$ be a 2-category. Then for all $X,Y \in \cat C$ the inclusion $\Fun^L(Y,X) \subseteq \Fun(Y,X)$ is closed under retracts.
\end{corollary}
\begin{proof}
This follows from \cref{rslt:adjunction-for-limit-of-morphisms} together with the characterization of retracts in terms of limits and colimits (see \cite[\href{https://kerodon.net/tag/0403}{Tag 0403}]{kerodon} and \cite[\href{https://kerodon.net/tag/0404}{Tag 0404}]{kerodon}).
\end{proof}

\begin{remark}
Similarly to \cref{rmk:symmetry-of-limits-and-adjunctions} one can deduce analogous versions of the above results by applying them to $\cat C^\op$, $\cat C^\co$ and $\cat C^{\co,\op}$.
\end{remark}

\end{document}